\documentclass[reqno,11pt,openany]{amsbook}
\usepackage{amsfonts,amsmath,amssymb}
\usepackage{graphicx}
\usepackage[top=3.2cm, bottom=3.3cm]{geometry}
\numberwithin{section}{chapter}
\numberwithin{equation}{section}

\def\overg{\overline{g}}

\def\a{{\alpha}}

\def\g{\mathbf{g}}
\def\be{{\beta}}
\def\ga{\gamma}

\def\eps{\epsilon}

\def\la{\lambda}

\def\om{\omega}

\def\ze{\zeta}

\def\varep{\varepsilon}

\def\al{\alpha}

\def\HH{{\mathcal H}}
\def\LL{{\mathcal L}}

\def\uF{{\underline{F}}}

\def\D{{\bf D}}

\def\R{{\bf R}}

\def\g{{\bf g}}

\def\f12{{\frac 1 2}}

\def\phii{\widetilde{\varphi}}

\def\va{\vartheta}

\def\f{\widetilde{f}}

\newcommand{\bea}{\begin{eqnarray}}
\newcommand{\eea}{\end{eqnarray}}
\def\beaa{\begin{eqnarray*}}
\def\eeaa{\end{eqnarray*}}

\def\ba{\begin{array}}
\def\ea{\end{array}}

\newtheorem{theorem}{Theorem}[chapter]
\newtheorem{lemma}[theorem]{Lemma}
\newtheorem{proposition}[theorem]{Proposition}
\newtheorem{corollary}[theorem]{Corollary}
\newtheorem{definition}[theorem]{Definition}
\newtheorem{remark}[theorem]{Remark}

\usepackage{etoolbox}
\makeatletter
\patchcmd{\@maketitle}{\newpage}{}{}{} 
\makeatother

\begin{document}
\title[The Einstein-Klein-Gordon coupled system]{The Einstein-Klein-Gordon coupled system: global stability of the Minkowski solution}

\author{Alexandru D. Ionescu}
\address{Princeton University}
\email{aionescu@math.princeton.edu}

\author{Beno\^{i}t Pausader}
\address{Brown University}
\email{benoit\_pausader@brown.edu}

\maketitle

\setcounter{tocdepth}{1}

\tableofcontents

\chapter{Introduction}\label{intro}

\section{The Einstein-Klein-Gordon coupled system}\label{firstSec}

The Einstein field equations  of General Relativity are a covariant geometric system that connect the Ricci tensor of a Lorentzian metric $\g$ on a manifold $M$ to the energy-momentum tensor of the matter fields in the spacetime, according to the equation
\begin{equation}\label{EFE}
G_{\al\be}=8\pi T_{\al\be}.
\end{equation}
Here $G_{\al\be}=\R_{\al\be}-(1/2)\R\g_{\al\be}$ is the Einstein tensor, where $\R_{\al\be}$ is the Ricci tensor, $\R$ is the scalar curvature, and $T_{\al\be}$ is the energy-momentum tensor of the matter in the space.

In this monograph we are concerned with the Einstein-Klein-Gordon coupled system, which describes the coupled evolution of an unknown Lorentzian metric $\g$ and a massive scalar field $\psi$. In this case the associated energy momentum tensor $T_{\al\be}$ is given by
\begin{equation}\label{EFE2}
T_{\al\be}:=\D_\al\psi\D_\be\psi-\frac{1}{2}\g_{\al\be}\big(\D_\mu\psi\D^\mu\psi+\psi^2\big),
\end{equation}
where $\D$ denotes covariant derivatives. 

Our goal is to prove definitive results on the global stability of the flat space among solutions of the Einstein-Klein-Gordon system. Our main theorems in this monograph include:

(1) A proof of global regularity (in wave coordinates) of solutions of the Einstein-Klein-Gordon coupled system, in the case of small, smooth, and localized perturbations of the stationary Minkowski solution $(\g,\psi)=(m,0)$;

(2) Precise asymptotics of the metric components and the Klein-Gordon field as the time goes to infinity, including the construction of modified (nonlinear) scattering profiles and quantitative bounds for convergence;

(3) Classical estimates on the solutions at null and timelike infinity, such as bounds on the metric components, weak peeling estimates of the Riemann curvature tensor, ADM and Bondi energy identities and estimates, and asymptotic description of null and timelike geodesics.

The general plan is to work in a standard gauge (the classical wave coordinates) and transform the geometric Einstein-Klein-Gordon system into a coupled system of quasilinear wave and Klein-Gordon equations. We then analyze this system in a framework inspired by the recent advances in the global existence theory for quasilinear dispersive models, such as plasma models and water waves.

More precisely, we rely on a combination of energy estimates and Fourier analysis. At a very general level one should think that energy estimates are used, in combination with vector-fields, to control high regularity norms of the solutions. The Fourier analysis is used, mostly in connection with normal forms, analysis of resonant sets, and a special ``designer" norm, to prove dispersion and decay in lower regularity norms. 

The method we present here incorporates Fourier analysis in a critical way. Its main advantage over the classical physical space methods is the ability to identify clearly resonant and non-resonant nonlinear quadratic interactions. We can then use normal forms to dispose of the non-resonant interactions, which leads to very precise estimates. In particular, some our asymptotic results appear to be new even in the case of the Einstein-vacuum equations (corresponding to $\psi=0$) mainly because we allow a large class of non-isotropic perturbations. Indeed, our assumptions on the metric on the initial slice are weak, essentially of the type 
\begin{equation*}
\g_{\al\be}=m_{\al\be}+\varep_0O(\langle x\rangle^{-1+}),\qquad \partial\g_{\al\be}=\varep_0O(\langle x\rangle^{-2+}),
\end{equation*}
consistent with non-isotropic decay, and we are still able to prove global regularity and derive suitable asymptotics, such as weak peeling estimates for the Riemann tensor and construct a Bondi energy function.

\subsection{Wave coordinates and PDE formulation of the problem} The system of equations \eqref{EFE}--\eqref{EFE2} is a geometric system, written in covariant form. To analyze it quantitatively and state our main theorems we need to fix a system of coordinates and reformulate our problem as a PDE problem.

We start by recalling some of the basic definitions and formulas of Lorentzian geometry. At this stage, all the formulas are completely analogous to the Riemannian case, hold in any dimension, and the computations can be performed in local coordinates. A standard reference is the book of Wald \cite{Wa}. Assume $\g$ is a sufficiently smooth Lorentzian metric in a $4$ dimensional open set $O$. We assume that we are working in a system of coordinates $x^0,x^1,x^2,x^3$ in $O$. We define the connection coefficients ${\bf \Gamma}$ and the covariant derivative $\D$ by
\begin{equation}\label{Ei3}
{\bf \Gamma}_{\mu\al\be}:=\g(\partial_\mu,\D_{\partial_\be}\partial_\al)=\frac{1}{2}(\partial_\al\g_{\be\mu}+\partial_\be\g_{\al\mu}-\partial_\mu\g_{\al\be}),
\end{equation}
where $\partial_\mu:=\partial_{x^\mu}$, $\mu\in\{0,1,2,3\}$. Thus
\begin{equation}\label{Ei4}
\D_{\partial_\alpha}\partial_\beta=\D_{\partial_\be}\partial_\al={\bf \Gamma}^\nu_{\,\,\,\,\al\be}\partial_\nu,\qquad {\bf \Gamma}^\nu_{\,\,\,\,\al\be}:=g^{\mu\nu}{\bf \Gamma}_{\mu\al\be},
\end{equation}
where $\g^{\al\be}$ is the inverse of the matrix $\g_{\al\be}$, i.e. $\g^{\al\be}\g_{\mu\be}=\delta_\mu^\al$. For $\mu,\nu\in\{0,1,2,3\}$ let
\begin{equation}\label{Ei2}
\begin{split}
&{\bf \Gamma}_\mu:=\g^{\al\be}{\bf \Gamma}_{\mu\al\be}=\g^{\al\be}\partial_\al \g_{\be\mu}-\frac{1}{2}\g^{\al\be}\partial_\mu\g_{\al\be},\qquad{\bf \Gamma}^\nu:=\g^{\mu\nu}{\bf \Gamma}_\mu.
\end{split}
\end{equation}
We record also the useful general identity
\begin{equation}\label{tr9}
\partial_\alpha\g^{\mu\nu}=-\g^{\mu\rho}\g^{\nu\lambda}\partial_\alpha \g_{\rho\lambda},
\end{equation}
and the Jacobi formula
\begin{equation}\label{Ei3.1}
\partial_\al(\log |\g|)=\g^{\mu\nu}\partial_\al\g_{\mu\nu},\qquad \al\in\{0,1,2,3\},
\end{equation}
where $|\g|$ denotes the determinant of the matrix $\g_{\al\be}$ in local coordinates.

Covariant derivatives can be calculated in local coordinates according to the general formula
\begin{equation}\label{Ei5}
\D_\al T_{\be_1\ldots\be_n}=\partial_\al T_{\be_1\ldots\be_n}-\sum_{j=1}^n{\bf \Gamma}^\nu_{\,\,\,\,\al\be_j}T_{\be_1\ldots\nu\ldots\be_n},
\end{equation}
for any covariant tensor $T$. In particular, for any scalar function $f$
\begin{equation}\label{Ei6}
\square_\g f=\g^{\al\be}\D_\al\D_\be f=\widetilde{\square}_\g f-{\bf \Gamma}^\nu\partial_\nu f,
\end{equation}
where $\widetilde{\square}_\g:=\g^{\al\be}\partial_\al\partial_\be$ denotes the reduced wave operator.

The Riemann curvature tensor measures commutation of covariant derivatives according to the covariant formula
\begin{equation}\label{Ei55}
\D_\al\D_\be\omega_\mu-\D_\be\D_\al\omega_\mu=\mathbf{R}_{\al\be\mu}^{\,\,\,\,\,\,\,\,\,\,\,\,\nu}\omega_\nu,
\end{equation}
for any form $\omega$. The Riemann tensor $\mathbf{R}$ satisfies the symmetry properties
\begin{equation}\label{Riem1}
\begin{split}
&\mathbf{R}_{\al\be\mu\nu}=-\mathbf{R}_{\be\al\mu\nu}=-\mathbf{R}_{\al\be\nu\mu}=\mathbf{R}_{\mu\nu\al\be},\\
&\mathbf{R}_{\al\be\mu\nu}+\mathbf{R}_{\be\mu\al\nu}+\mathbf{R}_{\mu\al\be\nu}=0,
\end{split}
\end{equation}
and the covariant Bianchi identities
\begin{equation}\label{Riem2}
\D_\rho \mathbf{R}_{\al\be\mu\nu}+\D_\al \mathbf{R}_{\be\rho\mu\nu}+\D_\be \mathbf{R}_{\rho\al\mu\nu}=0.
\end{equation}
Its components can be calculated in local coordinates in terms of the connection coefficients according to the formula
\begin{equation}\label{Ei6.5}
\mathbf{R}_{\al\be\mu}^{\,\,\,\,\,\,\,\,\,\,\,\,\rho}=-\partial_\alpha{\bf \Gamma}^\rho_{\,\,\,\,\be\mu}+\partial_\be{\bf\Gamma}^\rho_{\,\,\,\,\al\mu}-{\bf\Gamma}^\rho_{\,\,\,\,\al\nu}{\bf \Gamma}^\nu_{\,\,\,\,\be\mu}+{\bf \Gamma}^\rho_{\,\,\,\,\be\nu}{\bf \Gamma}^\nu_{\,\,\,\,\al\mu}.
\end{equation}
Therefore, the Ricci tensor $\R_{\al\mu}=g^{\be\rho}\mathbf{R}_{\al\be\mu\rho}$ is given by the formula
\begin{equation*}
\R_{\al\mu}=-\partial_\alpha{\bf\Gamma}^\rho_{\,\,\,\,\rho\mu}+\partial_\rho{\bf \Gamma}^\rho_{\,\,\,\,\al\mu}-{\bf \Gamma}^\rho_{\,\,\,\,\nu\al}{\bf \Gamma}^\nu_{\,\,\,\,\rho\mu}+{\bf \Gamma}^\rho_{\,\,\,\,\rho\nu}{\bf \Gamma}^\nu_{\,\,\,\,\al\mu}.
\end{equation*}
Simple calculations using \eqref{Ei3} and \eqref{Ei2} show that the Ricci tensor is given by
\begin{equation}\label{Ei20}
2\R_{\al\mu}=-\widetilde{\square}_\g \g_{\al\mu}+\partial_\al{\bf\Gamma}_\mu+\partial_\mu{\bf \Gamma}_\al+F^{\geq 2}_{\al\mu}(g,\partial g),
\end{equation}
where $F^{\geq 2}_{\al\be}(g,\partial g)$ is a quadratic semilinear expression, 
\begin{equation}\label{tr9.5}
\begin{split}
&F^{\geq 2}_{\al\be}(g,\partial g)=\frac{1}{2}\g^{\rho\mu}\g^{\nu\lambda}\big\{\partial_\nu \g_{\rho\mu}\partial_\be \g_{\al\lambda}+\partial_\nu \g_{\rho\mu}\partial_\al \g_{\be\lambda}-\partial_\nu \g_{\rho\mu}\partial_\lambda \g_{\al\be}\big\}\\
&+\g^{\rho\mu}\g^{\nu\lambda}\big\{-\partial_\rho \g_{\mu\lambda}\partial_\al \g_{\be\nu}-\partial_\rho \g_{\mu\lambda}\partial_\be \g_{\al\nu}+\partial_\rho \g_{\mu\lambda}\partial_\nu \g_{\al\be}+\partial_\al \g_{\rho\la}\partial_\mu \g_{\be\nu}+\partial_\be \g_{\rho\la}\partial_\mu \g_{\al\nu}\big\}\\
&-\frac{1}{2}\g^{\rho\mu}\g^{\nu\lambda}(\partial_\al \g_{\nu\mu}+\partial_\nu \g_{\al\mu}-\partial_\mu \g_{\al\nu})(\partial_\be \g_{\rho\lambda}+\partial_\rho \g_{\be\lambda}-\partial_\lambda \g_{\be\rho}).
\end{split}
\end{equation}

We consider the Einstein field equations \eqref{EFE}--\eqref{EFE2} for an unknown space-time $(M,\g)$; for simplicity, we drop the factor of $8\pi$ from of energy-momentum tensor. The covariant Bianchi identities $\D^\al G_{\a\be}=0$ can be used to derive an evolution equation for the massive scalar field $\psi$. The equation is
\begin{equation}\label{EFE3}
\square_\g\psi-\psi=0.
\end{equation}
Therefore the main unknowns in the problem are the metric tensor $\g$ and the scalar field $\psi$, which satisfy the covariant coupled equations \eqref{EFE} and \eqref{EFE3}. 

To construct solutions we need to fix a system of coordinates. In this paper we work in {\it{wave coordinates}}, which is the condition
\begin{equation}\label{wavecoo}
{\bf \Gamma}^\al=-\square_\g x^\al\equiv 0\,\,\text{ for }\,\,\alpha\in\{0,1,2,3\}.
\end{equation}
Our construction of global solutions of the Einstein-Klein-Gordon system is based on the following proposition, which can be proved by straightforward calculations.

\begin{proposition}\label{equivalence}
Assume $\g$ is a Lorentzian metric in a $4$ dimensional open set $O$, with induced covariant derivative $\D$ and Ricci curvature $\mathbf{R}_{\al\be}$, and $\psi:O\to\mathbb{R}$ is a scalar. Let $x^0,x^1,x^2,x^3$ denote a system of coordinates in $O$ and let ${\mathbf{\Gamma}}^\nu$ be defined as in \eqref{Ei2}.

(i) Assume that $(\g,\psi)$ satisfy the Einstein-Klein-Gordon system 
\begin{equation}\label{tr3}
\begin{split}
&\R_{\al\be}-\D_\al\psi\D_\be\psi-\frac{\psi^2}{2}\g_{\al\be}=0,\qquad\square_\g\psi-\psi=0,
\end{split}
\end{equation}
in $O$. Assume also that ${\bf\Gamma}^\mu\equiv 0$ in $O$, $\mu\in\{0,1,2,3\}$ (the harmonic gauge condition). Then
\begin{equation}\label{tr4}
\begin{split}
&\widetilde{\square}_\g \g_{\al\be}+2\partial_\al\psi\partial_\be\psi+\psi^2\g_{\al\be}-F^{\geq 2}_{\al\be}(g,\partial g)=0,\\
&\widetilde{\square}_\g\psi-\psi=0,
\end{split}
\end{equation}
where the quadratic semilinear terms $F^{\geq 2}_{\al\be}(g,\partial g)$ are defined in \eqref{tr9.5} and $\widetilde{\square}_\g:=\g^{\al\be}\partial_\al\partial_\be$ denotes the reduced wave operator .

(ii) Conversely, assume that the equations \eqref{tr4} (the reduced Einstein-Klein-Gordon system) hold in $O$. Then 
\begin{equation}\label{tr7}
\begin{split}
&\R_{\al\be}-\partial_\al\psi\partial_\be\psi-\frac{\psi^2}{2}\g_{\al\be}-\frac{1}{2}(\partial_\al {\bf\Gamma}_\be+\partial_\be{\bf \Gamma}_\al)=0,\\
&\square_\g\psi-\psi+{\bf\Gamma}^\mu\partial_\mu\psi=0,
\end{split}
\end{equation}
and the functions ${\bf\Gamma}_\be=g_{\be\nu}{\bf\Gamma}^\nu$ satisfy the reduced wave equations
\begin{equation}\label{tr6}
\begin{split}
\widetilde{\square}_\g{\bf\Gamma}_\be&=2{\bf\Gamma}^\nu\partial_\nu\psi\partial_\be\psi+\g^{\rho\al}[{\bf\Gamma}^\nu_{\,\,\,\,\rho\al}(\partial_\nu {\bf\Gamma}_\be+\partial_\be {\bf\Gamma}_\nu)+{\bf\Gamma}^\nu_{\,\,\,\,\rho\be}(\partial_\al {\bf\Gamma}_\nu+\partial_\nu {\bf\Gamma}_\al)]+\partial_\mu {\bf\Gamma}_\nu\partial_\be \g^{\mu\nu}.
\end{split}
\end{equation}
In particular, the pair $(\g,\psi)$ solves the Einstein-Klein-Gordon system \eqref{tr3} if ${\bf\Gamma}_\mu\equiv 0$ in $O$.
\end{proposition}

Our basic strategy to construct global solutions of the Einstein-Klein-Gordon system is to use Proposition \ref{equivalence}. We construct first the pair $(\g,\psi)$ by solving  the reduced Einstein-Klein-Gordon system \eqref{tr4} (regarded as a quasilinear wave-Klein-Gordon system) in the domain $\mathbb{R}^3\times[0,\infty)$. In addition, we arrange that ${\bf\Gamma}_\mu,\partial_t{\bf\Gamma}_\mu$ vanish on the initial hypersurface, so they vanish in the entire open domain, as a consequence of the wave equations \eqref{tr6}. Therefore the pair $(\g,\psi)$ solves the Einstein-Klein-Gordon system as desired. 

In other words, the problem is reduced to constructing global solutions of the quasilinear system \eqref{tr4} for initial data compatible with the wave coordinates condition.

\section{The global regularity theorem} 

To state our global regularity theorem we introduce first several spaces of functions on $\mathbb{R}^3$. 

\begin{definition}\label{spaces}
For $a\geq 0$ let $H^a$ denote the usual Sobolev spaces of index $a$ on $\mathbb{R}^3$. We define also the Banach spaces $H^{a,b}_\Omega$, $a,b\in\mathbb{Z}_+$, by the norms
\begin{equation}\label{qaz4}
\|f\|_{H^{a,b}_{\Omega}}:=\sum_{|\alpha|\leq b}\|\Omega^\alpha f\|_{H^a}.
\end{equation} 
We also define the weighted Sobolev spaces $H^{a,b}_{S,wa}$ and $H^{a,b}_{S,kg}$ by the norms
\begin{equation}\label{qaz4.2}
\|f\|_{H^{a,b}_{S,wa}}:=\sum_{|\beta'|\leq |\beta|\leq b}\|x^{\beta'}\partial^{\beta} f\|_{H^a},\qquad \|f\|_{H^{a,b}_{S,kg}}:=\sum_{|\beta|,|\beta'|\leq b}\|x^{\beta'}\partial^{\beta} f\|_{H^a},
\end{equation}
where $x^{\beta'}=x_1^{\beta'_1}x_2^{\beta'_2}x_3^{\beta'_3}$ and $\partial^\beta:=\partial_1^{\beta_1}\partial_2^{\beta_2}\partial_3^{\beta_3}$. Notice that $H^{a,b}_{S,kg}\hookrightarrow H^{a,b}_{S,wa}\hookrightarrow H^{a,b}_{\Omega}\hookrightarrow H^a$.
\end{definition}

To implement the strategy described above and use Proposition \ref{equivalence} we need to prescribe suitable initial data. Let $\Sigma_0=\{(x,t)\in \mathbb{R}^3\times[0,\infty):\,t=x^0=0\}$. We assume that $\overline{g},k$ are given symmetric tensors on $\Sigma_0$, such that $\overline{g}$ is a Riemannian metric on $\Sigma_0$. We assume also that $\psi_0,\psi_1:\Sigma_0\to\mathbb{R}$ are given initial data for the scalar field $\psi$. 

We start by prescribing the metric components on $\Sigma_0$
\begin{equation*}
\g_{ij}=\overline{g}_{ij},\qquad \g_{0i}=0,\qquad \g_{00}=-1.
\end{equation*}
The conditions $\g_{00}=-1$ and $\g_{0i}=0$ hold only on the initial hypersurface and are not propagated by the flow. They are imposed mostly for convenience and do not play a significant role in the analysis. We also prescribe the time derivative of the metric tensor
\begin{equation*}
\partial_t\g_{ij}=-2 k_{ij},
\end{equation*}
in such a way that $k$ is the second fundamental form of the surface $\Sigma_0$, $k(X,Y)=-\g(\D_Xn,Y)$, where $n=\partial_0$ is the future-oriented unit normal vector-field on $\Sigma_0$. The conditions ${\bf \Gamma}_\alpha=0$, $\alpha\in\{0,1,2,3\}$, can be used to determine the other components of the initial data for the pair $(\g,\psi)$ on the hypersurface $\Sigma_0$, which are
\begin{equation}\label{data}
\begin{split}
&\g_{ij}=\overline{g}_{ij},\qquad \g_{0i}=\g_{i0}=0,\qquad \g_{00}=-1,\\
&\partial_t\g_{ij}=-2k_{ij},\qquad \partial_t\g_{00}=2\overline{g}^{ij}k_{ij},\qquad \partial_t\g_{n0}=\overline{g}^{ij}\partial_i\overline{g}_{jn}-\frac{1}{2}\overline{g}^{ij}\partial_n\overline{g}_{ij},\\
&\psi=\psi_0,\qquad \partial_t\psi=\psi_1.
\end{split}
\end{equation} 

The remaining restrictions $\partial_t{\bf\Gamma}_\al=0$ lead to the constraint equations. In view of \eqref{tr7} the constraint equations are equivalent to the conditions $\R_{\al 0}-(1/2)\R\g_{\al 0}=T_{\al0}$, $\al\in\{0,1,2,3\}$, where $T_{\al\be}$ is as in \eqref{EFE2}. This leads to four constraint equations
\begin{equation}\label{asum2}
\begin{split}
\overline{D}_n(\overline{g}^{ij}k_{ij})-\overline{g}^{ij}\overline{D}_jk_{in}&=\psi_1\overline{D}_n\psi_0,\qquad n\in\{1,2,3\},\\
\overline{R}+\overline{g}^{ij}\overline{g}^{mn}(k_{ij}k_{mn}-k_{im}k_{jn})&=\psi_1^2+\overline{g}^{ij}\overline{D}_i\psi_0\overline{D}_j\psi_0+\psi_0^2,
\end{split}
\end{equation}
where $\overline{D}$ denotes the covariant derivative induced by the metric $\overline{g}$ on $\Sigma_0$, and $\overline{R}$ is the scalar curvature of the metric $\overline{g}$ on $\Sigma_0$.

We are now ready to state our first main theorem, which concerns global regularity of the system \eqref{tr4} for small initial data $(\overline{g}_{ij},k_{ij},\psi_0,\psi_1)$.

\begin{theorem}\label{Main}
Let $\Sigma_0:=\{(x,t)\in\mathbb{R}^4:t=0\}$ and assume that $(\overline{g}_{ij},k_{ij},\psi_0,\psi_1)$ is an initial  data set on $\Sigma_0$, satisfying the constraint equations \eqref{asum2} and the smallness conditions
\begin{equation}\label{asum1}
\begin{split}
\sum_{n=0}^{3}\sum_{i,j=1}^3\big\{\big\|\,|\nabla|^{1/2+\delta/4}&(\overline{g}_{ij}-\delta_{ij})\big\|_{H^{N(n),n}_{S,wa}}+\|\,|\nabla|^{-1/2+\delta/4}k_{ij}\|_{H^{N(n),n}_{S,wa}}\big\}\\
&+\sum_{n=0}^{3}\big\{\|\langle\nabla\rangle\psi_0\|_{H^{N(n),n}_{S,kg}}+\|\psi_1\|_{H^{N(n),n}_{S,kg}}\big\}\leq\varepsilon_0\leq\overline{\varepsilon}.
\end{split}
\end{equation}
Here $N_0:=40$, $d:=10$, $\delta:=10^{-10}$, $N(0):=N_0+16d$, $N(n)=N_0-nd$ for $n\geq 1$, $\overline{\varepsilon}$ is a small constant, and the operators $|\nabla|$ and $\langle\nabla\rangle$ are defined by the multipliers $|\xi|$ and $\langle\xi\rangle$.

(i) Then the reduced Einstein-Klein-Gordon system
\begin{equation}\label{asum3}
\begin{split}
\widetilde{\square}_\g \g_{\al\be}+2\partial_\al\psi\partial_\be\psi+\psi^2\g_{\al\be}-F^{\geq 2}_{\al\be}(g,\partial g)&=0,\\
\widetilde{\square}_\g\psi-\psi&=0,
\end{split}
\end{equation}
admits a unique global solution $(\g,\psi)$ in $M:=\{(x,t)\in\mathbb{R}^4:t\geq 0\}$, with initial data $(\overline{g}_{ij},k_{ij},\psi_0,\psi_1)$ on $\Sigma_0$ (as described in \eqref{data}). Here $F^{\geq 2}_{\al\be}(g,\partial g)$ are as in \eqref{tr9.5} and $\widetilde{\square}_\g=\g^{\mu\nu}\partial_\mu\partial_\nu$. The solution satisfies the harmonic gauge conditions
\begin{equation}\label{asum5}
0={\bf\Gamma}_\mu=\g^{\al\be}\partial_\al\g_{\be\mu}-\frac{1}{2}\g^{\al\be}\partial_\mu\g_{\al\be},\qquad\mu\in\{0,1,2,3\}
\end{equation} 
in $M$.  Moreover, the metric $\g$ stays close and converges as $t\to\infty$ to the Minkowski metric and $\psi$ stays small and converges to $0$ as $t\to\infty$ (in suitable norms). 

(ii) In view of Proposition \ref{equivalence}, the pair $(\g,\psi)$ is a global\footnote{In our geometric context, globality means that all future directed timelike and null geodesics starting from points in $M$ extend forever with respect to their affine parametrization.}  solution in $M$ of the Einstein-Klein-Gordon coupled system 
\begin{equation}\label{tr3.7}
\R_{\al\be}-\D_\al\psi\D_\be\psi-\frac{\psi^2}{2}\g_{\al\be}=0,\qquad \square_\g\psi-\psi=0.
\end{equation}
\end{theorem}

The proof of Theorem \ref{Main} is based on a complex bootstrap argument, involving energy estimates, vector-fields, Fourier analysis, and nonlinear scattering. We outline some of its main elements in subsection \ref{Remar2} below.

The global regularity conclusion of Theorem \ref{Main} is essentially a qualitative statement, which can only be proved by a precise quantitative analysis of the space-time. In Chapter \ref{Chapter6} we state and prove more precise theorems describing our space-time. These theorems include global quantitative control and nonlinear scattering of the metric tensor and the Klein-Gordon field (Theorem \ref{winr}), pointwise decay estimates in the physical space (Theorem \ref{Precisedmetric} and Lemma \ref{twoder0}), global control of timelike and null geodesics (Theorem \ref{winr100}), weak peeling estimates for the Riemann curvature tensor (Theorem \ref{PeelingR} and Proposition \ref{PeelingR2}), and ADM and Bondi energy formulas (Proposition \ref{ADM4}, Proposition \ref{linmo2+}, Theorem \ref{BondThm}, and Proposition \ref{Bond40}). We will discuss some of these more precise conclusions in section \ref{Remar2} below.

In the rest of this section we discuss previous related work and motivate some of the assumptions on the initial data.

\subsection{Global stability results in General Relativity}  Global stability of physical solutions is an important topic in General Relativity. For example, the global nonlinear stability of the Minkowski space-time among solutions of the Einstein-vacuum equation is a central theorem in the field, due to Christodoulou-Klainerman \cite{ChKl}. See also the more recent extensions of Klainerman-Nicol\`o \cite{KlNi}, Lindblad-Rodnianski \cite{LiRo}, Bieri and Zipser \cite{BiZi}, Speck \cite{Sp}. 

More recently, small data global regularity theorems have also been proved for other coupled Einstein field equations. The Einstein-Klein-Gordon system (the same system we analyze here) was considered recently by LeFloch-Ma \cite{LeMaEKG}, who proved small data global regularity for restricted data, which agree with a Schwarzschild solution with small mass outside a compact set. A similar result was announced by Wang \cite{Wa1}. 

Our main goals in this monograph are (1) to work with general unrestricted small initial data, and (2) develop the full asymptotic analysis of the space-time. A similar global regularity result for general small data was announced recently by LeFloch-Ma \cite{LeMaNew}.

We also refer to the recent work by Fajman-Joudioux-Smulevici \cite{FaJoSm} and Lindblad-Taylor \cite{LiTa} on the global stability of the Einstein-Vlasov system, at least for certain classes of ``restricted data" (see below for a longer discussion).

In a different direction, one can also raise the question of linear and nonlinear stability of other physical solutions of the Einstein equations. Stability of the Kerr family of solutions has been under intense study in recent years, first at the linearized level (see for example  \cite{DaRoSh,HaHiVa} and the references therein) and more recently at the full nonlinear level (see \cite{GiKlSz,HiVa,KlSz}). 

The stability of Kerr in the presence of a massive scalar field seems interesting as well. Solutions to the Klein-Gordon equation in Kerr can grow exponentially even from smooth initial data, as shown in \cite{Shl}, 
and this phenomenon was used by Chodosh--Shlapentokh-Rothman \cite{ChSh} to construct a curve of time-periodic solutions of the Einstein-Klein-Gordon system bifurcating from (empty) Kerr (see \cite{HeRa} for a prior numerical construction). Therefore a result on stability of Kerr similar to our main theorem could only be possible, if at all, in a stronger topology where this curve is not continuous (see also the discussion on the mini bosons in subsection \ref{minibos} below).

\subsubsection{Restricted initial data} One can often simplify considerably the global ana\-ly\-sis of wave and Klein-Gordon equations by considering initial data of compact support. The point is that the solutions have the finite speed of propagation, thus remain supported inside a light cone, and one can use the hyperbolic foliation method and its refinements (see \cite{LeMaHyp} for a recent account) to analyze the evolution.

However, to implement this method one needs to first control the solution on an initial hyperboloid (the ``initial data"), so the method is restricted to the case when one can establish such control. Due to the finite speed of propagation, this is possible for compactly supported data (for systems of wave or Klein-Gordon equations), or data that agrees with the Schwarzschild solution outside a compact set (in the case of the Einstein equations).

The use of ``restricted initial data" coupled with the hyperbolic foliation method, leads to significant simplifications of the global analysis, particularly at the level of proving decay. In the context of the Einstein equations these ideas have been used by many authors, such as Friedrich \cite{Fr}, Lindblad-Rodnianski \cite{LiRo}, Fajman-Joudioux-Smulevici \cite{FaJoSm}, Lindblad-Taylor \cite{LiTa}, LeFloch-Ma \cite{LeMaEKG}, Wang \cite{Wa1}, and Klainerman-Szeftel \cite{KlSz}. 

\subsection{Simplified Wave-Klein-Gordon models} The system \eqref{asum3} is complicated, but one can gain intuition by looking at simpler models. For example, one can consider the simplified  system
\begin{equation}\label{on1}
\begin{split}
-\square u&=A^{\alpha\beta}\partial_\alpha v\partial_\beta v+Dv^2,\\
(-\square +1)v&=uB^{\alpha\beta}\partial_\alpha\partial_\beta v+Euv,
\end{split}
\end{equation}
where $u,v$ are real-valued functions, and $A^{\alpha\beta}$, $B^{\alpha\beta}$, $D$, and $E$ are real constants. 

The system \eqref{on1} was derived by LeFloch-Ma \cite{LeMa} as a model for the full Einstein-Klein-Gordon system \eqref{asum3}. Intuitively, the deviation of the Lorentzian metric $\g$ from the Minkowski metric is replaced by a scalar function $u$, and the massive scalar field $\psi$ is replaced by $v$. The system \eqref{on1} has the same linear structure as the Einstein-Klein-Gordon system \eqref{asum3}, but only keeps, schematically, quadratic interactions that involve the massive scalar field.

Small data global regularity for the system \eqref{on1} was proved by LeFloch-Ma \cite{LeMa} in the case of compactly supported initial data (the restricted data case), using the hyperbolic foliation method. For general small initial data, global regularity was proved by the authors \cite{IoPa3}. 

A similar system, the massive Maxwell-Klein-Gordon system, was analyzed recently by Klainerman--Wang--Yang \cite{KlWaYa}, who also proved global regularity for general small initial data, using a different method. Coupled Wave-Klein-Gordon systems have also been considered in 2D, where the decay is slower and the global analysis requires nonlinearities with much more favorable structure (see, for example, Ifrim-Stingo \cite{IfSt} and the references therein).
\smallskip

\subsection{Small data global regularity results}\label{Remar5} The system \eqref{asum3} can be easily transformed into a quasilinear coupled system of wave and Klein-Gordon equations. Indeed, let $m$ denote the Minkowski metric and write
\begin{equation*}
\g_{\al\be}=m_{\al\be}+h_{\al\be},\qquad \g^{\al\be}=m^{\al\be}+g_{\geq 1}^{\al\be},\qquad\al,\be\in\{0,1,2,3\}.
\end{equation*}
It follows from \eqref{asum3} that the metric components $h_{\al\be}$ satisfy the nonlinear wave equations
\begin{equation}\label{laq11.1}
(\partial_0^2-\Delta)h_{\al\be}=\mathcal{N}^{h}_{\al\be}:=\mathcal{KG}_{\al\be}+g_{\geq 1}^{\mu\nu}\partial_\mu\partial_\nu h_{\al\be}-F^{\geq 2}_{\al\be}(g,\partial g)
\end{equation}
where $F^{\geq 2}_{\al\be}(g,\partial g)$ are the semilinear terms  in \eqref{tr9.5} and $\mathcal{KG}_{\al\be}:=2\partial_\al\psi\partial_\be\psi+\psi^2(m_{\al\be}+h_{\al\be})$. Moreover, the field $\psi$ satisfies the quasilinear Klein-Gordon equation
\begin{equation}\label{laq11.3}
(\partial_0^2-\Delta+1)\psi=\mathcal{N}^{\psi}:=g^{\mu\nu}_{\geq 1}\partial_\mu\partial_\nu\psi.
\end{equation}

Therefore Theorem \ref{Main} can be regarded as a small data global regularity result for a quasilinear evolution system. Several important techniques have been developed over the years in the study of such problems, starting with seminal contributions of John, Klainerman, Shatah, Simon, Christodoulou, Alinhac, and Delort \cite{Alin2,Alin3,Ch,ChKl,Del,DeFaXu,Jo,JoKl,Kl2,KlVf,Kl,Kl4,Sh,Si}. These include the vector-field method, normal forms, and the isolation of null structures. 

In the last few years new methods have emerged in the study of global solutions of quasilinear evolutions, inspired mainly by the advances in semilinear theory.
The basic idea is to combine the classical energy and vector-fields methods with refined analysis of the Duhamel formula, using the Fourier transform.
This is the essence of the ``method of space-time resonances'' of Germain-Masmoudi-Shatah \cite{GeMaSh,GeMaSh2} and Gustafson-Nakanishi-Tsai \cite{GuNaTs},
and the refinements by the authors and their collaborators in \cite{DeIoPa,DeIoPaPu,GuIoPa,GuPa,IoPa1,IoPa2,IoPu2,IoPu,IoPu3,KaPu}, using atomic decompositions and more sophisticated norms. 

According to this general philosophy, to prove Theorem \ref{Main} we work both in the physical space, mainly to prove energy estimates (including vector-fields), and in the Fourier space, mainly to prove decay of the solutions in time and nonlinear scattering. One of the main difficulties in the analysis of the system \eqref{laq11.1}--\eqref{laq11.3} comes from the fact that this is a multiple speed system, in the sense that the linear evolution admits different speeds of propagation, corresponding to the wave and the Klein-Gordon components. As a result the set of ``characteristics'' and the quadratic resonances are more involved, and one has a more limited set of geometric symmetries.

\subsection{Initial-data assumptions} The precise form of the smallness assumptions \eqref{asum1} on the metric initial data $\overline{g}_{ij}$ and $k_{ij}$ is important. Indeed, in view of the positive mass theorem of Schoen-Yau \cite{ScYa}, one expects the metric components $\overline{g}_{ij}-\delta_{ij}$ to decay like $M/\langle x\rangle$ and the second fundamental form $k$ to decay like $M/\langle x\rangle^2$, where $M\ll 1$ is the mass. Capturing this type of decay, using $L^2$ based norms, is precisely the role of the homogeneous multipliers $|\nabla|^{1/2+\delta/4}$ and $|\nabla|^{-1/2+\delta/4}$ in \eqref{asum1}. Notice that these multipliers are essentially sharp, up to the $\delta/4$ power.

Our assumptions on the metric are essentially of the type 
\begin{equation}\label{approxas}
g_{ij}=\delta_{ij}+\varep_0O(\langle x\rangle^{-1+\delta/4}),\qquad k_{ij}=\varep_0O(\langle x\rangle^{-2+\delta/4})
\end{equation}
at time $t=0$. These are less restrictive than the assumptions used sometimes even in the vacuum case $\psi\equiv 0$, see for example \cite{ChKl}, \cite{KlNi}, or \cite{LiRo2}, in the sense that the initial data is not assumed to agree with the Schwarzschild initial data up to lower order terms. For maximal time foliations, our assumptions are, however, more restrictive than the ones in Bieri work \cite{BiZi}, but we are able to prove more precise asymptotic bounds on the metric and the Riemann curvature tensor, see section \ref{Remar2} below.

We remark also that our assumptions \eqref{asum1} allow for non-isotropic initial data, possibly with different ``masses" in different directions. For the vacuum case, initial data of this type, satisfying the constraint equations, have been constructed recently by Carlotto-Schoen \cite{CaSh}.

\subsection{The mini-bosons}\label{minibos} A serious potential obstruction to small data global stability theorems is the presence of non-decaying ``small" solutions, such as small solitons. A remarkable fact is that there are such small non-decaying solutions for the Einstein-Klein-Gordon system, namely the so-called mini-boson stars. These are time-periodic (therefore non-decaying) and spherically symmetric exact solutions of the Einstein-Klein-Gordon system. They were discovered numerically by physicists, such as Kaup \cite{Kau}, Friedberg--Lee--Pang \cite{FrLePa} (see also \cite{LiPa}), and then constructed rigorously by Biz\'{o}n--Wasserman \cite{BiWa}.

These mini-bosons can be thought of as arbitrarily small (hence the name) in certain topologies, as explained in \cite{BiWa}. However, the mini-bosons (in particular the Klein-Gordon component) are not small in the stronger topology we use here, as described by \eqref{asum1}, so we can thankfully avoid them in our analysis.

\section{Main ideas and further asymptotic results}\label{Remar2}

The classical mechanism to establish small data global regularity for quasilinear equations has two main components:

(1) Propagate control of energy functionals (high order Sobolev norms and vector-fields);

(2) Prove dispersion/decay of the solution over time.

These are our basic goals here as well. In subsections \ref{hus1}--\ref{hus6} we describe some of the main difficulties and outline some of our arguments to prove global regularity.

In subsections \ref{hus7}--\ref{BonIntro} we summarize some of the additional results we prove in Chapter 6, concerning the global geometry of our space-time.

\subsection{Energy estimates}\label{hus1} A key component of our analysis is to prove energy estimates for solutions of our system \eqref{laq11.1}--\eqref{laq11.3}. These energy estimates involve vector-fields, corresponding to the natural symmetries of the linearized equations. In our case we use the Lorentz vector-fields $\Gamma_a$ and the rotation vector-fields $\Omega_{ab}$ (from the Minkowski space)
\begin{equation}\label{qaz2}
\Gamma_a:=x_a\partial_t+t\partial_a,\qquad \Omega_{ab}:=x_a\partial_b-x_b\partial_a,
\end{equation}
for $a,b\in\{1,2,3\}$. These vector-fields commute with both the wave operator and the Klein-Gordon operator in the flat Minkowski space.  We note that the scaling vector field, $S=t\partial_t+x\cdot\nabla_x$ does not satisfy nice commutation properties with the linearized system (due to the Klein-Gordon field), so we cannot use it in our analysis.

To state our main energy estimates we define the normalized solutions $U^{\LL h_{\al\be}}$ and $U^{\LL \psi}$ and the associated {\it{linear profiles}}  $V^{\LL h_{\al\be}}$ and $V^{\LL\psi}$ by
\begin{equation}\label{introdesc}
\begin{split}
U^{\LL h_{\al\be}}(t)&:=\partial_t (\LL h_{\al\be})(t)-i\Lambda_{wa} (\LL h_{\al\be})(t),\qquad V^{\LL h_{\al\be}}(t):=e^{it\Lambda_{wa}}U^{\LL h_{\al\be}}(t),\\
U^{\LL\psi}(t)&:=\partial_t(\LL \psi)(t)-i\Lambda_{kg} (\LL\psi)(t),\qquad\qquad\,\,\,\,\, V^{\LL\psi}(t):=e^{it\Lambda_{kg}}U^{\LL\psi}(t),
\end{split}
\end{equation}
where $\Lambda_{wa}=|\nabla|$, $\Lambda_{kg}=\langle\nabla\rangle=\sqrt{|\nabla|^2+1}$. Here $\LL$ denote differential operators obtained by applying up to 3 vector-fields $\Gamma_a$ or $\Omega_{ab}$, and these operators are applied to the metric components $h_{\al\be}$ and the field $\psi$. 

The complex-valued normalized solutions $U^{\LL h_{\al\be}}$ and $U^{\LL \psi}$ capture both the time derivatives (as the real part) and the spatial derivatives (as the imaginary part) of the variables $h_{\al\be}$ and $\psi$. The linear profiles $V^{\LL h_{\al\be}}$ and $V^{\LL \psi}$, which are defined by going forward in time along the nonlinear evolution, and then going backwards in time along the linear flow, capture the cumulative effect of the nonlinearity over time.

The main energy estimates we prove in this paper are
\begin{equation}\label{introdesc2}
\begin{split}
\big\|(\langle t\rangle|\nabla|_{\leq 1})^{\delta/4}|\nabla|^{-1/2}U^{\LL h_{\al\be}}(t)\big\|_{H^{n(\mathcal{L})}}+\big\|U^{\LL\psi}(t)\big\|_{H^{n(\mathcal{L})}}&\lesssim\varep_0\langle t\rangle^{H(\mathcal{L})\delta},
\end{split}
\end{equation}
for a suitable hierarchy of parameters $n(\mathcal{L})$ and $H(\mathcal{L})$ that depend on the differential operator $\mathcal{L}$. We remark that the energy estimates we prove for the metric variables $U^{\LL h_{\al\be}}$ also contain significant information at low frequencies, due to the factors $|\nabla|^{-1/2}$, which are connected to the natural $|x|^{-1+}$ decay of the metric components $h_{\al\be}$.

\subsection{Null structure and decomposition of the metric tensor} The proof of the energy estimates \eqref{introdesc2} is involved, mainly because the nonlinearities $\mathcal{N}^h_{\al\be}$ and $\mathcal{N}^\psi$ have complicated structure, both at the semilinear level (for $\mathcal{N}^h_{\al\be}$) and at the quasilinear level.

To begin with, the semilinear terms $F_{\al\be}^{\geq 2}(g,\partial g)$ do not have the classical null structure. They have, however, what is called a weak null structure still suitable for global analysis, as discovered by Lindblad-Rodnianski \cite{LiRo}. To identify and use this weak null structure we need to decompose the tensor $h_{\al\be}$. 

The standard way to decompose the metric tensor in General Relativity is based on null frames (see for instance \cite{ChKl} or \cite{LiRo}). Here we use a different decomposition of the metric tensor, reminiscent of the $\mathrm{div}$-$\mathrm{curl}$ decomposition of vector-fields in fluid models, which is connected to the classical work of Arnowitt--Deser--Misner \cite{ADM} on the Hamiltonian formulation of General Relativity. For us, this decomposition has the advantage of being more compatible with the Fourier transform and the vector-fields $\Omega_{ab},\Gamma_a$.

More precisely, let $R_j=|\nabla|^{-1}\partial_j$, $j\in\{1,2,3\}$, denote the Riesz transforms on $\mathbb{R}^3$, and let 
\begin{equation}\label{laq2}
\begin{split}
&F:=(1/2)[h_{00}+R_jR_kh_{jk}],\qquad \uF:=(1/2)[h_{00}-R_jR_kh_{jk}],\\
&\rho:=R_jh_{0j},\qquad \omega_j:=\in_{jkl}R_kh_{0l},\\
&\Omega_j:=\in_{jkl}R_kR_mh_{lm},\qquad \vartheta_{jk}:=\in_{jmp}\in_{knq}R_mR_nh_{pq}.
\end{split}
\end{equation} 
Geometrically, the variables $F$, $\underline{F}$, $\rho$, and $\omega$ are linked to the shift vector, while $\vartheta$ corresponds to the (linearized) coordinate free component of the spatial metric (or coincides at the linear level with the
spatial metric in spatially harmonic coordinates). The metric tensor $h$ can be recovered linearly from the components $F,\uF,\rho,\omega_j,\Omega_j,\va_{jk}$. 

Our analysis shows that the components $F,\omega_j,\Omega_j,\va_{jk}$ satisfy good wave equations, with all the quadratic semilinear terms having suitable null structure. On the other hand, the components $\uF$ and $\rho$ (which are related elliptically due to the harmonic gauge conditions) satisfy wave equations with some quadratic semilinear terms with no null structure. However, these non-null quadratic semilinear terms have the redeeming feature that they can be expressed only in terms of the good components $\vartheta_{jk}$. 

This algebraic structure suggests that we should aim to prove that the good components $F,\omega_j,\Omega_j,\va_{jk}$ do not grow during the evolution, in suitable norms to be made precise. On the other hand, the components $\uF,\rho$ as well as all the components $\LL h_{\al\be}$ and $\LL\psi$ which contain some weighted vector-fields $\Omega_{ab},\Gamma_a$ should be allowed to grow in time\footnote{Notice that these vector-fields are adapted to the Minkowski geometry, containing the coordinate functions $x_a$ and $t$, thus they can only be useful up to $\langle t\rangle^{0+}$ losses.} slowly, at suitable rates to be determined. At a qualitative level, this is precisely what our final conclusions are.

\subsection{Low frequencies and the quasilinear terms} To prove \eqref{introdesc2} we would like to apply the operators $\mathcal{L}$ to the equations \eqref{laq11.1}--\eqref{laq11.3} and perform energy estimates. The most difficult contributions are coming from the quasilinear terms $g_{\geq 1}^{\mu\nu}\partial_\mu\partial_\nu h_{\al\be}$ and $g_{\geq 1}^{\mu\nu}\partial_\mu\partial_\nu\psi$, when the entire operator $\LL$ hits the undifferentiated metric components $g_{\geq 1}^{\mu\nu}$, i.e. the terms
\begin{equation}\label{introdesc3}
\LL(g_{\geq 1}^{\mu\nu})\cdot \partial_\mu\partial_\nu h_{\al\be}\qquad\text{ and }\qquad \LL(g_{\geq 1}^{\mu\nu})\cdot \partial_\mu\partial_\nu \psi.
\end{equation}
These terms have their own quasilinear null structure, which helps with the analysis of medium frequencies, but the real difficulty is to bound the contributions of very low frequencies of the factor $\LL(g_{\geq 1}^{\mu\nu})$. The issue is that this factor does not have a spatial derivative to help at very low frequencies. When $\LL=Id$, we can symmetrize the system and obtain improved estimates, but no algebraic symmetrization is possible when $\mathcal{L}\neq Id$. This is the most difficult case of the analysis, and requires establishing a special hierarchy of growth and normal forms.

Some of the main choices we make in the proof, like the low frequency structure given by the multiplier $|\nabla|^{-1/2}$ in the energy estimates \eqref{introdesc2} or the precise hierarchy of energy growth (see \eqref{fvc1.0}) are motivated by the energy analysis of the terms in \eqref{introdesc3}. Even for moderately small frequencies, when the null structure of the quasilinear terms can be used to construct a normal form, one still has a loss of derivatives (at low frequencies) when iterating the equation. This is handled in a similar way to \cite{DeIoPa, DeIoPaPu}, using paradifferential calculus to extract the quasilinear terms and set up a second symmetrization.

\subsection{Wave-Klein-Gordon interactions} The main difficulty in the analysis of the Einstein-Klein-Gordon system is to understand the quadratic interactions of the metric tensor and the massive field. There are two such interactions: the semilinear quadratic terms $\langle\nabla\rangle\psi\ast\langle\nabla\rangle\psi$, which are part of the metric nonlinearities in $\mathcal{N}^h_{\al\be}$, and the quasilinear quadratic terms $h\ast\langle\nabla\rangle^2\psi$, which are part of nonlinearity $\mathcal{N}^\psi$

These bilinear interactions do not have natural null structure. In fact there are even undifferentiated terms, such as 
\begin{equation}\label{introdesc3.2}
\text{ terms like }\psi^2\text{ in }\mathcal{N}^h_{\al\be}\qquad\text{ and }\qquad \text{ terms like }h\psi\text{ in }\mathcal{N}^\psi,
\end{equation}
where the last terms arise by writing $g^{00}_{\geq 1}\partial_0^2\psi=-h_{00}(\Delta-1)\psi+\text{cubic expressions}$. 

Our main tool to estimate these interactions is to use normal forms. This requires that we understand their resonant structure, which is controlled by the quadratic phase
\begin{equation}\label{PhaseReso}
\Phi(\xi,\eta)=\Lambda_{kg}(\xi)\pm\Lambda_{kg}(\eta)\pm \Lambda_{wa}(\xi+\eta),
\end{equation} 
where $\Lambda_{kg}$ and $\Lambda_{wa}$ are defined as in \eqref{introdesc}. The main observation is that
\begin{equation}\label{PhaseReso2}
\vert\Phi(\xi,\eta)\vert\gtrsim (1+\vert\xi\vert+\vert\eta\vert)^{-2}\vert\xi+\eta\vert,
\end{equation}
see Lemma \ref{pha2} for more precise statement of this type. This means that the only resonant interactions involve the $0$ frequency of the metric components, which highlights again the importance of having precise optimal estimates on the low frequencies of the metric tensor. 

The undifferentiated terms in \eqref{introdesc2} have their own specific issues, mostly having to do with the low frequencies of the Klein-Gordon field $\psi$. These low frequencies disperse slowly, and one should think of these undifferentiated terms as the most difficult wave-Klein-Gordon interactions. In particular, the Klein-Gordon $Z$-norm defined in \eqref{introdesc5} below is designed at low frequencies mainly to be able to capture these terms.

The analysis of wave-Klein-Gordon bilinear interactions has similarities with the analysis in the earlier work of the authors \cite{IoPa3}. However, in \cite{IoPa3} the equation for the wave component $u$ was much simpler, and it was possible to prove sharp $1/t$ pointwise decay for this component. The situation here is substantially more difficult mainly because some of the metric components do not have this optimal pointwise decay, and we need additional normal form arguments to control the contribution of the Klein-Gordon nonlinearity $\mathcal{N}^\psi$.

\subsection{Weighted estimates on the profiles and pointwise decay} The energy estimates \eqref{introdesc2} can be used to prove weighted bounds on the profiles $V^{\LL h_{\al\be}}$ and $V^{\LL\psi}$. Indeed, it turns out that the application of a vector-field $\Gamma_a$ on a normalized solution is essentially equivalent to multiplication by $x_a$ of the associated profile, up to nonlinear effects (see Lemma \ref{ident} for the precise identity). Therefore, using \eqref{introdesc2}, we have
\begin{equation}\label{introdesc4}
\begin{split}
2^{k/2}(2^{k^-}\langle t\rangle)^{\delta/4}\|P_k(x_lV^{\mathcal{L}h_{\al\be}})(t)\|_{L^2}+2^{k^+}\|P_k(x_lV^{\mathcal{L}\psi})(t)\|_{L^2}\lesssim\varep_0\langle t\rangle^{H'(\LL)\delta}2^{-n'(\mathcal{L})k^+},
\end{split}
\end{equation}
for any $k\in\mathbb{Z}$, $l\in\{1,2,3\}$, and differential operator $\LL$ containing at most $2$ vector-fields $\Gamma_a$ or $\Omega_{ab}$. Here $P_k$ denote Littlewood-Paley projections to frequencies $\approx 2^k$.

The weighted profile bounds in \eqref{introdesc4} are an important component of our bootstrap argument. Using linear estimates (Lemma \ref{LinEstLem}) they imply almost optimal $\langle t\rangle^{-1+}$ pointwise decay estimates on the metric components and the Klein-Gordon field, with improved decay both at very low and very high frequency.

We emphasize, however, that weighted estimates on linear profiles are a lot stronger than pointwise decay estimates on solutions, and serve many other purposes. For example, space localization of the linear profiles gives us the main information we need to decompose the various nonlinear contributions both in frequency and space.
\smallskip

\subsection{Uniform bounds and the $Z$-norm}\label{hus5} We still need to deal with an important issue, namely we have to show that the good metric components $F,\omega_j,\Omega_j,\va_{jk}$ and the massive field $\psi$ do not grow in time in a suitable norm. For this we use what we call {\it{the Z-norm method}}: we define the norms
\begin{equation}\label{introdesc5}
\|f\|_{Z_{wa}}:=\sup_{k\in\mathbb{Z}}2^{N_0k^+}2^{k^-(1+\kappa)}\|\widehat{P_kf}\|_{L^\infty},\qquad \|f\|_{Z_{kg}}:=\sup_{k\in\mathbb{Z}}2^{N_0k^+}2^{k^-(1/2-\kappa)}\|\widehat{P_kf}\|_{L^\infty},
\end{equation}
where $x^+=\max(x,0)$ and $x^-=\min(x,0)$ for any $x\in\mathbb{R}$, $N_0=40$ as before, and $\kappa=10^{-3}$. Then we show that
\begin{equation}\label{introdesc6}
\|V^{F}(t)\|_{Z_{wa}}+\|V^{\omega_a}(t)\|_{Z_{wa}}+\|V^{\vartheta_{ab}}(t)\|_{Z_{wa}}+\|V^{\psi}(t)\|_{Z_{kg}}\lesssim\varep_0,
\end{equation}
for any $t\in[0,\infty)$ and $a,b\in\{1,2,3\}$, where the profiles $V^G$ are defined as in \eqref{introdesc},
\begin{equation}\label{introdesc7}
U^{G}(t):=\partial_t G(t)-i\Lambda_{wa}G(t),\qquad V^{G}(t):=e^{it\Lambda_{wa}}U^{G}(t),
\end{equation}
for $G\in\{F,\omega_a,\va_{ab}\}$. The components $\Omega_a$ can be recovered elliptically from $\omega_a$ and do not have to be considered separately. 

The choice of the $Z$-norm is very important, and one should think of it as analogous to the choice of the ``resolution norm" in the case of semilinear evolutions (the classical choices being Strichartz norms or $X^{s,b}$ norms). It has to complement well the information coming from energy estimates. The uniformity in time in \eqref{introdesc6} is the main point, in particular allowing us to prove sharp $\varep_0\langle t\rangle^{-1}$ pointwise decay on some components of the metric. 

The $Z$-norms defined in \eqref{introdesc5} measure the $L^\infty$ norm of solutions in the Fourier space, with weights that are particularly important at low frequencies. They cannot be propagated using energy estimates, since they are not $L^2$ based norms. We use instead the Duhamel formula, in the Fourier space, which leads to derivative loss. Because of this the $Z$-norm bounds \eqref{introdesc6} are weaker than the energy bounds \eqref{introdesc2} at very high frequencies. One should think of the $Z$-norm bounds as effective at middle frequencies, say $\langle t\rangle^{-1/2}\lesssim 2^k\lesssim \langle t\rangle^{1/2}$. 

The $Z$-norm method, with different choices of the norm itself, depending on the problem, was used recently by the authors and their collaborators in several small data global regularity problems, for water waves and plasmas, see \cite{DeIoPa, DeIoPaPu, GuIoPa, IoPa1, IoPa2, IoPa3, IoPu, IoPu3}.

\subsection{Modified scattering}\label{hus6} There is an additional difficulty in our problem, namely the solutions do not scatter linearly as $t\to\infty$. This is due to the low frequencies of the metric tensor in the quasilinear terms $g_{\geq 1}^{\mu\nu}\partial_\mu\partial_\nu h_{\al\be}$ and $g^{\mu\nu}_{\geq 1}\partial_\mu\partial_\nu\psi$, which create a long-range perturbation. 

To prove the uniform $Z$-norm bounds \eqref{introdesc6}, we need to renormalize the profiles. More precisely, we define the wave phase correction (related to optical functions)
\begin{equation}\label{introdesc10}
\begin{split}
\Theta_{wa}(\xi,t):=\int_0^t\Big\{&h_{00}^{low}(s\xi/\Lambda_{wa}(\xi),s)\frac{\Lambda_{wa}(\xi)}{2}\\
&+h_{0j}^{low}(s\xi/\Lambda_{wa}(\xi),s)\xi_j+h_{jk}^{low}(s\xi/\Lambda_{wa}(\xi),s)\frac{\xi_j\xi_k}{2\Lambda_{wa}(\xi)}\Big\}\,ds
\end{split}
\end{equation}
and the Klein-Gordon phase correction
\begin{equation}\label{introdesc11}
\begin{split}
\Theta_{kg}(\xi,t):=\int_0^t\Big\{&h_{00}^{low}(s\xi/\Lambda_{kg}(\xi),s)\frac{\Lambda_{kg}(\xi)}{2}\\
&+h_{0j}^{low}(s\xi/\Lambda_{kg}(\xi),s)\xi_j+h_{jk}^{low}(s\xi/\Lambda_{kg}(\xi),s)\frac{\xi_j\xi_k}{2\Lambda_{kg}(\xi)}\Big\}\,ds,
\end{split}
\end{equation}
where $h^{low}_{\al\be}$ are low frequency components of the metric tensor,
\begin{equation}\label{introdesc12}
\widehat{h^{low}_{\al\be}}(\rho,s):=\varphi_{\leq 0}(\langle s\rangle^{p_0}\rho)\widehat{h_{\al\be}}(\rho,s),\qquad p_0:=0.68.
\end{equation}  
The choice of $p_0$, slightly bigger than $2/3$, is important in the proof to justify the correction. Geometrically, the two phase corrections $\Theta_{wa}$ and $\Theta_{kg}$ are obtained by integrating suitable low frequency components of the metric tensor along the characteristics of the wave and the Klein-Gordon linear flows. 

The nonlinear profiles are then obtained by multiplication in the Fourier space,
\begin{equation}\label{introdesc13}
\widehat{V^{G}_\ast}(\xi,t):=e^{-i\Theta_{wa}(\xi,t)}\widehat{V^{G}}(\xi,t),\qquad \widehat{V^{\psi}_\ast}(\xi,t):=e^{-i\Theta_{kg}(\xi,t)}\widehat{V^{\psi}}(\xi,t),
\end{equation}
for $G\in\{F,\omega_a,\va_{ab}\}$. Notice that $\|V^{G}_\ast\|_{Z_{wa}}=\|V^{G}\|_{Z_{wa}}$ and $\|V^{\psi}_\ast\|_{Z_{kg}}=\|V^{\psi}\|_{Z_{kg}}$, since the phases $\Theta_{wa}$ and $\Theta_{kg}$ are real-valued. The point of this construction is that the new nonlinear profiles $V^{F}_\ast$, $V^{\omega_a}_\ast$, $V^{\vartheta_{ab}}_\ast$, and $V^{\psi}_\ast$ converge as the time goes to infinity, i.e.
\begin{equation}\label{introdesc14}
\begin{split}
\|V^F_\ast(t)-V^F_\infty\|_{Z_{wa}}+\|V^{\omega_a}_\ast(t)-V^{\omega_a}_\infty\|_{Z_{wa}}+\|V^{\vartheta_{ab}}_\ast(t)-V^{\vartheta_{ab}}_\infty\|_{Z_{wa}}&\lesssim\varep_0\langle t\rangle^{-\delta/2},\\
\|V^\psi_\ast(t)-V^\psi_\infty\|_{Z_{kg}}&\lesssim\varep_0\langle t\rangle^{-\delta/2},
\end{split}
\end{equation}
where $V^{F}_\infty, V^{\omega_a}_\infty, V^{\vartheta_{ab}}_\infty\in Z_{wa}$ and $V^{\psi}_\infty\in Z_{kg}$ are the {\it{nonlinear scattering data}}. These functions, in particular the components $V^{\vartheta_{ab}}_\infty$ and $V^\psi_\infty$ are important in the asymptotic analysis of our space-time. Chapter 5 is mainly concerned with the proofs of the bounds \eqref{introdesc14}.

\subsection{Asymptotic bounds and causal geodesics}\label{hus7} Our core bootstrap argument relies on controlling the solution both in the physical space and in the Fourier space, as summarized above. However, after closing the main bootstrap argument, we can derive classical bounds on the solutions in the physical space, without explicit use of the Fourier transform.

We start with decay estimates in the physical space. Let 
\begin{equation}\label{binr20}
L:=\partial_t+\partial_r,\qquad\underline{L}:=\partial_t-\partial_r,
\end{equation}
where $r:=|x|$ and $\partial_r:=|x|^{-1}x^j\partial_j$. Let
\begin{equation}\label{binr20.4}
\mathcal{T}:=\{L,r^{-1}\Omega_{12},r^{-1}\Omega_{23},r^{-1}\Omega_{31}\}
\end{equation}
denote the set of ``good" vector-fields, tangential to the (Minkowski) light cones.

In Theorem \ref{Precisedmetric} we prove that the metric components satisfy the bounds 
\begin{equation}\label{binr25.7}
\vert h(x,t)\vert+\langle t+r\rangle\vert \partial_Vh(x,t)\vert+\langle t-r\rangle\vert\partial_{\underline L}h(x,t)\vert\lesssim \varepsilon_0\langle t+r\rangle^{2\delta'-1},
\end{equation}
in the manifold $M:=\{(x,t)\in\mathbb{R}^3\times[0,\infty)\}$, where $r=|x|$, $V\in\mathcal{T}$, $h\in\{h_{\al\be}\}$, $\partial_W:=W^\al\partial_\al$, and $\delta'=2000\delta$. The scalar field decays faster but with no derivative improvement
\begin{equation}\label{binr25.4}
\begin{split}
\vert \psi(x,t)\vert+\vert\partial_0\psi(x,t)\vert&\lesssim\varepsilon_0\langle t+r\rangle^{\delta'/2-1}\langle r\rangle^{-1/2},\\
\vert \partial_b\psi(x,t)\vert&\lesssim\varepsilon_0\langle t+r\rangle^{\delta'/2-3/2},\qquad b\in\{1,2,3\}.
\end{split}
\end{equation}
Also, in Lemma \ref{twoder0} we show that the second order derivatives to the metric satisfy the bounds
\begin{equation}\label{binr25.8}
\langle r\rangle^2\vert \partial_{V_1}\partial_{V_2}h(x,t)\vert+\langle t-r\rangle^2\vert\partial_{\underline L}^2h(x,t)\vert+\langle t-r\rangle\langle r\rangle\vert\partial_{\underline L}\partial_{V_1}h(x,t)\vert\lesssim \varepsilon_0\langle r\rangle^{3\delta'-1},
\end{equation}
in the region $M':=\{(x,t)\in M:\,t\geq 1,\,|x|\geq 2^{-8}t\}$, where $V_1,V_2\in\mathcal{T}$ are good vector-fields.

The pointwise bounds \eqref{binr25.7}--\eqref{binr25.8} are as expected, including the small $\delta'$ losses which are due to our weak assumptions \eqref{approxas} on the initial data. These bounds follow mainly from the profile bounds \eqref{introdesc4} and linear estimates. 

As an application, we can describe precisely the future-directed causal geo\-de\-sics in our space-time $M$. Indeed, in Theorem \ref{winr100} we show that if $p=(p^0,p^1,p^2,p^3)$ is a point in $M$ and $v=v^\al\partial_\al$ is a null or time-like vector at $p$, normalized with $v^0=1$, then there is a unique affinely parametrized global geodesic curve $\gamma:[0,\infty)\to M$ with
\begin{equation*}
\gamma(0)=p=(p^0,p^1,p^2,p^3),\qquad\dot{\gamma}(0)=v=(v^0,v^1,v^2,v^3).
\end{equation*}
Moreover, the geodesic curve $\gamma$ becomes asymptotically parallel to a geodesic line of the Minkowski space, i.e. there is a vector $v_\infty=(v^0_\infty,v^1_\infty,v^2_\infty,v^3_\infty)$ such that, for any $s\in[0,\infty)$,
\begin{equation*}
|\dot{\gamma}(s)-v_\infty|\lesssim\varep_0(1+s)^{-1+6\delta'}\quad\text{ and }\quad |\gamma(s)-v_\infty s-p|\lesssim\varep_0(1+s)^{6\delta'}.
\end{equation*}

\subsection{Weak peeling estimates}\label{hus8} These are classical estimates on asymptotically flat space-times, which assert, essentially, that certain components of the Riemann curvature tensor have improved decay compared to the general estimate $|{\bf R}|\lesssim \varepsilon_0\langle t+r\rangle^{-1+}\langle t-r\rangle^{-2}$. The rate of decay is mainly determined by the {\it{signature}} of the component. 

More precisely, we use the Minkowski frames $(L,\underline{L},e_a)$, where $L,\underline{L}$ are as in \eqref{binr20} and $e_a\in\mathcal{T}_h:=\{r^{-1}\Omega_{12},r^{-1}\Omega_{23},r^{-1}\Omega_{31}\}$ and assign signature $+1$ to the vector-field $L$, $-1$ to the vector-field $\underline{L}$, and $0$ to the horizontal vector-fields in $\mathcal{T}_h$. With $e_1,e_2,e_3,e_4\in\mathcal{T}_h$, we define $\mathrm{Sig}(a)$ as the set of components of the Riemann tensor of total signature $a$, so
\begin{equation}\label{Qpeel3}
\begin{split}
\mathrm{Sig}(-2)&:=\{\mathbf{R}(\underline{L},e_1,\underline{L},e_2)\},\\
\mathrm{Sig}(2)&:=\{\mathbf{R}(L,e_1,L,e_2)\},\\
\mathrm{Sig}(-1)&:=\{\mathbf{R}(\underline{L},e_1,e_2,e_3),\mathbf{R}(\underline{L},L,\underline{L},e_1)\},\\
\mathrm{Sig}(1)&:=\{\mathbf{R}(L,e_1,e_2,e_3),\mathbf{R}(L,\underline{L},L,e_1)\},\\
\mathrm{Sig}(0)&:=\{\mathbf{R}(e_1,e_2,e_3,e_4),\mathbf{R}(L,\underline{L},e_1,e_2),\mathbf{R}(L,e_1,\underline{L},e_2),\mathbf{R}(L,\underline{L},L,\underline{L})\}.
\end{split}
\end{equation}
These components capture the entire curvature tensor, due to the symmetries \eqref{Riem1}. 

In Theorem \ref{PeelingR} we prove that if $\Psi_{(a)}\in\mathrm{Sig}(a)$, $a\in\{-2,-1,1,2\}$ then
\begin{equation}\label{Qpeel4}
\begin{split}
|\Psi_{(-2)}(x,t)|&\lesssim \varep_0\langle r\rangle^{7\delta'-1}\langle t-r\rangle^{-2},\\
|\Psi_{(-1)}(x,t)|&\lesssim \varep_0\langle r\rangle^{7\delta'-2}\langle t-r\rangle^{-1},\\
|\Psi_{(2)}(x,t)|+|\Psi_{(1)}(x,t)|+|\Psi_{(0)}(x,t)|&\lesssim \varep_0\langle r\rangle^{7\delta'-3},
\end{split}
\end{equation}
in the region $M'=\{(x,t)\in M:\,t\geq 1\text{ and }|x|\geq 2^{-8}t\}$. This holds in all cases except if $\Psi_{(0)}$ is of the form $\mathbf{R}(L,e_1,\underline{L},e_2)\in\mathrm{Sig}(0)$, in which case we can only prove the weaker bounds
\begin{equation}\label{Qpeel5}
|\mathbf{R}(L,e_1,\underline{L},e_2)(x,t)|\lesssim \varep_0\langle r\rangle^{7\delta'-2}\langle t-r\rangle^{-1}.
\end{equation}

Notice that we define our decomposition in terms of the Minkowski null pair $(L,\underline{L})$ instead of more canonical null frames (or tetrads) adapted to the metric $\g$ (see, for example, \cite{ChKl}, \cite{KlNi}, \cite{KlNi2}). This is not important however, since the weak peeling estimates are invariant under natural changes of the frame of the form $(L,\underline{L},e_a)\to (L',\underline{L}',e'_a)$, satisfying 
\begin{equation*}
|(L-L')(x,t)|+|(\underline{L}-\underline{L}')(x,t)|+|(e_a-e'_a)(x,t)|\lesssim r^{-1+2\delta'}\qquad\text{ in }M'.
\end{equation*}
As we show in Proposition \ref{PeelingR2}, one can in fact restore the full $\varep_0\langle r\rangle^{7\delta'-3}$ decay of the component $\mathbf{R}(L',e'_1,\underline{L}',e'_2)$, provided that $L'$ is almost null, i.e. $|\g(L',L')(x,t)|\lesssim \langle r\rangle^{-2+4\delta'}$ in $M'$.

The almost cubic decay we prove in \eqref{Qpeel4}--\eqref{Qpeel5} seems optimal in our problem, for two reasons. First, the Ricci components themselves involve squares of the massive field, and cannot decay better than $\langle r\rangle^{-3+}$ in $M'$. Moreover, the almost cubic decay is also formally consistent with the weak peeling estimates of Klainerman--Nicolo \cite[Theorem 1.2 (b)]{KlNi2} in the setting of our more general metrics (one would need to formally take $\gamma=-1/2-$ and $\delta=2+$ with the notation in \cite{KlNi2}, to match our decay assumptions \eqref{approxas} on the initial data; this range of parameters is not allowed however in \cite{KlNi2} as $\delta$ is assumed to be $<3/2$).

\subsection{The ADM energy and the linear momentum}\label{hus9} The ADM energy (or the ADM mass) measures the total deviation of our space-time from the Minkowski solution. It is calculated according to the standard formula (see for example \cite{Ba})
\begin{equation}\label{Imass1}
E_{ADM}(t):=\frac{1}{16\pi}\lim_{R\to\infty}\int_{S_{R,t}}(\partial_j\g_{nj}-\partial_n\g_{jj})\frac{x^n}{|x|}\,dx,
\end{equation}
where the integration is over large (Euclidean) spheres $S_{R,t}\subseteq\Sigma_t=\{(x,t):\,x\in\mathbb{R}^3\}$ of radius $R$. In our case we show in Proposition \ref{ADM4} that the energy $E_{ADM}(t)=E_{ADM}$ is well defined and constant in time. Moreover, it is non-negative and can be expressed in terms of the scattering profiles $V^\psi_\infty$ and $V^{\vartheta_{mn}}_\infty$ (see \eqref{introdesc14}) according to the formula
\begin{equation}\label{Imass2}
E_{ADM}=\frac{1}{16\pi}\|V^\psi_\infty\|_{L^2}^2+\frac{1}{64\pi}\sum_{m,n\in\{1,2,3\}}\|V^{\vartheta_{mn}}_\infty\|_{L^2}^2.
\end{equation}

We can also prove conservation of one other natural quantity, namely the linear momentum. Let $N$ denote the future unit normal vector-field to the hypersurface $\Sigma_t$, let $\overline{g}_{ab}=\mathbf{g}_{ab}$ denote the induced (Riemannian) metric on $\Sigma_t$, and define the second fundamental form 
\begin{equation*}
k_{ab}:=-\g(\D_{\partial_a}N,\partial_b)=\g(N,\D_{\partial_a}\partial_b)=N^\al\mathbf{\Gamma}_{\al a b},\qquad a,b\in\{1,2,3\}.
\end{equation*} 
Then we define the linear momentum $\mathbf{p}_a$, $a\in\{1,2,3\}$,
\begin{equation*}
\mathbf{p}_a(t):=\frac{1}{8\pi}\lim_{R\to\infty}\int_{S_{R,t}}\pi_{ab}\frac{x^b}{|x|}\,dx,\qquad \pi_{ab}:=k_{ab}-\mathrm{tr}k\overline{g}_{ab},
\end{equation*}
In Proposition \ref{linmo2+} we prove that the functions $\mathbf{p}_a$ are well defined and constant in time. Moreover, we show that  $\sum_{a\in\{1,2,3\}}\mathbf{p}_a^2\leq E_{ADM}^2$, so the ADM mass $M_{ADM}:=\big(E_{ADM}^2-\sum_{a\in\{1,2,3\}}\mathbf{p}_a^2\big)^{1/2}\geq 0$ is well defined.

We remark that the momentum $\mathbf{p}_a$ vanishes in the case of metrics $\g$ that agree with the Schwarzschild metric (including time derivatives) up to lower order terms. In particular, it vanishes in the case of metrics considered in earlier work on the stability for the Einstein vacuum equations, such as \cite{ChKl, KlNi, LiRo}. However, in our non-isotropic case the linear momentum does not necessarily vanish, and the quantities $\mathbf{p}_a$ defined above are natural conserved quantities of the evolution.

\subsection{The Bondi energy}\label{BonIntro} To define a Bondi energy we have to be more careful. We would like to compute integrals over large spheres as in \eqref{Imass1}, and then take the limit along outgoing null cones towards null infinity. But the limit exists only if we account properly for the geometry of the problem.  

First we need to understand the bending of the light cones caused by the long-range effect of the nonlinearity (i.e. the modified scattering). For this we construct (in Lemma \ref{AlmostOpticalLem}) an {\it{almost optical function}} $u:M'\to\mathbb{R}$, satisfying the properties
\begin{equation}\label{Imass3}
u(x,t)=\vert x\vert-t+u^{cor}(x,t),\qquad \g^{\alpha\beta}\partial_\alpha u\partial_\beta u=O(\varepsilon_0\langle r\rangle^{-2+6\delta'}).
\end{equation}
In addition, the correction $u^{cor}=O(\varepsilon_0\langle r\rangle^{3\delta'})$ is close to $\Theta_{wa}/|x|$ (see \eqref{introdesc10}) near the light cone,
\begin{equation}\label{Imass5}
\begin{split}
\Big\vert u^{cor}(x,t)-\frac{\Theta_{wa}(x,t)}{|x|}\Big\vert\lesssim \varep_0\langle r\rangle^{-1+3\delta'}(\langle r\rangle^{0.68}+\langle t-\vert x\vert \rangle),
\end{split}
\end{equation}
if $(x,t)\in M',\,\big|t-|x|\big|\leq t/10$. Notice that we work with an approximate optical condition $\g^{\alpha\beta}\partial_\alpha u\partial_\beta u=O(\varepsilon_0\langle r\rangle^{-2+6\delta'})$ instead of the classical optical condition $\g^{\alpha\beta}\partial_\alpha u\partial_\beta u=0$. This is mostly for convenience, since the weaker condition is still good enough for our analysis and almost optical functions are much easier to construct than exact optical functions.

For any $t\geq 1$ we define the hypersurface $\Sigma_t:=\{(x,t)\in M:\,x\in\mathbb{R}^3\}$, and let $\overline{g}_{jk}=\mathbf{g}_{jk}$ denote the induced (Riemannian) metric on $\Sigma_t$. With $u$ as above, we define the modified spheres $S^u_{R,t}:=\{x\in\Sigma_t:\,u(x,t)=R\}$ and let $\mathbf{n}_j:=\partial_ju(\overline{g}^{ab}\partial_au\partial_bu)^{-1/2}$ 
denote the unit vector-field normal to the spheres $S^u_{R,t}$. For $R\in\mathbb{R}$ and $t$ large (say $t\geq 2|R|+10$) we define
\begin{equation}\label{Imass7}
E_{Bondi}(R):=\frac{1}{16\pi}\lim_{t\to\infty}\int_{S^u_{R,t}} \overg^{ab}(\partial_a h_{jb}-\partial_j h_{ab})\mathbf{n}^j\,d\sigma,
\end{equation}
where $d\sigma= d\sigma(\overline{g})$ is the surface measure induced by the metric $\overline{g}$. Notice that this definition is a more geometric version of the definition \eqref{Imass1}, in the sense that the integration is with respect to the metric $\overg$. Geometrically, we fix $R$ and integrate on surfaces $S^u_{R,t}$ which live on the ``light cone" $\{u(x,t)=R\}$ 

In Theorem \ref{BondThm} we prove our main result: the limit in \eqref{Imass7} exists, and $E_{Bondi}:\mathbb{R}\to \mathbb{R}$ is a well-defined increasing and continuous function on $\mathbb{R}$, which increases from the Klein-Gordon energy $E_{KG}$ to the full ADM energy $E_{ADM}$, i.e.
\begin{equation}\label{Imass8}
\lim_{R\to-\infty}E_{Bondi}(R)=E_{KG}:=\frac{1}{16\pi}\|V^\psi_\infty\|_{L^2}^2,\qquad\lim_{R\to\infty}E_{Bondi}(R)=E_{ADM}.
\end{equation}

The definition \eqref{Imass7} of the Bondi energy is consistent with the general heuristics in \cite[Chapter 11]{Wa} and with the definition in \cite[Section 4.3.4]{Poi}. It also has expected properties, like monotonicity, continuity, and satisfies the limits \eqref{Imass8}. 

However, it is not clear to us if this definition is identical to the definition used by Klainerman-Nicolo \cite[Chapter 8.5]{KlNi}, starting from the Hawking mass. In fact, at the level of generality of our metrics \eqref{approxas}, it not even clear that one can prove sharp $r^{-3}$ pointwise decay on some of the signature $0$ components of the curvature tensor, which is one of the ingredients of the argument in \cite{KlNi}. 

We notice that the Klein-Gordon energy $E_{KG}$ is part of $E_{Bondi}(R)$, for all $R\in\mathbb{R}$. This is consistent with the geometric intuition, since the matter travels at speeds lower than the speed of light and accumulates at timelike infinity, not at null infinity. We can further measure its radiation by taking limits along timelike cones. Indeed, for $\alpha\in (0,1)$ let 
\begin{equation}\label{Imass9}
E_{i^+}(\alpha):=\frac{1}{16\pi}\lim_{t\to\infty}\int_{S_{\alpha t,t}}(\partial_j h_{nj}-\partial_n h_{jj})\frac{x^n}{|x|}\,dx,
\end{equation}
where the integration is over the Euclidean spheres $S_{\alpha t,t}\subseteq\Sigma_t$ of radius $\alpha t$. In Proposition \ref{Bond40} we prove that the limit in \eqref{Imass9} exists, and $E_{i^+}:(0,1)\to\mathbb{R}$ is a well-defined continuous and increasing function, satisfying
\begin{equation}\label{Imass10}
\lim_{\alpha\to 0}E_{i^+}(\alpha)=0,\qquad \lim_{\alpha\to 1}E_{i^+}(\alpha)=E_{KG}.
\end{equation}

\subsection{Organization} The rest of this monograph is organized as follows:

In Chapter 2 we introduce our main notations and definitions and describe precisely our bootstrap construction (Proposition \ref{bootstrap}). Then we prove several important lemmas that are being used in the rest of the analysis, such as Lemmas \ref{pha2} and \ref{PhaWave} on the structure and bounds on quadratic resonances, Lemma \ref{LinEstLem} concerning linear estimates for wave and Klein-Gordon evolutions, and Lemmas \ref{Lembil1}--\ref{LemBil4} concerning bilinear estimates. Finally, we use these lemmas and the bootstrap hypothesis to prove linear estimates on the solutions and the profiles, such as Lemmas \ref{dtv2} (localized $L^2$ bounds) and Lemma \ref{dtv3} (pointwise decay).

In Chapter 3 we analyze our main nonlinearities $\LL\mathcal{N}^h_{\al\be}$ and $\LL\mathcal{N}^\psi$ at a fixed time $t$. The main results in this chapter are Proposition \ref{plk1} (localized $L^2$, $L^\infty$, and weighted $L^2$ bounds on these nonlinearities), Lemmas \ref{nake12}--\ref{nake12.5} (identification of the energy disposable nonlinear components), and especially 
Proposition \ref{gooddec} (decomposition of the main nonlinearities). 

In Chapter 4 we prove the main bootstrap bounds \eqref{bootstrap3.1} on the energy functionals. We start from the decomposition in Proposition \ref{gooddec}, perform energy estimates, and prove bounds on all the resulting space-time integrals. The main space-time bounds are  stated in Proposition \ref{totbounds}, and are proved in the rest of the chapter, using normal forms, null structures, and paradifferential calculus in some of the harder cases.

In Chapter 5 we first prove the main bootstrap bounds \eqref{bootstrap3.2} (weighted estimates on profiles) in Proposition \ref{DiEs1}, as a consequence of the improved energy estimates and the nonlinear bounds in Proposition \ref{plk1}. Then we prove the main bootstrap bounds \eqref{bootstrap3.4} (the $Z$-norm estimates). This proof has several steps, such as the renormalization procedures in \eqref{Nor4}-\eqref{Nor5} and \eqref{bax3}--\eqref{bax4}, and the estimates \eqref{Nor40} and \eqref{bax12} showing boundedness and convergence of the nonlinear profiles in suitable norms in the Fourier space.

In Chapter 6 we prove a full, quantitative version of our main global regularity result (Theorem \ref{winr}) as well as all the other consequences on the asymptotic structure of our space-times, as described in detail in subsections \ref{hus7}--\ref{BonIntro} above.

\subsection{Acknowledgements} The first author was supported in part by NSF grant DMS-1600028 and NSF-FRG grant DMS-1463753. The second author was supported in part by NSF grant DMS-1362940 and by a Sloan Research fellowship.

\chapter{The main construction and preliminary estimates}

\section{Setup and the main bootstrap proposition}\label{prelset}

In this section we introduce most of our notations and definitions and state our main bootstrap proposition.

\subsection{The nonlinearities $\mathcal{N}^h_{\al\be}$ and $\mathcal{N}^\psi$}\label{decomp} Let $m$ denote the Min\-kow\-ski metric and write
\begin{equation}\label{zaq1}
\g_{\al\be}=m_{\al\be}+h_{\al\be},\qquad \g^{\al\be}=m^{\al\be}+g_{\geq 1}^{\al\be}.
\end{equation}
We start by rewriting our system as a wave-Klein-Gordon coupled quasilinear system:

\begin{proposition}\label{LemEKG}
Assume $(g,\psi)$ in a solution in $\mathbb{R}^3\times[0,T]$ of the reduced Einstein-Klein-Gordon system in harmonic gauge \eqref{asum3}-\eqref{asum5}. For  $\al,\be\in\{0,1,2,3\}$ we have
\begin{equation}\label{zaq11.1}
(\partial_0^2-\Delta)h_{\al\be}=\mathcal{N}^{h}_{\al\be}:=\mathcal{KG}_{\al\be}+g_{\geq 1}^{\mu\nu}\partial_\mu\partial_\nu h_{\al\be}-F^{\geq 2}_{\al\be}(g,\partial g)
\end{equation}
where 
\begin{equation}\label{zaq11.12}
\mathcal{KG}_{\al\be}:=2\partial_\al\psi\partial_\be\psi+\psi^2(m_{\al\be}+h_{\al\be}).
\end{equation}
Moreover
\begin{equation}\label{zaq11.3}
(\partial_0^2-\Delta+1)\psi=\mathcal{N}^{\psi}:=g^{\mu\nu}_{\geq 1}\partial_\mu\partial_\nu\psi.
\end{equation}
In addition, the nonlinearities $F^{\geq 2}_{\al\be}(g,\partial g)$ admit the decompositions
\begin{equation}\label{tr10}
F^{\geq 2}_{\al\be}(g,\partial g)=Q_{\al\be}+P_{\al\be},
\end{equation}
where
\begin{equation}\label{tr11}
\begin{split}
Q_{\al\be}:&=\g^{\rho\rho'}\g^{\lambda\la'}(\partial_\al h_{\rho'\la'}\partial_{\rho} h_{\be\la}-\partial_{\rho} h_{\rho'\lambda'}\partial_\al h_{\be\la})+\g^{\rho\rho'}\g^{\lambda\la'}(\partial_\be h_{\rho'\la'}\partial_{\rho} h_{\al\la}-\partial_{\rho} h_{\rho'\lambda'}\partial_\be h_{\al\la})\\
&+\frac{1}{2}\g^{\rho\rho'}\g^{\lambda\la'}(\partial_{\la'} h_{\rho\rho'}\partial_\be h_{\al\la}-\partial_\be h_{\rho\rho'}\partial_{\la'} h_{\al\la})+\frac{1}{2}\g^{\rho\rho'}\g^{\lambda\la'}(\partial_{\la'} h_{\rho\rho'}\partial_\al h_{\be\la}-\partial_\al h_{\rho\rho'}\partial_{\la'} h_{\be\la})\\
&-\g^{\rho\rho'}\g^{\lambda\la'}(\partial_{\la} h_{\al\rho'}\partial_\rho h_{\be\la'}-\partial_\rho h_{\al\rho'}\partial_{\la} h_{\be\la'})+\g^{\rho\rho'}\g^{\lambda\la'}\partial_{\rho'} h_{\al\la'}\partial_\rho h_{\be\la}
\end{split}
\end{equation}
and
\begin{equation}\label{tr12}
P_{\al\be}:=-\frac{1}{2}\g^{\rho\rho'}\g^{\lambda\lambda'}\partial_\al h_{\rho'\la'}\partial_\be h_{\rho\la}+\frac{1}{4}\g^{\rho\rho'}\g^{\lambda\la'}\partial_\al h_{\rho\rho'}\partial_\be h_{\la\la'}.
\end{equation}
\end{proposition}

\begin{proof} The identities \eqref{zaq11.1} and \eqref{zaq11.3} follow directly from the system \eqref{asum3}. The identities \eqref{tr10}, which allow us to extract the null components of the semilinear nonlinearities, follow by explicit calculations from the identities \eqref{tr9.5}, see for example \cite[Lemma 4.1]{LeMaEKG} and notice that $\partial_\rho h_{\mu\nu}=\partial_\rho g_{\mu\nu}$). 
\end{proof}

\subsection{Decomposition of the metric tensor} Let $R_j=|\nabla|^{-1}\partial_j$, $j\in\{1,2,3\}$, denote the Riesz transforms on $\mathbb{R}^3$, and notice that $\delta_{jk}R_jR_k=-I$. To identify null structures we use a double Hodge decomposition for the metric tensor, which is connected to the work of Arnowitt--Deser--Misner \cite{ADM} on the Hamiltonian formulation of General Relativity.  Let 
\begin{equation}\label{zaq2}
\begin{split}
&F:=(1/2)[h_{00}+R_jR_kh_{jk}],\qquad\uF:=(1/2)[h_{00}-R_jR_kh_{jk}],\\
&\rho:=R_jh_{0j},\qquad\omega_j:=\in_{jkl}R_kh_{0l},\\
&\Omega_j:=\in_{jkl}R_kR_mh_{lm},\qquad\vartheta_{jk}:=\in_{jmp}\in_{knq}R_mR_nh_{pq}.
\end{split}
\end{equation} 
Notice that $\omega$ and $\Omega$ are divergence-free vector-fields,
\begin{equation}\label{zaq3}
R_j\omega_j=0,\qquad R_j\Omega_j=0,
\end{equation}
and $\vartheta$ is a symmetric and divergence-free tensor-field,
\begin{equation}\label{zaq4}
\vartheta_{jk}=\vartheta_{kj},\qquad R_j\vartheta_{jk}=0,\qquad R_k\vartheta_{jk}=0.
\end{equation}
Moreover, using the general formula $\in_{mnk}\in_{pqk}=\delta_{mp}\delta_{nq}-\delta_{mq}\delta_{np}$ one can recover the tensor $h$ according to the formulas
\begin{equation}\label{zaq5}
\begin{split}
&h_{00}=F+\uF,\\
&h_{0j}=-R_j\rho+\in_{jkl}R_k\omega_l,\\
&h_{jk}=R_jR_k(F-\uF)-(\in_{klm}R_j+\in_{jlm}R_k)R_l\Omega_m+\in_{jpm}\in_{kqn}R_pR_q\vartheta_{mn}.
\end{split}
\end{equation}

As a consequence of the harmonic gauge condition, we will show in Lemma \ref{LemWaveGauge} that the main dynamical variables are $F,\uF,\omega_j$ and the traceless part of $ \vartheta_{jk}$, while the variables $\rho$ and $\Omega_j$, can be expressed elliptically in terms of these variables, up to quadratic remainders. This is important in identifying suitable null structures of the nonlinearities $\mathcal{N}^{h}_{\al\be}$ in section \ref{DecompositionNonlin}.

We often need to extract the linear part of the matrix $g_{\geq 1}^{\al\be}$. So we write $g_{\geq 1}^{\al\be}=g_1^{\al\be}+g_{\geq 2}^{\al\be}$ and use the identity 
\begin{equation*}
\begin{split}
\delta^\al_\be&=\g^{\al\rho}\g_{\be\rho}=(m^{\al\rho}+g_1^{\al\rho}+g_{\geq 2}^{\al\rho})(m_{\be\rho}+h_{\be\rho})\\
&=\delta^\al_\be+(m^{\al\rho}h_{\be\rho}+g_1^{\al\rho}m_{\be\rho})+(g_1^{\al\rho}h_{\be\rho}+g_{\geq 2}^{\al\rho}m_{\be\rho}+g_{\geq 2}^{\al\rho}h_{\be\rho}).
\end{split}
\end{equation*} 
Therefore, we can define $g_1^{\al\be}$ and $g_{\geq 2}^{\al\be}$ by 
\begin{equation}\label{zaq22}
\begin{split}
&g_1^{00}=-h_{00},\qquad g_1^{0j}=g_1^{j0}=h_{0j},\qquad g_1^{jk}=-h_{jk},\\
&g_{\geq 2}^{\al\rho}m_{\be\rho}+g_{\geq 2}^{\al\rho}h_{\be\rho}+g_1^{\al\rho}h_{\be\rho}=0.
\end{split}
\end{equation}

\subsection{Littlewood-Paley projections} We fix $\varphi:\mathbb{R}\to[0,1]$ an even smooth 
function supported in $[-8/5,8/5]$ and equal to $1$ in $[-5/4,5/4]$. For simplicity of notation, we also 
let $\varphi:\mathbb{R}^3\to[0,1]$ denote the corresponding radial function on $\mathbb{R}^3$. Let
\begin{equation*}
\varphi_{k}(x):=\varphi(|x|/2^k)-\varphi(|x|/2^{k-1})\text{ for any }k\in\mathbb{Z},\qquad \varphi_I:=\sum_{m\in I\cap\mathbb{Z}}\varphi_m\text{ for any }I\subseteq\mathbb{R}.
\end{equation*}
For any $B\in\mathbb{R}$ let 
\begin{equation*}
\varphi_{\leq B}:=\varphi_{(-\infty,B]},\quad\varphi_{\geq B}:=\varphi_{[B,\infty)},\quad\varphi_{<B}:=\varphi_{(-\infty,B)},\quad \varphi_{>B}:=\varphi_{(B,\infty)}.
\end{equation*}
For any $a<b\in\mathbb{Z}$ and $j\in[a,b]\cap\mathbb{Z}$ let
\begin{equation}\label{Alx80}
\varphi^{[a,b]}_j:=
\begin{cases}
\varphi_{j}\quad&\text{ if }a<j<b,\\
\varphi_{\leq a}\quad&\text{ if }j=a,\\
\varphi_{\geq b}\quad&\text{ if }j=b.
\end{cases}
\end{equation}

Let $P_k$, $k\in\mathbb{Z}$, (respectively $P_I$, $I\subseteq\mathbb{R}$) denote the operators on $\mathbb{R}^3$ defined by the Fourier multipliers $\xi\to \varphi_k(\xi)$ (respectively $\xi\to \varphi_I(\xi)$). For simplicity of notation let $P'_k=P_{[k-2,k+2]}$. 

For any $x\in\mathbb{Z}$ let $x^{+}=\max(x,0)$ and $x^-:=\min(x,0)$. Let
\begin{equation*}
\mathcal{J}:=\{(k,j)\in\mathbb{Z}\times\mathbb{Z}_+:\,k+j\geq 0\}.
\end{equation*}
For any $(k,j)\in\mathcal{J}$ let
\begin{equation*}
\phii^{(k)}_j(x):=
\begin{cases}
\varphi_{\leq -k}(x)\quad&\text{ if }k+j=0\text{ and }k\leq 0,\\
\varphi_{\leq 0}(x)\quad&\text{ if }j=0\text{ and }k\geq 0,\\
\varphi_j(x)\quad&\text{ if }k+j\geq 1\text{ and }j\geq 1,
\end{cases}
\end{equation*}
and notice that, for any $k\in\mathbb{Z}$ fixed, $\sum_{j\geq-\min(k,0)}\phii^{(k)}_j=1$. 

For $(k,j)\in\mathcal{J}$ let $Q_{j,k}$ denote the operator
\begin{equation}\label{qjk}
(Q_{j,k}f)(x):=\phii^{(k)}_j(x)\cdot P_kf(x).
\end{equation}
In view of the uncertainty principle the operators $Q_{j,k}$ are relevant only when $2^j2^k\gtrsim 1$, which explains the definitions above.

We will often estimate bilinear interactions, like products of two functions. For $k\in\mathbb{Z}$ let
\begin{equation}\label{gb7}
\mathcal{X}_k:=\{(k_1,k_2)\in\mathbb{Z}^2:\,|\max(k,k_1,k_2)-\mathrm{med}(k,k_1,k_2)|\leq 4\}.
\end{equation}
Notice that $P_k(P_{k_1}f\cdot P_{k_2}g)\equiv 0$ if $(k_1,k_2)\notin\mathcal{X}_k$.

\subsection{Vector-fields} Recall the vector-fields $\Gamma_j$ and $\Omega_{jk}$ defined in \eqref{qaz2}. These vector-fields satisfy simple commutation relations, which can be written schematically in the form
\begin{equation}\label{qaz2.1}
\begin{split}
&[\partial,\partial]=0,\qquad [\partial,\Omega]=\partial,\qquad [\partial,\Gamma]=\partial,\\
&[\Omega,\Omega]=\Omega,\qquad [\Omega,\Gamma]=\Gamma, \qquad [\Gamma,\Gamma]=\Omega,
\end{split}
\end{equation}
where $\partial$ denote generic coordinate vector-fields, $\Omega$ denote generic rotation vector-fields, and $\Gamma$ denote generic Lorentz vector-fields. For any $\alpha=(\alpha_1,\alpha_2,\alpha_3)\in(\mathbb{Z}_+)^3$ we define
\begin{equation}\label{qaz3}
\partial^\alpha:=\partial_1^{\alpha_1}\partial_2^{\alpha_2}\partial_3^{\alpha_3},\qquad \Omega^\alpha:=\Omega_{23}^{\alpha_1}\Omega_{31}^{\alpha_2}\Omega_{12}^{\alpha_3},\qquad \Gamma^\alpha:=\Gamma_1^{\alpha_1}\Gamma_2^{\alpha_2}\Gamma_3^{\alpha_3}.
\end{equation}

For any $n,q\in\mathbb{Z}_+$ we define $\mathcal{V}_n^q$ as the set of differential operators of the form
\begin{equation}\label{qaz3.1}
\mathcal{V}^q_n:=\big\{\mathcal{L}=\Gamma^a\Omega^b:\,|a|+|b|\leq n,\,q(\mathcal{L}):=|a|\leq q\big\}.
\end{equation}
Here $q(\mathcal{L})$ denotes the number of vector-fields transversal to the surfaces $\Sigma_a:=\{(x,t)\in\mathbb{R}^3\times\mathbb{R}:\,t=a\}$. We remark that in our proof we distinguish between the Lorentz vector-fields $\Gamma$ (which lead to slightly faster growth rates, see the definition \eqref{fvc1.0}) and the rotational vector-fields $\Omega$. Notice that, for any $\al\in\{0,1,2,3\}$ and any $\mathcal{L}_1\in\mathcal{V}_{n_1}^{q_1}$, $\mathcal{L}_2\in\mathcal{V}_{n_2}^{q_2}$, we have
\begin{equation}\label{qaz3.5}
\begin{split}
\mathcal{L}_1\mathcal{L}_2&=\text{ sum of operators in }\mathcal{V}^{q_1+q_2}_{n_1+n_2},\\
[\partial_\al,\mathcal{L}_1]&=\text{ sum of operators of the form }\partial_\be\mathcal{L}'\text{ for some }\beta\in\{0,1,2,3\},\,\mathcal{L}'\in \mathcal{V}_{n_1-1}^{q_1}.
\end{split}
\end{equation}

\subsection{Linear profiles and the $Z$-norms} We define the {\it{normalized solutions}} $U^{h_{\al\be}}$, $U^F$, $U^{\underline{F}}$, $U^\rho$, $U^{\omega_a}$, $U^{\Omega_a}$, $U^{\vartheta_{ab}}$, $U^\psi$ and their associated {\it{linear profiles}}  $V^{h_{\al\be}}$, $V^F$, $V^{\underline{F}}$, $V^\rho$, $V^{\omega_a}$, $V^{\Omega_a}$, $V^{\vartheta_{ab}}$, $V^\psi$, $\al,\be\in\{0,1,2,3\}$, $a,b\in\{1,2,3\}$, by
\begin{equation}\label{variables4}
\begin{split}
&U^{G}(t):=\partial_t G(t)-i\Lambda_{wa} G(t),\qquad V^{G}(t):=e^{it\Lambda_{wa}}U^{G}(t),\qquad G\in\{h_{\al\be},F,\underline{F},\rho,\omega_a,\Omega_a,\vartheta_{ab}\},\\
&U^{\psi}(t):=\partial_t\psi(t)-i\Lambda_{kg} \psi(t),\qquad\,\,\,V^{\psi}(t):=e^{it\Lambda_{kg}}U^{\psi}(t),
\end{split}
\end{equation}
where 
\begin{equation}\label{qaz1}
\Lambda_{wa}:=|\nabla|,\qquad \Lambda_{kg}:=\langle\nabla\rangle=\sqrt{|\nabla|^2+1}.
\end{equation}
More generally, for $\mathcal{L}\in\mathcal{V}_3^3$ (see definition \eqref{qaz3.1}) we define the {\it{weighted linear profiles}}{\footnote{Notice that we only apply the differential operators $\mathcal{L}$ to the metric components $h_{\al\be}$, but not to the variables $F,\underline{F},\rho,\omega_a,\Omega_a,\vartheta_{ab}$. In fact, the operators $\mathcal{L}$ that contain only the rotation vector-fields can be applied to these variables without any difficulty, but the combination of several Lorentz vector-fields $\Gamma_a$ and Riesz transforms leads to unwanted contributions at low frequency.}}
\begin{equation}\label{variables4L}
\begin{split}
&U^{\mathcal{L}h_{\al\be}}(t):=(\partial_t-i\Lambda_{wa})(\mathcal{L} h_{\al\be})(t),\qquad V^{\mathcal{L}h_{\al\be}}(t):=e^{it\Lambda_{wa}}U^{\mathcal{L}h_{\al\be}}(t),\\
&U^{\mathcal{L}\psi}(t):=(\partial_t-i\Lambda_{kg})(\mathcal{L} \psi)(t),\qquad\,\,V^{\mathcal{L}\psi}(t):=e^{it\Lambda_{kg}}U^{\mathcal{L}\psi}(t).
\end{split}
\end{equation}
Also, for $\ast\in\{F,\underline{F},\rho,\omega_a,\Omega_a,\vartheta_{ab},\mathcal{L}h_{\al\be},\mathcal{L}\psi\}$, $\mathcal{L}\in\mathcal{V}_{3}^3$, we define
\begin{equation}\label{notationL}
U^{\ast,+}:=U^{\ast},\qquad U^{\ast,-}:=\overline{U^{\ast}},\qquad V^{\ast,+}:=V^{\ast},\qquad V^{\ast,-}:=\overline{V^{\ast}}.
\end{equation}

The functions $F,\underline{F},\rho,\omega_a,\Omega_a,\vartheta_{ab},\mathcal{L}h_{\al\be},\mathcal{L}\psi$ can be recovered linearly from the normalized variables $U^{F}$, $U^{\underline{F}}$, $U^{\rho}$, $U^{\omega_a}$, $U^{\Omega_a}$, $U^{\vartheta_{ab}}$, $U^{\mathcal{L}h_{\al\be}}$, $U^{\mathcal{L}\psi}$ by the formulas
\begin{equation}\label{on5}
\begin{split}
&\partial_0 G=(U^{G}+\overline{U^{G}})/2,\qquad \Lambda_{wa}G=i(U^{G}-\overline{U^{G}})/2,\qquad G\in\{F,\underline{F},\rho,\omega_a,\Omega_a,\vartheta_{ab},\mathcal{L}h_{\al\be}\},\\
&\partial_0 \mathcal{L}\psi=(U^{\mathcal{L}\psi}+\overline{U^{\mathcal{L}\psi}})/2,\qquad\,\, \Lambda_{kg}\mathcal{L}\psi=i(U^{\mathcal{L}\psi}-\overline{U^{\mathcal{L}\psi}})/2.
\end{split}
\end{equation}

Let
\begin{equation}\label{on9.1}
\mathcal{P}:=\{(wa,+),(wa,-),(kg,+),(kg,-)\}.
\end{equation}
Let $\Lambda_{wa,+}(\xi)=\Lambda_{wa}(\xi)=|\xi|$, $\Lambda_{wa,-}(\xi)=-\Lambda_{wa,+}(\xi)$, $\Lambda_{kg,+}(\xi)=\Lambda_{kg}(\xi)=\sqrt{|\xi|^2+1}$, $\Lambda_{kg,-}(\xi)=-\Lambda_{kg,+}(\xi)$. For any $\sigma,\mu,\nu\in\mathcal{P}$ we define the {\it{quadratic phase function}}
\begin{equation}\label{on9.2}
\Phi_{\sigma\mu\nu}:\mathbb{R}^3\times\mathbb{R}^3\to\mathbb{R},\qquad \Phi_{\sigma\mu\nu}(\xi,\eta):=\Lambda_{\sigma}(\xi)-\Lambda_{\mu}(\xi-\eta)-\Lambda_\nu(\eta).
\end{equation}

In our analysis we will need a few parameters:
\begin{equation}\label{parameters}
\begin{split}
&N_0:=40,\quad d:=10,\quad\kappa:=10^{-3},\quad \delta:=10^{-10},\quad \delta':=2000\delta,\quad\gamma:=\delta/4.
\end{split}
\end{equation}
We define also the numbers $N(n)$ (which measure the number of Sobolev derivatives under control at the level of $n$ vector-fields),
\begin{equation}\label{fvc1.1}
N(0):=N_0+16d,\qquad N(n):=N_0-dn\,\text{ for }\,n\in\{1,2,3\}.
\end{equation}
Let $|\xi|_{\leq 1}$ denote a smooth increasing radial function on $\mathbb{R}^3$ equal to $|\xi|$ if $|\xi|\leq 1/2$ and equal to $1$ if $|\xi|\geq 2$. Let $|\nabla|_{\leq 1}^\gamma$ denote the associated operator defined by the multiplier $\xi\to |\xi|_{\leq 1}^\gamma$.

We are now ready to define the main $Z$-norms.

\begin{definition}\label{MainZDef}
For any $x\in\mathbb{R}$ let $x^+=\max(x,0)$ and $x^-=\min(x,0)$. We define the spaces $Z_{wa}$ and $Z_{kg}$ by the norms
\begin{equation}\label{sec5}
\|f\|_{Z_{wa}}:=\sup_{k\in\mathbb{Z}}2^{N_0k^+}2^{k^-(1+\kappa)}\|\widehat{P_kf}\|_{L^\infty}
\end{equation}
and
\begin{equation}\label{sec5.1}
\|f\|_{Z_{kg}}:=\sup_{k\in\mathbb{Z}}2^{N_0k^+}2^{k^-(1/2-\kappa)}\|\widehat{P_kf}\|_{L^\infty}.
\end{equation}
\end{definition}

\subsection{The main bootstrap proposition}\label{bootstrap0} Our main result is the following proposition:

\begin{proposition}\label{bootstrap} Assume that $(g,\psi)$ is a solution of the system \eqref{asum3}-\eqref{asum5} on the time interval $[0,T]$, $T\geq 1$. Assume that the initial data set $(\overline{g}_{ij},k_{ij},\psi_0,\psi_1)$, see \eqref{data}, satisfies the smallness conditions \eqref{asum1} and the constraint equations \eqref{asum2}. 

Define $U^G,U^{\mathcal{L}h_{\al\be}}, U^{\mathcal{L}\psi}$ as in \eqref{variables4}--\eqref{variables4L} and recall the definitions \eqref{qaz3.1}. Assume that, for any $t\in[0,T]$, the solution satisfies the bootstrap hypothesis
\begin{equation}\label{bootstrap2.1}
\sup_{q\leq n\leq 3,\,\mathcal{L}\in\mathcal{V}_{n}^q}\langle t\rangle^{-H(q,n)\delta}\big\{\|(\langle t\rangle|\nabla|_{\leq 1})^{\gamma}|\nabla|^{-1/2}U^{\mathcal{L}h_{\al\be}}(t)\|_{H^{N(n)}}+\|U^{\mathcal{L}\psi}(t)\|_{H^{N(n)}}\big\}\leq\varep_{1},
\end{equation}
\begin{equation}\label{bootstrap2.2}
\begin{split}
\sup_{q\leq n\leq 2,\,\mathcal{L}\in\mathcal{V}_{n}^q}&\sup_{k\in\mathbb{Z},\,l\in\{1,2,3\}}\,2^{N(n+1)k^+}\langle t\rangle^{-H(q+1,n+1)\delta}\\
&\big\{2^{k/2}(2^{k^-}\langle t\rangle)^{\gamma}\|P_k(x_lV^{\mathcal{L}h_{\al\be}})(t)\|_{L^2}+2^{k^+}\|P_k(x_lV^{\mathcal{L}\psi})(t)\|_{L^2}\big\}\leq\varep_{1},
\end{split}
\end{equation}
and
\begin{equation}\label{bootstrap2.4}
\|V^{F}(t)\|_{Z_{wa}}+\|V^{\omega_a}(t)\|_{Z_{wa}}+\|V^{\vartheta_{ab}}(t)\|_{Z_{wa}}+\langle t\rangle^{-\delta}\|V^{h_{\al\be}}(t)\|_{Z_{wa}}+\|V^{\psi}(t)\|_{Z_{kg}}\leq\varep_{1},
\end{equation}
for any $\al,\be\in\{0,1,2,3\}$ and $a,b\in\{1,2,3\}$. Here $\langle t\rangle:=\sqrt{1+t^2}$, $\varep_1:=\varep_0^{2/3}$, and
\begin{equation}\label{fvc1.0}
H(q,n):=
\begin{cases}
1\qquad&\text{ if }q=0\text{ and }n=0,\\
60(n-1)+20\qquad&\text{ if }q=0\text{ and }n\geq 1,\\
200(n-1)+30\qquad&\text{ if }q=1\text{ and }n\geq 1,\\
100(q+1)(n-1)\qquad&\text{ if }q\geq 2.
\end{cases}
\end{equation}
Then, for any $t\in[0,T]$, $\al,\be\in\{0,1,2,3\}$, and $a,b\in\{1,2,3\}$, one has the improved bounds
\begin{equation}\label{bootstrap3.1}
\sup_{q\leq n\leq 3,\,\mathcal{L}\in\mathcal{V}_{n}^q}\langle t\rangle^{-H(q,n)\delta}\big\{\|(\langle t\rangle|\nabla|_{\leq 1})^{\gamma}|\nabla|^{-1/2}U^{\mathcal{L}h_{\al\be}}(t)\|_{H^{N(n)}}+\|U^{\mathcal{L}\psi}(t)\|_{H^{N(n)}}\big\}\lesssim\varep_0,
\end{equation}
\begin{equation}\label{bootstrap3.2}
\begin{split}
\sup_{q\leq n\leq 2,\,\mathcal{L}\in\mathcal{V}_{n}^q}&\sup_{k\in\mathbb{Z},\,l\in\{1,2,3\}}\,2^{N(n+1)k^+}\langle t\rangle^{-H(q+1,n+1)\delta}\\
&\big\{2^{k/2}(2^{k^-}\langle t\rangle)^{\gamma}\|P_k(x_lV^{\mathcal{L}h_{\al\be}})(t)\|_{L^2}+2^{k^+}\|P_k(x_lV^{\mathcal{L}\psi})(t)\|_{L^2}\big\}\lesssim\varep_0,
\end{split}
\end{equation}
and
\begin{equation}\label{bootstrap3.4}
\|V^{F}(t)\|_{Z_{wa}}+\|V^{\omega_a}(t)\|_{Z_{wa}}+\|V^{\vartheta_{ab}}(t)\|_{Z_{wa}}+\langle t\rangle^{-\delta}\|V^{h_{\al\be}}(t)\|_{Z_{wa}}+\|V^{\psi}(t)\|_{Z_{kg}}\lesssim\varep_0.
\end{equation}
\end{proposition}

We will show in section \ref{Proof1} below that the smallness assumptions \eqref{asum1} on the initial data imply the bounds \eqref{bootstrap3.1}--\eqref{bootstrap3.4} at time $t=0$, and for all $t\in[0,2]$. Then we will show that Proposition \ref{bootstrap} implies the our main conclusions, in the quantitative form of Theorem \ref{winr}, and derive some additional asymptotic information about the solutions (in Chapter \ref{Chapter6}). Thus most of the work in this monograph, Chapters 2, 3, 4, and 5, is concerned with the proof of Proposition \ref{bootstrap}.

Notice that we aim to control, simultaneously, three types of norms: (i) energy norms involving up to 3 vector-fields $\Gamma_a$ and $\Omega_{ab}$, measured in Sobolev spaces, (ii) weighted norms on the linear profiles $V^{\LL h_{\al\be}}$ and $V^{\LL\psi}$, and (iii) the $Z$-norms on the undifferentiated profiles. See subsections \ref{hus1}--\ref{hus6} for a discussion and motivation on the roles of these different norms.

The function $H$ defined in \eqref{fvc1.0} is important, as it establishes a hierarchy of growth of the various energy norms. At the conceptual level this is needed because we define the weighted vector-fields $\Gamma_a,\Omega_{ab}$ in terms of the coordinate functions $x_j$ and $t$, thus we expect (at least logarithmic) losses as we apply more of these vector-fields. 

At the technical level, the function $H$ satisfies superlinear inequalities of the form
\begin{equation}\label{SuperlinearH1}
H(q_1,n_1)+H(q_2,n_2)\leq H(q_1+q_2,n_1+n_2)-40,
\end{equation}
if $n_1,n_2\geq 1$ and $n_1+n_2\leq 3$, and more refined versions. Such inequalities are helpful to estimate nonlinear interactions when the vector-fields split among the different components.

We notice also that we treat the two types of weighted vector-fields $\Gamma_a$ and $\Omega_{ab}$ differently, in the sense that the application of the non-tangential vector-fields $\Gamma_a$ leads to more loss in terms of time growth than the application of the tangential vector-fields $\Omega_{ab}$ (for example $H(0,1)=20<H(1,1)=30$). This is a subtle technical point to keep in mind, connected to a more general difficulty of estimating the effect of non-tangential vector-fields.

\section{Some lemmas}\label{SLSec} We are now ready to start the proof of the main bootstrap proposition \ref{bootstrap}. In this section we collect several results that are used in the rest of the paper. 

\subsection{General lemmas} We start with a lemma that is used often in integration by parts arguments. See \cite[Lemma 5.4]{IoPa2} for the proof.

\begin{lemma}\label{tech5} Assume that $0<\eps\leq 1/\eps\leq K$, $N\geq 1$ is an integer, and $f,g\in C^{N+1}(\mathbb{R}^3)$. Then
\begin{equation}\label{ln1}
\Big|\int_{\mathbb{R}^3}e^{iKf}g\,dx\Big|\lesssim_N (K\eps)^{-N}\big[\sum_{|\alpha|\leq N}\eps^{|\alpha|}\|D^\alpha_xg\|_{L^1}\big],
\end{equation}
provided that $f$ is real-valued,
\begin{equation}\label{ln2}
|\nabla_x f|\geq \mathbf{1}_{{\mathrm{supp}}\,g},\quad\text{ and }\quad\|D_x^\alpha f \cdot\mathbf{1}_{{\mathrm{supp}}\,g}\|_{L^\infty}\lesssim_N\eps^{1-|\alpha|},\,2\leq |\alpha|\leq N+1.
\end{equation}
\end{lemma}

To bound multilinear operators, we often use the following simple lemma.

\begin{lemma}\label{L1easy}
(i) Assume that $l\geq 2$, $f_1,\ldots, f_l,f_{l+1}\in L^2(\mathbb{R}^3)$, and $M:(\mathbb{R}^3)^l\to\mathbb{C}$ is a continuous compactly supported function. Then
\begin{equation}\label{ener62}
\begin{split}
\Big|\int_{(\mathbb{R}^3)^l}M(\xi_1,\ldots,\xi_l)\cdot\widehat{f_1}(\xi_1)\cdot\ldots\cdot \widehat{f_l}(\xi_l)\cdot\widehat{f_{l+1}}(-\xi_1-\ldots-\xi_l)\,d\xi_1\ldots d\xi_l\Big|\\
\lesssim \big\|\mathcal{F}^{-1}M\big\|_{L^1((\mathbb{R}^3)^l)}\|f_1\|_{L^{p_1}}\cdot\ldots\cdot\|f_{l+1}\|_{L^{p_{l+1}}},
\end{split}
\end{equation}
for any exponents $p_1,\ldots,p_{l+1}\in[1,\infty]$ satisfying $1/p_1+\ldots+1/p_{l+1}=1$. 

(ii) As a consequence, if $q,p_2,p_3\in[1,\infty]$ satisfy $1/p_2+1/p_3=1/q$ then
\begin{equation}\label{ener62.1}
\Big\|\mathcal{F}_{\xi}^{-1}\Big\{\int_{\mathbb{R}^3}M(\xi,\eta)\widehat{f}(\eta)\widehat{g}(-\xi-\eta)\,d\eta\Big\}\Big\|_{L^{q}}\lesssim \big\|\mathcal{F}^{-1}M\big\|_{L^1}\|f\|_{L^{p_2}}\|g\|_{L^{p_3}}.
\end{equation}
\end{lemma}

We recall also a Hardy-type estimate. See \cite[Lemma 3.5]{IoPa3} for the proof of (i), while part (ii) follows easily from definitions.

\begin{lemma}\label{hyt1}
(i) For $f\in L^2(\mathbb{R}^3)$ and $k\in\mathbb{Z}$ let
\begin{equation}\label{hyt2}
\begin{split}
&A_k:=\|P_kf\|_{L^2}+\sum_{l=1}^3\|\varphi_k(\xi)(\partial_{\xi_l}\widehat{f})(\xi)\|_{L^2_\xi},\qquad B_k:=\Big[\sum_{j\geq\max(-k,0)}2^{2j}\|Q_{j,k}f\|_{L^2}^2\Big]^{1/2}.
\end{split}
\end{equation} 
Then, for any $k\in\mathbb{Z}$,
\begin{equation}\label{hyt3}
A_k\lesssim\sum_{|k'-k|\leq 4}B_{k'}
\end{equation}
and
\begin{equation}\label{hyt3.1}
B_k\lesssim
\begin{cases}
\sum_{|k'-k|\leq 4}A_{k'}\qquad&\text{ if } k\geq 0,\\
\sum_{k'\in\mathbb{Z}}A_{k'}2^{-|k-k'|/2}\min(1,2^{k'-k})\qquad&\text{ if } k\leq 0.
\end{cases}
\end{equation}

(ii) If $m\in\mathcal{M}_0$, see \eqref{mults0}, then, for any $(k,j)\in\mathcal{J}$,
\begin{equation}\label{hyt4.4}
\|Q_{j,k}\{\mathcal{F}^{-1}(m\widehat{f}\,)\}\|_{L^2}\lesssim \sum_{j'\geq \max(-k,0)}\|Q_{j',k}f\|_{L^2}2^{-4|j-j'|}.
\end{equation}
\end{lemma}

\subsection{The phases $\Phi_{\sigma\mu\nu}$} Our normal form analysis relies on precise bounds on the phases $\Phi_{\sigma\mu\nu}$ defined in \eqref{on9.2}. We summarize the results we need in this subsection.

We consider first ${\mathrm{Wave}}\,\times{\mathrm{KG}}\to{\mathrm{KG}}$ and ${\mathrm{KG}}\,\times{\mathrm{KG}}\to{\mathrm{Wave}}$ interactions. These interactions are weakly elliptic, in the sense that the corresponding phases $\Phi_{\sigma\mu\nu}$ do not vanish, except when the wave frequency vanishes. More precisely, we have the following quantitative estimates.

\begin{lemma}\label{pha2} (i) Assume that $\Phi_{\sigma\mu\nu}$ is as in \eqref{on9.2}. If $|\xi|,|\xi-\eta|,|\eta|\in[0,b]$, $1\leq b$, then
\begin{equation}\label{pha3}
\begin{split}
&|\Phi_{\sigma\mu\nu}(\xi,\eta)|\geq |\xi|/(4b^2)\qquad \text{ if }\,\,(\sigma,\mu,\nu)=((wa,\iota),(kg,\iota_1),(kg,\iota_2)),\\
&|\Phi_{\sigma\mu\nu}(\xi,\eta)|\geq |\eta|/(4b^2)\qquad\text{ if }\,\,(\sigma,\mu,\nu)=((kg,\iota),(kg,\iota_1),(wa,\iota_2)).
\end{split}
\end{equation}

(ii) Assume that $k,k_1,k_2\in \mathbb{Z}$ and $n$ is a multiplier such that $\|\mathcal{F}^{-1}n\|_{L^1(\mathbb{R}^3\times\mathbb{R}^3)}\leq 1$. Let $\overline{k}=\max(k,k_1,k_2)$. If $(\sigma,\mu,\nu)=((wa,\iota),(kg,\iota_1),(kg,\iota_2))$ then
\begin{equation}\label{pha4}
\big\|\mathcal{F}^{-1}\{\Phi_{\sigma\mu\nu}(\xi,\eta)^{-1}n(\xi,\eta)\cdot\varphi_k(\xi)\varphi_{k_1}(\xi-\eta)\varphi_{k_2}(\eta)\}\big\|_{L^1(\mathbb{R}^3\times\mathbb{R}^3)}\lesssim 2^{-k}2^{4\overline{k}^+}.
\end{equation}
Moreover, if $(\sigma,\mu,\nu)=((kg,\iota),(kg,\iota_1),(wa,\iota_2))$ then
\begin{equation}\label{pha4.1}
\big\|\mathcal{F}^{-1}\{\Phi_{\sigma\mu\nu}(\xi,\eta)^{-1}n(\xi,\eta)\cdot\varphi_k(\xi)\varphi_{k_1}(\xi-\eta)\varphi_{k_2}(\eta)\}\big\|_{L^1(\mathbb{R}^3\times\mathbb{R}^3)}\lesssim 2^{-k_2}2^{4\overline{k}^+}.
\end{equation}
\end{lemma}

\begin{proof} The  bounds \eqref{pha3}, \eqref{pha4}, \eqref{pha4.1} were proved in Lemma 3.3 in \cite{IoPa3}. 
\end{proof}

We will also consider ${\mathrm{Wave}}\,\times{\mathrm{Wave}}\to{\mathrm{Wave}}$ interactions. In this case the corresponding bilinear phases $\Phi_{\sigma\mu\nu}$ can vanish on large sets, when the frequencies are parallel, and the strength of these interactions depends significantly on the angle between the frequencies of the inputs. To measure these angles, for $\iota_1,\iota_2\in\{+,-\}$ we define the functions
\begin{equation}\label{par1}
\Xi_{\iota_1\iota_2}:(\mathbb{R}^3\setminus\{0\})^2\to B_2,\qquad \Xi_{\iota_1\iota_2}(\theta,\eta):=\iota_1\frac{\theta}{|\theta|}-\iota_2\frac{\eta}{|\eta|},
\end{equation}
where $B_R:=\{x\in\mathbb{R}^3:\,|x|\leq R\}$. Let $\Xi_{\iota_1\iota_2,k}(\theta,\eta):=\iota_1\theta_k/|\theta|-\iota_2\eta_k/|\eta|$, $k\in\{1,2,3\}$.

In multilinear estimates we use sometimes the following elementary lemma:

\begin{lemma}\label{par70}
(i) By convention, let $++=--=+$ and $+-=-+=-$. Then
\begin{equation}\label{par2.1}
|\Xi_{\iota_1\iota_2}(\theta,\eta)|=\sqrt{\frac{2(|\theta||\eta|-\iota_1\iota_2\theta\cdot\eta)}{|\theta||\eta|}}.
\end{equation}

(ii) We define the functions
\begin{equation}\label{par71}
\widetilde{\Xi}:(\mathbb{R}^3\setminus\{0\})^2\to[0,2],\qquad \widetilde{\Xi}(\theta,\eta):=\Big|\frac{\theta}{|\theta|}-\frac{\eta}{|\eta|}\Big|\Big|\frac{\theta}{|\theta|}+\frac{\eta}{|\eta|}\Big|=2\sqrt{1-\frac{|\theta\cdot\eta|^2}{|\theta|^2|\eta|^2}}.
\end{equation}
Then, for any $\theta,\eta\in\mathbb{R}^3\setminus\{0\}$ we have
\begin{equation}\label{par72}
\min\{|\Xi_{+}(\theta,\eta)|,|\Xi_{-}(\theta,\eta)|\}\leq \widetilde{\Xi}(\theta,\eta)\leq 2\min\{|\Xi_{+}(\theta,\eta)|,|\Xi_{-}(\theta,\eta)|\}.
\end{equation}
Moreover, for any $x,y,z\in\mathbb{R}^3\setminus\{0\}$ we have
\begin{equation}\label{par72.4}
\widetilde{\Xi}(x,z)\leq 2[\widetilde{\Xi}(x,y)+\widetilde{\Xi}(y,z)].
\end{equation}
In addition, if $y+z\neq 0$ then
\begin{equation}\label{par73}
\widetilde{\Xi}(x,y+z)|y+z|\leq \widetilde{\Xi}(x,y)|y|+\widetilde{\Xi}(x,z)|z|.
\end{equation}
\end{lemma}

\begin{proof} The identities \eqref{par2.1} follow directly from definitions. The inequalities \eqref{par72} follow from the definitions as well, once we notice that $\widetilde{\Xi}(\theta,\eta)=|\Xi_{+}(\theta,\eta)|\cdot |\Xi_{-}(\theta,\eta)|$ and $|\Xi_{+}(\theta,\eta)|+|\Xi_{-}(\theta,\eta)|\geq 2$. To prove \eqref{par72.4} we notice that
\begin{equation*}
\min\{|\Xi_+(x,z)|,|\Xi_-(x,z)|\}\leq \min\{|\Xi_+(x,y)|,|\Xi_-(x,y)|\}+\min\{|\Xi_+(y,z)|,|\Xi_-(y,z)|\},
\end{equation*}
and then use \eqref{par72}. Finally, to prove \eqref{par73} we may assume that $x=(1,0,0)$, $y=(y_1,y')$, $z=(z_1,z')$, and estimate, using just \eqref{par71},
\begin{equation*}
\widetilde{\Xi}(x,y+z)|y+z|=2|y'+z'|\leq 2|y'|+2|z'|=\widetilde{\Xi}(x,y)|y|+\widetilde{\Xi}(x,z)|z|,
\end{equation*}
which gives the desired conclusion.
\end{proof}

We consider now trilinear expressions localized with respect to angular separation, as well as expressions resulting from normal form transformations.

\begin{lemma}\label{PhaWave} (i) Assume $\chi_1:\mathbb{R}^3\to[0,1]$ is a smooth function supported in the ball $B_2$, $\iota_1,\iota_2\in\{+,-\}$, $b\leq 2$, $f,f_1,f_2\in L^2(\mathbb{R}^3)$, and $k,k_1,k_2\in\mathbb{Z}$. Let
\begin{equation}\label{yip1}
L^b_{k,k_1,k_2}:=\int_{\mathbb{R}^3\times\mathbb{R}^3}m(\xi-\eta,\eta)\chi_1(2^{-b}\Xi_{\iota_1\iota_2}(\xi-\eta,\eta))\cdot \widehat{P_{k_1}f_1}(\xi-\eta)\widehat{P_{k_2}f_2}(\eta)\widehat{P_{k}f}(\xi)\,d\xi d\eta,
\end{equation}
where $m$ is a symbol satisfying $\|\mathcal{F}^{-1}(m)\|_{L^1(\mathbb{R}^6)}\leq 1$. Then
\begin{equation}\label{yip2}
|L^b_{k,k_1,k_2}|\lesssim \min\{2^{-b},2^{k_1-k}+1\}\|P_{k_1}f_1\|_{L^{\infty}}\|P_{k_2}f_2\|_{L^{2}}\|P_{k}f\|_{L^2}.
\end{equation}

(ii) Assume that $\chi_2:\mathbb{R}^3\to[0,1]$ is a smooth function supported in $B_2\setminus B_{1/2}$, $(\sigma,\mu,\nu)=((wa,\iota),(wa,\iota_1),(wa,\iota_2))$, $\iota,\iota_1,\iota_2\in\{+,-\}$. Then, for any $k,k_1,k_2\in\mathbb{Z}$,
\begin{equation}\label{par5}
|\Phi_{\sigma\mu\nu}(\xi,\eta)^{-1}\cdot\varphi_k(\xi)\varphi_{k_1}(\xi-\eta)\varphi_{k_2}(\eta)\chi_2(2^{-b}\Xi_{\iota_1\iota_2}(\xi-\eta,\eta))|\lesssim 2^{-2b}2^{-\min(k_1,k_2)},
\end{equation}
where $\Phi_{\sigma\mu\nu}(\xi,\eta)=\Lambda_{\sigma}(\xi)-\Lambda_{\mu}(\xi-\eta)-\Lambda_\nu(\eta)$ as in \eqref{on9.2}. In addition, if we define
\begin{equation}\label{yip5}
M^b_{k,k_1,k_2}:=\int_{\mathbb{R}^3\times\mathbb{R}^3}m(\xi-\eta,\eta)\frac{\chi_2(2^{-b}\Xi_{\iota_1\iota_2}(\xi-\eta,\eta))}{\Phi_{\sigma\mu\nu}(\xi,\eta)}\widehat{P_{k_1}f_1}(\xi-\eta)\widehat{P_{k_2}f_2}(\eta)\widehat{P_{k}f}(\xi)\,d\xi d\eta,
\end{equation}
where $\|\mathcal{F}^{-1}(m)\|_{L^1(\mathbb{R}^6)}\leq 1$, then
\begin{equation}\label{yip6}
|M^b_{k,k_1,k_2}|\lesssim 2^{-2b}2^{-\min(k_1,k_2)}\cdot\min\{2^{-b},2^{k_1-k}+1\}\|P_{k_1}f_1\|_{L^{\infty}}\|P_{k_2}f_2\|_{L^{2}}\|P_{k}f\|_{L^2}.
\end{equation}

(iii) If $\chi_3:\mathbb{R}\to[0,1]$ is a smooth function supported in $[-2,2]\setminus[-1/2,1/2]$ and $\Phi_{\sigma,\mu,\nu}$ is as above then, for any $k,k_1,k_2\in\mathbb{Z}$,
\begin{equation}\label{par6}
\sup_{\xi\in\mathbb{R}^3}\Big\|\int_{\mathbb{R}^3}e^{iy\cdot\eta}\frac{\chi_3(2^{-2b}|\Xi_{\iota_1\iota_2}(\xi-\eta,\eta)|^2)}{\Phi_{\sigma\mu\nu}(\xi,\eta)}\varphi_k(\xi)\varphi_{k_1}(\xi-\eta)\varphi_{k_2}(\eta)\,d\eta\Big\|_{L^1_y}\lesssim 2^{-2b}2^{-\min(k_1,k_2)}.
\end{equation}
\end{lemma}

\begin{remark} All estimates would follow from the ``natural" localization bounds
\begin{equation}\label{nag1der}
\big\|\mathcal{F}^{-1}\{\chi_1(2^{-b}\Xi_{\iota_1\iota_2}(\xi-\eta,\eta))\varphi_{0}(\xi-\eta)\varphi_{0}(\eta)\}\big\|_{L^1(\mathbb{R}^6)}\lesssim 1
\end{equation}
and the identities \eqref{par7}. We are not able to prove these bounds, but we prove the weaker bounds \eqref{nag1} and \eqref{nag10}, which still allow us to derive the conclusions of the lemma.
\end{remark}

\begin{proof} {\bf{Step 1.}} We show first that, for any $k_1,k_2\in\mathbb{Z}$ and $\iota_1,\iota_2\in\{+,-\}$, 
\begin{equation}\label{nag1}
\big\|\mathcal{F}^{-1}\{\chi_1(2^{-b}\Xi_{\iota_1\iota_2}(\xi-\eta,\eta))\varphi_{k_1}(\xi-\eta)\varphi_{k_2}(\eta)\}\big\|_{L^1(\mathbb{R}^6)}\lesssim 2^{-b}.
\end{equation}
After changes of variables we may assume that $\iota_1=\iota_2=+$, $k_1=k_2=0$. We have to prove that $\big\|F_b\big\|_{L^1}\lesssim 2^{-b}$ for any $b\leq 2$, where
\begin{equation}\label{nag2}
F_b(x,y):=\int_{\mathbb{R}^6}e^{ix\cdot\theta}e^{iy\cdot\eta}\chi_1(2^{-b}\Xi_{++}(\theta,\eta))\varphi_0(\theta)\varphi_0(\eta)\,d\theta d\eta.
\end{equation}
In fact, we will prove the stronger pointwise bounds
\begin{equation}\label{nag3}
|F_b(x,y)|\lesssim \frac{1}{\big[1+2^{2b}|x|^2+2^{2b}|y|^2\big]^8}\frac{1}{\big[1+\min\{|x|,|y|\}\widetilde{\Xi}(x,y)\big]^8}\frac{2^{2b}}{(1+|x|^2+|y|^2)^{1/2}},
\end{equation}
for any $x,y\in\mathbb{R}^3$, which would clearly imply the desired conclusion $\big\|F_b\big\|_{L^1}\lesssim 2^{-b}$.

In proving \eqref{nag3}, we may assume that $|y|\leq |x|$, $x=(x_1,0,0)$, $x_1\geq 0$, $y=(y_1,y_2,y_3)$, $|y_2|\geq |y_3|$. The desired bounds \eqref{nag3} are then equivalent to
\begin{equation}\label{nag4}
|F_b(x,y)|\lesssim \frac{1}{[1+2^{2b}x_1^2]^{8}}\frac{1}{[1+|y_2|]^8}\frac{2^{2b}}{(1+x_1^2)^{1/2}},
\end{equation}
where $y=(y_1,y')$. Notice that this follows easily if $x_1\lesssim 1$, since the $\theta,\eta$ integration is taken over a set of volume $\approx 2^{2b}$.

We will integrate by parts in $\theta$ and $\eta$ using the operators
\begin{equation}\label{nag5}
2^{2b}\Delta_\theta,\qquad 2^{2b}\Delta_\eta, \qquad L_\theta:=\theta_j\partial_{\theta_j},\qquad L_\eta:=\eta_j\partial_{\eta_j},\qquad S_{ij}:=\theta_i\partial_{\theta_j}+\eta_i\partial_{\eta_j}.
\end{equation}
The main point is that the vector-fields $L_\theta, L_\eta, S_{ij}$ act well on $\Xi_{++}(\theta,\eta)$. More precisely,
\begin{equation}\label{nag6}
\begin{split}
&(L_\theta \Xi_{++,k})(\theta,\eta)=0,\qquad (L_\eta\Xi_{++,k})(\theta,\eta)=0,\\
&(S_{ij}\Xi_{++,k})(\theta,\eta)=\delta_{jk}\big[\theta_i/|\theta|-\eta_i/|\eta|\big]-\big[\theta_i\theta_j\theta_k/|\theta|^3-\eta_i\eta_j\eta_k/|\eta|^3\big].
\end{split}
\end{equation}
In particular, it follows that if $\mathcal{O}$ is any combination of the operators $L_\theta, L_\eta, S_{ij}, 2^{2b}\Delta_\theta$ then $\mathcal{O}\big\{\chi_1(2^{-b}\Xi_{++}(\theta,\eta))\varphi_0(\theta)\varphi_0(\eta)\big\}$ is a function of the form $\widetilde{\chi}[2^{-b}\Xi_{++}(\theta,\eta),\theta,\eta]$, where $\widetilde{\chi}=\widetilde{\chi}_{\mathcal{O}}$ is a smooth function supported in the set $B_2\times(B_2\setminus B_{1/2})^2$.

We notice now that 
\begin{equation*}
S_{j1}\big\{e^{ix\cdot\theta}e^{i y\eta}\big\}=i\big(\theta_jx_1+\eta_jy_1\big)\{e^{ix\cdot\theta}e^{i y\eta}\}.
\end{equation*}
Therefore, if $|y_1|\leq 2^{-20}x_1$ then we can integrate by parts using only the vector-fields $S_{j1}$  and gain a factor of $x_1$ at every iteration. It follows that $|F_b(x,y)|\lesssim 2^{2b}(1+x_1^2)^{-20}$, which is better than the desired bounds \eqref{nag4}. 

On the other hand, if $|y_1|\geq 2^{-20}x_1$ then we first integrate by parts using $2^{2b}\Delta_\theta$ and $S_{j2}$, $j\in\{1,2,3\}$, to gain the factors $(1+2^{2b}x_1^2)^{-10}$ and $(1+|y_2|)^{-10}$ in \eqref{nag4}. Letting
\begin{equation}\label{nag9}
F'_b(x,y):=\int_{\mathbb{R}^6}e^{ix\cdot\theta}e^{iy\cdot\eta}\widetilde{\chi}[2^{-b}\Xi_{++}(\theta,\eta),\theta,\eta]\,d\theta d\eta,
\end{equation}
where $\widetilde{\chi}$ is a smooth function supported in the set $B_2\times(B_2\setminus B_{1/2})^2$, it remains to prove that 
\begin{equation}\label{nag8}
|F'_b(x,y)|\lesssim 2^{2b}(1+|x_1|)^{-1}.
\end{equation}

For this we use the scaling vector-field $L_\theta$. For integers $n\geq 0$ we define $F'_{b;n}$ by inserting cutoff functions of the form $\varphi_{n}^{[0,\infty)}(x\cdot \theta)$ in the integral in \eqref{nag9}. Clearly, $|F'_{b;0}(x,y)|\lesssim 2^{2b}(1+|x_1|)^{-1}$. We integrate by parts twice using the scaling vector-field $L_\theta$ to show that $|F'_{b;n}(x,y)|\lesssim 2^{-n}2^{2b}(1+|x_1|)^{-1}$, and the desired bounds \eqref{nag8} follow.

{\bf{Step 2.}} With $\chi_1$ as in (i), we show now that 
\begin{equation}\label{nag10}
\big\|\mathcal{F}^{-1}\{\chi_1(2^{-b}\Xi_{\iota_1\iota_2}(\xi-\eta,\eta))\chi'(2^{-b}(\eta/|\eta|-e))\varphi_{k_1}(\xi-\eta)\varphi_{k_2}(\eta)\}\big\|_{L^1(\mathbb{R}^6)}\lesssim 1,
\end{equation}
for any smooth function $\chi':\mathbb{R}^3\to[0,1]$  supported in the ball $B_2$ and any vector $e\in\mathbb{S}^2$. Indeed, to prove this we can rescale to $k_1=k_2=0$. Then we define
\begin{equation}\label{nag11}
G_b(x,y):=\int_{\mathbb{R}^6}e^{ix\cdot\theta}e^{iy\cdot\eta}\chi_1(2^{-b}\Xi_{++}(\theta,\eta))\chi'(2^{-b}(\eta/|\eta|-e))\varphi_0(\theta)\varphi_0(\eta)\,d\theta d\eta,
\end{equation}
and notice that it suffices to prove the pointwise bounds
\begin{equation}\label{nag12}
|G_b(x,y)|\lesssim \frac{1}{\big[1+2^{2b}|x|^2+2^{2b}|y|^2\big]^8}\frac{2^{4b}}{(1+|x\cdot e|^2+|y\cdot e|^2)^{8}},
\end{equation}
for any $x,y\in\mathbb{R}^3$. These bounds follow easily by integration by parts as before, using \eqref{nag6} and the operators $2^{2b}\Delta_\theta$, $2^{2b}\Delta_\eta$, $L_\theta$, $L_\eta$, and $S_{ij}$ defined in \eqref{nag5}. 

{\bf{Step 3.}} We prove now the bounds \eqref{yip2}. The estimates with the factor $2^{-b}$ follow from \eqref{nag1} and \eqref{ener62}. To prove the estimates with the factor $2^{k_1-k}+1$ we need to introduce an angular decomposition. Given $q\in\mathbb{Z}$, $q\leq 2$, we fix a $2^q$-net $\mathfrak{N}_q$ on $\mathbb{S}^2$ and define
\begin{equation}\label{nag15}
\varphi_{k;q,e}(\xi):=\varphi_k(\xi)\frac{\varphi_{\leq 0}(2^{-2q}|\xi/|\xi|-e|^2)}{\sum_{e'\in\mathfrak{N}_q}\varphi_{\leq 0}(2^{-2q}|\xi/|\xi|-e'|^2)},\qquad P_{k;q,e}f:=\mathcal{F}^{-1}\{\varphi_{k;q,e}\widehat{f}\},
\end{equation}
for any $k\in\mathbb{Z}$ and $e\in\mathfrak{N}_q$. We insert the partition of unity $\{\varphi_{k_2;b,e}(\eta)\}_{e\in\mathfrak{N}_b}$ in the integrals in \eqref{yip1}. Notice also that if $|\Xi_{\iota_1\iota_2}(\xi-\eta,\eta)|\lesssim 2^b$ and $|\eta/|\eta|-e|\lesssim 2^b$ then $\widetilde{\Xi}(\xi,e)\lesssim 2^b(2^{k_1-k}+1)$ in the support of the integral, as a consequence of \eqref{par72.4}--\eqref{par73}. Thus
\begin{equation*}
\begin{split}
&L^{b}_{k,k_1,k_2}=\sum_{e\in\mathfrak{N}_b}L^{b,e}_{k,k_1,k_2},\\
&L^{b,e}_{k,k_1,k_2}:=\int_{\mathbb{R}^3\times\mathbb{R}^3}m(\xi-\eta,\eta)\chi_1(2^{-b}\Xi_{\iota_1\iota_2}(\xi-\eta,\eta))\cdot \widehat{P_{k_1}f_1}(\xi-\eta)\widehat{P_{k_2;b,e}f_2}(\eta)\widehat{P_kA_{b',e}f}(\xi)\,d\xi d\eta,
\end{split}
\end{equation*}
where $b':=b+\max\{k_1-k,0\}+C$, and $\widehat{A_{b',e}f}(\xi)=\widehat{f}(\xi)\varphi_{\leq 0}(2^{-b'}\widetilde{\Xi}(\xi,e))$. Thus 
\begin{equation*}
\begin{split}
|L^{b}_{k,k_1,k_2}|^2&\lesssim \|P_{k_1}f_1\|_{L^\infty}^2\big\{\sum_{e\in\mathfrak{N}_b}\|P_{k_2;b,e}f_2\|_{L^2}^2\big\}\big\{\sum_{e\in\mathfrak{N}_b}\|P_kA_{b',e}f\|_{L^2}^2\big\}\\
&\lesssim \|P_{k_1}f_1\|_{L^\infty}^2\|P_{k_2}f_2\|_{L^2}^2\cdot 2^{2(b'-b)}\|P_kf\|_{L^2}^2,
\end{split}
\end{equation*}
using \eqref{nag10}, \eqref{ener62}, and orthogonality. The desired bounds \eqref{yip2} follow.

{\bf{Step 4.}} We prove now the bounds in (ii). With $\theta=\xi-\eta$ we write
\begin{equation}\label{par7}
\frac{1}{\Phi_{\sigma\mu\nu}(\xi,\eta)}=\frac{1}{\iota|\xi|-\iota_1|\xi-\eta|-\iota_2|\eta|}=\frac{\iota|\theta+\eta|+\iota_1|\theta|+\iota_2|\eta|}{2(\theta\cdot\eta-\iota_1\iota_2|\theta||\eta|)}=\frac{\iota|\theta+\eta|+\iota_1|\theta|+\iota_2|\eta|}{-\iota_1\iota_2|\theta||\eta||\Xi_{\iota_1\iota_2}(\theta,\eta)|^2}.
\end{equation}
The bounds \eqref{par5} and \eqref{yip6} follow, using also \eqref{yip2}. 

{\bf{Step 5.}} Since $\|\mathcal{F}^{-1}(m\cdot m')\|_{L^1}\lesssim \|\mathcal{F}^{-1}m\|_{L^1}\|\mathcal{F}^{-1}m'\|_{L^1}$, for \eqref{par6} it suffices to prove that
\begin{equation}\label{par9}
\Big\|\varphi_k(\xi)\int_{\mathbb{R}^3}e^{iy\cdot\eta}\chi_3\big(2^{-2b}|\Xi_{\iota_1\iota_2}(\xi-\eta,\eta)|^2\big)\cdot\varphi_{k_1}(\xi-\eta)\varphi_{k_2}(\eta)\,d\eta\Big\|_{L^1_y}\lesssim 1,
\end{equation}
for any $\xi\in\mathbb{R}^3$. By symmetry and rotation, we may assume that $k_2\leq k_1$ and $\xi=(\xi_1,0,0)$, $\xi_1>0$. The bounds \eqref{par9} follow by standard integration by parts if $b\geq -20$, since the function $H_{b,\xi}(\eta):=2^{-2b}|\Xi_{\iota_1\iota_2}(\xi-\eta,\eta)|^2$ satisfies differential bounds of the form $|D^\al_{\eta}H_{b;\xi}(\eta)|\lesssim_{|\alpha|}2^{-|\alpha|k_2}$ in the support of the integral, for all multi-indices $\alpha\in\mathbb{Z}_+^3$. 

When $b\leq -20$ we have to be slightly more careful. Notice that if $b\leq -20$ then the function $\chi_3(2^{-2b}|\Xi_{\iota_1\iota_2}(\xi-\eta,\eta)|^2)\cdot\varphi_k(\xi)\varphi_{k_1}(\xi-\eta)\varphi_{k_2}(\eta)$ is nontrivial only if $b\leq k-\max(k_1,k_2)+10$. This can be verified easily by considering the two cases $\iota_1=\iota_2$ and $\iota_1=-\iota_2$, and examining the definitions. Notice also that, using \eqref{par2.1},
\begin{equation}\label{par9.1}
\begin{split}
H_{b;\xi}(\eta)&=2^{-2b+1}\frac{|\xi-\eta||\eta|-\iota_1\iota_2(\xi-\eta)\cdot\eta}{|\xi-\eta||\eta|}=\frac{2^{-2b+1}[|\xi-\eta|^2|\eta|^2-((\xi-\eta)\cdot\eta)^2]}{|\xi-\eta||\eta|[|\xi-\eta||\eta|+\iota_1\iota_2(\xi-\eta)\cdot\eta]}\\
&=\frac{2^{-2b+1}\xi_1^2(\eta_2^2+\eta_3^2)}{|\xi-\eta||\eta|[|\xi-\eta||\eta|+\iota_1\iota_2(\xi-\eta)\cdot\eta]}.
\end{split}
\end{equation} 
Notice that the denominator of fraction in the right-hand side above is $\approx 2^{2k_1+2k_2}$ in the support of the integral. Thus the $\eta$ integral in \eqref{par9} is supported in the set $\mathcal{R}_{b;\xi}:=\big\{|\eta|\approx 2^{k_2},\,|\xi-\eta|\approx 2^{k_1},\,\sqrt{\eta_2^2+\eta_3^2}\approx 2^{b+k_1+k_2-k}\big\}$ (according to the remark above, we may assume that $2^{b+k_1+k_2-k}\lesssim 2^{\min(k_1,k_2)}=2^{k_2}$). Moreover, it is easy to see that
\begin{equation*}
|\partial_{\eta_1}^{\al}H_{b;\xi}(\eta)|\lesssim 2^{-k_2|\al|},\quad |\partial_{\eta_2}^{\al}H_{b;\xi}(\eta)|\lesssim 2^{-|\al|(b+k_1+k_2-k)},\quad |\partial_{\eta_3}^{\al}H_{b;\xi}(\eta)|\lesssim 2^{-|\al|(b+k_1+k_2-k)},
\end{equation*}
for $\eta\in\mathcal{R}_{b;\xi}$, and $\al\in[0,10]$. The bounds \eqref{par9} follow by integration by parts in $\eta$.
\end{proof}

\subsection{Linear estimates}\label{line} Linear estimates and decompositions in the Fourier space are the main building blocks to prove nonlinear estimates. In this subsection we summarize the estimates we need.

In some of the estimates in the paper we need to localize in the Fourier space to rotational invariant sets that are thinner than dyadic. For this we fix a smooth function $\chi:\mathbb{R}\to[0,1]$ supported in $[-2,2]$ with the property that $\sum_{n\in\mathbb{Z}}\chi(x-n)=1$ for all $x\in\mathbb{R}$. Then we define the operators $\mathcal{C}_{n,l}$, $n\geq 4$, $l\in\mathbb{Z}$, by
\begin{equation}\label{suploc}
\widehat{\mathcal{C}_{n,l}g}(\xi):=\chi(|\xi|2^{-l}-n)\widehat{g}(\xi).
\end{equation}

We record now several linear dispersive estimates (see \cite[Lemma 3.4]{IoPa3} for all the proofs). 

\begin{lemma}\label{LinEstLem} For any $f\in L^2(\mathbb{R}^3)$ and $(k,j)\in\mathcal{J}$ let
\begin{equation}\label{defin}
f_{j,k}:=P'_{k}Q_{j,k}f,\qquad Q_{\leq j,k}f:=\sum_{j'\in[\max(-k,0),j]}Q_{j',k}f,\qquad f_{\leq j,k}:=P'_{k}Q_{\leq j,k}f,
\end{equation}
where $P'_k=P_{[k-2,k+2]}$. For simplicity of notation, let
\begin{equation}\label{defin2}
f^\ast_{j,k}:=Q_{j,k}f,\qquad f^\ast_{\leq j,k}:=Q_{\leq j,k}f.
\end{equation}

(i) Then, for any $\alpha\in(\mathbb{Z}_+)^3$,
\begin{equation}\label{Linfty3.4}
\|D^\alpha_\xi\widehat{f_{j,k}}\|_{L^2}\lesssim 2^{|\alpha|j}\|\widehat{f^\ast_{j,k}}\|_{L^2},\qquad \|D^\alpha_\xi\widehat{f_{j,k}}\|_{L^\infty}\lesssim 2^{|\alpha|j}\|\widehat{f^\ast_{j,k}}\|_{L^\infty}.
\end{equation}
Moreover we have
\begin{equation}\label{Linfty3.3}
\|\widehat{f_{j,k}}\|_{L^\infty}\lesssim \min\big\{2^{3j/2} \|f^\ast_{j,k}\|_{L^2},2^{j/2-k}2^{\delta(j+k)/20}\|f^\ast_{j,k}\|_{H^{0,1}_\Omega}\big\},
\end{equation}
\begin{equation}\label{Linfty3.33}
\|\widehat{f_{j,k}}(r\theta)\|_{L^2(r^2dr)L^\infty_\theta}\lesssim 2^{j+k} \|f^\ast_{j,k}\|_{L^2},
\end{equation}
\begin{equation}\label{Linfty3.34}
\|\widehat{f_{j,k}}(r\theta)\|_{L^2(r^2dr)L^p_\theta}\lesssim_p \|f^\ast_{j,k}\|_{H^{0,1}_\Omega},\qquad p\in[2,\infty),
\end{equation}
and
\begin{equation}\label{Linfty3.31}
\|\widehat{f^\ast_{j,k}}-\widehat{f_{j,k}}\|_{L^\infty}\lesssim 2^{3j/2}2^{-4(j+k)}\|P_kf\|_{L^2}.
\end{equation}

(ii) For any $t\in\mathbb{R}$, $(k,j)\in\mathcal{J}$, and $f\in L^2(\mathbb{R}^3)$ we have
\begin{equation}\label{Linfty1}
\|e^{-it\Lambda_{wa}}f_{j,k}\|_{L^\infty}\lesssim 2^{3k/2}\min(1,2^j\langle t\rangle^{-1})\|f^\ast_{j,k}\|_{L^2}.
\end{equation}
In addition, if $|t|\geq 1$ and $j\geq\max(-k,0)$, then we have the stronger bounds
\begin{equation}\label{Linfty1.5}
\|\varphi_{[-80,80]}(\langle t\rangle^{-1}x)(e^{-it\Lambda_{wa}}f_{j,k})(x)\|_{L^\infty_x}\lesssim \langle t\rangle^{-1}2^{k/2}(1+\langle t\rangle 2^{k})^{\delta/20}\|f^\ast_{j,k}\|_{H^{0,1}_\Omega};
\end{equation}
\begin{equation}\label{Linfty1.6}
\|e^{-it\Lambda_{wa}}f_{j,k}\|_{L^\infty}\lesssim \langle t\rangle^{-1}2^{k/2}(1+\langle t\rangle 2^{k})^{\delta/20}\|f^\ast_{j,k}\|_{H^{0,1}_\Omega}\qquad \text{ if }\,2^{j}\leq 2^{-10}\langle t\rangle;
\end{equation}
\begin{equation}\label{Linfty1.1}
\|e^{-it\Lambda_{wa}}f_{\leq j,k}\|_{L^\infty}\lesssim 2^{2k}\langle t\rangle^{-1}\|\widehat{f^\ast_{\leq j,k}}\|_{L^\infty}\qquad \text{ if }\,2^{j}\lesssim\langle t\rangle^{1/2} 2^{-k/2}.
\end{equation}

(iii) For any $t\in\mathbb{R}$, $(k,j)\in\mathcal{J}$, and $f\in L^2(\mathbb{R}^3)$ we have
\begin{equation}\label{Linfty3}
\begin{split}
\|e^{-it\Lambda_{kg}}f_{j,k}\|_{L^\infty}\lesssim \min\big\{2^{3k/2},2^{3k^+}\langle t\rangle^{-3/2}2^{3j/2}\big\}\|f^\ast_{j,k}\|_{L^2}.
\end{split}
\end{equation}
Moreover, if $|t|\geq 1$ and $j\geq\max(-k,0)$, then we have the stronger bounds
\begin{equation}\label{Linfty3.6}
\|e^{-it\Lambda_{kg}}f_{j,k}\|_{L^\infty}\lesssim 2^{5k^+}\langle t\rangle^{-3/2}2^{j/2-k^-}(1+\langle t\rangle 2^{2k^-})^{\delta/20}\|f^\ast_{j,k}\|_{H^{0,1}_\Omega}\qquad \text{ if }\,2^{j}\leq 2^{k^--20}\langle t\rangle;
\end{equation}
\begin{equation}\label{Linfty3.1}
\|e^{-it\Lambda_{kg}}f_{\leq j,k}\|_{L^\infty}\lesssim 2^{5k^+}\langle t\rangle^{-3/2}\|\widehat{f^\ast_{\leq j,k}}\|_{L^\infty}\qquad \text{ if }\,2^{j}\lesssim \langle t\rangle^{1/2}.
\end{equation}

(iv) The bounds \eqref{Linfty1.6}, \eqref{Linfty1.1}, \eqref{Linfty3.6} can be improved by using the super-localization operators $\mathcal{C}_{n,l}$ defined in \eqref{suploc}. Indeed, assume that $|t|\geq 1$, $j\geq\max(-k,0)$, and $l\leq k-6$.  Then
\begin{equation}\label{Linfty1.6*}
\Big\{\sum_{n\geq 4}\|e^{-it\Lambda_{wa}}\mathcal{C}_{n,l}f_{j,k}\|_{L^\infty}^2\Big\}^{1/2}\lesssim \langle t\rangle^{-1}2^{l/2}(1+\langle t\rangle 2^{k})^{\delta/20}\|f^\ast_{j,k}\|_{H^{0,1}_\Omega}
\end{equation}
provided that $2^{j}+2^{-l}\lesssim \langle t\rangle(1+\langle t\rangle 2^{k})^{-\delta/20}$. Moreover, if $2^{j}+2^{-l}\lesssim\langle t\rangle^{1/2} 2^{-k/2}$ then
\begin{equation}\label{Linfty1.1*}
\sup_{n\geq 4}\|e^{-it\Lambda_{wa}}\mathcal{C}_{n,l}f_{\leq j,k}\|_{L^\infty}\lesssim 2^{k}2^l\langle t\rangle^{-1}\|\widehat{f^\ast_{\leq j,k}}\|_{L^\infty}.
\end{equation}
Finally, if $2^{j}+2^{-l}\lesssim \langle t\rangle 2^{k^-}(1+\langle t\rangle 2^{2k^-})^{-\delta/20}$ then
\begin{equation}\label{Linfty3.6*}
\Big\{\sum_{n\geq 4}\|e^{-it\Lambda_{kg}}\mathcal{C}_{n,l}f_{j,k}\|_{L^\infty}^2\Big\}^{1/2}\lesssim 2^{5k^+}\langle t\rangle^{-1}2^{l/2}2^{-k^-}(1+\langle t\rangle 2^{2k^-})^{\delta/20}\|f^\ast_{j,k}\|_{H^{0,1}_\Omega}.
\end{equation}
\end{lemma}

\subsection{Bilinear estimates} Linear estimates are insufficient to bound some of the quadratic terms in our nonlinearities. In this subsection we use $TT^\ast$ arguments to prove several additional bilinear estimates involving solutions of wave and Klein-Gordon equations: 

\begin{lemma}\label{Lembil1} Assume $k,k_2\in\mathbb{Z}$, $(k_1,j_1)\in\mathcal{J}$, $t\in\mathbb{R}$, and $f,g\in L^2(\mathbb{R}^3)$. Define $f_{j_1,k_1}$, $f_{\leq j_1,k_1}$, $f^\ast_{j,k}$, $f^\ast_{\leq j,k}$ as in \eqref{defin}--\eqref{defin2}. Assume that $m$ is a bounded multiplier satisfying $\|\mathcal{F}^{-1}m\|_{L^1}\leq 1$ and $I_m$ is a bilinear operator as in \eqref{abc36.1}. If $|t|\geq 1$, and 
\begin{equation}\label{bil1.01}
2^{j_1}\leq \langle t\rangle^{1/2}2^{-k_1/2}+2^{-\min(k,k_1,k_2)}
\end{equation}
then
\begin{equation}\label{bil1}
\|P_k I_m[e^{-it\Lambda_{wa}}f_{\leq j_1,k_1},P_{k_2}g]\|_{L^2}\lesssim 2^{\min(k,k_2)/2}\langle t\rangle^{-1}2^{3k_1/2}\|\widehat{f^\ast_{\leq j_1,k_1}}\|_{L^\infty}\|P_{k_2}g\|_{L^2}.
\end{equation}
\end{lemma}

\begin{proof} By duality we may assume that $k\leq k_2$, and the point is to gain both factors $\langle t\rangle^{-1}$ and $2^{k/2}$ in the right-hand side of \eqref{bil1}. By writing
\begin{equation*}
m(\xi-\eta,\eta)=C\int_{\mathbb{R}^6} K(x,y) e^{-ix\cdot\xi}e^{-i(y-x)\cdot\eta}\,dxdy,\qquad \|K\|_{L^1}\leq 1,
\end{equation*}
and combining the factors $e^{-ix\cdot\xi}$ and $e^{-i(y-x)\cdot\eta}$ with the $L^2$ functions, we may also assume $m\equiv 1$ and write, for simplicity, $I=I_m$. We estimate first 
\begin{equation*}
\|P_kI[e^{-it\Lambda_{wa}}f_{\leq j_1,k_1},P_{k_2}g]\|_{L^2}\lesssim 2^{3k/2}\|f_{\leq j_1,k_1}\|_{L^2}\|P_{k_2}g\|_{L^2}\lesssim 2^{3k/2}2^{3k_1/2}\|\widehat{f_{\leq j_1,k_1}}\|_{L^\infty}\|P_{k_2}g\|_{L^2}.
\end{equation*}
This suffices if $\langle t\rangle 2^k\lesssim 1$. On the other hand, if $\langle t\rangle 2^k\geq 2^{40}$ and $\langle t\rangle 2^{k_1}\leq 2^{40}$ then we estimate
\begin{equation*}
\|P_kI[e^{-it\Lambda_{wa}}f_{\leq j_1,k_1},P_{k_2}g]\|_{L^2}\lesssim 2^{3k_1}\|\widehat{f_{\leq j_1,k_1}}\|_{L^\infty}\|P_{k_2}g\|_{L^2},
\end{equation*}
which suffices. If $\langle t\rangle 2^k\geq 2^{40}$, $\langle t\rangle 2^{k_1}\geq 2^{40}$, and $k_1\leq k+10$ then $2^{-k}\leq 2^{-k_1+10}\leq \langle t\rangle^{1/2}2^{-k_1/2}$. Therefore $2^{j_1}\leq\langle t\rangle^{1/2}2^{-k_1/2+1}$ and the bound \eqref{Linfty1.1} gives
\begin{equation*}
\begin{split}
\|P_kI[e^{-it\Lambda_{wa}}f_{\leq j_1,k_1},P_{k_2}g]\|_{L^2}&\lesssim \|e^{-it\Lambda_{wa}}f_{\leq j_1,k_1}\|_{L^\infty}\|P_{k_2}g\|_{L^2}\\
&\lesssim 2^{2k_1}\langle t\rangle^{-1}\|\widehat{f^\ast_{\leq j_1,k_1}}\|_{L^\infty}\|P_{k_2}g\|_{L^2},
\end{split}
\end{equation*}
which suffices. It remains to prove \eqref{bil1} when 
\begin{equation}\label{bil8}
k\leq k_1-10,\qquad \langle t\rangle 2^{k}\geq 2^{40},\qquad 2^{j_1}\leq 2^{-k}+\langle t\rangle^{1/2}2^{-k_1/2}.
\end{equation}

{\bf{Case 1.}} Assume first that \eqref{bil8} holds and, in addition,
\begin{equation}\label{bil8.1}
2^{-k}\geq \langle t\rangle^{1/2}2^{-k_1/2}.
\end{equation}
In particular, $2^{j_1}\leq 2^{-k+1}$ and $k\leq 1$. We pass to the Fourier space and write
\begin{equation*}
\|P_kI[e^{-it\Lambda_{wa}}f_{\leq j_1,k_1},P_{k_2}g]\|_{L^2}^2=C\int_{\mathbb{R}^3\times\mathbb{R}^3}\widehat{P_{k_2}g}(\eta)\overline{\widehat{P_{k_2}g}(\rho)}L(\eta,\rho)\,d\eta d\rho
\end{equation*}
where 
\begin{equation}\label{bil9}
L(\eta,\rho):=\int_{\mathbb{R}^3}\varphi^2_k(\xi)e^{-it\Lambda_{wa}(\xi-\eta)}\widehat{f_{\leq j_1,k_1}}(\xi-\eta)e^{it\Lambda_{wa}(\xi-\rho)}\overline{\widehat{f_{\leq j_1,k_1}}(\xi-\rho)}\,d\xi.
\end{equation}
Using Schur's lemma, for \eqref{bil1} it suffices to prove that
\begin{equation}\label{bil10}
\begin{split}
&\sup_{\rho\in\mathbb{R}^3}\int_{\mathbb{R}^3}\varphi_{[k_1-4,k_1+4]}(\eta)\varphi_{[k_1-4,k_1+4]}(\rho)|L(\eta,\rho)|\,d\eta\lesssim 2^{k}\langle t\rangle^{-2}2^{3k_1}\|\widehat{f^\ast_{\leq j_1,k_1}}\|_{L^\infty}^2,\\
&\sup_{\eta\in\mathbb{R}^3}\int_{\mathbb{R}^3}\varphi_{[k_1-4,k_1+4]}(\eta)\varphi_{[k_1-2,k_1+4]}(\rho)|L(\eta,\rho)|\,d\rho\lesssim 2^{k}\langle t\rangle^{-2}2^{3k_1}\|\widehat{f^\ast_{\leq j_1,k_1}}\|_{L^\infty}^2.
\end{split}
\end{equation} 
Since $L(\rho,\eta)=\overline{L(\eta,\rho)}$, it suffices to prove the first bound in \eqref{bil10}. 

We would like to integrate by parts in $\xi$ in the integral definition of the kernel $L$. Let $n_0$ denote the smallest integer satisfying $2^{n_0}\geq (2^k\langle t\rangle)^{-1}$ and let
\begin{equation}\label{bil19.1}
\begin{split}
L_n(\eta,\rho):=\int_{\mathbb{R}^3}\varphi_n^{[n_0,\infty)}&(\Lambda'_{wa}(\xi-\eta)-\Lambda'_{wa}(\xi-\rho))\varphi^2_k(\xi)e^{-it(\Lambda_{wa}(\xi-\eta)-\Lambda_{wa}(\xi-\rho))}\\
&\times\widehat{f_{\leq j_1,k_1}}(\xi-\eta)\overline{\widehat{f_{\leq j_1,k_1}}(\xi-\rho)}\,d\xi.
\end{split}
\end{equation}
We may assume that $\|\widehat{f^\ast_{\leq j_1,k_1}}\|_{L^\infty}=1$. Assume $n=n_0+p$, $p\geq 0$.  If $p\geq 1$ then we integrate by parts in $\xi$, using Lemma \ref{tech5} with $K\approx \langle t\rangle 2^n$ and $\eps\approx 2^k$ (recall that $2^{j_1}\lesssim \eps^{-1}$ and $2^{-n-k_1}\lesssim 2^{-k_1}\langle t\rangle 2^k\lesssim \eps^{-1}$, due to \eqref{bil8.1}). It follows that, for all $p\geq 0$,
\begin{equation*}
\begin{split}
|L_n(\eta,\rho)|\lesssim 2^{-4p}\int_{\mathbb{R}^3}&\varphi_{\leq n+4}(\Lambda'_{wa}(\xi-\eta)-\Lambda'_{wa}(\xi-\rho))\varphi_{[k-4,k+4]}(\xi)\\
&\times\varphi_{[k_1-4,k_1+4]}(\xi-\eta)\varphi_{[k_1-4,k_1+4]}(\xi-\rho)\,d\xi.
\end{split}
\end{equation*} 
Therefore, after changes of variables,
\begin{equation*}
\begin{split}
\int_{\mathbb{R}^3}\varphi_{[k_1-4,k_1+4]}(\eta)&|L_n(\eta,\rho)|\,d\eta\lesssim 2^{-4p}\int_{\mathbb{R}^3\times\mathbb{R}^3}\varphi_{\leq n+4}(\Lambda'_{wa}(y)-\Lambda'_{wa}(x))\varphi_{[k-4,k+4]}(\rho+x)\\
&\times\varphi_{[k_1-4,k_1+4]}(\rho+x-y)\varphi_{[k_1-4,k_1+4]}(y)\varphi_{[k_1-4,k_1+4]}(x)\,dxdy,
\end{split}
\end{equation*}
for any $\rho\in\mathbb{R}^3$ with $|\rho|\in [2^{k_1-6},2^{k_1+6}]$. Since $\Lambda'_{wa}(z)=z/|z|$, the integration in $y$ in the expression above is essentially in a rectangle of sides smaller than $2^n2^{k_1}\times 2^n2^{k_1}\times 2^{k_1}$. Thus
\begin{equation*}
\int_{\mathbb{R}^3}\varphi_{[k_1-4,k_1+4]}(\eta)|L_n(\eta,\rho)|\,d\eta\lesssim 2^{-4p}2^{2n}2^{3k_1}2^{3k}\lesssim 2^{-2p}\langle t\rangle^{-2}2^{3k_1}2^{k}.
\end{equation*}
The desired conclusion \eqref{bil10} follows.

{\bf{Case 2.}} Assume now that \eqref{bil8} holds and, in addition,
\begin{equation}\label{bil9.1}
2^{-k}\leq \langle t\rangle^{1/2}2^{-k_1/2}.
\end{equation}
We fix a smooth function $\chi:\mathbb{R}\to[0,1]$ supported in $[-2,2]$ with the property that $\sum_{n\in\mathbb{Z}}\chi(x-n)=1$ for all $x\in\mathbb{R}$. Then we decompose
\begin{equation*}
\begin{split}
&f_{\leq j_1,k_1}=\sum_{n}f_{\leq j_1,k_1;n},\qquad \widehat{f_{\leq j_1,k_1;n}}(\xi):=\widehat{f_{\leq j_1,k_1}}(\xi)\chi(2^{-k}|\xi|-n).
\end{split}
\end{equation*}
Let $\widehat{g_{k_2;n}}(\xi):=\widehat{P_{k_2}g}(\xi)\varphi_{\leq 4}(2^{-k}|\xi|-n)$. Clearly 
\begin{equation*}
P_kI[e^{-it\Lambda_{wa}}f_{\leq j_1,k_1},P_{k_2}g]=\sum_{n}P_k I[e^{-it\Lambda_{wa}}f_{\leq j_1,k_1;n},g_{k_2;n}].
\end{equation*} 
The sum in $n$ has at most $C2^{k_1-k}$ nontrivial terms, so, by orthogonality,
\begin{equation*}
\begin{split}
\|P_k I[e^{-it\Lambda_{wa}}f_{\leq j_1,k_1},P_{k_2}g]\|_{L^2}&\lesssim \sum_{n}\|I[e^{-it\Lambda_{wa}}f_{\leq j_1,k_1;n},g_{k_2;n}]\|_{L^2}\\
&\lesssim 2^{(k_1-k)/2}\Big\{\sum_{n}\|e^{-it\Lambda_{wa}}f_{\leq j_1,k_1;n}\|_{L^\infty}^2\|g_{k_2;n}\|_{L^2}^2\Big\}^{1/2}\\
&\lesssim 2^{(k_1-k)/2}\|P_{k_2}g\|_{L^2}\sup_{n}\|e^{-it\Lambda_{wa}}f_{\leq j_1,k_1;n}\|_{L^\infty}.
\end{split}
\end{equation*} 
For \eqref{bil1} it suffices to prove that, for any $n\approx 2^{k_1-k}$,
\begin{equation*}
\|e^{-it\Lambda_{wa}}f_{\leq j_1,k_1;n}\|_{L^\infty}\lesssim 2^k2^{k_1}\langle t\rangle^{-1}\|\widehat{f^\ast_{\leq j_1,k_1}}\|_{L^\infty},
\end{equation*}
which follows from \eqref{Linfty1.1*}.
\end{proof}

We also have some variants using only rotational derivatives:

\begin{lemma}\label{Lembil2} Assume $k,k_1,k_2\in\mathbb{Z}$, $(k_1,j_1),(k_2,j_2)\in\mathcal{J}$, $t\in\mathbb{R}$, $|t|\geq 1$, and $f,g\in L^2(\mathbb{R}^3)$. Define $f_{j_1,k_1}, f^\ast_{j_1,k_1}, g_{j_2,k_2}, g^\ast_{j_2,k_2}$ as in \eqref{defin}--\eqref{defin2}. 

(i) If $m$ is a bounded multiplier satisfying $\|\mathcal{F}^{-1}m\|_{L^1}\leq 1$, $I_m$ is the associated bilinear operator as in \eqref{abc36.1}, and
\begin{equation}\label{ret1}
2^{j_1}\lesssim\langle t\rangle(1+2^{k_1}\langle t\rangle)^{-\delta/20}+2^{-\min(k,k_1,k_2)}
\end{equation}
then
\begin{equation}\label{Bil31}
\|P_k I_m[e^{-it\Lambda_{wa}}f_{j_1,k_1},P_{k_2}g]\|_{L^2}\lesssim 2^{\min(k,k_2)/2}\langle t\rangle^{-1}(1+\langle t\rangle2^{k_1})^{\delta/20}\|f^\ast_{j_1,k_1}\|_{H^{0,1}_\Omega}\|P_{k_2}g\|_{L^2}.
\end{equation}

(ii) If $m\in\mathcal{M}$ and $I_m$ is the associated bilinear operator as in \eqref{mults}--\eqref{abc36.1}, $\iota_2\in\{+,-\}$,
\begin{equation}\label{rex1}
2^{k_1},2^{k_2}\in[\langle t\rangle^{-1+\delta/2},\langle t\rangle^{2/\delta}]\qquad\text{ and }\qquad 2^{j_2}\leq \langle t\rangle^{1-\delta/2},
\end{equation}
then, with $g_{j_2,k_2}^+:=g_{j_2,k_2}$ and $g_{j_2,k_2}^-:=\overline{g_{j_2,k_2}}$,
\begin{equation}\label{rex2}
\|I_m[e^{-it\Lambda_{wa}}f_{j_1,k_1},e^{-it\Lambda_{wa,\iota_2}}g_{j_2,k_2}^{\iota_2}]\|_{L^2}\lesssim 2^{k_1/2}\langle t\rangle^{-1+\delta/2}\|g^\ast_{j_2,k_2}\|_{L^2}\|f^\ast_{j_1,k_1}\|_{H^{0,1}_\Omega}.
\end{equation}
\end{lemma}

\begin{proof} (i) As before, by duality we may assume that $k\leq k_2$ and the point is to gain both factors $\langle t\rangle^{-1}$ and $2^{k/2}$. We may also assume that $m\equiv 1$ and write $I=I_m$. 

We estimate first, using just the Cauchy-Schwarz inequality,
\begin{equation*}
\|P_kI[e^{-it\Lambda_{wa}}f_{j_1,k_1},P_{k_2}g]\|_{L^2}\lesssim 2^{3\min(k,k_1)/2}\|f_{j_1,k_1}\|_{L^2}\|P_{k_2}g\|_{L^2}.
\end{equation*}
This suffices to prove \eqref{Bil31} if $2^{\min(k,k_1)}\lesssim\langle t\rangle^{-1}(1+\langle t\rangle2^{k_1})^{\delta/20}$. On the other hand, if $2^{\min(k,k_1)}\gg\langle t\rangle^{-1}(1+\langle t\rangle2^{k_1})^{\delta/20}$ and $k_1\leq k+20$ then we use \eqref{Linfty1.6} to estimate
\begin{equation*}
\begin{split}
\|P_kI[e^{-it\Lambda_{wa}}f_{j_1,k_1},P_{k_2}g]\|_{L^2}&\lesssim \|e^{-it\Lambda_{wa}}f_{j_1,k_1}\|_{L^\infty}\|P_{k_2}g]\|_{L^2}\\
&\lesssim \langle t\rangle^{-1}2^{k_1/2}(1+\langle t\rangle2^{k_1})^{\delta/20}\|f^\ast_{j_1,k_1}\|_{H^{0,1}_\Omega}\|P_{k_2}g\|_{L^2},
\end{split}
\end{equation*}
which suffices. It remains to prove \eqref{Bil31} when 
\begin{equation}\label{Bil32}
k\leq k_1-20,\qquad \langle t\rangle 2^{k}\geq (1+\langle t\rangle2^{k_1})^{\delta/20},\qquad 2^{j_1}\lesssim\langle t\rangle(1+2^{k_1}\langle t\rangle)^{-\delta/20}.
\end{equation}

We decompose
\begin{equation*}
\begin{split}
&f_{j_1,k_1}=\sum_{n\geq 4}f_{j_1,k_1;n},\qquad f_{j_1,k_1;n}:=\mathcal{C}_{n,k}f_{j_1,k_1}.
\end{split}
\end{equation*}
Let $\widehat{g_{k_2;n}}(\xi):=\widehat{P_{k_2}g}(\xi)\varphi_{\leq 4}(2^{-k}|\xi|-n)$. Clearly 
\begin{equation*}
P_kI[e^{-it\Lambda_{wa}}f_{j_1,k_1},P_{k_2}g]=\sum_{n\geq 4}P_k I[e^{-it\Lambda_{wa}}f_{j_1,k_1;n},g_{k_2;n}].
\end{equation*} 
Therefore
\begin{equation*}
\begin{split}
\|P_k I[e^{-it\Lambda_{wa}}f_{j_1,k_1},&P_{k_2}g]\|_{L^2}\lesssim \sum_{n\geq 4}\big\|e^{-it\Lambda_{wa}}f_{j_1,k_1;n}\|_{L^\infty}\big\|g_{k_2;n}\big\|_{L^2}\\
&\lesssim \Big\{\sum_{n\geq 4}\|e^{-it\Lambda_{wa}}f_{j_1,k_1;n}\|_{L^\infty}^2\Big\}^{1/2}\Big\{\sum_{n\geq 4}\|g_{k_2;n}\|_{L^2}^2\Big\}^{1/2}\\
&\lesssim \langle t\rangle^{-1}2^{k/2}(1+\langle t\rangle 2^{k_1})^{\delta/20}\|f^\ast_{j_1,k_1}\|_{H^{0,1}_\Omega}\cdot\|P_{k_2}g\|_{L^2},
\end{split}
\end{equation*} 
where we used \eqref{Linfty1.6*} (see the restrictions \eqref{Bil32}) and orthogonality in the last inequality. This completes the proof of \eqref{Bil31}.

(ii) We decompose $e^{-it\Lambda_{wa}}f_{j_1,k_1}=P_{[k_1-4,k_1+4]}H_{k_1}^1+P_{[k_1-4,k_1+4]}H_{k_1}^2$, where
\begin{equation*}
H_{k_1}^1(x):=\varphi_{[-40,40]}(x/\langle t\rangle)e^{-it\Lambda_{wa}}f_{j_1,k_1}(x),\quad H_{k_1}^2(x):=(1-\varphi_{[-40,40]}(x/\langle t\rangle))e^{-it\Lambda_{wa}}f_{j_1,k_1}(x),
\end{equation*}
for $x\in\mathbb{R}^3$. In view of \eqref{Linfty1.5}, we have
\begin{equation*}
\|H_{k_1}^1\|_{L^\infty}\lesssim 2^{k_1/2}\langle t\rangle^{-1+\delta/2}\|f^\ast_{j_1,k_1}\|_{H^{0,1}_\Omega},
\end{equation*}
so the contribution of $H_{k_1}^1$ is bounded as claimed. 

On the other hand, we claim that the contribution of $H_{k_1}^2$ is negligible, 
\begin{equation}\label{rex4}
\|I_m[P_{[k_1-4,k_1+4]}H_{k_1}^2,e^{-it\Lambda_{wa,\iota_2}}g_{j_2,k_2}^{\iota_2}]\|_{L^2}\lesssim \langle t\rangle^{-2}\|g^\ast_{j_2,k_2}\|_{L^2}\|f^\ast_{j_1,k_1}\|_{L^2}.
\end{equation}
Indeed, the definitions \eqref{mults} and \eqref{abc36.1} show that
\begin{equation*}
|I_m[P_{l_1}F,P_{l_2}G](x)|\lesssim_N(|F|\ast K^N_{l_1})(x)\cdot (|G|\ast K^N_{l_2})(x),
\end{equation*}
for any $l_1,l_2\in\mathbb{Z}$ and $N\geq 10$, where $K_l^N(y):=2^{3l}(1+2^{l^-}|y|)^{-N}$. Moreover
\begin{equation*}
\|\varphi_{[-30,30]}(x/\langle t\rangle)\cdot (|H_{k_1}^2|\ast K_{l_1}^N)(x)\|_{L^\infty}\lesssim \langle t\rangle^{-6}\|f^\ast_{j_1,k_1}\|_{L^2},
\end{equation*}
for $l_1\in[k_1-4,k_1+4]$, in view of the support restriction on $H_{k_1}^2$ and the assumption $2^{k_1^-}\langle t\rangle\gtrsim  \langle t\rangle^{\delta/2}$. Also, using Lemma \ref{tech5} and the assumption $2^{j_2}\leq \langle t\rangle^{1-\delta/2}$, we have
\begin{equation*}
\|(1-\varphi_{[-30,30]}(x/\langle t\rangle))\cdot (|e^{-it\Lambda_{wa}}g_{j_2,k_2}|\ast K_{l_2}^N)(x)\|_{L^2}\lesssim \langle t\rangle^{-6}\|g^\ast_{j_2,k_2}\|_{L^2},
\end{equation*}
for $l_2\in[k_2-4,k_2+4]$. The desired estimates \eqref{rex4} follow from the last three bounds.
\end{proof}

We also need some bilinear estimates involving the Klein--Gordon flow:

\begin{lemma}\label{LemBil4} Assume $k,k_1,k_2\in\mathbb{Z}$, $(k_1,j_1)\in\mathcal{J}$, $t\in\mathbb{R}$, $|t|\geq 1$, and $f,g\in L^2(\mathbb{R}^3)$. Define $f_{j_1,k_1}, f^\ast_{j_1,k_1}$ as in \eqref{defin}--\eqref{defin2}. If $\|\mathcal{F}^{-1}m\|_{L^1}\leq 1$, $I_m$ is the bilinear operator as in \eqref{abc36.1},
\begin{equation}\label{LemBil42}
k\leq k^-_1-10\qquad\text{ and }\qquad 2^{j_1}\lesssim\langle t\rangle2^{k_1^-}(1+2^{2k_1^-}\langle t\rangle)^{-\delta/20}
\end{equation}
then
\begin{equation}\label{LemBil43}
\|P_k I_m[e^{-it\Lambda_{kg}}f_{j_1,k_1},P_{k_2}g]\|_{L^2}\lesssim 2^{k/2}\langle t\rangle^{-1}2^{5k_1^+}2^{-k_1^-}(1+\langle t\rangle2^{2k_1^-})^{\delta/20}\|f^\ast_{j_1,k_1}\|_{H^{0,1}_\Omega}\|P_{k_2}g\|_{L^2}.
\end{equation}
\end{lemma}

\begin{proof} This is similar to the proof of Lemma \ref{Lembil2} (i). We may assume $m\equiv 1$ and write $I=I_m$. We estimate first, using just the Cauchy-Schwarz inequality,
\begin{equation*}
\|P_kI[e^{-it\Lambda_{kg}}f_{j_1,k_1},P_{k_2}g]\|_{L^2}\lesssim 2^{3k/2}\|f_{j_1,k_1}\|_{L^2}\|P_{k_2}g\|_{L^2},
\end{equation*}
which suffices if $2^{k+k_1^-}\lesssim\langle t\rangle^{-1}(1+\langle t\rangle2^{2k_1^-})^{\delta/20}$. On the other hand, if $2^{k+k_1^-}\gg\langle t\rangle^{-1}(1+\langle t\rangle2^{2k_1^-})^{\delta/20}$ and $k\leq k_1^--10$ then $2^{2k_1^-}\gg\langle t\rangle^{-1}$ and we decompose
\begin{equation*}
\begin{split}
&f_{j_1,k_1}=\sum_{n\geq 4}f_{j_1,k_1;n},\qquad f_{j_1,k_1;n}:=\mathcal{C}_{n,k}f_{j_1,k_1}(\xi).
\end{split}
\end{equation*}
Let $\widehat{g_{k_2;n}}(\xi):=\widehat{P_{k_2}g}(\xi)\varphi_{\leq 4}(2^{-k}|\xi|-n)$. Clearly 
\begin{equation*}
P_kI[e^{-it\Lambda_{kg}}f_{j_1,k_1},P_{k_2}g]=\sum_{n\geq 4}P_k I[e^{-it\Lambda_{wa}}f_{j_1,k_1;n},g_{k_2;n}].
\end{equation*} 
Therefore, as in the proof of Lemma \ref{Lembil2} (i),
\begin{equation*}
\begin{split}
\|P_k I[e^{-it\Lambda_{kg}}f_{j_1,k_1},&P_{k_2}g]\|_{L^2}\lesssim \sum_{n\geq 4}\big\|e^{-it\Lambda_{kg}}f_{j_1,k_1;n}\|_{L^\infty}\big\|g_{k_2;n}\big\|_{L^2}\\
&\lesssim \Big\{\sum_{n\geq 4}\|e^{-it\Lambda_{kg}}f_{j_1,k_1;n}\|_{L^\infty}^2\Big\}^{1/2}\Big\{\sum_{n\geq 4}\|g_{k_2;n}\|_{L^2}^2\Big\}^{1/2}\\
&\lesssim \langle t\rangle^{-1}2^{k/2}2^{5k_1^+}2^{-k_1^-}(1+\langle t\rangle 2^{2k_1^-})^{\delta/20}\|f^\ast_{j_1,k_1}\|_{H^{0,1}_\Omega}\cdot\|P_{k_2}g\|_{L^2},
\end{split}
\end{equation*} 
where we used \eqref{Linfty3.6*} (see \eqref{LemBil42} and recall that $2^{k+k_1^-}\gg\langle t\rangle^{-1}(1+\langle t\rangle2^{2k_1^-})^{\delta/20}$) and ortho\-go\-nality in the last inequality. This completes the proof of \eqref{LemBil43}.
\end{proof}

\subsection{Interpolation inequalities}\label{inte0}

Finally, we need some interpolation bounds involving $L^p$ spaces and rotation vector-fields.

\begin{lemma}\label{inte1}
(i) Assume that $f\in H^{0,1}_\Omega$, $k\in\mathbb{Z}$, and $A_0\leq A_1\leq B\in [0,\infty)$. If 
\begin{equation}\label{consu}
\|Q_{j,k}f\|_{L^2}\leq A_0,\qquad \|Q_{j,k}f\|_{H^{0,1}_\Omega}\leq A_1,\qquad 2^{j+k}\|Q_{j,k}f\|_{H^{0,1}_\Omega}\leq B
\end{equation}
for all $j\geq -k^-$, then
\begin{equation}\label{consu2.1}
\|\widehat{P_k f}\|_{L^\infty}\lesssim 2^{-3k/2}A_0^{(1-\delta)/4}B^{(3+\delta)/4}
\end{equation}
and
\begin{equation}\label{consu2}
\|\widehat{P_k f}\|_{L^\infty}\lesssim 2^{-3k/2}A_1^{(1-\delta)/2}B^{(1+\delta)/2}.
\end{equation}

(ii) If $f\in H^{0,2}_\Omega$, $\Omega\in\{\Omega_{23},\,\Omega_{31},\,\Omega_{12}\}$, and $k\in\mathbb{Z}$ then
\begin{equation}\label{inte4}
\|P_k\Omega f\|_{L^4}\lesssim \|P_k f\|^{1/2}_{L^\infty}\|P_k f\|^{1/2}_{H^{0,2}_\Omega}.
\end{equation}
Similarly, if $f\in H^{0,3}_\Omega$ and $\Omega^2\in\{\Omega_{23}^{a_1}\Omega_{31}^{a_2}\Omega_{12}^{a_3}:\,a_1+a_2+a_3=2\}$ then
\begin{equation}\label{inte5}
\begin{split}
\|P_k\Omega f\|_{L^6}&\lesssim \|P_k f\|^{2/3}_{L^\infty}\|P_k f\|^{1/3}_{H^{0,3}_\Omega},\\
\|P_k\Omega^2 f\|_{L^3}&\lesssim \|P_k f\|^{1/3}_{L^\infty}\|P_k f\|^{2/3}_{H^{0,3}_\Omega}.
\end{split}
\end{equation}
Finally, we have the $L^2$ interpolation estimates
\begin{equation}\label{vcx22}
\Vert P_k\Omega f\Vert_{L^2}\lesssim \Vert P_kf\Vert_{L^2}^{1/2}\Vert P_kf\Vert_{H^{0,2}_\Omega}^{1/2}
\end{equation}
and
\begin{equation}\label{vcx22.5}
\begin{split}
\Vert P_k\Omega f\Vert_{L^2}&\lesssim \Vert P_kf\Vert_{L^2}^{2/3}\Vert P_kf\Vert_{H^{0,3}_\Omega}^{1/3},\\
\Vert P_k\Omega^2f\Vert_{L^2}&\lesssim \Vert P_kf\Vert_{L^2}^{1/3}\Vert P_kf\Vert_{H^{0,3}_\Omega}^{2/3}.
\end{split}
\end{equation}

\end{lemma}

\begin{proof} (i) The bounds follow from \eqref{Linfty3.3}: with $f_{j,k}=P_{k'}Q_{j,k}f$ we have
\begin{equation*}
\begin{split}
\|\widehat{f_{j,k}}\|_{L^\infty}&\lesssim \min\{2^{3j/2}\|Q_{j,k}f\|_{L^2},2^{j/2-k}2^{\delta(j+k)/20}\|Q_{j,k}f\|_{H^{0,1}_\Omega}\}\\
&\lesssim 2^{-3k/2}\min\{2^{3(j+k)/2}A_0,2^{-(j+k)/2}2^{\delta(j+k)/20}B\}.
\end{split}
\end{equation*}
The desired bounds \eqref{consu2.1} follow by summing over $j$ and considering the two cases $2^{j+k}\leq (B/A_0)^{1/2}$ and $2^{j+k}\geq (B/A_0)^{1/2}$. Similarly,
\begin{equation*}
\|\widehat{f_{j,k}}\|_{L^\infty}\lesssim 2^{j/2-k}2^{\delta(j+k)/20}\|Q_{j,k}f\|_{H^{0,1}_\Omega}\lesssim 2^{-3k/2}2^{(j+k)(1/2+\delta/20)}\min\{A_1,B2^{-(j+k)}\}
\end{equation*}
The desired bounds \eqref{consu2} follow again by summing over $j$.

(ii) For \eqref{inte4} we let $g:=P_kf$ and use integration by parts to write
\begin{equation*}
\|\Omega g\|_{L^4}^4=\int_{\mathbb{R}^3}\Omega g\,\Omega g\,\overline{\Omega g}\,\overline{\Omega g}\,dx=-\int_{\mathbb{R}^3}g\cdot \Omega\{\Omega g\,\overline{\Omega g}\,\overline{\Omega g}\}\,dx.
\end{equation*}
Therefore we can estimate 
\begin{equation*}
\|\Omega g\|_{L^4}^4\lesssim \int_{\mathbb{R}^3}|g|\,|\Omega g|^2\,|\Omega^2 g|\,dx\lesssim \|g\|_{L^\infty}\|\Omega^2 g\|_{L^2}\|\Omega g\|_{L^4}^2,
\end{equation*}
which gives \eqref{inte4}.

Similarly, to prove \eqref{inte5} we estimate as above
\begin{equation}\label{inte6}
\|\Omega g\|_{L^6}^6\lesssim \int_{\mathbb{R}^3}|g|\,|\Omega g|^4\,|\Omega^2 g|\,dx\lesssim \|g\|_{L^\infty}\|\Omega^2 g\|_{L^3}\|\Omega g\|_{L^6}^4
\end{equation}
and
\begin{equation*}
\|\Omega^2 g\|_{L^3}^3=\Big|\int_{\mathbb{R}^3}\Omega^2 g\,\overline{\Omega^2 g}(\Omega^2 g\overline{\Omega^2 g})^{1/2}\,dx\Big|\lesssim \int_{\mathbb{R}^3}|\Omega g|\,|\Omega^2 g|\,|\Omega^3 g|\,dx\lesssim \|g\|_{H^{0,3}_\Omega}\|\Omega^2 g\|_{L^3}\|\Omega g\|_{L^6},
\end{equation*}
where $\Omega^3g$ denote vector-fields of the form  $\Omega_{23}^{a_1}\Omega_{31}^{a_2}\Omega_{12}^{a_3}$ with $a_1+a_2+a_3=3$. Therefore, using also \eqref{inte6} and simplifying,
\begin{equation*}
\|\Omega^2 g\|_{L^3}^2\lesssim \|g\|_{H^{0,3}_\Omega}\|\Omega g\|_{L^6}\lesssim \|g\|_{H^{0,3}_\Omega}\big(\|g\|_{L^\infty}\|\Omega^2 g\|_{L^3}\big)^{1/2},
\end{equation*}
which gives the bounds in the second line of \eqref{inte5}. The bounds in the first line now follow from \eqref{inte6}.

The $L^2$ bounds \eqref{vcx22} and \eqref{vcx22.5} follow by similar arguments.
\end{proof}

We need also a bilinear interpolation lemma:

\begin{lemma}\label{LemBil5}
Assume $k,k_1,k_2\in\mathbb{Z}$, $k+4\leq K:=\min(k_1,k_2)$, $t\in\mathbb{R}$, and $f,g\in L^2(\mathbb{R}^3)$. Assume $J\geq \max(-K,0)$ satisfies 
\begin{equation}\label{Bil5.1}
2^{J}\leq \langle t\rangle^{1/2}2^{-K/2}+2^{-k},
\end{equation}
and define $f_{\leq J,k_1}$ and $g_{\leq J,k_2}$ as in \eqref{defin}. Then
\begin{equation}\label{Bil5}
\Vert P_kI[e^{-it\Lambda_{wa}}f_{\leq J,k_1},e^{-it\Lambda_{wa}}g_{\le J,k_2}]\Vert_{L^2}\lesssim 2^{k/2}2^{3K/2}\langle t\rangle^{-1}\Vert \widehat{P_{k_1}f}\Vert_{L^a}\Vert \widehat{P_{k_2}g}\Vert_{L^b},
\end{equation}
for all exponents $a,b\in[2,\infty]$ satisfying $1/a+1/b=1/2$.
\end{lemma}

\begin{proof} The conclusion follows directly from Lemma \ref{Lembil1} (which corresponds to the cases $(a,b)=(\infty,2)$ and $(a,b)=(2,\infty)$) and bilinear interpolation.
\end{proof}

\section{Analysis of the linear profiles}\label{dtv}

In this section we use the main bootstrap assumptions \eqref{bootstrap2.1}--\eqref{bootstrap2.4} and Lemma \ref{LinEstLem} to derive many linear bounds on the profiles $V^\ast$ and the normalized sultions $U^\ast$.

For $t\in[0,T]$, $(j,k)\in\mathcal{J}$, $J\geq\max(-k,0)$, and
\begin{equation*}
X\in\big\{F,\underline{F},\rho,\omega_a,\Omega_a,\vartheta_{ab},\mathcal{L}h_{\al\be},\mathcal{L}\psi:\,\mathcal{L}\in\mathcal{V}_{3}^3\big\},
\end{equation*}
we define the profiles $V^{X,\pm}$ as in \eqref{notationL}. If $\mathcal{L}\in\mathcal{V}_2^2$ we define also the {\it{space-localized profiles}}
\begin{equation}\label{on11.3}
V_{j,k}^{X,\pm}(t):=P'_{k}Q_{j,k}V^{X,\pm}(t),\quad V_{\leq J,k}^{X,\pm}(t):=\sum_{j\leq J}V_{j,k}^{X,\pm}(t),\quad V_{> J,k}^{X,\pm}(t):=\sum_{j> J}V_{j,k}^{X,\pm}(t).
\end{equation}
and the associated localized solutions
\begin{equation}\label{on11.36}
U_{j,k}^{X,\pm}(t):=e^{-it\Lambda_{\mu,\pm}}V_{j,k}^{X,\pm}(t),\quad U_{\leq J,k}^{X,\pm}(t):=\sum_{j\leq J}U_{j,k}^{X,\pm}(t),\quad U_{> J,k}^{X,\pm}(t):=\sum_{j> J}U_{j,k}^{X,\pm}(t),
\end{equation}
where $\mu=kg$ if $X=\mathcal{L}\psi$ and $\mu=wa$ otherwise. For simplicity of notation, we sometimes let $V^{X}_{\ast}:=V^{X,+}_{\ast}$ and $U^{X}_{\ast}:=U^{X,+}_{\ast}$, and notice that $V^{X,-}_{\ast}=\overline{V^{X}_{\ast}}$ and $U^{X,-}_{\ast}=\overline{U^{X}_{\ast}}$.

Given two pairs $(q,n),(q',n')\in\mathbb{Z}_+^2$ we write $(q,n)\leq (q',n')$ if $q\leq q'$ and $n\leq n'$. We collect first several bounds on our profiles. In the lemmas below we let $h$ denote generic components of the linearized metric, i.e.  $h\in\big\{h_{\al\be}:\al,\be\in\{0,1,2,3\}\big\}$.

\begin{lemma}\label{dtv2}
Assume that $(g,\psi)$ is a solution of the system \eqref{asum3}-\eqref{asum5} on some time interval $[0,T]$, $T\geq 1$, satisfying the bootstrap hypothesis  \eqref{bootstrap2.1}-\eqref{bootstrap2.4}. 

(i) For any $t\in[0,T]$ and $\mathcal{L}\in\mathcal{V}_n^q$, $n\leq 3$,  we have
\begin{equation}\label{vcx1}
\|(\langle t\rangle|\nabla|_{\leq 1})^{\gamma}|\nabla|^{-1/2}V^{\mathcal{L}h}(t)\|_{H^{N(n)}}+\|V^{\mathcal{L}\psi}(t)\|_{H^{N(n)}}\lesssim \varep_1\langle t\rangle^{H(q,n)\delta}.
\end{equation}
In addition, if $(k,j)\in\mathcal{J}$ and  $\mathcal{L}\in\mathcal{V}_n^q$, $n\leq 2$, then
\begin{equation}\label{vcx1.1}
2^{k/2}(2^{k^-}\langle t\rangle)^{\gamma}2^j\|Q_{j,k}V^{\mathcal{L}h}(t)\|_{L^2}+2^{k^+}2^j\|Q_{j,k}V^{\mathcal{L}\psi}(t)\|_{L^2}\lesssim \varep_1Y(k,t;q,n),
\end{equation}
and
\begin{equation}\label{cnb2}
2^{k^+}\|P_kV^{\mathcal{L}\psi}(t)\|_{L^2}+\|P_kV^{\mathcal{L}\psi}(t)\|_{H^{0,1}_\Omega}\lesssim \varep_12^{k^-}Y(k,t;q,n),
\end{equation}
where, for $q,n\in\{0,1,2\}$,
\begin{equation}\label{wws3}
Y(k,t;q,n):=\langle t\rangle^{H(q+1,n+1)\delta}2^{-N(n+1)k^+}.
\end{equation}

(ii) For any $H\in\big\{F,\omega_a,\vartheta_{ab}:\,a,b\in\{1,2,3\}\big\}$ and $k\in\mathbb{Z}$
\begin{equation}\label{vcx1.2}
\begin{split}
\|\widehat{P_kV^{H}}(t)\|_{L^\infty}+\langle t\rangle^{-\delta}\|\widehat{P_kV^{h}}(t)\|_{L^\infty}&\lesssim \varep_12^{-k^--\kappa k^-}2^{-N_0 k^+}\\
\|\widehat{P_kV^{\psi}}(t)\|_{L^\infty}&\lesssim \varep_12^{-k^-/2+\kappa k^-}2^{-N_0k^+}.
\end{split}
\end{equation}
Moreover, if $(k,j)\in\mathcal{J}$ and $H'\in\big\{F,\underline{F},\rho,\omega_a,\Omega_a,\vartheta_{ab}\}$ then
\begin{equation}\label{vcx1.2*}
\begin{split}
\|P_kV^{H'}(t)\|_{H^{0,a}_\Omega}&\lesssim \varep_1\langle t\rangle^{H(0,a)\delta}2^{-N(a)k^+}2^{k/2}(\langle t\rangle 2^{k^-})^{-\gamma},\quad\text{ if }a\leq 3,\\
2^j\|Q_{j,k}V^{H'}(t)\|_{H^{0,a}_\Omega}&\lesssim\varep_1\langle t\rangle^{H(1,a+1)\delta}2^{-N(a+1)k^+}2^{-k/2}(\langle t\rangle 2^{k^-})^{-\gamma},\quad\text{ if }a\leq 2.
\end{split}
\end{equation}
\end{lemma} 

\begin{proof} The bounds \eqref{vcx1} follow directly from \eqref{bootstrap2.1}, and the bounds \eqref{vcx1.2} follow from \eqref{bootstrap2.4} and Definition \ref{MainZDef}. 

It follows from \eqref{bootstrap2.2} that
\begin{equation}\label{cnb1}
2^{k/2}(2^{k^-}\langle t\rangle)^{\ga}\|\varphi_k(\xi)(\partial_{\xi_l}\widehat{V^{\mathcal{L}h}})(\xi,t)\|_{L^2_\xi}+2^{k^+}\|\varphi_k(\xi)(\partial_{\xi_l}\widehat{V^{\mathcal{L}\psi}})(\xi,t)\|_{L^2_\xi}\lesssim \varep_1Y(k,t;q,n),
\end{equation}
for any $k\in\mathbb{Z}$, $l\in\{1,2,3\}$, and $\mathcal{L}\in\mathcal{V}_n^q$, $(q,n)\leq (2,2)$.
The bounds in \eqref{vcx1.1} follow using also \eqref{vcx1} and Lemma \ref{hyt1} (i). The bounds \eqref{cnb2} follow from \eqref{cnb1} and the estimates 
\begin{equation*}
\Vert P_kV^{\mathcal{L}\psi}(t)\Vert_{H^{0,1}_\Omega}\lesssim 2^k\|\varphi_k(\xi)\nabla_\xi(\widehat{V^{\mathcal{L}\psi}})(\xi,t)\|_{L^2_\xi}\lesssim\varep_12^{k^-}Y(k,t;q,n).
\end{equation*}

To prove \eqref{vcx1.2*} we notice first that if $\Omega\in\{\Omega_{23},\Omega_{31},\Omega_{12}\}$ then
\begin{equation}\label{vcx1new2}
\Omega Q_{j,k}V^{X}=Q_{j,k}\Omega V^{X}=Q_{j,k}V^{\Omega X}\qquad\text{ and }\qquad\Omega P_kV^{X}=P_k\Omega V^{X}=P_kV^{\Omega X},
\end{equation}
for suitable profiles $X$. Recall that the functions $H'$ are defined by taking Riesz transforms of the metric components $h_{\al\be}$ (see \eqref{zaq2}). The bounds in the first line of \eqref{vcx1.2*} follow from \eqref{vcx1}, while the bounds in the second line follow from \eqref{vcx1.1} and Lemma \ref{hyt1} (ii). 
\end{proof}

We prove now several pointwise decay bounds on the normalized solutions.

\begin{lemma}\label{dtv3}
(i) For any $k\in\mathbb{Z}$, $t\in[0,T]$, and $\mathcal{L}\in\mathcal{V}_n^q$, $n\leq 2$, we have
\begin{equation}\label{wws1}
\sum_{j\geq -k^-}\|U_{j,k}^{\mathcal{L}h}(t)\|_{L^\infty}\lesssim \varep_1\langle t\rangle^{-1+H(q+1,n+1)\delta}2^{k^-}2^{-N(n+1)k^++2k^+}\min\{1,\langle t\rangle 2^{k^-}\}^{1-\delta},
\end{equation}
where $h\in\big\{h_{\al\be}:\al,\be\in\{0,1,2,3\}\big\}$ as before. In addition, if $2^{k^-}\langle t\rangle\geq 2^{20}$ then
\begin{equation}\label{wws11}
\sum_{2^j\in [2^{-k^-},2^{-20}\langle t\rangle]}\|U_{j,k}^{\mathcal{L}h}(t)\|_{L^\infty}\lesssim \varep_1\langle t\rangle^{-1+H(q,n+1)\delta}2^{k^-}2^{-N(n+1)k^++2k^+}.
\end{equation}
Moreover,
\begin{equation}\label{wws2}
\sum_{j\geq -k^-}\|U_{j,k}^{\mathcal{L}\psi}(t)\|_{L^\infty}\lesssim \varep_1\langle t\rangle^{-1+H(q+1,n+1)\delta}2^{k^-/2}2^{-N(n+1)k^++2k^+}\min\{1,2^{2k^-}\langle t\rangle\}
\end{equation}
and, if $j\geq -k^-$, 
\begin{equation}\label{wws2*}
\|U_{j,k}^{\mathcal{L}\psi}(t)\|_{L^\infty}\lesssim \varep_1\langle t\rangle^{-3/2+H(q+1,n+1)\delta}2^{j/2}2^{-N(n+1)k^++2k^+}.
\end{equation}
Finally, if $2^{2k^--20}\langle t\rangle\geq 1$ and $\mathcal{L}\in\mathcal{V}_n^q$, $n\leq 1$, then
\begin{equation}\label{wws2.5}
\begin{split}
\sum_{2^j\in [2^{-k^-},2^{k^--20}\langle t\rangle]}&\|U_{j,k}^{\mathcal{L}\psi}(t)\|_{L^\infty}\lesssim \varep_1\langle t\rangle^{-3/2+H(q+1,n+2)\delta}2^{-k^-/2}2^{-N(n+2)k^++5k^+}(\langle t\rangle 2^{2k^-})^{\delta/4}.
\end{split}
\end{equation}

(ii) In the case $n=0$ ($\mathcal{L}=Id$) these bounds can be improved slightly. More precisely, assume $k,J\in\mathbb{Z}$, $t\in[0,T]$, $2^{k^-}\langle t\rangle\geq 2^{20}$, and $2^J\in [2^{-k^-},2^{-10}\langle t\rangle]$. Then for any $H\in\big\{F,\omega_a,\vartheta_{ab}:\,a,b\in\{1,2,3\}\big\}$ and $\al,\be\in\{0,1,2,3\}$,
\begin{equation}\label{wws12x}
\|U_{\leq J,k}^{H}(t)\|_{L^\infty}+\langle t\rangle^{-\delta}\|U_{\leq J,k}^{h}(t)\|_{L^\infty}\lesssim \varep_1\langle t\rangle^{-1}2^{k^--\kappa k^-}2^{-N_0k^++5k^+}.
\end{equation}
Moreover, if $k,J\in\mathbb{Z}$, $t\in[0,T]$, $2^{2k^-}\langle t\rangle\geq 2^{20}$, and $2^J\in [2^{-k^-},2^{k^--20}\langle t\rangle]$, then
\begin{equation}\label{wws13x}
\|U_{\leq J,k}^{\psi}(t)\|_{L^\infty}\lesssim \varep_1\langle t\rangle^{-3/2}2^{-k^-/2+\kappa k^-/20}2^{-N_0k^++5k^+}.
\end{equation}
\end{lemma}

\begin{proof} (i) We prove first \eqref{wws1}. We estimate, using just \eqref{vcx1.1},
\begin{equation}\label{wws5}
\|e^{-it\Lambda_{wa}}V_{j,k}^{\mathcal{L}h}(t)\|_{L^\infty}\lesssim 2^{3k/2}\|V_{j,k}^{\mathcal{L}h}(t)\|_{L^2}\lesssim \varep_1Y(k,t;q,n)2^{k}2^{-j}(2^{k^-}\langle t\rangle)^{-\ga}.
\end{equation}
This suffices to prove \eqref{wws1} if $2^k\lesssim \langle t\rangle^{-1}$, by summing over $j\geq -k$. On the other hand, if $2^k\geq 2^{20}\langle t\rangle^{-1}$ then \eqref{wws5} still suffices to control the sum over $j$ with $2^j\geq 2^{-10}\langle t\rangle$. Finally, if $2^k\geq \langle t\rangle^{-1}$ and $2^j\leq 2^{-10}\langle t\rangle$ then we use \eqref{Linfty1} and \eqref{vcx1.1} to estimate
\begin{equation*}
\|e^{-it\Lambda_{wa}}V_{j,k}^{\mathcal{L}h}(t)\|_{L^\infty}\lesssim 2^{3k/2}\langle t\rangle^{-1}2^j\|Q_{j,k}V^{\mathcal{L}h}(t)\|_{L^2}\lesssim \langle t\rangle^{-1}2^{k}\varep_1Y(k,t;q,n)(\langle t\rangle 2^{k^-})^{-\ga}.
\end{equation*}
The desired conclusion \eqref{wws1} follows by summation over $j$ with $2^j\in [2^{-k^-}, 2^{-10}\langle t\rangle]$.

For \eqref{wws11} we use \eqref{vcx1} and \eqref{Linfty1.6}. Recalling \eqref{vcx1new2} we estimate the left-hand side by
\begin{equation*}
\begin{split}
C&\sum_{2^j\in [2^{-k^-},2^{-20}\langle t\rangle]}2^{k^+}\langle t\rangle^{-1}2^{k^-/2}(\langle t\rangle 2^{k^-})^{\delta/8}\|Q_{j,k}V^{\mathcal{L}h}(t)\|_{H^{0,1}_\Omega}\\
&\lesssim \sum_{2^j\in [2^{-k^-},2^{-20}\langle t\rangle]}2^{k^+}\langle t\rangle^{-1}2^{k^-/2}(\langle t\rangle 2^{k^-})^{\delta/20}\cdot \varep_1\langle t\rangle^{H(q,n+1)\delta}2^{-N(n+1)k^+}2^{k/2}(\langle t\rangle 2^{k^-})^{-\ga}.
\end{split}
\end{equation*}
The desired bound \eqref{wws11} follows. 

We prove now \eqref{wws2}. As in \eqref{wws5} we have
\begin{equation}\label{wws8}
\|e^{-it\Lambda_{kg}}V_{j,k}^{\mathcal{L}\psi}(t)\|_{L^\infty}\lesssim 2^{3k/2}\|V_{j,k}^{\mathcal{L}\psi}(t)\|_{L^2}\lesssim \varep_1Y(k,t;q,n)2^{3k/2}2^{-k^+}2^{-j}.
\end{equation}
This suffices to prove the desired bound when $2^{2k}\leq\langle t\rangle^{-1}$. This bound also suffices to control the sum over $j$ with $2^j\geq \langle t\rangle 2^{k^-}2^{-k^+}$ when $2^{2k}\geq\langle t\rangle^{-1}$. On the other hand, if $2^j\leq \langle t\rangle 2^{k^-}2^{-k^+}$ then we use \eqref{Linfty3} and \eqref{vcx1.1} to estimate
\begin{equation*}
\|e^{-it\Lambda_{kg}}V_{j,k}^{\mathcal{L}\psi}(t)\|_{L^\infty}\lesssim 2^{3k^+}\langle t\rangle^{-3/2}2^{3j/2}\|Q_{j,k}V^{\mathcal{L}\psi}(t)\|_{L^2}\lesssim \varep_1Y(k,t;q,n)2^{2k^+}\langle t\rangle^{-3/2}2^{j/2}.
\end{equation*}
The desired bound \eqref{wws2} follows by summation over $j$ with $2^j\leq \langle t\rangle 2^{k^-}2^{-k^+}$.

The bounds \eqref{wws2*} follow from \eqref{Linfty3} and \eqref{vcx1.1}. The bounds \eqref{wws2.5} follows from \eqref{Linfty3.6} and \eqref{vcx1.1}, once we notice that $\Omega_{ab} Q_{j,k}V^{\mathcal{L}\psi}=Q_{j,k}V^{\Omega_{ab}\mathcal{L}\psi}$, for any rotation vector-field $\Omega_{ab}$.

(ii)  To prove \eqref{wws12x} we define $J_0$ such that $2^{J_0}=2^{10}\langle t\rangle^{1/2}2^{-k/2}$ and estimate
\begin{equation*}
\|e^{-it\Lambda_{wa}}V_{\leq J,k}^{H}(t)\|_{L^\infty}\lesssim 2^{2k}\langle t\rangle^{-1}\|\widehat{Q_{\leq Jk}V^H}\|_{L^\infty}\lesssim 2^{2k}\langle t\rangle^{-1}\varep_12^{-k^--\kappa k^-}2^{-N_0k^+},
\end{equation*}
if $J\leq J_0$, using \eqref{Linfty1.1} and \eqref{vcx1.2}. If $J\geq J_0$ then we estimate the remaining contribution by
\begin{equation*}
\begin{split}
C\sum_{j\in[J_0,J]}\|e^{-it\Lambda_{wa}}V_{j,k}^{H}(t)\|_{L^\infty}&\lesssim \sum_{j\geq J_0}2^{k^+}\langle t\rangle^{-1}2^{k^-/2}(\langle t\rangle 2^{k^-})^{\delta/20}\|Q_{j,k}V^{H}(t)\|_{H^{0,1}_\Omega}\\
&\lesssim 2^{-N(2)k^++3k^+}\langle t\rangle^{-3/2+(H(1,2)+1)\delta}2^{k^-/2},
\end{split}
\end{equation*}
where we used \eqref{Linfty1.6} and \eqref{vcx1.2*}. These two bounds suffice to prove the estimates \eqref{wws12x} for $H$ when $2^k\lesssim \langle t\rangle^{1/(5d)}$; if $2^k\geq \langle t\rangle^{1/(5d)}$ then the desired estimates for $H$ follow by Sobolev embedding from \eqref{vcx1}. The bounds for the metric components $h$ follow by a similar argument.

The bounds \eqref{wws13x} follow in a similar way, using just \eqref{vcx1} if $2^k\gtrsim \langle t\rangle^{1/(10d)}$ and \eqref{vcx1.2} if $2^k\lesssim \langle t\rangle^{-1/2+\kappa/4}$.  On the other hand, if $2^k\in[\langle t\rangle^{-1/2+\kappa/4},\langle t\rangle^{1/(10d)}]$ then we use \eqref{Linfty3.1} and \eqref{vcx1.2} if $2^J\leq 2^{10}\langle t\rangle^{1/2}$, and \eqref{Linfty3.6} and \eqref{vcx1.1} to estimate the remaining contribution if $2^J\geq 2^{10}\langle t\rangle^{1/2}$.
\end{proof}

\begin{remark} We notice that the last two bounds \eqref{wws12x} and \eqref{wws13x} provide sharp pointwise decay at the rate of $\langle t\rangle^{-1}$ and $\langle t\rangle^{-3/2}$ for some parts of the metric tensor and of the Klein-Gordon field. In all the other pointwise bounds in Lemma \ref{dtv3} we allow small $\langle t\rangle^{C\delta}$ losses relative to these sharp decay rates. 
\end{remark}

We prove now several linear bounds on the profiles $V^{\Omega^ah}$ and $V^{\Omega^a\psi}$. These bounds are slight improvements in certain ranges of the bounds one could derive directly from the bootstrap assumptions. These improvements are important in several nonlinear estimates, and are possible because we use interpolation (Lemma \ref{inte1}) to take advantage of the stronger assumptions \eqref{vcx1.2} we have on the functions $V^h$ and $V^\psi$. 

\begin{lemma}\label{dtv6}
(i) For $a\in[0,3]$ let $\Omega^a$ denote vector-fields of the form $\Omega_{23}^{a_1}\Omega_{31}^{a_2}\Omega_{12}^{a_3}$ with $a_1+a_2+a_3\leq a$. If $t\in[0,T]$ and $2^k\gtrsim \langle t\rangle^{-1}$ then
\begin{equation}\label{UOmegaMoreReg}
\begin{split}
\Vert P_kV^{\Omega^1 h}(t)\Vert_{L^2}&\lesssim \varepsilon_12^{k/2}2^{-N(1)k^+}\langle t\rangle^{H(0,1)\delta}\cdot \langle t\rangle^{20\delta}2^{-6dk^+},\\
\Vert P_kV^{\Omega^2 h}(t)\Vert_{L^2}&\lesssim \varepsilon_12^{k/2}2^{-N(2)k^+}\langle t\rangle^{H(0,2)\delta}\cdot \langle t\rangle^{20\delta}2^{-4dk^+},
\end{split}
\end{equation}
where $h\in\big\{h_{\al\be}:\al,\be\in\{0,1,2,3\}\big\}$ as before. Moreover
\begin{equation}\label{UOmegaMoreReg6}
\begin{split}
\Vert P_kV^{\Omega^1 \psi}(t)\Vert_{L^2}&\lesssim \varepsilon_12^{-N(1)k^+}\langle t\rangle^{H(0,1)\delta}\cdot 2^{k^-}\langle t\rangle^{10\delta}2^{-dk^+},\\
\Vert P_kV^{\Omega^2 \psi}(t)\Vert_{L^2}&\lesssim \varepsilon_12^{-N(2)k^+}\langle t\rangle^{H(0,2)\delta}\cdot 2^{k^-/3}\langle t\rangle^{15\delta}2^{-2dk^+},
\end{split}
\end{equation}
and, for any $j\geq -k^-$,
\begin{equation}\label{UOmegaMoreReg5}
\begin{split}
\Vert V^{\Omega^1 h}_{j,k}(t)\Vert_{L^2}&\lesssim \varepsilon_12^{k/2}2^{-N(1)k^+}\langle t\rangle^{H(0,1)\delta}\cdot \langle t\rangle^{50\delta}2^{dk^+}2^{-(2/3)(j+k)},\\
\Vert V^{\Omega^2 h}_{j,k}(t)\Vert_{L^2}&\lesssim \varepsilon_12^{k/2}2^{-N(2)k^+}\langle t\rangle^{H(0,2)\delta}\cdot \langle t\rangle^{25\delta}2^{dk^+}2^{-(1/3)(j+k)}.
\end{split}
\end{equation}

(ii) In addition, for any $J\geq -k^-$,
\begin{equation}\label{UOmegaMoreReg2}
\begin{split}
\Vert U^{\Omega^1 h}_{\leq J,k}(t)\Vert_{L^4}&\lesssim \varepsilon_12^{3k^-/4}2^{-\frac{N(1)+N(2)}{2}k^+}\langle t\rangle^{-1/2+\frac{H(0,1)+H(0,2)}{2}\delta}\cdot \langle t\rangle^{6\delta}2^{2k^+},\\
\Vert U^{\Omega^1 h}_{\leq J,k}(t)\Vert_{L^6}&\lesssim \varepsilon_12^{5k^-/6}2^{-\frac{2N(1)+N(3)}{3}k^+}\langle t\rangle^{-2/3+\frac{2H(0,1)+H(0,3)}{3}\delta}\cdot \langle t\rangle^{8\delta}2^{2k^+},\\
\Vert U^{\Omega^2 h}_{\leq J,k}(t)\Vert_{L^3}&\lesssim \varepsilon_12^{2k^-/3}2^{-\frac{N(1)+2N(3)}{3}k^+}\langle t\rangle^{-1/3+\frac{H(0,1)+2H(0,3)}{3}\delta}\cdot \langle t\rangle^{4\delta}2^{2k^+}.
\end{split}
\end{equation}

Moreover, if $2^{2k^-}\langle t\rangle\geq 2^{20}$ and $2^J\in [2^{-k^-},2^{k^--20}\langle t\rangle]$ then
\begin{equation}\label{UOmegaMoreReg3}
\begin{split}
\Vert U^{\Omega^1 \psi}_{\leq J,k}(t)\Vert_{L^4}&\lesssim \varepsilon_12^{-\frac{N(1)+N(2)}{2}k^+}\langle t\rangle^{-3/4+\delta H(0,2)/2}2^{-k^-/4}\cdot 2^{\kappa k^-/60}2^{-2k^+},\\
\Vert U^{\Omega^1 \psi}_{\leq J,k}(t)\Vert_{L^6}&\lesssim \varepsilon_12^{-\frac{2N(1)+N(3)}{3}k^+}\langle t\rangle^{-1+\delta H(0,3)/3}2^{-k^-/3}\cdot 2^{\kappa k^-/60}2^{-2.5k^+},\\
\Vert U^{\Omega^2 \psi}_{\leq J,k}(t)\Vert_{L^3}&\lesssim \varepsilon_12^{-\frac{N(1)+2N(3)}{3}k^+}\langle t\rangle^{-1/2+2\delta H(0,3)/3}2^{-k^-/6}\cdot 2^{\kappa k^-/60}2^{-1.5k^+}.
\end{split}
\end{equation}
\end{lemma} 

\begin{proof} 

(i) We use the general interpolation inequalities in \eqref{vcx22}--\eqref{vcx22.5}. 
The bounds \eqref{UOmegaMoreReg} (which are relevant only when $2^{dk+}\geq 2^{10}\langle t\rangle^\delta$) follow from \eqref{vcx1}. 

The bounds in the first line of \eqref{UOmegaMoreReg6} follow directly from \eqref{cnb2} if $k\leq 0$ (notice that $H(1,1)-H(0,1)=10$) and from \eqref{vcx1} and \eqref{vcx22} if $k\geq 0$. Similarly, the bounds in the second line of \eqref{UOmegaMoreReg6} follow from \eqref{vcx22.5}, \eqref{vcx1.2} (if $k\leq 0$) and \eqref{vcx1} (if $k\geq 0$).

To prove the bounds in the first line of \eqref{UOmegaMoreReg5} we use \eqref{vcx22.5}, \eqref{vcx1}, and \eqref{vcx1.1}, 
\begin{equation*}
\begin{split}
\Vert V^{\Omega^1 h}_{j,k}(t)\Vert_{L^2}&\lesssim \Vert V^{h}_{j,k}(t)\Vert_{L^2}^{2/3}\Vert V^{h}_{j,k}(t)\Vert_{H^{0,3}_\Omega}^{1/3}\\
&\lesssim  [\varep_12^{k/2}\langle t\rangle^{H(1,1)}2^{-N(1)k^+}2^{-(j+k)}]^{2/3}[\varep_12^{k/2}\langle t\rangle^{H(0,3)}2^{-N(3)k^+}]^{1/3},
\end{split}
\end{equation*}
which gives the desired bounds. The estimates in the second line follow in a similar way.

(ii) To prove \eqref{UOmegaMoreReg2} we use \eqref{inte4}--\eqref{inte5}. Indeed, using also \eqref{wws1} and \eqref{vcx1},
\begin{equation*}
\begin{split}
\Vert U^{\Omega^1 h}_{\leq J,k}(t)\Vert_{L^4}&\lesssim \Vert U^{h}_{\leq J,k}(t)\Vert_{L^\infty}^{1/2}\Vert U^{h}_{\leq J,k}(t)\Vert_{H^{0,2}_\Omega}^{1/2}\\
&\lesssim \varep_1\big(\langle t\rangle^{-1+H(1,1)\delta}2^{k^-}2^{-N(1)k^++2k^+}\big)^{1/2}\big(\langle t\rangle^{H(0,2)\delta}2^{k/2}2^{-N(2)k^+}\big)^{1/2},
\end{split}
\end{equation*}
which gives the bounds in the first line of \eqref{UOmegaMoreReg2} (recall that $H(1,1)=H(0,1)+10$). Similarly,
\begin{equation*}
\begin{split}
\Vert U^{\Omega^1 h}_{\leq J,k}(t)\Vert_{L^6}&\lesssim \Vert U^{h}_{\leq J,k}(t)\Vert_{L^\infty}^{2/3}\Vert U^{h}_{\leq J,k}(t)\Vert_{H^{0,3}_\Omega}^{1/3}\\
&\lesssim \varep_1\big(\langle t\rangle^{-1+H(1,1)\delta}2^{k^-}2^{-N(1)k^++2k^+}\big)^{2/3}\big(\langle t\rangle^{H(0,3)\delta}2^{k/2}2^{-N(3)k^+}\big)^{1/3}
\end{split}
\end{equation*}
and
\begin{equation*}
\begin{split}
\Vert U^{\Omega^2 h}_{\leq J,k}(t)\Vert_{L^3}&\lesssim \Vert U^{h}_{\leq J,k}(t)\Vert_{L^\infty}^{1/3}\Vert U^{h}_{\leq J,k}(t)\Vert_{H^{0,3}_\Omega}^{2/3}\\
&\lesssim \varep_1\big(\langle t\rangle^{-1+H(1,1)\delta}2^{k^-}2^{-N(1)k^++2k^+}\big)^{1/3}\big(\langle t\rangle^{H(0,3)\delta}2^{k/2}2^{-N(3)k^+}\big)^{2/3}.
\end{split}
\end{equation*}
The remaining bounds in \eqref{UOmegaMoreReg2} follow.{\footnote{We notice that these bounds can be improved slightly if $2^J\leq \langle t\rangle 2^{-20}$, by using the $L^\infty$ bounds \eqref{wws11} instead of \eqref{wws1}. This is not useful for us, however, since we will mostly apply the bounds \eqref{UOmegaMoreReg2} to estimate the contributions of localized profiles corresponding to large $j$.}}

The bounds \eqref{UOmegaMoreReg3} follow in a similar way, using \eqref{inte4}--\eqref{inte5}, \eqref{wws13x} and \eqref{vcx1}:
\begin{equation*}
\begin{split}
\Vert U^{\Omega^1 \psi}_{\leq J,k}(t)\Vert_{L^4}\lesssim \varep_1\big(\langle t\rangle^{-3/2}2^{-k^-/2+\kappa k^-/20}2^{-N(1)k^+-5k^+}\big)^{1/2}\big(\langle t\rangle^{H(0,2)\delta}2^{-N(2)k^+}\big)^{1/2},
\end{split}
\end{equation*}
\begin{equation*}
\begin{split}
\Vert U^{\Omega^1 \psi}_{\leq J,k}(t)\Vert_{L^6}\lesssim \varep_1\big(\langle t\rangle^{-3/2}2^{-k^-/2+\kappa k^-/20}2^{-N(1)k^+-5k^+}\big)^{2/3}\big(\langle t\rangle^{H(0,3)\delta}2^{-N(3)k^+}\big)^{1/3},
\end{split}
\end{equation*}
\begin{equation*}
\begin{split}
\Vert U^{\Omega^2 \psi}_{\leq J,k}(t)\Vert_{L^3}\lesssim \varep_1\big(\langle t\rangle^{-3/2}2^{-k^-/2+\kappa k^-/20}2^{-N(1)k^+-5k^+}\big)^{1/3}\big(\langle t\rangle^{H(0,3)\delta}2^{-N(3)k^+}\big)^{2/3}.
\end{split}
\end{equation*}
The bounds in \eqref{UOmegaMoreReg3} follow.
\end{proof}

We also record a few additional $L^\infty$ bounds in the Fourier space.

\begin{lemma}\label{box100.5}
 If $k\in\mathbb{Z}$ and $\mathcal{L}\in\mathcal{V}_n^q$, $n\leq 1$, then
\begin{equation}\label{vcx1.3*}
\|\widehat{P_kV^{\mathcal{L} h_{\al\be}}}(t)\|_{L^\infty}\lesssim \varep_12^{-k^--\delta k^-/2}\langle t\rangle^{\frac{H(q,n+1)+H(q+1,n+2)}{2}\delta}2^{-N_0k^++(n+3/2)dk^+-3k^+/4}.
\end{equation}
Moreover, if $(k,j)\in\mathcal{J}$ and $\mathcal{L}\in\mathcal{V}_n^q$, $n\leq 1$, then
\begin{equation}\label{vcx1.15}
2^{3k/2}\big\|\widehat{Q_{j,k}V^{\mathcal{L}h_{\al\be}}}(t)\big\|_{L^\infty}+2^{k}\big\|\widehat{Q_{j,k}V^{\mathcal{L}\psi}}(t)\big\|_{L^\infty}\lesssim\varep_12^{-j/2+\delta j/4}Y(k,t;q,n+1)2^{\delta k^+/4}.
\end{equation}
\end{lemma}

\begin{proof} The bounds \eqref{vcx1.15} follow from \eqref{Linfty3.3}, \eqref{vcx1.1}, and \eqref{vcx1new2}. To prove the bounds \eqref{vcx1.3*} we use the estimates
\begin{equation*}
\begin{split}
\|Q_{j,k}V^{\mathcal{L}h_{\al\be}}(t)\|_{H^{0,1}_\Omega}&\lesssim \varep_1\langle t\rangle^{H(q,n+1)\delta}2^{k/2}2^{-N(n+1)k^+}(2^{k^-}\langle t\rangle)^{-\ga},\\
2^{j+k}\|Q_{j,k}V^{\mathcal{L}h_{\al\be}}(t)\|_{H^{0,1}_\Omega}&\lesssim \varep_1\langle t\rangle^{H(q+1,n+2)\delta}2^{k/2}2^{-N(n+2)k^+}(2^{k^-}\langle t\rangle)^{-\ga},
\end{split}
\end{equation*}
which follow from \eqref{vcx1}--\eqref{vcx1.1}, and the bounds \eqref{consu2}. 
\end{proof}

\chapter{The nonlinearities $\mathcal{N}^h_{\al\be}$ and $\mathcal{N}^\psi$}

\section{Localized $L^2$ bilinear estimates} \label{LocInter}

\subsection{The main nonlinearities} In this section we start our detailed analysis of the nonlinearities $\mathcal{N}_{\al\be}^{h}$ and $\mathcal{N}^\psi$. We start from the identities \eqref{zaq11.1}, which we write in the form
\begin{equation*}
\mathcal{N}_{\al\be}^{h}=g_{\geq 1}^{00}\partial_0^2h_{\al\be}+\sum_{(\mu,\nu)\neq (0,0)}g_{\geq 1}^{\mu\nu}\partial_\mu\partial_\nu h_{\al\be}+\mathcal{KG}_{\al\be}-F^{\geq 2}_{\al\be}(g,\partial g).
\end{equation*}
We would like to eliminate the terms that contain two time derivatives in both the metric nonlinearities, in order to express all the terms elliptically in terms of the normalized solutions $U^{h_{\al\be}}$ and $U^\psi$. Indeed, since $\partial_0^2h_{\al\be}=\Delta h_{\al\be}+\mathcal{N}_{\al\be}^{h}$, we have
\begin{equation}\label{sac1}
\mathcal{N}_{\al\be}^{h}=(1-g_{\geq 1}^{00})^{-1}\Big[\sum_{(\mu,\nu)\neq (0,0)}g_{\geq 1}^{\mu\nu}\partial_\mu\partial_\nu h_{\al\be}+g_{\geq 1}^{00}\Delta h_{\al\be}+\mathcal{KG}_{\al\be}-F^{\geq 2}_{\al\be}(g,\partial g)\Big].
\end{equation}

Recall the decomposition of the metric components,
\begin{equation}\label{box100}
g^{\al\be}=m^{\al\be}+g_{\geq 1}^{\al\be}=m^{\al\be}+g_{1}^{\al\be}+g_{\geq 2}^{\al\be},
\end{equation}
into the Minkowski metric, a linearized metric, and a quadratic metric. Using this decomposition we extract the quadratic components of the nonlinearity 
\begin{equation}\label{sac1.1}
\mathcal{N}_{\al\be}^{h,2}:=\mathcal{KG}^2_{\al\be}+\mathcal{Q}^2_{\al\be}+\mathcal{S}^2_{\al\be},
\end{equation}
where $\mathcal{KG}^2_{\al\be}$ are semilinear quadratic terms that involve the Klein--Gordon field,
\begin{equation}\label{sac1.5}
\mathcal{KG}^2_{\al\be}:=2\partial_\al\psi\partial_\be\psi+\psi^2m_{\al\be},
\end{equation}
$\mathcal{Q}^2_{\al\be}$ are quasilinear quadratic terms,
\begin{equation}\label{sac1.2}
\mathcal{Q}^2_{\al\be}:=\sum_{(\mu,\nu)\neq (0,0)}g_{1}^{\mu\nu}\partial_\mu\partial_\nu h_{\al\be}+g_{1}^{00}\Delta h_{\al\be}=-h_{00}\Delta h_{\al\be}+2h_{0j}\partial_0\partial_jh_{\al\be}-h_{jk}\partial_j\partial_kh_{\al\be},
\end{equation}
compare with the formulas \eqref{sac1}, and 
\begin{equation*}
\mathcal{S}^2_{\al\be}:=-(Q^2_{\al\be}+P^2_{\al\be}),
\end{equation*}
are semilinear quadratic terms that involve the metric components,  where
\begin{equation}\label{zaq20}
\begin{split}
Q^2_{\al\be}&:=m^{\rho\rho'}m^{\lambda\la'}(\partial_\al h_{\rho'\la'}\partial_{\rho} h_{\be\la}-\partial_{\rho} h_{\rho'\lambda'}\partial_\al h_{\be\la})+m^{\rho\rho'}m^{\lambda\la'}(\partial_\be h_{\rho'\la'}\partial_{\rho} h_{\al\la}-\partial_{\rho} h_{\rho'\lambda'}\partial_\be h_{\al\la})\\
&+\frac{1}{2}m^{\rho\rho'}m^{\lambda\la'}(\partial_{\la'} h_{\rho\rho'}\partial_\be h_{\al\la}-\partial_\be h_{\rho\rho'}\partial_{\la'} h_{\al\la})+\frac{1}{2}m^{\rho\rho'}m^{\lambda\la'}(\partial_{\la'} h_{\rho\rho'}\partial_\al h_{\be\la}-\partial_\al h_{\rho\rho'}\partial_{\la'} h_{\be\la})\\
&-m^{\rho\rho'}m^{\lambda\la'}(\partial_{\la} h_{\al\rho'}\partial_\rho h_{\be\la'}-\partial_\rho h_{\al\rho'}\partial_{\la} h_{\be\la'})+m^{\rho\rho'}m^{\lambda\la'}\partial_{\rho'} h_{\al\la'}\partial_\rho h_{\be\la}
\end{split}
\end{equation}
and
\begin{equation}\label{zaq21}
P^2_{\al\be}:=-\frac{1}{2}m^{\rho\rho'}m^{\lambda\lambda'}\partial_\al h_{\rho'\la'}\partial_\be h_{\rho\la}+\frac{1}{4}m^{\rho\rho'}m^{\lambda\la'}\partial_\al h_{\rho\rho'}\partial_\be h_{\la\la'}.
\end{equation}
Compare with the formulas \eqref{tr11}--\eqref{tr12}. Let $\mathcal{N}_{\al\be}^{h,\geq 3}:=\mathcal{N}_{\al\be}^{h}-\mathcal{N}_{\al\be}^{h,2}$ denote the cubic and higher order components of $\mathcal{N}_{\al\be}^{h}$.

Similarly, recall also the Klein--Gordon nonlinearities $\mathcal{N}^\psi$ defined in \eqref{zaq11.3}. We write
\begin{equation*}
\mathcal{N}^\psi=\sum_{(\mu,\nu)\neq (0,0)}g_{\geq 1}^{\mu\nu}\partial_\mu\partial_\nu\psi+g_{\geq 1}^{00}\partial_0^2\psi=\sum_{(\mu,\nu)\neq (0,0)}g_{\geq 1}^{\mu\nu}\partial_\mu\partial_\nu\psi+g_{\geq 1}^{00}(\Delta\psi-\psi+\mathcal{N}^\psi).
\end{equation*}
Therefore
\begin{equation}\label{gb10}
\mathcal{N}^\psi=(1-g_{\geq 1}^{00})^{-1}\Big[\sum_{(\mu,\nu)\neq (0,0)}g_{\geq 1}^{\mu\nu}\partial_\mu\partial_\nu\psi+g_{\geq 1}^{00}(\Delta\psi-\psi)\Big].
\end{equation}
Using also the identities \eqref{zaq22} we extract the quadratic component of $\mathcal{N}^\psi$,
\begin{equation}\label{npsi2}
\mathcal{N}^{\psi,2}:=\sum_{(\mu,\nu)\neq (0,0)}g_{1}^{\mu\nu}\partial_\mu\partial_\nu\psi+g_{1}^{00}(\Delta\psi-\psi)=-h_{00}(\Delta\psi-\psi)+2h_{0j}\partial_0\partial_j\psi-h_{jk}\partial_j\partial_k\psi,
\end{equation}
and let $\mathcal{N}^{\psi,\geq 3}:=\mathcal{N}^\psi-\mathcal{N}^{\psi,2}$ denote the cubic and higher order component of $\mathcal{N}^\psi$.

One of our main goals is to prove good bounds on the various components of the nonlinearities $\mathcal{N}^{h}_{\al\be}$ and $\mathcal{N}^\psi$. More precisely, we are looking to prove frequency localized $L^2$ bounds on $\mathcal{L}\mathcal{N}^{h}_{\al\be}$ and $\mathcal{L}\mathcal{N}^{\psi}$. 

Ideally, we would like to prove that these nonlinearities satisfy bounds of the form
\begin{equation}\label{cnb10}
\begin{split}
\sum_{\al,\be\in\{0,1,2,3\}}\|P_k(\mathcal{L}\mathcal{N}_{\al\be}^{h})(t)\|_{L^2}&\lesssim \varep_1\min(2^k,\langle t\rangle^{-1}) \sum_{\al,\be\in\{0,1,2,3\}}\|P_k(\nabla\mathcal{L}h_{\al\be})\|_{L^2},\\
\|P_k(\mathcal{L}\mathcal{N}^{\psi})(t)\|_{L^2}&\lesssim \varep_1\min(2^k,\langle t\rangle^{-1}) \|P_k(\mathcal{L}\psi)\|_{L^2}.
\end{split}
\end{equation}
Unfortunately, such optimal bounds do not hold for most of the important components of the nonlinearities. As we will see, we have both derivative loss, due to the quasilinear nature of the system, and loss of decay in time, due to the slower decay of the metric components $h_{\al\be}$. However, we can still prove estimates that are somewhat close to \eqref{cnb10}, but with certain losses. To quantify this, we define the acceptable loss function
\begin{equation}\label{AcceptableLoss}
\ell(0,0):=3,\qquad\ell(0,1):=13,\qquad\ell(1,1):=23,\qquad\ell(q,n):=33\,\,\text{ if }\,\,n\ge2.
\end{equation}
Notice that $\ell(q,n)+7 \le H(q,n)$ if $n\ge1$. 

\subsection{Multipliers and bilinear operators} We define two classes of multipliers $\mathcal{M}_0$ and $\mathcal{M}$ by
\begin{equation}\label{mults0}
\mathcal{M}_0:=\{m:\mathbb{R}^3\setminus\{0\}\to\mathbb{C}:\, |x|^{|\alpha|}\,|D^\al_x m(x)|\lesssim_{|\alpha|}1\text{ for any }\al\in\mathbb{Z}_+^3\text{ and }x\in \mathbb{R}^3\setminus\{0\}\},
\end{equation}
and
\begin{equation}\label{mults}
\begin{split}
\mathcal{M}:=\{&m:\mathbb{R}^6\to\mathbb{C}:\, m(x,y)=m_1(x,y)m'(x+y),\,\,m'\in\mathcal{M}_0,\\
&|x|^{|\alpha|}|y|^{\be}\,|D^\al_xD^\be_y m_1(x,y)|\lesssim_{|\alpha|,|\be|}1\text{ for any }\al,\be\in\mathbb{Z}_+^3\text{ and }x,y\in \mathbb{R}^3\setminus\{0\}\}.
\end{split}
\end{equation}
In most of our applications the multipliers in $\mathcal{M}$ will be of the form $m_1(x)m_2(y)$, where $m_1,m_2\in\mathcal{M}_0$. We will also need to allow sums of such multipliers in order to be able to define the important classes of null multipliers $\mathcal{M}^{null}_{\pm}\subseteq\mathcal{M}$, see Definition \ref{nullStr1}.

In some of our constructions, in particular connected to the application of normal forms and associated angular cutoffs, the class of multipliers $\mathcal{M}$ is too restrictive. To treat such situations we define a more general class of multipliers
\begin{equation}\label{multsstar}
\mathcal{M}^\ast:=\{m\in L^\infty(\mathbb{R}^6):\,\|\mathcal{F}^{-1}\{m\cdot\varphi_{k_1}(x)\varphi_{k_2}(y)\varphi_{k}(x+y)\}\|_{L^1(\mathbb{R}^6)}\lesssim 1\text{ for any }k_1,k_2,k\in\mathbb{Z}\}.
\end{equation}

Given a bounded multiplier $m$ let $I=I_m$ denote the bilinear operator
\begin{equation}\label{abc36.1}
\widehat{I[f,g]}(\xi)=\widehat{I_m[f,g]}(\xi):=\frac{1}{8\pi^3}\int_{\mathbb{R}^3}m(\xi-\eta,\eta)\widehat{f}(\xi-\eta)\widehat{g}(\eta)\,d\eta.
\end{equation}
We will often use the simple $L^2$ bounds
\begin{equation}\label{abc36.6}
\|P_kI_m[P_{k_1}f,P_{k_2}g]\|_{L^2}\lesssim 2^{3\min\{k,k_1,k_2\}/2}\|P_{k_1}f\|_{L^2}\|P_{k_2}g\|_{L^2}
\end{equation}
for any multiplier $m$ satisfying $\|m\|_{L^\infty}\leq 1$, $f,g\in L^2(\mathbb{R}^3)$, and $k,k_1,k_2\in\mathbb{Z}$.

\subsection{Frequency localized bilinear estimates} In this subsection we prove several bounds on localized bilinear interactions, which are the main building blocks for the estimates on the nonlinearities $\mathcal{N}_{\al\be}^{h}$ and $\mathcal{N}^\psi$ in the next section. Notice that, as a consequence of the definitions \eqref{fvc1.0}, we have the superlinear inequalities
\begin{equation}\label{SuperlinearH2}
\begin{split}
H(q_1,n_1+1)+H(q_2,n_2)&\le H(q_1+q_2,n_1+n_2)+20,\\
H(q_1,n_1+1)+H(q_2,n_2)&\le H(q_1+q_2,n_1+n_2)-40\qquad\hbox{ if }q_2\geq 1,\\
H(q_1+1,n_1+1)+H(q_2,n_2)&\le H(q_1+q_2,n_1+n_2)+30\qquad\hbox{ if }q_2\geq 1,
\end{split}
\end{equation}
which hold when $n_1\geq 1$, $n_2\geq 1$, $n_1+n_2\leq 3$. 

We start by proving $L^2$ bounds on localized bilinear interactions of the metric components. 

\begin{lemma}\label{box1}
Assume that $\mathcal{L}_1\in\mathcal{V}_{n_1}^{q_1}$, $\mathcal{L}_2\in\mathcal{V}_{n_2}^{q_2}$, $n_1+n_2\leq 3$, $h_1,h_2\in\{h_{\al\be}:\,\al,\be\in\{0,1,2,3\}\}$, and $\iota_1,\iota_2\in\{+,-\}$. Assume that $m\in\mathcal{M}$ (see \eqref{mults}), $I=I_m$ is defined as in \eqref{abc36.1}, and let 
\begin{equation}\label{box1.5}
I_{k,k_1,k_2}^{wa,1}(t):=2^{-k/2}\big\Vert P_kI[P_{k_1}U^{\mathcal{L}_1h_1,\iota_1},P_{k_2}U^{\mathcal{L}_2h_2,\iota_2}](t)\big\Vert_{L^2},
\end{equation}
for any $t\in[0,T]$ and $k,k_1,k_2\in\mathbb{Z}$. Then
\begin{equation}\label{box2}
\begin{split}
I_{k,k_1,k_2}^{wa,1}(t)\lesssim \varep_1^22^{-k/2}2^{3\min\{k,k_1,k_2\}/2}&\big(\langle t\rangle^22^{k_1^-+k_2^-}\big)^{-\ga}2^{k_1/2+k_2/2}\\
&\times\langle t\rangle^{[H(q_1,n_1)+H(q_2,n_2)]\delta}2^{-N(n_1)k_1^+-N(n_2)k_2^+}.
\end{split}
\end{equation}

In addition, assuming $n_1\leq n_2$ (in particular $n_1\leq 1$), we have:

(1) if $k=\min\{k,k_1,k_2\}$ and $n_1=1$ then
\begin{equation}\label{box3}
I_{k,k_1,k_2}^{wa,1}(t)\lesssim \varep_1^2\langle t\rangle^{-1+\delta[H(q_2,n_2)+H(q_1,n_1+1)+1]}2^{-2k_2^+}2^{k_2^-/4}2^{-N(n_2)k^+};
\end{equation}

(2) if $k=\min\{k,k_1,k_2\}$ and $n_1=0$ then
\begin{equation}\label{box3.5}
I_{k,k_1,k_2}^{wa,1}(t)\lesssim \varep_1^2\langle t\rangle^{-1+\delta[H(q_2,n_2)+\ell(q_2,n_2)]}2^{-2k_2^+}2^{k_2^-/4}2^{-N(n_2)k^+};
\end{equation}

(3) if $k_1=\min\{k,k_1,k_2\}$ and $n_1\in\{0,1\}$ then
\begin{equation}\label{box4}
I_{k,k_1,k_2}^{wa,1}(t)\lesssim \varep_1^2\langle t\rangle^{-1+ \delta[H(q_2,n_2)+H(q_1,n_1+1)+1]}2^{k_1}2^{-2k_1^+}2^{|k|/4}2^{-N(n_2)k^+};
\end{equation}

(4) if $k_2=\min\{k,k_1,k_2\}$ and $n_1=1$ then
\begin{equation}\label{box5}
I_{k,k_1,k_2}^{wa,1}(t)\lesssim \varep_1^2\langle t\rangle^{-1+ \delta[H(q_2,n_2)+H(q_1,n_1+1)+1]}2^{k_2}2^{-2k_2^+}2^{|k|/4}2^{-N(n_1+1)k^+};
\end{equation}

(5) if $k_2=\min\{k,k_1,k_2\}$ and $(n_1,n_2)\in\{(0,0),(0,1)\}$ then
\begin{equation}\label{box5.5}
I_{k,k_1,k_2}^{wa,1}(t)\lesssim  \varep_1^2\langle t\rangle^{-1+\delta[H(q_2,n_2)+1]}2^{k_2}2^{-2k_2^+}2^{|k|/4}2^{-N_0k^+};
\end{equation}

(6) if $k_2=\min\{k,k_1,k_2\}$ and $(n_1,n_2)\in\{(0,2),(0,3)\}$ then
\begin{equation}\label{box5.6}
I_{k,k_1,k_2}^{wa,1}(t)\lesssim  \varep_1^2\langle t\rangle^{-1+\delta[H(q_2,n_2)+\ell(q_2,n_2)]}2^{k_2}2^{-2k_2^+}2^{|k|/4}2^{-N(1)k^+}.
\end{equation}
\end{lemma}

\begin{proof} We remark first that bilinear $\rm{Wave}\times\mathrm{Wave}$ interactions appear in the metric nonlinearities $\mathcal{N}^{h_{\al\be}}$, both in semilinear and in quasilinear form. According to the general philosophy described in \eqref{cnb10}, we would like to have bounds on the form
\begin{equation}\label{box5.7}
I_{k,k_1,k_2}^{wa,1}(t)\lesssim \varep_1^22^{-|k_1-k_2|}\min(2^{k},\langle t\rangle^{-1})2^{-N(n_1+n_2)k^++k^+}.
\end{equation}
The factors $2^{-|k_1-k_2|}$ in the right-hand side are critical, in order to be able to estimate the quasilinear components of the nonlinearities $\mathcal{N}^{h}_{\al\be}$. We notice that the bounds \eqref{box3}--\eqref{box5.6} that we actually prove are variations of the ideal bounds \eqref{box5.7}, with small $\langle t\rangle^{C\delta}$ loss of decay and loss of derivative $2^{k^+}$ in some cases. For later use, it is very important to minimize the time decay loss as much as possible. 

The estimates follow from a case by case analysis, using Lemmas \ref{dtv2} and \ref{dtv3}, and the bilinear estimates in Lemmas \ref{Lembil1} and \ref{Lembil2}. We estimate first, using just \eqref{abc36.6}
\begin{equation}\label{GeneralQuadraticEst1}
I_{k,k_1,k_2}^{wa,1}(t)\lesssim  2^{-k/2}2^{3\min\{k,k_1,k_2\}/2}\|P_{k_1}U^{\mathcal{L}_1h_1,\iota_1}(t)\|_{L^2}\|P_{k_2}U^{\mathcal{L}_2h_2,\iota_2}(t)\|_{L^2},
\end{equation}
which gives \eqref{box2} in view of \eqref{vcx1}. To prove the rest of the bounds, we consider three cases.

{\bf{Step 1.}} We prove first \eqref{box3} and \eqref{box3.5}. Assume that $k=\min\{k,k_1,k_2\}$. We may also assume that $|k_1-k_2|\leq 4$, $2^k\gtrsim \langle t\rangle^{-1}$,  and $2^{k_2}\lesssim\langle t\rangle^{1/20}$ (otherwise the bounds follow from \eqref{box2}). Let $J_1$ be the largest integer such that $2^{J_1}\le \langle t\rangle(1+2^{k_1}\langle t\rangle)^{-\delta/20}$ and decompose
\begin{equation}\label{box15}
\begin{split}
P_{k_1}U^{\mathcal{L}_1h_1,\iota_1}(t)&=U_{\leq J_1,k_1}^{\mathcal{L}_1h_1,\iota_1}(t)+U_{>J_1,k_1}^{\mathcal{L}_1h_1,\iota_1}(t)\\
&=e^{-it\Lambda_{wa,\iota_1}}V_{\leq J_1,k_1}^{\mathcal{L}_1h_1,\iota_1}(t)+e^{-it\Lambda_{wa,\iota_1}}V_{>J_1,k_1}^{\mathcal{L}_1h_1,\iota_1}(t),
\end{split}
\end{equation}
 see \eqref{on11.3}--\eqref{on11.36}. Using \eqref{Bil31}, \eqref{vcx1}, and \eqref{vcx1.1} we estimate
\begin{equation}\label{box15.1}
\begin{split}
2^{-k/2}&\big\Vert P_kI[U_{\leq J_1,k_1}^{\mathcal{L}_1h_1,\iota_1}(t), P_{k_2}U^{\mathcal{L}_2h_2,\iota_2}(t)]\big\Vert_{L^2}\\
&\lesssim \langle t\rangle^{-1}(1+2^{k_1}\langle t\rangle)^{\delta/10}\big\|Q_{\leq J_1,k_1}V^{\mathcal{L}_1h_1}(t)\big\|_{H^{0,1}_\Omega}\|P_{k_2}U^{\mathcal{L}_2h_2}(t)\|_{L^2}\\
&\lesssim \varep_1^2\langle t\rangle^{-1+\delta[H(q_2,n_2)+H(q_1,n_1+1)+1]}2^{-N(n_2)k_2^+-5k_2^+}2^{k_2^-/2}
\end{split}
\end{equation}
and, using \eqref{abc36.6},
\begin{equation}\label{box15.2}
\begin{split}
2^{-k/2}&\big\Vert P_kI[U_{> J_1,k_1}^{\mathcal{L}_1h_1,\iota_1}(t), P_{k_2}U^{\mathcal{L}_2h_2,\iota_2}(t)]\big\Vert_{L^2}\lesssim  2^{k}\big\|U_{>J_1,k_1}^{\mathcal{L}_1h_1}(t)\big\|_{L^2}\|P_{k_2}U^{\mathcal{L}_2h_2}(t)\|_{L^2}\\
&\lesssim \varep_1^2\langle t\rangle^{-1+\delta[H(q_2,n_2)+H(q_1+1,n_1+1)+1]}2^{-N(n_2)k_2^+-5k_2^+}2^{k^-/2}.
\end{split}
\end{equation}
Moreover, if $n_2\leq 2$ then we can use \eqref{vcx1.1} and \eqref{wws1} to estimate
\begin{equation}\label{box15.3}
\begin{split}
2^{-k/2}&\big\Vert P_kI[U_{> J_1,k_1}^{\mathcal{L}_1h_1,\iota_1}(t), P_{k_2}U^{\mathcal{L}_2h_2,\iota_2}(t)]\big\Vert_{L^2}\lesssim  2^{-k/2}\big\|U_{>J_1,k_1}^{\mathcal{L}_1h_1}(t)\big\|_{L^2}\|P_{k_2}U^{\mathcal{L}_2h_2}(t)\|_{L^\infty}\\
&\lesssim \varep_1^22^{-k/2}\langle t\rangle^{-2+\delta'}2^{-N(n_1+1)k_1^+-N(n_2+1)k_2^++5k_2^+}.
\end{split}
\end{equation}

Since $H(1,1)=30$, the bounds \eqref{box3.5} follow from \eqref{box15.1} and \eqref{box15.2} if $n_1=0$ and $n_2\geq 2$. The bounds \eqref{box3} follow from \eqref{box15.1} and \eqref{box15.3} if $n_1\geq 1$ (in this case $n_2\leq 2$). 

It remains to prove the bounds \eqref{box3.5} when $n_1=0$ and $n_2\leq 1$. The estimates \eqref{box15.3} and \eqref{box15.1} still suffice if $2^{k_2}\leq \langle t\rangle^{-\delta'}$, but are slightly too weak when $2^{k_2}\geq \langle t\rangle^{-\delta'}$. In this case we need a different decomposition: let $J'_1$ be the largest integer such that $2^{J'_1}\le \langle t\rangle^{1/2}2^{-k_1/2}+2^{-k}$ and decompose $P_{k_1}U^{h_1,\iota_1}(t)=U_{\leq J'_1,k_1}^{h_1,\iota_1}(t)+U_{>J'_1,k_1}^{h_1,\iota_1}(t)$ as in \eqref{box15}. Using \eqref{bil1}, \eqref{vcx1}, and \eqref{vcx1.2} we estimate
\begin{equation}\label{box15.4}
\begin{split}
2^{-k/2}&\big\Vert P_kI[U_{\leq J'_1,k_1}^{h_1,\iota_1}(t), P_{k_2}U^{\mathcal{L}_2h_2,\iota_2}(t)]\big\Vert_{L^2}\\
&\lesssim \langle t\rangle^{-1}2^{3k_1/2}\big\|\mathcal{F}\{Q_{\leq J_1,k_1}V^{h_1}\}(t)\big\|_{L^\infty}\|P_{k_2}U^{\mathcal{L}_2h_2}(t)\|_{L^2}\\
&\lesssim \varep_1^2\langle t\rangle^{-1+\delta[H(q_2,n_2)+1]}2^{-N(n_2)k_2^+-5k_2^+}2^{k_2^-/2}.
\end{split}
\end{equation}
Moreover, using \eqref{vcx1.1} and \eqref{wws1} we estimate
\begin{equation}\label{box15.5}
\begin{split}
2^{-k/2}&\big\Vert P_kI[U_{> J'_1,k_1}^{h_1,\iota_1}(t), P_{k_2}U^{\mathcal{L}_2h_2,\iota_2}(t)]\big\Vert_{L^2}\lesssim  2^{-k/2}\big\|U_{>J'_1,k_1}^{h_1}(t)\big\|_{L^2}\|P_{k_2}U^{\mathcal{L}_2h_2}(t)\|_{L^\infty}\\
&\lesssim \varep_1^22^{-k/2}\langle t\rangle^{-1+\delta'}2^{-N(n_2+1)k_2^++5k_2^+}2^{k_2^-/2}\min\{2^{-J'_1}2^{-N(1)k_1^+},2^{-N(0)k_1^+}\}.
\end{split}
\end{equation}
By analyzing the cases $2^k\geq \langle t\rangle^{\delta'}$, $2^k\in[\langle t\rangle^{-1/2},\langle t\rangle^{\delta'}]$, and $2^k\leq\langle t\rangle^{-1/2}$ it is easy to see that the right-hand side of \eqref{box15.5} is suitably bounded, as claimed in the right-hand side of \eqref{box3.5}. The desired conclusion  follows using also \eqref{box15.4}.

{\bf{Step 2.}} We prove now \eqref{box4}. Assume that $k_1=\min\{k,k_1,k_2\}$. We may also assume $|k-k_2|\leq 4$ and $2^{k_1}\gtrsim\langle t\rangle^{-1}$  (otherwise the bounds follow from \eqref{box2}). We estimate first
\begin{equation}\label{box16.1}
\begin{split}
2^{-k/2}&\big\Vert P_kI[P_{k_1}U^{\mathcal{L}_1h_1,\iota_1},P_{k_2}U^{\mathcal{L}_2h_2,\iota_2}](t)\big\Vert_{L^2}\lesssim 2^{-k_2/2}\big\|P_{k_1}U^{\mathcal{L}_1h_1}(t)\big\|_{L^\infty}\|P_{k_2}U^{\mathcal{L}_2h_2}(t)\|_{L^2}\\
&\lesssim \varep_1^2\langle t\rangle^{-1+[H(q_2,n_2)+H(q_1+1,n_1+1)+1]\delta}2^{-N(n_2)k_2^+}2^{k_1}2^{-5k_1^+},
\end{split}
\end{equation}
using \eqref{vcx1} and \eqref{wws1}. This suffices to prove \eqref{box4} if $2^{k_2}\lesssim \langle t\rangle^{-\delta'}$ or if $2^{k_2}\gtrsim \langle t\rangle^{\delta'}$. Also, the desired bounds follow directly from \eqref{box2} if $2^{k_1}\lesssim \langle t\rangle^{-1+10\delta}$. 

It remains to consider the case $\langle t\rangle \gg 1$, $2^{k_2}\in[\langle t\rangle^{-\delta'},\langle t\rangle^{\delta'}]$, $2^{k_1}\geq\langle t\rangle^{-1+10\delta}$. Let $J_1$ be the largest integer such that $2^{J_1}\le \langle t\rangle 2^{-30}$ and decompose $P_{k_1}U^{\mathcal{L}_1h_1,\iota_1}(t)$ as in \eqref{box15}. Let $J_2$ be the largest integer such that $2^{J_2}\le \langle t\rangle^{1/2}$ and decompose $P_{k_2}U^{\mathcal{L}_2h_2,\iota_2}(t)$ in a similar way. Then
\begin{equation}\label{box16.2}
\begin{split}
2^{-k/2}&\big\Vert P_kI[U_{\leq J_1,k_1}^{\mathcal{L}_1h_1,\iota_1}(t), P_{k_2}U^{\mathcal{L}_2h_2,\iota_2}(t)]\big\Vert_{L^2}\lesssim 2^{-k_2/2}\big\|U_{\leq J_1,k_1}^{\mathcal{L}_1h_1}(t)\big\|_{L^\infty}\|P_{k_2}U^{\mathcal{L}_2h_2}(t)\|_{L^2}\\
&\lesssim \varep_1^2\langle t\rangle^{-1+\delta[H(q_2,n_2)+H(q_1,n_1+1)+1]}2^{-N(n_2)k_2^+}2^{k_1}2^{-5k_1^+},
\end{split}
\end{equation}
using \eqref{vcx1} and \eqref{wws11}. Moreover, 
\begin{equation}\label{box16.3}
\begin{split}
2^{-k/2}&\big\Vert P_kI[U_{>J_1,k_1}^{\mathcal{L}_1h_1,\iota_1}(t), U_{>J_2,k_2}^{\mathcal{L}_2h_2,\iota_2}(t)]\big\Vert_{L^2}\lesssim 2^{-k_2/2}\big\|U_{>J_1,k_1}^{\mathcal{L}_1h_1}(t)\big\|_{L^\infty}\|U_{>J_2,k_2}^{\mathcal{L}_2h_2}(t)\|_{L^2}\\
&\lesssim \varep_1^2\langle t\rangle^{-5/4}2^{-k_2}2^{-N(n_2+1)k_2^+}2^{k_1}2^{-5k_1^+},
\end{split}
\end{equation}
using \eqref{vcx1.1} and \eqref{wws1}. In addition, for any $j_1>J_1$ we use \eqref{rex2} and \eqref{vcx1} to estimate
\begin{equation}\label{box16.4}
\begin{split}
2^{-k/2}&\big\Vert P_kI[U_{j_1,k_1}^{\mathcal{L}_1h_1,\iota_1}(t), U_{\leq J_2,k_2}^{\mathcal{L}_2h_2,\iota_2}(t)]\big\Vert_{L^2}\\
&\lesssim 2^{-k_2/2}2^{k_1/2}\langle t\rangle^{-1+3\delta/4}\big\|Q_{j_1,k_1}V^{\mathcal{L}_1h_1}(t)\big\|_{H^{0,1}_\Omega}\|Q_{\leq J_2,k_2}V^{\mathcal{L}_2h_2}(t)\|_{L^2}\\
&\lesssim \varep_1^2\langle t\rangle^{-1+3\delta/4+H(q_1,n_1+1)\delta+H(q_2,n_2)\delta}2^{k_1}2^{-5k_1^+}2^{-N(n_2)k_2^+}.
\end{split}
\end{equation}
Finally, when $2^{j_1}\geq\langle t\rangle^2$ then we just use \eqref{vcx1} and \eqref{vcx1.1} to estimate
\begin{equation}\label{box16.5}
\begin{split}
2^{-k/2}&\big\Vert P_kI[U_{j_1,k_1}^{\mathcal{L}_1h_1,\iota_1}(t), U_{\leq J_2,k_2}^{\mathcal{L}_2h_2,\iota_2}(t)]\big\Vert_{L^2}\lesssim 2^{-k_2/2}2^{3k_1/2}\big\|V_{j_1,k_1}^{\mathcal{L}_1h_1}(t)\big\|_{L^2}\|V_{\leq J_2,k_2}^{\mathcal{L}_2h_2}(t)\|_{L^2}\\
&\lesssim \varep_1^22^{-3j_1/4}2^{k_1}2^{-5k_1^+}2^{-N(n_2)k_2^+}.
\end{split}
\end{equation}
The desired bounds \eqref{box4} follow from \eqref{box16.2}--\eqref{box16.5} when $2^{k_2}\in[\langle t\rangle^{-\delta'},\langle t\rangle^{\delta'}]$.

{\bf{Step 3.}} Finally, we prove \eqref{box5}--\eqref{box5.6}. Assume that $k_2=\min\{k,k_1,k_2\}$. We may also assume $|k-k_1|\leq 4$, $2^{k_2}\gtrsim\langle t\rangle^{-1}$, and $2^{k_1}\in[\langle t\rangle^{-4/5},\langle t\rangle]$ (otherwise the bounds follow from \eqref{box2}). Let $J_1$ be the largest integer such that $2^{J_1}\le \langle t\rangle(1+2^{k_1}\langle t\rangle)^{-\delta/20}$ and decompose $P_{k_1}U^{\mathcal{L}_1h_1,\iota_1}(t)$ as in \eqref{box15}. Using \eqref{Bil31}, \eqref{vcx1}, and \eqref{vcx1.1} we estimate
\begin{equation}\label{box21.1}
\begin{split}
2^{-k/2}&\big\Vert P_kI[U_{\leq J_1,k_1}^{\mathcal{L}_1h_1,\iota_1}(t), P_{k_2}U^{\mathcal{L}_2h_2,\iota_2}(t)]\big\Vert_{L^2}\\
&\lesssim 2^{-k_1/2}\langle t\rangle^{-1+\delta/4}2^{k_2/2}\big\|Q_{\leq J_1,k_1}V^{\mathcal{L}_1h_1}(t)\big\|_{H^{0,1}_\Omega}\|P_{k_2}U^{\mathcal{L}_2h_2}(t)\|_{L^2}\\
&\lesssim \varep_1^2\langle t\rangle^{-1+\delta[H(q_2,n_2)+H(q_1,n_1+1)+1]}2^{k_2}2^{-5k_2^+}2^{-N(n_1+1)k_1^+}
\end{split}
\end{equation}
and
\begin{equation}\label{box21.2}
\begin{split}
2^{-k/2}&\big\Vert P_kI[U_{> J_1,k_1}^{\mathcal{L}_1h_1,\iota_1}(t), P_{k_2}U^{\mathcal{L}_2h_2,\iota_2}(t)]\big\Vert_{L^2}\lesssim  2^{-k_1/2}2^{3k_2/2}\big\|U_{>J_1,k_1}^{\mathcal{L}_1h_1}(t)\big\|_{L^2}\|P_{k_2}U^{\mathcal{L}_2h_2}(t)\|_{L^2}\\
&\lesssim \varep_1^2\langle t\rangle^{-1+\delta[H(q_2,n_2)+H(q_1+1,n_1+1)+1]}2^{-5k_2^+}2^{2k_2}2^{-k_1}2^{-N(n_1+1)k_1^+}.
\end{split}
\end{equation}
Since $H(1,1)=30$, these bounds clearly suffice to prove \eqref{box5.6} when $(n_1,n_2)\in\{(0,2),(0,3)\}$.

It remains to consider the cases $(n_1,n_2)\in\{(0,0),(0,1),(1,1),(1,2)\}$. Assume first that $2^{k}\notin[\langle t\rangle^{-8\delta'},\langle t\rangle^{8\delta'}]$. Then we estimate, using just \eqref{vcx1} and \eqref{wws1},
\begin{equation}\label{box21.4}
\begin{split}
2^{-k/2}&\big\Vert P_kI[P_{k_1}U^{\mathcal{L}_1h_1,\iota_1}(t), P_{k_2}U^{\mathcal{L}_2h_2,\iota_2}(t)]\big\Vert_{L^2}\\
&\lesssim  2^{-k_1/2}\big\|P_{k_1}U^{\mathcal{L}_1h_1}(t)\big\|_{L^2}\|P_{k_2}U^{\mathcal{L}_2h_2}(t)\|_{L^\infty}\lesssim \varep_1^2\langle t\rangle^{-1+\delta'/2}2^{k_2}2^{-5k_2^+}2^{-N(n_1)k_1^+},
\end{split}
\end{equation}
which suffices to prove \eqref{box5}--\eqref{box5.5} when $2^{|k|}\geq\langle t\rangle^{8\delta'}$. 

On the other hand, if 
\begin{equation}\label{box21.6}
(n_1,n_2)\in\{(0,0),(0,1)\}\qquad\text{ and }\qquad 2^{k}\in [\langle t\rangle^{-8\delta'},\langle t\rangle^{8\delta'}],
\end{equation}
then we would like to use Lemma \ref{Lembil1}. Let $J'_1$ be the largest integer such that $2^{J'_1}\le \langle t\rangle^{1/2}2^{-k_1/2}$ and decompose $P_{k_1}U^{h_1,\iota_1}(t)$ as in \eqref{box15}. Using \eqref{bil1}, \eqref{vcx1}, and \eqref{vcx1.2} we estimate
\begin{equation}\label{box21.7}
\begin{split}
2^{-k/2}&\big\Vert P_kI[U_{\leq J'_1,k_1}^{h_1,\iota_1}(t), P_{k_2}U^{\mathcal{L}_2h_2,\iota_2}(t)]\big\Vert_{L^2}\\
&\lesssim 2^{-k_1/2+k_2/2}\langle t\rangle^{-1}2^{3k_1/2}\big\|\mathcal{F}\{Q_{\leq J'_1,k_1}V^{h_1}\}(t)\big\|_{L^\infty}\|P_{k_2}U^{\mathcal{L}_2h_2,\iota_2}(t)\|_{L^2}\\
&\lesssim \varep_1^2\langle t\rangle^{-1+\delta[H(q_2,n_2)+1]}2^{k_2}2^{-5k_2^+}2^{-\kappa k_1^-}2^{-N_0k_1^+}.
\end{split}
\end{equation}
Moreover, using \eqref{vcx1.1} and \eqref{wws1} we estimate
\begin{equation}\label{box21.8}
\begin{split}
2^{-k/2}&\big\Vert P_kI[U_{> J'_1,k_1}^{h_1,\iota_1}(t), P_{k_2}U^{\mathcal{L}_2h_2,\iota_2}(t)]\big\Vert_{L^2}\lesssim  2^{-k_1/2}\big\|Q_{>J'_1,k_1}V^{h_1}(t)\big\|_{L^2}\|P_{k_2}U^{\mathcal{L}_2h_2}(t)\|_{L^\infty}\\
&\lesssim \varep_1^2\langle t\rangle^{-3/2+\delta'}2^{-k_1/2}2^{k_2}2^{-N(1)k_1^+-5k_2^+}.
\end{split}
\end{equation}
These two bounds clearly suffice to prove \eqref{box5.5} when $2^{k}\in [\langle t\rangle^{-8\delta'},\langle t\rangle^{8\delta'}]$.

Finally, assume that 
\begin{equation}\label{box21.9}
(n_1,n_2)\in\{(1,1),(1,2)\}\qquad\text{ and }\qquad 2^{k}\in [\langle t\rangle^{-8\delta'},\langle t\rangle^{8\delta'}].
\end{equation}
Let $J_1$ be the largest integer such that $2^{J_1}\le \langle t\rangle(1+2^{k_1}\langle t\rangle)^{-\delta/20}$ as before and notice that \eqref{box21.1} give suitable bounds for the contributions of $U_{\leq J_1,k_1}^{\mathcal{L}_1h_1,\iota_1}(t)$. Moreover,
\begin{equation}\label{box21.10}
\begin{split}
2^{-k/2}&\big\Vert P_kI[U_{> J_1,k_1}^{\mathcal{L}_1h_1,\iota_1}(t), P_{k_2}U^{\mathcal{L}_2h_2,\iota_2}(t)]\big\Vert_{L^2}\lesssim  2^{-k_1/2}\big\|Q_{>J_1,k_1}V^{\mathcal{L}_1h_1}(t)\big\|_{L^2}\|P_{k_2}U^{\mathcal{L}_2h_2}(t)\|_{L^\infty}\\
&\lesssim \varep_1^2\langle t\rangle^{-2+\delta'}2^{-5k_2^+}2^{k_2}2^{-k_1}2^{-N(n_1+1)k_1^+},
\end{split}
\end{equation}
using \eqref{vcx1.1} and \eqref{wws1}. This suffices to complete the proof of \eqref{box5}.
\end{proof}

We prove now $L^2$ bounds on localized bilinear interactions of the Klein-Gordon field. 

\begin{lemma}\label{box31}
Assume that $\mathcal{L}_1\in\mathcal{V}_{n_1}^{q_1}$, $\mathcal{L}_2\in\mathcal{V}_{n_2}^{q_2}$, $n_1+n_2\leq 3$, $n_1\leq n_2$. Assume also that $m\in\mathcal{M}$ (see \eqref{mults}), $I=I_m$ is defined as in \eqref{abc36.1}, and let 
\begin{equation}\label{box31.5}
I_{k,k_1,k_2}^{wa,2}(t):=2^{-k/2}\big\Vert P_kI[P_{k_1}U^{\mathcal{L}_1\psi,\iota_1},P_{k_2}U^{\mathcal{L}_2\psi,\iota_2}](t)\big\Vert_{L^2},
\end{equation}
for any $t\in[0,T]$, $\iota_1,\iota_2\in\{+,-\}$,  and $k,k_1,k_2\in\mathbb{Z}$. Then
\begin{equation}\label{box32}
\begin{split}
&I_{k,k_1,k_2}^{wa,2}(t)\lesssim \varep_1^22^{-k/2}2^{3\min\{k,k_1,k_2\}/2}2^{-N(n_1)k_1^+-N(n_2)k_2^+}\\
&\times\min\big\{\langle t\rangle^{H(q_1,n_1)\delta},2^{k_1^-}\langle t\rangle^{H(q_1+1,n_1+1)\delta}\big\}\min\big\{\langle t\rangle^{H(q_2,n_2)\delta},2^{k_2^-}\langle t\rangle^{H(q_2+1,n_2+1)\delta}\big\},
\end{split}
\end{equation}
where the second factor in the right-hand side is, by definition, $\langle t\rangle^{H(q_2,n_2)\delta}$ where $n_2=3$. Moreover
\begin{equation}\label{box33}
\begin{split}
I_{k,k_1,k_2}^{wa,2}(t)\lesssim \varep_1^2\langle t\rangle^{-1+\delta[H(q_1+q_2,n_1+n_2)+\ell(q_1+q_2,n_1+n_2)]}2^{-N(n_1+n_2)k^+}2^{-k^+/4}.
\end{split}
\end{equation}
\end{lemma}

\begin{proof} As in the previous lemma, we remark that bilinear $\rm{KG}\times\mathrm{KG}$ interactions appear in the metric nonlinearities $\mathcal{N}^{h}_{\al\be}$, in semilinear form. According to the general philosophy described in \eqref{cnb10}, ideally we would like to have bounds on the form
\begin{equation}\label{box35.7}
I_{k,k_1,k_2}^{wa,2}(t)\lesssim \varep_1^2\min(2^{k},\langle t\rangle^{-1})2^{-N(n_1+n_2)k^+-k^+/2}.
\end{equation}
We notice that the estimates \eqref{box32}--\eqref{box33} we actually prove are variations of these optimal bounds, with small $\langle t\rangle^{C\delta}$ loss of decay in some cases.

The estimates \eqref{box32} follow using just $L^2$ bounds, see \eqref{abc36.6}, \eqref{vcx1}, and \eqref{cnb2}. To prove \eqref{box33}, we will have to consider several cases. We will sometimes use the following general bounds
\begin{equation}\label{sac34}
\begin{split}
I^{wa,2}_{k,k_1,k_2}(t)&\lesssim 2^{-k/2}\|P_{k_1}U^{\mathcal{L}_1\psi}(t)\|_{L^\infty}\|P_{k_2}U^{\mathcal{L}_2\psi}(t)\|_{L^2}\\
&\lesssim \varep_1^2\langle t\rangle^{-1+\delta(H(q_1+1,n_1+1)+H(q_2,n_2))}2^{-k/2}2^{k_1^-/2}2^{-N(n_2)k_2^+}2^{-N(n_1+1)k_1^++2k_1^+},
\end{split}
\end{equation} 
which follow from \eqref{vcx1} and \eqref{wws2}. Similarly, if $n_2\leq 2$,
\begin{equation}\label{sac34*}
\begin{split}
I^{wa,2}_{k,k_1,k_2}(t)&\lesssim 2^{-k/2}\|P_{k_1}U^{\mathcal{L}_1\psi}(t)\|_{L^2}\|P_{k_2}U^{\mathcal{L}_2\psi}(t)\|_{L^\infty}\\
&\lesssim \varep_1^2\langle t\rangle^{-1+\delta(H(q_1,n_1)+H(q_2+1,n_2+1))}2^{-k/2}2^{k_2^-/2}2^{-N(n_1)k_1^+}2^{-N(n_2+1)k_2^++2k_2^+}.
\end{split}
\end{equation}

We can prove one more general bound of this type when $n_2\leq 2$ by decomposing the profiles. Indeed, let
\begin{equation}\label{abc43.5}
P_{k_1}U^{\mathcal{L}_1\psi,\iota_1}=\sum_{j_1\geq\max(-k_1,0)}U^{\mathcal{L}_1\psi,\iota_1}_{j_1,k_1},\qquad P_{k_2}U^{\mathcal{L}_2\psi,\iota_2}=\sum_{j_2\geq\max(-k_2,0)}U^{\mathcal{L}_2\psi,\iota_2}_{j_2,k_2},
\end{equation}
as in \eqref{on11.3}--\eqref{on11.36}. In view of \eqref{vcx1.1} and \eqref{wws2*}, we have
\begin{equation}\label{abc43}
\begin{split}
\|U^{\mathcal{L}_l\psi,\iota_l}_{j_l,k_l}(t)\|_{L^2}&\lesssim 2^{-j_l}\varep_1Y(k_l,t;q_l,n_l),\\
\|U^{\mathcal{L}_l\psi,\iota_l}_{j_l,k_l}(t)\|_{L^\infty}&\lesssim 2^{3k_l^+}\langle t\rangle^{-3/2}2^{j_l/2}\varep_1Y(k_l,t;q_l,n_l),
\end{split}
\end{equation}
for $l\in\{1,2\}$, where $Y$ is defined as in \eqref{wws3}. We use the $L^2\times L^\infty$ estimate for each interaction (as in \eqref{sac34}--\eqref{sac34*}), and place the factor with the larger $j$ in $L^2$ (in order to gain $2^{-\max(j_1,j_2)}$) and the factor with the smaller $j$ in $L^\infty$. After summation over $j_1,j_2$, it follows that
\begin{equation}\label{box37.6}
I^{wa,2}_{k,k_1,k_2}(t)\lesssim \varep_1^2Y(k_1,t;q_1,n_1)Y(k_2,t;q_2,n_2)2^{3(k_1^++k_2^+)}2^{-k/2}\langle t\rangle^{-3/2}2^{\min(k_1^-,k_2^-)/2}.
\end{equation}

{\bf{Step 1.}} Assume first that $n_1=1$, thus $(n_1,n_2)\in\{(1,1),(1,2)\}$. The desired bounds \eqref{box33} follow from \eqref{box32} if $2^{k_1^-+k_2^-+k^-}\lesssim \langle t\rangle^{-1-\delta'}$. They also follow if $2^{\max\{k_1,k_2\}}\gtrsim\langle t\rangle^{\delta'}$, using from \eqref{sac34}--\eqref{sac34*} if $k\geq 0$ and \eqref{box37.6} if $k\leq 0$. If $2^{\max\{k_1,k_2\}}\leq \langle t\rangle^{\delta'}$ then the bounds \eqref{box33} follow from \eqref{box37.6} if $2^k\gtrsim \langle t\rangle^{-1+10\delta'}$. After these reductions, it remains to prove \eqref{box33} when $|t|\gg 1$ and
\begin{equation}\label{box58}
2^{k_1},2^{k_2}\in[\langle t\rangle^{-8\delta'},\langle t\rangle^{\delta'}],\qquad 2^{k}\leq \langle t\rangle^{-1+10\delta'}.
\end{equation}

This is a $\mathrm{High}\times\mathrm{High}\to\mathrm{Low}$ interaction and loss of derivatives is not an issue, so we only need to justify the $\langle t\rangle^{-1+\delta[H(q_1+q_2,n_1+n_2)+33]}$ time decay. The bounds still follow from \eqref{box32} if $2^k\lesssim \langle t\rangle^{-1+30\delta}$, see \eqref{SuperlinearH1}. On the other hand, if $2^k\geq\langle t\rangle^{-1+30\delta}$ then we let $J_1$ be the largest integer such that $2^{J_1}\le \langle t\rangle^{3/4}$. Using \eqref{LemBil43}, \eqref{vcx1}--\eqref{cnb2}, and \eqref{wws2} we estimate
\begin{equation}\label{box58.1}
\begin{split}
2^{-k/2}&\big\Vert P_kI[U_{\leq J_1,k_1}^{\mathcal{L}_1\psi,\iota_1}(t), P_{k_2}U^{\mathcal{L}_2\psi,\iota_2}(t)]\big\Vert_{L^2}\\
&\lesssim \langle t\rangle^{-1+\delta}2^{5k_1^+}2^{-k_1^-}\|Q_{\leq J_1,k_1}V^{\mathcal{L}_1\psi}(t)\|_{H^{0,1}_\Omega}\|P_{k_2}U^{\mathcal{L}_2\psi}(t)\|_{L^2}\\
&\lesssim \varep_1^2\langle t\rangle^{-1+\delta[H(q_1+1,n_1+1)+H(q_2,n_2)+1]}
\end{split}
\end{equation}
and
\begin{equation}\label{box58.2}
\begin{split}
2^{-k/2}&\big\Vert P_kI[U_{> J_1,k_1}^{\mathcal{L}_1\psi,\iota_1}(t), P_{k_2}U^{\mathcal{L}_2\psi,\iota_2}(t)]\big\Vert_{L^2}\\
&\lesssim 2^{-k/2}\big\|Q_{>J_1,k_1}V^{\mathcal{L}_1\psi}(t)\big\|_{L^2}\|P_{k_2}U^{\mathcal{L}_2\psi}(t)\|_{L^\infty}\lesssim \varep_1^22^{-k/2}\langle t\rangle^{-7/4+\delta'}.
\end{split}
\end{equation}
We use now the last inequality in \eqref{SuperlinearH2}. The desired bounds \eqref{box33} follow if $q_2\geq 1$. 

Moreover, we can repeat the estimates \eqref{box58.1}--\eqref{box58.2} with the roles of the functions $P_{k_1}U^{\mathcal{L}_1\psi,\iota_1}$ and $P_{k_2}U^{\mathcal{L}_2\psi,\iota_2}$ reversed. Thus $I^{wa,2}_{k,k_1,k_2}(t)\lesssim \varep_1^2\langle t\rangle^{-1+\delta[H(q_2+1,n_2+1)+H(q_1,n_1)+1]}$, and the desired conclusion follows if $q_1\geq 1$. 

Finally, if $q_1=q_2=0$ then we bound $2^{-k_1^-}\|P_{k_1}U^{\mathcal{L}_1\psi}\|_{L^2}\lesssim Y(k,t;0,0)\lesssim \langle t\rangle^{H(1,1)}2^{-N(1)k_1^+}$, as a consequence of \eqref{cnb2} and the assumption $n_1=1$. Then we estimate, as in \eqref{box58.1}
\begin{equation*}
\begin{split}
2^{-k/2}&\big\Vert P_kI[P_{k_1}U^{\mathcal{L}_1\psi,\iota_1}(t),U_{\leq J_2,k_2}^{\mathcal{L}_2\psi,\iota_2}(t)]\big\Vert_{L^2}\\
&\lesssim \langle t\rangle^{-1+\delta}2^{5k_2^+}2^{-k_2^-}\|P_{k_1}U^{\mathcal{L}_1\psi}(t)\|_{L^2}\|Q_{\leq J_2,k_2}V^{\mathcal{L}_2\psi}(t)\|_{H^{0,1}_\Omega}\lesssim \varep_1^2\langle t\rangle^{-1+\delta[H(1,1)+H(0,n_2+1)+1]},
\end{split}
\end{equation*}
 where $J_2$ is the largest integer such that $2^{J_2}\le \langle t\rangle^{3/4}$. Since $H(1,1)=30$, this is consistent with the desired estimates \eqref{box33}. The contribution of $U_{> J_2,k_2}^{\mathcal{L}_2\psi,\iota_2}(t)$ can be bounded as in \eqref{box58.2}, using an $L^\infty\times L^2$ estimate and \eqref{vcx1.1}. This completes the proof of \eqref{box33} when $n_1=1$.

{\bf{Step 2.}} Assume now that $n_1=0$. Notice that $\|P_{k_1}U^{\psi,\iota_1}(t)\|_{L^2}\lesssim 2^{k_1^-+\kappa k_1^-}2^{-N_0k_1^++2k_1^+}$, as a consequence of \eqref{vcx1.2}. Therefore
\begin{equation}\label{box37.1}
I_{k,k_1,k_2}^{wa,2}(t)\lesssim  \varep_1^22^{-k^+/2}2^{\min\{k^-,k_1^-,k_2^-\}}2^{k_1^-+\kappa k_1^-}2^{-N_0k_1^++4k_1^+}2^{-N(n_2)k_2^+}\langle t\rangle^{H(q_2,n_2)\delta},
\end{equation}
as a consequence of \eqref{vcx1} and \eqref{abc36.6}. The desired bounds \eqref{box33} follow if $2^{k_1^-}\lesssim \langle t\rangle^{-1/2+\kappa/8}$ or if ($2^{k_1^-}\geq \langle t\rangle^{-1/2+\kappa/8}2^{20}$ and $2^{k^-+k_1^-}\lesssim\langle t\rangle^{-1}$).

On the other hand, if $\langle t\rangle\geq 2^{1/\delta}$,
\begin{equation}\label{box59.1}
2^{k_1^-}\geq \langle t\rangle^{-1/2+\kappa/8}\qquad\text{ and }\qquad 2^{k^-+k_1^-}\geq \langle t\rangle^{-1}
\end{equation}
then we let $J_1$ be the largest integer such that $2^{J_1}\le 2^{k_1^--20}\langle t\rangle$ and decompose $P_{k_1}U^{\psi,\iota_1}(t)=e^{-it\Lambda_{kg,\iota_1}}V_{\leq J_1,k_1}^{\psi,\iota_1}(t)+e^{-it\Lambda_{kg,\iota_1}}V_{>J_1,k_1}^{\psi,\iota_1}(t)$ as in \eqref{box15}. Using \eqref{wws13x}, \eqref{vcx1}, and \eqref{vcx1.1} we estimate
\begin{equation}\label{box37.2}
\begin{split}
2^{-k/2}&\big\Vert P_kI[U_{\leq J_1,k_1}^{\psi,\iota_1}(t), P_{k_2}U^{\mathcal{L}_2\psi,\iota_2}(t)]\big\Vert_{L^2}\lesssim 2^{-k/2}\|U_{\leq J_1,k_1}^{\psi}(t)\|_{L^\infty}\|P_{k_2}U^{\mathcal{L}_2\psi}(t)\|_{L^2}\\
&\lesssim \varep_1^2\langle t\rangle^{-3/2+\delta[H(q_2,n_2)+1]}2^{-k/2}2^{-k_1^-/2}2^{-N_0k_1^++5k_1^+}2^{-N(n_2)k_2^+}
\end{split}
\end{equation}
and, with $\underline{k}=\min\{k,k_1,k_2\}$, 
\begin{equation}\label{box37.3}
\begin{split}
2^{-k/2}&\big\Vert P_kI[U_{> J_1,k_1}^{\psi,\iota_1}(t), P_{k_2}U^{\mathcal{L}_2\psi,\iota_2}(t)]\big\Vert_{L^2}\lesssim 2^{-k/2}2^{3\underline{k}/2}\big\|U_{>J_1,k_1}^{\psi}(t)\big\|_{L^2}\|P_{k_2}U^{\mathcal{L}_2\psi}(t)\|_{L^2}\\
&\lesssim \varep_1^2\langle t\rangle^{-1+\delta[H(q_2,n_2)+H(1,1)+1]}2^{-k_1-k/2}2^{3\underline{k}/2}2^{-N(1)k_1^+}2^{-N(n_2)k_2^+}.
\end{split}
\end{equation}

{\bf{Case 2.1.}} Assume first that $k=\min(k,k_1,k_2)$ and $2^{k}\geq\langle t\rangle^{-1}$. Then $|k_1-k_2|\leq 4$ and the bounds \eqref{box37.2}--\eqref{box37.3} show that
\begin{equation}\label{box37.4}
I^{wa,2}_{k,k_1,k_2}(t)\lesssim \varep_1^2\langle t\rangle^{-1+\delta[H(q_2,n_2)+1]}2^{-N(n_2)k_2^+-5k_2^+}\big\{2^{-k/2}2^{-k_1^-/2}\langle t\rangle^{-1/2}+2^{k-k_1}\langle t\rangle^{H(1,1)\delta}\big\}.
\end{equation}
Since $H(1,1)=30$ and $2^{-k/2}2^{-k_1^-/2}\langle t\rangle^{-1/2}\lesssim 1$ (see \eqref{box59.1}), this suffices when $n_2\geq 2$ or when $2^{k-k_1}2^{-4k_2^+}\leq\langle t\rangle^{-\delta'}$. In the remaining case ($n_2\in\{0,1\}$ and $2^{k-k_1}2^{-4k_2^+}\in[\langle t\rangle^{-\delta'},1]$) we need to improve the bounds \eqref{box37.3}. Using \eqref{vcx1.1} and \eqref{wws2} we estimate
\begin{equation*}
\begin{split}
2^{-k/2}\big\Vert P_kI[U_{> J_1,k_1}^{\psi,\iota_1}(t), P_{k_2}U^{\mathcal{L}_2\psi,\iota_2}(t)]\big\Vert_{L^2}&\lesssim 2^{-k/2}\big\|U_{>J_1,k_1}^{\psi}(t)\big\|_{L^2}\|P_{k_2}U^{\mathcal{L}_2\psi}(t)\|_{L^\infty}\\
&\lesssim \varep_1^2\langle t\rangle^{-2+\delta'}2^{-k_1^-/2-k/2}.
\end{split}
\end{equation*}
Since $2^{4k^+}\leq\langle t\rangle^{\delta'}$ this can be combined with \eqref{box37.2} to complete the proof of \eqref{box33}. 

{\bf{Case 2.2.}} Assume now that $k_1=\min(k,k_1,k_2)$. Then $|k-k_2|\leq 4$ and the bounds \eqref{box37.2}--\eqref{box37.3} show that
\begin{equation*}
I^{wa,2}_{k,k_1,k_2}(t)\lesssim \varep_1^2\langle t\rangle^{-1+\delta[H(q_2,n_2)+1]}2^{-N(n_2)k_2^+-5k_1^+}\big\{2^{-k/2}2^{-k_1^-/2}\langle t\rangle^{-1/2}+2^{(k_1^--k)/2}\langle t\rangle^{H(1,1)\delta}\big\}.
\end{equation*}
In view of \eqref{box59.1}, this suffices when $n_2\geq 2$ or when $2^{(k_1^--k^-)/2}2^{-k^+/4}\leq\langle t\rangle^{-\delta'}$. In the remaining case ($n_2\in\{0,1\}$ and $2^{(k_1^--k^-)/2}2^{-k^+/4}\in[\langle t\rangle^{-\delta'},1]$) we improve the bounds \eqref{box37.3} by estimating 
\begin{equation}\label{save1}
\begin{split}
2^{-k/2}&\big\Vert P_kI[U_{> J_1,k_1}^{\psi,\iota_1}(t), P_{k_2}U^{\mathcal{L}_2\psi,\iota_2}(t)]\big\Vert_{L^2}\\
&\lesssim 2^{-k/2}\big\|U_{>J_1,k_1}^{\psi}(t)\big\|_{L^2}\|P_{k_2}U^{\mathcal{L}_2\psi}(t)\|_{L^\infty}\lesssim \varep_1^2\langle t\rangle^{-2+\delta'}2^{-k_1^-},
\end{split}
\end{equation}
using \eqref{vcx1.1} and \eqref{wws2}. Since $2^{k^+}\lesssim \langle t\rangle^{4\delta'}$ and $2^{-k_1^-}\lesssim \langle t\rangle^{1/2}$, the desired bounds \eqref{box33} follow.

{\bf{Case 2.3.}} Finally, assume that $k_2=\min(k,k_1,k_2)$. In proving \eqref{box33} we may also assume that $n_2\geq 1$, since the case $n_2=0$ follows from the analysis in Case 2.2 by reversing the roles of the functions $P_{k_1}U^{\psi,\iota_1}$ and $P_{k_2}U^{\psi,\iota_2}$. The bounds \eqref{box37.2}--\eqref{box37.3} show that
\begin{equation*}
I^{wa,2}_{k,k_1,k_2}(t)\lesssim \varep_1^2\langle t\rangle^{-1+\delta[H(q_2,n_2)+1]}2^{-4k_2^+-N(1)k_1^+}\big\{2^{-k/2}2^{-k_1^-/2}\langle t\rangle^{-1/2}2^{-4k_1^+}+2^{k_2^--k}\langle t\rangle^{H(1,1)\delta}\big\}.
\end{equation*}
In view of \eqref{box59.1}, this suffices when $n_2\geq 2$ or when ($n_2=1$  and $2^{k_2^--k^-}2^{-k^+/2}\leq\langle t\rangle^{-\delta'}$). In the remaining case ($n_2=1$ and $2^{k_2^--k^-}2^{-k^+/2}\in[\langle t\rangle^{-\delta'},1]$) we use \eqref{vcx1.1} and \eqref{wws2} to prove bounds identical to \eqref{save1},
\begin{equation*}
\begin{split}
2^{-k/2}\big\Vert P_kI[U_{> J_1,k_1}^{\psi,\iota_1}(t), P_{k_2}U^{\mathcal{L}_2\psi,\iota_2}(t)]\big\Vert_{L^2}\lesssim \varep_1^2\langle t\rangle^{-2+\delta'}2^{-k_1^-}.
\end{split}
\end{equation*}
The desired bounds \eqref{box33} follow in this last case. This completes the proof of the lemma.
\end{proof}

Finally we prove $L^2$ bounds on localized bilinear interactions of the Klein-Gordon field and the metric components. 

\begin{lemma}\label{box80}
Assume that $\mathcal{L}_1\in\mathcal{V}_{n_1}^{q_1}$, $\mathcal{L}_2\in\mathcal{V}_{n_2}^{q_2}$, $n_1+n_2\leq 3$. Assume also that $m\in\mathcal{M}$ (see \eqref{mults}), $I=I_m$ is defined as in \eqref{abc36.1}, and let
\begin{equation}\label{box81}
I_{k,k_1,k_2}^{kg}(t):=2^{k_2^+-k_1}\big\Vert P_kI[P_{k_1}U^{\mathcal{L}_1h,\iota_1},P_{k_2}U^{\mathcal{L}_2\psi,\iota_2}](t)\big\Vert_{L^2},
\end{equation}
for any $t\in[0,T]$, $\iota_1,\iota_2\in\{+,-\}$,  $h\in\{h_{\al\be}\}$, and $k,k_1,k_2\in\mathbb{Z}$. Then
\begin{equation}\label{box82}
\begin{split}
I_{k,k_1,k_2}^{kg}&(t)\lesssim \varep_1^22^{-k_1/2}2^{3\min\{k,k_1,k_2\}/2}2^{-N(n_1)k_1^+-N(n_2)k_2^++k_2^+}\\
&\times \langle t\rangle^{H(q_1,n_1)\delta}(\langle t\rangle 2^{k_1^-})^{-\ga}\min\big\{\langle t\rangle^{H(q_2,n_2)\delta},2^{k_2^-}\langle t\rangle^{H(q_2+1,n_2+1)\delta}\big\},
\end{split}
\end{equation}
where the second factor in the right-hand side is, by definition, $\langle t\rangle^{H(q_2,n_2)\delta}$ where $n_2=3$. In addition, we have:

(1) if $n_1=0$ and $n_2\geq 0$ then
\begin{equation}\label{box83}
I_{k,k_1,k_2}^{kg}(t)\lesssim \varep_1^2\langle t\rangle^{-1+\delta[H(q_2,n_2)+33]}2^{-N(n_2)k^++5k^+/4}2^{-2\min\{k_1^+,k_2^+\}};
\end{equation}

(2) if $n_1=n_2=0$ then
\begin{equation}\label{box83.5}
I^{kg}_{k,k_1,k_2}(t)\lesssim 
\begin{cases}
\varep_1^2\langle t\rangle^{-1+4\delta}2^{-2k_2^+}2^{-N(0)k_1^+-k_1^+/4}&\quad\text{ if }k_1\geq k_2;\\
\varep_1^2\langle t\rangle^{-1+4\delta}2^{-2k_1^+}2^{-N_0k_2^++6k_2^+}&\quad\text{ if }k_1\leq k_2;
\end{cases}
\end{equation}

(3) if $n_1\geq 1$ and $n_2=0$ then
\begin{equation}\label{box84}
I_{k,k_1,k_2}^{kg}(t)\lesssim \varep_1^2\langle t\rangle^{-1+\delta(H(q_1,n_1)+\ell(q_1,n_1))}2^{-N(n_1)k^+-k^+/4}2^{-2\min\{k_1^+,k_2^+\}};
\end{equation}

(4) if $n_1\geq 1$ and $n_2\geq 1$ then
\begin{equation}\label{box84.5}
I_{k,k_1,k_2}^{kg}(t)\lesssim \varep_1^2\langle t\rangle^{-1+\delta(H(q_1+q_2,n_1+n_2)+33)}2^{-N(n_1+n_2)k^+-k^+/4}2^{-2\min\{k_1^+,k_2^+\}}.
\end{equation}

\end{lemma}

\begin{proof}
We notice that bilinear $\rm{Wave}\times\mathrm{KG}$ interactions appear in the nonlinearities $\mathcal{N}^{\psi}$, in quasilinear form. As before, ideally we would like to have bounds on the form
\begin{equation}\label{box86}
I_{k,k_1,k_2}^{kg}(t)\lesssim \varep_1^2\min(2^{k},\langle t\rangle^{-1})2^{-N(n_1+n_2)k^++k^+},
\end{equation}
but we are only able to prove variations of these optimal bounds, with small losses.

The estimates \eqref{box82} follow using just $L^2$ bounds (see \eqref{abc36.6}, \eqref{vcx1}, and \eqref{cnb2}).

{\bf{Step 1.}} We observe that we have the following general bounds, which follow from \eqref{vcx1}, \eqref{wws1}, and \eqref{wws2}: if $n_1\leq 2$ then
\begin{equation}\label{box88}
\begin{split}
I^{kg}_{k,k_1,k_2}(t)&\lesssim 2^{k_2^+-k_1}\|P_{k_1}U^{\mathcal{L}_1h}(t)\|_{L^\infty}\|P_{k_2}U^{\mathcal{L}_2\psi}(t)\|_{L^2}\\
&\lesssim \varep_1^2\langle t\rangle^{-1+\delta(H(q_1+1,n_1+1)+H(q_2,n_2))}2^{-N(n_2)k_2^++k_2^+}2^{-N(n_1+1)k_1^++k_1^+}.
\end{split}
\end{equation} 
Similarly, if $n_2\leq 2$ and $2^{k_1}\gtrsim\langle t\rangle^{-1}$ then
\begin{equation}\label{box89}
\begin{split}
I^{kg}_{k,k_1,k_2}(t)&\lesssim 2^{k_2^+-k_1}\|P_{k_1}U^{\mathcal{L}_1h}(t)\|_{L^2}\|P_{k_2}U^{\mathcal{L}_2\psi}(t)\|_{L^\infty}\\
&\lesssim \varep_1^2\langle t\rangle^{-1+\delta(H(q_1,n_1)+H(q_2+1,n_2+1))}2^{-k_1/2}2^{k_2^-/2}2^{-N(n_1)k_1^+}2^{-N(n_2+1)k_2^++3k_2^+}.
\end{split}
\end{equation}
Since $H(1,1)=30$, the bounds \eqref{box83} follow from \eqref{box88} if $n_2\geq 1$. They also follow from \eqref{box88} if $n_2=0$ and $k_1\leq k_2$ and from the stronger bounds \eqref{box83.5} if $n_2=0$ and $k_1\geq k_2$.

{\bf{Step 2.}} Assume now that $n_2=0$ and we prove the bounds \eqref{box83.5}--\eqref{box84}. We have
\begin{equation}\label{box89.2}
\begin{split}
I^{kg}_{k,k_1,k_2}(t)&\lesssim 2^{k_2^+-k_1}2^{3\min\{k_1,k_2\}/2}\|P_{k_1}U^{\mathcal{L}_1h}(t)\|_{L^2}\|P_{k_2}U^{\psi}(t)\|_{L^2}\\
&\lesssim \varep_1^2\langle t\rangle^{\delta H(q_1,n_1)}2^{\min\{k_1^-,k_2^-\}}2^{k_2^-+\kappa k_2^-}(\langle t\rangle 2^{k_1^-})^{-\ga}2^{-N(n_1)k_1^+-k_1^+/2}2^{-N_0k_2^++5k_2^+},
\end{split}
\end{equation}
using just $L^2$ estimates and the inequality $2^{3\min\{k_1,k_2\}/2}\lesssim 2^{k_1^-/2}2^{\min\{k_1^-,k_2^-\}}2^{3k_2^+/2}$. This suffices to prove \eqref{box83.5}--\eqref{box84} if $2^{\min\{k^-_1,k_2^-\}}2^{k_2^-}\lesssim \langle t\rangle^{-1}$. 

Assume now that $\langle t\rangle^{-1}\leq 2^{\min\{k^-_1,k_2^-\}}2^{k_2^--40}$. Let $J_2$ denote the largest integer such that $2^{J_2}\leq 2^{k_2^--20}\langle t\rangle$. Using \eqref{wws13x} and \eqref{vcx1} we estimate
\begin{equation}\label{box89.1}
\begin{split}
2^{k_2^+-k_1}\big\Vert P_kI[P_{k_1}U^{\mathcal{L}_1h,\iota_1},&U^{\psi,\iota_2}_{\leq J_2,k_2}](t)\big\Vert_{L^2}\lesssim 2^{k_2^+-k_1}\|P_{k_1}U^{\mathcal{L}_1h}(t)\|_{L^2}\|U^{\psi}_{\leq J_2,k_2}(t)\|_{L^\infty}\\
&\lesssim \varep_1^2\langle t\rangle^{-3/2+\delta H(q_1,n_1)}2^{-k^-_1/2}2^{-k_2^-/2}2^{-N(n_1)k_1^+-k_1^+/2}2^{-N_0k_2^++6k_2^+}\\
&\lesssim \varep_1^2\langle t\rangle^{-1+\delta H(q_1,n_1)}2^{-N(n_1)k_1^+-k_1^+/2}2^{-N_0k_2^++6k_2^+}.
\end{split}
\end{equation}
Using also \eqref{vcx1.1} and $L^2$ estimates, we also have, with $\underline{k}=\min\{k,k_1,k_2\}$ as before,
\begin{equation}\label{box89.4}
\begin{split}
2^{k_2^+-k_1}\big\Vert P_kI[P_{k_1}U^{\mathcal{L}_1h,\iota_1},&U^{\psi,\iota_2}_{>J_2,k_2}](t)\big\Vert_{L^2}\lesssim 2^{k_2^+-k_1}2^{3\underline{k}/2}\|P_{k_1}U^{\mathcal{L}_1h}(t)\|_{L^2}\|U^{\psi}_{>J_2,k_2}(t)\|_{L^2}\\
&\lesssim \varep_1^2\langle t\rangle^{-1+\delta (H(q_1,n_1)+30)}2^{-N(n_1)k_1^+-k_1^+/2}2^{-N(1)k_2^++4k_2^+}.
\end{split}
\end{equation}
The last two bounds suffice to prove \eqref{box84} when $n_1\geq 2$.

On the other hand, if $n_1\leq 1$ then we can also estimate, using \eqref{wws1},
\begin{equation}\label{box89.5}
\begin{split}
2^{k_2^+-k_1}\big\Vert P_kI[P_{k_1}U^{\mathcal{L}_1h,\iota_1},&U^{\psi,\iota_2}_{>J_2,k_2}](t)\big\Vert_{L^2}\lesssim 2^{k_2^+-k_1}\|P_{k_1}U^{\mathcal{L}_1h}(t)\|_{L^\infty}\|U^{\psi}_{>J_2,k_2}(t)\|_{L^2}\\
&\lesssim \varep_1^2\langle t\rangle^{-2+\delta'}2^{-k_2^-}2^{-N(n_1+1)k_1^++2k_1^+}2^{-N(1)k_2^++2k_2^+}.
\end{split}
\end{equation}
The bounds \eqref{box83.5} and \eqref{box84} with $n_1=1$ follow from \eqref{box89.1} and \eqref{box89.5} if $2^{k_1^+}+2^{k_2^+}\lesssim \langle t\rangle^{4\delta'}$. On the other hand, if $k_1\geq k_2$ and $2^{k_1^+}\geq \langle t\rangle^{4\delta'}$ then the bounds \eqref{box83.5} and \eqref{box84} with $n_1=1$ follow from \eqref{box89.1} and \eqref{box89.4}. Finally, if $k_1\leq k_2$ and $2^{k_2^+}\geq \langle t\rangle^{4\delta'}$ then the bounds \eqref{box83.5} and \eqref{box84} with $n_1=1$ follow from \eqref{box88}. 

This completes the proof of the bounds \eqref{box83.5} and \eqref{box84} in all cases.

{\bf{Step 3.}} We prove now the bounds \eqref{box84.5} when $k_2\leq k_1$. The bounds \eqref{box89} give the desired conclusion if $2^{k_1}\gtrsim\langle t\rangle^{\delta'}$. Also, the bounds \eqref{box82} give the desired conclusion if $2^{k_1}\lesssim\langle t\rangle^{-1/2-\delta'}$.

In the remaining case
\begin{equation}\label{box91.1}
k_2\leq k_1,\qquad 2^{k_1}\in [\langle t\rangle^{-1/2-\delta'},\langle t\rangle^{\delta'}],
\end{equation} 
we decompose $P_{k_2}U^{\mathcal{L}_2\psi,\iota_2}=\sum_{j_2\geq\max(-k_2,0)}U^{\mathcal{L}_2\psi,\iota_2}_{j_2,k_2}$ as in \eqref{abc43.5}. Then we estimate
\begin{equation*}
\begin{split}
2^{k_2^+-k_1}&\big\Vert P_kI[P_{k_1}U^{\mathcal{L}_1h,\iota_1}(t),U_{j_2,k_2}^{\mathcal{L}_2\psi,\iota_2}(t)]\big\Vert_{L^2}\\
&\lesssim 2^{k_2^+-k_1}\|P_{k_1}U^{\mathcal{L}_1h}(t)\|_{L^2}\|U_{j_2,k_2}^{\mathcal{L}_2\psi}(t)\Vert_{L^\infty}\lesssim \varep_1^2\langle t\rangle^{-3/2+2\delta'}2^{j_2/2}2^{-k_1/2},
\end{split}
\end{equation*}
using \eqref{vcx1} and \eqref{wws2*}. This suffices to estimate the contribution of the localized profiles for which $2^{j_2}\lesssim \langle t\rangle^{1/3}$. On the other hand, using \eqref{wws11} and \eqref{vcx1.1} we also estimate 
\begin{equation*}
\begin{split}
2^{k_2^+-k_1}&\big\Vert P_kI[P_{k_1}U^{\mathcal{L}_1h,\iota_1}(t),U_{j_2,k_2}^{\mathcal{L}_2\psi,\iota_2}(t)]\big\Vert_{L^2}\\
&\lesssim 2^{k_2^+-k_1}\|P_{k_1}U^{\mathcal{L}_1h}(t)\|_{L^\infty}\|U_{j_2,k_2}^{\mathcal{L}_2\psi}(t)\Vert_{L^2}\lesssim \varep_1^2\langle t\rangle^{-1+2\delta'}2^{-j_2},
\end{split}
\end{equation*}
which suffices to estimate the contribution of the localized profiles for which $2^{j_2}\geq\langle t\rangle^{1/3}$.

{\bf{Step 4.}} Finally, we prove the bounds \eqref{box84.5} when $k_1\leq k_2$. The bounds \eqref{box88} give the desired conclusion if $2^{k_2}\gtrsim\langle t\rangle^{\delta'}$. Also, the bounds \eqref{box82} and \eqref{SuperlinearH1} give the desired conclusion if $2^{k_1}\lesssim\langle t\rangle^{-1+40\delta}$. 

If $q_2\geq 1$, the bounds \eqref{box84.5} follow from \eqref{box88} and the last inequality in \eqref{SuperlinearH2}. If $q_2=0$ then we let $J_1$ be the largest integer such that $2^{J_1}\le \langle t\rangle 2^{-30}$ and decompose $P_{k_1}U^{\mathcal{L}_1h,\iota_1}(t)$ as in \eqref{box15}. We also decompose $P_{k_2}U^{\mathcal{L}_2h,\iota_2}(t)=U_{\leq J_2,k_2}^{\mathcal{L}_2h,\iota_2}(t)+U_{>J_2,k_2}^{\mathcal{L}_2h,\iota_2}(t)$, where $J_2$ is the largest integer such that $2^{J_2}\le \langle t\rangle^{4\delta'}$. Then
\begin{equation*}
\begin{split}
2^{k_2^+-k_1}&\big\Vert P_kI[U_{\leq J_1,k_1}^{\mathcal{L}_1h,\iota_1}(t), P_{k_2}U^{\mathcal{L}_2\psi,\iota_2}(t)]\big\Vert_{L^2}\lesssim 2^{k_2^+-k_1}\big\|U_{\leq J_1,k_1}^{\mathcal{L}_1h}(t)\big\|_{L^\infty}\|P_{k_2}U^{\mathcal{L}_2\psi}(t)\|_{L^2}\\
&\lesssim \varep_1^2\langle t\rangle^{-1+\delta[H(0,n_2)+H(q_1,n_1+1)+1]}2^{-N(n_2)k_2^++k_2^+}2^{-5k_1^+},
\end{split}
\end{equation*}
using \eqref{vcx1} and \eqref{wws11}. This suffices due to the first inequality in \eqref{SuperlinearH2}. Also
\begin{equation*}
\begin{split}
2^{k_2^+-k_1}&\big\Vert P_kI[U_{>J_1,k_1}^{\mathcal{L}_1h,\iota_1}(t), U_{>J_2,k_2}^{\mathcal{L}_2\psi,\iota_2}(t)]\big\Vert_{L^2}\lesssim 2^{k_2^+-k_1}\big\|U_{>J_1,k_1}^{\mathcal{L}_1h}(t)\big\|_{L^\infty}\|U_{>J_2,k_2}^{\mathcal{L}_2\psi}(t)\|_{L^2}\\
&\lesssim \varep_1^2\langle t\rangle^{-1+\delta[H(1,n_2+1)+H(q_1+1,n_1+1)+1]}2^{-J_2}2^{-N(n_2+1)k_2^++k_2^+}2^{-5k_1^+},
\end{split}
\end{equation*}
using \eqref{vcx1.1} and \eqref{wws1}. For \eqref{box84.5} it remains to prove that
\begin{equation}\label{box93.5}
\begin{split}
2^{k_2^+-k_1}&\big\Vert P_kI[U_{>J_1,k_1}^{\mathcal{L}_1h,\iota_1}(t), U_{\leq J_2,k_2}^{\mathcal{L}_2\psi,\iota_2}(t)]\big\Vert_{L^2}\\
&\lesssim \varep_1^2\langle t\rangle^{-1+\delta(H(q_1,n_1+n_2)+33)}2^{-N(n_1+n_2)k^+-k^+/4}2^{-2k_1^+},
\end{split}
\end{equation}
provided that $|t|\gg 1$ and 
\begin{equation}\label{box93}
q_2=0,\qquad \langle t\rangle^{-1+40\delta}\leq 2^{k_1}\leq 2^{k_2}\leq \langle t\rangle^{\delta'},\qquad 2^{J_2}\leq \langle t\rangle^{4\delta'}.
\end{equation}
Let $X:=2^{k_2^+-k_1}\big\Vert P_kI[U_{>J_1,k_1}^{\mathcal{L}_1h,\iota_1}(t), U_{\leq J_2,k_2}^{\mathcal{L}_2\psi,\iota_2}(t)]\big\Vert_{L^2}$ denote the expression in \eqref{box93.5}.

{\bf{Case 4.1.}} Assume first that $q_1=0$. We use the bounds in Lemma \ref{dtv6}. If $n_1=1$ and $n_2=1$ then we use the $L^4$ bounds in \eqref{UOmegaMoreReg2} and \eqref{UOmegaMoreReg3} to estimate
\begin{equation*}
\begin{split}
X&\lesssim 2^{k_2^+-k_1}\|U_{>J_1,k_1}^{\mathcal{L}_1h}(t)\|_{L^4}\|U_{\leq J_2,k_2}^{\mathcal{L}_2\psi}(t)\|_{L^4}\lesssim \varep_1^2\langle t\rangle^{-5/4+110\delta}2^{-k_1^-/4}2^{-k_2^-/4}2^{-N(2)k_2^+-2k_2^+}2^{-4k_1^+}.
\end{split}
\end{equation*}
Alternatively, we could use the $L^2$ bounds in \eqref{UOmegaMoreReg6}--\eqref{UOmegaMoreReg5} to estimate
\begin{equation*}
\begin{split}
X&\lesssim 2^{k_2^+-k_1}2^{3k_1/2}\|U_{>J_1,k_1}^{\mathcal{L}_1h}(t)\|_{L^2}\|U_{\leq J_2,k_2}^{\mathcal{L}_2\psi}(t)\|_{L^2}\lesssim \varep_1^2\langle t\rangle^{-2/3+110\delta}2^{k_1^-/3}2^{k_2^-}2^{-N(1)k_2^+}2^{-4k_1^+}.
\end{split}
\end{equation*}
We use the first estimates if $\langle t\rangle 2^{k_1^-}2^{k_2^-}\geq 1$ and the second estimates if $\langle t\rangle 2^{k_1^-}2^{k_2^-}\leq 1$. The desired bounds \eqref{box93.5} follow if $n_1=n_2=1$.

Similarly, if $n_1=1$ and $n_2=2$ then we use the $L^6$ bounds in \eqref{UOmegaMoreReg2} and the $L^3$ bounds in \eqref{UOmegaMoreReg3} to estimate
\begin{equation*}
\begin{split}
X&\lesssim 2^{k_2^+-k_1}\|U_{>J_1,k_1}^{\mathcal{L}_1h}(t)\|_{L^6}\|U_{\leq J_2,k_2}^{\mathcal{L}_2\psi}(t)\|_{L^3}\lesssim \varep_1^2\langle t\rangle^{-7/6+170\delta}2^{-k_1^-/6}2^{-k_2^-/6}2^{-N(3)k_2^+-2k_2^+}2^{-4k_1^+}.
\end{split}
\end{equation*}
Alternatively, using the $L^2$ bounds in \eqref{UOmegaMoreReg6}--\eqref{UOmegaMoreReg5} we estimate
\begin{equation*}
\begin{split}
X&\lesssim 2^{k_2^+-k_1}2^{3k_1/2}\|U_{>J_1,k_1}^{\mathcal{L}_1h}(t)\|_{L^2}\|U_{\leq J_2,k_2}^{\mathcal{L}_2\psi}(t)\|_{L^2}\lesssim \varep_1^2\langle t\rangle^{-2/3+170\delta}2^{k_1^-/3}2^{k_2^-/3}2^{-N(2)k_2^+}2^{-4k_1^+}.
\end{split}
\end{equation*}
As before, the desired bounds \eqref{box93.5} follow if $n_1=1$ and $n_2=2$ from these two estimates. 

Finally, if $n_1=2$ and $n_2=1$ then we use the $L^3$ bounds in \eqref{UOmegaMoreReg2} and the $L^6$ bounds in \eqref{UOmegaMoreReg3} to estimate
\begin{equation*}
\begin{split}
X&\lesssim 2^{k_2^+-k_1}\|U_{>J_1,k_1}^{\mathcal{L}_1h}(t)\|_{L^3}\|U_{\leq J_2,k_2}^{\mathcal{L}_2\psi}(t)\|_{L^6}\lesssim \varep_1^2\langle t\rangle^{-4/3+160\delta}2^{-k_1^-/3}2^{-k_2^-/3}2^{-N(3)k_2^+-2k_2^+}2^{-4k_1^+}.
\end{split}
\end{equation*}
Using the $L^2$ bounds in \eqref{UOmegaMoreReg6}--\eqref{UOmegaMoreReg5} we also estimate
\begin{equation*}
\begin{split}
X&\lesssim 2^{k_2^+-k_1}2^{3k_1/2}\|U_{>J_1,k_1}^{\mathcal{L}_1h}(t)\|_{L^2}\|U_{\leq J_2,k_2}^{\mathcal{L}_2\psi}(t)\|_{L^2}\lesssim \varep_1^2\langle t\rangle^{-1/3+160\delta}2^{2k_1^-/3}2^{k_2^-}2^{-N(2)k_2^+}2^{-4k_1^+}.
\end{split}
\end{equation*}
As before, the desired bounds \eqref{box93.5} follow if $n_1=2$ and $n_2=1$ from these two estimates.

{\bf{Case 4.2.}} Assume now that $q_1\geq 1$ and $q_2=0$. Recall that
\begin{equation}\label{box95}
\begin{split}
\|U_{>J_1,k_1}^{\mathcal{L}_1h}(t)\|_{L^2}&\lesssim\varep_12^{k_1/2}\min\{\langle t\rangle^{H(q_1,n_1)\delta},\langle t\rangle^{-1}2^{-k_1}\langle t\rangle^{H(q_1+1,n_1+1)\delta}\}2^{-5k_1^+},\\
\|U_{\leq J_2,k_2}^{\mathcal{L}_2\psi}(t)\|_{L^2}&\lesssim\varep_12^{k_2^-}\langle t\rangle^{H(1,n_2)\delta}2^{-N(n_2)k_2^+},
\end{split}
\end{equation}
as a consequence of \eqref{vcx1}--\eqref{cnb2}. Using just $L^2$ estimates we have
\begin{equation*}
\begin{split}
X&\lesssim 2^{k_2^++k_1/2}\|U_{>J_1,k_1}^{\mathcal{L}_1h}(t)\|_{L^2}\|U_{\leq J_2,k_2}^{\mathcal{L}_2\psi}(t)\|_{L^2}\\
&\lesssim \varep_1^22^{k_1^-+k_2^-}\langle t\rangle^{\delta(H(q_1,n_1)+H(1,n_2))}2^{-N(n_2)k_2^++2k_2^+}2^{-4k_1^+}.
\end{split}
\end{equation*}
Since $H(1,n_2)+H(q_1,n_1)\leq H(q_1,n_1+n_2)-150$ (due to the last inequality in \eqref{SuperlinearH2} and the assumption $q_1\geq 1$), this suffices to prove \eqref{box93.5} when $2^{k_1^-+k_2^-}\lesssim \langle t\rangle^{-1+180\delta}$. 

Finally, assume that $2^{k_1^-+k_2^-}\geq \langle t\rangle^{-1+180\delta}$. If $n_2=1$ then we estimate, using \eqref{box95}, Sobolev embedding, and the $L^6$ bounds in \eqref{UOmegaMoreReg3},
\begin{equation*}
\begin{split}
X&\lesssim 2^{k_2^+-k_1}\|U_{>J_1,k_1}^{\mathcal{L}_1h}(t)\|_{L^3}\|U_{\leq J_2,k_2}^{\mathcal{L}_2\psi}(t)\|_{L^6}\\
&\lesssim \varep_1^2\min\{\langle t\rangle^{H(q_1,n_1)\delta},\langle t\rangle^{-1}2^{-k_1}\langle t\rangle^{H(q_1+1,n_1+1)\delta}\}2^{-5k_1^+}\cdot \langle t\rangle^{-1+50\delta}2^{-k_2^-/3}2^{-N(2)k_2^+-k_2^+}\\
&\lesssim \varep_1^22^{-k_2^-/3}\langle t\rangle^{-1/3}2^{-k_1^-/3}\langle t\rangle^{\frac{H(q_1+1,n_1+1)+2H(q_1,n_1)}{3}\delta}2^{-4k_1^+}\cdot \langle t\rangle^{-1+50\delta}2^{-N(2)k_2^+-k_2^+}.
\end{split}
\end{equation*}
Since $H(q_1+1,n_1+1)+2H(q_1,n_1)\leq 3H(q_1,n_1+1)$ (see \eqref{fvc1.0}), this suffices to prove \eqref{box93.5} if $n_2=1$. 

On the other hand, if $n_2=2$ (so we necessarily have $(q_1,n_1)=(1,1)$) we estimate
\begin{equation*}
\begin{split}
X&\lesssim 2^{k_2^+-k_1}\|U_{>J_1,k_1}^{\mathcal{L}_1h}(t)\|_{L^6}\|U_{\leq J_2,k_2}^{\mathcal{L}_2\psi}(t)\|_{L^3}\\
&\lesssim \varep_1^22^{k_1/2}\min\{\langle t\rangle^{H(q_1,n_1)\delta},\langle t\rangle^{-1}2^{-k_1}\langle t\rangle^{H(q_1+1,n_1+1)\delta}\}2^{-5k_1^+}\cdot \langle t\rangle^{-1/2+100\delta}2^{-k_2^-/6}2^{-N(3)k_2^+-k_2^+}\\
&\lesssim \varep_1^22^{-k_2^-/6}\langle t\rangle^{-2/3}2^{-k_1^-/6}\langle t\rangle^{\frac{2H(q_1+1,n_1+1)+H(q_1,n_1)}{3}\delta}2^{-4k_1^+}\cdot \langle t\rangle^{-1/2+100\delta}2^{-N(3)k_2^+-k_2^+},
\end{split}
\end{equation*}
using \eqref{box95} and the $L^3$ bounds in \eqref{UOmegaMoreReg3}. Since $[2H(2,2)+H(1,1)]/3\leq H(1,3)-100$ (see \eqref{fvc1.0}), this suffices to prove \eqref{box93.5} if $n_2=2$. This completes the proof of the lemma.
\end{proof}

\subsection{The classes of functions $\mathcal{G}_a$} In most cases, the cubic and higher order nonlinearities can be treated perturbatively, and do not play a significant role in the analysis. To justify this, we need good bounds on the quadratic metric components $g_{\geq 2}^{\al\be}$.

The metric components $g_{\geq 2}^{\al\be}$ satisfy the identities \eqref{zaq22}. Therefore they can be represented as infinite sums of monomials of degree $\geq 2$ in the functions $h_{\mu\nu}$. More generally, for integers $a\geq 1$ we define the sets
\begin{equation}\label{gb0.1}
\mathcal{G}_a:=\big\{G_a=h_1\cdot\ldots\cdot h_a:\,h_1,\ldots,h_a\in\{h_{\mu\nu}\}\big\}.
\end{equation}
By convention, set $\mathcal{G}_0:=\{1\}$. In this subsection we prove the following bounds:

\begin{lemma}\label{Ga} Assume $\mathcal{L}\in \mathcal{V}_n^q$, $t\in[0,T]$, $l\in\{1,2,3\}$, and $G_a\in\mathcal{G}_a$, $a\geq 2$. Then there is a constant $C_0\geq 1$  such that
\begin{equation}\label{gb5}
\begin{split}
\|P_k\mathcal{L}G_a(t)\|_{L^2}&\leq (C_0\varep_1)^a\langle t\rangle^{-1+\delta'}2^{-N(n)k^+}2^{-k/2}\min\{1,2^{k^-}\langle t\rangle\}^{1-\delta'},\\
\|P_k\mathcal{L} G_a(t)\|_{L^\infty}&\leq (C_0\varep_1)^a \langle t\rangle^{-2+\delta'}2^{-N(n+1)k^++3k^+}\min\{1,2^{k^-}\langle t\rangle\}^{2-\delta'},\\
\|P_k(x_l\mathcal{L}G_a)(t)\|_{L^2}&\leq (C_0\varep_1)^a(\langle t\rangle+2^{-k^-})^{\delta'}2^{-N(n+1)k^+}2^{-k^-/2},
\end{split}
\end{equation}
where and $\delta'=2000\delta$ (see \eqref{parameters}), and the inequality in the first line holds for all pairs $(q,n)$ with $n\leq 3$, while the inequalities in the last two lines hold for pairs $(q,n)$ with $n\leq 2$. Thus
\begin{equation}\label{gb1}
\begin{split}
\|P_k\mathcal{L}g_{\geq 2}^{\al\be}(t)\|_{L^2}&\lesssim \varep_1^2\langle t\rangle^{-1+\delta'}2^{-N(n)k^+}2^{-k/2}\min\{1,2^{k^-}\langle t\rangle\}^{1-\delta'},\\
\|P_k\mathcal{L}g_{\geq 2}^{\al\be}(t)\|_{L^\infty}&\lesssim \varep_1^2 \langle t\rangle^{-2+\delta'}2^{-N(n+1)k^++3k^+}\min\{1,2^{k^-}\langle t\rangle\}^{2-\delta'},\\
\|P_k(x_l\mathcal{L}g_{\geq 2}^{\al\be})(t)\|_{L^2}&\lesssim \varep_1^2(\langle t\rangle+2^{-k^-})^{\delta'}2^{-N(n+1)k^+}2^{-k^-/2}.
\end{split}
\end{equation}
\end{lemma}

\begin{remark} We notice that we prove more than just frequency-localized $L^2$ bounds on the functions $\mathcal{L}g^{\al\be}_{\geq 2}$. In particular, we prove weighted  $L^2$ bounds which are important in our bootstrap scheme, see the key estimates \eqref{wer3.1} and \eqref{wer3.3} in Proposition \ref{DiEs1}.
\end{remark}

\begin{proof} The bounds \eqref{vcx1}, \eqref{wws1}, and \eqref{bootstrap2.2} show that, for any $\mathcal{L}\in\mathcal{V}_n^q$ and $k\in\mathbb{Z}$,
\begin{equation}\label{gb3}
\begin{split}
\|P_k\mathcal{L}h_{\al\be}(t)\|_{L^2}&\lesssim \varep_1\langle t\rangle^{H(q,n)\delta}2^{-N(n)k^+}2^{-k/2}(\langle t\rangle 2^{k^-})^{-\ga},\\
\|P_k\mathcal{L} h_{\al\be}(t)\|_{L^\infty}&\lesssim \varep_1\langle t\rangle^{-1+H(q+1,n+1)\delta}2^{-N(n+1)k^++3k^+}\min\{1,2^{k^-}\langle t\rangle\}^{1-\delta},\\
\|P_k(x_l\mathcal{L}h_{\al\be})(t)\|_{L^2}&\lesssim \varep_1\langle t\rangle^{H(q+1,n+1)\delta}2^{-N(n+1)k^+}2^{-k^-/2}(2^{-k^-}+\langle t\rangle)(\langle t\rangle 2^{k^-})^{-\ga},
\end{split}
\end{equation}
where the first inequality holds for all pairs $(q,n)\leq (3,3)$, while the last two inequalities hold for pairs $(q,n)\leq (2,2)$. Indeed, for the last bound we estimate first
\begin{equation*}
\begin{split}
\|P_k(x_l\mathcal{L}h_{\al\be})(t)\|_{L^2}\lesssim 2^{-2k^-}\| P_{[k-4,k+4]}U^{\mathcal{L}h_{\al\be}}(t)\|_{L^2}+2^{-k^-}\|\varphi_k(\xi)\partial_{\xi_l}\widehat{U^{\mathcal{L}h_{\al\be}}}(\xi,t)\|_{L^2_\xi}.
\end{split}
\end{equation*}
Then we recall that $\widehat{U^{\mathcal{L}h_{\al\be}}}(\xi,t)=e^{-it|\xi|}\widehat{V^{\mathcal{L}h_{\al\be}}}(\xi,t)$. Therefore, the right-hand side of the inequality above is bounded by
\begin{equation*}
C2^{-k^-}(2^{-k^-}+\langle t\rangle)\| P_{[k-4,k+4]}U^{\mathcal{L}h_{\al\be}}(t)\|_{L^2}+C2^{-k^-}\|\varphi_k(\xi)\partial_{\xi_l}\widehat{V^{\mathcal{L}h_{\al\be}}}(\xi,t)\|_{L^2_\xi}.
\end{equation*}
The estimate in the last line of \eqref{gb3} follows from \eqref{vcx1} and \eqref{cnb1}.

{\bf{Step 1.}} We consider first the case $a=2$. Assume $\mathcal{L}_1\in\mathcal{V}_{n_1}^{q_1}$, $\mathcal{L}_2\in\mathcal{V}_{n_2}^{q_2}$, $(q_1,n_1)+(q_2,n_2)\leq (q,n)$. Assume also $h_1,h_2\in\{h_{\mu\nu}\}$. If $2^{k^-}\leq\langle t\rangle^{-1}$ then we bound, using \eqref{gb3} and \eqref{abc36.6},
\begin{equation}\label{gb7.1}
\begin{split}
\|P_k(\mathcal{L}_1h_1\cdot\mathcal{L}_2h_2)(t)\|_{L^2}&\lesssim \sum_{(k_1,k_2)\in\mathcal{X}_k}2^{3\min(k,k_1,k_2)/2}\|P_{k_1}\mathcal{L}_1h_1(t)\|_{L^2}\|P_{k_2}\mathcal{L}_2h_2(t)\|_{L^2}\\
&\lesssim \varep_1^2\langle t\rangle ^{\delta'}2^{k^-/2},
\end{split}
\end{equation}
where the sets $\mathcal{X}_k$ are as in \eqref{gb7}. Moreover, if $n_1,n_2\leq 2$ and $2^{k^-}\geq\langle t\rangle^{-1}$ then
\begin{equation}\label{gb7.4}
\begin{split}
\|P_k(\mathcal{L}_1h_1\cdot\mathcal{L}_2h_2)(t)\|_{L^\infty}&\lesssim \sum_{(k_1,k_2)\in\mathcal{X}_k}\|P_{k_1}\mathcal{L}_1h_1(t)\|_{L^\infty}\|P_{k_2}\mathcal{L}_2h_2(t)\|_{L^\infty}\\
&\lesssim \varep_1^2\langle t\rangle ^{-2+\delta'}2^{-N(n+1)k^++3k^+}.
\end{split}
\end{equation}

To estimate $\|P_k(x_l\mathcal{L}_1h_1\cdot\mathcal{L}_2h_2)(t)\|_{L^2}$ we combine the factor $x_l$ with the higher frequency term and estimate it in $L^2$, and estimate the lower frequency term in $L^\infty$. Using \eqref{gb3} as before, it follows that, for $l\in\{1,2,3\}$,
\begin{equation}\label{gb7.5}
\|P_k(x_l\mathcal{L}_1h_1\cdot\mathcal{L}_2h_2)(t)\|_{L^2}\lesssim \varep_1^2\langle t\rangle^{-1+\delta'}2^{-N(n+1)k^+}2^{-k^-/2}\min(1,2^{k^-}\langle t\rangle)^{1-\delta'}\cdot (2^{-k^-}+\langle t\rangle).
\end{equation} 

The $L^\infty$ bounds in \eqref{gb5} follow from \eqref{gb7.4} if $2^{k^-}\langle t\rangle\gtrsim 1$ and from \eqref{gb7.1} if  $2^{k^-}\langle t\rangle\lesssim 1$. The weighted $L^2$ bounds in the last line of \eqref{gb5} follow from \eqref{gb7.5}. The $L^2$ bounds in the first line of \eqref{gb5} follow from \eqref{gb7.1} if $2^{k^-}\langle t\rangle\lesssim 1$. It remains to prove that
\begin{equation}\label{gb8.2}
\|P_{k}(\mathcal{L}_1h_1\cdot \mathcal{L}_2h_2)(t)\|_{L^2}\lesssim\varep_1^2\langle t\rangle^{-1+\delta'}2^{-k/2}2^{-N(n)k^+},
\end{equation}
for any $k\in\mathbb{Z}$ with $2^k\geq \langle t\rangle^{-1}$, and any $\mathcal{L}_1\in\mathcal{V}_{n_1}^{q_1}$, $\mathcal{L}_2\in\mathcal{V}_{n_2}^{q_2}$ with $(q_1,n_1)+(q_2,n_2)\leq (q,n)$.

To prove \eqref{gb8.2}, we assume first that $n_1,n_2\leq 2$ and estimate
\begin{equation}\label{gb8.3}
\|P_{k}(\mathcal{L}_1h_1\cdot \mathcal{L}_2h_2)(t)\|_{L^2}\leq S_1+S_2+S_3,
\end{equation}
where, using \eqref{wws1} and the inequality in the first line of \eqref{gb3},
\begin{equation*}
S_1:=\sum_{k_1\leq k-6,\,|k_2-k|\leq 4}\|P_{k_1}\mathcal{L}_1h_1(t)\|_{L^\infty}\|P_{k_2}\mathcal{L}_2h_2(t)\|_{L^2}\lesssim \varep_1^2\langle t\rangle^{-1+\delta'/2}2^{-k/2}2^{-N(n_2)k^+},
\end{equation*}
\begin{equation*}
S_2:=\sum_{k_2\leq k-6,\,|k_1-k|\leq 4}\|P_{k_1}\mathcal{L}_1h_1(t)\|_{L^2}\|P_{k_2}\mathcal{L}_2h_2(t)\|_{L^\infty}\lesssim \varep_1^2\langle t\rangle^{-1+\delta'/2}2^{-k/2}2^{-N(n_1)k^+},
\end{equation*}
and
\begin{equation*}
\begin{split}
S_3:=\sum_{k_1,k_2\geq k-6,\,|k_1-k_2|\leq 10}&\|P_{k_1}\mathcal{L}_1h_1(t)\|_{L^\infty}\|P_{k_2}\mathcal{L}_2h_2(t)\|_{L^2}\lesssim \varep_1^2\langle t\rangle^{-1+\delta'/2}2^{-k/2}2^{-N(n_2)k^+}.
\end{split}
\end{equation*}
These bounds clearly suffice to prove \eqref{gb8.2}.

On the other hand if $\max(n_1,n_2)=3$, then we may assume that $n_2=3$ and $n_1=0$. The bounds on $S_1$ and $S_3$ above still hold, but the bounds on $S_2$ fail, because we do not have suitable $L^\infty$ bounds on $P_{k_2}\mathcal{L}_2h_2(t)$. However, we can still use $L^2$ bounds as in \eqref{gb7.1} to control the contribution of small frequencies $k_2$, i.e. $2^{k_2}\lesssim \langle t\rangle^{-1}$. For \eqref{gb8.2} it remains to show that
\begin{equation}\label{gb8.4}
\|P_{k}(P_{k_1}h_1\cdot P_{k_2}\mathcal{L}_2h_2)(t)\|_{L^2}\lesssim\varep_1^2\langle t\rangle^{-1+2\delta'/3}2^{-k/2}2^{-N(3)k^+-k^+},
\end{equation}
provided that $|k_1-k|\leq 4$, $\langle t \rangle^{-1}\leq 2^{k_2}$, $k_2\leq k-6$, and $2^k\leq\langle t\rangle ^{1/10}$.

To prove \eqref{gb8.4} we would like to use Lemma \ref{Lembil2} (i) (the simple idea of directly estimating $P_{k_1}h_1$ and $L^\infty$ and $P_{k_2}\mathcal{L}_2h_2$ in $L^2$ does not work when $k_2$ is small, due to the factor $2^{-k_2/2}$ coming from \eqref{gb3}). For this we write first
\begin{equation*}
h_1(t)=-|\nabla|^{-1}\Im(U^{h_1}(t))=-|\nabla|^{-1}\Im(e^{-it\Lambda_{wa}}V^{h_1}(t)),
\end{equation*}
see \eqref{variables4L}). Then we decompose $P_{k_1}V^{h_1,\iota_1}=V^{h_1,\iota_1}_{\leq J_1,k_1}+V^{h_1,\iota_1}_{>J_1,k_1}$, where $J_1$ is the largest integer such that $2^{J_1}\leq\langle t\rangle(1+2^{k_1}\langle t\rangle)^{-\delta}$ and $\iota_1\in\{+,-\}$. With $I$ defined as in \eqref{mults}--\eqref{abc36.1}, we estimate
\begin{equation*}
\begin{split}
\|P_{k}I[e^{-it\Lambda_{wa,\iota_1}}V^{h_1,\iota_1}_{\leq J_1,k_1}(t),P_{k_2}\mathcal{L}_2h_2)(t)]\|_{L^2}&\lesssim 2^{k_2/2}\langle t\rangle^{-1}(2^{k_1}\langle t\rangle)^{\delta}\|P_{k_2}\mathcal{L}_2h_2(t)\|_{L^2}\|P_{k_1}V^{h_1}\|_{H^{0,1}_\Omega}\\
&\lesssim\varep_1^2\langle t\rangle^{-1+\delta'/2}2^{k_1/2}2^{-N(1)k_1^++2k_1^+},
\end{split}
\end{equation*}
for $\iota_1\in\{+,-\}$, using \eqref{Bil31}, \eqref{gb3}, and \eqref{vcx1}. We also estimate
\begin{equation*}
\begin{split}
\|P_{k}I[e^{-it\Lambda_{wa,\iota_1}}V^{h_1,\iota_1}_{> J_1,k_1}(t),P_{k_2}\mathcal{L}_2h_2)(t)]\|_{L^2}&\lesssim 2^{3k_2/2}\|P_{k_2}\mathcal{L}_2h_2(t)\|_{L^2}\|Q_{>J_1,k_1}V^{h_1}(t)\|_{L^2}\\
&\lesssim\varep_1^2\langle t\rangle^{-1+\delta'/2}2^{k_1/2}2^{-N(1)k_1^++2k_1^+},
\end{split}
\end{equation*}
using \eqref{gb3} and \eqref{vcx1.1}. The estimates \eqref{gb8.4} follow from these two bounds, which completes the proof of \eqref{gb8.2}.

{\bf{Step 2.}} We prove now the bounds \eqref{gb5} for $a\geq 3$, by induction over $a$. The last two bounds in \eqref{gb5} follow as in {\bf{Step 1}}, since the bounds satisfied by $(C\varep_1)^{-a}G_a$ in \eqref{gb5} are stronger than the bounds satisfied by $\varep_1^{-1}h_{\al\be}$ in \eqref{gb3}. As before, it remains to prove that if $h\in\{h_{\al\be}\}$ and $G\in\mathcal{G}_a$ then
\begin{equation}\label{gb16.4}
\|P_{k}(\mathcal{L}_1h\cdot \mathcal{L}_2G)(t)\|_{L^2}\lesssim (C_0\varep_1)^a\varep_1\langle t\rangle^{-1+\delta'/2}2^{-k/2}2^{-N(n)k^+},
\end{equation}
for any $k\in\mathbb{Z}$ with $2^k\geq \langle t\rangle^{-1}$, and any $\mathcal{L}_1\in\mathcal{V}_{n_1}^{q_1}$, $\mathcal{L}_2\in\mathcal{V}_{n_2}^{q_2}$ with $(q_1,n_1)+(q_2,n_2)\leq (q,n)$.

The bounds \eqref{gb16.4} follow in the same way as the bounds \eqref{gb8.2} if $\max\{n_1,n_2\}\leq 2$, using $L^\infty\times L^2$ estimates with the lower frequency measured in $L^\infty$. On the other hand, if $(n_1,n_2)=(0,3)$, then we have to prove the analogue of the bounds \eqref{gb8.4}, which is  
\begin{equation*}
\|P_{k}(P_{k_1}h\cdot P_{k_2}\mathcal{L}_2G)(t)\|_{L^2}\lesssim\varep_1(C_0\varep_1)^a\langle t\rangle^{-1}2^{-k/2}2^{-N(3)k^+-k^+},
\end{equation*}
provided that $|k_1-k|\leq 4$, $\langle t \rangle^{-1}\leq 2^{k_2}$, $k_2\leq k-6$, and $2^k\leq\langle t\rangle ^{1/10}$. This is easier now, since we can just use the $L^2$ estimates $\|P_{k_2}\mathcal{L}_2G(t)\|_{L^2}\lesssim (C_0\varep_1)^a2^{-k_2/2}\langle t\rangle ^{-1+\delta'}$ from the induction hypothesis, and combine them with $L^\infty$ estimates on $P_{k_1}h$. The loss of the factor $2^{-k_2/2}$ is mitigated in this case by the gain of time decay. The proof in the case $(n_1,n_2)=(3,0)$ is similar, which completes the proof of the lemma.
\end{proof}

In some of the analysis in sections \ref{NullStruc} and \ref{DecompositionNonlin} we need a slightly different bound:

\begin{lemma}\label{Ga4}
Assume $\mathcal{L}\in \mathcal{V}_n^q$, $t\in[0,T]$, $k\in\mathbb{Z}$, and $G_{\geq 1}=\sum_{d\geq 1}a_dg_d$ for some functions $g_d\in\mathcal{G}_d$ (see \eqref{gb0.1}) and some coefficients $a_d\in\mathbb{R}$ with $|a_d|\leq C^d$. Then
\begin{equation}\label{nyu1}
\begin{split}
\|P_k|\nabla|^{-1}\mathcal{L}(G_{\geq 1}\partial_\rho h)(t)\|_{L^2}&\lesssim\varep_1^2\langle t\rangle^{-1+\delta'}2^{-N(n)k^+}2^{-k/2}\min\{1,2^{k^-}\langle t\rangle\}^{1-\delta'},\\
\|P_k|\nabla|^{-1}\mathcal{L}(G_{\geq 1}\partial_\rho h)(t)\|_{L^\infty}&\lesssim\varep_1^2\langle t\rangle^{-2+\delta'}2^{-N(n+1)k^++3k^+}\min\{1,2^{k^-}\langle t\rangle\}^{2-\delta'},
\end{split}
\end{equation}
where $h\in\{h_{\al\be}\}$, $\rho\in\{0,1,2,3\}$, and the bounds in the second line hold only if $n\leq 2$. 
\end{lemma}

\begin{proof} Notice that the bounds \eqref{nyu1} are slightly stronger than the $L^2$ and the $L^\infty$ bounds in \eqref{gb5}, but we do not prove weighted $L^2$ bounds. $\mathrm{Low}\times\mathrm{High}\to\mathrm{High}$ interactions can still be estimated in the same way, but some care is needed to estimate $\mathrm{High}\times\mathrm{High}\to\mathrm{Low}$ interactions, due to the factor $|\nabla|^{-1}$. More precisely, we prove that for any $k\in\mathbb{Z}$, and $\LL_1\in\mathcal{V}_{n_1}^{q_1}$, $\LL_2\in\mathcal{V}_{n_2}^{q_2}$ with $(q_1,n_1)+(q_2,n_2)\leq (q,n)$, we have
\begin{equation}\label{nyu2}
\begin{split}
\sum_{k_1,k_2\geq k}\|P_k(P_{k_1}\LL_1 G_{\geq 1}\cdot P_{k_2}\partial_\rho\LL_2h)(t)\|_{L^2}&\lesssim\varep_1^2\langle t\rangle^{-1+\delta'}2^{-N(n)k^+}2^{k/2}\min\{1,2^{k^-}\langle t\rangle\}^{1-\delta'},\\
\sum_{k_1,k_2\geq k}\|P_k(P_{k_1}\LL_1 G_{\geq 1}\cdot P_{k_2}\partial_\rho\LL_2h)(t)\|_{L^\infty}&\lesssim\varep_1^2\langle t\rangle^{-2+\delta'}2^{-N(n+1)k^++3k^+}2^k\min\{1,2^{k^-}\langle t\rangle\}^{2-\delta'}.
\end{split}
\end{equation}

These bounds follow easily when $2^k\langle t\rangle^{1-2\delta}\lesssim 1$ using just $L^2$ estimates as in \eqref{gb7.1}. They also follow easily when $k\geq 0$, since there is no potential derivative loss in this case. Finally, if $2^k\in[2^{10}\langle t\rangle^{2\delta-1},1]$ then the contribution of the components $a_dG_d$, $d\geq 2$, in $G_{\geq 1}$ can be bounded in the same way, using \eqref{gb5}. After these reductions, for \eqref{nyu2} it remains to show that
 \begin{equation}\label{nyu3}
\begin{split}
\sum_{k_1,k_2\geq k}2^{-k_1}\|P_k I[P_{k_1}U^{\LL_1 h_1,\iota_1},P_{k_2}U^{\LL_2h_2,\iota_2}](t)\|_{L^2}&\lesssim\varep_1^2\langle t\rangle^{-1+\delta'}2^{k/2},\\
\sum_{k_1,k_2\geq k}2^{-k_1}\|P_k I[P_{k_1}U^{\LL_1 h_1,\iota_1},P_{k_2}U^{\LL_2h_2,\iota_2}](t)\|_{L^\infty}&\lesssim\varep_1^2\langle t\rangle^{-2+\delta'}2^k,
\end{split}
\end{equation}
provided that $h_1,h_2\in\{h_{\al\be}\}$, $\iota_1,\iota_2\in\{+,-\}$, $I$ is as in \eqref{abc36.1}, and $2^k\in[2^{10}\langle t\rangle^{2\delta-1},1]$.

To prove \eqref{nyu3} we use the super-localized bounds \eqref{Linfty1.6*}. By symmetry, we may assume $n_1=\min(n_1,n_2)\leq 1$ and decompose $P_{k_1}U^{\LL_1 h_1,\iota_1}(t)=U^{\LL_1 h_1,\iota_1}_{\leq J_1,k_1}(t)+U_{>J_1,k_1}^{\LL_1 h_1,\iota_1}(t)$ as in \eqref{box15}, where $J_1$ is the largest integer such that $2^{J_1}\leq \langle t\rangle (1+2^{k_1}\langle t\rangle)^{-\delta}$. The $L^2$ bounds in \eqref{nyu3} follow using \eqref{Bil31} to estimate the contribution of $U^{\LL_1 h_1,\iota_1}_{\leq J_1,k_1}(t)$, and using \eqref{abc36.6} and \eqref{vcx1.1} to bound the contribution of $U^{\LL_1 h_1,\iota_1}_{>J_1,k_1}(t)$. For the $L^\infty$ bounds, it suffices to show that 
\begin{equation}\label{plk10.3}
\|P_kI[P_{k_1}U^{\mathcal{L}_1h_1,\iota_1},P_{k_2}U^{\mathcal{L}_2h_2,\iota_2}](t)\|_{L^\infty}\lesssim\varep_1^22^{k_1}2^{k}\langle t\rangle^{-2+3\delta'/4}2^{-2k_1^+}2^{-2k_2^+},
\end{equation}
provided that $|k_1-k_2|\leq 4$, $2^k\in[2^{10}\langle t\rangle^{-1},1]$, and $2^{k_1},2^{k_2}\in [2^k,\langle t\rangle]$. 

We decompose also $P_{k_2}U^{\mathcal{L}_2h_2,\iota_2}=U^{\mathcal{L}_2h_2,\iota_2}_{\leq J_2,k_2}+U^{\mathcal{L}_2h_2,\iota_2}_{> J_2,k_2}$, where $2^{J_1}=2^{J_2}$, and estimate
\begin{equation*}
\begin{split}
\big\|P_kI[U^{\mathcal{L}_1h_1,\iota_1}_{> J_1,k_1},&P_{k_2}U^{\mathcal{L}_2h_2,\iota_2}](t)\big\|_{L^\infty}+\big\|P_kI[U^{\mathcal{L}_1h_1,\iota_1}_{\leq J_1,k_1},U^{\mathcal{L}_2h_2,\iota_2}_{>J_2,k_2}](t)\big\|_{L^\infty}\\
&\lesssim\varep_1^22^{3k/2}2^{k_1/2}\langle t\rangle^{-2+\delta'/2}2^{-4k_1^+}2^{-4k_2^+},
\end{split}
\end{equation*}
using \eqref{vcx1.1} and \eqref{wws1}. Finally, using \eqref{Linfty1.6*} with $l=k$ we estimate
\begin{equation*}
\begin{split}
\big\|P_kI[U^{\mathcal{L}_1h_1,\iota_1}_{\leq J_1,k_1},U^{\mathcal{L}_2h_2,\iota_2}_{\leq J_2,k_2}](t)\big\|_{L^\infty}&\lesssim\sum_{|n_1-n_2|\leq 4}\big\|\mathcal{C}_{n_1,l}U^{\mathcal{L}_1h_1,\iota_1}_{\leq J_1,k_1}\big\|_{L^\infty}\big\|\mathcal{C}_{n_2,l}U^{\mathcal{L}_2h_2,\iota_2}_{\leq J_2,k_2}(t)\big\|_{L^\infty}\\
&\lesssim \varep_1^2\langle t\rangle^{-2+\delta'/2}2^k2^{k_1/2}2^{k_2/2}2^{-4k_1^+}2^{-4k_2^+},
\end{split}
\end{equation*}
using also \eqref{vcx1} in the estimate in the second line. This completes the proof of \eqref{plk10.3}.
\end{proof}

\section{Bounds on the nonlinearities $\mathcal{N}_{\al\be}^{h}$ and $\mathcal{N}^\psi$} \label{NullStruc}

Recall the decompositions, see \eqref{sac1}--\eqref{npsi2},
\begin{equation}\label{bom1}
\begin{split}
\mathcal{N}_{\al\be}^{h}&=\mathcal{KG}^2_{\al\be}+\mathcal{Q}^2_{\al\be}+\mathcal{S}^2_{\al\be}+\mathcal{N}_{\al\be}^{h,\geq 3},\\
\mathcal{N}^{\psi}&=\mathcal{N}^{\psi,2}+\mathcal{N}^{\psi,\geq 3}.
\end{split}
\end{equation}
In this subsection we prove several frequency-localized bounds on the nonlinear terms $\mathcal{L}\mathcal{KG}^2_{\al\be}$, $\mathcal{L}\mathcal{Q}^2_{\al\be}$, $\mathcal{L}\mathcal{S}^2_{\al\be}$, $\mathcal{L}\mathcal{N}_{\al\be}^{h,\geq 3}$, $\mathcal{L}\mathcal{N}^{\psi,2}$, and $\mathcal{L}\mathcal{N}^{\psi,\geq 3}$. One should think of these as rather general bounds; we could improve some of them slightly, both in terms of differentiability and time decay, but we do not pursue all the possible improvements at this stage. 

For $n\geq 0$ we define $\widetilde{\ell}(n)$ (a slightly worse loss function) and $\widetilde{N}(n)$ by 
\begin{equation}\label{plk2.5}
\widetilde{\ell}(n):=35\text{ if }n\geq 1,\quad\widetilde{\ell}(0):=25,\qquad \widetilde{N}(n):=N(n)\text{ if }n\geq 1,\quad\widetilde{N}(0):=N_0.
\end{equation}
The following proposition is our main result in this section.

\begin{proposition}\label{plk1}
For any $k\in\mathbb{Z}$, $t\in[0,T]$, $\mathcal{L}\in\mathcal{V}^q_n$, $q\leq n\leq 3$, and $\al,\be\in\{0,1,2,3\}$
\begin{equation}\label{wer4.0}
\|P_k(\LL\mathcal{N}_{\al\be}^{h})(t)\|_{L^2}\lesssim\varep_1^22^{k/2}\langle t\rangle^{H(q,n)\delta}2^{-\widetilde{N}(n)k^++7k^+}\min(2^{k}\langle t\rangle^{3\delta/2},\langle t\rangle^{-1+\widetilde{\ell}(n)\delta}),
\end{equation}
and
\begin{equation}\label{wer4.1}
\|P_k(\LL\mathcal{N}^{\psi})(t)\|_{L^2}\lesssim\varep_1^2\langle t\rangle^{H(q,n)\delta}2^{-\widetilde{N}(n)k^++7k^+}\min(2^{k}\langle t\rangle^{3\delta/2},\langle t\rangle^{-1+\widetilde{\ell}(n)\delta}).
\end{equation}
Moreover if $n\leq 2$ and $l\in\{1,2,3\}$ then
\begin{equation}\label{plk2}
\begin{split}
\|P_k(\LL\mathcal{N}_{\al\be}^{h})(t)\|_{L^\infty}&\lesssim\varep_1^22^{k}\langle t\rangle^{-1+4\delta'}2^{-N(n+1)k^++5k^+}\min(2^{k},\langle t\rangle^{-1}),\\
\|P_k(x_l\LL\mathcal{N}_{\al\be}^{h})(t)\|_{L^2}&\lesssim\varep_1^22^{k/2}\langle t\rangle^{H(q,n)\delta+\widetilde{\ell}(n)\delta}2^{-N(n+1)k^+-2k^+},
\end{split}
\end{equation}
and
\begin{equation}\label{plk3}
\begin{split}
\|P_k(\LL\mathcal{N}^{\psi})(t)\|_{L^\infty}&\lesssim\varep_1^22^{k/2}\langle t\rangle^{-1+4\delta'}2^{-N(n+1)k^++5k^+}\min(2^{k},\langle t\rangle^{-1}),\\
\|P_k(x_l\LL\mathcal{N}^{\psi})(t)\|_{L^2}&\lesssim\varep_1^2\langle t\rangle^{H(q,n)\delta+\widetilde{\ell}(n)\delta}2^{-N(n+1)k^+-2k^+}.
\end{split}
\end{equation}
\end{proposition}

The proposition follows from Lemmas \ref{abc2.2}--\ref{dtv8.0} below.

\subsection{The quadratic nonlinearities} We consider first the nonlinearities $\mathcal{KG}^2_{\al\be}$.

\begin{lemma}\label{abc2.2}
Assume that $t\in[0,T]$, $\mathcal{L}\in\mathcal{V}_n^q$, $q\leq n\leq 3$, and $k\in\mathbb{Z}$. Then
\begin{equation}\label{gb17.1}
\|P_k\{\mathcal{L}\mathcal{KG}^2_{\al\be}\}(t)\|_{L^2}\lesssim \varep_1^22^{k/2}2^{-N(n)k^+}\langle t\rangle^{H(q,n)\delta }\min(2^{k}\langle t\rangle^{3\delta/2},\langle t\rangle^{-1+\ell(q,n)\delta+\delta/2}).
\end{equation}
Moreover, if $l\in\{1,2,3\}$ and $n\leq 2$ then we also have the bounds
\begin{equation}\label{gb17.2}
\begin{split}
\|P_k\{\mathcal{L}\mathcal{KG}^2_{\al\be}\}(t)\|_{L^\infty}&\lesssim\varep_1^22^{k}\langle t\rangle^{-1+2\delta'}2^{-N(n+1)k^++3k^+}\min\{2^{k},\langle t\rangle^{-1}\},\\
\|P_k\{x_l\mathcal{L}\mathcal{KG}^2_{\al\be}\}(t)\|_{L^2}&\lesssim \varep_1^22^{k/2}2^{-N(n+1)k^+-5k^+}\langle t\rangle^{\delta[H(q,n)+\widetilde{\ell}(n)]}.
\end{split}
\end{equation}
\end{lemma}

\begin{proof} {\bf{Step 1.}} Recall the bilinear operators $I$ defined in \eqref{abc36.1} and the loss function $\ell$ defined in \eqref{AcceptableLoss}. In view of \eqref{on5}, for \eqref{gb17.1} it suffices to prove that
\begin{equation}\label{sac32}
\begin{split}
\sum_{(k_1,k_2)\in\mathcal{X}_k}&\|P_kI[P_{k_1}U^{\mathcal{L}_1\psi,\iota_1},P_{k_2}U^{\mathcal{L}_2\psi,\iota_2}](t)\|_{L^2}\\
&\lesssim \varep_1^22^{k/2}2^{-N(n)k^+}\langle t\rangle^{H(q,n)\delta }\min(2^{k}\langle t\rangle^{3\delta/2},\langle t\rangle^{-1+\ell(q,n)\delta+\delta/2})
\end{split}
\end{equation}
for any $\iota_1,\iota_2\in\{+,-\}$, $\mathcal{L}_1\in\mathcal{V}_{n_1}^{q_1}$, $\mathcal{L}_2\in\mathcal{V}_{n_2}^{q_2}$, $(q_1,n_1)+(q_2,n_2)\leq (q,n)$. This follows easily from the bounds \eqref{box32}--\eqref{box33}.

{\bf{Step 2.}} To prove the $L^\infty$ bounds in \eqref{gb17.2} it suffices to show that
\begin{equation}\label{plk15}
\begin{split}
\sum_{(k_1,k_2)\in\mathcal{X}_k}&\|P_kI[P_{k_1}U^{\mathcal{L}_1\psi,\iota_1},P_{k_2}U^{\mathcal{L}_2\psi,\iota_2}](t)\|_{L^\infty}\\
&\lesssim \varep_1^22^{k}\langle t\rangle^{-1+2\delta'}2^{-N(n+1)k^++3k^+}\min(2^{k},\langle t\rangle^{-1}).
\end{split}
\end{equation}
This follows from \eqref{sac32} if $2^k\leq\langle t\rangle^{-1}$, using just the Cauchy-Schwarz inequality. On the other hand, if $2^k\geq\langle t\rangle^{-1}$ then the contribution of the pairs $(k_1,k_2)$ for which $\min\{k_1,k_2\}\leq k+10$ can be bounded easily, using just \eqref{wws2}. These bounds also suffice when $\min\{k_1,k_2\}\geq k+10$, but only  if $k\geq 0$. Moreover, using just \eqref{cnb2},
\begin{equation*}
\begin{split}
\|P_kI[P_{k_1}U^{\mathcal{L}_1\psi,\iota_1},P_{k_2}U^{\mathcal{L}_2\psi,\iota_2}](t)\|_{L^\infty}&\lesssim 2^{3k}\|P_{k_1}U^{\mathcal{L}_1\psi}(t)\|_{L^2}\|P_{k_2}U^{\mathcal{L}_2\psi}(t)\|_{L^2}\\
&\lesssim \varep_1^22^{3k}2^{k_1}2^{k_2}\langle t\rangle^{\delta'}2^{-4k_1^+}2^{-4k_2^+},
\end{split}
\end{equation*}
and this suffices to bound the contribution of the pairs $(k_1,k_2)$ for which $2^{2k}2^{k_1^-}2^{k_2^-}\lesssim \langle t\rangle^{-2+\delta'}$.

In the remaining case 
\begin{equation}\label{plk15.5}
2^k\in[\langle t\rangle^{-1},1],\qquad \min\{k_1,k_2\}\geq k+10,\qquad 2^{-k}\leq 2^{k_1^-}\langle t\rangle^{1-\delta'/2},
\end{equation}
we use the super-localized $L^\infty$ bounds \eqref{Linfty3.6*}. For this we decompose $P_{k_1}U^{\mathcal{L}_1\psi,\iota_1}=U^{\mathcal{L}_1\psi,\iota_1}_{\leq J_1,k_1}+U^{\mathcal{L}_1\psi,\iota_1}_{> J_1,k_1}$ and $P_{k_2}U^{\mathcal{L}_2\psi,\iota_2}=U^{\mathcal{L}_2\psi,\iota_2}_{\leq J_2,k_2}+U^{\mathcal{L}_2\psi,\iota_2}_{> J_2,k_2}$, where $2^{J_1}=2^{J_2}\approx 2^{k_1^-}\langle t\rangle^{1-\delta'/2}$. Then we estimate
\begin{equation*}
\begin{split}
\|P_kI[U^{\mathcal{L}_1\psi,\iota_1}_{> J_1,k_1},&P_{k_2}U^{\mathcal{L}_2\psi,\iota_2}](t)\|_{L^\infty}\lesssim 2^{3k/2}\|U^{\mathcal{L}_1\psi,\iota_1}_{> J_1,k_1}(t)\|_{L^2}\|P_{k_2}U^{\mathcal{L}_2\psi,\iota_2}(t)\|_{L^\infty}\\
&\lesssim \varep_1^22^{3k/2}2^{-J_1}\langle t\rangle^{-1+\delta'}2^{k_2^-/2}2^{-4k_1^+}2^{-4k_2^+}\lesssim\varep_1^22^{3k/2}\langle t\rangle^{-2+3\delta'/2}2^{-k_2^-/2}2^{-4k_1^+}2^{-4k_2^+},
\end{split}
\end{equation*}
using \eqref{vcx1.1} and \eqref{wws2}. Similarly,
\begin{equation*}
\begin{split}
\|P_kI[U^{\mathcal{L}_1\psi,\iota_1}_{\leq J_1,k_1},&U^{\mathcal{L}_2\psi,\iota_2}_{>J_2,k_2}](t)\|_{L^\infty}\lesssim\varep_1^22^{3k/2}\langle t\rangle^{-2+3\delta'/2}2^{-k_2^-/2}2^{-4k_1^+}2^{-4k_2^+}.
\end{split}
\end{equation*}
Finally, using \eqref{Linfty3.6*} with $l=k$ and recalling the definition \eqref{suploc}, we estimate
\begin{equation*}
\begin{split}
\big\|P_kI[U^{\mathcal{L}_1\psi,\iota_1}_{\leq J_1,k_1},U^{\mathcal{L}_2\psi,\iota_2}_{\leq J_2,k_2}](t)\big\|_{L^\infty}&\lesssim\sum_{|n_1-n_2|\leq 4}\big\|\mathcal{C}_{n_1,l}U^{\mathcal{L}_1\psi,\iota_1}_{\leq J_1,k_1}\big\|_{L^\infty}\big\|\mathcal{C}_{n_2,l}U^{\mathcal{L}_2\psi,\iota_2}_{\leq J_2,k_2}(t)\big\|_{L^\infty}\\
&\lesssim \varep_1^2\langle t\rangle^{-2+\delta'}2^k2^{-4k_1^+}2^{-4k_2^+},
\end{split}
\end{equation*}
where we also used \eqref{cnb2} in the estimate in the second line. Therefore
\begin{equation}\label{plk16}
\|P_kI[P_{k_1}U^{\mathcal{L}_1\psi,\iota_1},P_{k_2}U^{\mathcal{L}_2\psi,\iota_2}](t)\|_{L^\infty}\lesssim \varep_1^22^{k}\langle t\rangle^{-2+3\delta'/2}2^{-4k_1^+}2^{-4k_2^+},
\end{equation}
which suffices to bound the contribution of the remaining pairs $(k_1,k_2)$ as in \eqref{plk15.5}.

{\bf{Step 3.}} To prove the weighted $L^2$ bounds in \eqref{gb17.2} it suffices to show that
\begin{equation}\label{sac33}
\begin{split}
\sum_{(k_1,k_2)\in\mathcal{X}_k,\,k_1\geq k_2}\|P_k\{x_lI[P_{k_1}U^{\mathcal{L}_1\psi,\iota_1},&P_{k_2}U^{\mathcal{L}_2\psi,\iota_2}]\}(t)\|_{L^2}\\
&\lesssim \varep_1^22^{k/2}2^{-N(n+1)k^+-5k^+}\langle t\rangle^{\delta[H(q,n)+\widetilde{\ell}(n)]}.
\end{split}
\end{equation}
Notice that we may assume that $k_1\geq k_2$ in \eqref{sac33}, due to symmetry. We write $U^{\mathcal{L}_j\psi,\iota_j}(t)=e^{-it\Lambda_{kg,\iota_j}}V^{\mathcal{L}_j\psi,\iota_j}(t)$, $j\in\{1,2\}$. 
When proving \eqref{sac33}, the $\partial_{\xi_l}$ derivative can hit either the multiplier $m(\xi-\eta,\eta)$, or the phase $e^{-it\Lambda_{kg,\iota_1}(\xi-\eta)}$, or the profile $\widehat{P_{k_1}V^{\mathcal{L}_1\psi,\iota_1}}(\xi-\eta)$. In the first two cases, the derivative effectively corresponds to multiplying by factors $\lesssim 2^{-k_1^-}$ or $\lesssim \langle t\rangle$, and the desired bounds are again consequences of Lemma \ref{box31}.

Let $\widehat{U^{\mathcal{L}_j\psi,\iota_j}_{\ast l,k_j}}(\xi,t):=e^{-it\Lambda_{kg,\iota_j}(\xi)}\partial_{\xi_l}\{\varphi_{k_j}\cdot \widehat{V^{\mathcal{L}_j\psi,\iota_j}}\}(\xi,t)$. In view of \eqref{cnb1} and \eqref{cnb2} we have
\begin{equation}\label{sac51}
\|U^{\mathcal{L}_j\psi,\iota_j}_{\ast l,k_j}(t)\|_{L^2}\lesssim \varep_1\langle t\rangle^{H(q_j+1,n_j+1)\delta}2^{-N(n_j+1)k_j^+}.
\end{equation}
Using these $L^2$ bounds, the $L^\infty$ bounds \eqref{wws2}, and the $L^2$ bounds \eqref{vcx1} we estimate
\begin{equation}\label{sac51.1}
\|P_k\{I[U^{\mathcal{L}_1\psi,\iota_1}_{\ast l,k_1},P_{k_2}U^{\mathcal{L}_2\psi,\iota_2}]\}(t)\|_{L^2}\lesssim \varep_1^2\langle t\rangle^{\delta'}\min(\langle t\rangle^{-1},2^{k_2^-},2^{k^-})2^{-8k_2^+}2^{-N(n_1+1)k_1^+}.
\end{equation}
 This suffices to prove \eqref{sac33}, except if there is derivative loss, which happens when 
\begin{equation}\label{sac33.1}
n_1=n\quad \text{ and }\quad 2^{k}\geq\langle t\rangle^{1/100}2^{10}.
\end{equation}
In this case the estimates \eqref{sac51.1} still suffice to bound the contribution of the pairs $(k_1,k_2)$ as in \eqref{sac33}, unless $k_2\in[-10k,k-10]$. In this last case we make the change of variables $\eta\to\xi-\eta$ in order to move the $\partial_{\xi_l}$ derivative on the low frequency factor, and estimate
\begin{equation*}
\|P_k\{I[P_{k_1}U^{\mathcal{L}_1\psi,\iota_1},U^{\mathcal{L}_2\psi,\iota_2}_{\ast l,k_2}]\}(t)\|_{L^2}\lesssim \varep_1^2\langle t\rangle^{\delta'}2^{-N(n)k_1^+}2^{3k_2^-/2},
\end{equation*}
for any $(k_1,k_2)\in\mathcal{X}_k$ with $k_2\leq k_1$. We remark that the loss of the factor of $\langle t\rangle ^{\delta'}$ is mitigated by the gain of derivative and the assumption $2^{k}\gtrsim\langle t\rangle^{1/100}$. This suffices to bound the remaining contributions as claimed in \eqref{sac33}.
\end{proof}

We consider now the quadratic nonlinearities $\mathcal{Q}^2_{\al\be}$ and $\mathcal{S}^2_{\al\be}$.

\begin{lemma}\label{abc2.3}
Assume that $t\in[0,T]$, $\mathcal{L}\in\mathcal{V}_n^q$, $q\leq n\leq 3$, and $k\in\mathbb{Z}$. Then
\begin{equation}\label{gb19.1}
\begin{split}
\|P_k\{\mathcal{L}\mathcal{Q}^2_{\al\be}\}(t)\|_{L^2}&+\|P_k\{\mathcal{L}\mathcal{S}^2_{\al\be}\}(t)\|_{L^2}\\
&\lesssim \varep_1^22^{k/2}2^{-\widetilde{N}(n)k^++2k^+}\langle t\rangle^{H(q,n)\delta }\min(2^{k}\langle t\rangle^{3\delta/2},\langle t\rangle^{-1+\widetilde{\ell}(n)\delta}).
\end{split}
\end{equation}
Moreover, if $n\leq 2$ and $l\in\{1,2,3\}$ then we also have the bounds
\begin{equation}\label{gb19.2}
\begin{split}
\|P_k\{\mathcal{L}\mathcal{Q}^2_{\al\be}\}(t)\|_{L^\infty}+\|P_k\{\mathcal{L}\mathcal{S}^2_{\al\be}\}(t)\|_{L^\infty}&\lesssim\varep_1^22^{k}\langle t\rangle^{-1+2\delta'}2^{-N(n+1)k^++3k^+}\min\{2^{k},\langle t\rangle^{-1}\},\\
\|P_k\{x_l\mathcal{L}\mathcal{Q}^2_{\al\be}\}(t)\|_{L^2}+\|P_k\{x_l\mathcal{L}\mathcal{S}^2_{\al\be}\}(t)\|_{L^2}&\lesssim \varep_1^22^{k/2}2^{-N(n+1)k^+-5k^+}\langle t\rangle^{\delta[H(q,n)+\widetilde{\ell}(n)]}.
\end{split}
\end{equation}
\end{lemma}

\begin{proof} We notice that $\mathcal{S}^2_{\al\be}$ is a sum of quadratic expressions of the form $\partial_\mu h_1\cdot\partial_\nu h_2$, where $h_1,h_2\in\{h_{\al\be}\}$ and $\mathcal{Q}^2_{\al\be}$ is a sum of quadratic expressions of the form $h_1\cdot\partial_\mu\partial_\nu h_2$, $(\mu,\nu)\neq (0,0)$. The differential operator $\LL$ can split between these two factors. Notice that commutation with one vector-field $\partial_\mu$ generates similar terms of the form $\partial_\nu\mathcal{L}'$ (see \eqref{qaz3.5}), while commutation with two vector-fields $\partial_\mu\partial_\nu$ (which appear in $\mathcal{Q}^2_{\al\be}$) leads to terms of the form $\partial_\al\partial_\be\mathcal{L'}$. If $\al=\be=0$ then we have to further replace $\partial_t^2\mathcal{L'}h$ with $\Delta\mathcal{L'}h+\mathcal{L'}\mathcal{N}^h$ (as in \eqref{sac1}), in order to have access to elliptic estimates.

{\bf{Step 1.}} In view of these considerations, for \eqref{gb19.1} it suffices to prove that
\begin{equation}\label{sac62}
\begin{split}
\sum_{(k_1,k_2)\in\mathcal{X}_k}&2^{|k_2-k_1|}\|P_kI[P_{k_1}U^{\mathcal{L}_1h_1,\iota_1},P_{k_2}U^{\mathcal{L}_2h_2,\iota_2}](t)\|_{L^2}\\
&\lesssim \varep_1^22^{k/2}2^{-\widetilde{N}(n)k^++2k^+}\langle t\rangle^{H(q,n)\delta }\min(2^{k}\langle t\rangle^{3\delta/2},\langle t\rangle^{-1+\widetilde{\ell}(n)\delta}),
\end{split}
\end{equation}
and
\begin{equation}\label{sac62.5}
\begin{split}
\sum_{(k_1,k_2)\in\mathcal{X}_k}&2^{-k_1}\|P_kI[P_{k_1}U^{\mathcal{L}_1h_1,\iota_1},P_{k_2}\mathcal{L}'\mathcal{N}^{h}](t)\|_{L^2}\\
&\lesssim \varep_1^22^{k/2}2^{-\widetilde{N}(n)k^++2k^+}\langle t\rangle^{H(q,n)\delta }\min(2^{k}\langle t\rangle^{3\delta/2},\langle t\rangle^{-1+\widetilde{\ell}(n)\delta}),
\end{split}
\end{equation}
for any $\iota_1,\iota_2\in\{+,-\}$, $h_1,h_2\in\{h_{\al\be}\}$, $\mathcal{N}^h\in\{\mathcal{N}^h_{\al\be}\}$, $\mathcal{L}_1\in\mathcal{V}_{n_1}^{q_1}$, $\mathcal{L}_2\in\mathcal{V}_{n_2}^{q_2}$, $\mathcal{L}'\in\mathcal{V}_{n'}^{q'}$, $(q_1,n_1)+(q_2,n_2)\leq (q,n)$,  $(q_1,n_1)+(q',n')\leq (q,n-1)$. 

Without loss of generality, in proving \eqref{sac62} we may assume that $n_1\leq n_2$. We use Lemma \ref{box1}. If $2^k\langle t\rangle^{3\delta/2}\lesssim\langle t\rangle^{-1+\widetilde{\ell}(n)\delta}$ then the bounds \eqref{sac62} follow from \eqref{box2} and \eqref{SuperlinearH1}. On the other hand, if $2^k\langle t\rangle^{3\delta/2}\lesssim\langle t\rangle^{-1+\widetilde{\ell}(n)\delta}$ then the contribution of the pairs $(k_1,k_2)\in\mathcal{X}_k$ with $\min(k_1,k_2)\geq k$ (thus $|k_1-k_2|\leq 10$) is bounded as claimed due to \eqref{box3}, \eqref{box3.5}, and the first inequality in \eqref{SuperlinearH2}. The contribution of the pairs $(k_1,k_2)\in\mathcal{X}_k$ with $k_1=\min(k,k_1,k_2)$ (thus $|k_2-k|\leq 4$) is bounded as claimed due to \eqref{box4} (this case gives the worst contribution to the growth function $\langle t\rangle^{\delta[H(q,n)+\widetilde{\ell}(n)]}$, when $n_1=0$). Finally, the contribution of the pairs $(k_1,k_2)\in\mathcal{X}_k$ with $k_2=\min(k,k_1,k_2)$ is bounded as claimed due to \eqref{box5} if $n_1=1$ and \eqref{box5.5}--\eqref{box5.6} if $n_1=0$. 

To prove \eqref{sac62.5} we use only $L^2$ bounds. We may assume, by induction over $n$ (in the case $n=0$ the right-hand side of \eqref{sac62.5} is trivial), that the bounds \eqref{wer4.0} hold for $P_{k_2}\mathcal{N}^{\mathcal{L}'h'}(t)$ and estimate the left-hand side of \eqref{sac62.5}, using also \eqref{vcx1}, by
\begin{equation}\label{plk10}
\begin{split}
\sum_{(k_1,k_2)\in\mathcal{X}_k}2^{-k_1}&2^{3\min\{k,k_1,k_2\}/2}\|P_{k_1}U^{\mathcal{L}_1h_1,\iota_1}(t)\|_{L^2}\|P_{k_2}\LL'\mathcal{N}^{h}(t)\|_{L^2}\\
&\lesssim \varep_1^22^{3k/2}2^{-\widetilde{N}(n)k^+}\langle t\rangle^{-1+\delta[H(q_1,n_1)+H(q',n')+\widetilde{\ell}(n')+2]}.
\end{split}
\end{equation}
This suffices to prove \eqref{sac62.5}, using also \eqref{SuperlinearH1}.

{\bf{Step 2.}} Similarly, to prove the $L^\infty$ bounds in \eqref{gb19.2} it suffices to show that
\begin{equation}\label{plk10.1}
\begin{split}
\sum_{(k_1,k_2)\in\mathcal{X}_k}&2^{|k_2-k_1|}\|P_kI[P_{k_1}U^{\mathcal{L}_1h_1,\iota_1},P_{k_2}U^{\mathcal{L}_2h_2,\iota_2}](t)\|_{L^\infty}\\
&\lesssim\varep_1^22^{k}\langle t\rangle^{-1+2\delta'}2^{-N(n+1)k^++3k^+}\min\{2^{k},\langle t\rangle^{-1}\}
\end{split}
\end{equation}
and
\begin{equation}\label{plk10.2}
\begin{split}
\sum_{(k_1,k_2)\in\mathcal{X}_k}&2^{-k_1}\|P_kI[P_{k_1}U^{\mathcal{L}_1h_1,\iota_1},P_{k_2}\LL'\mathcal{N}^{h}](t)\|_{L^\infty}\\
&\lesssim\varep_1^22^{k}\langle t\rangle^{-1+2\delta'}2^{-N(n+1)k^++3k^+}\min\{2^{k},\langle t\rangle^{-1}\}.
\end{split}
\end{equation}

These bounds follow easily from \eqref{sac62}--\eqref{sac62.5} if $2^k\lesssim \langle t\rangle^{-1+\delta'}$. Also, the bounds \eqref{plk10.2} follow from \eqref{wws1} and using Proposition \ref{plk1} inductively.

We prove now the bounds \eqref{plk10.1}. The contribution of the pairs $(k_1,k_2)$ with $\min\{k_1,k_2\}\leq k+10$ or $\max\{k_1,k_2\}\geq \langle t\rangle$ can be bounded easily, using just \eqref{wws2}. This also suffices to bound all the contributions if $k\geq 0$. In the remaining case, the desired bounds follow from \eqref{plk10.3}.

{\bf{Step 3.}} Similarly, to prove the weighted $L^2$ bounds in \eqref{gb19.2} it suffices to show that
\begin{equation}\label{sac63}
\begin{split}
\sum_{(k_1,k_2)\in\mathcal{X}_k,\,k_2\geq k_1}&2^{k_2-k_1}\|P_k\{x_lI[P_{k_1}U^{\mathcal{L}_1h_1,\iota_1},P_{k_2}U^{\mathcal{L}_2h_2,\iota_2}]\}(t)\|_{L^2}\\
&\lesssim \varep_1^22^{k/2}2^{-N(n+1)k^+-5k^+}\langle t\rangle^{\delta[H(q,n)+\widetilde{\ell}(n)]},
\end{split}
\end{equation}
and
\begin{equation}\label{sac63.5}
\begin{split}
\sum_{(k_1,k_2)\in\mathcal{X}_k}&2^{-k_1}\|P_k\{x_lI[P_{k_1}U^{\mathcal{L}_1h_1,\iota_1},P_{k_2}\LL'\mathcal{N}^{h}]\}(t)\|_{L^2}\\
&\lesssim \varep_1^22^{k/2}2^{-N(n+1)k^+-5k^+}\langle t\rangle^{\delta[H(q,n)+\widetilde{\ell}(n)]}.
\end{split}
\end{equation}
Notice that in \eqref{sac63} it is more convenient to assume that $k_1\leq k_2$ instead of $n_1\leq n_2$. 

To prove the bounds \eqref{sac63} we write $U^{\mathcal{L}_j h_j,\iota_j}(t)=e^{-it\Lambda_{wa,\iota_j}}V^{\mathcal{L}_j h_j,\iota_j}(t)$, $j\in\{1,2\}$ and let $\widehat{U^{\mathcal{L}_jh_j,\iota_j}_{\ast l,k_j}}(\xi,t):=e^{-it\Lambda_{wa,\iota_j}(\xi)}\partial_{\xi_l}\{\varphi_{k_j}\cdot \widehat{V^{\mathcal{L}_jh_j,\iota_j}}\}(\xi,t)$, as in Lemma \ref{abc2.2}. In view of \eqref{cnb1} and \eqref{vcx1} it follows that, for $j\in\{1,2\}$,
\begin{equation}\label{sac81}
\|U^{\mathcal{L}_jh_j,\iota_j}_{\ast l,k_j}(t)\|_{L^2}\lesssim \varep_1\langle t\rangle^{H(q_j+1,n_j+1)\delta}2^{-k_j/2}2^{-N(n_j+1)k_j^+}(\langle t\rangle 2^{k_j^-})^{-\ga}.
\end{equation}

To prove \eqref{sac63} we make the change of variables $\eta\to\xi-\eta$ and notice that the $\partial_{\xi_l}$ derivative can hit either the multiplier $m(\eta,\xi-\eta)$, or the phase $e^{-it\Lambda_{wa,\iota_2}(\xi-\eta)}$, or the higher frequency profile $\widehat{P_{k_2}V^{\mathcal{L}_2h_2,\iota_2}}(\xi-\eta)$. In the first two cases, the derivative effectively corresponds to multiplying by factors $\lesssim \langle t\rangle+2^{-k_2^-}$, and the desired bounds are then consequences of \eqref{sac62}. In the last case we estimate 
\begin{equation}\label{sac81.5}
\begin{split}
2^{k_2-k_1}&\|P_k\{I[P_{k_1}U^{\mathcal{L}_1h_1,\iota_1},U^{\mathcal{L}_2h_2,\iota_2}_{\ast l,k_2}]\}(t)\|_{L^2}\\
&\lesssim \varep_1^2\langle t\rangle^{\delta'}\min(\langle t\rangle^{-1},2^{k_1^-},2^{k^-})^{9/10}2^{-8k_1^+}2^{-N(n_2+1)k_2^++k_2^+},
\end{split}
\end{equation}
using the $L^2$ bounds \eqref{sac81} and the bounds \eqref{wws1} and \eqref{vcx1}. As in Lemma \ref{abc2.2}, this suffices except if there is derivative loss, which can happen only when 
\begin{equation}\label{sac82}
n_2=n\qquad\text{ and }\qquad 2^{k}\geq\langle t\rangle^{1/100}2^{10}.
\end{equation}
 
Assuming that \eqref{sac82} holds, the sum over $k_1\leq -10k$ or $k_1\geq k-10$ can still be bounded as before, using \eqref{sac81.5}. To bound the sum over $k_1\in[-10k,k-10]$ we return to the original formula for $I$ (without making the the change of variables $\eta\to\xi-\eta$) and notice that the $\partial_{\xi_l}$ derivative now hits the low frequency factor. Then we estimate, using just \eqref{sac81} and \eqref{vcx1},
\begin{equation*}
2^{k_2-k_1}\|P_k\{I[U^{\mathcal{L}_1h_1,\iota_1}_{\ast l,k_1},P_{k_2}U^{\mathcal{L}_2h_2,\iota_2}]\}(t)\|_{L^2}\lesssim \varep_1^2\langle t\rangle^{\delta'}2^{-N(n)k_2^++2k_2^+}2^{-\delta k_1^-}.
\end{equation*}
As before, the loss of the factor of $\langle t\rangle ^{\delta'}$ is compensated by the gain of derivative and the assumption $2^{k}\gtrsim\langle t\rangle^{1/100}$. This suffices to complete the proof of \eqref{sac63}. 

The bounds \eqref{sac63.5} are easier: as in \eqref{plk10} we estimate the left-hand side by
\begin{equation*}
\begin{split}
\sum_{(k_1,k_2)\in\mathcal{X}_k}2^{-k_1}&2^{3\min\{k,k_1,k_2\}/2}\|P_{k_1}U^{\mathcal{L}_1h_1,\iota_1}(t)\|_{L^2}\big[\|P_{k_2}\{x_l\LL'\mathcal{N}^{h}\}(t)\|_{L^2}+2^{-k_2^-}\|P'_{k_2}\LL'\mathcal{N}^{h}(t)\|_{L^2}\big]\\
&\lesssim \varep_1^22^{k}2^{-N(n+1)k^+-8k^+}\langle t\rangle^{\delta[H(q_1,n_1)+H(q',n')+\widetilde{\ell}(n')+2]},
\end{split}
\end{equation*}
where $P'_{k_2}=P_{[k-2,k+2]}$, using \eqref{vcx1} and Proposition \ref{plk1} inductively. This clearly suffices.
\end{proof}

Finally, we consider the quadratic nonlinearities $\mathcal{N}^{\psi,2}$.

\begin{lemma}\label{dtv8.2}
Assume that $t\in[0,T]$, $\mathcal{L}\in\mathcal{V}_n^q$, $n\leq 3$, and $k\in\mathbb{Z}$. Then
\begin{equation}\label{gb21.1}
\|P_k\{\mathcal{L}\mathcal{N}^{\psi,2}\}(t)\|_{L^2}\lesssim \varep_1^22^{-\widetilde{N}(n)k^++7k^+}\langle t\rangle^{H(q,n)\delta }\min(2^{k}\langle t\rangle^{3\delta/2},\langle t\rangle^{-1+\widetilde{\ell}(n)\delta}).
\end{equation}
Moreover, if $l\in\{1,2,3\}$ and $n\leq 2$ then we also have the bounds
\begin{equation}\label{gb21.2}
\begin{split}
\|P_k\{\mathcal{L}\mathcal{N}^{\psi,2}\}(t)\|_{L^\infty}&\lesssim\varep_1^22^{k/2}\langle t\rangle^{-1+2\delta'}2^{-N(n+1)k^++3k^+}\min\{2^{k},\langle t\rangle^{-1}\},\\
\|P_k\{x_l\mathcal{L}\mathcal{N}^{\psi,2}\}(t)\|_{L^2}&\lesssim \varep_1^22^{-N(n+1)k^+-3k^+}\langle t\rangle^{\delta[H(q,n)+\widetilde{\ell}(n)]}.
\end{split}
\end{equation}
\end{lemma} 

\begin{proof} We examine the definition \eqref{npsi2} and notice that $\mathcal{N}^{\psi,2}$ is a sum of terms of the form $h\cdot\partial_\mu\partial_\nu \psi$, $(\mu,\nu)\neq (0,0)$, or $h\cdot\psi$, where $h\in\{h_{\mu\nu}\}$. As in Lemma \ref{abc2.3}, we distributed the  vector-field $\LL$, commute with the derivatives $\partial_\mu\partial_\nu$, and  further replace $\partial_t^2\mathcal{L'}\psi$ with $(\Delta-1)\mathcal{L'}\psi+\mathcal{L'}\mathcal{N}^\psi$.

{\bf{Step 1.}} For \eqref{gb21.1} it suffices to prove that
\begin{equation}\label{gb21.3}
\begin{split}
\sum_{(k_1,k_2)\in\mathcal{X}_k}2^{k_2^+-k_1}&\|P_kI[P_{k_1}U^{\mathcal{L}_1h,\iota_1},P_{k_2}U^{\mathcal{L}_2\psi,\iota_2}](t)\|_{L^2}\\
&\lesssim\varep_1^22^{-\widetilde{N}(n)k^++7k^+}\langle t\rangle^{H(q,n)\delta }\min(2^{k}\langle t\rangle^{3\delta/2},\langle t\rangle^{-1+\widetilde{\ell}(n)\delta}),
\end{split}
\end{equation}
and
\begin{equation}\label{plk17}
\begin{split}
\sum_{(k_1,k_2)\in\mathcal{X}_k}2^{-k_1}&\|P_kI[P_{k_1}U^{\mathcal{L}_1h,\iota_1},P_{k_2}\mathcal{L}'\mathcal{N}^{\psi}](t)\|_{L^2}\\
&\lesssim\varep_1^22^{-\widetilde{N}(n)k^++7k^+}\langle t\rangle^{H(q,n)\delta }\min(2^{k}\langle t\rangle^{3\delta/2},\langle t\rangle^{-1+\widetilde{\ell}(n)\delta}),
\end{split}
\end{equation}
for any $\iota_1,\iota_2\in\{+,-\}$, $h\in\{h_{\al\be}\}$, $\mathcal{L}_1\in\mathcal{V}_{n_1}^{q_1}$, $\mathcal{L}_2\in\mathcal{V}_{n_2}^{q_2}$, $\mathcal{L}'\in\mathcal{V}_{n'}^{q'}$, $(q_1,n_1)+(q_2,n_2)\leq (q,n)$, $(q_1,n_1)+(q',n')\leq (q,n-1)$. 

To prove \eqref{gb21.3} we use Lemma \ref{box80}. The bounds follow from \eqref{box82} if $2^{k}\langle t\rangle^{3\delta/2}\lesssim\langle t\rangle^{-1+\widetilde{\ell}(n)\delta}$. If $2^{k}\langle t\rangle^{3\delta/2}\geq\langle t\rangle^{-1+\widetilde{\ell}(n)\delta}$ then \eqref{gb21.3} follow from \eqref{box83.5} when $n_1=n_2=0$, or from \eqref{box83} if $n_1=0$ and $n_2\geq 1$, or from \eqref{box84} if $n_1\geq 1$ and $n_2=0$, or from \eqref{box84.5} if $n_1\geq 1$ and $n_2\geq 1$.

The bounds \eqref{plk17} follow easily, using just $L^2$ estimates as in \eqref{plk10}.

{\bf{Step 2.}} To prove the $L^\infty$ bounds in \eqref{gb21.2} it suffices to show that
\begin{equation}\label{plk18}
\begin{split}
\sum_{(k_1,k_2)\in\mathcal{X}_k}2^{k_2^+-k_1}&\|P_kI[P_{k_1}U^{\mathcal{L}_1h,\iota_1},P_{k_2}U^{\mathcal{L}_2\psi,\iota_2}](t)\|_{L^\infty}\\
&\lesssim\varep_1^22^{k/2}\langle t\rangle^{-1+2\delta'}2^{-N(n+1)k^++3k^+}\min\{2^{k},\langle t\rangle^{-1}\},
\end{split}
\end{equation}
and
\begin{equation}\label{plk19}
\begin{split}
\sum_{(k_1,k_2)\in\mathcal{X}_k}2^{-k_1}&\|P_kI[P_{k_1}U^{\mathcal{L}_1h,\iota_1},P_{k_2}\mathcal{L}'\mathcal{N}^{\psi}](t)\|_{L^\infty}\\
&\lesssim\varep_1^22^{k/2}\langle t\rangle^{-1+2\delta'}2^{-N(n+1)k^++3k^+}\min\{2^{k},\langle t\rangle^{-1}\}.
\end{split}
\end{equation}

These bounds follow easily from \eqref{gb21.3}--\eqref{plk17} if $2^k\lesssim \langle t\rangle^{-1+\delta'}$, using the Cauchy-Schwarz inequality. Also, the bounds \eqref{plk19} follow from \eqref{wws1} and using Proposition \ref{plk1} inductively.

We prove now the bounds \eqref{plk18}. The contribution of the pairs $(k_1,k_2)$ with $\min\{k_1,k_2\}\leq k+10$ or $\max\{k_1,k_2\}\geq \langle t\rangle$ can be bounded easily, using just \eqref{wws1} and \eqref{wws2}. This also suffices to bound all the contributions if $k\geq 0$. To summarize, it remains to show that 
\begin{equation}\label{plk19.5}
2^{-k_1}\|P_kI[P_{k_1}U^{\mathcal{L}_1h_1,\iota_1},P_{k_2}U^{\mathcal{L}_2\psi,\iota_2}](t)\|_{L^\infty}\lesssim\varep_1^22^{k/2}\langle t\rangle^{-2+3\delta'/2}2^{-2k_1^+}2^{-2k_2^+},
\end{equation}
provided that $2^k\in[\langle t\rangle^{-1+\delta'},1]$, $|k_1-k_2|\leq 4$, and $k_1,k_2\in[k+10,\langle t\rangle]$. 

This is similar to the bounds \eqref{plk16}. First we estimate the left-hand side of \eqref{plk19.5} by $2^{3k/2}\varep_1^2\langle t\rangle^{-1+\delta'/2}2^{k/2}2^{-4k_1^+}2^{-4k_2^+}$, using \eqref{cnb2} and \eqref{wws1}. This suffices if $2^{k}2^{k_2^-}\lesssim \langle t\rangle^{-1+\delta'}$. On the other hand, if $2^k\in[\langle t\rangle^{-1+\delta'},1]$ then we decompose $P_{k_1}U^{\mathcal{L}_1h_1,\iota_1}=U^{\mathcal{L}_1h_1,\iota_1}_{\leq J_1,k_1}+U^{\mathcal{L}_1h_1,\iota_1}_{> J_1,k_1}$ and $P_{k_2}U^{\mathcal{L}_2\psi,\iota_2}=U^{\mathcal{L}_2\psi,\iota_2}_{\leq J_2,k_2}+U^{\mathcal{L}_2\psi,\iota_2}_{> J_2,k_2}$, where $2^{J_1}\approx \langle t\rangle^{1-\delta'/2}$ and $2^{J_2}\approx 2^{k_2^-}\langle t\rangle^{1-\delta'/2}$. Then we estimate
\begin{equation*}
\begin{split}
2^{-k_1}\big\|P_kI[U^{\mathcal{L}_1h_1,\iota_1}_{> J_1,k_1},&P_{k_2}U^{\mathcal{L}_2\psi,\iota_2}](t)\big\|_{L^\infty}+2^{-k_1}\big\|P_kI[U^{\mathcal{L}_1h_1,\iota_1}_{\leq J_1,k_1},U^{\mathcal{L}_2\psi,\iota_2}_{>J_2,k_2}](t)\big\|_{L^\infty}\\
&\lesssim\varep_1^22^{3k/2}2^{-k_2}\langle t\rangle^{-2+3\delta'/2}2^{-4k_1^+}2^{-4k_2^+},
\end{split}
\end{equation*}
using \eqref{vcx1.1}, \eqref{wws1}, and \eqref{wws2}. Morever, using \eqref{Linfty1.6*} and \eqref{Linfty3.6*} with $l=k$ we estimate
\begin{equation*}
\begin{split}
2^{-k_1}\big\|P_kI[U^{\mathcal{L}_1h_1,\iota_1}_{\leq J_1,k_1},U^{\mathcal{L}_2\psi,\iota_2}_{\leq J_2,k_2}](t)\big\|_{L^\infty}&\lesssim 2^{-k_1}\sum_{|n_1-n_2|\leq 4}\big\|\mathcal{C}_{n_1,l}U^{\mathcal{L}_1h_1,\iota_1}_{\leq J_1,k_1}\big\|_{L^\infty}\big\|\mathcal{C}_{n_2,l}U^{\mathcal{L}_2\psi,\iota_2}_{\leq J_2,k_2}(t)\big\|_{L^\infty}\\
&\lesssim \varep_1^2\langle t\rangle^{-2+\delta'}2^k2^{-k_1/2}2^{-4k_1^+}2^{-4k_2^+},
\end{split}
\end{equation*}
using also \eqref{vcx1} and \eqref{cnb2} in the estimate in the second line. This completes the proof of \eqref{plk19.5}.

{\bf{Step 3.}} As before, to prove the weighted $L^2$ bounds in \eqref{gb21.2} it suffices to show that
\begin{equation}\label{abc99.1}
\sum_{(k_1,k_2)\in\mathcal{X}_k}2^{k_2^+-k_1}\|P_k\{x_lI[P_{k_1}U^{\mathcal{L}_1h,\iota_1},P_{k_2}U^{\mathcal{L}_2\psi,\iota_2}]\}(t)\|_{L^2}\lesssim\varep_1^22^{-N(n+1)k^+-3k^+}\langle t\rangle^{\delta[H(q,n)+\widetilde{\ell}(n)]}
\end{equation}
and
\begin{equation}\label{plk22}
\sum_{(k_1,k_2)\in\mathcal{X}_k}2^{-k_1}\|P_k\{x_lI[P_{k_1}U^{\mathcal{L}_1h,\iota_1},P_{k_2}\mathcal{L}'\mathcal{N}^{\psi}]\}(t)\|_{L^2}\lesssim\varep_1^22^{-N(n+1)k^+-3k^+}\langle t\rangle^{\delta[H(q,n)+\widetilde{\ell}(n)]}.
\end{equation}

The bounds \eqref{plk22} are easy. As in the proof of \eqref{sac63.5} we estimate the left-hand side by
\begin{equation*}
\begin{split}
\sum_{(k_1,k_2)\in\mathcal{X}_k}2^{-k_1}&2^{3\min\{k,k_1,k_2\}/2}\|P_{k_1}U^{\mathcal{L}_1h_1,\iota_1}(t)\|_{L^2}\big[\|P_{k_2}\{x_l\mathcal{L}'\mathcal{N}^{\psi}\}(t)\|_{L^2}+2^{-k_2^-}\|P'_{k_2}\mathcal{L}'\mathcal{N}^{\psi}(t)\|_{L^2}\big]\\
&\lesssim \varep_1^22^{k/2}2^{-N(n+1)k^+-8k^+}\langle t\rangle^{\delta[H(q_1,n_1)+H(q',n')+\widetilde{\ell}(n')+2]},
\end{split}
\end{equation*}
using \eqref{vcx1} and Proposition \ref{plk1} inductively. This clearly suffices.

To prove \eqref{abc99.1} we write $U^{\mathcal{L}_2\psi,\iota_2}=e^{-it\Lambda_{kg,\iota_2}}V^{\mathcal{L}_2\psi,\iota_2}$ and examine the formula \eqref{abc36.1}. We make the change of variables $\eta\to\xi-\eta$ and notice that the $\partial_{\xi_l}$ derivative can hit either the multiplier $m(\eta,\xi-\eta)$, or the phase $e^{-it\Lambda_{kg,\iota_2}(\xi-\eta)}$, or the profile $\widehat{P_{k_2}V^{\mathcal{L}_2\psi,\iota_2}}(\xi-\eta)$. In the first two cases, the derivative effectively corresponds to multiplying by factors $\lesssim \langle t\rangle$ or $\lesssim 2^{-k_2^-}$, and changing the multiplier $m_1$, in a way that still satisfies \eqref{mults}. The desired bounds are then consequences of the bounds \eqref{gb21.3} (in the case $2^{k_2^-}\lesssim \langle t\rangle^{-1}$ we need to apply \eqref{box82} again to control the corresponding contributions).

It remains to consider the case when the $\partial_{\xi_l}$ derivative hits the profile $\widehat{P_{k_2}V^{\mathcal{L}_2\psi,\iota_2}}(\xi-\eta)$. It suffices to prove that 
\begin{equation}\label{abc99.2}
\sum_{(k_1,k_2)\in\mathcal{X}_k}2^{k_2^+-k_1}\|P_kI[P_{k_1}U^{\mathcal{L}_1h,\iota_1},U^{\mathcal{L}_2\psi,\iota_2}_{\ast l,k_2}](t)\|_{L^2}\lesssim\varep_1^2\langle t\rangle^{\delta[H(q,n)+\widetilde{\ell}(n)]}2^{-N(n+1)k^+-4k^+},
\end{equation}
where, as in Lemma \ref{abc2.2}, $\widehat{U^{\mathcal{L}_2\psi,\iota_2}_{\ast l,k_2}}(\xi,t)=e^{-it\Lambda_{kg,\iota_2}(\xi)}\partial_{\xi_l}\{\varphi_{k_2}\cdot \widehat{V^{\mathcal{L}_2\psi,\iota_2}}\}(\xi,t)$. We have, see \eqref{sac51},
\begin{equation}\label{abc99.3}
\|U^{\mathcal{L}_2\psi,\iota_2}_{\ast l,k_2}(t)\|_{L^2}\lesssim \varep_1\langle t\rangle^{H(q_2+1,n_2+1)\delta}2^{-N(n_2+1)k_2^+}.
\end{equation}

Assume first that $n_1\geq 1$. The contribution of the pairs $(k_1,k_2)$ in \eqref{abc99.2} for which $k_1\leq k_2+10$ is bounded as claimed, using the $L^\infty$ estimates \eqref{wws1} and the $L^2$ bounds \eqref{abc99.3}. Similarly, the contribution of the pairs $(k_1,k_2)$ for  which $k_2\leq k_1-10$  (thus $|k-k_1|\leq 4$) is bounded as claimed if $2^k\lesssim \langle t\rangle^{1/100}$. Finally, if $2^k\geq \langle t\rangle^{1/100}$ then the contribution of the pairs $(k_1,k_2)$ for  which $k_2\leq k_1-10$ is bounded as claimed using \eqref{vcx1} and \eqref{abc99.3}.

Assume now that 
\begin{equation}\label{abc99.4}
n_1=0\qquad\text{ and }\qquad n_2=n.
\end{equation}
We need to be slightly more careful than before. If $2^{k}\lesssim \langle t\rangle^{1/100}$ then we can just use the $L^\infty$ bounds in \eqref{wws1} and the $L^2$ bounds \eqref{abc99.3} to prove \eqref{abc99.2}. On the other hand, if $2^{k}\geq \langle t\rangle^{1/100}$ then the contribution of the pairs $(k_1,k_2)$ with $k_1\geq k-10$ or $k_1\leq -10k$ can be estimated as before, using just \eqref{abc99.3} and the $L^2$ bounds in \eqref{vcx1}. 

To estimate the contribution of the remaining pairs, we need to avoid the derivative loss. We go back to \eqref{abc99.1}. It remains to prove that if $2^{k}\geq \langle t\rangle^{1/100}$ then
\begin{equation}\label{abc99.6}
\begin{split}
\sum_{(k_1,k_2)\in\mathcal{X}_k,\,k_1\in[-10k,k-10]}2^{k_2^+-k_1}&\|P_k\{\partial_{\xi_l}\mathcal{F}\{I[P_{k_1}U^{\mathcal{L}_1h,\iota_1},P_{k_2}U^{\mathcal{L}_2\psi,\iota_2}]\}\}(t)\|_{L^2_\xi}\\
&\lesssim\varep_1^22^{-N(n+1)k^+-3k^+}.
\end{split}
\end{equation}

We do not make the change of variables $\eta\to\xi-\eta$ now, so the $\partial_{\xi_l}$ derivative hits the low frequency factor or the multiplier in the definition of the operator $I$. The contribution when the derivative $\partial_{\xi_l}$ hits the multiplier can be bounded easily using \eqref{wws1}. Using the $L^2$ bounds \eqref{vcx1} on $\psi$, and $L^2\times L^\infty$ estimates as before, for \eqref{abc99.6} it suffices to prove that
\begin{equation}\label{abc99.7}
2^{-k_1}\big\|\mathcal{F}^{-1}\{\partial_{\xi_l}[\widehat{P_{k_1}U^{\mathcal{L}_1h}}(\xi,t)]\}\big\|_{L^\infty}\lesssim \varep_1\langle t\rangle^{\delta'}2^{-2k_1^+}2^{-2\delta k_1^-}.
\end{equation}

To prove \eqref{abc99.7} we replace $U^{\mathcal{L}_1h}(t)$ by $e^{-it|\nabla|}V^{\mathcal{L}_1h}(t)$ and use either \eqref{wws1} when the $\partial_{\xi_l}$ derivative hits the factor $\varphi_{k_1}(\xi)e^{-it|\nabla|}$ or  \eqref{cnb1} when the derivative hits the factor profile. The desired bounds \eqref{abc99.7} follow. This completes the proof of the lemma. 
\end{proof}

\subsection{The cubic and higher order nonlinearities} We now prove bounds on the nonlinea\-ri\-ties $\mathcal{L}\mathcal{N}^{h,\geq 3}$ and $\mathcal{L}\mathcal{N}^{\psi,\geq 3}$.

\begin{lemma}\label{abc2.1}
Assume that $t\in[0,T]$, $\mathcal{L}\in\mathcal{V}_n^q$, $n\leq 2$, $l\in\{1,2,3\}$, and $k\in\mathbb{Z}$. Then
\begin{equation}\label{abc3.00}
\|P_k(\mathcal{L}\mathcal{N}^{h,\geq 3}_{\al\be})(t)\|_{L^2}\lesssim\varep_1^2\langle t\rangle^{-1/2+4\delta'}2^{-\widetilde{N}(n)k^++6k^+}\min(2^{k},\langle t\rangle^{-1})^{3/2},
\end{equation}
\begin{equation}\label{abc6.00}
\|P_k(x_l\mathcal{L}\mathcal{N}^{h,\geq 3}_{\al\be})(t)\|_{L^2}\lesssim\varep_1^22^{k/2}\langle t\rangle^{-1/2+4\delta'}2^{-N(n+1)k^+-3k^+}.
\end{equation}
Moreover, if $n=3$ and $\mathcal{L}\in\mathcal{V}_n^q$ then
\begin{equation}\label{plk28}
\|P_k(\mathcal{L}\mathcal{N}^{h,\geq 3}_{\al\be})(t)\|_{L^2}\lesssim\varep_1^22^{k/2}\langle t\rangle^{-1/2+4\delta'}2^{-\widetilde{N}(n)k^++6k^+}\min(2^{k},\langle t\rangle^{-1}).
\end{equation}
\end{lemma}

\begin{proof} We notice that $\mathcal{N}^{h,\geq 3}_{\al\be}$ is a sum of terms of the form $G_{\geq 1}\cdot \mathcal{N}_1$, where $G_{\geq 1}=\sum_{d\geq 1}a_dg_d$ for some functions $g_d\in\mathcal{G}_d$ (see \eqref{gb0.1}) and some coefficients $a_d\in\mathbb{R}$ with $|a_d|\leq C^d$, and $\mathcal{N}_1$ is a quadratic term similar to those in $\mathcal{KG}^2_{\al\be}$, $\mathcal{Q}^2_{\al\be}$, $\mathcal{S}^2_{\al\be}$. It suffices to prove that if $\LL\in\mathcal{V}_n^q$ then
\begin{equation}\label{sac3.00}
\|P_k(\mathcal{L}(G_{\geq 1}\cdot \mathcal{N}_1))(t)\|_{L^2}\lesssim\varep_1^2\langle t\rangle^{-1/2+4\delta'}2^{-\widetilde{N}(n)k^++6k^+}\min(2^{k},\langle t\rangle^{-1})^{3/2},
\end{equation}
\begin{equation}\label{sac4.00}
\|P_k(x_l\mathcal{L}(G_{\geq 1}\cdot \mathcal{N}_1))(t)\|_{L^2}\lesssim\varep_1^22^{k/2}\langle t\rangle^{-1/2+4\delta'}2^{-N(n+1)k^+-3k^+},
\end{equation}
provided that $n\leq 2$. Moreover, if $n=3$ then
\begin{equation}\label{plk28.1}
\|P_k(\mathcal{L}(G_{\geq 1}\cdot \mathcal{N}_1))(t)\|_{L^2}\lesssim\varep_1^22^{k/2}\langle t\rangle^{-1/2+4\delta'}2^{-\widetilde{N}(n)k^++6k^+}\min(2^{k},\langle t\rangle^{-1}).
\end{equation}

Concerning the functions $G_{\geq 1}$ and $\mathcal{N}_1$, we may assume that (see \eqref{gb3} and \eqref{gb1}), 
\begin{equation}\label{abc3.01}
\begin{split}
\|P_k\mathcal{D}G_{\geq 1}(t)\|_{L^2}&\lesssim \varep_1\langle t\rangle^{\delta'}2^{-N(m)k^+}2^{-k/2} 2^{-\ga k^-},\\
\|P_k\mathcal{D} G_{\geq 1}(t)\|_{L^\infty}&\lesssim \varep_1\langle t\rangle^{-1+\delta'}2^{-N(m+1)k^++3k^+}\min\{1,2^{k^-}\langle t\rangle\}^{1-\delta},\\
\|P_k(x_l\mathcal{D}G_{\geq 1})(t)\|_{L^2}&\lesssim \varep_1\langle t\rangle^{\delta'}2^{-N(m+1)k^+}2^{-k^-/2}(2^{-k^-}+\langle t\rangle) 2^{-\ga k^-},
\end{split}
\end{equation}
where $l\in\{1,2,3\}$ and $\mathcal{D}\in\mathcal{V}_m^q$, where the inequalities in the first line hold for all pairs $(q,m)\leq (3,3)$, while the inequalities in the last two lines hold only for pairs $(q,n)\leq (2,2)$. We may also assume that $\mathcal{N}_1$ satisfies the bounds (see Lemmas \ref{abc2.2} and \ref{abc2.3})
\begin{equation}\label{abc3.02}
\begin{split}
\|P_k\mathcal{D}\mathcal{N}_1(t)\|_{L^2}&\lesssim \varep_1^22^{k/2}2^{-\widetilde{N}(m)k^++2k^+}\langle t\rangle^{-1+\delta'}\min\{1,2^k\langle t\rangle\},\\
\|P_k\mathcal{D} \mathcal{N}_1(t)\|_{L^\infty}&\lesssim \varep_1^22^{k}\langle t\rangle^{-2+2\delta'}2^{-N(m+1)k^++3k^+}\min\{1,2^{k}\langle t\rangle\},\\
\|P_k(x_l\mathcal{D}\mathcal{N}_1)(t)\|_{L^2}&\lesssim \varep_1^22^{k/2}2^{-N(m+1)k^+-5k^+}\langle t\rangle^{\delta'}.
\end{split}
\end{equation}

The bounds \eqref{sac3.00} follow easily from \eqref{abc3.01}--\eqref{abc3.02}, using $L^2\times L^\infty$ estimates similar to those in Lemma \ref{Ga}, with the higher frequency factor placed in $L^2$ and the lower frequency factor in $L^\infty$. The proof of \eqref{sac4.00} also follows from \eqref{abc3.01}--\eqref{abc3.02}. The case $2^k\leq \langle t\rangle^{-1}$ follows by $L^2$ estimates as before. If $2^k\geq \langle t\rangle^{-1}$ then we first combine the weight $x_l$ with the higher frequency and place it in $L^2$. Using \eqref{abc3.01}--\eqref{abc3.02}, we have
\begin{equation*}
\begin{split}
\sum_{(k_1,k_2)\in\mathcal{X}_k,\,k_1\leq k_2+10}\|P_{k_1}\mathcal{L}_1G_{\geq 1}(t)\|_{L^\infty}\cdot\|P_{k_2}(x_l\mathcal{L}_2 \mathcal{N}_1)(t)\|_{L^2}\lesssim \varep_1^3\langle t\rangle^{-1+3\delta'}2^{-N(n+1)k^+-4k^+}
\end{split}
\end{equation*}
and
\begin{equation}\label{sho3}
\sum_{(k_1,k_2)\in\mathcal{X}_k,\,k_2\leq k_1-10}\|P_{k_1}(x_l\mathcal{L}_1G_{\geq 1})(t)\|_{L^2}\cdot\|P_{k_2}(\mathcal{L}_2 \mathcal{N}_1)(t)\|_{L^\infty}\lesssim \varep_1^3\langle t\rangle^{-1+3\delta'}2^{-N(n_1+1)k^+}2^{k^-/2}.
\end{equation}
These bounds suffice in most cases, except when ($n_1=n$ and $2^k\geq \langle t\rangle^{1/8}$), because of the loss of derivative in \eqref{sho3}. In this case, however, we combine the weight $x_l$ with the lower frequency and use the estimate
\begin{equation*}
\sum_{(k_1,k_2)\in\mathcal{X}_k,\,k_2\leq k_1-10}\|P_{k_1}(\mathcal{L}_1 G_{\geq 1})(t)\|_{L^2}\cdot 2^{3k_2/2}\|P_{k_2}(x_l\mathcal{L}_2 \mathcal{N}_1)(t)\|_{L^2}\lesssim \varep_1^3\langle t\rangle^{3\delta'}2^{-N(n)k^+},
\end{equation*}
which suffices. 

The bounds \eqref{plk28.1} also follow from \eqref{abc3.01}--\eqref{abc3.02}. We use first $L^2\times L^\infty$ estimates with the lower frequency placed in $L^\infty$. This suffices in most cases, except when all the vector-fields apply to the low frequency factor (so we do not have $L^\infty$ control of this factor). In this case, however, we can reverse the two norms, and still prove \eqref{plk28.1} after some simple analysis.
\end{proof}

\begin{lemma}\label{dtv8.0}
Assume that $t\in[0,T]$, $\mathcal{L}\in\mathcal{V}_n^q$, $n\leq 3$, and $k\in\mathbb{Z}$. Then
\begin{equation}\label{abc31.00}
\|P_k(\mathcal{L}\mathcal{N}^{\psi,\geq 3})(t)\|_{L^2}\lesssim\varep_1^2\langle t\rangle^{-0.6}2^{-\widetilde{N}(n)k^++6k^+}\min\{2^{k},\langle t\rangle^{-1}\}.
\end{equation}
Moreover, if $l\in\{1,2,3\}$ and $n\leq 2$ then we also have the  bounds
\begin{equation}\label{abc31.300}
\begin{split}
\|P_k(\mathcal{L}\mathcal{N}^{\psi,\geq 3})(t)\|_{L^\infty}&\lesssim\varep_1^22^{k/2}\langle t\rangle^{-1.6}2^{-N(n+1)k^++5k^+}\min\{2^{k},\langle t\rangle^{-1}\},\\
\|P_k(x_l\mathcal{L}\mathcal{N}^{\psi,\geq 3})(t)\|_{L^2}&\lesssim\varep_1^2\langle t\rangle^{-0.6}2^{-N(n+1)k^+-2k^+}.
\end{split}
\end{equation}
\end{lemma}

\begin{proof} We examine the identity \eqref{gb10}. It suffices to prove that if $\Psi'\in\{\psi,\partial_\mu\psi,\partial_\mu\partial_\nu\psi:\,\mu,\nu\in\{0,1,2,3\},\,(\mu,\nu)\neq (0,0)\}$ and $G_{\geq 2}$ is a sum of the form $\sum_{d\geq 2}a_dg_d$ for some functions $g_d\in\mathcal{G}_d$ (see \eqref{gb0.1}) and some coefficients $a_d\in\mathbb{R}$ with $|a_d|\leq C^d$, then
\begin{equation}\label{gb11.00}
\|P_k\{\mathcal{L}(G_{\geq 2}\cdot \Psi')\}(t)\|_{L^2}\lesssim \varep_1^2\langle t\rangle^{-0.6}2^{-\widetilde{N}(n)k^++6k^+}\min\{2^{k},\langle t\rangle^{-1}\}
\end{equation}
and, assuming that $n\leq 2$,
\begin{equation}\label{gb12.00}
\begin{split}
\|P_k\{\mathcal{L}(G_{\geq 2}\cdot \Psi')\}(t)\|_{L^\infty}&\lesssim\varep_1^22^{k/2}\langle t\rangle^{-1.6}2^{-N(n+1)k^++5k^+}\min\{2^{k},\langle t\rangle^{-1}\},\\
\|P_k\{x_l\mathcal{L}(G_{\geq 2}\cdot \Psi')\}(t)\|_{L^2}&\lesssim\varep_1^2\langle t\rangle^{-0.6}2^{-N(n+1)k^+-2k^+}.
\end{split}
\end{equation}

Concerning the functions $G_{\geq 2}$ and $\Psi'$, we may assume that (see \eqref{gb1})
\begin{equation}\label{abc31.01}
\begin{split}
\|P_k\mathcal{D}G_{\geq 2}(t)\|_{L^2}&\lesssim \varep_1^2\langle t\rangle^{-1+\delta'}2^{-N(m)k^+}2^{-k/2} \min\{1,2^{k^-}\langle t\rangle\}^{1-\delta'},\\
\|P_k\mathcal{D} G_{\geq 2}(t)\|_{L^\infty}&\lesssim \varep_1^2\langle t\rangle^{-2+\delta'}2^{-N(m+1)k^++3k^+}\min\{1,2^{k^-}\langle t\rangle\}^{2-\delta'},\\
\|P_k(x_l\mathcal{D}G_{\geq 2})(t)\|_{L^2}&\lesssim \varep_1^2(2^{-k^-}+\langle t\rangle)^{\delta'}2^{-N(m+1)k^+}2^{-k^-/2},
\end{split}
\end{equation}
where $l\in\{1,2,3\}$ and $\mathcal{D}\in\mathcal{V}_m^q$, where the inequalities in the first line hold for all pairs $(q,m)\leq (3,3)$, while the inequalities in the last two lines hold only for pairs $(q,m)\leq (2,2)$. Also, as a consequence of \eqref{vcx1}, \eqref{wws2}, and \eqref{cnb1}
\begin{equation}\label{abc31.02}
\begin{split}
\|P_k\mathcal{D}\Psi'(t)\|_{L^2}&\lesssim \varep_1\langle t\rangle^{\delta'}2^{-N(m)k^++k^+},\\
\|P_k\mathcal{D} \Psi'(t)\|_{L^\infty}&\lesssim \varep_1\langle t\rangle^{-1+\delta'}2^{-N(m+1)k^++3k^+}2^{k^-/2}\min\{1,2^{2k^-}\langle t\rangle\},\\
\|P_k(x_l\mathcal{D}\Psi')(t)\|_{L^2}&\lesssim \varep_1\langle t\rangle^{\delta'}(1+2^{k^-}\langle t\rangle)2^{-N(m+1)k^++k^+},
\end{split}
\end{equation}

As before, we start with $L^2$ estimates,
\begin{equation}\label{abc31.03}
\begin{split}
\|P_k(\mathcal{L}_1G_{\geq 2}\cdot\mathcal{L}_2 \Psi')(t)\|_{L^2}&\lesssim\sum_{(k_1,k_2)\in\mathcal{X}_k}2^{3\min(k,k_1,k_2)/2}\|P_{k_1}\mathcal{L}_1G_{\geq 2}(t)\|_{L^2}\|P_{k_2}\mathcal{L}_2 \Psi'(t)\|_{L^2}\\
&\lesssim \varep_1^32^{k^-}\langle t\rangle^{-1+3\delta'} 2^{-N(n)k^++2k^+},
\end{split}
\end{equation}
provided that $\mathcal{L}_1\in\mathcal{V}_{n_1}^{q_1}$, $\mathcal{L}_2\in\mathcal{V}_{n_2}^{q_2}$, $n_1+n_2=n$. This suffices to prove \eqref{gb11.00} if $2^k\lesssim \langle t\rangle^{-2/3}$ or if $2^k\gtrsim \langle t\rangle^{1/5}$. Moreover, if $2^k\in[\langle t\rangle^{-2/3},\langle t\rangle^{1/5}]$ then we use $L^2\times L^\infty$ estimates and \eqref{abc31.01}--\eqref{abc31.02}, with the lower frequency factor estimated in $L^\infty$ and the higher frequency factor estimated in $L^2$. This suffices in most cases, except when $n=3$ and all the vector-fields apply to the low frequency factor. In this case, however, we can reverse the two norms, and still prove the desired conclusion \eqref{gb11.00}. 

The $L^\infty$ bounds in \eqref{gb12.00} follow from the $L^2$ bounds \eqref{abc31.03} if $2^k\lesssim \langle t\rangle^{-0.9}$. On the other hand, if $2^k\geq\langle t\rangle^{-0.9}$ then we just use the $L^2$ and the $L^\infty$ bounds in \eqref{abc31.01}--\eqref{abc31.02}. 

The proof of the weighted $L^2$ bounds in \eqref{gb12.00} is similar. The case $2^k\leq \langle t\rangle^{-1}$ follows by $L^2$ estimates as before. As in the proof of Lemma \ref{abc2.1}, if $2^k\geq \langle t\rangle^{-1}$ then we first combine the weight $x_l$ with the higher frequency and place it in $L^2$. This gives the desired bounds \eqref{gb12.00} in most cases, except when $2^k\geq \langle t\rangle^{1/10}$ and ($n_1=n$ or $n_2=n$), because of the loss of derivative. In these cases, however, we combine the weight $x_l$ with the lower frequency factor and use the $L^2$ estimates in the first and third lines of \eqref{abc31.01}--\eqref{abc31.02}. The desired bounds \eqref{gb12.00} follow in these cases as well.
\end{proof}

\subsection{Additional low frequency bounds} We prove now some additional linear bounds on the solutions $P_kU^{\LL h}$ when $k$ is very small. These bounds are important in some of the energy estimates in the next sections, when the vector-fields hit very low frequency factors.

\begin{lemma}\label{TrivialVF}

Assume that $t\in[0,T]$, $k\le 0$, and $J+k\geq 0$. If $\mathcal{L}=\Gamma_a\mathcal{L}^\prime$, $a\in\{1,2,3\}$, $\mathcal{L}^\prime\in\mathcal{V}^q_n$, $n\leq 2$, and $h\in\{h_{\al\be}\}$ then
\begin{equation}\label{TrivialVF2}
(2^{k^-}\langle t\rangle)^\ga 2^{-k/2}\Vert \varphi_{\leq J}(x)\cdot P_kU^{\mathcal{L}h}(t)\Vert_{L^2}\lesssim \varepsilon_12^k(2^J+\langle t\rangle)\cdot\langle t\rangle^{H(q,n)\delta+2\delta}.
\end{equation}
In addition, if $|\alpha|\leq 3$, $\mathcal{L}^\prime\in\mathcal{V}^q_n$, $n+|\alpha|\leq 3$, then
\begin{equation}\label{TrivialVF3}
(2^{k^-}\langle t\rangle)^\ga 2^{-k/2}\Vert \varphi_{\leq J}(x)\cdot P_kU^{\Omega^{\alpha}\LL'h}(t)\Vert_{L^2}\lesssim \varepsilon_1 2^{\vert\alpha\vert(J+k)} \cdot\langle t\rangle^{H(q,n)\delta}.
\end{equation}
\end{lemma}

\begin{proof} These bounds, which should be compared with \eqref{vcx1}, are used only when $2^k\lesssim \langle t\rangle^{-1+\delta'}$. Recall the formula $U^{\mathcal{L}h}(t)=(\partial_t-i\Lambda_{wa})(\mathcal{L}h)(t)$. To prove \eqref{TrivialVF2} we write
\begin{equation}\label{TrivialVF5}
U^{\mathcal{L}h}(t)=(\partial_t-i\Lambda_{wa})(\Gamma_a\mathcal{L}'h)(t)=\Gamma_aU^{\mathcal{L'}h}(t)+[\partial_t-i\Lambda_{wa},\Gamma_a](\mathcal{L}'h)(t).
\end{equation}
The commutator can be bounded easily, without needing spacial localization,
\begin{equation*}
\big\|P_k[\partial_t-i\Lambda_{wa},\Gamma_a](\mathcal{L}'h)(t)\big\|_{L^2}\lesssim \big\|P_kU^{\mathcal{L}'h}(t)\big\|_{L^2}\lesssim \varep_1(2^{k^-}\langle t\rangle)^{-\ga} 2^{k/2}\cdot\langle t\rangle^{H(q,n)\delta},
\end{equation*}
see \eqref{vcx1}. For the main term we write $\Gamma_aU^{\mathcal{L'}h}(t)=t\partial_aU^{\mathcal{L'}h}(t)+x_a\partial_tU^{\mathcal{L'}h}(t)$, and estimate
\begin{equation*}
\big\|P_kt\partial_aU^{\mathcal{L'}h}(t)\big\|_{L^2}\lesssim \langle t\rangle2^k\big\|P_kU^{\mathcal{L}'h}(t)\big\|_{L^2}\lesssim \varep_1(2^{k^-}\langle t\rangle)^{-\ga} 2^{k/2}\langle t\rangle^{H(q,n)\delta}\cdot\langle t\rangle 2^k.
\end{equation*}
Using also the identity $\partial_tU^{\mathcal{L'}h}(t)=-i\Lambda_{wa}U^{\mathcal{L'}h}(t)+\mathcal{L'}\mathcal{N}^{h}(t)$
and spacial localization,
\begin{equation}\label{TrivialVF9}
\begin{split}
\Vert \varphi_{\leq J}(x)\cdot P_k(x_a\partial_t&U^{\mathcal{L}'h}(t))\Vert_{L^2}\lesssim 2^J\Vert  P'_k(\partial_tU^{\mathcal{L}'h}(t))\Vert_{L^2}\\
&\lesssim 2^{J+k}\Vert  P'_kU^{\mathcal{L}'h}(t)\Vert_{L^2}+2^{J}\Vert  P'_k\mathcal{L}'\mathcal{N}^{h}(t)\Vert_{L^2}\\
&\lesssim 2^{J+k}\varep_1(2^{k^-}\langle t\rangle)^{-\ga} 2^{k/2}\langle t\rangle^{H(q,n)\delta}+2^J\varep_1^22^{3k/2}\langle t\rangle^{H(q,n)\delta+3\delta/2}
\end{split}
\end{equation}
where $P'_k=P_{[k-2,k+2]}$ and we use \eqref{wer4.0} in the last line. The desired conclusion \eqref{TrivialVF2} follows from the last three bounds.

The proof of \eqref{TrivialVF3} is easier, since $P_kU^{\Omega^{\alpha}\LL'h}(t)=\Omega^{\alpha}P_kU^{\LL'h}(t)$, and each $\Omega$ vector-field generates a factor of $2^{J+k}$, as in \eqref{TrivialVF9} above.
\end{proof}

\subsection{Additional bounds on some quadratic nonlinearities} We will also need some slightly stronger bounds on some of the components of the nonlinearities $\mathcal{L}\mathcal{N}^h_{\al\be}$ when $\mathcal{L}\in\mathcal{V}_1^q$.

\begin{lemma}\label{Onev1}
For any $k\in\mathbb{Z}$, $t\in[0,T]$, $\mathcal{L}\in\mathcal{V}^q_n$, $n\leq 1$, and $\al,\be\in\{0,1,2,3\}$
\begin{equation}\label{onev2}
\|P_k(\LL\mathcal{KG}^2_{\al\be})(t)\|_{L^2}+\|P_k(\LL\mathcal{S}^2_{\al\be})(t)\|_{L^2}\lesssim\varep_1^22^{k/2}\langle t\rangle^{-1+H(q,n)\delta+3\delta/2}2^{-N(n)k^++k^+}.
\end{equation}
Moreover, if $2^k\leq\langle t\rangle^{-1/10}$ then
\begin{equation}\label{onev3}
\|P_k(\LL\mathcal{Q}^2_{\al\be})(t)\|_{L^2}\lesssim\varep_1^22^{k/2}\langle t\rangle^{-1+H(q,n)\delta+3\delta/2}.
\end{equation}
\end{lemma}

\begin{proof} The quadratic nonlinearities $\mathcal{KG}^2_{\al\be}$, $\mathcal{S}^2_{\al\be}$, and $\mathcal{Q}^2_{\al\be}$ are defined in \eqref{sac1.5}--\eqref{zaq21}. The point of the lemma is the slightly better estimates in terms of powers of $\langle t\rangle$, when at most one vector-field acts on nonlinearities. We prove the desired bounds in several steps.

{\bf{Step 1.}} We consider first bilinear interactions of the metric components.  For later use we prove slightly stronger frequency-localized estimates. Assume $k,k_1,k_2\in\mathbb{Z}$, $t\in[0,T]$, $\mathcal{L}_2\in\mathcal{V}^{q_2}_{n_2}$, $n_2\leq 1$, $h_1,h_2\in\{h_{\al\be}\}$, and $\iota_1,\iota_2\in\{+,-\}$. Assume that $m$ is a multiplier satisfying $\|\mathcal{F}^{-1}m\|_{L^1(\mathbb{R}^6)}\leq 1$ and define $I_m$ as in \eqref{abc36.1}. Using just $L^2$ estimates we have
\begin{equation}\label{onev4}
\begin{split}
\big\Vert P_kI_m[P_{k_1}U^{h_1,\iota_1}(t), P_{k_2}U^{\mathcal{L}_2h_2,\iota_2}(t)]\big\Vert_{L^2}&\lesssim \varep_1^22^{3\min\{k,k_1,k_2\}/2}\big(\langle t\rangle^22^{k_1^-+k_2^-}\big)^{-\ga}2^{k_1/2+k_2/2}\\
&\times\langle t\rangle^{[H(q_2,n_2)+1]\delta}2^{-N(0)k_1^+-N(n_2)k_2^+}.
\end{split}
\end{equation}

These bounds suffice if $\langle t\rangle \lesssim 1$, or if $2^k\lesssim \langle t\rangle^{-1}$. Moreover, if $\langle t\rangle \gg 1$, $2^k\geq\langle t\rangle^{-1}$, and $k=\min\{k,k_1,k_2\}$ then we also have the estimates (see \eqref{box15.4}--\eqref{box15.5})
\begin{equation}\label{onev5}
\big\Vert P_kI_m[P_{k_1}U^{h_1,\iota_1}(t), P_{k_2}U^{\mathcal{L}_2h_2,\iota_2}(t)]\big\Vert_{L^2}\lesssim \varep_1^22^{k/2}\langle t\rangle^{-1+\delta[H(q_2,n_2)+1]}2^{-N(n_2)k_2^+-5k_2^+}2^{k_2^-/4}.
\end{equation}

On the other hand, if $\langle t\rangle \gg 1$, $2^k\geq\langle t\rangle^{-1}$, and $k_1=\min\{k,k_1,k_2\}$ then $2^k\approx 2^{k_2}$ and
\begin{equation}\label{onev6}
\big\Vert P_kI_m[P_{k_1}U^{h_1,\iota_1}(t), P_{k_2}U^{\mathcal{L}_2h_2,\iota_2}(t)]\big\Vert_{L^2}\lesssim \varep_1^22^{-|k_1|}\langle t\rangle^{-1+\delta'}2^{k/2}2^{-N(n_2)k^+},
\end{equation}
using the bounds \eqref{wws1} and \eqref{vcx1}. If, in addition, $\langle t\rangle^{-8\delta'}\leq 2^{k_1}\leq 2^{k_2}\leq\langle t\rangle^{8\delta'}$ then we let $J_1$ be the largest integer such that $2^{J_1}\le \langle t\rangle 2^{-30}$, decompose $P_{k_1}U^{h_1,\iota_1}(t)=U^{h_1,\iota_1}_{\leq J_1,k_1}(t)+U^{h_1,\iota_1}_{>J_1,k_1}(t)$ as in \eqref{box15}, and use \eqref{vcx1.1}, \eqref{wws1}, together with the stronger bounds \eqref{wws12x} on $\|U^{h_1,\iota_1}_{\leq J_1,k_1}(t)\|_{L^\infty}$, to show that
\begin{equation}\label{onev7}
\big\Vert P_kI_m[P_{k_1}U^{h_1,\iota_1}(t), P_{k_2}U^{\mathcal{L}_2h_2,\iota_2}(t)]\big\Vert_{L^2}\lesssim \varep_1^22^{-|k_1|/2}\langle t\rangle^{-1+\delta[H(q_2,n_2)+1]}2^{k/2}2^{-N(n_2)k^+}.
\end{equation}

Finally, if $\langle t\rangle \gg 1$, $2^k\geq\langle t\rangle^{-1}$, and $k_2=\min\{k,k_1,k_2\}$ then $2^k\approx 2^{k_1}$ and, using \eqref{box5.5},
\begin{equation}\label{onev8}
\big\Vert P_kI_m[P_{k_1}U^{h_1,\iota_1}(t), P_{k_2}U^{\mathcal{L}_2h_2,\iota_2}(t)]\big\Vert_{L^2}\lesssim  \varep_1^2\langle t\rangle^{-1+\delta[H(q_2,n_2)+1]}2^{-|k_2|}2^{|k|/4}2^{-N_0k^++k^+}.
\end{equation}

We examine the identities \eqref{zaq20}--\eqref{zaq21} and estimate
\begin{equation}\label{onev9}
\|P_k(\LL\mathcal{S}^2_{\al\be})(t)\|_{L^2}\lesssim\sum_{k_1,k_2,\iota_1,\iota_2,h_1,h_2,\mathcal{L}_2\in\mathcal{V}^{q}_{n}}\big\Vert P_kI_m[P_{k_1}U^{h_1,\iota_1}(t), P_{k_2}U^{\mathcal{L}_2h_2,\iota_2}(t)]\big\Vert_{L^2}.
\end{equation} 
The desired bounds in \eqref{onev2} follow from \eqref{onev4} if $\langle t\rangle \lesssim 1$ or if $2^k\lesssim \langle t\rangle^{-1}$. On the other hand, if $\langle t\rangle \gg 1$ and $2^k\geq\langle t\rangle^{-1}$ then the $\mathrm{High}\times\mathrm{High}\to\mathrm{Low}$ interactions in \eqref{onev9} can be estimated as claimed using \eqref{onev5}. The  $\mathrm{Low}\times\mathrm{High}\to\mathrm{High}$ interactions in \eqref{onev9} can be estimated using \eqref{onev6}--\eqref{onev7} if $k_1\leq k_2$ or if $n=0$ (thus $\LL_2=Id$); they can also be estimated using \eqref{onev8} if $k_2\leq k_1$ and $n=1$. This completes the proof of the bounds on $\|P_k(\LL\mathcal{S}^2_{\al\be})(t)\|_{L^2}$ in \eqref{onev2}.

Similarly, using the identities \eqref{sac1.2} we estimate
\begin{equation}\label{onev10}
\begin{split}
\|P_k(\LL\mathcal{Q}^2_{\al\be})(t)\|_{L^2}\lesssim&\sum_{k_1,k_2,\iota_1,\iota_2,h_1,h_2,\mathcal{L}_2\in\mathcal{V}^{q}_{n}}\Big\{2^{|k_1-k_2|}\big\Vert P_kI_m[P_{k_1}U^{h_1,\iota_1}(t), P_{k_2}U^{\mathcal{L}_2h_2,\iota_2}(t)]\big\Vert_{L^2}\\
&+\sum_{k_1,k_2,\iota_1,h_1,\mathcal{N}^h_2}2^{-k_1}\big\Vert P_kI_m[P_{k_1}U^{h_1,\iota_1}(t), P_{k_2}\mathcal{N}^{h}_2(t)]\big\Vert_{L^2}\Big\},
\end{split}
\end{equation} 
where $\mathcal{N}^h_2\in\big\{\mathcal{N}^h_{\al\be},\,\al,\be\in\{0,1,2,3\}\big\}$. The terms in the second line of \eqref{onev10} are generated by commuting $\LL$ and derivatives, and then replacing $\partial_0^2 h$ with $\Delta h+\mathcal{N}^h$, as in the proof of Lemma \ref{abc2.3}. The desired bounds \eqref{onev3} follow as before, using \eqref{onev4}, \eqref{onev5}, \eqref{onev6}, and \eqref{onev8} to control the terms in the first line of \eqref{onev10} (recall that $2^k\leq\langle t\rangle^{-1/10}$), and the $L^2$ estimates \eqref{wer4.0} and \eqref{vcx1} to control the terms in the second line.

{\bf{Step 2.}} We consider now bilinear interactions of the Klein-Gordon field, and prove again slightly stronger frequency-localized estimates. Assume $k,k_1,k_2\in\mathbb{Z}$, $t\in[0,T]$, $\mathcal{L}_2\in\mathcal{V}^{q_2}_{n_2}$, $n_2\leq 1$, $\iota_1,\iota_2\in\{+,-\}$, and $m$ satisfies $\|\mathcal{F}^{-1}m\|_{L^1(\mathbb{R}^6)}\leq 1$. Using $L^2\times L^\infty$ estimates (with the higher frequency in $L^2$), and the bounds \eqref{wws2}, \eqref{vcx1}, and \eqref{cnb2}, we have
\begin{equation}\label{onev12}
\begin{split}
\big\Vert P_kI_m[P_{k_1}U^{\psi,\iota_1}(t), P_{k_2}U^{\mathcal{L}_2\psi,\iota_2}(t)]\big\Vert_{L^2}&\lesssim \varep_1^2\langle t\rangle^{-1+\delta'}2^{\min\{k_1^-,k_2^-\}/2}2^{\max\{k_1^-,k_2^-\}}2^{-N(n_2)\max\{k_1^+,k_2^+\}}.
\end{split}
\end{equation}

Assume that $2^{k_1^-}\geq \langle t\rangle^{-1/2}2^{20}$. Then we decompose $P_{k_1}U^{\psi,\iota_1}(t)=U^{\psi,\iota_1}_{\leq J_1,k_1}(t)+U^{\psi,\iota_1}_{>J_1,k_1}(t)$ as in \eqref{on11.3}--\eqref{on11.36}, where $J_1$ is the largest integer satisfying $2^{J_1}\leq 2^{k_1^--20}\langle t\rangle$. Using \eqref{wws13x}, \eqref{wws2}, and \eqref{vcx1}--\eqref{vcx1.1} we estimate
\begin{equation}\label{onev13}
\begin{split}
\big\Vert P_k&I_m[P_{k_1}U^{\psi,\iota_1}(t), P_{k_2}U^{\mathcal{L}_2\psi,\iota_2}(t)]\big\Vert_{L^2}\\
&\lesssim \|U^{\psi,\iota_1}_{\leq J_1,k_1}(t)\|_{L^\infty}\|P_{k_2}U^{\mathcal{L}_2\psi,\iota_2}(t)\|_{L^2}+\|U^{\psi,\iota_1}_{>J_1,k_1}(t)\|_{L^2}\|P_{k_2}U^{\mathcal{L}_2\psi,\iota_2}(t)\|_{L^\infty}\\
&\lesssim \varep_1^2\langle t\rangle^{-3/2+H(q_2,n_2)\delta}2^{-k_1^-/2+\kappa k_1^-/20}2^{-N(1)k_1^+}2^{-5k_2^+}.
\end{split}
\end{equation}

Finally, since $\|P_{k_1}U^{\psi,\iota_1}(t)\|_{L^2}\lesssim 2^{k_1^-+\kappa k_1^-}2^{-N_0k_1^++2k_1^+}$ (see \eqref{vcx1.2}) we can use just $L^2$ bounds to estimate
\begin{equation}\label{onev14}
\begin{split}
\big\Vert P_kI_m[P_{k_1}U^{\psi,\iota_1}(t),&P_{k_2}U^{\mathcal{L}_2\psi,\iota_2}(t)]\big\Vert_{L^2}\\
&\lesssim  \varep_1^22^{3\min\{k^-,k_1^-,k_2^-\}/2}2^{k_1^-+\kappa k_1^-}2^{-N(1)k_1^+}2^{-N(n_2)k_2^+}\langle t\rangle^{H(q_2,n_2)\delta}.
\end{split}
\end{equation}

We can now complete the proof of \eqref{onev2}. We estimate first 
\begin{equation}\label{onev15}
\|P_k(\LL\mathcal{KG}^2_{\al\be})(t)\|_{L^2}\lesssim\sum_{k_1,k_2,\iota_1,\iota_2,\mathcal{L}_2\in\mathcal{V}^{q}_{n}}\big\Vert P_kI_m[P_{k_1}U^{\psi,\iota_1}(t), P_{k_2}U^{\mathcal{L}_2\psi,\iota_2}(t)]\big\Vert_{L^2}.
\end{equation} 
The bounds claimed in \eqref{onev2} follow using only \eqref{onev12} if $2^k\gtrsim\langle t\rangle^{\delta'}$. If $2^k\in[1,\langle t\rangle^{\delta'}]$ then the desired bounds follow using \eqref{onev13} for the contribution of the pairs $(k_1,k_2)\in\mathcal{X}_k$ with $2^{k_1^-}\geq\langle t\rangle^{-1/2}2^{20}$, and \eqref{onev14} for the other pairs. 

Assume now that $2^k\leq 1$. First we use \eqref{onev14} to bound the contribution of the pairs $(k_1,k_2)$ for which $2^{k_1^-+k^-}\lesssim\langle t\rangle^{-1}$. On the other hand if $2^{k_1^-+k^-}\geq \langle t\rangle^{-1}2^{100}$ then we use \eqref{onev13} if $2^{k_1^-}\geq \langle t\rangle^{-1/2}2^{20}$ and \eqref{onev12} for the remaining pairs with $2^{k_1^-}\leq \langle t\rangle^{-1/2}2^{20}\leq 2^{k^-}$. This completes the proof of \eqref{onev2}.
\end{proof}

\section{Decompositions of the main nonlinearities}\label{DecompositionNonlin}

\subsection{The metric tensors $\LL h$} Recall the variables $F,\uF,\rho,\omega_j,\Omega_j,\vartheta_{jk}$ defined in \eqref{zaq2}, and the identities \eqref{zaq5} that recover the metric components. To include vector-fields we need to expand this definition. More precisely, given a symmetric covariant $2$-tensor $H_{\al\be}$ we define
\begin{equation}\label{zaq2l}
\begin{split}
F=F[H]&:=(1/2)[H_{00}+R_jR_kH_{jk}],\\
\uF=\uF[H]&:=(1/2)[H_{00}-R_jR_kH_{jk}],\\
\rho=\rho[H]&:=R_jH_{0j},\\
\omega_j=\omega_j[H]&:=\in_{jkl}R_kH_{0l},\\
\Omega_j=\Omega_j[H]&:=\in_{jkl}R_kR_mH_{lm},\\
\vartheta_{jk}=\vartheta_{jk}[H]&:=\in_{jmp}\in_{knq}R_mR_nH_{pq}.
\end{split}
\end{equation} 
We often apply this decomposition to $H_{\al\be}=\mathcal{L}h_{\al\be}$, $\mathcal{L}\in\mathcal{V}^3_3$. As a general rule, here and in other places, we first apply all the vector-fields to the components $h_{\al\be}$, and then take the Riesz transforms. So we define the variables
\begin{equation}\label{zaq2l.1}
G^{\mathcal{L}}:=G[\mathcal{L}h],\qquad U^{G^\LL,\pm}(t):=\partial_tG^{\LL}(t)\mp i|\nabla|G^{\LL}(t),\qquad G\in\{F,\uF,\rho,\omega_j,\Omega_j,\vartheta_{jk}\}.
\end{equation}
Compare with the definitions \eqref{variables4}-\eqref{variables4L}.

As in \eqref{zaq5}, we can recover the tensor $H$ according to the formulas 
\begin{equation}\label{zaq5l}
\begin{split}
&H_{00}=F[H]+\uF[H],\\
&H_{0j}=-R_j\rho[H]+\in_{jkl}R_k\omega_l[H],\\
&H_{jk}=R_jR_k(F[H]-\uF[H])-(\in_{klm}R_j+\in_{jlm}R_k)R_l\Omega_m[H]+\in_{jpm}\in_{kqn}R_pR_q\vartheta_{mn}[H].
\end{split}
\end{equation} 

The harmonic gauge condition \eqref{asum5} gives
\begin{equation}\label{zaq5.1}
m^{\al\be}\partial_\al h_{\be\mu}-\frac{1}{2}m^{\al\be}\partial_\mu h_{\al\be}=E^{\geq 2}_\mu:=-g_{\geq 1}^{\al\be}\partial_\al h_{\be\mu}+\frac{1}{2}g_{\geq 1}^{\al\be}\partial_\mu h_{\al\be}.
\end{equation}
These identities and the identities \eqref{zaq11.1} can be used to derive a set of elliptic equations for the variables $F,\uF,\rho,\omega_j,\Omega_j,\vartheta_{jk}$. More precisely, let $R_0:=|\nabla|^{-1}\partial_t$ and 
\begin{equation}\label{zaq21.1}
\tau=\tau[H]:=(1/2)[\delta_{jk}H_{jk}+R_jR_kH_{jk}]=-(1/2)\delta_{jk}\vartheta[H]_{jk}.
\end{equation}
We show below that the variables $\rho$ and $\Omega_j$ can be expressed in terms of the other variables, up to quadratic errors, while the variables $\tau$ are in fact quadratic. More precisely:

\begin{lemma}\label{LemWaveGauge} Assume that $\mathcal{L}\in\mathcal{V}_3^3$ and define 
\begin{equation}\label{zaq7l}
E^{com}_{\mathcal{L},\mu}:=m^{\al\be}[\partial_\al,\LL]h_{\be\mu}-\frac{1}{2}m^{\al\be}[\partial_\mu,\LL]h_{\al\be}.
\end{equation}
Then the  variables $\rho^\LL,\Omega_j^\LL$ (defined in \eqref{zaq2l.1}) satisfy the elliptic-type identities
\begin{equation}\label{zaq11.2}
\begin{split}
\rho^\LL&=R_0\uF^\LL+R_0\tau^\LL+|\nabla|^{-1}\LL E_0^{\geq 2}+|\nabla|^{-1}E^{com}_{\mathcal{L},0},\\
\Omega_j^\LL&=R_0\omega_j^\LL+|\nabla|^{-1}\in_{jlk}R_l\LL E_k^{\geq 2}+|\nabla|^{-1}\in_{jlk}R_lE^{com}_{\mathcal{L},k},
\end{split}
\end{equation}
where $\tau^\LL:=\tau[\LL h]$. The variables $\tau^\LL$ satisfies the identities
\begin{equation}\label{zaq11}
\begin{split}
2|\nabla|^2\tau^\LL&=\partial_\al\LL E_\al^{\geq 2}+\partial_\al E^{com}_{\mathcal{L},\al}+\underline{F}[\LL \mathcal{N}^h]+\tau[\LL \mathcal{N}^h],\\
2|\nabla|\partial_0\tau^\LL&=-|\nabla|\LL E_0^{\geq 2}+R_k\partial_0\LL E_k^{\geq 2}-|\nabla|E^{com}_{\mathcal{L},0}+\partial_0R_kE^{com}_{\mathcal{L},k}+\rho[\LL \mathcal{N}^h].
\end{split}
\end{equation}
\end{lemma}

\begin{proof} As a consequence of \eqref{zaq5.1} we have
\begin{equation}\label{zaq5.1l}
m^{\al\be}\partial_\al \LL h_{\be\mu}-\frac{1}{2}m^{\al\be}\partial_\mu \LL h_{\al\be}=\LL E^{\geq 2}_\mu+m^{\al\be}[\partial_\al,\LL] h_{\be\mu}-\frac{1}{2}m^{\al\be}[\partial_\mu, \LL] h_{\al\be}=\LL E^{\geq 2}_\mu+E^{com}_{\mathcal{L},\mu}.
\end{equation}
Notice that $m^{\al\be} \LL h_{\al\be}=2\tau^\LL-2F^\LL$. Therefore \eqref{zaq5.1l} with $\mu=0$ gives
\begin{equation*}
-\partial_0(F^\LL+\uF^\LL)+|\nabla|\rho^\LL-\partial_0(\tau^\LL-F^\LL)=\LL E^{\geq 2}_0+E^{com}_{\mathcal{L},0}.
\end{equation*}
This simplifies to 
\begin{equation}\label{zaq8.1}
-\partial_0\uF^\LL+|\nabla|\rho^\LL-\partial_0\tau^\LL=\LL E^{\geq 2}_0+E^{com}_{\mathcal{L},0}.
\end{equation}
Similarly, using \eqref{zaq5.1l} with $\mu=k\in\{1,2,3\}$ gives
\begin{equation*}
-\partial_0\LL h_{0k}+\partial_j\LL h_{jk}-\frac{1}{2}\partial_k(m^{\al\be} \LL h_{\al\be})=\LL E^{\geq 2}_k+E^{com}_{\mathcal{L},k}.
\end{equation*}
Taking the divergence and the curl, and using \eqref{zaq2l}, this gives
\begin{equation}\label{zaq8.2}
-\partial_0\rho^\LL-|\nabla|\uF^\LL+|\nabla|\tau^\LL=R_k\LL E^{\geq 2}_k+R_kE^{com}_{\mathcal{L},k},
\end{equation}
and
\begin{equation}\label{zaq8.3}
-\partial_0\omega_j^\LL+|\nabla|\Omega_j^\LL=\in_{jlk}R_l\LL E^{\geq 2}_k+\in_{jlk}R_lE^{com}_{\mathcal{L},k}.
\end{equation}
The identities \eqref{zaq11.2} follow from \eqref{zaq8.1} and \eqref{zaq8.3}.

We can now use \eqref{zaq8.1} and \eqref{zaq8.2} to derive the identities \eqref{zaq11} for $\tau^\LL$. Indeed, applying $\partial_0$ to the first equation and $|\nabla|$ to the second equation and adding up we have
\begin{equation*}
-(\partial_0^2+|\nabla|^2)\underline{F}^\LL-(\partial_0^2+|\nabla|^2)\tau^\LL+2|\nabla|^2\tau^\LL=\partial_0\LL E_0^{\geq 2}+|\nabla|R_kE_k^{\geq 2}+\partial_0E^{com}_{\mathcal{L},0}+|\nabla|R_kE^{com}_{\mathcal{L},k}.
\end{equation*}
Similarly, applying $-|\nabla|$ to the first equation and $\partial_0$ to the second equation and adding up,
\begin{equation*}
-(\partial_0^2+|\nabla|^2)\rho^\LL+2|\nabla|\partial_0\tau^\LL=-|\nabla|\LL E_0^{\geq 2}+R_k\partial_0\LL E_k^{\geq 2}-|\nabla|E^{com}_{\mathcal{L},0}+\partial_0R_kE^{com}_{\mathcal{L},k}.
\end{equation*}
The desired identities \eqref{zaq11} follow since $(\partial_0^2+|\nabla|^2)G^\LL=G[\LL \mathcal{N}^h]$, $G\in\{\uF,\tau,\rho\}$.
\end{proof}

\subsection{Null structures} In the analysis of the wave nonlinearities $\mathcal{N}^{h}$ it is important to identify null components, for which we prove better estimates. We start with some definitions.

\begin{definition}\label{nullStr1}
We define two classes of null multipliers $\mathcal{M}_+^{null}$ and $\mathcal{M}_-^{null}$,
\begin{equation}\label{mults1}
\begin{split}
\mathcal{M}^{null}_\pm:=\big\{n:(\mathbb{R}^3\setminus 0)^2&\to\mathbb{C}:\,n(x,y)=(x_i/|x|\mp y_i/|y|)m(x,y)\\
&\text{ for some }m\in\mathcal{M}\text{ and }i\in\{1,2,3\}\big\}.
\end{split}
\end{equation}
For any $(q,n)$ with $0\leq q\leq n\leq 3$ we define the set of "semilinear null forms of order $(q,n)$" as the set of finite sums of expressions of the form
\begin{equation}\label{nake1}
I^{null}_{n_{\iota_1\iota_2}}[U^{\LL_1h_1,\iota_1},U^{\LL_2h_2,\iota_2}]:=\frac{1}{8\pi^3}\int_{\mathbb{R}^3}n_{\iota_1\iota_2}(\xi-\eta,\eta)\widehat{U^{\LL_1h_1,\iota_1}}(\xi-\eta)\widehat{U^{\LL_2h_2,\iota_2}}(\eta)\,d\eta,
\end{equation}
where $\iota_1,\iota_2\in\{+,-\}$, $h_1,h_2\in\{h_{\al\be}\}$, $\mathcal{L}_1\in\mathcal{V}_{n_1}^{q_1}$, $\mathcal{L}_2\in\mathcal{V}_{n_2}^{q_2}$, $(q_1,n_1)+(q_2,n_2)\leq (q,n)$ and $n_{\iota_1\iota_2}\in\mathcal{M}^{null}_{\iota_1\iota_2}$. Here, by convention, $++=--=+$ and $+-=-+=-$.
\end{definition}

Notice that our definition of semilinear null forms contains the classical null forms
\begin{equation}\label{nake2.5}
\partial_\al h_1\partial_\be h_2-\partial_\be h_1\partial_\al h_2\qquad\text{ and }\qquad m^{\al\be}\partial_\al h_1\partial_\be h_2,
\end{equation}
for $h_1,h_2\in\{h_{\mu\nu}\}$, $\al,\be,\mu,\nu\in\{0,1,2,3\}$. Indeed, since 
\begin{equation*}
\partial_0h=(1/2)[U^{h,+}+U^{h,-}], \qquad \partial_jh=(i/2)[R_jU^{h,+}-R_jU^{h,-}],
\end{equation*}
for any $h\in\{h_{\al\be}\}$ (see \eqref{variables4L}), we have, for $a,b\in\{1,2,3\}$,
\begin{equation}\label{nake2}
\partial_ah_1\partial_bh_2-\partial_bh_1\partial_ah_2=\sum_{\iota_1,\iota_2\in\{+,-\}}I^{null}_{n^{a,b}_{\iota_1\iota_2}}[U^{h_1,\iota_1},U^{h_2,\iota_2}],\quad n^{a,b}_{\iota}(\theta,\eta):=\frac{\iota}{4}\frac{\theta_a\eta_b-\theta_b\eta_a}{|\theta||\eta|},
\end{equation}
\begin{equation}\label{nake2.1}
\partial_0h_1\partial_bh_2-\partial_bh_1\partial_0h_2=\sum_{\iota_1,\iota_2\in\{+,-\}}\iota_1I^{null}_{n^{0,b}_{\iota_1\iota_2}}[U^{h_1,\iota_1},U^{h_2,\iota_2}],\quad n^{0,b}_{\iota}(\theta,\eta):=\frac{1}{4}\Big[\frac{\theta_b}{|\theta|}-\iota\frac{\eta_b}{|\eta|}\Big],
\end{equation}
\begin{equation}\label{nake2.2}
m^{\al\be}\partial_\al h_1\partial_\be h_2=\sum_{\iota_1,\iota_2\in\{+,-\}}I^{null}_{\widetilde{n}_{\iota_1\iota_2}}[U^{h_1,\iota_1},U^{h_2,\iota_2}],\quad \widetilde{n}_{\iota}(\theta,\eta):=\frac{1}{4}\Big[-1+\iota\frac{\theta\cdot\eta}{|\theta||\eta|}\Big].
\end{equation}
It is easy to verify that the symbols $n^{a,b}_{\iota},n^{0,b}_{\iota},\widetilde{n}_\iota$ are in $\mathcal{M}^{null}_\iota$, as defined in \eqref{mults1} (in fact $n^{a,b}_{\iota}\in \mathcal{M}^{null}_+\cap\mathcal{M}^{null}_-$), therefore the classical null forms in \eqref{nake2.5} are all semilinear null forms of order $(0,0)$. The vector-fields $\LL$ can be incorporated as well, without any difficulty. 

We remark that our definition \eqref{nullStr1} of semilinear null forms is slightly more general because we would like to allow forms expressed in terms of the variables $F,\uF,\rho,\omega_j,\Omega_j,\vartheta_{jk}$, which involve the Riesz transforms. For example, an expression of the form
\begin{equation}\label{zaq8.55}
fg-R_jR_kf\cdot R_jR_kg
\end{equation}
is a semilinear null form of order $(0,0)$ if $f,g\in\{R_1^{a_1}R_2^{a_2}R_3^{a_3}\partial_\al h_{\mu\nu}\}$.

We also need to define suitable classes of cubic error terms. More precisely:

\begin{definition}\label{nake4}
We define two sets of quadratic and higher order expressions $\mathcal{QU}_0\subseteq\mathcal{QU}$
\begin{equation}\label{sho3.9}
\begin{split}
\mathcal{QU}&:=\big\{\partial_\al h_1\partial_\be h_2,\,G_{\geq 1}\cdot\partial_\al h_1\partial_\be h_2,\,\mathcal{KG}_{\al\be}^2,\,G_{\geq 1}\cdot\mathcal{KG}_{\al\be}^2,\,G_{\geq 1}\cdot\partial_\al\partial_\be h_2\big\},\\
\mathcal{QU}_0&:=\big\{\partial_\al h_1\partial_\be h_2,\,G_{\geq 1}\cdot\partial_\al h_1\partial_\be h_2,\,\mathcal{KG}_{\al\be}^2,\,G_{\geq 1}\cdot\mathcal{KG}_{\al\be}^2\big\},\\
\end{split}
\end{equation}
where $\al,\be\in\{0,1,2,3\}$, $h_1,h_2\in\{h_{\mu\nu}\}$, $\mathcal{KG}^2_{\al\be}$ are defined in \eqref{sac1.5}, and $G_{\geq 1}=\sum_{d\geq 1}a_dg_d$ for some functions $g_d\in\mathcal{G}_d$ (see \eqref{gb0.1}) and some coefficients $a_d\in\mathbb{R}$ with $|a_d|\leq C^d$.

A "semilinear cubic remainder of order $(q,n)$" is a finite sum of expressions of the form 
\begin{equation}\label{zaq8.6}
I[|\nabla|^{-1}\LL_1\mathcal{N},\partial_\al\LL_2h_{\mu\nu}]\quad\text {or }\quad I[|\nabla|^{-1}\LL_1\mathcal{N},|\nabla|^{-1}\LL_2\mathcal{N}']\quad\text{ or }\quad I[\LL_1\mathcal{N}_0,\LL_2h_{\mu\nu}]
\end{equation}
where $I$ is as in \eqref{abc36.1}, $\mathcal{N},\mathcal{N'}\in\mathcal{QU}$, $\mathcal{N}_0\in\mathcal{QU}_0$, $\mathcal{L}_i\in\mathcal{V}_{n_i}^{q_i}$, and $(q_1,n_1)+(q_2,n_2)\leq (q,n)$.
\end{definition}

For example, the functions $\mathcal{N}_{\mu\nu}^{h}$, $\partial_\mu E_{\nu}^{\geq 2}, |\nabla|\partial_\mu\tau$ can be written as sums of terms of the form $R^a\mathcal{N}$ for $\mathcal{N}\in\mathcal{QU}$, where $R^a=R_1^{a_1}R_2^{a_2}R_3^{a_3}$ (see \eqref{zaq5.1} and \eqref{zaq11}). 

We notice that functions in $\mathcal{QU}$ satisfy quadratic-type bounds similar to $\mathcal{N}_{\al\be}^{h}$:

\begin{lemma} \label{SHO4} If $\mathcal{N}\in\mathcal{QU}$, $\mathcal{L}\in\mathcal{V}^q_n$, $n\leq 3$, $t\in[0,T]$, and $k\in\mathbb{Z}$ then
\begin{equation}\label{sho4}
\|P_k\{\mathcal{L}\mathcal{N}\}(t)\|_{L^2}\lesssim \varep_1^22^{k/2}\langle t\rangle^{H(q,n)\delta}2^{-\widetilde{N}(n)k^++7k^+}\min(2^{k}\langle t\rangle^{3\delta/2},\langle t\rangle^{-1+\widetilde{\ell}(n)\delta}).
\end{equation}
Moreover, if $n\leq 2$ and $l\in\{1,2,3\}$ then we also have the bounds
\begin{equation}\label{sho4.2}
\begin{split}
\|P_k\{\mathcal{L}\mathcal{N}\}(t)\|_{L^\infty}&\lesssim\varep_1^22^{k}\langle t\rangle^{-1+4\delta'}2^{-N(n+1)k^++5k^+}\min\{2^{k},\langle t\rangle^{-1}\},\\
\|P_k\{x_l\mathcal{L}\mathcal{N}\}(t)\|_{L^2}&\lesssim \varep_1^22^{k/2}\langle t\rangle^{\delta'}2^{-N(n+1)k^+-2k^+}.
\end{split}
\end{equation}
Finally, if $\mathcal{L}\in\mathcal{V}^q_n$, $n\leq 3$, $t\in[0,T]$, $k\in\mathbb{Z}$, and $2^k\lesssim\langle t\rangle^{-\delta'}$ then
\begin{equation}\label{sho4.3}
\|P_k\{\mathcal{L}\mathcal{N}\}(t)\|_{L^2}\lesssim \varep_1^22^{k/2}\langle t\rangle^{H(q,n)\delta}\langle t\rangle^{-1+\ell(q,n)\delta+\delta/2}.
\end{equation}
\end{lemma}

\begin{proof} The bounds \eqref{sho4}--\eqref{sho4.2} follow from the proofs of Lemmas \ref{abc2.2}, \ref{abc2.3}, \ref{abc2.1}. The bounds \eqref{sho4.3} follow from \eqref{gb17.1} if $\mathcal{N}=\mathcal{KG}^2_{\al\be}$ and from \eqref{sac3.00} if $\mathcal{N}$ is a cubic term. It remains to prove that
\begin{equation*}
\begin{split}
\|P_k\{\mathcal{L}\mathcal{Q}^2_{\al\be}\}(t)\|_{L^2}+\|P_k\{\mathcal{L}\mathcal{S}^2_{\al\be}\}(t)\|_{L^2}\lesssim \varep_1^22^{k/2}\langle t\rangle^{H(q,n)\delta}\langle t\rangle^{-1+\ell(q,n)\delta+\delta/2}
\end{split}
\end{equation*}
if $2^k\lesssim\langle t\rangle^{-\delta'}$, where $\mathcal{Q}^2_{\al\be}$ and $\mathcal{S}^2_{\al\be}$ are as in \eqref{sac1.2}--\eqref{zaq21}. This can be proved as in Lemma \ref{abc2.3}, using Lemma \ref{box1}. Most of the terms gain a factor of $2^{k/2}$, which is more than enough to give the additional time decay. The only exception are the $\mathrm{High}\times\mathrm{High}\to\mathrm{Low}$ interactions which are controlled using \eqref{box2}--\eqref{box3.5}; however, in these interactions we already have the time decay factor $\langle t\rangle^{H(q,n)\delta}\langle t\rangle^{-1+\ell(q,n)\delta+\delta/2}$ as claimed.
\end{proof}

\subsection{Energy disposable nonlinearities} To prove energy estimate in the next section we need to bound the contribution of space-time  integrals. Many resulting terms can be estimated easily, without normal form analysis, using just decay properties. In this subsection we identify these terms. We start with a definition.

\begin{definition}\label{nake10} Assume that $(q,n)\leq (3,3)$. A function $L:\mathbb{R}^3\times[0,T]\to\mathbb{C}$ will be called "wave-disposable of order $(q,n)$" if, for any $t\in[0,T]$,
\begin{equation}\label{nake11}
\big\||\nabla|^{-1/2}L(t)\big\|_{H^{N(n)}}\lesssim \varep_1^2\langle t\rangle^{-1+H(q,n)\delta-\delta/2}.
\end{equation}
Similarly, a function $L:\mathbb{R}^3\times[0,T]$ will be called "KG-disposable of order $(q,n)$" if
\begin{equation}\label{nake11.1}
\big\|L(t)\big\|_{H^{N(n)}}\lesssim \varep_1^2\langle t\rangle^{-1+H(q,n)\delta-\delta/2}.
\end{equation}
\end{definition}

We show first that most cubic and higher order terms in $\LL\mathcal{N}_{\al\be}^{h}$ and $\LL\mathcal{N}^{\psi}$ are disposable. 

\begin{lemma}\label{nake12}
If $(q,n)\leq (3,3)$ then any semilinear cubic remainder of order $(q,n)$ (see Definition \ref{nake4}) is wave-disposable. In addition, if $h_1,h_2\in \{h_{\mu\nu}\}$ then terms of the form 
\begin{equation}\label{nake13.3}
I[|\nabla|^{-1}\LL_1 (G_{\geq 1}\partial_\rho h_1),\LL_2\psi]
\end{equation}
are KG-disposable of order $(q,n)$ for all $\LL_1\in\mathcal{V}_{n_1}^{q_1}$, $\LL_2\in\mathcal{V}_{n_2}^{q_2}$ with $(q_1,n_1)+(q_2,n_2)\leq (q,n)$.

Moreover, if $n_2<n$, then terms of the form
\begin{equation}\label{nake13}
I[|\nabla|^{-1}\LL_1 (G_{\geq 1}\partial_\rho h_1),\partial_\al\partial_\be\LL_2h_2]
\end{equation}
are wave-disposable of order $(q,n)$, while terms of the form
\begin{equation}\label{nake13.1}
I[|\nabla|^{-1}\LL_1 (G_{\geq 1}\partial_\rho h_1),\partial_\al\partial_\be\LL_2\psi]
\end{equation}
are KG-disposable of order $(q,n)$.
\end{lemma}

\begin{proof} Time decay is not an issue in this lemma, since we are considering cubic and higher order terms, but we need to be careful to avoid possible derivative loss. For any $(q,n)$ with $(q,n)\leq (3,3)$ and $t\in[0,T]$ we define the frequency envelopes $\{b_k\}_{k\in\mathbb{Z}}=\{b_k(q,n;t)\}_{k\in\mathbb{Z}}$ by
\begin{equation}\label{nake15}
\begin{split}
b^0_k(q,n;t)&:=\sup_{\mathcal{K}\in\mathcal{V}_n^q,\,\al,\be\in\{0,1,2,3\}}\langle t\rangle^{-H(q,n)\delta}\big\{\big\|P_k\{(\langle t\rangle|\nabla|_{\leq 1})^{\ga}|\nabla|^{-1/2}U^{\mathcal{K}h_{\al\be}}\}(t)\big\|_{H^{N(n)}}\\
&\qquad\qquad\qquad\qquad\qquad\qquad\quad\,\,\,\,+\|P_kU^{\mathcal{K}\psi}(t)\|_{H^{N(n)}}\big\},\\
b_k(q,n;t)&:=\varep_12^{-\ga|k|/4}+\sum _{k'\in\mathbb{Z}}2^{-\ga|k-k'|/4}b^0_{k'}(q,n;t).
\end{split}
\end{equation}
In view of the bootstrap assumption \eqref{bootstrap2.1} we have
\begin{equation}\label{nake15.1}
\sum_{k\in\mathbb{Z}}(b_k(q,n;t))^2\lesssim \varep_1^2\quad \text{ and }\quad b_k(q,n;t)\leq b_{k'}(q,n;t)2^{\ga|k-k'|/4}\text{ for any }k,k'\in\mathbb{Z}.
\end{equation}
The main point of the definition is that we have the slightly better $L^2$ bounds
\begin{equation}\label{nake15.2}
\begin{split}
\big\|P_kU^{\mathcal{K}h_{\al\be}}(t)\big\|_{L^2}&\lesssim b_k(q,n;t)\langle t\rangle^{H(q,n)\delta}2^{k/2}2^{-N(n)k^+}(\langle t\rangle 2^{k^-})^{-\ga},\\
\big\|P_kU^{\mathcal{K}\psi}(t)\big\|_{L^2}&\lesssim b_k(q,n;t)\langle t\rangle^{H(q,n)\delta}2^{-N(n)k^+},
\end{split}
\end{equation}
for any $k\in\mathbb{Z}$, $\mathcal{K}\in\mathcal{V}^q_n$, and $\al,\be\in\{0,1,2,3\}$ (compare with \eqref{vcx1}).

To prove the conclusions we need two more quadratic bounds: if $\mathcal{N}_0\in\mathcal{QU}_0$ and $\mathcal{N}\in\mathcal{QU}$ (see \eqref{sho3.9}), $k\geq 0$, $t\in[0,T]$, and $\mathcal{K}\in\mathcal{V}_{n'}^{q'}$ then
\begin{equation}\label{nake15.5}
\begin{split}
\|P_k\{\mathcal{K}\mathcal{N}_0\}(t)\|_{L^2}&\lesssim \varep_1b_k(q',n';t)2^{k/2}2^{-N(n')k}\langle t\rangle^{-1+\delta'/2},\\
\|P_k\{\mathcal{K}\mathcal{N}\}(t)\|_{L^2}&\lesssim \varep_1b_k(q',n';t)2^{3k/2}2^{-N(n')k}\langle t\rangle^{-1+\delta'/2}.
\end{split}
\end{equation}
Notice that these bounds are improvements over the general bounds \eqref{sho4} when $k\geq 0$. 

{\bf{Step 1.}} We assume first that the bounds \eqref{nake15.5} hold and show how to prove the conclusions of the lemma. Recall also the bounds \eqref{nyu1}. In view of the definitions it suffices to show that
\begin{equation}\label{nake15.6}
\begin{split}
\|P_kI[|\nabla|^{-1}\LL_1\mathcal{N},U^{\LL_2h_2,\iota_2}](t)\|_{L^2}&\lesssim \varep_1b_k(q,n;t)2^{k/2}2^{-N(n)k^+}\langle t\rangle^{-1},\\
\|P_kI[|\nabla|^{-1}\LL_1\mathcal{N},|\nabla|^{-1}\LL_2\mathcal{N}'](t)\|_{L^2}&\lesssim \varep_1b_k(q,n;t)2^{k/2}2^{-N(n)k^+}\langle t\rangle^{-1},\\
\|P_kI[\LL_1\mathcal{N}_0,|\nabla|^{-1}U^{\LL_2h_2,\iota_2}](t)\|_{L^2}&\lesssim \varep_1b_k(q,n;t)2^{k/2}2^{-N(n)k^+}\langle t\rangle^{-1},
\end{split}
\end{equation}
and
\begin{equation}\label{nake15.66}
\|P_kI[|\nabla|^{-1}\LL_1 (G_{\geq 1}\partial_\rho h_1),\langle\nabla\rangle^{-1}U^{\LL_2\psi,\iota_2}](t)\|_{L^2}\lesssim \varep_1^22^{-\delta|k|}2^{-N(n)k^+}\langle t\rangle^{-1},
\end{equation}
for any $k\in\mathbb{Z}$, $t\in[0,T]$, $\mathcal{N},\mathcal{N'}\in\mathcal{QU}$, $\mathcal{N}_0\in\mathcal{QU}_0$, $\LL_1\in\mathcal{V}_{n_1}^{q_1}$, $\LL_2\in\mathcal{V}_{n_2}^{q_2}$, $(q_1,n_1)+(q_2,n_2)\leq (q,n)$. Moreover, if $n_2<n$ and $\mathcal{N}^{h}\in\{\mathcal{N}^h_{\al\be}\}$, then we also have to prove the bounds
\begin{equation}\label{nake15.7}
\begin{split}
\|P_kI[|\nabla|^{-1}\LL_1 (G_{\geq 1}\partial_\rho h_1),|\nabla|U^{\LL_2h_2,\iota_2}](t)\|_{L^2}&\lesssim \varep_1^22^{-\delta|k|}2^{k/2}2^{-N(n)k^+}\langle t\rangle^{-1},\\
\|P_kI[|\nabla|^{-1}\LL_1 (G_{\geq 1}\partial_\rho h_1),\LL_2\mathcal{N}^{h}](t)\|_{L^2}&\lesssim \varep_1^22^{-\delta|k|}2^{k/2}2^{-N(n)k^+}\langle t\rangle^{-1},
\end{split}
\end{equation}
and 
\begin{equation}\label{nake15.8}
\begin{split}
\|P_kI[|\nabla|^{-1}\LL_1 (G_{\geq 1}\partial_\rho h_1),\langle\nabla\rangle U^{\LL_2\psi,\iota_2}](t)\|_{L^2}&\lesssim \varep_1^22^{-\delta|k|}2^{-N(n)k^+}\langle t\rangle^{-1},\\
\|P_kI[|\nabla|^{-1}\LL_1 (G_{\geq 1}\partial_\rho h_1),\LL\mathcal{N}^{\psi}](t)\|_{L^2}&\lesssim \varep_1^22^{-\delta|k|}2^{-N(n)k^+}\langle t\rangle^{-1},
\end{split}
\end{equation}
in order to show that the expressions in \eqref{nake13}--\eqref{nake13.1} are also disposable.

The proofs of \eqref{nake15.6}--\eqref{nake15.8} rely on $L^2\times L^\infty$ estimates, as in Lemmas \ref{Ga}-\ref{Ga4}. In most cases we place the high frequency factor in $L^2$ and the low frequency factor in $L^\infty$, except when the low frequency factor carries all the 3 vector-fields. We provide all the details only for the proof of the harder estimates in the first line of \eqref{nake15.6}, which require the frequency envelopes. 

The functions $|\nabla|^{-1}\LL_1\mathcal{N}$ and $U^{\LL_2h_2,\iota_2}$ satisfy the bounds
\begin{equation}\label{nake16.0}
\begin{split}
\|P_{k_1}(|\nabla|^{-1}\mathcal{L}_1\mathcal{N})(t)\|_{L^2}&\lesssim \varep_1b_{k_1}(q,n;t)2^{-N(n_1)k_1^++k_1^+/2}\min(2^{k_1},\langle t\rangle^{-1})^{1/2-\delta'},\\
\|P_{k_1}(|\nabla|^{-1}\mathcal{L}_1\mathcal{N})(t)\|_{L^\infty}&\lesssim \varep_1^2\langle t\rangle^{-1+4\delta'}2^{-N(n_1+1)k_1^++5k_1^+}\min(2^{k_1},\langle t\rangle^{-1}),
\end{split}
\end{equation}
see \eqref{nake15.5}, \eqref{sho4}, and \eqref{sho4.2}, and
\begin{equation}\label{nake16.1}
\begin{split}
\|P_{k_2}(U^{\LL_2h_2,\iota_2})(t)\|_{L^2}&\lesssim b_{k_2}(q,n;t)2^{-N(n_2)k_2^++k_2^+/2}2^{ k_2^-/2-\delta k_2^-}\langle t\rangle^{\delta'},\\
\|P_{k_2}(U^{\LL_2h_2,\iota_2})(t)\|_{L^\infty}&\lesssim \varep_1\langle t\rangle^{-1+\delta'}2^{k_2^-}2^{-N(n_2+1)k_2^++2k_2^+},
\end{split}
\end{equation}
see \eqref{nake15.2} and \eqref{wws1}. As before, the $L^\infty$ bounds in the second lines of \eqref{nake16.0} and \eqref{nake16.1} hold only if $n_1\leq 2$ and $n_2\leq 2$ respectively. We estimate first
\begin{equation*}
\begin{split}
\|P_kI[|\nabla|^{-1}\LL_1\mathcal{N},U^{\LL_2h_2,\iota_2}](t)\|_{L^2}&\lesssim \sum_{(k_1,k_2)\in\mathcal{X}_k}2^{3k/2}\|P_{k_1}(|\nabla|^{-1}\mathcal{L}_1\mathcal{N})(t)\|_{L^2}\|P_{k_2}(U^{\LL_2h_2,\iota_2})(t)\|_{L^2}\\
&\lesssim \varep_1^3\langle t\rangle^{-1/2+4\delta'}2^{3k/2},
\end{split}
\end{equation*}
which suffices to prove \eqref{nake15.6} if $2^k\lesssim \langle t\rangle^{-0.51}$. Moreover, if $2^k\geq \langle t\rangle^{-0.51}$ and $n_1,n_2\in[0,2]$ then we estimate, using also \eqref{nake15.1},
\begin{equation*}
\begin{split}
\|P_kI[|\nabla|^{-1}\LL_1\mathcal{N},U^{\LL_2h_2,\iota_2}](t)\|_{L^2}&\lesssim S_1+S_2,
\end{split}
\end{equation*}
where
\begin{equation*}
\begin{split}
S_1:=\sum_{(k_1,k_2)\in\mathcal{X}_k,\,k_1\geq k-8}&\|P_{k_1}(|\nabla|^{-1}\mathcal{L}_1\mathcal{N})(t)\|_{L^2}\|P_{k_2}(U^{\LL_2h_2,\iota_2})(t)\|_{L^\infty}\\
&\lesssim \varep_1^2b_{k}(q,n;t)2^{-N(n_1)k^++k^+/2}\langle t\rangle^{-3/2+4\delta'},
\end{split}
\end{equation*}
\begin{equation*}
\begin{split}
S_2:=\sum_{(k_1,k_2)\in\mathcal{X}_k,\,k_1\leq k-8}&\|P_{k_1}(|\nabla|^{-1}\mathcal{L}_1\mathcal{N})(t)\|_{L^\infty}\|P_{k_2}(U^{\LL_2h_2,\iota_2})(t)\|_{L^2}\\
&\lesssim \varep_1^2b_{k}(q,n;t)2^{-N(n_2)k^++k^+/2}\langle t\rangle^{-5/3}.
\end{split}
\end{equation*}
Finally, if $n_1=3$ (thus $n_2=0$, $n=3$), then we estimate
\begin{equation*}
\begin{split}
\|P_kI[|\nabla|^{-1}\LL_1\mathcal{N},U^{\LL_2h_2,\iota_2}](t)\|_{L^2}&\lesssim \sum_{(k_1,k_2)\in\mathcal{X}_k}\|P_{k_1}(|\nabla|^{-1}\mathcal{L}_1\mathcal{N})(t)\|_{L^2}\|P_{k_2}(U^{\LL_2h_2,\iota_2})(t)\|_{L^\infty}\\
&\lesssim \varep_1^2b_{k}(q,n;t)2^{-N(n)k^++k^+/2}\langle t\rangle^{-3/2+4\delta'},
\end{split}
\end{equation*}
while if  $n_2=3$ (thus $n_1=0$, $n=3$), then
\begin{equation*}
\begin{split}
\|P_kI[|\nabla|^{-1}\LL_1\mathcal{N},U^{\LL_2h_2,\iota_2}](t)\|_{L^2}&\lesssim \sum_{(k_1,k_2)\in\mathcal{X}_k}\|P_{k_1}(|\nabla|^{-1}\mathcal{L}_1\mathcal{N})(t)\|_{L^\infty}\|P_{k_2}(U^{\LL_2h_2,\iota_2})(t)\|_{L^2}\\
&\lesssim \varep_1^2b_{k}(q,n;t)2^{-N(n)k^++k^+/2}\langle t\rangle^{-5/3}.
\end{split}
\end{equation*}
These bounds suffice to prove \eqref{nake15.6} when $2^k\geq \langle t\rangle^{-0.51}$.

The estimates \eqref{nake15.7}--\eqref{nake15.8} are slightly easier, because we do not need to carry the frequency envelopes. The functions $|\nabla|^{-1}\LL_1 (G_{\geq 1}\partial_\rho h_1)$ satisfy the bounds in Lemma \ref{Ga4}, while the functions $\mathcal{N}^{\LL_2h_2}$ and $\mathcal{N}^{\LL_2\psi}$ satisfy the bounds in Proposition \ref{plk1}. Since $n_2<n$ there is no derivative loss, and the estimates \eqref{nake15.7}--\eqref{nake15.8} follow in the same way as \eqref{nake15.6}.

{\bf{Step 2.}} We prove now the bounds \eqref{nake15.5}. We examine the definition \eqref{sho3.9}. If $\mathcal{N}^1_0:=\partial_\al h_1\cdot\partial_\be h_2$ and $\mathcal{N}^2_0:=G_{\geq 1}\partial_\al h_1\cdot\partial_\be h_2$, $\al,\be\in\{0,1,2,3\}$, $h_1,h_2\in\{h_{\mu\nu}\}$  then 
\begin{equation}\label{nake23}
\|P_k\{\mathcal{K}\mathcal{N}_0^1\}(t)\|_{L^2}+\|P_k\{\mathcal{K}\mathcal{N}_0^2\}(t)\|_{L^2}\lesssim \varep_1b_k(q',n';t)2^{k/2}2^{-N(n')k}\langle t\rangle^{-1+\delta'},
\end{equation}
for any $k\geq 0$, $t\in[0,T]$, and $\mathcal{K}\in\mathcal{V}_{n'}^{q'}$. This is easy to see using just $L^2\times L^\infty$ estimates as before, and the bounds \eqref{nake15.2}, \eqref{wws1}, and \eqref{abc3.01}. Similarly, if $\mathcal{N}^3_0:=\mathcal{KG}^2_{\al\be}$ (see \eqref{sac1.5}) for some $\al,\be\in\{0,1,2,3\}$ and $\mathcal{N}_0^4:=G_{\geq 1}\cdot \mathcal{KG}^2_{\al\be}$ then
\begin{equation}\label{nake23.2}
\|P_k\{\mathcal{K}\mathcal{N}_0^3\}(t)\|_{L^2}+\|P_k\{\mathcal{K}\mathcal{N}_0^4\}(t)\|_{L^2}\lesssim \varep_1^22^{-N(n')k}\langle t\rangle^{-1+\delta'},
\end{equation}
for any $k\geq 0$, $t\in[0,T]$, and $\mathcal{K}\in\mathcal{V}_{n'}^{q'}$. Finally, if $(\al,\be)\neq (0,0)$, then
\begin{equation}\label{nake23.3}
\|P_k\{\mathcal{K}(G_{\geq 1}\cdot \partial_\al\partial_\be h_2)\}(t)\|_{L^2}\lesssim \varep_1b_k(q',n';t)2^{-N(n')k+3k/2}\langle t\rangle^{-1+\delta'},
\end{equation}
for any $k\geq 0$, $t\in[0,T]$, and $\mathcal{K}\in\mathcal{V}_{n'}^{q'}$, using again $L^2\times L^\infty$ estimates as before (and Proposition \ref{plk1} to bound the commutator term $[\mathcal{K},\partial_\al\partial_\be] h_2$). As a consequence of the last three bounds and the identitites \eqref{sac1}, the metric nonlinearities $\mathcal{N}_{\mu\nu}^{h}$ satisfy bounds similar to \eqref{nake23.3}, 
\begin{equation*}
\|P_k\{\mathcal{K}\mathcal{N}_{\mu\nu}^h\}(t)\|_{L^2}\lesssim \varep_1b_k(q',n',t)2^{-N(n')k+3k/2}\langle t\rangle^{-1+\delta'}.
\end{equation*}
Therefore, we can use the equation $\partial_0^2h=\Delta h+\mathcal{N}^h$ to prove the bounds \eqref{nake23.3} for $\al=\be=0$ as well. This completes the proof of \eqref{nake15.5}.
\end{proof}

We show now that all the quadratic terms arising as commutators are also energy disposable. 
\begin{lemma}\label{nake12.5}
Assume that $(q,n)\leq (3,3)$, $\LL_1\in\mathcal{V}_{n_1}^{q_1}$, $\LL_2\in\mathcal{V}_{n_2}^{q_2}$, $(q_1,n_1)+(q_2,n_2)\leq (q,n-1)$. If  $h_1,h_2\in\{h_{\mu\nu}\}$, $\rho,\al,\be\in\{0,1,2,3\}$ then quadratic terms of the form
\begin{equation}\label{nake13.6}
I[U^{\LL_1\psi,\iota_1},U^{\LL_2\psi,\iota_2}]\quad\text{ or }\quad I[R_\rho\LL_1 h_1,\partial_\al\partial_\be\LL_2h_2]\quad\text{ or }\quad I[\partial_\rho\LL_1 h_1,\partial_\be\LL_2h_2]
\end{equation}
are wave-disposable of order $(q,n)$,  where $\iota_1,\iota_2\in\{+,-\}$. Moreover, terms of the form
\begin{equation}\label{nake13.7}
I[R_{\rho}\LL_1 h_1,\partial_\al\partial_\be\LL_2\psi]
\end{equation}
are KG-disposable of order $(q,n)$.
\end{lemma}

Terms such as those in \eqref{nake13.6} and \eqref{nake13.7} will be called "wave (respectively KG) commutator remainders" of order $(q,n)$.

\begin{proof} Derivative loss is not an issue in this lemma, since $n_1+n_2\leq n-1$, but we need to be careful with the time decay.  We show first for any $k,k_1,k_2\in\mathbb{Z}$ and $t\in[0,T]$
\begin{equation}\label{nake18}
\begin{split}
2^{-k/2}\|P_kI&[P_{k_1}U^{\LL_1\psi,\iota_1},P_{k_2}U^{\LL_2\psi,\iota_2}](t)\|_{L^2}\\
&\lesssim \varep_1^2\langle t\rangle ^{-1+\delta[H(q,n-1)+\ell(q,n-1)+1]}2^{-N(n-1)k^+}2^{-\ga(|k|+|k_1|+|k_2|)/4}.
\end{split}
\end{equation}
Indeed, these bounds follow from \eqref{box32} and \eqref{SuperlinearH1} if $2^{\min\{k,k_1,k_2\}}\lesssim \langle t\rangle^{-1}$ (see also the definitions \eqref{AcceptableLoss}). On the other hand, if  $2^{\min \{k,k_1,k_2\}}\geq \langle t\rangle^{-1}$ then \eqref{nake18} follows from the bounds \eqref{box33}. Since 
\begin{equation}\label{nake18.1}
H(q,n-1)+\ell(q,n-1)+12\leq H(q,n),
\end{equation}
the bounds \eqref{nake18} suffice to show that $I[U^{\LL_1\psi,\iota_1},U^{\LL_2\psi,\iota_2}]$ is wave-disposable. 

Moreover, we also have
\begin{equation}\label{nake18.01}
\begin{split}
2^{-k/2}2^{|k_1-k_2|}\|P_kI&[P_{k_1}U^{\LL_1h_1,\iota_1},P_{k_2}U^{\LL_2h_2,\iota_2}](t)\|_{L^2}\\
&\lesssim \varep_1^2\langle t\rangle ^{-1+\delta[H(q,n)-4]}2^{-N(n)k^+-2k^+}2^{-\ga(|k|+|k_1|+|k_2|)/4},
\end{split}
\end{equation}
for any $k,k_1,k_2\in\mathbb{Z}$, $t\in[0,T]$, $h_1,h_2\in\{h_{\al\be}\}$, $\iota_1,\iota_2\in\{+,-\}$. These bounds follow from \eqref{box2} if $2^{\min\{k,k_1,k_2\}}\lesssim \langle t\rangle^{-1}$ and from \eqref{box3}--\eqref{box5.6} (see also \eqref{SuperlinearH2} and use \eqref{box5.5} instead of \eqref{box4} when $n_1=n_2=0$) if $2^{\min \{k,k_1,k_2\}}\geq \langle t\rangle^{-1}$. Thus terms of the form $I[R_\rho\LL_1 h_1,\partial_\al\partial_\be\LL_2h_2]$ and $I[\partial_\rho\LL_1 h_1,\partial_\mu\LL_2h_2]$ are also wave-disposable of order $(q,n)$ (in the case $(\al,\be)= (0,0)$ we replace first $\partial_0^2\LL_2 h_2$ with $\Delta \LL_2h_2+\LL_2\mathcal{N}^h$ and use \eqref{wer4.0} and \eqref{plk2} to bound the nonlinear contribution).

Finally, we can use Lemma \ref{box80} in a similar way to show that 
\begin{equation}\label{nake18.03}
\begin{split}
2^{k_2^+-k_1}\|P_kI&[P_{k_1}U^{\LL_1h_1,\iota_1},P_{k_2}U^{\LL_2\psi,\iota_2}](t)\|_{L^2}\\
&\lesssim \varep_1^2\langle t\rangle ^{-1+\delta[H(q,n)-4]}2^{-N(n)k^+-2k^+}2^{-\ga(|k|+|k_1|+|k_2|)/4},
\end{split}
\end{equation}
for any $k,k_1,k_2\in\mathbb{Z}$, $t\in[0,T]$, $h_1\in\{h_{\al\be}\}$, $\iota_1,\iota_2\in\{+,-\}$, thus expressions of the form $I[R_\rho\LL_1 h_1,\partial_\al\partial_\be\LL_2\psi]$ are KG-disposable of order $(q,n)$. This completes the proof.
\end{proof}

\subsection{The main decomposition} Given three operators $\mathcal{L}_1=\Gamma^{a'}\Omega^{b'}$, $\mathcal{L}_2=\Gamma^{a''}\Omega^{b''}$, $\mathcal{L}=\Gamma^{a}\Omega^{b}$ (see \eqref{qaz3}), we say that $\LL_1+\LL_2=\LL$ if $a'+a''=a$ and $b'+b''=b$. Therefore
\begin{equation}\label{sumL}
\LL(fg)=\sum_{\LL_1+\LL_2=\LL}c_{\LL_1,\LL_2}\LL_1f\cdot\LL_2g,
\end{equation}
for some coefficients $c_{\LL_1,\LL_2}\in[0,\infty)$. We are now ready to prove an important proposition concerning the decomposition of the nonlinearities $\mathcal{N}^{\LL h}_{\al\be}$ and $\mathcal{N}^{\LL\psi}$. 

\begin{proposition}\label{gooddec}
If $\al,\be\in\{0,1,2,3\}$, $(q,n)\leq (3,3)$, and $\mathcal{L}\in \mathcal{V}_n^q$ then
\begin{equation}\label{BoxWave}
\LL\mathcal{N}^{h}_{\al\be}=\sum_{\mu,\nu\in\{0,1,2,3\}}\widetilde{g}^{\mu\nu}_{\geq 1}\partial_\mu\partial_\nu (\LL h_{\al\be})+\mathcal{Q}_{wa}^{\LL}(h_{\al\be})+\mathcal{S}_{\al\be}^{\LL,1}+\mathcal{S}_{\al\be}^{\LL,2}+\mathcal{KG}^{\LL}_{\al\be}+\mathcal{R}^{\LL h}_{\al\be}
\end{equation}
and
\begin{equation}\label{BoxKG}
\LL\mathcal{N}^{\psi}=\sum_{\mu,\nu\in\{0,1,2,3\}}\widetilde{g}^{\mu\nu}_{\geq 1}\partial_\mu\partial_\nu (\LL\psi)+h_{00}\LL\psi+\mathcal{Q}_{kg}^{\LL}(\psi)+\mathcal{R}^{\LL\psi},
\end{equation}
where the remainders $\mathcal{R}^{\LL h}_{\al\be}$ and $\mathcal{R}^{\LL\psi}$ are wave-disposable (respectively KG-disposable) of order $(q,n)$, and the reduced metric components $\widetilde{g}^{\mu\nu}_{\geq 1}$ are defined by
\begin{equation}\label{BoxWave1}
\widetilde{g}^{00}_{\geq 1}:=0,\qquad \widetilde{g}^{0j}_{\geq 1}:=(1-g_{\geq 1}^{00})^{-1}g_{\geq 1}^{0j},\qquad \widetilde{g}^{jk}_{\geq 1}:=(1-g_{\geq 1}^{00})^{-1}\big[g_{\geq 1}^{jk} +g_{\geq 1}^{00}\delta^{jk}\big].
\end{equation}

$\bullet\,\,$ The terms $\mathcal{Q}_{wa}^{\LL}(h_{\al\be})$ and $\mathcal{Q}_{kg}^{\LL}(\psi)$ are given by 
\begin{equation}\label{BoxWave1.2}
\mathcal{Q}^{\LL}_{wa}(h_{\al\be}):=\sum_{G\in\{F,\uF,\omega_n,\vartheta_{mn}\}}\sum_{\iota_1,\iota_2\in\{+,-\}}\sum_{\LL_1+\LL_2=\LL,\,\LL_2\neq \LL}c_{\LL_1,\LL_2}I^{null}_{\mathfrak{q}^{G,wa}_{\iota_1\iota_2}}[|\nabla|^{-1}U^{G^{\LL_1},\iota_1},|\nabla|U^{\LL_2h_{\al\be},\iota_2}],
\end{equation}
\begin{equation}\label{BoxWave1.3}
\mathcal{Q}^{\LL}_{kg}(\psi):=\sum_{G\in\{F,\uF,\omega_n,\vartheta_{mn}\}}\sum_{\iota_1,\iota_2\in\{+,-\}} \sum_{\LL_1+\LL_2=\LL,\,\LL_2\neq \LL}c_{\LL_1,\LL_2}I_{\mathfrak{q}^{G,kg}_{\iota_1\iota_2}}[|\nabla|^{-1}U^{G^{\LL_1},\iota_1},\langle\nabla\rangle U^{\LL_2\psi,\iota_2}],
\end{equation}
The multipliers $\mathfrak{q}^{F,wa}_{\iota_1\iota_2},\,\mathfrak{q}^{\uF,wa}_{\iota_1\iota_2},\,\mathfrak{q}^{\omega_n,wa}_{\iota_1\iota_2},\,\mathfrak{q}^{\vartheta_{mn},wa}_{\iota_1\iota_2}\in\mathcal{M}_{\iota_1\iota_2}^0$ and $\mathfrak{q}^{F,kg}_{\iota_1\iota_2},\,\mathfrak{q}^{\uF,kg}_{\iota_1\iota_2},\,\mathfrak{q}^{\omega_n,kg}_{\iota_1\iota_2},\,\mathfrak{q}^{\vartheta_{mn},kg}_{\iota_1\iota_2}\in\mathcal{M}$ are given explicitly in \eqref{BoxWave7}-\eqref{BoxWave7.5}.

$\bullet\,\,$ The semilinear terms $\mathcal{S}_{\al\be}^{\LL,1}$, $\mathcal{S}_{\al\be}^{\LL,2}$, and $\mathcal{KG}^{\LL}_{\al\be}$ are given by 
\begin{equation}\label{BoxWave2.3}
\mathcal{S}_{\al\be}^{\LL,1}:=\sum_{h_1,h_2\in\{h_{\mu\nu}\}}\sum_{\iota_1,\iota_2\in\{+,-\}} \sum_{\LL_1+\LL_2=\LL}I^{null}_{\mathfrak{n}_{\iota_1\iota_2}}[U^{\LL_1h_1,\iota_1},U^{\LL_2h_2,\iota_2}],
\end{equation}
\begin{equation}\label{BoxWave2.4}
\mathcal{S}^{\LL,2}_{\al\be}:=\sum_{\LL_1+\LL_2=\LL}(1/2)c_{\LL_1,\LL_2}R_pR_q\partial_\al\va^{\LL_1}_{mn}\cdot R_pR_q\partial_\be\va^{\LL_2}_{mn},
\end{equation}
\begin{equation}\label{BoxWave2.2}
\mathcal{KG}_{\al\be}^{\LL}:=\sum_{\LL_1+\LL_2=\LL}c_{\LL_1,\LL_2}(2\partial_\al\LL_1\psi\cdot\partial_\be\LL_2\psi+m_{\al\be}\LL_1\psi\LL_2\psi),
\end{equation}
where $\mathfrak{n}_{\iota_1\iota_2}=\mathfrak{n}_{\iota_1\iota_2}(\LL_1,\LL_2,h_1,h_2)$ are null multipliers in $\mathcal{M}^0_{\iota_1\iota_2}$ (see \eqref{mults1}).
\end{proposition}

\begin{remark} (1) The nonlinearities $\LL\mathcal{N}^{h}_{\al\be}$ contain five types of components:

(i) The top order terms $\widetilde{g}^{\mu\nu}_{\geq 1}\partial_\mu\partial_\nu (\LL h_{\al\be})$ which lead to derivative loss in energy estimates and normal form analysis;

(ii) The terms $\mathcal{Q}^{\LL}_{wa}(h_{\al\be})$, which are sums of bilinear interactions of the metric components, with null multipliers $\mathfrak{q}_{\iota_1\iota_2}^{G,wa}$. These interactions are all defined by null multipliers, but are not of the same type as the semilinear null forms. The issue is the additional anti-derivative on the first factor $|\nabla|^{-1}U^{G^{\LL_1},\pm}$, which leads to significant difficulties at very low frequencies, particularly when these first factors carry all the vector-fields $\LL_1=\LL$;

(iii) Generic semilinear null terms $\mathcal{S}_{\al\be}^{\LL,1}$;

(iv) The special terms $\mathcal{S}^{\LL,2}_{\al\be}$, which involve non-null bilinear interactions of the "good" metric components $\vartheta$;

(v) The Klein-Gordon nonlinearities $\mathcal{KG}_{\al\be}^{\LL}$. These nonlinearities do not have null structure and, in fact, involve the massive field $\psi$ in undifferentiated form.

(2) The Klein-Gordon nonlinearities $\LL\mathcal{N}^{\psi}$ contain two types of quasilinear components, which are somewhat similar to the first two types of metric nonlinearities described above. 
\end{remark}

\begin{proof} {\bf{Step 1.}} We start with the quasilinear components of the metric nonlinearities, which can be written in the form $\sum_{\mu,\nu\in\{0,1,2,3\}}\widetilde{g}_{\geq 1}^{\mu\nu}\partial_\mu\partial_\nu h_{\al\be}$ (see \eqref{sac1} and \eqref{BoxWave1}). Thus
\begin{equation}\label{BoxWave2}
\LL\big\{\widetilde{g}_{\geq 1}^{\mu\nu}\partial_\mu\partial_\nu h_{\al\be}\big\}=\widetilde{g}_{\geq 1}^{\mu\nu}\LL(\partial_\mu\partial_\nu h_{\al\be})+\sum_{\LL_1+\LL_2=\LL,\,\LL_2\neq \LL}c_{\LL_1,\LL_2}\LL_1\widetilde{g}_{\geq 1}^{\mu\nu}\LL_2(\partial_\mu\partial_\nu h_{\al\be}).
\end{equation}
The first term in the right-hand side can be replaced by $\widetilde{g}_{\geq 1}^{\mu\nu}\partial_\mu\partial_\nu (\LL h_{\al\be})$, while the second term can be replaced by $\sum_{\LL_1+\LL_2=\LL,\,\LL_2\neq \LL}c_{\LL_1,\LL_2}\LL_1\widetilde{g}^{\mu\nu}_1\partial_\mu\partial_\nu (\LL_2 h_{\al\be})$ up to wave-disposable errors (due to Lemmas \ref{nake12} and \ref{nake12.5}), where $\widetilde{g}^{\mu\nu}_1$ is the linear part of $\widetilde{g}_{\geq 1}^{\mu\nu}$. If view of the formulas \eqref{zaq22} and \eqref{BoxWave1}, we have
\begin{equation}\label{BoxWave2.9}
\widetilde{g}^{00}_1=0,\qquad \widetilde{g}^{0j}_1=h_{0j},\qquad\widetilde{g}^{jk}_1=-h_{jk}-h_{00}\delta_{jk}.
\end{equation}
Therefore
\begin{equation}\label{BoxWave5}
\begin{split}
\LL\big\{\widetilde{g}^{\mu\nu}_{\geq 1}\partial_\mu\partial_\nu h_{\al\be}\big\}-\widetilde{g}_{\geq 1}^{\mu\nu}\partial_\mu\partial_\nu (\LL h_{\al\be})&=\sum_{\LL_1+\LL_2=\LL,\,\LL_2\neq \LL}c_{\LL_1,\LL_2}\big[2\LL_1h_{0j}\cdot\partial_0\partial_j \LL_2h_{\al\be}\\
&-\LL_1(h_{jk}+h_{00}\delta_{jk})\cdot\partial_j\partial_k \LL_2h_{\al\be}\big]+\mathcal{R}_{\al\be}^{\LL h,1},
\end{split}
\end{equation}
where $\mathcal{R}_{\al\be}^{\LL h,1}$ is wave-disposable of order $(q,n)$. Using \eqref{zaq5l} with $H=\LL_1 h$, we have
\begin{equation*}
\begin{split}
&2\LL_1h_{0j}\cdot\partial_0\partial_j\LL_2h_{\al\be}-\LL_1(h_{jk}+h_{00}\delta_{jk})\cdot\partial_j\partial_k\LL_2 h_{\al\be}\\
&=[-2R_j\rho^{\LL_1}+2\in_{jmn}R_m\omega_n^{\LL_1}]\cdot\partial_j\partial_0\LL_2h_{\al\be}+\big[-\delta_{jk}(F^{\LL_1}+\uF^{\LL_1})-R_jR_k(F^{\LL_1}-\uF^{\LL_1})\\
&+(\in_{jlm}R_k+\in_{klm}R_j)R_l\Omega_m^{\LL_1}-\in_{jpm}\in_{kqn}R_pR_q\va_{mn}^{\LL_1}\big]\cdot\partial_j\partial_k\LL_2h_{\al\be}.
\end{split}
\end{equation*}
We use now the formulas \eqref{zaq11.2}. After reorganizing the terms, the expression above becomes
\begin{equation}\label{BoxWave5.2}
\begin{split}
&-[\delta_{jk}+R_jR_k]F^{\LL_1}\cdot \partial_j\partial_k\LL_2h_{\al\be}-\big\{2R_jR_0\uF^{\LL_1}\cdot\partial_j\partial_0\LL_2h_{\al\be}+(\delta_{jk}-R_jR_k)\uF^{\LL_1}\cdot \partial_j\partial_k\LL_2h_{\al\be}\big\}\\
&+\big\{2\in_{jmn}R_m\omega_n^{\LL_1}\cdot\partial_j\partial_0\LL_2h_{\al\be}+(\in_{jln}R_k+\in_{kln}R_j)R_lR_0\omega_n^{\LL_1}\cdot\partial_j\partial_k\LL_2h_{\al\be}\big\}\\
&-\in_{jpm}\in_{kqn}R_pR_q\va_{mn}^{\LL_1}\cdot\partial_j\partial_k\LL_2h_{\al\be}-2R_jR_0\tau^{\LL_1}\cdot\partial_j\partial_0\LL_2h_{\al\be}+\mathcal{R}_{\al\be}^{\LL h,2}.
\end{split}
\end{equation} 
Here $\mathcal{R}_{\al\be}^{\LL h,2}$ corresponds to the contribution of the error terms in \eqref{zaq11.2}, and is wave-disposable of order $(q,n)$, in view of Lemma \ref{nake12} (see \eqref{nake13}) and Lemma \ref{nake12.5}. The terms $-2R_jR_0\tau^{\LL_1}\cdot\partial_j\partial_0\LL_2h_{\al\be}$ are wave-disposable errors (due to Lemma \ref{ControlSC} below), while the other terms can be re\-written as claimed in \eqref{BoxWave1.2}, using the identities 
\begin{equation*}
\partial_0\LL_2h_{\al\be}=(1/2)(U^{\LL_2h_{\al\be},+}+U^{\LL_2h_{\al\be},-}),\qquad |\nabla|\LL_2h_{\al\be}=(i/2)(U^{\LL_2h_{\al\be},+}-U^{\LL_2h_{\al\be},-}),
\end{equation*}
see \eqref{on5}, and similar identities for $F^{\LL_1},\uF^{\LL_1},\omega_n^{\LL_1},\va_{mn}^{\LL_1}$.  The symbols $\mathfrak{q}^{G,wa}_{\iota_1\iota_2}$ are given by
\begin{equation}\label{BoxWave7}
\begin{split}
\mathfrak{q}^{F,wa}_{\iota_1\iota_2}:=\frac{-\iota_1\iota_2}{4}\Big[1-\frac{(\theta\cdot\eta)^2}{|\theta|^2|\eta|^2}\Big],&\qquad \mathfrak{q}^{\uF,wa}_{\iota_1\iota_2}:=\frac{-\iota_1\iota_2}{4}\Big[1-\iota_1\iota_2\frac{\theta\cdot\eta}{|\theta||\eta|}\Big]^2,\\
\mathfrak{q}^{\omega_n,wa}_{\iota_1\iota_2}:=\frac{-i\iota_1}{2}\frac{\in_{jpn}\theta_p\eta_j}{|\theta||\eta|}\Big[1-\iota_1\iota_2\frac{\theta\cdot\eta}{|\theta||\eta|}\Big],&\qquad \mathfrak{q}^{\va_{mn},wa}_{\iota_1\iota_2}:=\frac{\iota_1\iota_2}{4}\frac{\in_{jpm}\theta_p\eta_j}{|\theta||\eta|}\frac{\in_{kqn}\theta_q\eta_k}{|\theta||\eta|}.
\end{split}
\end{equation}
These multipliers are similar to the classical null multipliers in \eqref{nake2}--\eqref{nake2.2}, thus belong to $\mathcal{M}^0_{\iota_1\iota_2}$ as claimed.

{\bf{Step 2.}} We consider now the Klein-Gordon nonlinearity $\mathcal{N}^\psi$ defined in \eqref{gb10}. The analysis is similar to the analysis of the quasilinear wave nonlinearities. One can first place most of the cubic and higher order terms and the commutator terms in the KG-disposable remainder, due to Lemmas \ref{nake12} and \ref{nake12.5}, thus
\begin{equation*}
\LL\mathcal{N}^\psi=\widetilde{g}_{\geq 1}^{\mu\nu}\LL(\partial_\mu\partial_\nu \psi)+\LL(h_{00}\psi)+\sum_{\LL_1+\LL_2=\LL,\,\LL_2\neq \LL}c_{\LL_1,\LL_2}\LL_1\widetilde{g}^{\mu\nu}_1\LL_2(\partial_\mu\partial_\nu \psi)+\mathcal{R}^{\LL\psi,1},
\end{equation*}
where $\mathcal{R}^{\LL\psi,1}$ is KG-disposable, and the additional term $\LL(h_{00}\psi)$ comes from the term $(1-g^{00}_{\geq 1})^{-1}g^{00}_{\geq 1}\psi$ in \eqref{gb10}. The first two terms in the right-hand side are as claimed in \eqref{BoxKG}, and we can decompose the remaining term as in \eqref{BoxWave5}--\eqref{BoxWave5.2}, with $h_{\al\be}$ replaced by $\psi$. The terms $-2R_jR_0\tau^{\LL_1}\cdot\partial_j\partial_0\LL_2\psi$ are KG-disposable errors (due to Lemma \ref{ControlSC} below), while all the other terms, including $\LL(h_{00}\psi)$,  are accounted for in $\mathcal{Q}^{\LL}_{kg}(\psi)$. The resulting symbols $\mathfrak{q}^{G,kg}_{\iota_1\iota_2}$ can be calculated as before, and are given explicitly by
\begin{equation}\label{BoxWave7.5}
\begin{split}
\mathfrak{q}^{F,kg}_{\iota_1\iota_2}:=\frac{-\iota_1\iota_2}{4}\Big[1-\frac{(\theta\cdot\eta)^2}{|\theta|^2\langle\eta\rangle^2}\Big],&\qquad \mathfrak{q}^{\uF,kg}_{\iota_1\iota_2}:=\frac{-\iota_1\iota_2}{4}\Big[1-\iota_1\iota_2\frac{\theta\cdot\eta}{|\theta|\langle\eta\rangle}\Big]^2,\\
\mathfrak{q}^{\omega_n,kg}_{\iota_1\iota_2}:=\frac{-i\iota_1}{2}\frac{\in_{jpn}\theta_p\eta_j}{|\theta|\langle\eta\rangle}\Big[1-\iota_1\iota_2\frac{\theta\cdot\eta}{|\theta|\langle\eta\rangle}\Big],&\qquad \mathfrak{q}^{\va_{mn},kg}_{\iota_1\iota_2}:=\frac{\iota_1\iota_2}{4}\frac{\in_{jpm}\theta_p\eta_j}{|\theta|\langle\eta\rangle}\frac{\in_{kqn}\theta_q\eta_k}{|\theta|\langle\eta\rangle}.
\end{split}
\end{equation}

{\bf{Step 3.}} Finally, we consider the semilinear terms coming from the last two terms in \eqref{sac1}. The cubic terms and the commutators can be safely included in the wave-disposable remainders, due to Lemmaa \ref{nake12}--\ref{nake12.5}. The Klein-Gordon contributions coming from $\mathcal{KG}^2_{\al\be}$ are included in the terms $\mathcal{KG}^\LL_{\al\be}$ in \eqref{BoxWave2.2}. The null contributions coming from $Q^2_{\al\be}$ are included in the terms $\mathcal{S}^{\LL,1}_{\al\be}$ in \eqref{BoxWave2.3} (see \eqref{nake2}--\eqref{nake2.2}). The contributions of the terms $P^2_{\al\be}$ in \eqref{zaq21} are recovered in the terms $\mathcal{S}^{\LL,1}_{\al\be}$ and $\mathcal{S}^{\LL,2}_{\al\be}$, due to Lemma \ref{Pstructure} below.
\end{proof}

In the analysis of the quasilinear terms in Lemma \ref{gooddec} we used that certain quadratic expression involving $\tau$ are energy disposable. We prove this below:

\begin{lemma}\label{ControlSC}
Assume that $(q,n)\leq (3,3)$, $\LL_1\in\mathcal{V}_{n_1}^{q_1}$, $\LL_2\in\mathcal{V}_{n_2}^{q_2}$, $(q_1,n_1)+(q_2,n_2)\leq (q,n)$, and $n_2<n$ (so $n\geq 1$). If  $h\in\{h_{\mu\nu}\}$ and $\al\in\{0,1,2,3\}$ then quadratic terms of the form
\begin{equation}\label{saur00}
I[R_0\tau^{\mathcal{L}_1},|\nabla|\partial_\alpha\mathcal{L}_2h]\qquad\text{ and }\qquad I[\tau^{\mathcal{L}_1},|\nabla|\partial_\alpha \mathcal{L}_2h]
\end{equation}
are wave-disposable of order $(q,n)$. Similarly, quadratic terms of the form
\begin{equation}\label{saur01}
I[R_0\tau^{\mathcal{L}_1},\langle\nabla\rangle\partial_\alpha\mathcal{L}_2\psi]\qquad\text{ and }\qquad I[\tau^{\mathcal{L}_1},\langle\nabla\rangle\partial_\alpha\mathcal{L}_2\psi]
\end{equation}
are $KG$-disposable of order $(q,n)$.
\end{lemma}

\begin{proof} The main point is that the metric components $\tau^{\LL_1}$ have quadratic character, up to lower order terms, due to the identities \eqref{zaq11}. At low frequencies, however, the resulting quadratic bounds are not effective, and we need to trivialize one-vector field using Lemma \ref{TrivialVF}. 

More precisely, with $\overline{|k|}:=\max\{|k|,|k_1|,|k_2|\}$ it suffices to prove that
\begin{equation}\label{saur1}
2^{k_2}2^{-k/2}\|P_kI[P_{k_1}R_\mu\tau^{\mathcal{L}_1},P_{k_2}U^{\mathcal{L}_2h,\iota_2}](t)\|_{L^2}\lesssim\varep_1^2\langle t\rangle^{-1+H(q,n)\delta-\delta/2}2^{-N(n)k^+}2^{-\ga\overline{|k|}/4}
\end{equation}
and
\begin{equation}\label{saur2}
2^{k_2^+}\|P_kI[P_{k_1}R_\mu\tau^{\mathcal{L}_1},P_{k_2}U^{\mathcal{L}_2\psi,\iota_2}](t)\|_{L^2}\lesssim\varep_1^2\langle t\rangle^{-1+H(q,n)\delta-\delta/2}2^{-N(n)k^+}2^{-\ga\overline{|k|}/4},
\end{equation}
for any $\mu\in\{0,1,2,3\}$, $\iota_2\in\{+,-\}$, $k,k_1,k_2\in\mathbb{Z}$, and $t\in[0,T]$. These bounds follow from \eqref{box2} and \eqref{box82} if $2^{\min\{k,k_1,k_2\}}\lesssim\langle t\rangle^{-1-4\delta}$ or if $2^{\max\{k,k_1,k_2\}}\gtrsim\langle t\rangle^{2.1}$. In the remaining range
\begin{equation}\label{saur3}
\langle t\rangle\geq 2^{\delta^{-1}}\quad\text{ and }\quad 2^k,\,2^{k_1},\,2^{k_2}\in[\langle t\rangle^{-1-4\delta},\langle t\rangle^{2.1}]
\end{equation}
we divide the proof into several steps.

{\bf{Step 1.}} Assume first that \eqref{saur3} holds and $k_2=\min\{k,k_1,k_2\}$. Then
\begin{equation}\label{saur5}
\begin{split}
2^{k_2}2^{-k/2}&\|P_kI[P_{k_1}R_\mu\tau^{\mathcal{L}_1},P_{k_2}U^{\mathcal{L}_2h,\iota_2}](t)\|_{L^2}\lesssim 2^{k_2-k/2}\Vert P_{k_1}R_\mu\tau^{\mathcal{L}_1}(t)\Vert_{L^2}\Vert P_{k_2}U^{\mathcal{L}_2h}(t)\Vert_{L^\infty}\\
&\lesssim\varepsilon_1^22^{-k/2}\langle t\rangle^{-1+\delta'}2^{2k_2^-}2^{-8k_2^+}\cdot 2^{-k_1/2}2^{\ga|k_1|}2^{-N(n_1)k_1^+}
\end{split}
\end{equation}
and
\begin{equation}\label{saur6}
\begin{split}
2^{k_2^+}&\|P_kI[P_{k_1}R_\mu\tau^{\mathcal{L}_1},P_{k_2}U^{\mathcal{L}_2\psi,\iota_2}](t)\|_{L^2}\lesssim 2^{k_2^+}\Vert P_{k_1}R_\mu\tau^{\mathcal{L}_1}(t)\Vert_{L^2}\Vert P_{k_2}U^{\mathcal{L}_2\psi}(t)\Vert_{L^\infty}\\
&\lesssim\varepsilon_1^2\langle t\rangle^{-1+\delta'}2^{k_2^-/2}\min\{1,\langle t\rangle 2^{2k_2^-}\}2^{-8k_2^+}\cdot 2^{-k_1/2}2^{\ga|k_1|}2^{-N(n_1)k_1^+},
\end{split}
\end{equation}
using the $L^\infty$ estimates \eqref{wws1} and \eqref{wws2} on the second factor. The desired bounds \eqref{saur1}--\eqref{saur2} follow unless $\langle t\rangle^{-1/2}\lesssim 2^{k_2}\leq 2^{k_1}\lesssim \langle t\rangle^{2\delta'}$. In this case, however, we use the identities \eqref{zaq11} and replace $R_\mu\tau^{\mathcal{L}_1}$ with lower order linear terms of the forms $R^a|\nabla|^{-2}\partial_\mu\partial_\nu\LL'_1h_{\al\be}$ (coming from the linear commutator terms) and nonlinear terms of the form $R^a|\nabla|^{-2}\LL''_1\mathcal{N}$, where $\mathcal{N}\in\mathcal{QU}$ (see \eqref{sho3.9}), $\LL'_1\in\mathcal{V}^{q_1}_{n_1-1}$, $\LL''_1\in\mathcal{V}^{q_1}_{n_1}$, $\al,\be,\mu,\nu\in\{0,1,2,3\}$, and $R^a=R_1^{a_1}R_2^{a_2}R_3^{a_3}$.

The contributions of the linear commutators $R^a|\nabla|^{-2}\partial_\mu\partial_\nu\LL'_1h_{\al\be}$ are bounded as claimed, due to the more general estimates \eqref{nake18.01} and \eqref{nake18.03}. The contributions of the nonlinear terms $R^a|\nabla|^{-2}\LL''_1\mathcal{N}$ can be bounded using $L^2\times L^\infty$ estimates, as in \eqref{saur5}--\eqref{saur6}, and recalling the $L^2$ bounds \eqref{sho4} and the assumption $\langle t\rangle^{-1/2}\lesssim 2^{k_2}\leq 2^{k_1}\lesssim \langle t\rangle^{2\delta'}$.

{\bf{Step 2.}} Assume now that \eqref{saur3} holds and $k=\min\{k,k_1,k_2\}$. The bounds \eqref{saur1}--\eqref{saur2} follow from \eqref{box3}--\eqref{box3.5} and \eqref{saur6} unless $\langle t\rangle^{-1/2-2\delta'}\lesssim 2^{k_2}\approx 2^{k_1}\lesssim \langle t\rangle^{2\delta'}$. In this case, however, we can again use the identities \eqref{zaq11}. As before, we estimate the contributions of the linear commutators using \eqref{nake18.01} and \eqref{nake18.03}, and the contributions of the nonlinear terms $R^a|\nabla|^{-2}\LL''_1\mathcal{N}$ using the $L^2$ bounds \eqref{sho4}. The desired bounds \eqref{saur1}--\eqref{saur2} follow.

{\bf{Step 3.}} Assume now that $k_1=\min\{k,k_1,k_2\}$ and $n_1<n$. Recall that $n_2<n$, thus $H(q_1,n_1)+H(q_2,n_2)\leq H(q,n)-40$ (see \eqref{SuperlinearH1}). The desired bounds follow from \eqref{box2} and \eqref{box82} if $2^{k_1}\lesssim \langle t\rangle^{-1+39\delta}$ or $2^{k_2}\gtrsim \langle t\rangle$. On the other hand, if $\langle t\rangle^{-1+39\delta}\leq 2^{k_1}\leq 2^{k_2}\leq \langle t\rangle$ then we use the identities \eqref{zaq11} as before, and replace $R_\mu\tau^{\mathcal{L}_1}$ with lower order linear terms of the forms $R^a|\nabla|^{-2}\partial_\mu\partial_\nu\LL'_1h_{\al\be}$ and nonlinear terms $R^a|\nabla|^{-2}\LL''_1\mathcal{N}$. The contributions of the lower order linear terms can be suitably controlled using \eqref{nake18.01}--\eqref{nake18.03}. Moreover, we estimate
\begin{equation}\label{saur8}
\begin{split}
2^{k_2}2^{-k/2}&\|P_kI[P_{k_1}R^a|\nabla|^{-2}\LL''_1\mathcal{N},P_{k_2}U^{\mathcal{L}_2h,\iota_2}](t)\|_{L^2}\\
&\lesssim 2^{k/2}2^{3k_1/2}\Vert P_{k_1}|\nabla|^{-2}\LL''_1\mathcal{N}(t)\Vert_{L^2}\Vert P_{k_2}U^{\mathcal{L}_2h}(t)\Vert_{L^2}\\
&\lesssim\varepsilon_1^22^k\langle t\rangle^{-1+H(q_1,n_1)\delta+36\delta}\langle t\rangle^{H(q_2,n_2)\delta}2^{-N(n_2)k^+}
\end{split}
\end{equation}
and
\begin{equation}\label{saur9}
\begin{split}
2^{k_2^+}&\|P_kI[P_{k_1}R^a|\nabla|^{-2}\LL''_1\mathcal{N},P_{k_2}U^{\mathcal{L}_2\psi,\iota_2}](t)\|_{L^2}\\
&\lesssim 2^{k^+}2^{3k_1/2}\Vert P_{k_1}|\nabla|^{-2}\LL''_1\mathcal{N}(t)\Vert_{L^2}\Vert P_{k_2}U^{\mathcal{L}_2\psi}(t)\Vert_{L^2}\\
&\lesssim\varepsilon_1^22^{k^+}\langle t\rangle^{-1+H(q_1,n_1)\delta+36\delta}\langle t\rangle^{H(q_2,n_2)\delta}2^{-N(n_2)k^+},
\end{split}
\end{equation}
using \eqref{sho4} and \eqref{vcx1}. The desired bounds \eqref{saur1}--\eqref{saur2} follow.

{\bf{Step 4.}} Finally, assume that $k_1=\min\{k,k_1,k_2\}$ and $(q_1,n_1)=(q,n), (q_2,n_2)=(0,0)$. The proof is harder in this case mainly because the decomposition \eqref{zaq11} is not effective when $k_1$ is very small, say $2^{k_1}\approx\langle t\rangle^{-1}$. We consider two cases:

{\bf{Case 4.1.}} Let
\begin{equation}\label{saur10}
Y(0,1):=16,\qquad Y(1,1):=26,\qquad Y(q,n):=56\,\,\text{ if }\,\,n\geq 2.
\end{equation}
If $2^{k_1}\in[\langle t\rangle^{-1-4\delta},\langle t\rangle^{-1+Y(q,n)\delta}]$ then we prove the more general bounds 
\begin{equation}\label{saur11}
\begin{split}
2^{k_2-k_1}2^{-k/2}\|P_kI[P_{k_1}&U^{\mathcal{L}_1h_1,\iota_2},P_{k_2}U^{h_2,\iota_2}](t)\|_{L^2}\lesssim\varep_1^2\langle t\rangle^{-1+H(q,n)\delta-\delta/2}2^{-N(n)k^+}2^{-\ga\overline{|k|}/4}
\end{split}
\end{equation}
and
\begin{equation}\label{saur12}
\begin{split}
2^{k_2^+-k_1}\|P_kI[P_{k_1}&U^{\mathcal{L}_1h_1,\iota_1},P_{k_2}U^{\psi,\iota_2}](t)\|_{L^2}\lesssim\varep_1^2\langle t\rangle^{-1+H(q,n)\delta-\delta/2}2^{-N(n)k^+}2^{-\ga\overline{|k|}/4},
\end{split}
\end{equation}
for any $\mathcal{L}_1\in\mathcal{V}^q_n$, $n\geq 1$, $\iota_1,\iota_2\in\{+,-\}$, $h_1,h_2\in\{h_{\al\be}\}$, $k,k_1,k_2\in\mathbb{Z}$, and $t\in[0,T]$.

The bounds \eqref{saur11}--\eqref{saur12} follow easily using just $L^2$ estimates unless $2^{|k_2|}\lesssim \langle t\rangle^{\delta'}$ and $\langle t\rangle\gg 1$. In this case we decompose 
$P_{k_2}U^{h_2,\iota_2}(t)=U^{h_2,\iota_2}_{\leq J_2,k_2}(t)+U^{h_2,\iota_2}_{> J_2,k_2}(t)$ and $P_{k_2}U^{\psi,\iota_2}(t)=U^{\psi,\iota_2}_{\leq J_2,k_2}(t)+U^{\psi,\iota_2}_{> J_2,k_2}(t)$ as in \eqref{on11.3}--\eqref{on11.36}, where $J_2$ is the largest integer satisfying $2^{J_2}\leq\langle t\rangle^{1/10}$. The contributions of the functions $U^{h_2,\iota_2}_{> J_2,k_2}(t)$ and $U^{\psi,\iota_2}_{> J_2,k_2}(t)$ to \eqref{saur11} and \eqref{saur12} respectively can be bounded easily, using again just $L^2$ estimates.

We now consider the main terms. Let $J_1$ denote the largest integer satisfying $2^{J_1}\leq (2^{-k_1}+\langle t\rangle)\langle t\rangle^{\delta/4}$ and decompose
\begin{equation}\label{saur12.2}
U^\ast_{1,\leq J_1}:=P'_{k_1}(\varphi_{\leq J_1}\cdot P_{k_1}U^{\mathcal{L}_1h_1,\iota_1}),\,\,\,\,\,\,U^\ast_{1,> J_1}:=P'_{k_1}(\varphi_{> J_1}\cdot P_{k_1}U^{\mathcal{L}_1h_1,\iota_1}).
\end{equation}
Notice that in this case we decompose in the physical space the normalized solutions  $P_{k_1}U^{\mathcal{L}_1h_1,\iota_1}$, not the profiles  $P_{k_1}V^{\mathcal{L}_1h_1,\iota_1}$. The point is that the functions $U^\ast_{1,\leq J_1}$ satisfy the $L^2$ bounds
\begin{equation}\label{saur14}
\|U^\ast_{1,\leq J_1}(t)\|_{L^2}\lesssim\varep_12^{k_1/2}\langle t\rangle^{H(q,n)\delta-3\delta/4}.
\end{equation}
These bounds are stronger than \eqref{vcx1} (notice the gain of $\langle t\rangle^{-3\delta/4}$) and follow from Lemma \ref{TrivialVF}. Indeed, if $q=0$ then we use \eqref{TrivialVF3} and the assumption $2^{k_1}\in[\langle t\rangle^{-1-4\delta},\langle t\rangle^{-1+Y(q,n)\delta}]$, so
\begin{equation*}
\|U^\ast_{1,\leq J_1}(t)\|_{L^2}\lesssim (2^{k_1}\langle t\rangle)^{-\ga}2^{k_1/2}(1+2^{k_1}\langle t\rangle)\langle t\rangle^{\delta/4}\langle t\rangle^{H(0,n-1)\delta}.
\end{equation*}
The desired bounds \eqref{saur14} follow when $q=0$ since $H(0,n-1)+Y(0,n)\leq H(0,n)-2$, see \eqref{saur10}. The proof is similar in the case $q\geq 1$, using \eqref{TrivialVF2} instead of \eqref{TrivialVF3}.

Using \eqref{bil1}, \eqref{saur14}, and \eqref{vcx1.2}, we find that
\begin{equation}\label{saur16}
\begin{split}
2^{k_2/2-k_1}\|P_kI[U^\ast_{1,\leq J_1},U^{h_2,\iota_2}_{\leq J_2,k_2}](t)\|_{L^2}&\lesssim 2^{-k_1/2}\langle t\rangle^{-1}2^{2k_2}\Vert U^\ast_{1,\leq J_1}(t)\Vert_{L^2}\Vert \widehat{P_{k_2}V^{h_2}}(t)\Vert_{L^\infty}\\
&\lesssim \varepsilon_1^2\langle t\rangle^{-1+H(q,n)\delta-3\delta/4}2^{-N(n)k^+}2^{k^-/2}2^{-3k^+}.
\end{split}
\end{equation}
Similarly, using \eqref{wws13x} and \eqref{vcx1.2} we have
\begin{equation}\label{saur17}
\begin{split}
2^{k_2^+-k_1}\|P_kI[&U^\ast_{1,\leq J_1},U^{\psi,\iota_2}_{\leq J_2,k_2}](t)\|_{L^2}\\
&\lesssim 2^{k_2^+-k_1}\Vert U^\ast_{1,\leq J_1}(t)\Vert_{L^2}\min\{\Vert U^{\psi}_{\leq J_2,k_2}(t)\Vert_{L^\infty},2^{3k_1/2}\Vert U^{\psi}_{\leq J_2,k_2}(t)\Vert_{L^2}\}\\
&\lesssim \varepsilon_1^2\langle t\rangle^{-1+H(q,n)\delta-3\delta/4}2^{-N(n)k^+}2^{\kappa k^-/20}2^{-3k^+}.
\end{split}
\end{equation}

Finally, we claim that the remaining contributions are negligible,
\begin{equation}\label{bnm37}
2^{k_2/2-k_1}\|P_kI[U^\ast_{1,> J_1},U^{h_2,\iota_2}_{\leq J_2,k_2}](t)\|_{L^2}+2^{k_2^+-k_1}\|P_kI[U^\ast_{1,> J_1},U^{\psi,\iota_2}_{\leq J_2,k_2}](t)\|_{L^2}\lesssim \varepsilon_1^2\langle t\rangle^{-2}2^{-N(n)k^+}.
\end{equation}
To see this we use approximate finite speed of propagation arguments. Indeed, we observe that
\begin{equation*}
\begin{split}
P_kI[&U^\ast_{1,> J_1},U^{h_2,\iota_2}_{\leq J_2,k_2}](x,t)\\
&=C\int_{\mathbb{R}^6}\varphi_{>J_1}(z)P_{k_1}U^{\mathcal{L}h_1,\iota_1}(z)\cdot \varphi_{\le J_2}(y)P_{k_2}V^{h_2,\iota_2}(y)K(x-y,x-z)\,dydz,\\
\end{split}
\end{equation*}
where the kernel is given by
\begin{equation*}
\begin{split}
K(y',z')&:=\int_{\mathbb{R}^6}e^{-i\iota_2 t\vert\eta\vert}e^{iy'\cdot \eta}e^{iz'\cdot\theta}m(\theta,\eta)\varphi_k(\theta+\eta)\varphi_{[k_1-2,k_1+2]}(\theta)\varphi_{[k_2-2,k_2+2]}(\eta)\,d\theta d\eta.
\end{split}
\end{equation*}
The point is that $|K(y',z')|$ is small where $|y'-z'|\gtrsim 2^{J_1}$. Indeed, for any $M\geq 1$ we have
\begin{equation*}
\begin{split}
\vert K(y',z')\vert&\lesssim_M 2^{3k_1}2^{3k}\left[1+\frac{|y'|+|z'|}{\langle t\rangle+2^{-k_1}}\right]^{-M},
\end{split}
\end{equation*}
using integration by parts either in $\eta$ or in $\theta$. Since $|y-z|\gtrsim 2^{J_1}\approx \langle t\rangle^{\delta/4}(\langle t\rangle+2^{-k_1})$ in the support of the integral, it follows that the first expression in \eqref{bnm37} is negligible as claimed. The second expression can be bounded in the same way. In view of \eqref{saur16}--\eqref{saur17}, this completes the proof of the bounds \eqref{saur11}--\eqref{saur12}.

{\bf{Case 4.2.}} Finally we prove the bounds \eqref{saur1}--\eqref{saur2} when 
\begin{equation}\label{saur25}
k_1=\min\{k,k_1,k_2\},\qquad (q_1,n_1)=(q,n), (q_2,n_2)=(0,0),\qquad 2^{k_1}\geq\langle t\rangle^{-1+Y(q,n)\delta}.
\end{equation}
We use the identities \eqref{zaq11}. The contributions of the linear commutators $R^a|\nabla|^{-2}\partial_\mu\partial_\nu\LL'_1h_{\al\be}$ are bounded as claimed, due to the estimates \eqref{nake18.01}--\eqref{nake18.03}. It remains to prove that
\begin{equation}\label{saur31}
2^{k_2/2}\|P_kI[P_{k_1}|\nabla|^{-2}\LL_1\mathcal{N},P_{k_2}U^{h,\iota_2}](t)\|_{L^2}\lesssim\varep_1^2\langle t\rangle^{-1+H(q,n)\delta-\delta/2}2^{-N(n)k^+}2^{-\ga\overline{|k|}/4}
\end{equation}
and
\begin{equation}\label{saur32}
2^{k_2^+}\|P_kI[P_{k_1}|\nabla|^{-2}\LL_1\mathcal{N},P_{k_2}U^{\psi,\iota_2}](t)\|_{L^2}\lesssim\varep_1^2\langle t\rangle^{-1+H(q,n)\delta-\delta/2}2^{-N(n)k^+}2^{-\ga\overline{|k|}/4},
\end{equation}
for any $\mathcal{N}\in\mathcal{QU}$. Using \eqref{sho4}, \eqref{wws1}, and \eqref{wws2}, these bounds follow if $2^{k_1}\gtrsim \langle t\rangle^{-1/2}$ or if $2^{k_2}\notin[\langle t\rangle^{-\delta'},\langle t\rangle^{\delta'}]$. 

On the other hand, if $2^{k_1}\leq\langle t\rangle^{-1/2}$ and $2^{k_2}\in[\langle t\rangle^{-\delta'},\langle t\rangle^{\delta'}]$ then we would like to use \eqref{sho4.3} and estimate as in \eqref{saur16}--\eqref{saur17}. For this we replace first $P_{k_2}U^{h,\iota_2}$ with $U^{h,\iota_2}_{\leq J_2,k_2}$ and $P_{k_2}U^{\psi,\iota_2}$ with $U^{\psi,\iota_2}_{\leq J_2,k_2}$ at the expense of acceptable errors, where $J_2$ is the largest integer satisfying $2^{J_2}\leq\langle t\rangle^{1/10}$. Then we estimate, using \eqref{bil1}, \eqref{sho4.3}, and \eqref{vcx1.2},
\begin{equation*}
\begin{split}
2^{k_2/2}\|P_kI[P_{k_1}|\nabla|^{-2}&\LL_1\mathcal{N},U^{h_2,\iota_2}_{\leq J_2,k_2}](t)\|_{L^2}\lesssim 2^{k_1/2}\langle t\rangle^{-1}2^{2k_2}\Vert P_{k_1}|\nabla|^{-2}\LL_1\mathcal{N}(t)\Vert_{L^2}\Vert \widehat{P_{k_2}V^{h_2}}(t)\Vert_{L^\infty}\\
&\lesssim \varepsilon_1^22^{-k_1}\langle t\rangle^{-2+H(q,n)\delta+\ell(q,n)\delta+\delta/2}2^{-N(n)k^+}2^{k^--\kappa k^-}2^{-3k^+}.
\end{split}
\end{equation*}
Similarly, using \eqref{wws13x} and \eqref{vcx1.2} we have
\begin{equation*}
\begin{split}
2^{k_2^+}\|P_kI[&P_{k_1}|\nabla|^{-2}\LL_1\mathcal{N},U^{\psi,\iota_2}_{\leq J_2,k_2}](t)\|_{L^2}\\
&\lesssim 2^{k_2^+-2k_1}\Vert P_{k_1}\LL_1\mathcal{N} (t)\Vert_{L^2}\min\{\Vert U^{\psi}_{\leq J_2,k_2}(t)\Vert_{L^\infty},2^{3k_1/2}\Vert U^{\psi}_{\leq J_2,k_2}(t)\Vert_{L^2}\}\\
&\lesssim \varepsilon_1^22^{-k_1}\langle t\rangle^{-2+H(q,n)\delta+\ell(q,n)\delta+\delta/2}2^{-N(n)k^+}2^{\kappa k^-/20}2^{-3k^+}.
\end{split}
\end{equation*}
Since $2^{-k_1}\langle t\rangle^{-1}\lesssim \langle t\rangle^{-Y(q,n)\delta}$ (see \eqref{saur25}) and $Y(q,n)\geq \ell(q,n)+2$, these bounds suffice to prove \eqref{saur31}--\eqref{saur32}. This completes the proof of the lemma.
\end{proof}

We show now that the nonlinearities $\LL P_{\al\be}^2$ can be written as sums of bilinear expressions involving only the good metric components $\vartheta$, null semilinear forms, and acceptable errors.

\begin{lemma}\label{Pstructure}
If $\LL\in\mathcal{V}^q_n$ then, with $\mathcal{S}^{\LL,2}_{\al\be}$ as in \eqref{BoxWave2.4},
\begin{equation}\label{Pstruc1}
\LL P_{\al\be}^2=-\mathcal{S}^{\LL,2}_{\al\be}+\Pi^{\LL}_{\al\be}+\mathcal{C}^{\LL,\geq 3}_{\al\be}+\mathcal{T}^{\LL}_{\al\be},
\end{equation}
where $\Pi^\LL_{\al\be}$ are semilinear null forms of order $(q,n)$, $\mathcal{C}^{\LL,\geq 3}_{\al\be}$ are semilinear cubic remainders of order $(q,n)$ (of the first and second type described in \eqref{zaq8.6}), and $\mathcal{T}^{\LL}_{\al\be}$ are wave commutator remainders of order $(q,n)$ (see \eqref{nake13.6}).
\end{lemma}

\begin{proof}
Using \eqref{zaq21} we write
\begin{equation*}
\begin{split}
\LL P^2_{\al\be}&=\sum_{\LL_1+\LL_2=\LL}c_{\LL_1,\LL_2}\big\{-(1/2)\partial_\al \LL_1h_{00}\partial_\be \LL_2h_{00}+\partial_\al \LL_1h_{0j}\partial_\be \LL_2h_{0j}-(1/2)\partial_\al \LL_1h_{jk}\partial_\be \LL_2h_{jk}\\
&+(1/4)\partial_\al\LL_1(-h_{00}+\delta_{jk}h_{jk})\partial_\be\LL_2(-h_{00}+\delta_{j'k'}h_{j'k'})\Big\}+\mathcal{T}^{\LL,1}_{\al\be}
\end{split}
\end{equation*}
for some wave commutator remainders $\mathcal{T}^{\LL,1}_{\al\be}$ of order $(q,n)$. Using now \eqref{zaq5l} and assuming $\LL_i\in\mathcal{V}^{q_i}_{n_i}$, $(q_1,n_1)+(q_2,n_2)=(q,n)$, we rewrite the expression between the brackets as
\begin{equation*}
\begin{split}
&-(1/2)\partial_\al (F^{\LL_1}+\uF^{\LL_1})\partial_\be (F^{\LL_2}+\uF^{\LL_2})+\partial_\al (-R_j\rho^{\LL_1}+\in_{jkl}R_k\omega_l^{\LL_1})\partial_\be(-R_j\rho^{\LL_2}+\in_{jk'l'}R_{k'}\omega_{l'}^{\LL_2})\\
&-(1/2)\partial_\al [R_jR_k(F^{\LL_1}-\uF^{\LL_1})-(\in_{klm}R_j+\in_{jlm}R_k)R_l\Omega_m^{\LL_1}+\in_{jpm}\in_{kqn}R_pR_q\vartheta^{\LL_1}_{mn}]\\
&\times\partial_\be [R_jR_k(F^{\LL_2}-\uF^{\LL_2})-(\in_{kl'm'}R_j+\in_{jl'm'}R_k)R_{l'}\Omega_{m'}^{\LL_2}+\in_{jp'm'}\in_{kq'n'}R_{p'}R_{q'}\vartheta_{m'n'}^{\LL_2}]\\
&+\partial_\al(\tau^{\LL_1}-F^{\LL_1})\partial_\be(\tau^{\LL_2}-F^{\LL_2}).
\end{split}
\end{equation*}
We replace also $\rho^{\LL_i}$ with $R_0\uF^{\LL_i}+R_0\tau^{\LL_i}+|\nabla|^{-1}\LL_iE_0^{\geq 2}+|\nabla|^{-1}E^{com}_{\LL_i,0}$ and $\Omega_j^{\LL_i}$ with $R_0\omega_j^{\LL_i}+|\nabla|^{-1}\in_{jlk}R_l{\LL_i}E_k^{\geq 2}+|\nabla|^{-1}\in_{jlk}R_lE^{com}_{\LL_i,k}$, according to the harmonic gauge identities \eqref{zaq11.2}. The contributions of the error terms lead to either semilinear cubic remainders of order $(q,n)$ (of the first type described in \eqref{zaq8.6}) or wave commutator remainders of order $(q,n)$ (of the last type in \eqref{nake13.6}). Then we regroup and expand the remaining terms, in the form
\begin{equation}\label{Pstruc3.5}
A^{FF}_{\al\be}+A^{F\uF}_{\al\be}+A^{F\om}_{\al\be}+A^{F\vartheta}_{\al\be}+A^{\uF\,\uF}_{\al\be}+A^{\uF\om}_{\al\be}+A^{\uF\vartheta}_{\al\be}+A^{\om\om}_{\al\be}+A^{\om\vartheta}_{\al\be}+A^{\vartheta\va}_{\al\be}+\widetilde{A}^\tau_{\al\be},
\end{equation}
where
\begin{equation*}
A^{FF}_{\al\be}=(1/2)(\partial_\al F^{\LL_1}\cdot \partial_\be F^{\LL_2}-R_jR_k\partial_\al F^{\LL_1}\cdot R_jR_k\partial_\be F^{\LL_2}),
\end{equation*}
\begin{equation*}
\begin{split}
A^{F\uF}_{\al\be}&=-(1/2)(\partial_\al F^{\LL_1}\cdot \partial_\be \uF^{\LL_2}-R_jR_k\partial_\al F^{\LL_1}\cdot R_jR_k\partial_\be \uF^{\LL_2})\\
&-(1/2)(\partial_\be F^{\LL_2}\cdot \partial_\al \uF^{\LL_1}-R_jR_k\partial_\be F^{\LL_2}\cdot R_jR_k\partial_\al \uF^{\LL_1}),
\end{split}
\end{equation*}
\begin{equation*}
A^{F\om}_{\al\be}=\in_{klm}R_kR_j\partial_\al F^{\LL_1}\cdot R_l R_jR_0\partial_\be\om_m^{\LL_2}+\in_{klm}R_kR_j\partial_\be F^{\LL_2}\cdot R_l R_jR_0\partial_\al\om_m^{\LL_1},
\end{equation*}
\begin{equation*}
A^{F\va}_{\al\be}=-\frac{1}{2}\in_{jpm}R_jR_k\partial_\al F^{\LL_1}\cdot R_p\in_{kqn} R_q\partial_\be\va^{\LL_2}_{mn}-\frac{1}{2}\in_{jpm}R_jR_k\partial_\be F^{\LL_2}\cdot R_p\in_{kqn} R_q\partial_\al\va^{\LL_1}_{mn},
\end{equation*}
\begin{equation*}
A^{\uF\,\uF}_{\al\be}=-\frac{1}{2}(\partial_\al \uF^{\LL_1}\cdot \partial_\be \uF^{\LL_2}+R_jR_k\partial_\al \uF^{\LL_1}\cdot R_jR_k\partial_\be \uF^{\LL_2}-2R_0R_k\partial_\al \uF^{\LL_1}\cdot R_0R_k\partial_\be \uF^{\LL_2}),
\end{equation*}
\begin{equation*}
\begin{split}
A^{\uF\,\om}_{\al\be}&=-\in_{jkl}R_jR_0\partial_\al\uF^{\LL_1}\cdot  R_k\partial_\be\omega_l^{\LL_2}-\in_{jkl}R_jR_0\partial_\be\uF^{\LL_2}\cdot  R_k\partial_\al\omega_l^{\LL_1}\\
&-\in_{klm}R_kR_j\partial_\al\uF^{\LL_1}\cdot R_lR_jR_0\partial_\be\omega_m^{\LL_2}-\in_{klm}R_kR_j\partial_\be\uF^{\LL_2}\cdot R_lR_jR_0\partial_\al\omega_m^{\LL_1},
\end{split}
\end{equation*}
\begin{equation*}
A^{\uF\va}_{\al\be}=\frac{1}{2}\in_{jpm}R_jR_k\partial_\al\uF^{\LL_1}\cdot R_p\in_{kqn} R_q\partial_\be\va_{mn}^{\LL_2}+\frac{1}{2}\in_{jpm}R_jR_k\partial_\be\uF^{\LL_2}\cdot R_p\in_{kqn} R_q\partial_\al\va_{mn}^{\LL_1},
\end{equation*}
\begin{equation*}
\begin{split}
A^{\om\om}_{\al\be}&=R_l\partial_\al \omega_m^{\LL_1}\cdot  R_l\partial_\be\omega_m^{\LL_2}-R_jR_0 R_l\partial_\al\omega_m^{\LL_1}\cdot R_jR_0 R_l\partial_\be\omega_m^{\LL_2}\\
&-( R_l\partial_\al\omega_m^{\LL_1}\cdot  R_m\partial_\be\omega_l^{\LL_2}-R_jR_0 R_l\partial_\al\omega_m^{\LL_1}\cdot R_jR_0 R_m\partial_\be  \omega_l^{\LL_2})\\
&-(1/2)\in_{klm}\in_{jl'm'}R_l R_jR_0\partial_\al\om_m^{\LL_1}\cdot R_k R_{l'}R_0\partial_\be\om_{m'}^{\LL_2}\\
&-(1/2)\in_{jlm}\in_{kl'm'}R_l R_kR_0\partial_\al\om_m^{\LL_1}\cdot R_j R_{l'}R_0\partial_\be\om_{m'}^{\LL_2},
\end{split}
\end{equation*}
\begin{equation*}
\begin{split}
A^{\om\va}_{\al\be}&=\in_{jkq}R_j R_lR_0\partial_\al\om_m^{\LL_1}\cdot R_kR_l\partial_\be\va_{qm}^{\LL_2}-\in_{jkq}R_j R_lR_0\partial_\al\om_m^{\LL_1}\cdot R_kR_m\partial_\be\va_{ql}^{\LL_2}\\
&+\in_{jkq}R_j R_lR_0\partial_\be\om_m^{\LL_2}\cdot R_kR_l\partial_\al\va_{qm}^{\LL_1}-\in_{jkq}R_j R_lR_0\partial_\be\om_m^{\LL_2}\cdot R_kR_m\partial_\al\va_{ql}^{\LL_1},
\end{split}
\end{equation*}
\begin{equation*}
\begin{split}
A^{\va\va}_{\al\be}=R_pR_q\partial_\al\va_{mn}^{\LL_1}\cdot R_pR_n\partial_\be\va_{mq}^{\LL_2}-\frac{1}{2}R_pR_q\partial_\al\va_{mn}^{\LL_1}\cdot R_pR_q\partial_\be\va_{mn}^{\LL_2}-\frac{1}{2}R_pR_q\partial_\al\va_{mn}^{\LL_1}\cdot R_mR_n\partial_\be\va_{pq}^{\LL_2},
\end{split}
\end{equation*}
\begin{equation*}
\begin{split}
\widetilde{A}^{\tau}_{\al\be}&=\partial_\al\tau^{\LL_1}\cdot \partial_\be\tau^{\LL_2}+R_0R_j\partial_\al\tau^{\LL_1}\cdot R_0R_j\partial_\be\tau^{\LL_2}-\partial_\al\tau^{\LL_1}\cdot \partial_\be F^{\LL_2}-\partial_\al F^{\LL_1}\cdot \partial_\be\tau^{\LL_2}\\
&+R_0R_j\partial_\al\tau^{\LL_1}\cdot R_0R_j\partial_\be\uF^{\LL_2}+R_0R_j\partial_\al\uF^{\LL_1}\cdot R_0R_j\partial_\be\tau^{\LL_2}\\
&-R_0R_j\partial_\al\tau^{\LL_1}\cdot \in_{jpq}R_p\partial_\be\om_q^{\LL_2}-R_0R_j\partial_\be\tau^{\LL_2}\cdot \in_{jpq}R_p\partial_\al\om_q^{\LL_1}.
\end{split}
\end{equation*}

All the terms in $A^{FF}_{\al\be}$, $A^{F\uF}_{\al\be}$, $A^{F\vartheta}_{\al\be}$, $A^{\uF\vartheta}_{\al\be}$ are clearly semilinear null forms (see also \eqref{zaq8.55}). Since $\va$ is divergence free (see \eqref{zaq4}), the terms $R_pR_q\partial_\al\va_{mn}^{\LL_1}\cdot R_pR_n\partial_\be\va_{mq}^{\LL_2}$ and $R_pR_q\partial_\al\va_{mn}^{\LL_1}\cdot R_mR_n\partial_\be\va_{pq}^{\LL_2}$ in $A^{\va\va}_{\al\be}$ are also semilinear null forms. The remaining term in $A^{\va\va}_{\al\be}$ generates the terms $\mathcal{S}_{\al\be}^{\LL,2}$ in \eqref{Pstruc1}. 

The terms that contain $R_0$ require a little more care. We notice first that if $\alpha\neq 0$ then $R_0\partial_\alpha=R_\alpha\partial_0$; moreover, if $\alpha=0$ then we use the identities
\begin{equation}\label{Pstruc4}
R_0\partial_0G^{\LL_i}=|\nabla|^{-1}\{\Delta G^{\LL_i}+G[\LL_i\mathcal{N}^h]\}=R_m\partial_m G^{\LL_i}+|\nabla|^{-1}G[\LL_i\mathcal{N}^h],
\end{equation}
for any $G\in\{F,\uF,\om_m,\vartheta_{mn},\tau\}$. Thus we can express all the terms in $A^{F\om}_{\al\be}$, $A^{\uF\om}_{\al\be}$, and $A^{\om\va}_{\al\be}$ as sums of semilinear null forms and semilinear cubic remainders (of the first type in \eqref{zaq8.6}). 

To deal with the remaining terms, we claim that
\begin{equation}\label{Pstruc5}
R^aR_0\partial_\al G_1^{\LL_1}\cdot R^bR_0\partial_\be G_2^{\LL_2}\approx R^aR_j\partial_\al G_1^{\LL_1}\cdot R^bR_j\partial_\be G_2^{\LL_2}
\end{equation}
for any $G_1,G_2\in\{F,\uF,\om_m,\vartheta_{mn},\tau\}$, $\al,\be\in\{0,1,2,3\}$, and $R^a=R_1^{a_1}R_2^{a_2}R_3^{a_3}$, $R^b=R_1^{b_1}R_2^{b_2}R_3^{b_3}$, where the identity holds up to sums of semilinear null forms and semilinear cubic remainders (of the first two types in \eqref{zaq8.6}). Indeed, if $\al,\be\in\{1,2,3\}$ then this follows from \eqref{nake2.2} and the identities $R_0\partial_\rho=R_\rho\partial_0$, $R_j\partial_\rho=R_\rho\partial_j$, $\rho\in\{\al,\be\}$. If $\alpha=0$ or $\beta=0$ then \eqref{Pstruc5} follows using first \eqref{Pstruc4}, to extract the cubic remainders, and combining with \eqref{nake2}--\eqref{nake2.2}.

We examine now all the terms in $A^{\uF\,\uF}_{\al\be}$ and $A^{\om\om}_{\al\be}$; as a consequence of \eqref{Pstruc5} and \eqref{zaq8.55}, they can all be written as sums of semilinear null forms and semilinear cubic remainders. 

Similarly, all the terms in $\widetilde{A}^{\tau}_{\al\be}$ are sums of semilinear null forms, semilinear cubic remainders, and wave commutator remainders of order $(q,n)$. Indeed, the main point is that $\partial_\mu\tau^{\LL_i}$ can be written as a sum of expressions of the form $R^a|\nabla|^{-1}\LL'\mathcal{N}$ and $R^a\partial_\rho\LL''h$, $\mathcal{N}\in\mathcal{QU}$, $h\in\{h_{\al\be}\}$, $\LL'\in\mathcal{V}_{n_i}^{q_i}$, $\LL''\in\mathcal{V}_{n_i-1}^{q_i}$, $\rho\in\{0,1,2,3\}$ (due to \eqref{zaq11}). We replace also $R_0XR_0Y$ with $R_kXR_kY$, at the expense of acceptable errors, in three of the terms in $\widetilde{A}^{\tau}_{\al\be}$, and use \eqref{Pstruc4} for the terms in the last line of $\widetilde{A}^{\tau}_{\al\be}$ if $\alpha=0$ or $\beta=0$. It follows that all terms in $\widetilde{A}^{\tau}_{\al\be}$ can be written as a sum of acceptable errors, and the lemma follows.
\end{proof}

\chapter{Improved energy estimates}

\section{Setup and preliminary reductions}\label{EnergyEstimates1}

In this chapter we prove the main energy estimates \eqref{bootstrap3.1}. More precisely:

\begin{proposition}\label{EnergEst2}
With the hypothesis in Proposition \ref{bootstrap}, for any $t\in[0,T]$ and $\mathcal{L}\in\mathcal{V}^q_n$, $n\leq 3$, we have
\begin{equation}\label{EnergEst3}
\|(\langle t\rangle|\nabla|_{\leq 1})^{\ga}|\nabla|^{-1/2}U^{\mathcal{L}h_{\al\be}}(t)\|_{H^{N(n)}}+\|U^{\mathcal{L}\psi}(t)\|_{H^{N(n)}}\lesssim\varep_0\langle t\rangle^{H(q,n)\delta}.
\end{equation}
\end{proposition}

\subsection{Energy increments}

To prove Proposition \ref{EnergEst2} we start by defining suitable energy functionals. We consider the modified metric $\widetilde{g}^{\mu\nu}:=m^{\mu\nu}+\widetilde{g}^{\mu\nu}_{\geq 1}$, where $\widetilde{g}^{\mu\nu}_{\geq 1}$ are as in \eqref{BoxWave1}, and the modified wave operator $\widetilde{\square}_{\widetilde{g}}:=\widetilde{g}^{\mu\nu}\partial_\mu\partial_\nu$. Suppose that $\phi$ solves, for $\lambda\in\{0,1\}$
\begin{equation*}
-\widetilde{\Box}_{\widetilde{g}}\phi+\lambda(1-h_{00})\phi=N.
\end{equation*}
We define
\begin{equation}\label{EnergyIncrement1}
\begin{split}
\mathcal{E}_{\lambda}(\phi)&:=\frac{1}{2}\int_{\mathbb{R}^3}\left\{-\widetilde{g}^{00}(\partial_t\phi)^2+\widetilde{g}^{jk}\partial_j\phi\partial_k\phi+\lambda(1-h_{00})(\phi)^2\right\}dx,\\
\mathcal{B}(\phi)&:=\int_{\mathbb{R}^3}\left(\partial_j\widetilde{g}^{0j}(\partial_t\phi)^2-(1/2)\partial_t\widetilde{g}^{jk}\partial_j\phi\partial_k\phi+\partial_j\widetilde{g}^{jk}\partial_t\phi\partial_k\phi+(\lambda/2) \partial_th_{00}(\phi)^2\right)dx.
\end{split}
\end{equation}
Since $\widetilde{g}^{00}=-1$, we have
\begin{equation}\label{EnergyIncrement1.1}
\frac{d}{dt}\mathcal{E}_\lambda(\phi)=-\mathcal{B}(\phi)+\int_{\mathbb{R}^3}\partial_t\phi\cdot  N dx.
\end{equation}
Using \eqref{wws1}, one easily observes that, provided $\varepsilon_1$ is small enough,
\begin{equation}\label{equivEnergy}
\begin{split}
\mathcal{E}_0(\phi)\le \Vert \nabla_{x,t}\phi\Vert_{L^2}^2\le 4\mathcal{E}_0(\phi),\qquad
\mathcal{E}_1(\phi)\le \Vert \nabla_{x,t}\phi\Vert_{L^2}^2+\Vert \phi\Vert_{L^2}^2\le 4 \mathcal{E}_1(\phi).
\end{split}
\end{equation}
This will form the basis of our energy estimates for both $h$ and $\psi$, as we apply these identities to $\LL h_{\al\be}$ (with $\lambda=0$) and $\LL\psi$ (with $\lambda=1$), with suitable multipliers,  and use \eqref{BoxWave}--\eqref{BoxKG}. 

The bulk terms $\mathcal{B}$ have similar regularity as the energy, and have additional "null structure". Indeed, using the formulas \eqref{BoxWave2.9}, \eqref{zaq5l} (with $H=h$), and \eqref{zaq11.2} (notice that the commutator errors are trivial in this case), we have
\begin{equation}\label{Bphi}
\begin{split}
&\partial_j\widetilde{g}^{0j}(\partial_t\phi)^2-(1/2)\partial_t\widetilde{g}^{jk}\partial_j\phi\partial_k\phi+\partial_j\widetilde{g}^{jk}\partial_t\phi\partial_k\phi\\
&=(\partial_t\phi)^2\partial_jh_{0j}+(1/2)\partial_j\phi\partial_k\phi\partial_t(h_{jk}+\delta_{jk}h_{00})-\partial_t\phi\partial_k\phi\partial_j(h_{jk}+\delta_{jk}h_{00})+E^{\al\be}_{\geq 2}\partial_{\al}\phi\partial_\be\phi\\
&=\frac{1}{2}\left(\partial_j\phi\cdot\partial_j\phi\cdot\partial_tF+\partial_j\phi\cdot\partial_k\phi\cdot R_jR_k\partial_tF\right)-\in_{klm}\partial_k\phi\cdot\left[\partial_j\phi\cdot R_jR_l\partial_t\Omega_m+\partial_t\phi\cdot\partial_l\Omega_m\right]\\
&+\frac{1}{2}\left(2(\partial_t\phi)^2\cdot\partial_t\underline{F}+\partial_j\phi\cdot\partial_j\phi\cdot\partial_t\uF-4\partial_t\phi\cdot\partial_k\phi\cdot \partial_k\underline{F}-\partial_j\phi\cdot\partial_k\phi\cdot R_jR_k\partial_t\underline{F}\right)\\
&+\frac{1}{2}\in_{jpm}\in_{kqn}\partial_j\phi\cdot\partial_k\phi\cdot R_pR_q\partial_t\vartheta_{mn}+(\partial_t\phi)^2\cdot\partial_t\tau+E_{\ge 2}^{\alpha\beta}\cdot\partial_\alpha\phi\cdot\partial_\beta\phi.
\end{split}
\end{equation}
The coefficients of the last two terms in the last line above are of the form $R^a|\nabla|^{-1}\mathcal{N}$, for some $\mathcal{N}\in\mathcal{QU}$ (see \eqref{zaq11} and Definition \ref{nake4}). 

We claim now that the main terms in \eqref{Bphi} can be expressed in terms of semilinear null forms. Indeed, if $\phi$ is a wave unknown then we let $U^{\phi,\iota}:=\left(\partial_t-\iota i\Lambda_{wa}\right)\phi$ and rewrite
\begin{equation}\label{Bphi2}
\frac{1}{2}\int_{\mathbb{R}^3}\left(\partial_j\phi\cdot\partial_j\phi\cdot\partial_tF+\partial_j\phi\cdot\partial_k\phi\cdot R_jR_k\partial_tF\right)\,dx=\sum_{\iota_1,\iota_2,\iota\in\{\pm\}} \mathcal{I}_{p^{F,wa}_{\iota_1,\iota_2,\iota}}[U^{F,\iota_1},U^{\phi,\iota_2},U^{\phi,\iota}]
\end{equation}
where
\begin{equation}\label{NonlinearMultW1}
\mathcal{I}_{p}[F,G,H]:=\frac{1}{(8\pi^3)^2}\int_{\mathbb{R}^3\times\mathbb{R}^3}p(\xi,\eta)\widehat{F}(\xi-\eta)\widehat{G}(\eta) \overline{\widehat{H}(\xi)}\,d\xi d\eta,
\end{equation}
and
\begin{equation}\label{symp1}
p^{F,wa}_{\iota_1,\iota_2,\iota}(\xi,\eta):=\frac{1}{16}\Big[\iota\iota_2\frac{\eta}{|\eta|}\cdot\frac{\xi}{|\xi|}-\iota\iota_2\Big(\frac{\xi-\eta}{|\xi-\eta|}\cdot\frac{\eta}{|\eta|}\Big)\Big(\frac{\xi-\eta}{|\xi-\eta|}\cdot\frac{\xi}{|\xi|}\Big)\Big].
\end{equation}

The other terms in the right-hand side of \eqref{Bphi} can be written in a similar way, as sums of integrals of the form $\mathcal{I}_{p^{G,wa}_{\iota_1,\iota_2,\iota}}[U^{G,\iota_1},U^{\phi,\iota_2},U^{\phi,\iota}]$, $G\in\{\uF,\Omega_m,\va_{mn}\}$, with symbols
\begin{equation}\label{symp2}
\begin{split}
p^{\uF,wa}_{\iota_1,\iota_2,\iota}(\xi,\eta):=\frac{1}{16}\Big[2+\iota\iota_2\frac{\eta}{|\eta|}\cdot\frac{\xi}{|\xi|}-2\iota_1\iota_2\frac{\xi-\eta}{|\xi-\eta|}\cdot\frac{\eta}{|\eta|}-2\iota_1\iota\frac{\xi-\eta}{|\xi-\eta|}\cdot\frac{\xi}{|\xi|}\\
+\iota\iota_2\Big(\frac{\xi-\eta}{|\xi-\eta|}\cdot\frac{\eta}{|\eta|}\Big)\Big(\frac{\xi-\eta}{|\xi-\eta|}\cdot\frac{\xi}{|\xi|}\Big)\Big],
\end{split}
\end{equation}
\begin{equation}\label{symp3}
p^{\Omega_{m},wa}_{\iota_1,\iota_2,\iota}(\xi,\eta):=\frac{1}{8}\iota_2\in_{klm}\frac{(\xi-\eta)_l}{|\xi-\eta|}\frac{\eta_k}{|\eta|}\Big(\iota\frac{\xi-\eta}{|\xi-\eta|}\cdot\frac{\xi}{|\xi|}-\iota_1\Big),
\end{equation}
\begin{equation}\label{symp4}
p^{\va_{mn},wa}_{\iota_1,\iota_2,\iota}(\xi,\eta):=-\frac{1}{16}\iota\iota_2\Big(\in_{jpm}\frac{(\xi-\eta)_p}{|\xi-\eta|}\frac{\eta_j}{|\eta|}\Big)\Big(\in_{kqn}\frac{(\xi-\eta)_q}{|\xi-\eta|}\frac{\xi_k}{|\xi|}\Big).
\end{equation}

The calculation is similar if $\phi$ is a Klein-Gordon variable. In this case we define $U^{\phi,\iota}:=\left(\partial_t-\iota i\Lambda_{kg}\right)\phi$. We add up the term $\partial_th_{00}\phi^2/2$ from \eqref{EnergyIncrement1}, and rewrite the main terms in the right-hand side of \eqref{Bphi} as sums of integrals of the form $\mathcal{I}_{p^{G,kg}_{\iota_1,\iota_2,\iota}}[U^{G,\iota_1},U^{\phi,\iota_2},U^{\phi,\iota}]$, $G\in\{F,\uF,\Omega_m,\va_{mn}\}$, $\iota_1,\iota_2,\iota\in\{+,-\}$. The resulting symbols are 
\begin{equation}\label{symp5}
p^{F,kg}_{\iota_1,\iota_2,\iota}(\xi,\eta):=\frac{1}{16}\Big[\iota\iota_2\frac{\eta\cdot\xi-1}{\langle\eta\rangle\langle\xi\rangle}-\iota\iota_2\Big(\frac{\xi-\eta}{|\xi-\eta|}\cdot\frac{\eta}{\langle\eta\rangle}\Big)\Big(\frac{\xi-\eta}{|\xi-\eta|}\cdot\frac{\xi}{\langle\xi\rangle}\Big)\Big],
\end{equation}
\begin{equation}\label{symp6}
\begin{split}
p^{\uF,kg}_{\iota_1,\iota_2,\iota}(\xi,\eta):=\frac{1}{16}\Big[2+\iota\iota_2\frac{\eta\cdot\xi-1}{\langle\eta\rangle\langle\xi\rangle}-2\iota_1\iota_2\frac{\xi-\eta}{|\xi-\eta|}\cdot\frac{\eta}{\langle\eta\rangle}-2\iota_1\iota\frac{\xi-\eta}{|\xi-\eta|}\cdot\frac{\xi}{\langle\xi\rangle}\\
+\iota\iota_2\Big(\frac{\xi-\eta}{|\xi-\eta|}\cdot\frac{\eta}{\langle\eta\rangle}\Big)\Big(\frac{\xi-\eta}{|\xi-\eta|}\cdot\frac{\xi}{\langle\xi\rangle}\Big)\Big],
\end{split}
\end{equation}
\begin{equation}\label{symp7}
p^{\Omega_{m},kg}_{\iota_1,\iota_2,\iota}(\xi,\eta):=\frac{1}{8}\iota_2\in_{klm}\frac{(\xi-\eta)_l}{|\xi-\eta|}\frac{\eta_k}{\langle\eta\rangle}\Big(\iota\frac{\xi-\eta}{|\xi-\eta|}\cdot\frac{\xi}{\langle\xi\rangle}-\iota_1\Big),
\end{equation}
\begin{equation}\label{symp8}
p^{\va_{mn},kg}_{\iota_1,\iota_2,\iota}(\xi,\eta):=-\frac{1}{16}\iota\iota_2\Big(\in_{jpm}\frac{(\xi-\eta)_p}{|\xi-\eta|}\frac{\eta_j}{\langle\eta\rangle}\Big)\Big(\in_{kqn}\frac{(\xi-\eta)_q}{|\xi-\eta|}\frac{\xi_k}{\langle\xi\rangle}\Big).
\end{equation}

\subsection{The main space-time bounds} We would like now to use the calculations in the previous subsection to start the proof of Proposition \ref{EnergEst2}. Assume $\LL\in \mathcal{V}^q_n$ and define
\begin{equation}\label{ban7}
P^n_{wa}:=\langle \nabla\rangle^{N(n)}\vert\nabla\vert^{-1/2}\vert\nabla\vert_{\le 1}^\ga,\qquad P^n_{kg}:=\langle \nabla\rangle^{N(n)}.
\end{equation}
We define also the associated multipliers $P^n_{wa}(\xi):=\langle \xi\rangle^{N(n)}\vert\xi\vert^{-1/2}\vert\xi\vert_{\le 1}^\ga$ and $P^n_{kg}(\xi):=\langle \xi\rangle^{N(n)}$, $\xi\in\mathbb{R}^3$. We first let $\phi:=P^n_{wa}(\LL h)$, $h\in\{h_{\al\be}\}$, and write
\begin{equation*}
\begin{split}
\left\vert \mathcal{E}_0(P_{wa}^n(\mathcal{L}h))(t)\right\vert
&\lesssim \left\vert \mathcal{E}_0(P_{wa}^n(\mathcal{L}h))(0)\right\vert+\left\vert\int_{0}^t \mathcal{B}(P_{wa}^n(\mathcal{L}h))(s)ds\right\vert\\
&+\left\vert \int_{0}^t\int_{\mathbb{R}^3}\partial_s(P_{wa}^n(\mathcal{L}h))(s)\cdot\widetilde{\Box}_{\widetilde{g}}(P_{wa}^n(\mathcal{L}h))(s)ds \right\vert,
\end{split}
\end{equation*}
for any $t\in[0,T]$, as a consequence of \eqref{EnergyIncrement1.1}. Since $\left\vert \mathcal{E}_0(P_{wa}^n(\mathcal{L}h))(s)\right\vert\approx \|P_{wa}^n(U^{\LL h}(s))\|^2_{L^2}$ for any $s\in [0,t]$ (due to \eqref{equivEnergy}), to prove the first inequality in \eqref{EnergEst3} it suffices to show that
\begin{equation}\label{ban9}
\begin{split}
\left\vert\int_{0}^t \mathcal{B}(P_{wa}^n(\mathcal{L}h))(s)\,ds\right\vert+\left\vert \int_{0}^t\int_{\mathbb{R}^3}P_{wa}^n(\partial_s(\mathcal{L}h))(s)\cdot\widetilde{\Box}_{\widetilde{g}}(P_{wa}^n(\mathcal{L}h))(s)ds \right\vert\lesssim \varep_1^3\langle t\rangle^{2H(q,n)\delta-2\ga},
\end{split}
\end{equation}
for any $t\in[0,T]$. Similarly, to prove the second inequality in \eqref{EnergEst3} it suffices to show that
\begin{equation}\label{ban9.5}
\begin{split}
\Big\vert\int_{0}^t &\mathcal{B}(P_{kg}^n(\mathcal{L}\psi))(s)\,ds\Big\vert+\Big\vert \int_{0}^t\int_{\mathbb{R}^3}P_{kg}^n(\partial_s(\mathcal{L}\psi))(s)\\
&\times\big\{-\widetilde{\Box}_{\widetilde{g}}(P_{kg}^n(\mathcal{L}\psi))(s)+(1-h_{00})P_{kg}^n(\mathcal{L}\psi)(s)\big\}ds \Big\vert\lesssim \varep_1^3\langle t\rangle^{2H(q,n)\delta}.
\end{split}
\end{equation}

We decompose the time integrals into dyadic pieces. More precisely, given $t\in[0,T]$, we fix a suitable decomposition of the function $\mathbf{1}_{[0,t]}$, i.e. we fix functions $q_0,\ldots,q_{L+1}:\mathbb{R}\to[0,1]$, $|L-\log_2(2+t)|\leq 2$, with the properties
\begin{equation}\label{nh2}
\begin{split}
&\mathrm{supp}\,q_0\subseteq [0,2], \quad \mathrm{supp}\,q_{L+1}\subseteq [t-2,t],\quad\mathrm{supp}\,q_m\subseteq [2^{m-1},2^{m+1}]\text{ for } m\in\{1,\ldots,L\},\\
&\sum_{m=0}^{L+1}q_m(s)=\mathbf{1}_{[0,t]}(s),\qquad q_m\in C^1(\mathbb{R})\text{ and }\int_0^t|q'_m(s)|\,ds\lesssim 1\text{ for }m\in \{1,\ldots,L\}.
\end{split}
\end{equation}

The following proposition provides our main estimates on space-time contributions.

\begin{proposition}\label{totbounds}
Assume that $t\in[0,T]$, $\iota,\iota_1,\iota_2\in\{+,-\}$, $h,h_1,h_2\in\{h_{\al\be}\}$, $\LL\in\mathcal{V}^q_n$, $\LL_i\in\mathcal{V}^{q_i}_{n_i}$, $i\in\{1,2\}$, $\LL_1+\LL_2\leq\LL$.  For any bounded symbol $q$ let $I_q$ denote the associated operator as in \eqref{abc36.1}. Let $J_m$ denote the supports of the functions $q_m$ defined above. Then:

(1) if $\LL_2\neq \LL$ and $\mathfrak{q}^{\ast,wa}_{\iota_1\iota_2}$ is a symbol of the form $\mathfrak{q}^{G,wa}_{\iota_1\iota_2}(\theta,\eta)m_0(\theta)m_1(\eta)m_2(\theta+\eta)$, where $G\in \{F,\uF,\omega_n,\va_{mn}\}$, $\mathfrak{q}^{G,wa}_{\iota_1\iota_2}$ are as in \eqref{BoxWave7}, and $m_0,m_1,m_2\in\mathcal{M}_0$, then
\begin{equation}\label{EEEst3}
\begin{split}
\Big|\int_{J_m}\int_{\mathbb{R}^3}q_m(s)[P_{wa}^n(\xi)]^2\mathcal{F}\{I_{\mathfrak{q}^{\ast,wa}_{\iota_1\iota_2}}[\vert\nabla\vert^{-1}U^{\mathcal{L}_1h_1,\iota_1},\vert \nabla\vert U^{\mathcal{L}_2h_2,\iota_2}]\}(\xi,s)\cdot\overline{\widehat{U^{\mathcal{L}h,\iota}}(\xi,s)}\,d\xi ds\Big|\\
\lesssim \varep_1^32^{2H(q,n)\delta m-2\ga m};
\end{split}
\end{equation}

(2) if $\mathfrak{n}_{\iota_1\iota_2}\in\mathcal{M}^{null}_{\iota_1\iota_2}$ then
\begin{equation}\label{EEEst4}
\begin{split}
\Big|\int_{J_m}\int_{\mathbb{R}^3}q_m(s)\mathcal{F}\{P_{wa}^nI_{\mathfrak{n}_{\iota_1\iota_2}}[U^{\mathcal{L}_1h_1,\iota_1},U^{\mathcal{L}_2h_2,\iota_2}]\}(\xi,s)\cdot\overline{\widehat{P_{wa}^nU^{\mathcal{L}h,\iota}}(\xi,s)}\,d\xi ds\Big|\lesssim \varep_1^32^{-\delta m}
\end{split}
\end{equation}
and
\begin{equation}\label{EEEst4.4}
\begin{split}
\Big|\int_{J_m}\int_{\mathbb{R}^3}q_m(s)\mathcal{F}\{I_{\mathfrak{n}_{\iota_1\iota_2}}[U^{h_1,\iota_1},P_{wa}^nU^{\mathcal{L}h_2,\iota_2}]\}(\xi,s)\cdot\overline{\widehat{P_{wa}^nU^{\mathcal{L}h,\iota}}(\xi,s)}\,d\xi ds\Big|\lesssim \varep_1^32^{-\delta m};
\end{split}
\end{equation}

(3) if $\mathfrak{m}\in\mathcal{M}$ and $\vartheta\in\{\va_{mn}\}$, then
\begin{equation}\label{EEEst2}
\begin{split}
\Big|\int_{J_m}\int_{\mathbb{R}^3}q_m(s)[P_{wa}^n(\xi)]^2\mathcal{F}\{I_{\mathfrak{m}}[U^{\va^{\mathcal{L}_1},\iota_1},U^{\va^{\mathcal{L}_2},\iota_2}]\}(\xi,s)\cdot\overline{\widehat{U^{\mathcal{L}h,\iota}}(\xi,s)}\,d\xi ds\Big|\\
\lesssim \varep_1^32^{2H(q,n)\delta m-2\ga m};
\end{split}
\end{equation}

(4) if $\mu,\nu\in\{0,1,2,3\}$ and $\widetilde{g}^{\mu\nu}_{\geq 1}$ are as in \eqref{BoxWave1}, then
\begin{equation}\label{EEcom1}
\begin{split}
\Big|\int_{J_m}\int_{\mathbb{R}^3}q_m(s)\mathcal{F}\big\{P_{wa}^n[\widetilde{g}^{\mu\nu}_{\geq 1}\partial_\mu\partial_\nu(\LL h_2)]-\widetilde{g}^{\mu\nu}_{\geq 1}\partial_\mu\partial_\nu(P_{wa}^n\LL h_2)\big\}(\xi,s)\cdot\overline{\widehat{P_{wa}^nU^{\mathcal{L}h,\iota}}(\xi,s)}\,d\xi ds\Big|\\
\lesssim \varep_1^32^{2H(q,n)\delta m-2\ga m};
\end{split}
\end{equation}

(5) if $\mathfrak{m}\in\mathcal{M}$ then
\begin{equation}\label{EEEst1}
\begin{split}
\Big|\int_{J_m}\int_{\mathbb{R}^3}q_m(s)[P_{wa}^n(\xi)]^2\mathcal{F}\{I_{\mathfrak{m}}[U^{\mathcal{L}_1\psi,\iota_1},U^{\mathcal{L}_2\psi,\iota_2}]\}(\xi,s)\cdot\overline{\widehat{U^{\mathcal{L}h,\iota}}(\xi,s)}\,d\xi ds\Big|\\
\lesssim \varep_1^32^{2H(q,n)\delta m-2\ga m};
\end{split}
\end{equation}

(6) if $\mathfrak{m}\in\mathcal{M}$ and $n_2<n$ then
\begin{equation}\label{EkgEst3}
\begin{split}
\Big|\int_{J_m}\int_{\mathbb{R}^3}q_m(s)[P_{kg}^n(\xi)]^2\mathcal{F}\{I_{\mathfrak{m}}[\vert\nabla\vert^{-1}U^{\mathcal{L}_1h_1,\iota_1},\langle\nabla\rangle U^{\mathcal{L}_2\psi,\iota_2}]\}(\xi,s)\cdot\overline{\widehat{U^{\mathcal{L}\psi,\iota}}(\xi,s)}\,d\xi ds\Big|\\
\lesssim \varep_1^32^{2H(q,n)\delta m};
\end{split}
\end{equation}

(7) if $\mu,\nu\in\{0,1,2,3\}$ and $\widetilde{g}^{\mu\nu}_{\geq 1}$ are as in \eqref{BoxWave1}, then
\begin{equation}\label{EEcom2}
\begin{split}
\Big|\int_{J_m}\int_{\mathbb{R}^3}q_m(s)\mathcal{F}\big\{P_{kg}^n[\widetilde{g}^{\mu\nu}_{\geq 1}\partial_\mu\partial_\nu(\LL\psi)+h_{00}\LL\psi]-[\widetilde{g}^{\mu\nu}_{\geq 1}\partial_\mu\partial_\nu(P_{kg}^n\LL\psi)+h_{00}P_{kg}^n\LL\psi]\big\}(\xi,s)\\
\times\overline{\widehat{P_{kg}^nU^{\mathcal{L}\psi,\iota}}(\xi,s)}\,d\xi ds\Big|\lesssim \varep_1^32^{2H(q,n)\delta m};
\end{split}
\end{equation}

(8) if $\mathfrak{p}^{\ast,kg}_{\iota_1,\iota_2,\iota}$ is a symbol of the form $p^{G,kg}_{\iota_1,\iota_2,\iota}(\xi,\eta)m_0(\xi-\eta)m_1(\eta)m_2(\xi)$, where $m_0,m_1,m_2\in\mathcal{M}_0$, $G\in \{F,\uF,\Omega_n,\va_{mn}\}$,  and $p^{G,kg}_{\iota_1,\iota_2,\iota}$ are as in \eqref{symp5}--\eqref{symp8}, then
\begin{equation}\label{BulkKG}
\begin{split}
\Big|\int_{J_m}\int_{\mathbb{R}^3\times\mathbb{R}^3}q_m(s){\mathfrak{p}}^{\ast,kg}_{\iota_1,\iota_2,\iota}(\xi,\eta)\widehat{U^{h_1,\iota_1}}(\xi-\eta,s)\widehat{P^n_{kg}U^{\LL \psi,\iota_2}}(\eta,s) \overline{\widehat{P^n_{kg}U^{\LL\psi,\iota}}(\xi,s)}\,d\xi d\eta ds\Big|\\
\lesssim \varep_1^32^{2H(q,n)\delta m}.
\end{split}
\end{equation}
\end{proposition}

\begin{proof}[Proof of Proposition \ref{EnergEst2}] It is easy to see that the bounds \eqref{ban9}--\eqref{ban9.5} would follow from Proposition \ref{totbounds}. Indeed, we start from the identity $-\square(\LL h_{\al\be})=\LL\mathcal{N}^h_{\al\be}$, and use \eqref{BoxWave}. Thus
\begin{equation}\label{ban11}
\begin{split}
-\square&(P_{wa}^n\LL h_{\al\be})-\sum_{\mu,\nu\in\{0,1,2,3\}}\widetilde{g}^{\mu\nu}_{\geq 1}\partial_\mu\partial_\nu(P_{wa}^n\LL h_{\al\be})\\
&=\sum_{\mu,\nu\in\{0,1,2,3\}}\big\{P^n_{wa}[\widetilde{g}^{\mu\nu}_{\geq 1}\partial_\mu\partial_\nu (\LL h_{\al\be})]-\widetilde{g}^{\mu\nu}_{\geq 1}\partial_\mu\partial_\nu(P_{wa}^n\LL h_{\al\be})\big\}+P_{wa}^n\mathcal{Q}_{wa}^{\LL}(h_{\al\be})\\
&+P_{wa}^n\mathcal{S}_{\al\be}^{\LL,1}+P_{wa}^n\mathcal{S}_{\al\be}^{\LL,2}+P_{wa}^n\mathcal{KG}^{\LL}_{\al\be}+P_{wa}^n\mathcal{R}^{\LL h}_{\al\be}.
\end{split}
\end{equation}
Therefore the contribution of the second term in the left-hand side of \eqref{ban9} can be bounded as claimed, as a consequence of the estimates \eqref{EEEst3}--\eqref{EEEst1} and the definitions. 

To estimate the contribution of the first term in the left-hand side of \eqref{ban9} we examine the formulas \eqref{Bphi}--\eqref{symp4}. The contribution of the last two terms in \eqref{Bphi} is bounded as claimed, due to the $L^\infty$ bounds in \eqref{sho4.2}. Moreover, we claim that all the other terms are null forms that can be estimated using \eqref{EEEst4.4}. Indeed, the symbols $p_{\iota_1,\iota_2,\iota}^{\Omega_m,wa}$ and $p_{\iota_1,\iota_2,\iota}^{\va_{mn},wa}$ are clearly null in the variables $\xi-\eta$ and $\eta$, due to \eqref{nake2}. Moreover, by examining \eqref{symp1}--\eqref{symp2} we can write
\begin{equation*}
16p_{\iota_1,\iota_2,\iota}^{F,wa}(\xi,\eta)=\Big(\iota_1\frac{\xi-\eta}{|\xi-\eta|}-\iota_2\frac{\eta}{|\eta|}\Big)\cdot\Big(\iota_1\frac{\xi-\eta}{|\xi-\eta|}-\iota\frac{\xi}{|\xi|}\Big)-\frac{1}{4}\Big|\iota_1\frac{\xi-\eta}{|\xi-\eta|}-\iota_2\frac{\eta}{|\eta|}\Big|^2\Big|\iota_1\frac{\xi-\eta}{|\xi-\eta|}-\iota\frac{\xi}{|\xi|}\Big|^2,
\end{equation*}
and
\begin{equation*}
16p_{\iota_1,\iota_2,\iota}^{\uF,wa}(\xi,\eta)=\Big(\iota_1\frac{\xi-\eta}{|\xi-\eta|}-\iota_2\frac{\eta}{|\eta|}\Big)\cdot\Big(\iota_1\frac{\xi-\eta}{|\xi-\eta|}-\iota\frac{\xi}{|\xi|}\Big)+\frac{1}{4}\Big|\iota_1\frac{\xi-\eta}{|\xi-\eta|}-\iota_2\frac{\eta}{|\eta|}\Big|^2\Big|\iota_1\frac{\xi-\eta}{|\xi-\eta|}-\iota\frac{\xi}{|\xi|}\Big|^2,
\end{equation*}
thus $p_{\iota_1,\iota_2,\iota}^{F,wa}$ and $p_{\iota_1,\iota_2,\iota}^{\uF,wa}$ are sums of acceptable null symbols in the variables $\xi-\eta$ and $\eta$ multiplied by symbols of $\xi$, as desired.

The estimates \eqref{ban9.5} follow in a similar way, using \eqref{EkgEst3}-\eqref{EEcom2} to bound the contributions of the nonlinearities, and using \eqref{BulkKG} for the space-time integral of $\mathcal{B}(P_{kg}^n(\mathcal{L}\psi))(s)$.
\end{proof}

\subsection{Poincar\'e normal forms}

We examine now the bounds in Proposition \ref{totbounds}, and notice that all the space-time integrals in the left-hand sides do not have derivative loss. In most cases, however, the time decay we have is not enough to allow direct estimates. In such a situation we would like to integrate by parts in time (the method of normal forms) to prove the desired space-time estimates.

More precisely, assume that we are given a trilinear form
\begin{equation}\label{DefI}
\mathcal{G}_{\mathfrak{m}}[f,g,h]=\int_{\mathbb{R}^3\times\mathbb{R}^3}\mathfrak{m}(\xi-\eta,\eta)\widehat{f}(\xi-\eta)\widehat{g}(\eta)\overline{\widehat{h}(\xi)}\,d\xi d\eta,
\end{equation}
for a suitable multiplier $\mathfrak{m}$, and consider its associated quadratic phase
\begin{equation*}
\Phi_{\sigma\mu\nu}(\xi,\eta)=\Lambda_\sigma(\xi)-\Lambda_\mu(\xi-\eta)-\Lambda_\nu(\eta),
\end{equation*}
where $\sigma,\mu,\nu\in\{(wa,+),(wa,-),(kg,+),(kg,-)\}$ as in \eqref{on9.1}. Define
\begin{equation}\label{DefJ}
\begin{split}
\mathcal{H}_{\mathfrak{m}}[f,g,h]&:=\int_{\mathbb{R}^3\times\mathbb{R}^3}\frac{\mathfrak{m}(\xi-\eta,\eta)}{\Phi_{\sigma\mu\nu}(\xi,\eta)}\widehat{f}(\xi-\eta)\widehat{g}(\eta)\overline{\widehat{h}(\xi)}\,d\xi d\eta.
\end{split}
\end{equation}
Using integration by parts in time, for $m\in\{1,\ldots,L\}$ we have
\begin{equation}\label{NF2}
\begin{split}
-i\int_{\mathbb{R}}q_m(s)\mathcal{G}_{\mathfrak{m}}[f(s),g(s),h(s)]ds=\mathcal{H}^1+\mathcal{H}^2+\mathcal{H}^3+\mathcal{H}^4
\end{split}
\end{equation}
where the functions $q_m$ are defined as in \eqref{nh2}, and
\begin{equation}\label{NF3}
\begin{split}
\mathcal{H}^1&=\int_{\mathbb{R}}q^\prime_m(s)\mathcal{H}_{\mathfrak{m}}[f(s),g(s),h(s)]ds,\\
\mathcal{H}^2&=\int_{\mathbb{R}}q_m(s)\mathcal{H}_{\mathfrak{m}}[\left(\partial_s+i\Lambda_{\mu}\right)f(s),g(s),h(s)]ds,\\
\mathcal{H}^3&=\int_{\mathbb{R}}q_m(s)\mathcal{H}_{\mathfrak{m}}[f(s),\left(\partial_s+i\Lambda_\nu\right)g(s),h(s)]ds,\\
\mathcal{H}^4&=\int_{\mathbb{R}}q_m(s)\mathcal{H}_{\mathfrak{m}}[f(s),g(s),\left(\partial_s+i\Lambda_\sigma\right)h(s)]ds.\\
\end{split}
\end{equation}

In other words, we can estimate integrals like those in the left-hand side of \eqref{NF2} in terms of integrals such as those in \eqref{NF3}. The point is to gain time integrability. Indeed, the identities \eqref{zaq11.1} and \eqref{zaq11.3} show that 
\begin{equation}\label{roco1}
(\partial_t^2+\Lambda_{wa}^2)(\mathcal{L} h_{\al\be})=\mathcal{L} \mathcal{N}^{h}_{\al\be},\qquad (\partial_t^2+\Lambda_{kg}^2)(\mathcal{L} \psi)=\mathcal{L} \mathcal{N}^{\psi},
\end{equation}
for any $\al,\be\in\{0,1,2,3\}$ and $\mathcal{L}\in\mathcal{V}_n^q$, $n\in[0,3]$. Using the definitions \eqref{variables4L}, it follows that
\begin{equation}\label{roco2}
(\partial_t+i\Lambda_{wa})U^{\mathcal{L}h_{\al\be}}=\mathcal{L}\mathcal{N}^{h}_{\al\be},\qquad (\partial_t+i\Lambda_{kg})U^{\mathcal{L}\psi}=\LL\mathcal{N}^{\psi}.
\end{equation}

The main difficulty when applying \eqref{NF2} is the presence of time-resonances, which are frequencies $(\xi,\eta)$ for which $\Phi_{\sigma\mu\nu}(\xi,\eta)=0$, and which lead to problems in estimating the terms $\mathcal{H}_{\mathfrak{m}}[f,g,h]$ in \eqref{DefJ}. These resonances  are stronger in the case of trilinear wave interactions, where parallel frequencies $(\xi,\eta)$ lead to resonances; in the case of mixed interactions involving two Klein-Gordon fields and one wave component, resonances only occur when the frequency of the wave component vanishes (see Lemmas \ref{PhaWave} and \ref{pha2}).

In some cases we can use simple estimates to control trilinear expressions, such as 
\begin{equation}\label{EstimJ}
\begin{split}
\left\vert \mathcal{G}_{\mathfrak{m}}[P_{k_1}F,P_{k_2}G,P_kH]\right\vert&\lesssim \Vert \mathcal{F}^{-1}\{\mathfrak{m}\cdot \varphi_{kk_1k_2}\}\Vert_{L^1}\Vert P_{k_1}F\Vert_{L^{p_1}}\Vert P_{k_2}G\Vert_{L^{p_2}}\Vert P_kH\Vert_{L^{p}},\\
\end{split}
\end{equation}
provided that $k,k_1,k_2\in\mathbb{Z}$, $p,p_1,p_2\in[1,\infty]$, $1/p_1+1/p_2+1/p=1$, where
\begin{equation}\label{phikkk}
\varphi_{kk_1k_2}(\theta,\eta):=\varphi_{[k-2,k+2]}(\eta+\theta)\varphi_{[k_1-2,k_1+2]}(\theta)\varphi_{[k_2-2,k_2+2]}(\eta).
\end{equation}
These bounds follow from Lemma \ref{L1easy} (i). We also have pure $L^2$ bounds
\begin{equation}\label{EstimJ2}
\begin{split}
\left\vert \mathcal{G}_{\mathfrak{m}}[P_{k_1}F,P_{k_2}G,P_kH]\right\vert\lesssim 2^{3\min\{k,k_1,k_2\}/2}\Vert \mathfrak{m}\cdot \varphi_{kk_1k_2}\Vert_{L^\infty}\Vert P_{k_1}F\Vert_{L^{2}}\Vert P_{k_2}G\Vert_{L^{2}}\Vert P_kH\Vert_{L^{2}}.
\end{split}
\end{equation}
Trilinear expressions like $\mathcal{H}_{\mathfrak{m}}$ are often estimated using Lemmas \ref{PhaWave} and \ref{pha2}.

\subsection{Elements of paradifferential calculus}\label{ParaCalc}

There is one other issue when applying normal forms, namely the possible loss of derivatives coming from the quasilinear nonlinearities. As in some of our earlier work \cite{DeIoPa} and \cite{DeIoPaPu}, we deal with this using paradifferential calculus. 

In this subsection, we summarize some tools of basic paradifferential calculus and prove a proposition about paralinearization of our wave and Klein-Gordon operators. We recall first the definition of paradifferential operators (Weyl quantization): given a symbol $a=a(x,\zeta):\mathbb{R}^3\times\mathbb{R}^3\to\mathbb{C}$, we define the operator $T_a$ by
\begin{equation}\label{Tsigmaf2}
\begin{split}
\mathcal{F}\left\{T_{a}f\right\}(\xi)=\frac{1}{8\pi^3}\int_{\mathbb{R}^3}\chi_0\Big(\frac{\vert\xi-\eta\vert}{\vert\xi+\eta\vert}\Big)\widetilde{a}(\xi-\eta,(\xi+\eta)/2)\widehat{f}(\eta)d\eta,
\end{split}
\end{equation}
where $\widetilde{a}$ denotes the partial Fourier transform of $a$ in the first coordinate and $\chi_0=\varphi_{\le-100}$.

We will use a simple norm to estimate symbols: for $q\in [1,\infty]$ and $l\in \mathbb{R}$ we define
\begin{equation}\label{nor1}
\begin{split}
&\|a\|_{\mathcal{L}^q_l}:=\sup_{\zeta\in\mathbb{R}^3}(1+|\zeta|^2)^{-l/2}\|\,|a|(.,\zeta)\|_{L^q_x},\\
&|a|(x,\zeta):=\sum_{|\beta|\leq 20,\,|\alpha|\leq 3}|\zeta|^{|\beta|}|(D^\beta_\zeta D^\alpha_x a)(x,\zeta)|.
\end{split}
\end{equation}
The index $l$ is called the {\it{order}} of the symbol, and it measures the contribution of the symbol in terms of derivatives on $f$. Notice that we have the simple product rule
\begin{equation}\label{nor2}
\|ab\,\|_{\mathcal{L}^p_{l_1+l_2}}\lesssim \|a\|_{\mathcal{L}^q_{l_1}}\|b\|_{\mathcal{L}^r_{l_2}},\qquad 1/p=1/q+1/r.
\end{equation}

We will use a simple lemma concerning paradifferential operators. See, for example, \cite[Section 3]{DeIoPa} or \cite[Appendix A]{DeIoPaPu} for a longer discussion and proofs. 

We start with a lemma relating products of symbols and composition of operators.

\begin{lemma}\label{PropProd} (i) If $1/p=1/q+1/r$ and $k\in\mathbb{Z}$, and $l\in [-10,10]$ then
\begin{equation}
\label{LqBdTa}
\Vert P_kT_af\Vert_{L^p}\lesssim 2^{lk_+}\Vert a\Vert_{\mathcal{L}^q_l}\Vert P_{[k-2,k+2]}f\Vert_{L^r}.
\end{equation}

(ii) For any symbols $a,b$ let $E(a,b):=T_aT_b-T_{ab}$. If $1/p=1/q_1+1/q_2+1/r$, $k\in\mathbb{Z}$, and $l_1,l_2\in [-4,4]$ then 
\begin{equation}
\label{LqBdTa2}
2^{k_+}\Vert P_kE(a,b)f\Vert_{L^p}\lesssim (2^{l_1k_+}\Vert a\Vert_{\mathcal{L}^{q_1}_{l_1}})(2^{l_2k_+}\Vert b\Vert_{\mathcal{L}^{q_2}_{l_2}})\cdot\Vert P_{[k-4,k+4]}f\Vert_{L^r}.
\end{equation}
\end{lemma}

The main point of the lemma the gain of one derivative in the left-hand side for the operator $T_aT_b-T_{ab}$. We remark that one could in fact gain two derivatives by subtracting the contribution of the Poisson bracket between the symbols $a$ and $b$, defined by
\begin{align*}
\{ a,b \} := \nabla_x a \nabla_\zeta b - \nabla_\zeta a \nabla_x b.
\end{align*}
We do not need a refinement of this type in this paper.

An important feature of the Weyl paradifferential calculus is self-adjointness of operators defined by real-valued symbols. More precisely:

\begin{lemma}\label{PropSym} If $\|a\|_{\mathcal{L}_0^\infty}<\infty$ is real-valued then $T_a$ is a bounded self-adjoint operator on $L^2$. Moreover, we have
\begin{equation}\label{Alu2}
\overline{T_af}=T_{a'}\overline{f},\quad\text{ where }\quad a'(y,\zeta):=\overline{a(y,-\zeta)}.
\end{equation} 
\end{lemma}

\subsection{Paralinearization of the reduced wave operator} Recall the modified metric $\widetilde{g}^{\mu\nu}:=m^{\mu\nu}+\widetilde{g}^{\mu\nu}_{\geq 1}$, where $\widetilde{g}^{\mu\nu}_{\geq 1}$ are as in \eqref{BoxWave1}, and the modified wave operator $\widetilde{\square}_{\widetilde{g}}:=\widetilde{g}^{\mu\nu}\partial_\mu\partial_\nu$. 

Notice that $\widetilde{g}^{00}=-1$. We first define the main symbols
\begin{equation}\label{ParaSymbolsW}
\mathcal{D}_{wa}:=(\widetilde{g}^{0j}\zeta_j)^2+\widetilde{g}^{jk}\zeta_j\zeta_k,\qquad \Sigma_{wa}:=\sqrt{\mathcal{D}_{wa}}-\widetilde{g}^{0j}\zeta_j,\qquad \sigma_{wa}:=\sqrt{\mathcal{D}_{wa}}+\widetilde{g}^{0j}\zeta_j,
\end{equation}
and the main symbols for the Klein-Gordon components,
\begin{equation}\label{ParaSymbolsKG}
\mathcal{D}_{kg}:=(\widetilde{g}^{0j}\zeta_j)^2+1+\widetilde{g}^{jk}\zeta_j\zeta_k,\qquad \Sigma_{kg}:=\sqrt{\mathcal{D}_{kg}}-\widetilde{g}^{0j}\zeta_j,\qquad \sigma_{kg}:=\sqrt{\mathcal{D}_{kg}}+\widetilde{g}^{0j}\zeta_j.
\end{equation}
These definitions are related to the paradifferential identities in Proposition \ref{ParaBox}.

Using the formulas \eqref{BoxWave2.9} we can extract the linear and higher order components of the symbols $\Sigma_{wa}$ and $\Sigma_{kg}$. Indeed, we have
\begin{equation*}
\mathcal{D}_{wa}=|\zeta|^2-h_{00}|\zeta|^2-h_{jk}\zeta_j\zeta_k+|\zeta|^2\mathcal{D}_{wa}^{\geq 2},
\end{equation*}
where $\mathcal{D}_{wa}^{\geq 2}$ is a quadratic symbol of order $0$. Thus, in view of Lemma \ref{Ga},
\begin{equation}\label{ParaSymb1}
\begin{split}
&\Sigma_{wa}=\vert\zeta\vert\big(1+\Sigma_{wa}^1+\Sigma_{wa}^{\ge2}\big),\\
&\Sigma_{wa}^1:=-(1/2)[h_{00}+2h_{0j}\widehat{\zeta_j}+h_{jk}\widehat{\zeta_j}\widehat{\zeta_k}\big],\\
&\Vert \Sigma_{wa}^{\ge 2}\Vert_{\mathcal{L}^\infty_0}+\Vert \partial_t\Sigma_{wa}^{\ge 2}\Vert_{\mathcal{L}^\infty_0}\lesssim \varepsilon_1^2\langle t\rangle^{2\delta^\prime -2},
\end{split}
\end{equation}
where $\widehat{\zeta_j}:=\zeta_j/|\zeta|$. Similarly, we have the decomposition of the Klein-Gordon symbol $\Sigma_{kg}$,
\begin{equation}\label{ParaSymb2}
\begin{split}
&\Sigma_{kg}=\langle\zeta\rangle\big(1+\Sigma_{kg}^1+\Sigma_{kg}^{\ge2}\big),\\
&\Sigma_{kg}^1:=-(1/2)\Big[h_{00}\frac{|\zeta|^2}{\langle\zeta\rangle^2}+2h_{0j}\widehat{\zeta_j}+h_{jk}\frac{\zeta_j\zeta_k}{\langle\zeta\rangle^2}\Big],\\
&\Vert \Sigma_{kg}^{\ge 2}\Vert_{\mathcal{L}^\infty_0}+\Vert \partial_t\Sigma_{kg}^{\ge 2}\Vert_{\mathcal{L}^\infty_0}\lesssim \varepsilon_1^2\langle t\rangle^{2\delta^\prime -2}.
\end{split}
\end{equation}

 The main result of this subsection is the following:

\begin{proposition}\label{ParaBox} For $\LL\in\mathcal{V}_n^q$, $(q,n)\leq (3,3)$, we define the "quasilinear variables"
\begin{equation}\label{MathcalU}
\begin{split}
\mathcal{U}^{\mathcal{L}h_{\al\be}}:=\left(\partial_t-iT_{\sigma_{wa}}\right)(\mathcal{L}h_{\al\be}),\qquad \mathcal{U}^{\mathcal{L}\psi}:=\left(\partial_t-iT_{\sigma_{kg}}\right)(\mathcal{L}\psi).
\end{split}
\end{equation}

(i) Then, for any $t\in[0,T]$, $k\in\mathbb{Z}$, and $f\in\{\LL h_{\al\be},\LL\psi\}$,
\begin{equation}\label{DefParaUnknown}
\Vert P_k(\mathcal{U}^f-U^f)(t)\Vert_{L^2}\lesssim \varepsilon_1\min\{1,2^k\langle t\rangle\}^{1-\delta}\langle t\rangle^{-1+31\delta}\cdot \Vert P_{[k-2,k+2]}U^f(t)\Vert_{L^2}.
\end{equation}

(ii) Moreover, we have
\begin{equation}\label{ParaBox1}
\begin{split}
\left(\partial_t+iT_{\Sigma_{wa}}\right)\mathcal{U}^f&=-\widetilde{g}^{\mu\nu}\partial_\mu\partial_\nu f+\Pi_{wa}[U^f],\qquad\text{ for }f\in\{\LL h_{\al\be}\},\\
\left(\partial_t+iT_{\Sigma_{kg}}\right)\mathcal{U}^f&=(-\widetilde{g}^{\mu\nu}\partial_\mu\partial_\nu+1) f+\Pi_{kg}[U^f],\qquad\text{ for }f\in\{\LL\psi\},\\
\end{split}
\end{equation}
where the remainder terms $\Pi_\ast[U^f]$ satisfy, for $\ast\in\{wa,kg\}$, 
\begin{equation}\label{PerturbPara}
\begin{split}
\Pi_\ast[U^f]&=\sum_{h\in\{h_{\al\be}\},\,\iota_1,\iota_2\in\{+,-\}}I_{\mathfrak{m}^{h,\ast}_{\iota_1\iota_2}}[U^{h,\iota_1},U^{f,\iota_2}]+\mathcal{C}_\ast[U^f].
\end{split}
\end{equation}
Here $\mathfrak{m}^{h,\ast}_{\iota_1\iota_2}$ are multipliers in $\mathcal{M}$ satisfying the additional bounds
\begin{equation}\label{ParaBox4.5}
\|\mathcal{F}^{-1}\{\mathfrak{m}^{h,\ast}_{\iota_1\iota_2}(\theta,\eta)\varphi_{k_1}(\theta)\varphi_{k_2}(\eta)\}\|_{L^1}\lesssim \min\{1,2^{k_2-k_1}\}
\end{equation}
for any $k_1,k_2\in\mathbb{Z}$, and the cubic remainders satisfy
\begin{equation}\label{PerturbPara2}
\begin{split}
\Vert P_k\mathcal{C}_{wa}[U^f](t)\Vert_{L^2}&\lesssim \varepsilon_1^2\langle t\rangle^{-5/4}2^{k/2}2^{-N(n)k^+}b_k(q,n;t)\qquad\text{ for }f\in\{\LL h_{\al\be}\},\\
\Vert P_k\mathcal{C}_{kg}[U^f](t)\Vert_{L^2}&\lesssim \varepsilon_1^2\langle t\rangle^{-5/4}2^{-N(n)k^+}b_k(q,n;t)\qquad\text{ for }f\in\{\LL\psi\},
\end{split}
\end{equation}
where $b_k(q,n;t)$ are the frequency envelope coefficients defined in \eqref{nake15}.

(iii) As a consequence, with $Y'(0)=Y'(1)=2$ and $Y'(2)=Y'(3)=35$, we have\footnote{These bounds illustrate the main gains one can expect by using quasilinear profiles instead of linear profiles: there are no derivative losses in \eqref{ParaBox5} and \eqref{ParaBox5.02}, and a smaller loss in terms of time decay when $n\leq 1$ in \eqref{ParaBox5.01}, compared to the bounds on $(\partial_t+i\Lambda_{wa})P_kU^{\LL h_{\al\be}}=\LL\mathcal{N}^h_{\al\be}$ and $(\partial_t+i\Lambda_{wa})P_kU^{\LL \psi}=\LL\mathcal{N}^\psi$.}
\begin{equation}\label{ParaBox5}
\big\|\left(\partial_t+iT_{\Sigma_{wa}}\right)P_k\mathcal{U}^{\LL h_{\al\be}}(t)\big\|_{L^2}\lesssim \varepsilon_1\langle t\rangle^{-1+\delta'}2^{k/2}2^{-N(n)k^+}b_k(q,n;t),
\end{equation}
\begin{equation}\label{ParaBox5.01}
\begin{split}
\big\|\left(\partial_t+iT_{\Sigma_{wa}}\right)P_k\mathcal{U}^{\LL h_{\al\be}}(t)\big\|_{L^2}\lesssim \varepsilon_1^2\langle t\rangle^{-1+H(q,n)\delta+Y'(n)\delta}2^{k/2}2^{-\widetilde{N}(n)k^++7k^+},
\end{split}
\end{equation}
\begin{equation}\label{ParaBox5.02}
\big\|\left(\partial_t+iT_{\Sigma_{kg}}\right)P_k\mathcal{U}^{\LL\psi}(t)\big\|_{L^2}\lesssim \varepsilon_1\langle t\rangle^{-1+\delta'}2^{-N(n)k^+}b_k(q,n;t).
\end{equation}
\end{proposition}

\begin{proof} (i) Using \eqref{ParaSymb1} we have 
\begin{equation*}
\sigma_{wa}=\Sigma_{wa}+2\widetilde{g}^{0j}_{\geq 1}\zeta_j=|\zeta|+|\zeta|\big(\Sigma_{wa}^1+\Sigma_{wa}^{\geq 2}+2\widetilde{g}^{0j}_{\geq 1}\widehat{\zeta_j}\big).
\end{equation*}
The bounds \eqref{DefParaUnknown} follow using \eqref{LqBdTa} and the $L^\infty$ estimates \eqref{wws1} when $f\in\{\LL h_{\al\be}\}$. The proof is similar when $f\in\{\LL \psi\}$ using \eqref{ParaSymb2} instead of \eqref{ParaSymb1}.

(ii) We compute, using the definitions,
\begin{equation}\label{sda1}
\begin{split}
\left(\partial_t+iT_{\Sigma_{wa}}\right)\left(\partial_t-iT_{\sigma_{wa}}\right)f&=\left[\partial_t\partial_t+iT_{\Sigma_{wa}-\sigma_{wa}}\partial_t+T_{\Sigma_{wa}\sigma_{wa}}\right]f\\
&-iT_{\partial_t\sigma_{wa}}f+\left(T_{\Sigma_{wa}}T_{\sigma_{wa}}-T_{\Sigma_{wa}\sigma_{wa}}\right)f,
\end{split}
\end{equation}
and, recalling that $\widetilde{g}^{00}=-1$,
\begin{equation}\label{sda2}
-\widetilde{g}^{\mu\nu}\partial_\mu\partial_\nu f=\big[\partial_t\partial_t-2iT_{\widetilde{g}^{0j}\zeta_j}\partial_t+T_{\widetilde{g}^{jk}\zeta_j\zeta_k}\big]f+E_1+E_2,
\end{equation}
where
\begin{equation*}
\begin{split}
\widehat{E_1}(\xi)&:=\frac{2i}{8\pi^3}\int_{\mathbb{R}^3}\Big[\frac{\xi_j+\eta_j}{2}\chi_0\Big(\frac{\vert\xi-\eta\vert}{\vert\xi+\eta\vert}\Big)-\eta_j\Big]\widehat{\widetilde{g}^{0j}_{\geq 1}}(\xi-\eta)\widehat{\partial_tf}(\eta)d\eta,\\
\widehat{E_2}(\xi)&:=\frac{-1}{8\pi^3}\int_{\mathbb{R}^3}\Big[\frac{(\xi_j+\eta_j)(\xi_k+\eta_k)}{4}\chi_0\Big(\frac{\vert\xi-\eta\vert}{\vert\xi+\eta\vert}\Big)-\eta_j\eta_k\Big]\widehat{\widetilde{g}^{jk}_{\geq 1}}(\xi-\eta)\widehat{f}(\eta)d\eta.
\end{split}
\end{equation*}
In view of the definitions \eqref{ParaSymbolsW}, we have $\Sigma_{wa}-\sigma_{wa}=-2\widetilde{g}^{0j}\zeta_j$ and $\Sigma_{wa}\cdot\sigma_{wa}=\widetilde{g}^{jk}\zeta_j\zeta_k$. Thus the main terms in \eqref{sda1}--\eqref{sda2} are the same, and the identities in the first line of \eqref{ParaBox1} follow easily by extracting the quadratic terms (which do not have derivative loss and can be written as claimed in \eqref{PerturbPara}). The cubic and higher order remainders can be bounded as claimed in \eqref{PerturbPara2}, using the bounds in Lemma \ref{Ga}.

The analysis of the Klein-Gordon terms is similar, using the identities $\Sigma_{kg}-\sigma_{kg}=-2\widetilde{g}^{0j}\zeta_j$ and $\Sigma_{kg}\cdot\sigma_{kg}=1+\widetilde{g}^{jk}\zeta_j\zeta_k$, see \eqref{ParaSymbolsKG}.

We remark that we could obtain more information on the remainder terms, which consist mostly of null quadratic terms and cubic terms, but this is not necessary for our purpose.

(iii) In view of \eqref{sac1} and the definitions \eqref{BoxWave1} we have
\begin{equation}\label{ParaBox6}
\begin{split}
-\widetilde{g}^{\mu\nu}&\partial_\mu\partial_\nu \LL h_{\al\be}=(\partial_t^2-\Delta)\LL h_{\al\be}-\widetilde{g}^{\mu\nu}_{\geq 1}\partial_\mu\partial_\nu \LL h_{\al\be}\\
&=\LL(\widetilde{g}^{\mu\nu}_{\geq 1}\partial_\mu\partial_\nu h_{\al\be})-\widetilde{g}^{\mu\nu}_{\geq 1}\partial_\mu\partial_\nu \LL h_{\al\be}+\LL\big\{(1-g_{\geq 1}^{00})^{-1}[\mathcal{KG}_{\al\be}-F^{\geq 2}_{\al\be}(g,\partial g)]\big\}.
\end{split}
\end{equation}
To prove the estimate \eqref{ParaBox5} we use the identity in the first line of \eqref{ParaBox1}. Notice that the expression $(1-g_{\geq 1}^{00})^{-1}[\mathcal{KG}_{\al\be}-F^{\geq 2}_{\al\be}(g,\partial g)]$ in the right-hand side or \eqref{ParaBox6} is a sum of terms in $\mathcal{QU}_0$, see \eqref{sho3.9}. The contribution of these terms is therefore bounded as claimed, due to \eqref{nake15.5}. The cubic terms $\mathcal{C}[U^{\LL h_{\al\be}}]$ also satisfy acceptable estimates, due to \eqref{PerturbPara2}. Therefore it remains to prove that
\begin{equation}\label{ParaBox7}
\begin{split}
\|P_kI[U^{h_1,\iota_1},U^{\LL h_2,\iota_2}](t)\|_{L^2}&\lesssim \varepsilon_1\langle t\rangle^{-1+\delta'}2^{k/2}2^{-N(n)k^+}b_k(q,n;t),\\
\|P_kI[\LL_1\widetilde{g}^{\mu\nu}_{\geq 1},\partial_\al\partial_\be \LL_2h_2](t)\|_{L^2}&\lesssim \varepsilon_1^2\langle t\rangle^{-1+\delta'}2^{k/2}2^{-N(n)k^+}2^{-\delta|k|},
\end{split}
\end{equation}
for any $h_1,h_2\in\{h_{\al\be}\}$, $\iota_1,\iota_2\in\{+,-\}$, $\mu,\nu,\al,\be\in\{0,1,2,3\}$, $\LL_1\in\mathcal{V}^{q_1}_{n_1}$, $\LL_2\in\mathcal{V}^{q_2}_{n_2}$, $(q_1,n_1)+(q_2,n_2)\leq (q,n)$, and $n_2<n$. 

For $n\geq 2$, the estimates \eqref{ParaBox5.01} follow from the formulas \eqref{ParaBox1}, \eqref{ParaBox6}, and the proofs of the $L^2$ estimates in section \ref{NullStruc} such as \eqref{sac32}, \eqref{sac62}--\eqref{sac62.5}, and \eqref{sac3.00}. For $n\in\{0,1\}$ the estimates \eqref{ParaBox5.01} follow from the improved estimates in Lemma \ref{Onev1} and cubic estimates such as \eqref{sac3.00}. If $n=1$ the term $\LL(\widetilde{g}^{\mu\nu}_{1}\partial_\mu\partial_\nu h_{\al\be})-\widetilde{g}^{\mu\nu}_{1}\partial_\mu\partial_\nu \LL h_{\al\be}$ in \eqref{ParaBox6} can be estimated using the bounds \eqref{onev8} when the undifferentiated factor carries the vector-field and its frequency is small.

The proofs of bounds \eqref{ParaBox7} are similar to some of the proofs in section \ref{DecompositionNonlin}, such as the proofs of \eqref{nake15.6} and \eqref{nake15.7}. Indeed, for the bounds in the first line, we decompose the input functions dyadically in frequency, then use $L^2\times L^\infty$ estimates for the $\mathrm{Low}\times\mathrm{High}\to\mathrm{High}$ interactions and for the $\mathrm{High}\times\mathrm{High}\to\mathrm{Low}$ interactions when $k\geq 0$; the remaining $\mathrm{High}\times\mathrm{High}\to\mathrm{Low}$ interactions when $k\leq 0$ can be bounded using \eqref{box2} and \eqref{box3.5}. The proof of the bounds in the second line of \eqref{ParaBox7} is similar (compare with \eqref{nake15.7}).

The estimates \eqref{ParaBox5.02} are similar. We start from the identities
\begin{equation*}
\begin{split}
(-\widetilde{g}^{\mu\nu}\partial_\mu\partial_\nu+1) \LL\psi&=(\partial_t^2-\Delta+1)\LL \psi-\widetilde{g}^{\mu\nu}_{\geq 1}\partial_\mu\partial_\nu \LL\psi\\
&=\LL(\widetilde{g}^{\mu\nu}_{\geq 1}\partial_\mu\partial_\nu\psi)-\widetilde{g}^{\mu\nu}_{\geq 1}\partial_\mu\partial_\nu \LL\psi-\LL\big\{(1-g_{\geq 1}^{00})^{-1}g_{\geq 1}^{00}\psi\big\},
\end{split}
\end{equation*}
which follow from \eqref{gb10}. Then we use the identities in the second line of \eqref{ParaBox1}, and notice that all the resulting terms that need to be estimated are quadratic or higher order, and do not lose derivatives. The desired bounds are similar to some of the bounds we proved in section \ref{DecompositionNonlin}, such as \eqref{nake15.8} and \eqref{nake15.66}, using also the inequalities \eqref{ParaBox4.5} when $\LL=Id$. to compensate for the lower regularity of $U^h$ compared to $U^\psi$.
\end{proof}

\section{Pure wave interactions}\label{En2purewave}

In this section we consider Wave$\times$Wave$\times$Wave interactions, and prove the bounds \eqref{EEEst3}--\eqref{EEcom1} in Proposition \ref{totbounds}.

\subsection{Null interactions} We start with a lemma concerning null interactions:

\begin{lemma}\label{NullLHHWWWS}
With the assumptions of Proposition \ref{totbounds}, we have, for any $m\in\{0,\ldots,L+1\}$,
\begin{equation}\label{NullLHHWWWF2}
\begin{split}
\sum_{k,k_1,k_2\in\mathbb{Z}}&2^{N(n)k^+-k/2}(2^{N(n)k_1^+-k_1/2}+2^{N(n)k_2^+-k_2/2}+2^{N(n)k^+-k/2})\\
&\Big|\int_{J_m}q_m(s)\mathcal{G}_{\mathfrak{n}_{\iota_1\iota_2}}\big[P_{k_1}U^{\mathcal{L}_1h_1,\iota_1}(s),P_{k_2}U^{\mathcal{L}_2h_2,\iota_2}(s),P_kU^{\mathcal{L}h,\iota}(s)\big]\,ds\Big|\lesssim \varepsilon_1^32^{-\delta m}, 
\end{split}
\end{equation}
where $\mathfrak{n}_{\iota_1\iota_2}\in\mathcal{M}^{null}_{\iota_1\iota_2}$ is a null symbol and the operators $\mathcal{G}_{n_{\iota_1\iota_2}}$ are defined as in \eqref{DefI}. 

Moreover, if $\mathfrak{q}_{\iota_1\iota_2}$ is a double-null symbol of the form
\begin{equation}\label{doubnull}
\mathfrak{q}_{\iota_1\iota_2}(\theta,\eta)=(\iota_1\theta_i/|\theta|-\iota_2\eta_i/|\eta|)(\iota_1\theta_j/|\theta|-\iota_2\eta_j/|\eta|)m_1(\theta,\eta),
\end{equation}
where $i,j\in\{1,2,3\}$, $m_1\in\mathcal{M}$, and if $n_2<n$ then we also have
\begin{equation}\label{NullLHHWWWF}
\begin{split}
\sum_{k,k_1,k_2\in\mathbb{Z}}2^{2N(n)k^+-k+k_2-k_1}2^{2\ga(m+k^-)}&\\
\Big|\int_{J_m}q_m(s)\mathcal{G}_{\mathfrak{q}_{\iota_1\iota_2}}\big[P_{k_1}U^{\mathcal{L}_1h_1,\iota_1}(s),&P_{k_2}U^{\mathcal{L}_2h_2,\iota_2}(s),P_kU^{\mathcal{L}h,\iota}(s)\big]\,ds\Big|\lesssim \varepsilon_1^32^{2H(q,n)\delta m}.
\end{split}
\end{equation}
\end{lemma}

It is easy to see that this lemma gives four of the bounds in Proposition \ref{totbounds}:

\begin{corollary}\label{thrbo}
The estimates \eqref{EEEst3}, \eqref{EEEst4}, \eqref{EEEst4.4}, and \eqref{EEcom1} hold.
\end{corollary}

\begin{proof}[Proof of Corollary \ref{thrbo}] The bounds \eqref{EEEst4} and \eqref{EEEst4.4} clearly follow from \eqref{NullLHHWWWF2}. Notice also that all the symbols in \eqref{BoxWave7}, thus all the symbols $\mathfrak{q}_{\iota_1\iota_2}^{\ast,wa}$ in (1) of Proposition \ref{totbounds}, are double-null, so the bounds \eqref{EEEst3} follow from \eqref{NullLHHWWWF}.

To prove \eqref{EEcom1} we calculate, as in \eqref{BoxWave5}--\eqref{BoxWave5.2},
\begin{equation*}
\begin{split}
\widetilde{g}^{\mu\nu}_{\geq 1}\partial_\mu\partial_\nu H&=-[\delta_{jk}+R_jR_k]F\cdot \partial_j\partial_kH-\big\{2R_jR_0\uF\cdot\partial_j\partial_0H+(\delta_{jk}-R_jR_k)\uF\cdot \partial_j\partial_kH\big\}\\
&+\big\{2\in_{jmn}R_m\omega_n\cdot\partial_j\partial_0H+(\in_{jln}R_k+\in_{kln}R_j)R_lR_0\omega_n\cdot\partial_j\partial_kH\big\}\\
&-\in_{jpm}\in_{kqn}R_pR_q\va_{mn}\cdot\partial_j\partial_kH-2R_jR_0\tau\cdot\partial_j\partial_0H+\sum_{(\mu,\nu)\neq (0,0)}\widetilde{G}^{\mu\nu}_{\geq 2}\cdot \partial_\mu\partial_\nu H,
\end{split}
\end{equation*}
where $H\in\{\LL h_2,\,P^n_{wa}\LL h_2\}$. The quadratic coefficients $\widetilde{G}^{\mu\nu}_{\geq 2}$ are linear combinations of expressions of the form $R^a|\nabla|^{-1}(G_{\geq 1}\partial_\rho h)$ where $R^a=R_1^{a_1}R_2^{a_2}R_3^{a_3}$ and $G_{\geq 1}$ are as in Definition \ref{nake4}. Therefore, as in Proposition \ref{gooddec}, we have
\begin{equation}\label{sori1}
\begin{split}
P_{wa}^n&[\widetilde{g}^{\mu\nu}_{\geq 1}\partial_\mu\partial_\nu \LL h_2]-\widetilde{g}^{\mu\nu}_{\geq 1}\partial_\mu\partial_\nu (P_{wa}^n\LL h_2)=\sum_{G\in\{F,\uF,\omega_n,\va_{mn}\}}\sum_{\iota_1,\iota_2\in\{+,-\}}\\
&\big\{P_{wa}^n I_{\mathfrak{q}^{G,wa}_{\iota_1\iota_2}}[|\nabla|^{-1}U^{G,\iota_1},|\nabla|U^{\LL h_{2},\iota_2}]-I_{\mathfrak{q}^{G,wa}_{\iota_1\iota_2}}[|\nabla|^{-1}U^{G,\iota_1},P^n_{wa}|\nabla|U^{\LL h_{2},\iota_2}]\big\}\\
&-2\{P_{wa}^n[R_jR_0\tau\cdot\partial_j\partial_0\LL h_2]-R_jR_0\tau\cdot\partial_j\partial_0(P^n_{wa}\LL h_2)\}\\
&+\sum_{(\mu,\nu)\neq (0,0)}\{P^n_{wa}[\widetilde{G}^{\mu\nu}_{\geq 2}\cdot \partial_\mu\partial_\nu \LL h_2]-\widetilde{G}^{\mu\nu}_{\geq 2}\cdot \partial_\mu\partial_\nu (P^n_{wa}\LL h_2)\},
\end{split}
\end{equation}
where the null multipliers $\mathfrak{q}^{G,wa}_{\iota_1\iota_2}$ are defined in \eqref{BoxWave7}. We notice that
\begin{equation}\label{sori1.5}
\begin{split}
&P_{wa}^n I^0_{\mathfrak{q}^{G,wa}_{\iota_1\iota_2}}[|\nabla|^{-1}U^{G,\iota_1},|\nabla|U^{\LL h_{2},\iota_2}]-I^0_{\mathfrak{q}^{G,wa}_{\iota_1\iota_2}}[|\nabla|^{-1}U^{G,\iota_1},P^n_{wa}|\nabla|U^{\LL h_{2},\iota_2}]\\
&=P_{wa}^n I^0_{\mathfrak{n}^{1}_{\iota_1\iota_2}}[U^{G,\iota_1},U^{\LL h_{2},\iota_2}]+I^0_{\mathfrak{n}^{2}_{\iota_1\iota_2}}[U^{G,\iota_1},P_{wa}^n U^{\LL h_{2},\iota_2}],
\end{split}
\end{equation}
where
\begin{equation}\label{sori1.7}
\begin{split}
&\mathfrak{n}^{1}_{\iota_1\iota_2}(\theta,\eta):=\mathfrak{q}^{G,wa}_{\iota_1\iota_2}(\theta,\eta)\Big\{\varphi_{\leq -4}(|\theta|/|\eta|)\frac{P^n_{wa}(\eta+\theta)-P^n_{wa}(\eta)}{P^n_{wa}(\eta+\theta)}\frac{|\eta|}{|\theta|}+\varphi_{>-4}(|\theta|/|\eta|)\frac{|\eta|}{|\theta|}\Big\},\\
&\mathfrak{n}^{2}_{\iota_1\iota_2}(\theta,\eta):=-\mathfrak{q}^{G,wa}_{\iota_1\iota_2}(\theta,\eta)\varphi_{>-4}(|\theta|/|\eta|)\frac{|\eta|}{|\theta|}.
\end{split}
\end{equation}
These are null multipliers as in (2) of Proposition \ref{totbounds}, thus the resulting integrals are similar to the integrals in \eqref{EEEst4}--\eqref{EEEst4.4} and are suitably bounded due to \eqref{NullLHHWWWF2}. 

We claim now that the remaining terms in \eqref{sori1} are cubic-type acceptable errors satisfying
\begin{equation}\label{sori2}
\big\|P_{wa}^n[H\cdot\partial_\mu\partial_\nu\LL h_2](s)-H(s)\cdot\partial_\mu\partial_\nu(P^n_{wa}\LL h_2)(s)\big\|_{L^2}\lesssim \varep_1^2\langle s\rangle^{-1}
\end{equation}
for $(\mu,\nu)\neq (0,0)$ and $s\in[0,T]$, where $H$ is either $R_jR_0\tau$ or $\widetilde{G}^{\mu\nu}_{\geq 2}$. This would clearly suffice to complete the proof of \eqref{EEcom1}. To prove the bounds \eqref{sori2} we notice that $|\nabla| H$ is of the form $R^a|\nabla|^{-1}\mathcal{N}$, for some $\mathcal{N}\in\mathcal{QU}$ (see \eqref{sho3.9} and use \eqref{zaq11}). We decompose as in \eqref{sori1.7}, thus for \eqref{sori2} it suffices to prove that 
\begin{equation}\label{sori3}
\big\|P_{wa}^nI[|\nabla|^{-1}\mathcal{N},U^{\LL h_2,\iota_2}](s)\big\|_{L^2}+\big\|I[|\nabla|^{-1}\mathcal{N},P^n_{wa}U^{\LL h_2,\iota_2}](s)\big\|_{L^2}\lesssim \varep_1^2\langle s\rangle^{-1},
\end{equation}
for any $\iota_2\in\{+,-\}$ and bilinear operators $I$ as in \eqref{abc36.1}. The bound on the first term follows from the first inequality in \eqref{nake15.6}, while the second term can be bounded easily using \eqref{nake16.0} and the assumption $\|P^n_{wa}U^{\LL h_2}(s)\|_{L^2}\lesssim \varep_1\langle s\rangle^{\delta'}$. This completes the proof of \eqref{sori2} and \eqref{EEcom1}.
\end{proof}

\begin{proof}[Proof of Lemma \ref{NullLHHWWWS}] We remark first that the two estimates are somewhat similar, except that the bounds \eqref{NullLHHWWWF} involve stronger null multipliers (double-null), but we have to deal with an additional anti-derivative which causes significant difficulties at low frequencies.  

The contributions of very small frequencies can sometimes be bounded using just $L^2$ norms,
\begin{equation}\label{TrivialEstSmallFreq1}
\begin{split}
&\Big|\mathcal{G}_{\mathfrak{m}}\big[P_{k_1}U^{\mathcal{L}_1h_1,\iota_1}(s),P_{k_2}U^{\mathcal{L}_2h_2,\iota_2}(s),P_kU^{\mathcal{L}h,\iota}(s)\big]\Big|\\
&\lesssim 2^{\underline{k}}2^{\frac{\underline{k}+k_2+k+k_1}{2}}\cdot 2^{-\frac{k_1}{2}}\Vert P_{k_1}U^{\mathcal{L}_1h_1}(s)\Vert_{L^2}\cdot 2^{-\frac{k_2}{2}}\Vert P_{k_2}U^{\mathcal{L}_2h_2}(s)\Vert_{L^2}\cdot 2^{-\frac{k}{2}}\Vert P_{k}U^{\mathcal{L}h}(s)\Vert_{L^2}\\
&\lesssim\varepsilon_1^3 2^{2\underline{k}+\overline{k}}\cdot 2^{-\ga(k_1^-+k_2^-+k^-)}2^{-N(n_1)k_1^+-N(n_2)k_2^+-N(n)k^+}\cdot\langle s\rangle^{[H(\mathcal{L}_1)+H(\mathcal{L}_2)+H(\mathcal{L})]\delta-3\ga},
\end{split}
\end{equation}
provided that $\Vert \mathfrak{m}\Vert_{L^\infty}\lesssim 1$, where $\overline{k}:=\max\{k,k_1,k_2\}$, $\underline{k}:=\min\{k,k_1,k_2\}$, and $H(\LL)=H(q,n)$, $H(\LL_i)=H(q_i,n_i)$, $i\in\{1,2\}$. These bounds suffice to prove the desired conclusions if $2^m\lesssim 1$, so in the analysis below we may assume that $m\geq \delta^{-1}$.

{\bf Step 1.} We consider first the contribution of the triplets $(k,k_1,k_2)$ for which $k\leq \underline{k}+8$. In this case $2^{k_2-k_1}\lesssim 1$ and it suffices to focus on the harder estimates \eqref{NullLHHWWWF2}. In fact, we will prove the stronger bounds, for any $k,k_1,k_2\in\mathbb{Z}$ satisfying $k\leq \underline{k}+8$,
\begin{equation}\label{bnm12}
\begin{split}
&2^{N(n)k^+-k/2}(2^{N(n)k_1^+-k_1/2}+2^{N(n)k_2^+-k_2/2}+2^{N(n)k^+-k/2})\Big|\int_{J_m}q_m(s)\\
&\times\mathcal{G}_{\mathfrak{n}_{\iota_1\iota_2}}\big[P_{k_1}U^{\mathcal{L}_1h_1,\iota_1}(s),P_{k_2}U^{\mathcal{L}_2h_2,\iota_2}(s),P_kU^{\mathcal{L}h,\iota}(s)\big]\,ds\Big|\lesssim \varepsilon_1^32^{-\delta m}2^{-\delta(|k|+|k_1|+|k_2|)}.
\end{split}
\end{equation}

{\bf{Case 1.1.}} Assume first that $\iota_1=-\iota_2$. In this case, we notice that
\begin{equation}\label{trf1}
\Vert \mathcal{F}^{-1}\{\mathfrak{n}_{\iota_1\iota_2}\cdot\varphi_{kk_1k_2}\}\Vert_{L^1}\lesssim 2^{k-k_1}.
\end{equation}
If $\underline{k}\le -\delta^\prime m$ or if $\underline{k}\ge-\delta^\prime m$ and $\overline{k}\geq \delta'm$, then we can just use \eqref{EstimJ} with $(p_1,p_2,p)=(2,\infty,2)$  (if $n_2\leq n_1$) or $(p_1,p_2,p)=(\infty,2,2)$ (if $n_1\leq n_2$), and the estimates \eqref{vcx1} and \eqref{wws1}. On the other hand, if $k,k_1,k_2\in[-\delta^\prime m,\delta'm]$ then we may assume $m\leq L$ (the case $m=L+1$ is easier because $|J_m|\lesssim 1$) and decompose the multiplier into resonant and non-resonant contributions. More precisely, let $q_0=-8\delta^\prime m$ and
\begin{equation}\label{trf2}
\mathfrak{n}_{\iota_1\iota_2}=\mathfrak{n}_{\iota_1\iota_2}^r+\mathfrak{n}^{nr}_{\iota_1\iota_2},\qquad\mathfrak{n}_{\iota_1\iota_2}^r(\theta,\eta)=\varphi_{\le q_0}(\Xi_{\iota_1\iota_2}(\theta,\eta))\mathfrak{n}_{\iota_1\iota_2}(\theta,\eta),
\end{equation}
where $\Xi_{\iota_1\iota_2}$ is defined as in \eqref{par1}. 

We bound the contributions of the resonant parts $\mathfrak{n}_{\iota_1\iota_2}^r$ using the null structure, while for the non-resonant parts $\mathfrak{n}_{\iota_1\iota_2}^{nr}$ we use normal forms (the identity \eqref{NF2}). It follows from \eqref{vcx1}, \eqref{wws1}, \eqref{wer4.0} (or \eqref{nake15.5} if $l\geq 0$), and \eqref{plk2} (recalling also \eqref{roco2}) that
\begin{equation}\label{bnm1}
\begin{split}
\|P_l U^{\mathcal{K}h}(t)\|_{L^2}&\lesssim \varep_1\langle t\rangle^{\delta'/2}2^{l/2}2^{-\ga l^-}2^{-N(n')l^+},\\
\|P_l U^{\mathcal{K}h}(t)\|_{L^\infty}&\lesssim \varep_1\langle t\rangle^{\delta'/2-1}2^{l^-}2^{-N(n'+1)l^++2l^+}\min\{1,2^{l^-}\langle t\rangle\}^{1-\delta},
\end{split}
\end{equation}
and
\begin{equation}\label{bnm1.5}
\begin{split}
\|P_l (\partial_t+i\Lambda_{wa})U^{\mathcal{K}h}(t)\|_{L^2}&\lesssim \varep_1\langle t\rangle^{\delta'/2-1}2^{l/2}2^{-N(n')l^++l^+}\min\{1,2^{l^-}\langle t\rangle\},\\
\|P_l (\partial_t+i\Lambda_{wa})U^{\mathcal{K}h}(t)\|_{L^\infty}&\lesssim \varep_1\langle t\rangle^{4\delta'-2}2^{l^-}2^{-N(n'+1)l^++6l^+}\min\{1,2^{l^-}\langle t\rangle\},
\end{split}
\end{equation}
for any $l\in\mathbb{Z}$, $h\in\{h_{\al\be}\}$, $t\in[0,T]$, and $\mathcal{K}\in\mathcal{V}^{q'}_{n'}$, where the inequalities in the first lines hold for all $n'\leq 3$, while the inequalities in the second lines only hold for $n'\leq 2$. Thus, using \eqref{yip2} and \eqref{bnm1}, and recalling the null structure of the symbols $\mathfrak{n}_{\iota_1\iota_2}$, we have
\begin{equation}\label{bnm2}
\begin{split}
\Big|\mathcal{G}_{\mathfrak{n}_{\iota_1\iota_2}^r}&\big[P_{k_1}U^{\mathcal{L}_1h_1,\iota_1}(s),P_{k_2}U^{\mathcal{L}_2h_2,\iota_2}(s),P_kU^{\mathcal{L}h,\iota}(s)\big]\Big|\\
&\lesssim 2^{q_0}2^{\max(k_1,k_2)-k}\|P_{k_1}U^{\mathcal{L}_1h_1}(s)\|_{L^{p_1}}\|P_{k_2}U^{\mathcal{L}_2h_2}(s)\|_{L^{p_2}}\|P_kU^{\mathcal{L}h}(s)\|_{L^2}\\
&\lesssim \varep_1^32^{q_0}2^{2\delta'm}\langle s\rangle^{-1+2\delta'}2^{-N(n)k^++k/2}2^{-N(n)k_2^++k_2/2},
\end{split}
\end{equation}
assuming that $2^{-\delta'm}\lesssim 2^k\lesssim 2^{k_1}\approx 2^{k_2}\lesssim 2^{\delta'm}$, where $(p_1,p_2)=(2,\infty)$ if $n_1\geq n_2$ and $(p_1,p_2)=(\infty,2)\}$ if $n_1\leq n_2$. This suffices to estimate the resonant contributions as in \eqref{bnm12}.

For the non-resonant symbols $\mathfrak{n}_{\iota_1\iota_2}^{nr}$, we can use the normal form formulas \eqref{NF2}-\eqref{NF3} and the bounds \eqref{yip6}. Using \eqref{bnm1}--\eqref{bnm1.5} and estimating as in \eqref{bnm2} we have
\begin{equation}\label{bnm16}
\begin{split}
\Big|\mathcal{H}_{\mathfrak{n}_{\iota_1\iota_2}^{nr}}\big[P_{k_1}U^{\mathcal{L}_1h_1,\iota_1}(s),P_{k_2}U^{\mathcal{L}_2h_2,\iota_2}(s),P_kU^{\mathcal{L}h,\iota}(s)\big]\Big|\lesssim \varep_1^32^{-3q_0}2^{-0.9m},\\
\Big|\mathcal{H}_{\mathfrak{n}_{\iota_1\iota_2}^{nr}}\big[P_{k_1}U^{\mathcal{L}_1h_1,\iota_1}(s),P_{k_2}U^{\mathcal{L}_2h_2,\iota_2}(s),P_k(\partial_s+i\Lambda_{wa,\iota})U^{\mathcal{L}h,\iota}(s)\big]\Big|\lesssim \varep_1^32^{-3q_0}2^{-1.9m},\\
\Big|\mathcal{H}_{\mathfrak{n}_{\iota_1\iota_2}^{nr}}\big[P_{k_1}(\partial_s+i\Lambda_{wa,\iota_1})U^{\mathcal{L}_1h_1,\iota_1}(s),P_{k_2}U^{\mathcal{L}_2h_2,\iota_2}(s),P_kU^{\mathcal{L}h,\iota}(s)\big]\Big|\lesssim \varep_1^32^{-3q_0}2^{-1.9m},\\
\Big|\mathcal{H}_{\mathfrak{n}_{\iota_1\iota_2}^{nr}}\big[P_{k_1}U^{\mathcal{L}_1h_1,\iota_1}(s),P_{k_2}(\partial_s+i\Lambda_{wa,\iota_2})U^{\mathcal{L}_2h_2,\iota_2}(s),P_kU^{\mathcal{L}h,\iota}(s)\big]\Big|\lesssim \varep_1^32^{-3q_0}2^{-1.9m},
\end{split}
\end{equation}
for any $s\in J_m$ (recall that $2^{|k|}+2^{|k_1|}+2^{|k_2|}\lesssim 2^{\delta'm}$). This completes the proof of \eqref{bnm12}.

{\bf{Case 1.2.}} Assume now that $\iota_1=\iota_2$. If $k\geq \min\{k_1,k_2\}-10$ then $2^k\approx 2^{k_1}\approx 2^{k_2}$, the bounds \eqref{trf1} still hold, and the same proof as before gives the desired bounds \eqref{bnm12}.

On the other hand, if $k\leq \min\{k_1,k_2\}-10$ then we have no resonant contributions, i.e. $\Xi_{\iota_1\iota_2}(\theta,\eta)\gtrsim 1$ in the support of the integral. So we can integrate by parts directly using \eqref{NF2}. If $2^k\leq 2^{-m/4}$ or if $2^{k_1}\gtrsim 2^{\delta'm}$ then we use only the $L^2$ bounds in the first lines of \eqref{bnm1} and \eqref{bnm1.5}, and estimate as in \eqref{TrivialEstSmallFreq1} using \eqref{par5}. On the other hand, if $-m/4\leq k\leq \min\{k_1,k_2\}-10\leq \delta'm$ then we can estimate the resulting terms in \eqref{NF3} as in \eqref{bnm16}, using $L^2\times L^2\times L^\infty$ bounds and \eqref{bnm1}--\eqref{bnm1.5}. The desired bounds \eqref{bnm12} follow.

{\bf Step 2.} We complete now the proof of \eqref{NullLHHWWWF2} by analyzing the contribution of the triplets $(k,k_1,k_2)$ for which $k\ge\underline{k}+8$. By symmetry, we may assume that $n_1\le n_2$. 

{\bf{Case 2.1.}} We assume first that $n_1\geq 1$. In this case $1\leq n_1\leq n_2\leq 2$, $n_2<n$, and we can still prove the strong bounds \eqref{bnm12}. Indeed, if $\underline{k}\le -2\delta^\prime m$ or if $\underline{k}\ge-2\delta^\prime m$ and $\overline{k}\geq 2\delta'm$, then we can use \eqref{EstimJ} (with the lowest frequency placed in $L^\infty$) and the estimates \eqref{bnm1}. On the other hand, if $k,k_1,k_2\in[-2\delta'm,2\delta'm]$ then we still decompose $\mathfrak{n}_{\iota_1\iota_2}=\mathfrak{n}_{\iota_1\iota_2}^r+\mathfrak{n}^{nr}_{\iota_1\iota_2}$ as in \eqref{trf2}, with $q_0=-8\delta'm$, and apply \eqref{yip2} and \eqref{yip6}, together with \eqref{bnm1}--\eqref{bnm1.5}. The desired bounds follow easily by estimating the resonant and the non-resonant contributions as in \eqref{bnm2}--\eqref{bnm16}.

{\bf{Case 2.2.}} Assume now that $n_1=0$, $n_2\geq 1$, and consider the contribution of the triplets $(k,k_1,k_2)$ for which $k_1\geq k_2$. We can still prove the strong bounds \eqref{bnm12} because there is no derivative loss in this case. If $n_2\leq 2$ then we can still estimate as in \eqref{bnm2}--\eqref{bnm16}, using $L^2\times L^2\times L^\infty$ bounds with the lowest frequency placed in $L^\infty$. On the other hand, if $n_2=3$ (thus $n=3$), then \eqref{bnm12} follows from \eqref{TrivialEstSmallFreq1} if $k_2\leq -3m/4$; if $k_2\geq -3m/4$ then we decompose the symbol as in \eqref{trf2} and estimate as in \eqref{bnm2}--\eqref{bnm16}, with the terms corresponding to the frequency $\approx 2^{k_1}$ always placed in $L^\infty$.

{\bf{Case 2.3.}} Finally, assume that $n_1=0$ and consider the contribution of the triplets $(k,k_1,k_2)$ for which $k_1\leq k_2$. The main issue here is the loss of derivative in normal form arguments. Let
\begin{equation}\label{bnm21}
b_{k,m}(q,n):=\Big\{\frac{1}{|J_m|}\int_{J_m}(b_k(q,n;s))^2\,ds\Big\}^{1/2},
\end{equation}
where the frequency envelope coefficients $b_k(q,n;s)$ are defined in \eqref{nake15}. We will show that
\begin{equation}\label{bnm20}
\begin{split}
&2^{N(n)k^+-k/2}(2^{N(n)k_1^+-k_1/2}+2^{N(n)k^+-k/2})\Big|\int_{J_m}q_m(s)\\
&\times\mathcal{G}_{\mathfrak{n}_{\iota_1\iota_2}}\big[P_{k_1}U^{h_1,\iota_1}(s),P_{k_2}U^{\mathcal{L}h_2,\iota_2}(s),P_kU^{\mathcal{L}h,\iota}(s)\big]\,ds\Big|\lesssim \varepsilon_12^{-\delta m}2^{-\delta |k_1|} (b_{k,m}(q,n))^2,
\end{split}
\end{equation}
provided that $|k-k_2|\leq 4$ and $\mathcal{L}\in\mathcal{V}^q_n$, $(q,n)\leq (3,3)$. This is slightly weaker than the bounds \eqref{bnm12} (since $b_{k,m}(q,n)\gtrsim \varep_12^{-\delta|k|}$ due to \eqref{nake15}), but still suffices to complete the proof of \eqref{NullLHHWWWF2} due to the square summability of the coefficients $b_{k,m}(q,n)$ (see \eqref{nake15.1}). 

To prove \eqref{bnm20} we notice that we can still use the normal form argument as in {\bf{Case 2.1}} if $2^k\lesssim 2^{m/10}$, or to control the contribution of the factor $2^{N(n)k_1^+-k_1/2}$ in the left-hand side. If $2^k\gtrsim 2^{m/10}$ and $|k_1|\geq \delta'm$ then we can use just $L^\infty\times L^2\times L^2$ estimates (with the $L^2$ bounds on the two high frequency terms coming from \eqref{nake15.2}) to prove \eqref{bnm20}. 

In the remaining case when $|k_1|\leq \delta'm$ the resonant contribution can be estimated as in \eqref{bnm2}. The bound on the non-resonant contribution requires an additional idea (the use of paradifferential calculus) to avoid the derivative loss in the application of the Poincar\'{e} normal form; the desired bounds follow from Lemma \ref{ParalinearSym} below.

{\bf{Step 3.}} We turn now to the proof of \eqref{NullLHHWWWF}. In view of \eqref{NullLHHWWWF2} it suffices to prove that
\begin{equation}\label{bnm30}
\begin{split}
2^{-k_1}\Big|\int_{J_m}q_m(s)\mathcal{G}_{\mathfrak{q}_{\iota_1\iota_2}}\big[P_{k_1}U^{\mathcal{L}_1h_1,\iota_1}(s),P_{k_2}U^{\mathcal{L}_2h_2,\iota_2}(s),P_kU^{\mathcal{L}h,\iota}(s)\big]\,ds\Big|\\
\lesssim \varepsilon_1^32^{-2N(n)k^+}2^{2H(\LL)\delta m}2^{-2\ga m}2^{-\ga(|k|+|k_1+m|)/4},
\end{split}
\end{equation}
for any triplet $(k,k_1,k_2)$ for which $k_1=\underline{k}\leq\overline{k}-10$. Using first \eqref{TrivialEstSmallFreq1}, the left-hand side of \eqref{bnm30} is bounded by
\begin{equation}\label{bnm31}
C\varep_1^32^{k_1+k}2^{-\ga k_1^--2\ga k^-}2^{-8k_1^+}2^{-N(n)k^+-N(n_2)k^+}2^{m+\delta m(H(\LL_1)+H(\LL_2)+H(\LL))-3\ga m}.
\end{equation}
Since $N(n_2)\geq N(n)+10$, this implies \eqref{bnm30} unless 
\begin{equation}\label{bnm32}
2^{k_1+m}2^{-2\ga(k_1+m)^-}2^{-8k^+}2^{(1-\delta)k^-}\gtrsim 2^{\delta m(H(\LL)-H(\LL_1)-H(\LL_2))}.
\end{equation}
It remains to prove \eqref{bnm30} for triplets $(k,k_1,k_2)$ for which \eqref{bnm32} holds. In particular, we may assume that $m\geq \delta^{-2}$ in the analysis below.

{\bf{Case 3.1.}} Assume first that $n_1< n$, thus $H(\LL)-H(\LL_1)-H(\LL_2)\geq 40$ (due to \eqref{SuperlinearH1}) and $k_1+m\geq 40\delta m$ (due to \eqref{bnm32}). We may also assume that $m\leq L$, and use \eqref{NF2}. We observe that $|\mathfrak{q}_{\iota_1\iota_2}(\xi-\eta,\eta)(\Phi_{\sigma\mu\nu}(\xi,\eta))^{-1}|\lesssim 2^{-k_1}$, where $(\sigma,\mu,\nu)=((wa,\iota),(wa,\iota_1),(wa,\iota_2))$ in the support of the integral, due to the double-null assumption \eqref{doubnull} and the bounds \eqref{par5}. Using just $L^2$ estimates, the left-hand side of \eqref{bnm30} is bounded by
\begin{equation*}
\begin{split}
C2^{3k_1/2}2^{-2k_1}\sup_{s\in J_m}\big\{&(\|P_{k_1}U^{\mathcal{L}_1h_1}(s)\|_{L^2}+2^m\|P_{k_1}\mathcal{L}_1\mathcal{N}^h_1(s)\|_{L^2})\|P_{k_2}U^{\mathcal{L}_2h_2}(s)\|_{L^2}\|P_kU^{\mathcal{L}h}(s)\|_{L^2}\\
&+2^m\|P_{k_1}U^{\mathcal{L}_1h_1}(s)\|_{L^2}\|P_{k_2}\mathcal{L}_2\mathcal{N}^h_2(s)\|_{L^2}\|P_kU^{\mathcal{L}h}(s)\|_{L^2}\\
&+2^m\|P_{k_1}U^{\mathcal{L}_1h_1}(s)\|_{L^2}\|P_{k_2}U^{\mathcal{L}_2h_2}(s)\|_{L^2}\|P_k\mathcal{L}\mathcal{N}^h(s)\|_{L^2}\big\},
\end{split}
\end{equation*}
where $\mathcal{N}^h_1,\mathcal{N}^h_2,\mathcal{N}^h$ are suitable components of the metric nonlinearities $\mathcal{N}^h_{\al\be}$. In view of \eqref{vcx1} and \eqref{wer4.0}, all the terms in the expression above are dominated by
\begin{equation*}
C\varep_1^32^{k}2^{\delta m(H(\LL_1)+H(\LL_2)+H(\LL))+36\delta m}2^{-2N(n)k^+-2k^+},
\end{equation*}
which suffices to prove \eqref{bnm30} in this case (due to \eqref{SuperlinearH1}).

The same argument also proves the desired bounds \eqref{bnm30} when $n_2=0$, $(q_1,n_1)=(q,n)$, and $|k|\geq\delta' m$.

{\bf{Case 3.2.}} Assume now that \eqref{bnm32} holds, $n_2=0$, $(q_1,n_1)=(q,n)$, $|k|\leq \delta'm$, and 
\begin{equation}\label{bnm35}
k_1+m\leq Y(q,n)\delta m-\delta m,
\end{equation}
where $Y(q,n)$ is defined as in \eqref{saur10}. In particular, we may assume that $-2\delta m\leq k_1+m$, due to \eqref{bnm32}. To prove \eqref{bnm30} in this case, we trivialize one vector-field using Lemma \ref{TrivialVF}. As in \eqref{saur12.2} we define
\begin{equation*}
U^\ast_{1,\leq J_1}:=P'_{k_1}(\varphi_{\leq J_1}\cdot P_{k_1}U^{\mathcal{L}_1h_1,\iota_1}),\,\,\,\,\,\,U^\ast_{1,> J_1}:=P'_{k_1}(\varphi_{> J_1}\cdot P_{k_1}U^{\mathcal{L}_1h_1,\iota_1}).
\end{equation*}
where $J_1$ denotes the largest integer satisfying $J_1\leq \max\{-k_1,m\}+\delta m/4$. Then
\begin{equation}\label{saur14rep}
\|U^\ast_{1,\leq J_1}(s)\|_{L^2}\lesssim\varep_12^{k_1/2}2^{H(q,n)\delta m-7\delta m/4},
\end{equation}
see Lemma \ref{TrivialVF}. We decompose also $P_{k_2}U^{h_2,\iota_2}(s)=U^{h_2,\iota_2}_{\leq J_2,k_2}(s)+U^{h_2,\iota_2}_{> J_2,k_2}(s)$ as in \eqref{on11.3}--\eqref{on11.36}, where $J_2$ is the largest integer satisfying $J_2\leq m/4$.

With $I_{\mathfrak{q}_{\iota_1\iota_2}}$ defined as in \eqref{abc36.1}, the left-hand side of \eqref{bnm30} is bounded by
\begin{equation*}
C2^{-k_1}2^m\sup_{s\in J_m}\big\|P_kI_{\mathfrak{q}_{\iota_1\iota_2}}[P_{k_1}U^{\mathcal{L}_1h_1,\iota_1}(s),P_{k_2}U^{h_2,\iota_2}(s)]\big\|_{L^2}\|P_{[k-2,k+2]}U^{\mathcal{L}h,\iota}(s)\|_{L^2}.
\end{equation*}
In view of \eqref{vcx1}, for \eqref{bnm30} it suffices to prove that, for any $s\in J_m$,
\begin{equation}\label{jik1}
2^{-k_1}\big\|P_kI_{\mathfrak{q}_{\iota_1\iota_2}}[P_{k_1}U^{\mathcal{L}_1h_1,\iota_1},P_{k_2}U^{h_2,\iota_2}](s)\big\|_{L^2}\lesssim \varepsilon_1^22^{-N(n)k^+-k^+}2^{-m+H(q,n)\delta m}2^{-\ga m}2^{-\ga|k_1+m|/4}.
\end{equation}

In view of \eqref{bnm32}, we may assume that $|k|\leq\delta' m$. Using just $L^2$ estimates, we have
\begin{equation}\label{jik2}
2^{-k_1}\big\|P_kI_{\mathfrak{q}_{\iota_1\iota_2}}[P_{k_1}U^{\mathcal{L}_1h_1,\iota_1},U^{h_2,\iota_2}_{> J_2,k_2}](s)\big\|_{L^2}\lesssim \varepsilon_1^22^{-N_0k^+}2^{-m/5}2^{k_1}.
\end{equation}
To bound the contribution of the profile $U^{h_2,\iota_2}_{\leq J_2,k_2}(s)$ we use \eqref{bil1}, \eqref{vcx1.2}, and \eqref{saur14rep}, thus
\begin{equation}\label{jik3}
\begin{split}
2^{-k_1}\big\|P_kI_{\mathfrak{q}_{\iota_1\iota_2}}[U^\ast_{1,\leq J_1},U^{h_2,\iota_2}_{\leq J_2,k_2}](s)\big\|_{L^2}&\lesssim 2^{-k_1/2}2^{-m}\|U^\ast_{1,\leq J_1}(s)\|_{L^2}2^{3k_2/2}\|\widehat{P_{k_2}U^{h_2}}\|_{L^\infty}\\
&\lesssim\varep_1^22^{-m}2^{H(q,n)\delta m-7\delta m/4}2^{\delta m}2^{-N_0k^++2k_+}.
\end{split}
\end{equation}
Finally, the bilinear interaction of $U^\ast_{1,\leq J_1}(s)$ and $U^{h_2,\iota_2}_{\leq J_2,k_2}(s)$ is negligible
\begin{equation}\label{jik4}
2^{-k_1}\big\|P_kI_{\mathfrak{q}_{\iota_1\iota_2}}[U^\ast_{1,\leq J_1},U^{h_2,\iota_2}_{\leq J_2,k_2}](s)\big\|_{L^2}\lesssim\varep_1^22^{-2m},
\end{equation}
using an approximate finite speed of propagation argument as in the proof of \eqref{bnm37}. The desired bounds \eqref{jik1} follow from \eqref{jik2}--\eqref{jik4}.

{\bf Step 4.} To prove \eqref{bnm30} in the remaining cases we need to use angular localization and integrate by parts in time. More precisely, we decompose
\begin{equation}\label{jik6}
\mathfrak{q}_{\iota_1\iota_2}=\sum_{b\leq 4}\mathfrak{q}^b_{\iota_1\iota_2},\qquad\mathfrak{q}^b_{\iota_1\iota_2}(\theta,\eta)=\varphi_{b}(\Xi_{\iota_1\iota_2}(\theta,\eta))\mathfrak{q}_{\iota_1\iota_2}(\theta,\eta),
\end{equation}
and define the associated operators $\mathcal{G}_{\mathfrak{q}_{\iota_1\iota_2}^b}$ as in \eqref{DefI}. For \eqref{bnm30} it suffices to prove that
\begin{equation}\label{jik7}
\begin{split}
2^{-k_1}\Big|\int_{J_m}q_m(s)\mathcal{G}_{\mathfrak{q}_{\iota_1\iota_2}^b}\big[P_{k_1}U^{\mathcal{L}_1h_1,\iota_1}(s),P_{k_2}U^{h_2,\iota_2}(s),P_kU^{\mathcal{L}h,\iota}(s)\big]\,ds\Big|\\
\lesssim \varepsilon_1^32^{-2N(n)k^+}2^{2H(q,n)\delta m}2^{-3\ga m}2^{\delta b},
\end{split}
\end{equation}
provided that $\mathcal{L},\mathcal{L}_1\in\mathcal{V}^q_n$, $b\leq 4$, $m\geq \delta^{-2}$, and
\begin{equation}\label{jik8}
k_1\geq -m+(Y(q,n)-1)\delta m,\qquad |k|\leq\delta'm,\qquad k_1\leq\min\{k,k_2\}-6.
\end{equation}

{\bf{Case 4.1.}} Assume first that $b\leq-3\delta'm$. We decompose $P_{k_2}U^{h_2,\iota_2}(s)=U^{h_2,\iota_2}_{\leq J_2,k_2}(s)+U^{h_2,\iota_2}_{> J_2,k_2}(s)$ as in \eqref{on11.3}--\eqref{on11.36}, where $J_2$ is the largest integer satisfying $J_2\leq m-\delta'm$. Using \eqref{nag1} and the double null assumption \eqref{doubnull} we have $\|\mathcal{F}^{-1}(q_{\iota_1\iota_2}^b)\|_{L^1}\lesssim 2^b$. Then we estimate, using \eqref{Bil31} and \eqref{vcx1},
\begin{equation}\label{jik11}
\begin{split}
&2^{-k_1}\Big|\int_{J_m}q_m(s)\mathcal{G}_{\mathfrak{q}_{\iota_1\iota_2}^b}\big[P_{k_1}U^{\mathcal{L}_1h_1,\iota_1}(s),U^{h_2,\iota_2}_{\leq J_2,k_2}(s),P_kU^{\mathcal{L}h,\iota}(s)\big]\,ds\Big|\\
&\lesssim 2^{-k_1}2^m\sup_{s\in J_m}\big\|P_kI_{\mathfrak{q}_{\iota_1\iota_2}^b}\big[P_{k_1}U^{\mathcal{L}_1h_1,\iota_1}(s),U^{h_2,\iota_2}_{\leq J_2,k_2}(s)\big]\big\|_{L^2}\big\|P_{[k-2,k+2]}U^{\mathcal{L}h,\iota}(s)\big\|_{L^2}\\
&\lesssim 2^{-k_1}2^m\cdot 2^b2^{k_1/2}2^{-m+\delta m}\sup_{s\in J_m}\|P_{k_1}U^{\mathcal{L}_1h_1,\iota_1}(s)\|_{L^2}\|U^{h_2,\iota_2}_{\leq J_2,k_2}(s)\|_{H^{0,1}_\Omega}\|P'_{k}U^{\mathcal{L}h,\iota}(s)\|_{L^2}\\
&\lesssim \varep_1^32^b2^{\delta'm}2^k2^{-2N(n)k^+}.
\end{split}
\end{equation}
Moreover, using $L^2$ estimates, the double null assumption \eqref{doubnull}, and \eqref{vcx1}--\eqref{vcx1.1}, we have
\begin{equation}\label{jik12}
\begin{split}
&2^{-k_1}\Big|\int_{J_m}q_m(s)\mathcal{G}_{\mathfrak{q}_{\iota_1\iota_2}^b}\big[P_{k_1}U^{\mathcal{L}_1h_1,\iota_1}(s),U^{h_2,\iota_2}_{> J_2,k_2}(s),P_kU^{\mathcal{L}h,\iota}(s)\big]\,ds\Big|\\
&\lesssim 2^{-k_1}2^m2^{3k_1/2}2^{2b}\sup_{s\in J_m}\|P_{k_1}U^{\mathcal{L}_1h_1,\iota_1}(s)\|_{L^2}\|U^{h_2,\iota_2}_{> J_2,k_2}(s)\|_{L^2}\|P'_{k}U^{\mathcal{L}h,\iota}(s)\|_{L^2}\\
&\lesssim \varep_1^32^{2b}2^{k_1^-}2^{2\delta'm}2^{-2N(n)k^+}.
\end{split}
\end{equation}
The desired bounds \eqref{jik7} follow from \eqref{jik11}--\eqref{jik12} if $b\leq-3\delta'm$.

{\bf{Case 4.2.}} Assume now that the inequalities \eqref{jik8} hold, and, in addition, 
\begin{equation}\label{jik16}
b\geq-3\delta'm\qquad\text{ and }\qquad k_1\geq -0.6m.
\end{equation}
If $m=L+1$ then \eqref{jik7} follows easily using \eqref{EstimJ}. On the other hand, if $m\leq L$ then we integrate by parts in time and use \eqref{NF2}--\eqref{NF3}. Using \eqref{yip6}, the left-hand side of \eqref{jik7} is bounded by
\begin{equation}\label{jik17.5}
\begin{split}
C2^{-2k_1}2^{-3b}\sup_{s\in J_m}\big\{&(\|P_{k_1}U^{\mathcal{L}_1h_1}(s)\|_{L^2}+2^m\|P_{k_1}\mathcal{L}_1\mathcal{N}^h_1(s)\|_{L^2})\|P_{k_2}U^{h_2}(s)\|_{L^\infty}\|P_kU^{\mathcal{L}h}(s)\|_{L^2}\\
&+2^m\|P_{k_1}U^{\mathcal{L}_1h_1}(s)\|_{L^2}\|P_{k_2}\mathcal{N}^h_2(s)\|_{L^\infty}\|P_kU^{\mathcal{L}h}(s)\|_{L^2}\\
&+2^m\|P_{k_1}U^{\mathcal{L}_1h_1}(s)\|_{L^2}\|P_{k_2}U^{h_2}(s)\|_{L^\infty}\|P_k\mathcal{L}\mathcal{N}^h(s)\|_{L^2}\big\},
\end{split}
\end{equation}
where $\mathcal{N}^h_1,\mathcal{N}^h_2,\mathcal{N}^h$ are suitable components of the metric nonlinearities $\mathcal{N}^h_{\al\be}$. In view of \eqref{vcx1}, \eqref{wws1}, \eqref{wer4.0}, and \eqref{plk2}, all the terms in the expression above are dominated by
\begin{equation*}
C\varep_1^32^{-3k_1/2}2^{-3b}2^{-m+2\delta' m}\lesssim \varep_1^32^{-m/20},
\end{equation*}
where we used the assumptions \eqref{jik16} in the last inequality. This suffices to prove \eqref{jik7}.

{\bf{Case 4.3.}} Assume now that the inequalities \eqref{jik8} hold, and, in addition, 
\begin{equation}\label{jik20}
k_1\leq -0.6m\qquad\text{ and }\qquad b\in[-3\delta'm,-2\delta m].
\end{equation}
We can use the condition $k_1\leq -0.6 m$ to improve the argument in {\bf{Case 4.1.}} Indeed, let $J'_2:=-k_1$ and decompose $P_{k_2}U^{h_2,\iota_2}(s)=U^{h_2,\iota_2}_{\leq J'_2,k_2}(s)+U^{h_2,\iota_2}_{> J'_2,k_2}(s)$ as in \eqref{on11.3}--\eqref{on11.36}. Using \eqref{bil1} (instead of \eqref{Bil31}), \eqref{nag1}, the assumption \eqref{doubnull}, \eqref{vcx1}, and \eqref{vcx1.2}, we estimate
\begin{equation}\label{jik21}
\begin{split}
&2^{-k_1}\Big|\int_{J_m}q_m(s)\mathcal{G}_{\mathfrak{q}_{\iota_1\iota_2}^b}\big[P_{k_1}U^{\mathcal{L}_1h_1,\iota_1}(s),U^{h_2,\iota_2}_{\leq J'_2,k_2}(s),P_kU^{\mathcal{L}h,\iota}(s)\big]\,ds\Big|\\
&\lesssim 2^{-k_1}2^m\sup_{s\in J_m}\big\|P_kI_{\mathfrak{q}_{\iota_1\iota_2}^b}\big[P_{k_1}U^{\mathcal{L}_1h_1,\iota_1}(s),U^{h_2,\iota_2}_{\leq J'_2,k_2}(s)\big]\big\|_{L^2}\big\|P_{[k-2,k+2]}U^{\mathcal{L}h,\iota}(s)\big\|_{L^2}\\
&\lesssim 2^{-k_1}2^m\cdot 2^b2^{k_1/2}2^{-m}2^{3k_2/2}\sup_{s\in J_m}\|P_{k_1}U^{\mathcal{L}_1h_1,\iota_1}(s)\|_{L^2}\|\widehat{P_{k_2}U^{h_2,\iota_2}}(s)\|_{L^\infty}\|P'_{k}U^{\mathcal{L}h,\iota}(s)\|_{L^2}\\
&\lesssim \varep_1^32^b2^{2H(q,n)\delta m}2^{\delta m}2^{-2N(n)k^+-2k^+},
\end{split}
\end{equation}
compare with \eqref{jik11}. The contribution of the profile $U^{h_2,\iota_2}_{> J'_2,k_2}(s)$ can be estimated as in \eqref{jik12}, and the bounds \eqref{jik7} follow if $b\leq-2\delta m$.

{\bf{Case 4.4.}} Finally assume that the inequalities \eqref{jik8} hold, and, in addition, 
\begin{equation}\label{jik30}
k_1\leq -0.6m\qquad\text{ and }\qquad b\in[-2\delta m,4].
\end{equation} 
If $\iota\neq\iota_2$ then we can still use normal forms and estimate as in \eqref{jik17.5}, but with a factor of $2^{-k_1}$ replaced by $2^{-k}$, due to the better lower bound $|\Phi(\xi,\eta)|\gtrsim 2^{k}$. The desired bounds follow. Also in the case $m=L+1$ there is no loss of $2^m$ and the desired bounds follow easily using \eqref{EstimJ}.

On the other hand, if $\iota=\iota_2$ and $m\leq L$ then we may assume that $\iota=\iota_2=+$, by taking complex conjugates. The proof in this case is more complicated, as it requires switching to the quasilinear variables $\mathcal{U}^{h_2}$ and $\mathcal{U}^{\LL h}$, and is provided in Lemma \ref{ParalinearSym2} below.
\end{proof}

\subsection{Non-null semilinear terms} In this subsection we show how to estimate the remaining Wave$\times$Wave$\times$Wave interactions in Proposition \ref{totbounds}.

\begin{lemma}\label{NonNullTerm1}
With the assumptions of Proposition \ref{totbounds} we have
\begin{equation}\label{GoodGGFTerms2}
\begin{split}
\sum_{k,k_1,k_2}2^{2\ga k^-}2^{-k}2^{2N(n)k^+}\Big\vert\int_{J_m}q_m(s)\mathcal{G}_{\mathfrak{m}}[P_{k_1}U^{\vartheta^{\mathcal{L}_1},\iota_1}(s),P_{k_2}U^{\vartheta^{\mathcal{L}_2},\iota_2}(s),P_kU^{\mathcal{L}h,\iota}(s)]\,ds\Big\vert\\
\lesssim\varepsilon_1^32^{2H(q,n)\delta m-2\ga m},
\end{split}
\end{equation}
where $\mathfrak{m}\in\mathcal{M}$ and the operators $\mathcal{G}_{\mathfrak{m}}$ are defined as in \eqref{DefI}. Thus the bounds \eqref{EEEst2} hold.
\end{lemma}

\begin{proof} In proving \eqref{GoodGGFTerms2} it is important to keep in mind that $\vartheta_{ab}$ are among the "good" components of the metric, which satisfy strong bounds like \eqref{vcx1.2} uniformly in time. 

For suitable values $J\geq \max\{-k_a,0\}$ we will sometimes decompose $P_{k_a}U^{\vartheta^{\LL_a},\iota_a}=U^{\vartheta^{\LL_a},\iota_a}_{\leq J,k_a}+U^{\vartheta^{\LL_a},\iota_a}_{>J,k_a}$, $a\in\{1,2\}$, where 
\begin{equation}\label{mik0}
U^{\vartheta^{\LL_a},\iota_a}_{\leq J,k_a}:=e^{-it\Lambda_{wa,\iota_a}}P'_{k_a}(\varphi_{\leq J}\cdot P_{k_a}V^{\vartheta^{\LL_a},\iota_a}),\qquad U^{\vartheta^{\LL_a},\iota_a}_{> J,k_a}:=e^{-it\Lambda_{wa,\iota_a}}P'_{k_a}(\varphi_{> J}\cdot P_{k_a}V^{\vartheta^{\LL_a},\iota_a}),
\end{equation}
and $V^{\vartheta^{\LL_a},\iota_a}:=e^{it\Lambda_{wa,\iota_a}}U^{\vartheta^{\LL_a},\iota_a}$. We consider several cases.

{\bf{Step 1.}} We assume first that $\min\{n_1,n_2\}=0$. By symmetry, we may assume that $n_1=0$. Let $J_1:=\infty$ if $|k_1|>\delta'm-10$ and $J_1:=m-\delta'm$ if $|k_1|\leq\delta'm-10$, and decompose $P_{k_1}U^{\vartheta,\iota_1}=U^{\vartheta,\iota_1}_{\leq J_1,k_1}+U^{\vartheta,\iota_1}_{>J_1,k_1}$ as in \eqref{mik0} (in particular, $P_{k_1}U^{\vartheta,\iota_1}=U^{\vartheta,\iota_1}_{\leq J_1,k_1}$ if $|k_1|>\delta'm-10$).

In this case we prove the stronger bounds
\begin{equation}\label{mik4}
\sum_{k,k_1,k_2}2^{2\ga k^-}2^{-k}2^{2N(n)k^+}\Big\vert\int_{J_m}q_m(s)\mathcal{G}_{\mathfrak{m}}[U^{\vartheta,\iota_1}_{>J_1,k_1},P_{k_2}U^{\vartheta^{\mathcal{L}},\iota_2}(s),P_kU^{\mathcal{L}h,\iota}(s)]\,ds\Big\vert\lesssim\varepsilon_1^3,
\end{equation}
and, for any $s\in J_m$ and $k\in\mathbb{Z}$, 
\begin{equation}\label{mik1}
\begin{split}
2^{\ga(m+k^-)}{\sum_{k_1,k_2}}^\ast\|P_kI[&U^{\vartheta,\iota_1}_{\leq J_1,k_1}(s),P_{k_2}U^{\vartheta^{\mathcal{L}},\iota_2}(s)]\|_{L^2}\\
&\lesssim \varep_1b_k(q,n;s)2^{k/2}2^{-N(n)k^+}2^{-m+H(q,n)\delta m},
\end{split}
\end{equation}
where $I=I_{\mathfrak{m}}$ is as in \eqref{abc36.1}, the coefficients $b_k(q,n;t)$ are defined in \eqref{nake15}, and ${\sum_{k_1,k_2}}^\ast$ denotes the sum over pairs $(k_1,k_2)\in\mathcal{X}_k$ with the additional assumption $k_1\leq k_2$ if $\LL=\mathrm{Id}$ (thus $n_1=n_2=n=0$). It is clear that \eqref{mik4}--\eqref{mik1} would suffice to prove \eqref{GoodGGFTerms2}. 

Using \eqref{nake15.2}, \eqref{vcx1.2}, and $L^2$ estimates as in \eqref{abc36.6}, we have, for $\ast\in\{\leq,>\}$,
\begin{equation}\label{mik2}
\begin{split}
2^{\ga(m+k^-)}&\|P_kI[U^{\vartheta,\iota_1}_{\ast,k_1}(s),P_{k_2}U^{\vartheta^{\mathcal{L}},\iota_2}(s)]\|_{L^2}\\
&\lesssim \varep_12^{3\underline{k}/2}b_{k_2}(q,n;s)2^{\ga(k^--k_2^-)}2^{k_2/2}2^{-N(n)k_2^+}2^{H(q,n)\delta m}\cdot 2^{k_1^-/2-\kappa k_1^-}2^{-N_0k_1^++2k_1^+},
\end{split}
\end{equation}
for any $k,k_1,k_2\in\mathbb{Z}$. This suffices to prove both bounds \eqref{mik4} and \eqref{mik1} when $2^m\lesssim 1$, or when $2^k\lesssim 2^{-m+\ga m/2}$.

{\bf{Substep 1.1.}} We prove now the bounds \eqref{mik1} when $m\geq \delta^{-1}$ and $k\geq -m+\ga m/2+10$.  The bounds \eqref{mik2} suffice to control the contribution of the pairs $(k_1,k_2)$ for which $k_1\leq -3m/4$ or $k_1\geq m/10$. 

To control the contribution of the pairs $(k_1,k_2)$ with $|k_1|>\delta'm-10$ (thus $U^{\vartheta,\iota_1}_{\leq J_1,k_1}=P_{k_1}U^{\vartheta,\iota_1}$) we use $L^2\times L^\infty$ estimates. More precisely, in view of \eqref{nake15.2} and \eqref{wws1} we have
\begin{equation}\label{mik7}
\begin{split}
&2^{\ga(m+k^-)}\|P_kI[P_{k_1}U^{\vartheta,\iota_1}(s),P_{k_2}U^{\vartheta^{\mathcal{L}},\iota_2}(s)]\|_{L^2}\\
&\lesssim \varep_1b_{k_2}(q,n;s)2^{\ga(k^--k_2^-)}2^{k_2/2}2^{-N(n)k_2^+}2^{H(q,n)\delta m}\cdot 2^{k_1^-}2^{-m+\delta'm/2}2^{-N(1)k_1^++2k_1^+},
\end{split}
\end{equation}
which suffices if $2^{|k_1|}\gtrsim 2^{\delta'm}$ and $k_1\leq k-4$ or if $|k-k_1|\leq 8$ and $n\geq 2$. The contribution of the pairs $(k_1,k_2)$ with $|k-k_1|\leq 8$ and $n\leq 1$ can be controlled in a similar way, by estimating $P_{k_1}U^{\vartheta,\iota_1}(s)$ in $L^2$ and $P_{k_2}U^{\vartheta^{\mathcal{L}},\iota_2}(s)$ in $L^\infty$. Finally, the contribution of the pairs $(k_1,k_2)\in \mathcal{X}_k$ with $k_1,k_2\geq k+4$ can be estimated using the $L^\infty$ super-localized estimates \eqref{Linfty1.6*},
\begin{equation}\label{mik7.5}
\begin{split}
2^{\ga(m+k^-)}&\|P_kI[P_{k_1}U^{\vartheta,\iota_1}(s),P_{k_2}U^{\vartheta^{\mathcal{L}},\iota_2}(s)]\|_{L^2}\\
&\lesssim \sum_{|n_1-n_2|\leq 4} 2^{\ga(m+k^-)}\|\mathcal{C}_{n_1,k}P_{k_1}U^{\vartheta}(s)\|_{L^\infty}\|\mathcal{C}_{n_2,k}P_{k_2}U^{\vartheta^{\mathcal{L}}}(s)\|_{L^2}\\
&\lesssim \varep_1^22^{k_2/2}2^{-N(n)k_2^+}2^{H(q,n)\delta m}\cdot 2^{k_1/2}2^{k/2}2^{-m+\delta'm/2}2^{-8k_1^+}.
\end{split}
\end{equation}
This suffices to bound the contribution of all the pairs $(k_1,k_2)$ with $|k_1|>\delta'm-10$ in \eqref{mik1}.

We consider now the pairs $(k_1,k_2)$ with $|k_1|\leq \delta'm-10$ and use the more precise bilinear  estimates \eqref{bil1}. Let $J'_1$ denote the largest integer such that 
\begin{equation}\label{mik7.6}
2^{J'_1}\leq 2^{-10}[2^{m/2-k_1/2}+2^{-k_2}+2^{-k}]
\end{equation}
and apply first \eqref{bil1} (or \eqref{Linfty1.1} if $k_1=\min\{k,k_1,k_2\}$), and then \eqref{vcx1.2} to estimate
\begin{equation}\label{SimpleLinftyThetaTheta}
\begin{split}
&2^{\ga(m+k^-)}\Vert P_kI[U^{\vartheta,\iota_1}_{\leq J'_1,k_1}(s),P_{k_2}U^{\vartheta^\mathcal{L},\iota_2}(s)]\Vert_{L^2}\\
&\lesssim 2^{-m}2^{\min\{k,k_1,k_2\}/2}2^{3k_1/2}\Vert \widehat{P_{k_1}V^{\vartheta}(s)}\Vert_{L^\infty}\cdot 2^{\ga(m+k^-)}\Vert P_{k_2}U^{\vartheta^\mathcal{L}}(s)\Vert_{L^2}\\
&\lesssim\varepsilon_12^{-m+H(q,n)\delta m}\cdot 2^{k_1^-/4}2^{-(N_0-2)k_1^+}2^{\min\{k,k_1,k_2\}/2}2^{k_2/2}2^{\ga(k^--k_2^-)}b_{k_2}(q,n;s)2^{-N(n_2)k_2^+},
\end{split}
\end{equation}
where $U^{\vartheta,\iota_1}_{\leq J'_1,k_1}:=e^{-it\Lambda_{wa,\iota_1}}P'_{k_1}(\varphi_{\leq J'_1}\cdot P_{k_1}V^{\vartheta,\iota_1})$ as in \eqref{mik0}. Moreover, using \eqref{Linfty1.6} and \eqref{vcx1.1} instead, we estimate
\begin{equation}\label{ThetaThetaFCase2}
\begin{split}
\sum_{j_1\in[J'_1,J_1]}2^{\ga(m+k^-)}\Vert &P_kI[U^{\vartheta,\iota_1}_{j_1,k_1}(s),P_{k_2}U^{\vartheta^\mathcal{L},\iota_2}(s)]\Vert_{L^2}\\
&\lesssim \varepsilon_1^22^{-m+\delta'm}2^{-J'_1}2^{k_2/2}2^{-8k_1^+}2^{-N(n_2)k_2^+}.
\end{split}
\end{equation}
It is easy to see that these two bounds can be summed over $k_1,k_2\in\mathcal{X}_k$ with $|k_1|\leq \delta'm-10$ to complete the proof of \eqref{mik1}. In the (harder) case of pairs $k_1,k_2$ with $k_1,k_2\geq k+4$ one can use \eqref{nake15.1} to sum \eqref{SimpleLinftyThetaTheta}, and the definition \eqref{mik7.6} to sum \eqref{ThetaThetaFCase2}.

{\bf{Substep 1.2.}} We prove now the bounds \eqref{mik4} when $m\geq \delta^{-1}$ and $k\geq -m+\ga m/2+10$. We may assume that $|k_1|\leq \delta'm$ (due to the definition of $J_1$). We notice that we can still prove strong bounds similar to \eqref{mik1} if $k\leq -2\delta'm$, using just $L^2$ estimates and \eqref{vcx1}--\eqref{vcx1.1}.

Assume now that $k\geq -2\delta'm$. We would like to integrate by parts in time as in \eqref{NF2}. For this we need to decompose into resonant and non-resonant contributions. As in \eqref{trf2}, with $q_0=-8\delta^\prime m$ and $\Xi_{\iota_1\iota_2}$ as in \eqref{par1},  we decompose
\begin{equation}\label{trf2.1}
\begin{split}
\mathfrak{m}&=\mathfrak{m}^r+\mathfrak{m}^{nr},\qquad\mathfrak{m}^r(\theta,\eta)=\varphi_{\le q_0}(\Xi_{\iota_1\iota_2}(\theta,\eta))\mathfrak{m}(\theta,\eta).
\end{split}
\end{equation}

To bound the resonant contribution we use the smallness of Fourier support and Schur's test. For this we notice that, with $\varphi_{kk_1k_2}$ defined as in \eqref{phikkk},
\begin{equation*}
\begin{split}
&\sup_{\eta\in\mathbb{R}^3}\int_{\mathbb{R}^3}|\widehat{U^{\vartheta,\iota_1}_{>J_1,k_1}}(\xi-\eta,s)||\mathfrak{m}^r(\xi-\eta,\eta)|\varphi_{kk_1k_2}(\xi-\eta,\eta)\,d\xi\lesssim \big\|\widehat{U^{\vartheta,\iota_1}_{>J_1,k_1}}(s)\big\|_{L^2}\cdot 2^{3k_1/2}2^{q_0},\\
&\sup_{\xi\in\mathbb{R}^3}\int_{\mathbb{R}^3}|\widehat{U^{\vartheta,\iota_1}_{>J_1,k_1}}(\xi-\eta,s)||\mathfrak{m}^r(\xi-\eta,\eta)|\varphi_{kk_1k_2}(\xi-\eta,\eta)\,d\eta\lesssim \big\|\widehat{U^{\vartheta,\iota_1}_{>J_1,k_1}}(s)\big\|_{L^2}\cdot 2^{3k_1/2}2^{q_0}2^{k_2-k},
\end{split}
\end{equation*}
where in the second line we used that $\widetilde{\Xi}(\xi-\eta,\xi)\lesssim \widetilde{\Xi}(\xi-\eta,\eta)2^{k_2-k}\lesssim 2^{q_0}2^{k_2-k}$ in the support of the integral (see \eqref{par73}). Therefore, according to Schur's test,
\begin{equation*}
\begin{split}
2^{\ga(m+k^-)}\|P_kI_{\mathfrak{m}^r}&[U^{\vartheta,\iota_1}_{>J_1,k_1}(s),P_{k_2}U^{\vartheta^{\mathcal{L}},\iota_2}(s)]\|_{L^2}\\
&\lesssim 2^{\ga(m+k^-)}\|P_{k_2}U^{\vartheta^{\mathcal{L}}}(s)\|_{L^2}\big\|U^{\vartheta,\iota_1}_{>J_1,k_1}(s)\big\|_{L^2}\cdot 2^{3k_1/2}2^{q_0}(1+2^{k_2-k}),\\
&\lesssim \varep_12^{k_1}2^{-N(1)k_1^+}2^{-J_1}2^{q_0}2^{\delta'm}\cdot 2^{k_2/2}2^{-N(n)k_2^+}b_{k_2}(q,n;s)(1+2^{k_2-k})2^{\ga(k^--k_2^-)},
\end{split}
\end{equation*}
where $I_{\mathfrak{m}}$ is as in \eqref{abc36.1}. Since $|k_1|\leq \delta'm$ and $2^{-J_1}\lesssim 2^{-m+\delta'm}$, this suffices to show that
\begin{equation}\label{mik9}
\begin{split}
2^{\ga(m+k^-)}{\sum_{k_1,k_2}}^\ast\|P_kI_{\mathfrak{m}^r}[&U^{\vartheta,\iota_1}_{>J_1,k_1}(s),P_{k_2}U^{\vartheta^{\mathcal{L}},\iota_2}(s)]\|_{L^2}\\
&\lesssim \varep_1b_k(q,n;s)2^{k/2}2^{-N(n)k^+}2^{-m-\delta'm},
\end{split}
\end{equation}
if $2^k\gtrsim 2^{-2\delta'm}$. This is similar to \eqref{mik1} and implies the required bounds on the resonant contributions.

On the other hand, the non-resonant contributions corresponding to the symbol $\mathfrak{m}^{nr}$ can be treated as in the proof of \eqref{NullLHHWWWF2} in Lemma \ref{NullLHHWWWS}. We may assume first that $m\leq L$, integrate by parts in time, and notice that estimates like \eqref{bnm16} still hold; such estimates suffice to prove  \eqref{mik4} when $k\leq m/10$. In the remaining case when $k$ is very large and $|k_1|\leq\delta'm$, the desired conclusion follows using paradifferential calculus from Lemma \ref{ParalinearSym} below.

{\bf{Step 2.}} We prove now \eqref{GoodGGFTerms2} when $n_1,n_2\geq 1$. In this case we prove the strong bounds
\begin{equation}\label{mik11}
\begin{split}
2^{\ga(m+k^-)}\sum_{k_1,k_2\in\mathcal{X}_k}\|P_kI[&P_{k_1}U^{\vartheta^{\mathcal{L}_1},\iota_1}(s),P_{k_2}U^{\vartheta^{\mathcal{L}_2},\iota_2}(s)]\|_{L^2}\\
&\lesssim \varep_1^22^{-\ga|k|/4}2^{k/2}2^{-N(n)k^+}2^{-m+H(q,n)\delta m},
\end{split}
\end{equation}
for any $s\in J_m$ and $k\in\mathbb{Z}$ (compare with \eqref{mik1} and recall that $b_k(q,n;s)\geq\varep_12^{-\ga|k|/4}$). 

As in \eqref{mik2}, we recall that $\underline{k}=\min\{k,k_1,k_2\}$ and start with $L^2$ estimates,
\begin{equation}\label{mik12}
\begin{split}
2^{\ga(m+k^-)}&\|P_kI[P_{k_1}U^{\vartheta^{\mathcal{L}_1},\iota_1}(s),P_{k_2}U^{\vartheta^{\mathcal{L}_2},\iota_2}(s)]\|_{L^2}\\
&\lesssim \varep_1^22^{3\underline{k}/2}2^{k_1/2}2^{k_2/2}2^{-\ga k_1^--\ga k_2^-}2^{-N(n_1)k_1^+}2^{-N(n_2)k_2^+}2^{H(q_1,n_1)\delta m+H(q_2,n_2)\delta m},
\end{split}
\end{equation}
for any $k,k_1,k_2\in\mathbb{Z}$. This suffices to prove both bounds \eqref{mik11} when $2^m\lesssim 1$, or when $k\leq -m+35\delta m$ or when $k\geq m/9$ (using \eqref{SuperlinearH1} and $N(n)+10\leq\min \{N(n_1),N(n_2)\}$). 

On the other hand, if $m\geq \delta^{-1}$ and $k\in[-m+35\delta m,m/9]$ then the bounds \eqref{mik12} suffice to control the contribution of the pairs $(k_1,k_2)\in\mathcal{X}_k$ for which either $\min\{k_1,k_2\}\leq -3m/4$ or $\max\{k_1,k_2\}\geq m/9$. Using $L^2\times L^\infty$ estimates as in \eqref{mik7}--\eqref{mik7.5} we can also control the contribution of the pairs $(k_1,k_2)$ for which $\max\{|k_1|,|k_2|\}\geq\delta'm$. In the remaining range $k_1,k_2\in[-\delta'm,\delta'm]$ we consider two cases.

{\bf{Case 2.1.}} Assume first that $q_1=q_2=0$. Since rotations and Riesz transforms essentially commute (up to multipliers that are accounted for in the multiplier $\mathfrak{m}$), we may replace $U^{\vartheta^{\mathcal{L}_a},\iota_a}$ by $\mathcal{L}_a U^{\vartheta,\iota_a}$, $a\in\{1,2\}$, and use interpolation inequalities. Let 
\begin{equation}\label{mik18.2}
J^\ast=\max\big\{\frac{m-\max\{k_1,k_2\}}{2},-k\big\}.
\end{equation}

Assuming that $k\ge\max\{k_1,k_2\}-10$, we use \eqref{inte4}--\eqref{inte5} to estimate
\begin{equation*}
\begin{split}
\|P_kI[U^{\vartheta^{\mathcal{L}_1},\iota_1}_{\leq J^\ast,k_1}(s),&U^{\vartheta^{\mathcal{L}_2},\iota_2}_{\leq J^\ast,k_2}(s)]\|_{L^2}\lesssim\Vert \LL_1U^{\vartheta,\iota_1}_{\leq J^\ast,k_1}(s)\Vert_{L^a}\Vert \LL_2U^{\vartheta,\iota_2}_{\leq J^\ast,k_2}(s)\Vert_{L^b}\\
&\lesssim \Vert U^{\vartheta,\iota_1}_{\leq J^\ast,k_1}(s)\Vert_{L^\infty}^{\frac{n-n_1}{n}}\Vert P_{k_1}V^{\vartheta,\iota_1}\Vert_{H^{0,n}_\Omega}^{\frac{n_1}{n}}\Vert U^{\vartheta,\iota_2}_{\leq J^\ast,k_2}(s) \Vert_{L^\infty}^{\frac{n-n_2}{n}}\Vert P_{k_2}V^{\vartheta,\iota_2}\Vert_{H^{0,n}_\Omega}^\frac{n_2}{n},
\end{split}
\end{equation*}
where $U^{\vartheta^{\mathcal{L}_a},\iota_a}_{\leq J^\ast,k_a}$ are as in \eqref{mik0}, $n=n_1+n_2$, and
\begin{equation}\label{mik18.1}
\begin{split}
\frac{1}{a}=\frac{n_1}{n}\frac{1}{2},\qquad\frac{1}{b}=\frac{n_2}{n}\frac{1}{2}.
\end{split}
\end{equation}
Using now \eqref{wws12x} and \eqref{vcx1} we have
\begin{equation}\label{mik18}
\begin{split}
2^{\ga(m+k^-)}\|P_kI[&U^{\vartheta^{\mathcal{L}_1},\iota_1}_{\leq J^\ast,k_1}(s),U^{\vartheta^{\mathcal{L}_2},\iota_2}_{\leq J^\ast,k_2}(s)]\|_{L^2}\\
&\lesssim  \varep_1^22^{-m+H(0,n)\delta m}2^{k_1/2+k_2/2}2^{-N(n)k_1^+-N(n)k_2^+}2^{-|k_1|/8-|k_2|/8}.
\end{split}
\end{equation}

On the other hand, if $k\le\max\{k_1,k_2\}-10$ then we need a bilinear estimate to bring in the small factor $2^{k/2}$. We use Lemma \ref{LemBil5} followed by \eqref{inte4}--\eqref{inte5} to estimate
\begin{equation*}
\begin{split}
\|P_kI[&U^{\vartheta^{\mathcal{L}_1},\iota_1}_{\leq J^\ast,k_1}(s),U^{\vartheta^{\mathcal{L}_2},\iota_2}_{\leq J^\ast,k_2}(s)]\|_{L^2}\\
&\lesssim 2^{-m}2^{k/2}2^{3k_1/2}\Vert \mathcal{F}\{\LL_1P_{k_1}U^{\vartheta,\iota_1}\}(s)\Vert_{L^a}\Vert \mathcal{F}\{\LL_2P_{k_2}U^{\vartheta,\iota_2}\}(s)\Vert_{L^b}\\
&\lesssim 2^{-m}2^{k/2}2^{3k_1/2}\Vert \mathcal{F}\{P_{k_1}V^{\vartheta}\}\Vert_{L^\infty}^{\frac{n-n_1}{n}}\Vert P_{k_1}V^{\vartheta}\Vert_{H^{0,n}_\Omega}^\frac{n_1}{n}\Vert\mathcal{F}\{P_{k_2}V^{\vartheta}\}\Vert_{L^\infty}^{\frac{n-n_2}{n}}\Vert P_{k_2}V^{\vartheta}\Vert_{H^{0,n}_\Omega}^\frac{n_2}{n},
\end{split}
\end{equation*}
where $a,b$ are as in \eqref{mik18.1}. Using \eqref{vcx1.2} and \eqref{vcx1}, it follows that
\begin{equation}\label{mik18.5}
2^{\ga(m+k^-)}\|P_kI[U^{\vartheta^{\mathcal{L}_1},\iota_1}_{\leq J^\ast,k_1}(s),U^{\vartheta^{\mathcal{L}_2},\iota_2}_{\leq J^\ast,k_2}(s)]\|_{L^2}\lesssim  \varep_1^22^{-m+H(0,n)\delta m}2^{\ga k^-}2^{k/2}2^{-N_0k_1^+}2^{3k_1^-/4}.
\end{equation}

For the remaining contributions we just use \eqref{wws1} and \eqref{vcx1.1} and $L^2\times L^\infty$ estimates,
\begin{equation}\label{mik18.7}
\begin{split}
2^{\ga(m+k^-)}\big\{\|P_kI[&U^{\vartheta^{\mathcal{L}_1},\iota_1}_{\leq J^\ast,k_1}(s),U^{\vartheta^{\mathcal{L}_2},\iota_2}_{>J^\ast,k_2}(s)]\|_{L^2}+\|P_kI[U^{\vartheta^{\mathcal{L}_1},\iota_1}_{>J^\ast,k_1}(s),P_{k_2}U^{\vartheta^{\mathcal{L}_2},\iota_2}(s)]\|_{L^2}\big\}\\
\lesssim &\varepsilon_1^22^{-m+\delta'm}2^{-J^\ast}2^{|k_1|/2+|k_2|/2}2^{-N(n_1+1)k_1^+}2^{-N(n_2+1)k_2^+}2^{2(k_1^++k_2^+)}.
\end{split}
\end{equation}
We combine \eqref{mik18} and \eqref{mik18.7} to control the contribution of the pairs $(k_1,k_2)$ for which $k\geq \max\{k_1,k_2\}-10$; we also combine \eqref{mik18.5} and \eqref{mik18.7} to control the contribution of the pairs $(k_1,k_2)$ for which $k\leq \max\{k_1,k_2\}-10$. This completes the proof of \eqref{mik11}.

{\bf{Case 2.2.}} The case $\max\{q_1,q_2\}\ge 1$ is comparatively easier. Without loss of generality we may assume that $q_2\geq\max\{q_1,1\}$. We use Lemma \ref{Lembil2} and \eqref{vcx1} to estimate
\begin{equation*}
\begin{split}
2^{\ga(m+k^-)}&\|P_kI[U^{\vartheta^{\mathcal{L}_1},\iota_1}_{\leq m(1-\delta),k_1}(s),P_{k_2}U^{\vartheta^{\mathcal{L}_2},\iota_2}(s)]\|_{L^2}\\
&\lesssim 2^{\ga m}2^{-m+\delta m}2^{\min\{k,k_2\}/2}\Vert P_{k_1}U^{\vartheta^{\mathcal{L}_1}}(s)\Vert_{H^{0,1}_\Omega}\Vert P_{k_2}U^{\vartheta^{\LL_2}}(s)\Vert_{L^2}\\
&\lesssim 2^{-m+\delta m}2^{\min\{k,k_2\}/2}2^{k_1/2}2^{k_2/2}2^{-N(n_1+1)k_1^+}2^{-N(n_2)k_2^+}2^{H(q_1,n_1+1)\delta m+H(q_2,n_2)\delta m}.
\end{split}
\end{equation*}
This gives acceptable contributions using \eqref{SuperlinearH2}. On the other hand, using just $L^2\times L^\infty$ estimates, \eqref{vcx1.1}, and \eqref{wws1}, we have
\begin{equation*}
\begin{split}
2^{\ga(m+k^-)}&\|P_kI[U^{\vartheta^{\mathcal{L}_1},\iota_1}_{>m(1-\delta),k_1}(s),P_{k_2}U^{\vartheta^{\mathcal{L}_2},\iota_2}(s)]\|_{L^2}\lesssim \varep_1^22^{-1.9m},
\end{split}
\end{equation*}
for any $k_1,k_2\in [-\delta'm,\delta'm]$. The bounds \eqref{mik11} follow in this case as well, which completes the proof of the lemma.
\end{proof}

\subsection{Second symmetrization and paradifferential calculus}

In the case of interactions of vastly different frequencies $\vert\xi\vert\approx\vert\eta\vert\gg\vert\xi-\eta\vert$, $\iota=\iota_2$, the application of normal forms as in \eqref{NF2} leads to a loss of derivatives due to the quasilinear nature of the nonlinearities. To avoid this loss this we use paradifferential calculus and perform a second symmetrization.

Our first lemma concerns the contribution of very high frequencies and applies to conclude the analysis in {\bf{Case 2.3.}} in Lemma \ref{NullLHHWWWS} and {\bf{Substep 1.2.}} in Lemma \ref{NonNullTerm1}.

\begin{lemma}\label{ParalinearSym} If $\iota,\iota_1,\iota_2\in\{+,-\}$, $h,h_1,h_2\in\{h_{\al\be}\}$, $t\in[0,T]$, $m\in[\delta^{-1},L+1]$, $\vert k_1\vert\le \delta^\prime m$, $k,k_2\ge m/10$, $J_1\geq m-2\delta'm$, $(q,n)\leq (3,3)$, and $\mathcal{L},\LL_2\in\mathcal{V}^q_n$ then
\begin{equation}\label{bnm25}
\begin{split}
\Big|\int_{J_m}q_m(s)\mathcal{G}_{\mathfrak{m}^{nr}}&\big[U^{h_1,\iota_1}_{\leq J_1,k_1}(s),P_{k_2}U^{\mathcal{L}_2h_2,\iota_2}(s),P_kU^{\mathcal{L}h,\iota}(s)\big]\,ds\Big|\lesssim \varepsilon_12^{-2N(n)k+k}2^{-2\delta m} (b_{k,m}(q,n))^2,
\end{split}
\end{equation}
where $\mathfrak{m}^{nr}(\theta,\eta)=\mathfrak{m}(\theta,\eta)\varphi_{>q_0}(\Xi_{\iota_1\iota_2}(\theta,\eta))$ is the non-resonant part of a symbol $\mathfrak{m}\in\mathcal{M}$, $q_0=-8\delta'm$, and $b_{k,m}(q,n)$ are defined in \eqref{bnm21}.
\end{lemma}

\begin{proof} Notice that, as a consequence of \eqref{nake15.2}-\eqref{nake15.5},
\begin{equation}\label{sori8}
\|P_{k_2}U^{\mathcal{L}_2h_2,\iota_2}(s)\|_{L^2}+\|P_kU^{\mathcal{L}h,\iota}(s)\|_{L^2}\lesssim b_k(s)2^{\delta'm}2^{k/2-N(n)k}
\end{equation}
and
\begin{equation}\label{sori9}
\begin{split}
\|(\partial_s+i\Lambda_{wa,\iota_2})&P_{k_2}U^{\mathcal{L}_2h_2,\iota_2}(s)\|_{L^2}\\
&+\|(\partial_s+i\Lambda_{wa,\iota})P_kU^{\mathcal{L}h,\iota}(s)\|_{L^2}\lesssim \varep_1b_k(s)2^{-m+\delta'm}2^{3k/2-N(n)k},
\end{split}
\end{equation}
for any $s\in J_m$, where, for simplicity of notation, we let $b_k(s):=b_k(q,n;s)$. Moreover,
\begin{equation}\label{sori10}
\|U^{h_1,\iota_1}_{\leq J_1,k_1}(s)\|_{L^\infty}\lesssim \varep_12^{-m+\delta'm}2^{k_1^-}2^{-8k_1^+}
\end{equation}
and
\begin{equation}\label{sori11}
\|(\partial_s+i\Lambda_{wa,\iota_1})U^{h_1,\iota_1}_{\leq J_1,k_1}(s)\|_{L^\infty}\lesssim \varep_1^22^{-1.9m}.
\end{equation}
The bounds \eqref{sori10} follow from \eqref{wws1}. The bounds \eqref{sori11} are slightly harder because of the spacial cutoffs. To prove them we write, for any $\al,\be\{0,1,2,3\}$,
\begin{equation*}
\begin{split}
(\partial_s+i\Lambda_{wa})U^{h_{\al\be}}_{\leq J_1,k_1}(s)&=(\partial_s+i\Lambda_{wa})[e^{-is\Lambda_{wa}}P'_{k_1}(\varphi_{\leq J_1}\cdot P_{k_1}V^{h_{\al\be}}(s))]\\
&=e^{-is\Lambda_{wa}}P'_{k_1}(\varphi_{\leq J_1}\cdot \partial_sP_{k_1}V^{h_{\al\be}}(s))\\
&=e^{-is\Lambda_{wa}}P'_{k_1}(\varphi_{\leq J_1}\cdot P_{k_1}e^{is\Lambda_{wa}}\mathcal{N}^h_{\al\be}(s)).
\end{split}
\end{equation*}
Therefore, for any $x\in\mathbb{R}^3$,
\begin{equation*}
\begin{split}
(\partial_s+i\Lambda_{wa})&U^{h_{\al\be}}_{\leq J_1,k_1}(x,s)=C\int_{\mathbb{R}^6}e^{ix\cdot\xi}e^{-is|\xi|}\varphi'_{k_1}(\xi)\widehat{\varphi_{\leq J_1}}(\xi-\eta)e^{is|\eta|}\varphi_{k_1}(\eta)\widehat{\mathcal{N}^h_{\al\be}}(\eta,s)\,d\xi d\eta\\
&=C\int_{\mathbb{R}^6}e^{ix\cdot\eta}e^{ix\cdot\theta}e^{-is[|\eta+\theta|-|\eta|]}\varphi'_{k_1}(\eta+\theta)\widehat{\varphi_{\leq J_1}}(\theta)\varphi_{k_1}(\eta)\widehat{\mathcal{N}^h_{\al\be}}(\eta,s)\,d\theta d\eta,
\end{split}
\end{equation*}
where $\varphi'_{k_1}=\varphi_{[k_1-2,k_1+2]}$. Let
\begin{equation}\label{sori12}
L_{x,s}(\eta):=\int_{\mathbb{R}^3}e^{ix\cdot\theta}e^{-is[|\eta+\theta|-|\eta|]}\varphi'_{k_1}(\eta+\theta)(2^{3J_1}\widehat{\varphi}(2^{J_1}\theta))\,d\theta,
\end{equation}
so
\begin{equation}\label{sori13}
(\partial_s+i\Lambda_{wa})U^{h_{\al\be}}_{\leq J_1,k_1}(x,s)=C\int_{\mathbb{R}^6}e^{i(x-y)\cdot\eta}L_{x,s}(\eta)\varphi_{k_1}(\eta)\mathcal{N}^h_{\al\be}(y,s)\,d\eta dy.
\end{equation}
Since $s2^{-J_1}+2^{|k_1|}\lesssim 2^{2\delta'm}$ (the hypothesis of the lemma), it is easy to see that $|D^\alpha_{\eta}L_{x,s}(\eta)|\lesssim 2^{4\delta'm|\alpha|}$ for any $x\in\mathbb{R}^3$, $|\eta|\approx 2^{k_1}$, and multi-indices $\alpha$ with $|\alpha|\leq 10$. Therefore
\begin{equation*}
|(\partial_s+i\Lambda_{wa})U^{h_{\al\be}}_{\leq J_1,k_1}(x,s)|\lesssim \int_{\mathbb{R}^3}2^{4\delta'm}(1+|x-y|2^{-4\delta'm})^{-4}|\mathcal{N}^h_{\al\be}(y,s)|\,dy,
\end{equation*}
using integration by parts in $\eta$ in \eqref{sori13}. The desired conclusion \eqref{sori11} follows from \eqref{sho4.2}.

We divide the rest of the proof of the lemma into several steps.

{\bf{Step 1.}} We start with some preliminary reductions. First, we may assume that $m\leq L$ since in the case $m=L+1$ we have $|J_m|\lesssim 1$ and the desired bounds follow using \eqref{sori8}, \eqref{sori10}, and \eqref{yip2}. We may also assume that $\iota=\iota_2$; otherwise
\begin{equation*}
|\Phi_{\sigma\mu\nu}(\xi,\eta)|=|\Lambda_{wa,\iota}(\xi)-\Lambda_{wa,\iota_2}(\eta)-\Lambda_{wa,\iota_1}(\xi-\eta)|\gtrsim 2^k
\end{equation*}
in the support of the integral, so the normal form argument \eqref{NF2} still gives the desired conclusion since the loss of derivatives in \eqref{sori9} is compensated by the large denominator. By taking complex conjugates, we may actually assume that $\iota_2=\iota=+$. 

To continue we switch to the quasilinear variables defined in \eqref{MathcalU},
\begin{equation*}
\mathcal{U}^{\LL_2h_2}=(\partial_t-iT_{\sigma_{wa}})(\LL_2 h_{2})\quad\text{ and }\quad\mathcal{U}^{\LL h}=(\partial_t-iT_{\sigma_{wa}})(\LL h).
\end{equation*}
See subsection \ref{ParaCalc} for the definition of the paradifferential operators $T_a$ and the symbols $\sigma_{wa}$. In view of \eqref{DefParaUnknown} we have, for any $s\in J_m$,
\begin{equation}\label{sori19} 
\|P_{k_2}(\mathcal{U}^{\LL_2h_2}-U^{\LL_2h_2})(s)\|_{L^2}+\|P_{k}(\mathcal{U}^{\LL h}-U^{\LL h})(s)\|_{L^2}\lesssim 2^{-m+\delta'm}b_k(s)2^{k/2-N(n)k}.
\end{equation}
Thus we may replace $P_{k_2}U^{\LL_2h_2}(s)$ with $P_{k_2}\mathcal{U}^{\LL_2h_2}(s)$ and $P_{k}U^{\LL h}(s)$ with $P_{k}\mathcal{U}^{\LL h}(s)$ in the integral in \eqref{bnm25}, at the expense of acceptable errors. To summarize, it remains to prove that
\begin{equation}\label{sori20}
\Big|\int_{J_m}q_m(s)\mathcal{G}_{\mathfrak{m}^{nr}}\big[U^{h_1,\iota_1}_{\leq J_1,k_1}(s),P_{k_2}\mathcal{U}^{\mathcal{L}_2h_2}(s),P_k\mathcal{U}^{\mathcal{L}h}(s)\big]\,ds\Big|\lesssim \varepsilon_12^{-2N(n)k+k}2^{-2\delta m} b_{k,m}^2,
\end{equation}
provided that $m\leq L$ and, for simplicity, $b_{k,m}:=b_{k,m}(q,n)$. 

{\bf{Step 2.}} We integrate by parts in $s$ using \eqref{NF2}. The contributions of the first two terms, when the $d/ds$ derivative hits either the function $q_m(s)$ or the first term $U^{h_1,\iota_1}_{\leq J_1,k_1}(s)$ can be bounded easily, using the $L^\infty$ bounds \eqref{sori10}--\eqref{sori11} and the $L^2$ bounds \eqref{sori8} and \eqref{sori19} (there are no derivative losses in this case). So it remains to prove that
\begin{equation}\label{sori21}
\begin{split}
&\Big|\mathcal{H}_{\mathfrak{m}^{nr}}\big[U^{h_1,\iota_1}_{\leq J_1,k_1}(s),(\partial_s+i\Lambda_{wa})P_{k_2}\mathcal{U}^{\mathcal{L}_2h_2}(s),P_k\mathcal{U}^{\mathcal{L}h}(s)\big]\\
&+\mathcal{H}_{\mathfrak{m}^{nr}}\big[U^{h_1,\iota_1}_{\leq J_1,k_1}(s),P_{k_2}\mathcal{U}^{\mathcal{L}_2h_2}(s),(\partial_s+i\Lambda_{wa})P_k\mathcal{U}^{\mathcal{L}h}(s)\big]\Big|\lesssim \varepsilon_12^{-2N(n)k+k}2^{-m-2\delta m} b_{k}(s)^2,
\end{split}
\end{equation}
for any $s\in J_m$, where the operators $\mathcal{H}_{\mathfrak{m}^{nr}}$ are defined in \eqref{DefJ}. Notice that $$(\partial_s+i\Lambda_{wa})P_k\mathcal{U}^{\mathcal{L}h}=(\partial_s+iT_{\Sigma_{wa}})P_k\mathcal{U}^{\mathcal{L}h}-iT_{\Sigma_{wa}-|\ze|}P_k\mathcal{U}^{\mathcal{L}h},$$ 
where the symbols $\Sigma_{wa}$ are defined in \eqref{ParaSymbolsW}, and a similar identity holds for $P_{k_2}\mathcal{U}^{\mathcal{L}_2h_2}$. Therefore we use the bounds in \eqref{ParaBox5} to replace $(\partial_s+i\Lambda_{wa})P_{k_2}\mathcal{U}^{\mathcal{L}_2h_2}$ and $(\partial_s+i\Lambda_{wa})P_{k}\mathcal{U}^{\mathcal{L}h}$ with $-iT_{\Sigma_{wa}-|\ze|}P_{k_2}\mathcal{U}^{\mathcal{L}_2h_2}$ and $-iT_{\Sigma_{wa}-|\ze|}P_k\mathcal{U}^{\mathcal{L}h}$ respectively in \eqref{sori21}, at the expense of acceptable errors. For \eqref{sori21} it remains to prove that
\begin{equation}\label{sori24}
\begin{split}
&\Big|\mathcal{H}_{\mathfrak{m}^{nr}}\big[U^{h_1,\iota_1}_{\leq J_1,k_1}(s),iT_{\Sigma^{\geq 1}_{wa}}P_{k_2}\mathcal{U}^{\mathcal{L}_2h_2}(s),P_k\mathcal{U}^{\mathcal{L}h}(s)\big]\\
&+\mathcal{H}_{\mathfrak{m}^{nr}}\big[U^{h_1,\iota_1}_{\leq J_1,k_1}(s),P_{k_2}\mathcal{U}^{\mathcal{L}_2h_2}(s),iT_{\Sigma^{\geq 1}_{wa}}P_k\mathcal{U}^{\mathcal{L}h}(s)\big]\Big|\lesssim \varepsilon_12^{-2N(n)k+k}2^{-m-2\delta m} b_{k}(s)^2,
\end{split}
\end{equation}
for any $s\in J_m$, where $\Sigma^{\geq 1}_{wa}(x,\ze):=\Sigma_{wa}(x,\ze)-|\ze|$.

{\bf{Step 3.}} We now write explicitly the expression in the left-hand side of \eqref{sori24} and exploit the cancellation between the two terms to avoid derivative loss. Using the definitions we write
\begin{equation*}
\begin{split}
\mathcal{H}_{\mathfrak{m}^{nr}}\big[U^{h_1,\iota_1}_{\leq J_1,k_1},&iT_{\Sigma^{\geq 1}_{wa}}P_{k_2}\mathcal{U}^{\mathcal{L}_2h_2},P_k\mathcal{U}^{\mathcal{L}h}\big]=\frac{i}{8\pi^3}\int_{\mathbb{R}^9}\frac{\mathfrak{m}^{nr}(\xi-\eta,\eta)}{|\xi|-|\eta|-\iota_1|\xi-\eta|}\widetilde{\Sigma^{\geq 1}_{wa}}\Big(\eta-\rho,\frac{\eta+\rho}{2}\Big)\\
&\times\chi_0\Big(\frac{|\eta-\rho|}{|\eta+\rho|}\Big)\widehat{U^{h_1,\iota_1}_{\leq J_1,k_1}}(\xi-\eta)\widehat{P_{k_2}\mathcal{U}^{\mathcal{L}_2h_2}}(\rho)\overline{\widehat{P_k\mathcal{U}^{\mathcal{L}h}}(\xi)}\, d\xi d\eta d\rho,
\end{split}
\end{equation*}
\begin{equation*}
\begin{split}
\mathcal{H}_{\mathfrak{m}^{nr}}\big[U^{h_1,\iota_1}_{\leq J_1,k_1},&P_{k_2}\mathcal{U}^{\mathcal{L}_2h_2},iT_{\Sigma^{\geq 1}_{wa}}P_k\mathcal{U}^{\mathcal{L}h}\big]=\frac{-i}{8\pi^3}\int_{\mathbb{R}^9}\frac{\mathfrak{m}^{nr}(\xi-\eta,\eta)}{|\xi|-|\eta|-\iota_1|\xi-\eta|}\overline{\widetilde{\Sigma^{\geq 1}_{wa}}}\Big(\xi-\rho,\frac{\xi+\rho}{2}\Big)\\
&\times\chi_0\Big(\frac{|\xi-\rho|}{|\xi+\rho|}\Big)\widehat{U^{h_1,\iota_1}_{\leq J_1,k_1}}(\xi-\eta)\widehat{P_{k_2}\mathcal{U}^{\mathcal{L}_2h_2}}(\eta)\overline{\widehat{P_k\mathcal{U}^{\mathcal{L}h}}(\rho)}\, d\xi d\eta d\rho.
\end{split}
\end{equation*}
The key property that allows symmetrization is the reality of the symbol $\Sigma^{\geq 1}_{wa}$, which shows that $\overline{\widetilde{\Sigma^{\geq 1}_{wa}}}\Big(\xi-\rho,\frac{\xi+\rho}{2}\Big)=\widetilde{\Sigma^{\geq 1}_{wa}}\Big(\rho-\xi,\frac{\xi+\rho}{2}\Big)$. Therefore, after changes of variables, we have
\begin{equation}\label{sori28}
\begin{split}
&\mathcal{H}_{\mathfrak{m}^{nr}}\big[U^{h_1,\iota_1}_{\leq J_1,k_1}(s),iT_{\Sigma^{\geq 1}_{wa}}P_{k_2}\mathcal{U}^{\mathcal{L}_2h_2}(s),P_k\mathcal{U}^{\mathcal{L}h}(s)\big]\\
&+\mathcal{H}_{\mathfrak{m}^{nr}}\big[U^{h_1,\iota_1}_{\leq J_1,k_1}(s),P_{k_2}\mathcal{U}^{\mathcal{L}_2h_2}(s),iT_{\Sigma^{\geq 1}_{wa}}P_k\mathcal{U}^{\mathcal{L}h}(s)\big]\\
&=C\int_{\mathbb{R}^9}K_{\mathfrak{m}}(\xi,\eta,\rho)\widehat{U^{h_1,\iota_1}_{\leq J_1,k_1}}(\xi-\eta-\rho,s)\widehat{P_{k_2}\mathcal{U}^{\mathcal{L}_2h_2}}(\eta,s)\overline{\widehat{P_k\mathcal{U}^{\mathcal{L}h}}(\xi,s)}\, d\xi d\eta d\rho,
\end{split}
\end{equation}
where
\begin{equation}\label{sori29}
\begin{split}
K_{\mathfrak{m}}(\xi,\eta,\rho):&=\frac{\mathfrak{m}^{nr}(\xi-\eta-\rho,\eta+\rho)}{|\xi|-|\eta+\rho|-\iota_1|\xi-\eta-\rho|}\widetilde{\Sigma^{\geq 1}_{wa}}\Big(\rho,\frac{2\eta+\rho}{2}\Big)\chi_0\Big(\frac{|\rho|}{|2\eta+\rho|}\Big)\\
&-\frac{\mathfrak{m}^{nr}(\xi-\eta-\rho,\eta)}{|\xi-\rho|-|\eta|-\iota_1|\xi-\eta-\rho|}\widetilde{\Sigma^{\geq 1}_{wa}}\Big(\rho,\frac{2\xi-\rho}{2}\Big)\chi_0\Big(\frac{|\rho|}{|2\xi-\rho|}\Big).
\end{split}
\end{equation}

{\bf{Step 4.}} To prove \eqref{sori24} it suffices to show that
\begin{equation}\label{pint2}
\begin{split}
\Big|\int_{\mathbb{R}^9}K_{\mathfrak{m}}(\xi,\eta,\rho)\widehat{U^{h_1,\iota_1}_{\leq J_1,k_1}}(\xi-\eta-\rho,s)\widehat{P_{k_2}f_2}(\eta)\overline{\widehat{P_kf}(\xi)}\, d\xi d\eta d\rho\Big|\lesssim \varep_12^{-3m/2}\|P_kf\|_{L^2}\|P_{k_2}f_2\|_{L^2},
\end{split}
\end{equation}
for any $s\in J_m$ and $f,f_2\in L^2$. The main issue is the possible loss of derivative, so one should think of $|\xi|, |\eta|\in[2^{k-4},2^{k+4}]$ as large and $|\rho|, |\xi-\eta-\rho|\leq 2^{k-20}$ as small. For $a,b\in[0,1]$ let
\begin{equation}\label{pint2.5}
\begin{split}
S(\xi,\eta,\rho,s;a,b):=&\frac{\mathfrak{m}(\xi-\eta-\rho,\eta+a\rho)\varphi_{>q_0}(\Xi_{\iota_1+}(\xi-\eta-\rho,\eta+a\rho))}{|\xi-\rho+a\rho|-|\eta+a\rho|-\iota_1|\xi-\eta-\rho|}\\
&\times\widetilde{\Sigma^{\geq 1}_{wa}}\Big(\rho,\frac{2\xi-\rho}{2}-b(\xi-\eta-\rho),s\Big)\chi_0\Big(\frac{|\rho|}{|2\xi-\rho-2b(\xi-\eta-\rho)|}\Big),
\end{split}
\end{equation}
such that $K_{\mathfrak{m}}(\xi,\eta,\rho,s)=S(\xi,\eta,\rho,s;1,1)-S(\xi,\eta,\rho,s;0,0)$. It suffices to prove that
\begin{equation}\label{pint3}
\begin{split}
\Big|\int_{\mathbb{R}^9}\nabla_{a,b}S(\xi,\eta,\rho,s;a,b)\widehat{U^{h_1,\iota_1}_{\leq J_1,k_1}}(\xi-\eta-\rho,s)\widehat{P_{k_2}f_2}(\eta)\overline{\widehat{P_kf}(\xi)}\, d\xi d\eta d\rho\Big|\\
\lesssim \varep_12^{-3m/2}\|P_kf\|_{L^2}\|P_{k_2}f_2\|_{L^2},
\end{split}
\end{equation}
for any $a,b\in[0,1]$. 

We notice that the symbol $\Sigma^{\geq 1}_{wa}(x,\zeta)=\Sigma_{wa}(x,\zeta)-|\zeta|$ can be decomposed as a sum of symbols of the form $A_{d,l}G_{d,l}(x)\mu_{d,l}(\zeta)$, $d\geq 1$, $l\in\{1,\ldots, L(d)\}$, where $G_{d,l}\in\mathcal{G}_d$ (see definition \eqref{gb0.1}), $\mu_{d,l}$ are smooth homogeneous multipliers of order $1$, and the constants $A_{d,l}$ and $L(d)$ are bounded by $C^d$. Using the $L^\infty$ norms in \eqref{gb5} and \eqref{sori10}, together with the general estimate \eqref{ener62}, for \eqref{pint3} it suffices to bound 
\begin{equation}\label{pint4}
\|\mathcal{F}^{-1} M^{l,l_1,l_3}_{a,b}\|_{L^1(\mathbb{R}^9)}\lesssim 2^{m/4}2^{6l_3^+}
\end{equation}
for any $a,b\in[0,1]$ and integers $l\geq m/20$, $l_1\in[-2\delta' m,2\delta'm]$, and $l_3\leq l-20$, where
\begin{equation}\label{pint5}
\begin{split}
M_{a,b}^{l,l_1,l_3}(\eta,\rho,\theta):=\nabla_{a,b}\Big\{&\frac{\mathfrak{m}(\theta,\eta+a\rho)\varphi_{>q_0}(\Xi_{\iota_1+}(\theta,\eta+a\rho))}{|\theta+\eta+a\rho|-|\eta+a\rho|-\iota_1|\theta|}A(\eta+(1-b)\theta+\rho/2)\\
&\times\chi_0\Big(\frac{|\rho|}{|2\eta+2(1-b)\theta+\rho|}\Big)\Big\}\varphi_{l}(\eta)\varphi_{l_1}(\theta)\varphi_{l_3}(\rho).
\end{split}
\end{equation}
This multiplier is obtained from the expression in \eqref{pint2.5} by making the change of variables $\xi=\theta+\eta+\rho$, and $A:\mathbb{R}^3\to\mathbb{R}$ is a smooth homogeneous function of order $1$. 

To prove \eqref{pint4} we use first \eqref{par7}, thus
\begin{equation*}
\frac{\varphi_{>q_0}(\Xi_{\iota_1+}(\theta,\eta+a\rho))}{|\theta+\eta+a\rho|-|\eta+a\rho|-\iota_1|\theta|}=\frac{|\theta+\eta+a\rho|+|\eta+a\rho|+\iota_1|\theta|}{-\iota_1|\theta||\eta+a\rho|}\frac{\varphi_{>q_0}(\Xi_{\iota_1+}(\theta,\eta+a\rho))}{|\Xi_{\iota_1+}(\theta,\eta+a\rho)|^2}.
\end{equation*}
We use now \eqref{nag1} and recall the assumption $2^{-q_0}\lesssim 2^{8\delta'm}$. The desired bounds \eqref{pint4} follow by examining the terms resulting from taking the derivatives in $a$ or $b$, and recalling the algebra property $\|\mathcal{F}^{-1}(m\cdot m')\|_{L^1}\lesssim \|\mathcal{F}^{-1}(m)\|_{L^1}\|\mathcal{F}^{-1}(m')\|_{L^1}$. 
\end{proof}

Our last lemma concludes the analysis in {\bf{Case 4.4.}} in the proof of Lemma \ref{NullLHHWWWS}.

\begin{lemma}\label{ParalinearSym2}
If $\iota_1\in\{+,-\}$, $t\in[0,T]$, $m\in[\delta^{-2},L]$, $\mathcal{L},\LL_1\in\mathcal{V}^q_n$, $n\geq 1$, and 
\begin{equation}\label{jik40}
k_1\in[-m+(Y(q,n)-1)\delta m,-0.6m],\qquad |k|,|k_2|\leq\delta'm,\qquad b\in[-2\delta m,4],
\end{equation}
then, with $h,h_1,h_2\in\{h_{\al\be}\}$ and $\mathfrak{q}_{\iota_1\iota_2}^b$ defined as in \eqref{jik6} and \eqref{doubnull},
\begin{equation*}
\begin{split}
2^{-k_1}\Big|\int_{J_m}q_m(s)\mathcal{G}_{\mathfrak{q}_{\iota_1+}^b}\big[P_{k_1}U^{\mathcal{L}_1h_1,\iota_1}(s),P_{k_2}U^{h_2}(s),P_kU^{\mathcal{L}h}(s)\big]\,ds\Big|\lesssim \varepsilon_1^32^{-2N(n)k^+}2^{2H(q,n)\delta m}2^{-\delta m}.
\end{split}
\end{equation*}
\end{lemma}

\begin{proof} As in the proof of Lemma \ref{ParalinearSym} we replace first the solutions $U^{h_2}$ and $U^{\LL h}$ with the quasilinear variables $\mathcal{U}^{h_2}=(\partial_t-iT_{\sigma_{wa}})h_{2}$ and $\mathcal{U}^{\LL h}=(\partial_t-iT_{\sigma_{wa}})(\LL h)$ defined in \eqref{MathcalU}, at the expense of acceptable errors that can be estimated as in \eqref{sori19}. It remains to prove that
\begin{equation}\label{jik41}
\begin{split}
2^{-k_1}\Big|\int_{J_m}q_m(s)\mathcal{G}_{\mathfrak{q}_{\iota_1+}^b}\big[P_{k_1}U^{\mathcal{L}_1h_1,\iota_1}(s),P_{k_2}\mathcal{U}^{h_2}(s),P_k\mathcal{U}^{\mathcal{L}h}(s)\big]\,ds\Big|\lesssim \varepsilon_1^32^{-2N(n)k^+}2^{2H(q,n)\delta m-\delta m}.
\end{split}
\end{equation}
Then we apply the integration by parts identity \eqref{NF2}. With $\mathcal{H}_{\mathfrak{q}_{\iota_1+}^b}$ defined as in \eqref{DefJ}, for \eqref{jik41} it suffices to prove that
\begin{equation}\label{jik42}
\big|\mathcal{H}_{\mathfrak{q}_{\iota_1+}^b}\big[P_{k_1}U^{\mathcal{L}_1h_1,\iota_1}(s),P_{k_2}\mathcal{U}^{h_2}(s),P_k\mathcal{U}^{\mathcal{L}h}(s)\big]\big|\lesssim \varepsilon_1^32^{k_1}2^{-2N(n)k^+}2^{2H(q,n)\delta m-\delta m},
\end{equation}
\begin{equation}\label{jik43}
2^m\big|\mathcal{H}_{\mathfrak{q}_{\iota_1+}^b}\big[(\partial_s+i\Lambda_{wa,\iota_1})P_{k_1}U^{\mathcal{L}_1h_1,\iota_1}(s),P_{k_2}\mathcal{U}^{h_2}(s),P_k\mathcal{U}^{\mathcal{L}h}(s)\big]\big|\lesssim \varepsilon_1^32^{k_1}2^{-2N(n)k^+}2^{2H(q,n)\delta m-\delta m},
\end{equation}
and
\begin{equation}\label{jik44}
\begin{split}
&2^m\big|\mathcal{H}_{\mathfrak{q}_{\iota_1+}^b}\big[P_{k_1}U^{\mathcal{L}_1h_1,\iota_1}(s),(\partial_s+i\Lambda_{wa})P_{k_2}\mathcal{U}^{h_2}(s),P_k\mathcal{U}^{\mathcal{L}h}(s)\big]\\
&+\mathcal{H}_{\mathfrak{q}_{\iota_1+}^b}\big[P_{k_1}U^{\mathcal{L}_1h_1,\iota_1}(s),P_{k_2}\mathcal{U}^{h_2}(s),(\partial_s+i\Lambda_{wa})P_k\mathcal{U}^{\mathcal{L}h}(s)\big]\big|\lesssim \varepsilon_1^32^{k_1}2^{-2N(n)k^+}2^{2H(q,n)\delta m-\delta m},
\end{split}
\end{equation}
for any $s\in J_m$. We prove these estimate in several steps.

{\bf{Step 1.}} We start with the easier estimates \eqref{jik42} and \eqref{jik43}. The main point is that
\begin{equation}\label{jik45}
\|P_{k_1}U^{\mathcal{L}_1h_1,\iota_1}(s)\|_{L^2}+2^m\|(\partial_s+i\Lambda_{wa,\iota_1})P_{k_1}U^{\mathcal{L}_1h_1,\iota_1}(s)\|_{L^2}\lesssim \eps_12^{k_1/2}2^{H(q,n)\delta m+Y'(n)\delta m},
\end{equation}
where $Y'(1):=2$ and $Y'(2)=Y'(3):=35$. Indeed, these bounds follow from \eqref{vcx1} and \eqref{wer4.0} when $n\geq 2$. If $n=1$ they follow from \eqref{vcx1}, Lemma \ref{Onev1} (recall $k_1\leq -0.6m$), and \eqref{abc3.00}.

Using \eqref{nag1}, \eqref{par7}, and \eqref{doubnull} we have
\begin{equation}\label{jik47}
\Big\|\mathcal{F}^{-1}\Big\{\frac{\mathfrak{q}_{\iota_1+}^b(\xi-\eta,\eta)}{|\xi|-|\eta|-\iota_1|\xi-\eta|}\varphi_{kk_1k_2}(\xi-\eta,\eta)\Big\}\Big\|_{L^1(\mathbb{R}^6)}\lesssim 2^{-k_1}2^{-b}.
\end{equation}
With $J_2=-k_1$ we decompose, as in \eqref{on11.3}--\eqref{on11.36},
\begin{equation}\label{jik48}
P_{k_2}\mathcal{U}^{h_2}=P_{k_2}U^{h_2}+P_{k_2}(\mathcal{U}^{h_2}-U^{h_2})=U^{h_2,+}_{\leq J_2,k_2}+U^{h_2,+}_{>J_2,k_2}+P_{k_2}(\mathcal{U}^{h_2}-U^{h_2}).
\end{equation}
For $G(s)\in \{U^{\mathcal{L}_1h_1,\iota_1}(s),2^m(\partial_s+i\Lambda_{wa,\iota_1})P_{k_1}U^{\mathcal{L}_1h_1,\iota_1}(s)\}$, we estimate, using \eqref{bil1} and \eqref{jik47},
\begin{equation}\label{jik50}
\begin{split}
\big|\mathcal{H}_{\mathfrak{q}_{\iota_1+}^b}&\big[G(s),U^{h_2,+}_{\leq J_2,k_2}(s),P_k\mathcal{U}^{\mathcal{L}h}(s)\big]\big|\\
&\lesssim 2^{k_1/2}2^{-m}2^{3k/2}2^{-k_1}2^{-b}\|G(s)\|_{L^2}\|\widehat{P_{k_2}U^{h_2}}(s)\|_{L^\infty}\|P_k\mathcal{U}^{\mathcal{L}h}(s)\|_{L^2}\\
&\lesssim \varepsilon_1^32^{k^-/2}2^{-2N(n)k^+-4k^+}2^{2H(q,n)\delta m}2^{-m+Y'(n)\delta m+\delta m}2^{-b},
\end{split}
\end{equation}
using also \eqref{vcx1}, \eqref{vcx1.2}, \eqref{jik45} in the last line. We also estimate, using just $L^2$ bounds,
\begin{equation}\label{jik51}
\begin{split}
&\big|\mathcal{H}_{\mathfrak{q}_{\iota_1+}^b}\big[G(s),U^{h_2,+}_{>J_2,k_2}(s)+P_{k_2}(\mathcal{U}^{h_2}-U^{h_2})(s),P_k\mathcal{U}^{\mathcal{L}h}(s)\big]\big|\\
&\lesssim 2^{3k_1/2}2^{-k_1}2^{-b}\|G(s)\|_{L^2}\big\{\|U^{h_2,+}_{>J_2,k_2}(s)\|_{L^2}+\|P_{k_2}(\mathcal{U}^{h_2}-U^{h_2})(s)\|_{L^2}\big\}\|P_k\mathcal{U}^{\mathcal{L}h}(s)\|_{L^2}\\
&\lesssim \varepsilon_1^32^{2k_1}2^{2\delta'm},
\end{split}
\end{equation}
using also \eqref{vcx1}--\eqref{vcx1.1}, \eqref{DefParaUnknown}, \eqref{jik45} in the last line. The desired bounds \eqref{jik42}--\eqref{jik43} follow once we notice that $-m+Y'(n)\delta m+\delta m-b\leq k_1-6\delta m$, due to \eqref{jik40} and \eqref{saur10}.

{\bf{Step 2.}} We consider now the estimates \eqref{jik44}. We decompose
\begin{equation*}
\begin{split}
(\partial_s+i\Lambda_{wa})P_{k_2}\mathcal{U}^{h_2}=P_{k_2}(\partial_s+iT_{\Sigma_{wa}})\mathcal{U}^{h_2}-iP_{k_2}T_{\Sigma_{wa}-|\ze|}\mathcal{U}^{h_2},\\
P_k(\partial_s+i\Lambda_{wa})\mathcal{U}^{\mathcal{L}h}=P_k(\partial_s+iT_{\Sigma_{wa}})\mathcal{U}^{\mathcal{L}h}-iP_kT_{\Sigma_{wa}-|\ze|}\mathcal{U}^{\mathcal{L}h}.
\end{split}
\end{equation*}
For \eqref{jik44} it suffices to prove that 
\begin{equation}\label{jik55}
2^m\big|\mathcal{H}_{\mathfrak{q}_{\iota_1+}^b}\big[P_{k_1}U^{\mathcal{L}_1h_1,\iota_1}(s),P_{k_2}(\partial_s+iT_{\Sigma_{wa}})\mathcal{U}^{h_2}(s),P_k\mathcal{U}^{\mathcal{L}h}(s)\big]\big|\lesssim \varepsilon_1^32^{k_1}2^{-2N(n)k^+}2^{2H(q,n)\delta m-\delta m},
\end{equation}
\begin{equation}\label{jik56}
2^m\big|\mathcal{H}_{\mathfrak{q}_{\iota_1+}^b}\big[P_{k_1}U^{\mathcal{L}_1h_1,\iota_1}(s),P_{k_2}\mathcal{U}^{h_2}(s),P_k(\partial_s+iT_{\Sigma_{wa}})\mathcal{U}^{\mathcal{L}h}(s)\big]\big|\lesssim \varepsilon_1^32^{k_1}2^{-2N(n)k^+}2^{2H(q,n)\delta m-\delta m},
\end{equation}
and
\begin{equation}\label{jik57}
\begin{split}
&2^m\big|\mathcal{H}_{\mathfrak{q}_{\iota_1+}^b}\big[P_{k_1}U^{\mathcal{L}_1h_1,\iota_1}(s),P_{k_2}T_{\Sigma_{wa}-|\ze|}\mathcal{U}^{h_2}(s),P_k\mathcal{U}^{\mathcal{L}h}(s)\big]\\
&-\mathcal{H}_{\mathfrak{q}_{\iota_1+}^b}\big[P_{k_1}U^{\mathcal{L}_1h_1,\iota_1}(s),P_{k_2}\mathcal{U}^{h_2}(s),P_kT_{\Sigma_{wa}-|\ze|}\mathcal{U}^{\mathcal{L}h}(s)\big]\big|\lesssim \varepsilon_1^32^{k_1}2^{-2N(n)k^+}2^{2H(q,n)\delta m-\delta m}.
\end{split}
\end{equation}

We notice that the bounds \eqref{jik56} follow from \eqref{ParaBox5.01}, using the decomposition \eqref{jik48} and estimating as in \eqref{jik50}--\eqref{jik51}.

{\bf{Step 3.}} We prove now the bounds \eqref{jik55}. We use the formulas \eqref{ParaBox1} and \eqref{ParaBox6} with $\LL=Id$. The contribution of the cubic and higher order terms can be bounded easily, proceeding as in \eqref{jik51}. To control the main contributions we will prove that
\begin{equation}\label{jik58}
\begin{split}
2^m\big|\mathcal{Q}_\mathfrak{p}[P_{k_1}U^{\mathcal{L}_1h_1,\iota_1}(s),&P_{k_3}U^{h_3,\iota_3}(s),P_{k_4}U^{h_4,\iota_4}(s),P_k\mathcal{U}^{\mathcal{L}h}(s)]\big|\\
&\lesssim \varepsilon_1^42^{2k_1}2^{-2N(n)k^+}2^{2H(q,n)\delta m-10\delta m}2^{k^-/2}2^{k_3^-/4}2^{k_4^-/2-4k_4^+},
\end{split}
\end{equation}
and
\begin{equation}\label{jik59}
\begin{split}
2^m\big|\mathcal{Q}_\mathfrak{p}[P_{k_1}U^{\mathcal{L}_1h_1,\iota_1}(s),&P_{k_3}U^{\psi,\iota_3}(s),P_{k_4}U^{\psi,\iota_4}(s),P_k\mathcal{U}^{\mathcal{L}h}(s)]\big|\\
&\lesssim \varepsilon_1^42^{2k_1}2^{-2N(n)k^+}2^{2H(q,n)\delta m-10\delta m}2^{k^-/2}2^{k_3^-/4}2^{k_4^-/2-4k_4^+},
\end{split}
\end{equation}
for any $s\in J_m$, $\iota_3,\iota_4\in\{+,-\}$, $h_3,h_4\in\{h_{\al\be}\}$, $k_3\leq k_4\in\mathbb{Z}$. Here
\begin{equation}\label{jik60}
\mathcal{Q}_\mathfrak{p}[f_1,f_3,f_4,f]:=\int_{(\mathbb{R}^3)^3}\mathfrak{p}(\xi-\eta,\eta-\rho,\rho)\widehat{f_1}(\xi-\eta)\widehat{f_3}(\eta-\rho)\widehat{f_4}(\rho)\overline{\widehat{f}(\xi)}\,d\xi d\eta d\rho,
\end{equation}
and $\mathfrak{p}$ is a multiplier satisfying $\|\mathcal{F}^{-1}\mathfrak{p}\|_{L^1(\mathbb{R}^9)}\leq 1$. These bounds clearly suffice to prove \eqref{jik55}; they are in fact stronger than needed because we would like to apply them in the proof of the estimates \eqref{jik57} as well.

{\bf{Substep 3.1.}} We prove first the bounds \eqref{jik58}. Since $k_3\leq k_4$ we may assume that $k_4\geq k-8$. Using \eqref{ener62} we estimate first the left-hand side of \eqref{jik58} by
\begin{equation}\label{jik61}
\begin{split}
C2^m2^{3k_1/2}\|P_{k_1}U^{\mathcal{L}_1h_1,\iota_1}(s)\|_{L^2}&\|P_{k_3}U^{h_3,\iota_3}(s)\|_{L^\infty}\|P_{k_4}U^{h_4,\iota_4}(s)\|_{L^2}\|P_k\mathcal{U}^{\mathcal{L}h}(s)\|_{L^2}\\
&\lesssim \varep_1^42^{2k_1}2^{\delta'm}2^{k_3^-}2^{-N(0)k_4^++2k_4^+}2^{k^-/2}2^{-N(n)k^+},
\end{split}
\end{equation}
where we used \eqref{vcx1} and \eqref{wws1} in the second line. This suffices if either $k_3\leq -8\delta'm$ or $k_4\geq 8\delta'm$. On the other hand, if $k_3,k_4\in[-8\delta'm,8\delta'm]$ then we fix $J_3=J_4$ the largest integer smaller than $m/4$ and decompose $P_{k_3}U^{h_3,\iota_3}=U^{h_3,\iota_3}_{\leq J_3,k_3}+U^{h_3,\iota_3}_{>J_3,k_3}$ and $P_{k_4}U^{h_4,\iota_4}=U^{h_4,\iota_4}_{\leq J_4,k_4}+U^{h_4,\iota_4}_{>J_4,k_4}$ as in \eqref{on11.3}--\eqref{on11.36}. The contributions of the functions $U^{h_3,\iota_3}_{>J_3,k_3}$ and $U^{h_4,\iota_4}_{>J_4,k_4}$ can be estimated easily, using \eqref{vcx1.1}. After these reductions it remains to prove that
\begin{equation}\label{jik63}
\begin{split}
2^m\big|\mathcal{Q}_\mathfrak{p}[P_{k_1}U^{\mathcal{L}_1h_1,\iota_1}(s),&U^{h_3,\iota_3}_{\leq J_3,k_3}(s),U^{h_4,\iota_4}_{\leq J_4,k_4}(s),P_k\mathcal{U}^{\mathcal{L}h}(s)]\big|\\
&\lesssim \varepsilon_1^42^{2k_1}2^{-2N(n)k^+}2^{2H(q,n)\delta m-10\delta m}2^{k^-/2}2^{k_3^-/4}2^{k_4^-/2-4k_4^+},
\end{split}
\end{equation}
for any $s\in J_m$ and $k_3\leq k_4\in[-8\delta'm,8\delta'm]$.

To prove \eqref{jik63} we examine the formula \eqref{jik60} and write
\begin{equation*}
\begin{split}
\mathcal{Q}_\mathfrak{p}[f_1,f_3,f_4,f]=C\int_{\mathbb{R}^9}\int_{\mathbb{R}^9}K(x,y,z)e^{-ix\cdot \rho}e^{-iy\cdot(\xi-\rho)}e^{-iz\cdot (\eta-\xi)}\\
\times\widehat{f_1}(\rho)\widehat{f_3}(\xi-\rho)\widehat{f_4}(\eta-\xi)\overline{\widehat{f}(\eta)}\,d\xi d\eta d\rho\, dx dy dz.
\end{split}
\end{equation*}
after suitable changes of variables, where $K=\mathcal{F}^{-1}(\mathfrak{p})$. Since $\|K\|_{L^1}\lesssim 1$, we have
\begin{equation}\label{jik65}
\begin{split}
\big|\mathcal{Q}_\mathfrak{p}[f_1,f_3,f_4,f]\big|&\lesssim \sup_{x,y,z\in\mathbb{R}^3}\Big|\int_{\mathbb{R}^9}e^{ix\cdot\xi}(\widehat{f_1}(\rho)e^{iy\rho})\widehat{f_3}(\xi-\rho)\widehat{f_4}(\eta-\xi)\overline{(\widehat{f}(\eta)e^{iz\cdot\eta})}\,d\xi d\eta d\rho\Big|\\
&\lesssim  \sup_{y,z\in\mathbb{R}^3}\big\|I[f_1(.-y),f_3]\big\|_{L^2}\big\|I[\overline{f_4},f(.-z)]\big\|_{L^2},
\end{split}
\end{equation}
where $I$ is defined as in \eqref{abc36.1} with the multiplier $m$ equal to $1$. 

In our case, we apply \eqref{jik65} and \eqref{bil1} to estimate the left-hand side of \eqref{jik63} by
\begin{equation*}
\begin{split}
C2^m(2^{k_1/2}2^{-m}&2^{3k_3/2}\|P_{k_1}U^{\mathcal{L}_1h_1}(s)\|_{L^2}\|\widehat{P_{k_3}U^{h_3}}(s)\|_{L^\infty})\\
&\times(2^{k/2}2^{-m}2^{3k_4/2}\|P_k\mathcal{U}^{\mathcal{L}h}(s)\|_{L^2}\|\widehat{P_{k_4}U^{h_4}}(s)\|_{L^\infty})\\
&\lesssim \varep_1^42^{-m}2^{2H(q,n)\delta m+2\delta m}2^{k_1}2^{-N(n)k^++2k^+}2^{-N_0k_4^++2k_4^+}2^{k^-/2}2^{k_3^-/4}2^{k_4^-/2},
\end{split}
\end{equation*}
where we used \eqref{vcx1} and \eqref{vcx1.2} in the last line. This gives the claimed bounds \eqref{jik63} once we recall that $k_4\geq k-8$ and $2^{-m+(Y(q,n)-1)\delta m}\leq 2^{k_1}$, see \eqref{jik40} and recall \eqref{saur10}.

{\bf{Substep 3.2.}} We prove now the bounds \eqref{jik59}. Estimating as in \eqref{jik61}, this is easy using the $L^\infty$ estimates \eqref{wws2} unless $k_3,k_4\in[-8\delta'm,8\delta'm]$. In this case we fix $J_3=J_4$ the largest integer smaller than $m/4$, as before, and reduce to proving that 
\begin{equation*}
\begin{split}
2^m\big|\mathcal{Q}_\mathfrak{p}[&P_{k_1}U^{\mathcal{L}_1h_1,\iota_1}(s),U^{\psi,\iota_3}_{\leq J_3,k_3}(s),U^{\psi,\iota_4}_{\leq J_4,k_4}(s),P_k\mathcal{U}^{\mathcal{L}h}(s)]\big|\\
&\lesssim \varepsilon_1^42^{2k_1}2^{-2N(n)k^+}2^{2H(q,n)\delta m-10\delta m}2^{k^-/2}2^{k_3^-/4}2^{k_4^-/2-4k_4^+}
\end{split}
\end{equation*}
for any $s\in J_m$ and $k_3\leq k_4\in[-8\delta'm,8\delta'm]$. These bounds follow easily using \eqref{wws13x}.

{\bf{Step 4.}} Finally we prove the bounds \eqref{jik57}. We write $\Sigma_{wa}-|\zeta|=|\zeta|\Sigma_{wa}^1+|\zeta|\Sigma_{wa}^{\geq 2}$, as in \eqref{ParaSymb2}. The contribution of the symbol $\Sigma_{wa}^{\geq 2}$ leads to higher order terms that can be estimated using just $L^2$ bounds. To bound the main term we write, as in \eqref{sori28}--\eqref{sori29},
\begin{equation}\label{jik70}
\begin{split}
&\mathcal{H}_{\mathfrak{q}_{\iota_1+}^b}\big[P_{k_1}U^{\mathcal{L}_1h_1,\iota_1}(s),P_{k_2}T_{|\ze|\Sigma_{wa}^1}\mathcal{U}^{h_2}(s),P_k\mathcal{U}^{\mathcal{L}h}(s)\big]\\
&-\mathcal{H}_{\mathfrak{q}_{\iota_1+}^b}\big[P_{k_1}U^{\mathcal{L}_1h_1,\iota_1}(s),P_{k_2}\mathcal{U}^{h_2}(s),P_kT_{|\ze|\Sigma^1_{wa}}\mathcal{U}^{\mathcal{L}h}(s)\big]\\
&=C\int_{\mathbb{R}^9}A(\xi,\eta,\rho)\widehat{P_{k_1}U^{\mathcal{L}_1h_1,\iota_1}}(\xi-\eta-\rho,s)\widehat{P'_{k_2}\mathcal{U}^{h_2}}(\eta,s)\overline{\widehat{P'_k\mathcal{U}^{\mathcal{L}h}}(\xi,s)}\, d\xi d\eta d\rho,
\end{split}
\end{equation}
where $P'_k=P_{[k-2,k+2]}$, $P'_{k_2}=P_{[k_2-2,k_2+2]}$, and
\begin{equation}\label{jik71}
\begin{split}
A(\xi,\eta,\rho):&=\frac{\mathfrak{q}^{b}_{\iota_1+}(\xi-\eta-\rho,\eta+\rho)}{|\xi|-|\eta+\rho|-\iota_1|\xi-\eta-\rho|}\varphi_{k_2}(\eta+\rho)\widetilde{(|\ze|\Sigma^{1}_{wa})}\Big(\rho,\frac{2\eta+\rho}{2}\Big)\chi_0\Big(\frac{|\rho|}{|2\eta+\rho|}\Big)\varphi_k(\xi)\\
&-\frac{\mathfrak{q}_{\iota_1+}^{b}(\xi-\eta-\rho,\eta)}{|\xi-\rho|-|\eta|-\iota_1|\xi-\eta-\rho|}\varphi_{k_2}(\eta)\widetilde{(|\ze|\Sigma^{1}_{wa})}\Big(\rho,\frac{2\xi-\rho}{2}\Big)\chi_0\Big(\frac{|\rho|}{|2\xi-\rho|}\Big)\varphi_k(\xi-\rho).
\end{split}
\end{equation}

The formula \eqref{ParaSymb1} shows that $\widetilde{(|\ze|\Sigma^{1}_{wa})}\big(\rho,v\big)$ is a sum of expressions of the form $\widehat{h_3}(\rho)g(v)$, where $h_3\in\{h_{\al\be}\}$ and $g(v)$ is either $|v|$ or $|v|\widehat{v_j}$, or $|v|\widehat{v_j}\widehat{v_k}$. As in \eqref{pint2.5}, for $x,y\in[0,1]$ let
\begin{equation}\label{jik73}
\begin{split}
B(\xi,\eta,&\rho;x,y):=\frac{\mathfrak{q}_{\iota_1+}(\xi-\eta-\rho,\eta+x\rho)\varphi_{b}(\Xi_{\iota_1+}(\xi-\eta-\rho,\eta+x\rho))}{|\xi-\rho+x\rho|-|\eta+x\rho|-\iota_1|\xi-\eta-\rho|}\varphi_{k_2}(\eta+x\rho)\\
&\times g\Big(\frac{2\xi-\rho}{2}-y(\xi-\eta-\rho)\Big)\chi_0\Big(\frac{|\rho|}{|2\xi-\rho-2y(\xi-\eta-\rho)|}\Big)\varphi_k(\xi-\rho+x\rho)
\end{split}
\end{equation}
such that $A(\xi,\eta,\rho)$ is a sum over $h_3$ and $g$ of expressions of the form $\widehat{h_3}(\rho)[B(\xi,\eta,\rho;1,1)-B(\xi,\eta,\rho;0,0)]$. In view of these identities, for \eqref{jik57} it suffices to show that
\begin{equation}\label{jik75}
\begin{split}
2^m\sum_{k_3\leq k-10}&\Big|\int_{\mathbb{R}^9}\nabla_{x,y}B(\xi,\eta,\rho;x,y)|\rho|^{-1}\widehat{P_{k_1}U^{\mathcal{L}_1h_1,\iota_1}}(\xi-\eta-\rho,s)\widehat{P_{k_3}U^{h_3,\iota_3}}(\rho,s)\\
&\times\widehat{P'_{k_2}\mathcal{U}^{h_2}}(\eta,s)\overline{\widehat{P'_k\mathcal{U}^{\mathcal{L}h}}(\xi,s)}\, d\xi d\eta d\rho\Big|\lesssim \varepsilon_1^32^{k_1}2^{-2N(n)k^+}2^{2H(q,n)\delta m-\delta m},
\end{split}
\end{equation}
for any $x,y\in[0,1]$ and $\iota_3\in\{+,-\}$.

Let $B'(\xi,\eta,\rho):=\nabla_{x,y}B(\xi,\eta,\rho;x,y)$, for $x,y\in[0,1]$. To prove \eqref{jik75} we notice that 
\begin{equation}\label{jik76}
\big\|\mathcal{F}^{-1}\{B'(\xi,\eta,\rho)|\rho|^{-1}\varphi_{l_3}(\rho)\varphi_{l_1}(\xi-\eta-\rho)\varphi_{l_2}(\eta)\varphi_l(\xi)\}\big\|_{L^1(\mathbb{R}^9)}\lesssim 2^{-2b}(2^{-l_3}+2^{-l_1})
\end{equation}
for any $x,y\in[0,1]$, where $l,l_1,l_3,l_2\in\mathbb{Z}$, $l,l_2\in[-2\delta'm,2\delta'm]$, $|l-l_2|\leq 4$, $l_1\leq -0.6m+4$, $l_3\leq l-10$. Indeed, one can think of $2^l$ and $2^{l_2}$ as large and comparable, and $2^{l_1},2^{l_3}$ as small. Using \eqref{par7} we rewrite
\begin{equation}\label{jik78}
\begin{split}
B(\xi,\eta,&\rho;x,y)|\rho|^{-1}=\frac{\mathfrak{q}_{\iota_1+}(\xi-\eta-\rho,\eta+x\rho)\varphi_{b}(\Xi_{\iota_1+}(\xi-\eta-\rho,\eta+x\rho))}{-\iota_1|\Xi_{\iota_1+}(\xi-\eta-\rho,\eta+x\rho)|^2} \\
&\times\frac{[|\xi-\rho+x\rho|+|\eta+x\rho|+\iota_1|\xi-\eta-\rho|]\varphi_{k_2}(\eta+x\rho)\varphi_k(\xi-\rho+x\rho)}{|\rho||\eta+x\rho|\cdot|\xi-\eta-\rho|}\\
&\times g\Big(\frac{2\xi-\rho}{2}-y(\xi-\eta-\rho)\Big)\chi_0\Big(\frac{|\rho|}{|2\xi-\rho-2y(\xi-\eta-\rho)|}\Big).
\end{split}
\end{equation}
Using \eqref{nag1} and the double-null assumption \eqref{doubnull}, it is easy to see that
\begin{equation}\label{jik79}
\big\|\mathcal{F}^{-1}\{B(\xi,\eta,\rho;x,y)|\rho|^{-1}\varphi_{l_3}(\rho)\varphi_{l_1}(\xi-\eta-\rho)\varphi_{l_2}(\eta)\varphi_l(\xi)\}\big\|_{L^1(\mathbb{R}^9)}\lesssim 2^{-b}2^l2^{-l_1-l_3}
\end{equation}
for any $x,y\in[0,1]$ and $l,l_1,l_3,l_2$ as above. Taking $x$ derivatives generates factors $\lesssim 2^{-b}2^{l_3-l}$ from the terms in the first two lines of \eqref{jik78}, while taking $y$ derivatives generates factors $\lesssim 2^{l_1-l}$ from the terms in the third line of \eqref{jik78}. Combining these estimates yields \eqref{jik76}. 

We can now complete the proof of \eqref{jik75}. In view of \eqref{jik76} it suffices to show that
\begin{equation}\label{jik80}
\begin{split}
&2^m2^{-2b}\sum_{k_3\leq k-10}(2^{-k_1}+2^{-k_3})\Big|\int_{\mathbb{R}^9}\mathfrak{p}'(\rho,\xi-\eta-\rho,\eta)\widehat{P_{k_1}U^{\mathcal{L}_1h_1,\iota_1}}(\xi-\eta-\rho,s)\\
&\times \widehat{P_{k_3}U^{h_3,\iota_3}}(\rho,s)\widehat{P'_{k_2}\mathcal{U}^{h_2}}(\eta,s)\overline{\widehat{P'_k\mathcal{U}^{\mathcal{L}h}}(\xi,s)}\, d\xi d\eta d\rho\Big|\lesssim \varepsilon_1^32^{k_1}2^{-2N(n)k^+}2^{2H(q,n)\delta m-\delta m},
\end{split}
\end{equation}
provided that $\mathfrak{p}'$ is a multiplier satisfying $\|\mathcal{F}^{-1}\mathfrak{p}'\|_{L^1(\mathbb{R}^9)}\leq 1$. The sum over $k_3\geq k_1$ is bounded as claimed due to \eqref{jik58}, while the sum over $k_3\leq k_1$ can be estimated easily as in \eqref{jik61}. This completes the proof of the lemma.
\end{proof}

\section{Mixed Wave-Klein-Gordon interactions}\label{En3mixed}

We consider now the interactions of the metric components and the Klein-Gordon field, and prove the remaining bounds \eqref{EEEst1}--\eqref{BulkKG} in Proposition \ref{totbounds}. 

We start with the semilinear estimates.

\begin{lemma}\label{KKWINt}
With the assumptions of Proposition \ref{totbounds}, for any $m\in\{0,\ldots,L+1\}$ we have
\begin{equation}\label{KKFINt3}
\begin{split}
\sum_{k,k_1,k_2\in\mathbb{Z}}&2^{2N(n)k^+-k}2^{2\ga(m+k^-)}\Big|\int_{J_m}q_m(s)\\
&\times\mathcal{G}_\mathfrak{m}[P_{k_1}U^{\mathcal{L}_1\psi,\iota_1}(s),P_{k_2}U^{\mathcal{L}_2\psi,\iota_2}(s),P_kU^{\mathcal{L}h,\iota}(s)]\,ds\Big|\lesssim \varepsilon_1^32^{2H(q,n)\delta m}
\end{split}
\end{equation}
where $\mathfrak{m}\in\mathcal{M}^\ast$ (see \eqref{multsstar}) and the operators $\mathcal{G}_{\mathfrak{m}}$ are defined as in \eqref{DefI}. Moreover
\begin{equation}\label{KKFINt1}
\begin{split}
\sum_{k,k_1,k_2\in\mathbb{Z}}&2^{2N(n)k^+}2^{k_2^+-k_1}\Big|\int_{J_m}q_m(s)\\
&\times\mathcal{G}_\mathfrak{m}[P_{k_1}U^{\mathcal{L}_1h_1,\iota_1}(s),P_{k_2}U^{\mathcal{L}_2\psi,\iota_2}(s),P_kU^{\mathcal{L}\psi,\iota}(s)]\,ds\Big|\lesssim \varepsilon_1^32^{2H(q,n)\delta m}
\end{split}
\end{equation}
if $n_2<n$. Therefore the bounds \eqref{EEEst1}--\eqref{EkgEst3} hold.
\end{lemma}

\begin{proof} The proofs are similar to some of the proofs in section \ref{En2purewave}, using mainly $L^2$ or $L^\infty$ estimates on the frequency-localized solutions. In some cases we integrate by parts in time, using \eqref{NF2}--\eqref{NF3} and the bounds in Lemma \ref{pha2} on the resulting multipliers.

Recall some of the $L^2$ estimates we proved earlier,
\begin{equation}\label{mnb60.1}
\begin{split}
(2^{l^-}2^m)^{\ga}2^{-l/2}\|P_lU^{\mathcal{K}h}(s)\|_{L^2}&\lesssim \varep_1 2^{H(\mathcal{K})\delta m}2^{-N(n')l^+},\\
\|P_lU^{\mathcal{K}\psi}(s)\|_{L^2}&\lesssim \varep_12^{H(\mathcal{K})\delta m}2^{-N(n')l^+},
\end{split}
\end{equation}
and
\begin{equation}\label{mnb60.2}
\begin{split}
2^m 2^{-l/2}\|P_l(\partial_s+i\Lambda_{wa})U^{\mathcal{K}h}(s)\|_{L^2}&\lesssim \varep_1^22^{H(\mathcal{K})\delta m} 2^{-\widetilde{N}(n')l^++7l^+}\cdot 2^{35\delta m},\\
2^m\|P_l(\partial_s+i\Lambda_{kg})U^{\mathcal{K}\psi}(s)\|_{L^2}&\lesssim \varep_1^22^{H(\mathcal{K})\delta m} 2^{-\widetilde{N}(n')l^++7l^+}\cdot 2^{35\delta m},
\end{split}
\end{equation}
for any $\mathcal{K}\in\mathcal{V}_{n'}^{q'}$, $l\in\mathbb{Z}$, and $s\in J_m$. See \eqref{vcx1} and \eqref{wer4.0}--\eqref{wer4.1}. 

Letting $\underline{k}=\min\{k_1,k_2,k_3\}$ and using just the $L^2$ bounds \eqref{mnb60.1} we find that
\begin{equation}\label{kkf1}
\begin{split}
&2^{2N(n)k^+-k}2^{2\ga(m+k^-)}\big|\mathcal{G}_\mathfrak{m}[P_{k_1}U^{\mathcal{L}_1\psi,\iota_1}(s),P_{k_2}U^{\mathcal{L}_2\psi,\iota_2}(s),P_kU^{\mathcal{L}h,\iota}(s)]\big|\\
&\lesssim \varepsilon_1^3 2^{\underline{k}}2^{\ga(m+k^-)}2^{(\underline{k}-k)/2}2^{[H(\mathcal{L}_1)+H(\mathcal{L}_2)+H(\mathcal{L})]\delta m}2^{-N(n_1)k_1^+-N(n_2)k_2^++N(n)k^+},\\
\end{split}
\end{equation}
and
\begin{equation}\label{kkf2}
\begin{split}
&2^{2N(n)k^++k_2^+-k_1}\big|\mathcal{G}_\mathfrak{m}[P_{k_1}U^{\mathcal{L}_1h_1,\iota_1}(s),P_{k_2}U^{\mathcal{L}_2\psi,\iota_2}(s),P_kU^{\mathcal{L}\psi,\iota}(s)]\big|\\
&\lesssim \varepsilon_1^3 2^{\underline{k}}2^{-\ga(k_1^-+m)}2^{(\underline{k}-k_1)/2}2^{[H(\mathcal{L}_1)+H(\mathcal{L}_2)+H(\mathcal{L})]\delta m}2^{-N(n_1)k_1^+-N(n_2)k_2^++k_2^++N(n)k^+}.
\end{split}
\end{equation}
We prove the main estimates in several steps.

{\bf{Step 1.}} We consider first the case when $\min\{n_1,n_2,n\}\ge1$. The bounds \eqref{kkf1}--\eqref{kkf2} already suffice in this case to bound the contributions of the sums over the triples $(k,k_1,k_2)$ with $\min\{k,k_1,k_2\}\leq -m$ or $\max\{k_1,k_2,k_3\}\geq m/4$, due to the bounds \eqref{SuperlinearH1}. They also suffice to bound the entire sums when $|J_m|\lesssim 1$.

For the remaining contributions we integrate by parts using \eqref{NF2}--\eqref{NF3}. Using also Lemma \ref{pha2} (i) each term in the sum in the left-hand side of \eqref{KKFINt3} is bounded by
\begin{equation*}
\begin{split}
C2^{2N(n)k^+-k}&2^{2\ga(m+k^-)}2^{3\underline{k}/2}2^{-k}2^{2\max\{k^+,k_1^+,k_2^+\}}\\
\sup_{s\in J_m}\big\{&\big[\|P_{k_1}U^{\mathcal{L}_1\psi}(s)\|_{L^2}+2^m\|P_{k_1}(\partial_s+i\Lambda_{kg})U^{\mathcal{L}_1\psi}(s)\|_{L^2}\big]\|P_{k_2}U^{\LL_2 \psi}(s)\|_{L^2}\|P_kU^{\mathcal{L}h}(s)\|_{L^2}\\
&+2^m\|P_{k_1}U^{\mathcal{L}_1\psi}(s)\|_{L^2}\|P_{k_2}(\partial_s+i\Lambda_{kg})U^{\mathcal{L}_2\psi}(s)\|_{L^2}\|P_kU^{\mathcal{L}h}(s)\|_{L^2}\\
&+2^m\|P_{k_1}U^{\mathcal{L}_1\psi}(s)\|_{L^2}\|P_{k_2}U^{\LL_2\psi}(s)\|_{L^2}\|P_k(\partial_s+i\Lambda_{wa})U^{\mathcal{L}h}(s)\|_{L^2}\big\}.
\end{split}
\end{equation*} 
In view of \eqref{mnb60.1}--\eqref{mnb60.2} this is bounded by
\begin{equation*}
C\varepsilon_1^3 2^{\ga(m+k^-)}2^{3(\underline{k}-k)/2}2^{[H(\mathcal{L}_1)+H(\mathcal{L}_2)+H(\mathcal{L})]\delta m}2^{-N(n_1)k_1^+-N(n_2)k_2^++N(n)k^+}\cdot 2^{35\delta m}2^{9\max\{k^+,k_1^+,k_2^+\}}.
\end{equation*}
Since $n>\max\{n_1,n_2\}$ and recalling  the bounds \eqref{SuperlinearH1}, this suffices to prove \eqref{KKFINt3} (compare with \eqref{kkf1}). Similarly, each term in the sum in the left-hand side of \eqref{KKFINt1} is bounded by
\begin{equation*}
C\varepsilon_1^3 2^{3(\underline{k}-k_1)/2}2^{[H(\mathcal{L}_1)+H(\mathcal{L}_2)+H(\mathcal{L})]\delta m}2^{-N(n_1)k_1^+-N(n_2)k_2^++k_2^++N(n)k^+}\cdot 2^{35\delta m}2^{9\max\{k^+,k_1^+,k_2^+\}},
\end{equation*}
compare with \eqref{kkf2}, and the desired bounds \eqref{KKFINt1} follow.

{\bf{Step 2.}} We consider now the case when $\min\{n_1,n_2,n\}=0$. By symmetry, the possibilities are $(n_1=0,\,n\geq n_2\geq 0)$ in \eqref{KKFINt3} or $(n_2=0,\,n\geq n_1\geq 1)$ in \eqref{KKFINt1}. The two possibilities are similar, by changes of variables. More precisely, assume that $0\le q\le n\le 3$ and $\mathcal{L},\mathcal{L}_2\in\mathcal{V}_n^q$, $t\in[0,T]$, and $m\in\{0,\ldots,L+1\}$. For any $k,k_1,k_2\in\mathbb{Z}$ and $\iota,\iota_1,\iota_2\in\{+,-\}$ let
\begin{equation}\label{tat0}
\begin{split}
\mathcal{I}_{m;k,k_1,k_2}:=\int_{J_m}&q_m(s)\mathcal{G}_\mathfrak{m}[P_{k_1}U^{\psi,\iota_1}(s),P_{k_2}U^{\mathcal{L}_2\psi,\iota_2}(s),P_kU^{\mathcal{L}h,\iota}(s)]\,ds
\end{split}
\end{equation}
where $\|\mathcal{F}^{-1}(\mathfrak{m})\|_{L^1}\leq 1$. We will show that
\begin{equation}\label{tat1}
\sum_{k,k_1,k_2\in\mathbb{Z}}2^{2\ga(m+k^-)}2^{2N(n)k^+}2^{-k}|\mathcal{I}_{m;k,k_1,k_2}|\lesssim \varep_1^32^{2H(q,n)\delta m}
\end{equation}
and, if $n\geq 1$,
\begin{equation}\label{tat1.1}
\begin{split}
\sum_{k,k_1,k_2\in\mathbb{Z}}2^{k_1^+}2^{2N(n)k_2^+}2^{-k}|\mathcal{I}_{m;k,k_1,k_2}|\lesssim \varep_1^32^{2H(q,n)\delta m}.
\end{split}
\end{equation}
These two bounds would clearly suffice to complete the proof of \eqref{KKFINt3}--\eqref{KKFINt1}.

Using just the $L^2$ bounds \eqref{vcx1} and \eqref{vcx1.2}, we have
\begin{equation}\label{wer10.2}
\begin{split}
2^{-k}|\mathcal{I}_{m;k,k_1,k_2}|\lesssim \varep_1^3|J_m|2^{2H(q,n)\delta m}&2^{\underline{k}+k_1^-}2^{(\underline{k}-k)/2} \cdot 2^{-\ga(m+k^-)}2^{\kappa k_1^-}\\
&\times 2^{-N(n)k^+-N(n)k_2^+}2^{-(N_0-2)k_1^+}.
\end{split}
\end{equation} 
Using also the $L^\infty$ bounds \eqref{wws2} on the $P_{k_1}U^{\psi,\iota_1}$ component we also have
\begin{equation}\label{wer10.3}
\begin{split}
2^{-k}|\mathcal{I}_{m;k,k_1,k_2}|\lesssim \varep_1^3|J_m|2^{2H(q,n)\delta m}&2^{-k/2}2^{-m+\delta'm/4}2^{k_1^-/2}2^{-\ga(m+k^-)}\\
&\times 2^{-N(n)k^+-N(n)k_2^+}2^{-(N(1)-2)k_1^+}.
\end{split}
\end{equation}

The bounds \eqref{tat1}--\eqref{tat1.1} follow if $|J_m|\lesssim 1$. Indeed, the bounds \eqref{wer10.2} suffice to estimate the contribution of the triplets $(k,k_1,k_2)$ with $\vert k\vert+\vert k_1\vert+\vert k_2\vert\ge\delta^\prime m$ (in the case $n=0$ we use also a similar bound with the roles of $k_1$ and $k_2$ reversed). On the other hand, the contribution of the triplets $(k,k_1,k_2)$ with $\vert k\vert+\vert k_1\vert+\vert k_2\vert\le\delta^\prime m$ can be bounded using \eqref{wer10.3}.

On the other hand, if $|J_m|\approx 2^m\gg 1$ (thus $m\in[\delta^{-1},L]$) the bounds \eqref{wer10.2} still suffice to bound the contribution of triplets $(k,k_1,k_2)$ for which either $\underline{k}\leq -m$ or $\max\{k,k_1,k_2\}\geq 4m$. For the remaining contributions, we consider several cases.

{\bf{Step 3.}} We show first that if $m\in[1/\delta,L]$ then
\begin{equation}\label{njm1}
\begin{split}
\sum_{k,k_1,k_2\in\mathbb{Z},\,\underline{k}\geq -m,\,k\leq -0.6m}2^{2\ga (m+k)}2^{-k}|\mathcal{I}_{m;k,k_1,k_2}|&\lesssim \varep_1^32^{2H(q,n)\delta m},\\
\sum_{\substack{k,k_1,k_2\in\mathbb{Z},\,\underline{k}\geq -m,\, k\leq -0.6m}}2^{-k}2^{2N(n)k_2^+}2^{k_1^+}|\mathcal{I}_{m;k,k_1,k_2}|&\lesssim \varep_1^32^{2H(q,n)\delta m},\qquad n\ge1.
\end{split}
\end{equation}
This is the case of small frequencies $2^k$. The estimates \eqref{wer10.2} cleary suffice to control the contribution of the triplets $(k,k_1,k_2)$ for which $k\leq -0.6m$ and $k_1\le -0.4m$. They also suffice to control the contribution of the triplets $(k,k_1,k_2)$ for which $k\leq-0.6m$ and $(1+\ga)(m+k)+k_1^-(1+\kappa)-6.5k_1^+\leq 0$.

It remains to bound the contribution of the triplets $(k,k_1,k_2)$ for which 
\begin{equation}\label{njm2}
k\in[-m,-0.6m]\qquad\text{ and }\qquad (m+k)+k_1^-\geq 6k_1^+.
\end{equation}
In particular, $k_1\geq -m/2+100$. Let $J_1:=k_1^-+m-40$ and decompose $P_{k_1}U^{\psi,\iota_1}=U^{\psi,\iota_1}_{\leq J_1,k_1}+U^{\psi,\iota_1}_{>J_1,k_1}$ as in \eqref{on11.3}--\eqref{on11.36}. Let
\begin{equation}\label{njm2.2}
\begin{split}
&\mathcal{I}^1_{m;k,k_1,k_2}:=\int_{J_m}q_m(s)\mathcal{G}_\mathfrak{m}[U^{\psi,\iota_1}_{\leq J_1,k_1}(s),P_{k_2}U^{\mathcal{L}_2\psi,\iota_2}(s),P_kU^{\mathcal{L}h,\iota}(s)]\,ds,\\
&\mathcal{I}^2_{m;k,k_1,k_2}:=\int_{J_m}q_m(s)\mathcal{G}_\mathfrak{m}[U^{\psi,\iota_1}_{>J_1,k_1}(s),P_{k_2}U^{\mathcal{L}_2\psi,\iota_2}(s),P_kU^{\mathcal{L}h,\iota}(s)]\,ds.
\end{split}
\end{equation}  

Using \eqref{vcx1} and \eqref{wws13x} we estimate
\begin{equation*}
\begin{split}
2^{-k}|\mathcal{I}^1_{m;k,k_1,k_2}|&\lesssim 2^m2^{-k}\sup_{s\in J_m}\|P_kU^{\mathcal{L}h,\iota}(s)\|_{L^2}\|U^{\psi,\iota_1}_{\leq J_1,k_1}(s)\|_{L^\infty}\|P_{k_2}U^{\mathcal{L}_2\psi,\iota_2}(s)\|_{L^2}\\
&\lesssim\varep_1^32^{2H(q,n)\delta m}(2^m2^k2^{k_1^-})^{-1/2}2^{\kappa k_1^-/20}2^{-\ga(m+k)}2^{-N(n)k_2^+}2^{-(N_0-5)k_1^+}.
\end{split}
\end{equation*}
Therefore, for $(k,k_1,k_2)$ as in \eqref{njm2},
\begin{equation}\label{njm4}
\begin{split}
2^{-k}2^{2\ga(m+k)}|\mathcal{I}^1_{m;k,k_1,k_2}|&\lesssim \varep_1^32^{2H(q,n)\delta m}(2^m2^k2^{k_1^-})^{-1/2+\ga}2^{-\delta'\vert k_1\vert},\\
2^{-k}2^{2N(n)k_2^+}2^{k_1^+}|\mathcal{I}^1_{m;k,k_1,k_2}|
&\lesssim \varep_1^32^{2H(q,n)\delta m}(2^m2^k2^{k_1^-})^{-1/2}2^{-\ga (m+k)}2^{N(n)k_2^+-(N_0-6)k_1^+}.
\end{split}
\end{equation}

Similarly, using \eqref{vcx1}, \eqref{vcx1.1}, and $L^2$ bounds we estimate
\begin{equation*}
\begin{split}
2^{-k}|\mathcal{I}^2_{m;k,k_1,k_2}|&\lesssim 2^m2^{k/2}\sup_{s\in I_m}\|P_kU^{\mathcal{L}h,\iota}(s)\|_{L^2}\|U^{\psi,\iota_1}_{>J_1,k_1}(s)\|_{L^2}\|P_{k_2}U^{\mathcal{L}_2\psi,\iota_2}(s)\|_{L^2}\\
&\lesssim\varep_1^32^{ 2H(q,n)\delta m}2^{-\ga(m+k)}2^{30\delta m}2^{k-k_1^-}2^{-N(n)k_2^+}2^{-N(1)k_1^+}.
\end{split}
\end{equation*}
Therefore, for $(k,k_1,k_2)$ as in \eqref{njm2},
\begin{equation}\label{njm5}
\begin{split}
2^{-k}2^{2\ga(m+k)}|\mathcal{I}^2_{m;k,k_1,k_2}|&\lesssim \varep_1^32^{2H(q,n)\delta m}2^{\ga(m+k)}2^{30\delta m}2^{k-k_1^-}2^{-k_1^+},\\
2^{-k}2^{2N(n)k_2^+}2^{k_1^+}|\mathcal{I}^2_{m;k,k_1,k_2}|&\lesssim \varep_1^32^{2H(q,n)\delta m}2^{-\ga(m+k)}2^{30\delta m}2^{k-k_1^-}2^{N(n)k_2^+}2^{-N(1)k_1^++k_1^+}.
\end{split}
\end{equation} 
It is easy to see that \eqref{njm4}-\eqref{njm5} suffice to bound the remaining contribution of the triplets $(k,k_1,k_2)$ as in \eqref{njm2}. The desired estimates \eqref{njm1} follow.

{\bf{Step 4.}} With $\overline{k}=\max\{k,k_1,k_2\}$ we show now that if $m\in[1/\delta,L]$ then
\begin{equation}\label{njm6}
\begin{split}
\sum_{k,k_1,k_2\in\mathbb{Z},\,k\geq -0.6m,\,\overline{k}\leq 8\delta'm}2^{2\ga(m+k^-)}2^{2N(n)k^+}2^{-k}|\mathcal{I}_{m;k,k_1,k_2}|&\lesssim \varep_1^32^{2H(q,n)\delta m},\\
\sum_{k,k_1,k_2\in\mathbb{Z},\,k\geq -0.6m,\,\overline{k}\leq 8\delta'm}2^{2N(n)k_2^++k_1^+}2^{-k}|\mathcal{I}_{m;k,k_1,k_2}|&\lesssim \varep_1^32^{2H(q,n)\delta m}.
\end{split}
\end{equation}
For this we use normal forms, see \eqref{NF2}--\eqref{NF3}. Using also Lemma \ref{pha2} (ii) we estimate
\begin{equation*}
\begin{split}
2^{-k}|\mathcal{I}_{m;k,k_1,k_2}|\lesssim 2^{-2k}&2^{4\overline{k}^+}\sup_{s\in J_m}\big\{\|P_{k_1}U^{\psi}(s)\|_{L^\infty}\|P_{k_2}U^{\LL_2 \psi}(s)\|_{L^2}\|P_kU^{\mathcal{L}h}(s)\|_{L^2}\\
&+2^m\|P_{k_1}(\partial_s+i\Lambda_{kg})U^{\psi}(s)\|_{L^\infty}\|P_{k_2}U^{\LL_2 \psi}(s)\|_{L^2}\|P_kU^{\mathcal{L}h}(s)\|_{L^2}\\
&+2^m\|P_{k_1}U^{\psi}(s)\|_{L^\infty}\|P_{k_2}(\partial_s+i\Lambda_{kg})U^{\mathcal{L}_2\psi}(s)\|_{L^2}\|P_kU^{\mathcal{L}h}(s)\|_{L^2}\\
&+2^m\|P_{k_1}U^{\psi}(s)\|_{L^\infty}\|P_{k_2}U^{\LL_2\psi}(s)\|_{L^2}\|P_k(\partial_s+i\Lambda_{wa})U^{\mathcal{L}h}(s)\|_{L^2}\big\}.
\end{split}
\end{equation*} 
Using the $L^2$ bounds \eqref{mnb60.1}--\eqref{mnb60.2} and the $L^\infty$ bounds \eqref{wws2} and \eqref{plk3} we then estimate
\begin{equation*}
\begin{split}
2^{-k}|\mathcal{I}_{m;k,k_1,k_2}|\lesssim \varep_1^32^{-|k_1|/4}2^{-3k/2}2^{-m+8\delta'm}2^{-4\overline{k}^+},
\end{split}
\end{equation*} 
and the bounds \eqref{njm6} follow.

{\bf{Step 5.}} We now consider the case of high frequencies, and show that if $m\in[1/\delta,L]$ then
\begin{equation}\label{njm12}
\begin{split}
\sum_{k,k_1,k_2\in\mathbb{Z},\,k\geq -0.6m,\,\overline{k}\geq 8\delta' m}2^{2N(n)k^+}2^{2\ga(m+k^-)}2^{-k}|\mathcal{I}_{m;k,k_1,k_2}|&\lesssim \varep_1^32^{2H(n)\delta m},\\
\sum_{k,k_1,k_2\in\mathbb{Z},\,k\geq -0.6m,\,\overline{k}\geq 8\delta' m}2^{2N(n)k_2^++k_1^+}2^{-k}\vert \mathcal{I}_{m;k,k_1,k_2}|&\lesssim \varep_1^32^{2H(n)\delta m},\qquad n\geq 1.
\end{split}
\end{equation}

Assume first that $n\leq 2$. Using \eqref{vcx1}, \eqref{wws1}, and \eqref{wws2}, with the lowest frequency placed in $L^\infty$, for triplets $(k,k_1,k_2)$ as in \eqref{njm12} we have
\begin{equation*}
\begin{split}
2^{-k}|\mathcal{I}_{m;k,k_1,k_2}|&\lesssim\varep_1^32^{-m+2\delta'm}2^{-N(n)k_2^+}2^{-N(0)k_1^+}\qquad\text{ if }k=\underline{k},\\
2^{-k}|\mathcal{I}_{m;k,k_1,k_2}|&\lesssim\varep_1^32^{-m+2\delta'm}2^{k^-_1/2-4k_1^+}2^{-N(n)k_2^+}2^{-N(n)k^+-k/2}\qquad\text{ if }k_1=\underline{k},\\
2^{-k}|\mathcal{I}_{m;k,k_1,k_2}|&\lesssim\varep_1^32^{-m+2\delta'm}2^{k^-_2/2-4k_2^+}2^{-N(n)k^+-k/2}2^{-N(0)k_1^+}\qquad\text{ if }k_2=\underline{k}.
\end{split}
\end{equation*}
These bounds suffice to prove \eqref{njm12}, due to the gain of high derivative in all cases.

Assume now that $n=3$. The bounds \eqref{wer10.3} suffice to bound the contribution of the triplets $(k,k_1,k_2)$ as in \eqref{njm12} for which $k\geq 0$. On the other hand, if $k\leq 0$ (thus $k=\underline{k}$, $k_1,k_2\geq 8\delta'm-6$), then we let $J_1=m-40$ and decompose $P_{k_1}U^{\psi,\iota_1}=U^{\psi,\iota_1}_{\leq J_1,k_1}+U^{\psi,\iota_1}_{>J_1,k_1}$ and $\mathcal{I}_{m;k,k_1,k_2}=\mathcal{I}_{m;k,k_1,k_2}+\mathcal{I}_{m;k,k_1,k_2}$ as in \eqref{njm2.2}. Then we estimate
\begin{equation*}
\begin{split}
2^{-k}|\mathcal{I}^1_{m;k,k_1,k_2}|&\lesssim 2^m2^{-k}\sup_{s\in J_m}\|P_kU^{\mathcal{L}h,\iota}(s)\|_{L^2}\|U^{\psi,\iota_1}_{\leq J_1,k_1}(s)\|_{L^\infty}\|P_{k_2}U^{\mathcal{L}_2\psi,\iota_2}(s)\|_{L^2}\\
&\lesssim\varep_1^32^{2\delta' m}2^{-k/2}2^{-m/2}2^{-N(n)k_2^+}2^{-(N_0-5)k_1^+},
\end{split}
\end{equation*}
using \eqref{vcx1} and \eqref{wws13x}. Also
\begin{equation*}
\begin{split}
2^{-k}|\mathcal{I}^2_{m;k,k_1,k_2}|&\lesssim 2^m2^{k/2}\sup_{s\in I_m}\|P_kU^{\mathcal{L}h,\iota}(s)\|_{L^2}\|U^{\psi,\iota_1}_{>J_1,k_1}(s)\|_{L^2}\|P_{k_2}U^{\mathcal{L}_2\psi,\iota_2}(s)\|_{L^2}\\
&\lesssim\varep_1^32^{2\delta' m}2^{k}2^{-N(n)k_2^+}2^{-N(1)k_1^+},
\end{split}
\end{equation*}
using \eqref{vcx1}--\eqref{vcx1.1}. Therefore, if $k\in[-0.6m,0]$ we have
\begin{equation*}
2^{-k}|\mathcal{I}_{m;k,k_1,k_2}|\lesssim\varep_1^32^{2\delta' m}2^{-N(n)k_2^+}2^{-N(1)k_1^+}.
\end{equation*}
This suffices to bound the remaining contributions over $k\leq 0$ in \eqref{njm12}.
\end{proof}

We can now finally complete the proof of Proposition \ref{totbounds}.

\begin{lemma}\label{Pag1}
With the assumptions of Proposition \ref{totbounds}, for any $m\in\{0,\ldots,L+1\}$ we have
\begin{equation}\label{Pag2}
\begin{split}
\sum_{k,k_1,k_2\in\mathbb{Z}}&2^{N(n)k^+}(2^{N(n)k_2^+}+2^{k_2^+-k_1^+}2^{N(n)k_1^+})(2^{-k^+}+2^{-k_2^+})\Big|\int_{J_m}q_m(s)\\
&\times\mathcal{G}_\mathfrak{m}[P_{k_1}U^{h_1,\iota_1}(s),P_{k_2}U^{\mathcal{L}_2\psi,\iota_2}(s),P_kU^{\mathcal{L}\psi,\iota}(s)]\,ds\Big|\lesssim \varepsilon_1^32^{2H(q,n)\delta m},
\end{split}
\end{equation}
where $\mathcal{L},\mathcal{L}_2\in\mathcal{V}_n^q$, $n\leq 3$, and $\mathfrak{m}\in\mathcal{M}^\ast$. Moreover, if $\mathfrak{n}_{\iota_1\iota_2}\in\mathcal{M}^0_{\iota_1\iota_2}$ is a null multiplier then
\begin{equation}\label{Pag3}
\begin{split}
\sum_{k,k_1,k_2\in\mathbb{Z}}&2^{N(n)k^+}(2^{N(n)k_2^+}+2^{k_2^+-k_1^+}2^{N(n)k_1^+})\Big|\int_{J_m}q_m(s)\\
&\times\mathcal{G}_{\mathfrak{n}_{\iota_1\iota_2}}[P_{k_1}U^{h_1,\iota_1}(s),P_{k_2}U^{\mathcal{L}_2\psi,\iota_2}(s),P_kU^{\mathcal{L}\psi,\iota}(s)]\,ds\Big|\lesssim \varepsilon_1^32^{2H(q,n)\delta m},
\end{split}
\end{equation}
As a consequence, the bounds \eqref{EEcom2}--\eqref{BulkKG} hold.
\end{lemma}

\begin{proof} Using the $L^\infty$ bounds \eqref{wws1} and the $L^2$ bounds \eqref{vcx1} we have the general estimates
\begin{equation}\label{Pag4}
\begin{split}
\Big|\int_{J_m}q_m(s)&\mathcal{G}_\mathfrak{m}[P_{k_1}U^{h_1,\iota_1}(s),P_{k_2}U^{\mathcal{L}_2\psi,\iota_2}(s),P_kU^{\mathcal{L}\psi,\iota}(s)]\,ds\Big|\\
&\lesssim \varep_1^3|J_m|2^{2H(q,n)\delta m}2^{-m+30\delta m}2^{k_1^-}2^{-N(n)k^+-N(n)k_2^+}2^{-(N(1)-2)k_1^+}.
\end{split}
\end{equation}
As before, we divide the proof into several steps.

{\bf{Step 1.}} We consider first the contribution of large frequencies $k_1$, and show that
\begin{equation}\label{Pag5}
\begin{split}
\sum_{k,k_1,k_2\in\mathbb{Z},\,k_1\geq\overline{k}-8}&2^{N(n)k^+}(2^{N(n)k_2^+}+2^{k_2^+-k_1^+}2^{N(n)k_1^+})\Big|\int_{J_m}q_m(s)\\
&\times\mathcal{G}_\mathfrak{m}[P_{k_1}U^{h_1,\iota_1}(s),P_{k_2}U^{\mathcal{L}_2\psi,\iota_2}(s),P_kU^{\mathcal{L}\psi,\iota}(s)]\,ds\Big|\lesssim \varepsilon_1^32^{2H(q,n)\delta m},
\end{split}
\end{equation}
where $\overline{k}=\max\{k,k_1,k_2\}$. Clearly, $2^{N(n)k_2^+}\lesssim 2^{k_2^+-k_1^+}2^{N(n)k_1^+}$ for triplets $(k,k_1,k_2)$ as in \eqref{Pag5}, and by symmetry we may assume $k_2\leq k$. 

If $n\leq 2$ then we can use the $L^\infty$ bounds \eqref{wws2} and the $L^2$ bounds \eqref{vcx1} to estimate
\begin{equation}\label{Pag6}
\begin{split}
\Big|\int_{J_m}q_m(s)&\mathcal{G}_\mathfrak{m}[P_{k_1}U^{h_1,\iota_1}(s),P_{k_2}U^{\mathcal{L}_2\psi,\iota_2}(s),P_kU^{\mathcal{L}\psi,\iota}(s)]\,ds\Big|\\
&\lesssim \varep_1^3|J_m|2^{2H(q,n)\delta m}2^{k_2^-/2}2^{-4k_2^+}2^{-m+\delta' m/2}2^{k_1^+/2}2^{-N(0)k_1^+-N(n)k^+}.
\end{split}
\end{equation}
The bounds \eqref{Pag4} (used for $n=3$) and \eqref{Pag6} (used for $n\leq 2$) show that the contribution of the triples $(k,k_1,k_2)$ in \eqref{Pag5} for which $\max\{|k|,|k_1|,|k_2|\}\geq 2\delta'm$ is bounded as claimed.

Assume now that $\max\{|k|,|k_1,|k_2|\}\leq 2\delta'm$. If $|J_m|\lesssim 1$ then the bounds \eqref{Pag4} and \eqref{Pag6} still suffice to control the contribution of these triples. On the other hand, if $m\in[\delta^{-1},L]$ then we use normal forms, see \eqref{NF2}--\eqref{NF3}, and Lemma \ref{pha2} (ii) to estimate
\begin{equation}\label{Pag8}
\begin{split}
\Big|\int_{J_m}&q_m(s)\mathcal{G}_\mathfrak{m}[P_{k_1}U^{h_1,\iota_1}(s),P_{k_2}U^{\mathcal{L}_2\psi,\iota_2}(s),P_kU^{\mathcal{L}\psi,\iota}(s)]\,ds\Big|\lesssim 2^{-k_1}2^{4\overline{k}^+}\sup_{s\in J_m}\\
&\big\{[\|P_{k_1}U^{h_1}(s)\|_{L^\infty}+2^m\|P_{k_1}(\partial_s+i\Lambda_{wa})U^{h_1}(s)\|_{L^\infty}]\|P_{k_2}U^{\LL_2 \psi}(s)\|_{L^2}\|P_kU^{\mathcal{L}\psi}(s)\|_{L^2}\\
&+2^m\|P_{k_1}U^{h_1}(s)\|_{L^\infty}\|P_{k_2}(\partial_s+i\Lambda_{kg})U^{\mathcal{L}_2\psi}(s)\|_{L^2}\|P_kU^{\mathcal{L}\psi}(s)\|_{L^2}\\
&+2^m\|P_{k_1}U^{h_1}(s)\|_{L^\infty}\|P_{k_2}U^{\LL_2\psi}(s)\|_{L^2}\|P_k(\partial_s+i\Lambda_{kg})U^{\mathcal{L}\psi}(s)\|_{L^2}\big\}.
\end{split}
\end{equation}
In view of the $L^2$ bounds \eqref{vcx1} and \eqref{wer4.1} and the $L^\infty$ bounds \eqref{wws1} and \eqref{plk2}, all the terms in the right-hand side of \eqref{Pag8} are bounded by $C2^{-m+8\delta'm}$. This suffices to bound the remaining contributions in \eqref{Pag5}.

{\bf{Step 2.}} We complete now the proof of \eqref{Pag2} by showing that
\begin{equation}\label{Pag10}
\begin{split}
\sum_{k,k_1,k_2\in\mathbb{Z},\,k_1\leq\overline{k}-8}&2^{N(n)k^+}2^{N(n)k_2^+}(2^{-k^+}+2^{-k_2^+})\Big|\int_{J_m}q_m(s)\\
&\times\mathcal{G}_\mathfrak{m}[P_{k_1}U^{h_1,\iota_1}(s),P_{k_2}U^{\mathcal{L}_2\psi,\iota_2}(s),P_kU^{\mathcal{L}\psi,\iota}(s)]\,ds\Big|\lesssim \varepsilon_1^32^{2H(q,n)\delta m}.
\end{split}
\end{equation}

Indeed, we may assume that $2^k\approx 2^{k_2}$ and use \eqref{Pag4} to bound the contribution of the triples $(k,k_1,k_2)$ in \eqref{Pag10} for which $\max\{|k|,|k_1,|k_2|\}\geq 2\delta'm$. For the remaining triples $(k,k_1|,k_2)$ with $\max\{|k|,|k_1|,|k_2|\}\leq 2\delta'm$ we use normal forms and estimate as in \eqref{Pag8}.

{\bf{Step 3.}} Finally, we complete the proof of \eqref{Pag3} by showing that
\begin{equation}\label{Pag15}
\begin{split}
\sum_{k,k_1,k_2\in\mathbb{Z},\,k_1\leq\overline{k}-8}&2^{N(n)k^+}2^{N(n)k_2^+}\Big|\int_{J_m}q_m(s)\\
&\times\mathcal{G}_{\mathfrak{n}_{\iota_1\iota_2}}[P_{k_1}U^{h_1,\iota_1}(s),P_{k_2}U^{\mathcal{L}_2\psi,\iota_2}(s),P_kU^{\mathcal{L}\psi,\iota}(s)]\,ds\Big|\lesssim \varepsilon_1^32^{2H(q,n)\delta m}.
\end{split}
\end{equation}

This is more subtle, essentially due to the possible loss of high derivative, and we need to exploit the null structure of the symbol $\mathfrak{n}_{\iota_1\iota_2}$. We start by estimating some of the easier contributions. Recall the coefficients $b_k(s)=b_k(q,n;s)$ defined in \eqref{nake15} and the time averages $b_{k,m}=b_{k,m}(q,n)$ defined in \eqref{bnm21}. Using \eqref{wws1} and \eqref{nake15.2} we have
\begin{equation}\label{Pag16}
\begin{split}
\Big|\int_{J_m}q_m(s)&\mathcal{G}_{\mathfrak{n}_{\iota_1\iota_2}}[P_{k_1}U^{h_1,\iota_1}(s),P_{k_2}U^{\mathcal{L}_2\psi,\iota_2}(s),P_kU^{\mathcal{L}\psi,\iota}(s)]\,ds\Big|\\
&\lesssim \varep_1|J_m|2^{2H(q,n)\delta m}2^{-m+30\delta m}2^{k_1^-}2^{-4k_1^+}2^{-N(n)k^+-N(n)k_2^+}b_{k,m}^2
\end{split}
\end{equation}
if $|k-k_2|\leq 8$. In view of \eqref{nake15.1} this suffices to control the contribution of the triplets $(k,k_1,k_2)$ in \eqref{Pag15} for which $|k_1|\geq\delta'm$.  

Assume that $|k_1|\leq\delta'm$. The bounds \eqref{Pag16} still suffice to control the sum if $|J_m|\lesssim 1$.  On the other hand, if $m\in[\delta^{-1},L]$ then we can use normal forms as in \eqref{Pag8} to control the contribution of the triplets $(k,k_1,k_2)$ for which $\max\{k,k_2\}\leq 20\delta'm+8$. 

After these reductions, we may assume that 
\begin{equation}\label{Pag18}
m\in[\delta^{-1},L],\qquad |k_1|\leq\delta'm,\qquad k,k_2\geq 20\delta'm.
\end{equation}
We can further dispose of the resonant part of the multipliers. As in \eqref{trf2} we decompose
\begin{equation}\label{Pag19}
\mathfrak{n}_{\iota_1\iota_2}=\mathfrak{n}_{\iota_1\iota_2}^r+\mathfrak{n}^{nr}_{\iota_1\iota_2},\qquad\mathfrak{n}_{\iota_1\iota_2}^r(\theta,\eta)=\varphi_{\le q_0}(\Xi_{\iota_1\iota_2}(\theta,\eta))\mathfrak{n}_{\iota_1\iota_2}(\theta,\eta),
\end{equation}
where $q_0:=-2\delta'm$ and $\Xi_{\iota_1\iota_2}$ is defined as in \eqref{par1}. In view of \eqref{yip2},
\begin{equation}\label{Pag20}
\begin{split}
\Big|\int_{J_m}q_m(s)&\mathcal{G}_{\mathfrak{n}_{\iota_1\iota_2}^r}[P_{k_1}U^{h_1,\iota_1}(s),P_{k_2}U^{\mathcal{L}_2\psi,\iota_2}(s),P_kU^{\mathcal{L}\psi,\iota}(s)]\,ds\Big|\\
&\lesssim \varep_12^{q_0}2^{2H(q,n)\delta m}2^{-m+30\delta m}2^{k_1^-}2^{-4k_1^+}2^{-N(n)k^+-N(n)k_2^+}b_{k,m}^2,
\end{split}
\end{equation}
for $k,k_1,k_2$ as in \eqref{Pag18}. The key factor $2^{q_0}$ in the right-hand side is due to the nullness assumption on the multiplier $\mathfrak{n}_{\iota_1\iota_2}$. The bounds \eqref{Pag20} suffice to control the contributions of the resonant interactions. To summarize, for \eqref{Pag15} it remains to show that
\begin{equation}\label{Pag22}
\Big|\int_{J_m}q_m(s)\mathcal{G}_{\mathfrak{n}_{\iota_1\iota_2}^{nr}}[P_{k_1}U^{h_1,\iota_1}(s),P_{k_2}U^{\mathcal{L}_2\psi,\iota_2}(s),P_kU^{\mathcal{L}\psi,\iota}(s)]\,ds\Big|\lesssim \varepsilon_12^{-2N(n)k^+}b_{k,m}^2,
\end{equation}
provided that $m,k,k_1,k_2$ satisfy \eqref{Pag18}.

{\bf{Step 4.}} The bounds \eqref{Pag22} are similar to the bounds \eqref{bnm25} in Lemma \ref{ParalinearSym}. Recall that 
\begin{equation}\label{Pag23}
\begin{split}
\|P_{k_2}U^{\mathcal{L}_2\psi,\iota_2}(s)\|_{L^2}+\|P_kU^{\mathcal{L}\psi,\iota}(s)\|_{L^2}\lesssim b_k(s)2^{\delta'm}2^{-N(n)k},\\
\|P_{k_1}U^{h_1,\iota_1}(s)\|_{L^\infty}\lesssim \varep_12^{-m+\delta'm}2^{k_1-4k_1^+},
\end{split}
\end{equation}
see \eqref{nake15.2}, and
\begin{equation}\label{Pag24}
\begin{split}
\|(\partial_s+i\Lambda_{kg,\iota_2})&P_{k_2}U^{\mathcal{L}_2\psi,\iota_2}(s)\|_{L^2}+\|(\partial_s+i\Lambda_{kg,\iota})P_kU^{\mathcal{L}\psi,\iota}(s)\|_{L^2}\lesssim \varep_1b_k(s)2^{-m+\delta'm}2^{-N(n)k+k}.
\end{split}
\end{equation}
These bounds follow easily using $L^2\times L^\infty$ estimates from the formulas \eqref{roco2} and \eqref{gb10}. 

To prove \eqref{Pag22} it is convenient for us to apply the wave operator $\Lambda_{wa}$ to the Klein-Gordon variables $U^{\mathcal{L}_2\psi}$ and $U^{\mathcal{L}\psi}$, instead of the more natural Klein-Gordon operator $\Lambda_{kg}$. The reason is to be able to reuse some of the estimates in Lemma \ref{ParalinearSym}, such as \eqref{pint2}, and also use some of the bounds in Lemma \ref{PhaWave} that do not involve derivative loss, such as \eqref{nag1}. We notice that the two operators are not too different at high frequency, 
\begin{equation}\label{Pag24.5}
\|(\Lambda_{wa}-\Lambda_{kg})P_{l}f\|_{L^2}\lesssim 2^{-l}\|P_lf\|_{L^2},\qquad l\in\{k,k_2\}.
\end{equation}

As in {\bf{Step 1}} in the proof of Lemma \ref{ParalinearSym}, the estimates \eqref{Pag22} follow easily if $\iota=-\iota_2$, using normal forms and the ellipticity of the phase $\Lambda_{wa,\iota}(\xi)-\Lambda_{wa,\iota_2}(\eta)-\Lambda_{\iota_1}(\xi-\eta)$. On the other hand, if $\iota=\iota_2$ then we may assume that $\iota=\iota=+$ and replace the variables $U^{\mathcal{L}_2\psi,\iota_2}$ and $U^{\mathcal{L}\psi,\iota}$ with the quasilinear variables $\mathcal{U}^{\mathcal{L}_2\psi}$ and $\mathcal{U}^{\mathcal{L}\psi}$, at the expense of acceptable errors (as bounded in \eqref{DefParaUnknown}). We then apply normal forms and notice that 
\begin{equation}\label{Pag25}
\Big\|\mathcal{F}^{-1}\Big\{\frac{\mathfrak{n}_{\iota_1+}^{nr}(\xi-\eta,\eta)}{|\xi|-|\eta|-\iota_1|\xi-\eta|}\varphi_{kk_1k_2}(\xi-\eta,\eta)\Big\}\Big\|_{L^1(\mathbb{R}^6)}\lesssim 2^{-k_1}2^{-3q_0}\lesssim 2^{8\delta'm},
\end{equation}
as a consequence of \eqref{nag1} and \eqref{par7} (notice that there is no high derivative loss in these bounds).  Therefore we can estimate the first two terms in \eqref{NF3} using just \eqref{Pag23}. 

To control the remaining two terms and prove \eqref{Pag22} it remains to show that
\begin{equation}\label{Pag26}
\begin{split}
&\Big|\mathcal{H}_{\mathfrak{n}^{nr}_{\iota_1+}}\big[P_{k_1}U^{h_1,\iota_1}(s),(\partial_s+i\Lambda_{wa})P_{k_2}\mathcal{U}^{\mathcal{L}_2\psi}(s),P_k\mathcal{U}^{\mathcal{L}\psi}(s)\big]\\
&+\mathcal{H}_{\mathfrak{n}^{nr}_{\iota_1+}}\big[P_{k_1}U^{h_1,\iota_1}(s),P_{k_2}\mathcal{U}^{\mathcal{L}_2\psi}(s),(\partial_s+i\Lambda_{wa})P_k\mathcal{U}^{\mathcal{L}\psi}(s)\big]\Big|\lesssim \varepsilon_12^{-2N(n)k}2^{-m} b_{k}(s)^2,
\end{split}
\end{equation}
for any $s\in J_m$, where the operators $\mathcal{H}_{\mathfrak{n}^{nr}_{\iota_1+}}$ are defined in \eqref{DefJ}. We decompose 
\begin{equation*}
(\partial_s+i\Lambda_{wa})P_k\mathcal{U}^{\mathcal{L}\psi}=(\partial_s+iT_{\Sigma_{kg}})P_k\mathcal{U}^{\mathcal{L}\psi}-iT_{\Sigma_{kg}-\Sigma_{wa}}P_k\mathcal{U}^{\mathcal{L}\psi}-iT_{\Sigma_{wa}^{\geq 1}}P_k\mathcal{U}^{\mathcal{L}\psi},
\end{equation*}
where $\Sigma_{wa},\Sigma_{kg}$ are defined in \eqref{ParaSymbolsW}-\eqref{ParaSymbolsKG} and $\Sigma^{\geq 1}_{wa}(x,\ze)=\Sigma_{wa}(x,\ze)-|\ze|$. Notice that $\|\Sigma_{kg}-\Sigma_{wa}\|_{\mathcal{L}^\infty_{-1}}\lesssim 1$, see the definition \eqref{nor1}. Therefore, using \eqref{LqBdTa} and \eqref{ParaBox5.02} we have
\begin{equation*}
\begin{split}
\|(\partial_s+iT_{\Sigma_{kg}})P_k\mathcal{U}^{\mathcal{L}\psi}(s)\|_{L^2}&\lesssim\varep_12^{-m+\delta'm}2^{-N(n)k}b_k(s),\\
\|T_{\Sigma_{kg}-\Sigma_{wa}}P_k\mathcal{U}^{\mathcal{L}\psi}\|_{L^2}&\lesssim 2^{\delta'm}2^{-N(n)k-k}b_k(s),
\end{split}
\end{equation*}
for any $s\in J_m$. Therefore, recalling that $k\geq 20\delta'm$ and the bounds \eqref{Pag24}, \eqref{Pag25}, these contributions to the second term in the left-hand side of \eqref{Pag26} can be bounded as claimed. Similarly, we can replace $(\partial_s+i\Lambda_{wa})P_{k_2}\mathcal{U}^{\mathcal{L}_2\psi}(s)$ with $-iT_{\Sigma_{wa}^{\geq 1}}P_{k_2}\mathcal{U}^{\mathcal{L}_2\psi}$ in the first term, at the expense of acceptable errors. For \eqref{Pag26} it remains to prove that, for any $s\in J_m$,
\begin{equation}\label{Pag27}
\begin{split}
&\Big|\mathcal{H}_{\mathfrak{n}^{nr}_{\iota_1+}}\big[P_{k_1}U^{h_1,\iota_1}(s),iT_{\Sigma^{\geq 1}_{wa}}P_{k_2}\mathcal{U}^{\mathcal{L}_2\psi}(s),P_k\mathcal{U}^{\mathcal{L}\psi}(s)\big]\\
&+\mathcal{H}_{\mathfrak{n}^{nr}_{\iota_1+}}\big[P_{k_1}U^{h_1,\iota_1}(s),P_{k_2}\mathcal{U}^{\mathcal{L}_2\psi}(s),iT_{\Sigma^{\geq 1}_{wa}}P_k\mathcal{U}^{\mathcal{L}\psi}(s)\big]\Big|\lesssim \varepsilon_12^{-2N(n)k}2^{-m} b_{k}(s)^2.
\end{split}
\end{equation}

The bounds \eqref{Pag27} are similar to the bounds \eqref{sori24} (with $J_1=\infty$). They follow by the same argument, first rewriting the expression as in \eqref{sori28}--\eqref{sori29}, and then using the general bounds \eqref{pint2}. This completes the proof of \eqref{Pag27}, and the desired bounds \eqref{Pag3} follow.

{\bf{Step 5.}} Finally we prove the bounds \eqref{EEcom2}--\eqref{BulkKG}. Notice that
\begin{equation}\label{Pag30}
\frac{\rho_j}{|\rho|}-\frac{\rho_j}{\langle\rho\rangle}=\frac{\rho_j}{|\rho|}\frac{1}{\langle\rho\rangle(|\rho|+\langle\rho\rangle)},
\end{equation}
for any $\rho\in\mathbb{R}^3$ and $j\in\{1,2,3\}$. We examine now the multipliers $p^{G,kg}_{\iota_1,\iota_2,\iota}$ in the formulas \eqref{symp5}--\eqref{symp8}, and compare them with the multipliers $p^{G,wa}_{\iota_1,\iota_2,\iota}$ in \eqref{symp1}--\eqref{symp4}. The multipliers $p^{G,wa}_{\iota_1,\iota_2,\iota}$ are null in the variables $\xi-\eta$ and $\eta$, as shown in the proof of Proposition \ref{EnergEst2} in section \ref{EnergyEstimates1}, and the difference between $p^{G,kg}_{\iota_1,\iota_2,\iota}-p^{G,wa}_{\iota_1,\iota_2,\iota}$ gains at least one derivative either in $\xi$ or in $\eta$, due to \eqref{Pag30}. The bounds \eqref{BulkKG} follow from \eqref{Pag2}--\eqref{Pag3}. 

To prove \eqref{EEcom2} we use identities similar to \eqref{sori1}--\eqref{sori1.7}, 
\begin{equation*}\label{Pag31}
\begin{split}
P_{kg}^n&[\widetilde{g}^{\mu\nu}_{\geq 1}\partial_\mu\partial_\nu \LL \psi+h_{00}\LL\psi]-[\widetilde{g}^{\mu\nu}_{\geq 1}\partial_\mu\partial_\nu (P_{kg}^n\LL \psi)+h_{00}P^n_{kg}\LL\psi]=\sum_{G\in\{F,\uF,\omega_n,\va_{mn}\}}\sum_{\iota_1,\iota_2\in\{+,-\}}\\
&\big\{P_{kg}^n I_{\mathfrak{q}^{G,kg}_{\iota_1\iota_2}}[|\nabla|^{-1}U^{G,\iota_1},\langle\nabla\rangle U^{\LL\psi,\iota_2}]-I_{\mathfrak{q}^{G,kg}_{\iota_1\iota_2}}[|\nabla|^{-1}U^{G,\iota_1},P^n_{kg}\langle\nabla\rangle U^{\LL\psi,\iota_2}]\big\}\\
&-2\{P_{kg}^n[R_jR_0\tau\cdot\partial_j\partial_0\LL \psi]-R_jR_0\tau\cdot\partial_j\partial_0(P^n_{kg}\LL \psi)\}\\
&+\sum_{(\mu,\nu)\neq (0,0)}\{P^n_{kg}[\widetilde{G}^{\mu\nu}_{\geq 2}\cdot \partial_\mu\partial_\nu \LL \psi]-\widetilde{G}^{\mu\nu}_{\geq 2}\cdot \partial_\mu\partial_\nu (P^n_{kg}\LL\psi)\},
\end{split}
\end{equation*}
where the null multipliers $\mathfrak{q}^{G,kg}_{\iota_1\iota_2}$ are defined in \eqref{BoxWave7.5} and, as in \eqref{sori1},  $\widetilde{G}^{\mu\nu}_{\geq 2}$ are linear combinations of expressions of the form $R^a|\nabla|^{-1}(G_{\geq 1}\partial_\rho h)$. The contributions of the cubic and higher order terms in the last two lines can be bounded as in \eqref{sori2}. For the main terms we decompose
\begin{equation}\label{Pag32}
\begin{split}
&P_{kg}^n I_{\mathfrak{q}^{G,kg}_{\iota_1\iota_2}}[|\nabla|^{-1}U^{G,\iota_1},\langle\nabla\rangle U^{\LL\psi,\iota_2}]-I_{\mathfrak{q}^{G,kg}_{\iota_1\iota_2}}[|\nabla|^{-1}U^{G,\iota_1},P^n_{kg}\langle\nabla\rangle U^{\LL\psi,\iota_2}]\\
&=P_{kg}^n I_{\mathfrak{m}^{1}_{\iota_1\iota_2}}[U^{G,\iota_1},U^{\LL\psi,\iota_2}]+P_{kg}^n I_{\mathfrak{m}^{2}_{\iota_1\iota_2}}[U^{G,\iota_1},U^{\LL\psi,\iota_2}]+I_{\mathfrak{m}^{3}_{\iota_1\iota_2}}[U^{G,\iota_1},P_{kg}^n U^{\LL\psi,\iota_2}],
\end{split}
\end{equation}
where
\begin{equation}\label{Pag33}
\begin{split}
&\mathfrak{m}^{1}_{\iota_1\iota_2}(\theta,\eta):=\varphi_{\leq 0}(\eta)\mathfrak{q}^{G,kg}_{\iota_1\iota_2}(\theta,\eta)\frac{P^n_{kg}(\eta+\theta)-P^n_{kg}(\eta)}{P^n_{kg}(\eta+\theta)}\frac{\langle\eta\rangle}{|\theta|},\\
&\mathfrak{m}^{2}_{\iota_1\iota_2}(\theta,\eta):=\varphi_{>1}(\eta)\mathfrak{q}^{G,kg}_{\iota_1\iota_2}(\theta,\eta)\Big\{\varphi_{\leq -4}(|\theta|/|\eta|)\frac{P^n_{kg}(\eta+\theta)-P^n_{kg}(\eta)}{P^n_{kg}(\eta+\theta)}\frac{\langle\eta\rangle}{|\theta|}+\varphi_{>-4}(|\theta|/|\eta|)\frac{\langle\eta\rangle}{|\theta|}\Big\},\\
&\mathfrak{m}^{3}_{\iota_1\iota_2}(\theta,\eta):=-\varphi_{>1}(\eta)\mathfrak{q}^{G,kg}_{\iota_1\iota_2}(\theta,\eta)\varphi_{>-4}(|\theta|/|\eta|)\frac{\langle\eta\rangle}{|\theta|}.
\end{split}
\end{equation}

It is easy to see that, for $a\in\{1,2,3\}$
\begin{equation*}
\|\mathcal{F}^{-1}\{\mathfrak{m}^{a}_{\iota_1\iota_2}(\theta,\eta)\varphi_{kk_1k_2}(\theta,\eta)\}\|_{L^1(\mathbb{R}^6)}\lesssim \min\{1,2^{k_2^+-k_1^+}\}\qquad\text{ for any }k,k_1,k_2\in\mathbb{Z}.
\end{equation*}
Therefore the contribution of the multiplier $\mathfrak{m}^{1}_{\iota_1\iota_2}$ can be controlled using just \eqref{Pag2}. To bound the contributions of the multipliers $\mathfrak{m}^{2}_{\iota_1\iota_2}$ and $\mathfrak{m}^{3}_{\iota_1\iota_2}$ we substitute first the symbols $\mathfrak{q}^{G,kg}_{\iota_1\iota_2}$ with $\mathfrak{q}^{G,wa}_{\iota_1\iota_2}$, and gain a factor of $\langle\eta\rangle^{-1}$ (compare \eqref{BoxWave7.5} and \eqref{BoxWave7}, and use \eqref{Pag30}). After this substitution we can thus use \eqref{Pag2} to estimate the contributions of the smoothing components, and \eqref{Pag3} for the contributions of the null components associated to the null symbols  $\mathfrak{q}^{G,wa}_{\iota_1\iota_2}$. This completes the proof of the bounds \eqref{EEcom2}--\eqref{BulkKG}.
\end{proof}

\chapter{Improved profile bounds}

\section{Weighted bounds} \label{LocProf}

In this section we prove the main bounds \eqref{bootstrap3.2}. Recall the identities \eqref{roco2}.
We also need a key identity that connects the vector-fields $\Gamma_l$ with weighted norms on profiles.

\begin{lemma}\label{ident}
Assume $\mu\in\{wa,kg\}$ and
\begin{equation}\label{wer0}
(\partial_t+i\Lambda_\mu)U=\mathcal{N},
\end{equation}
on $\mathbb{R}^3\times[0,T]$. If $V(t)=e^{it\Lambda_\mu}U(t)$ and $l\in\{1,2,3\}$ then, for any $t\in[0,T]$,
\begin{equation}\label{wer1}
\widehat{\Gamma_lU}(\xi,t)=i(\partial_{\xi_l}\widehat{\mathcal{N}})(\xi,t)+e^{-it\Lambda_{\mu}(\xi)}\partial_{\xi_l}[\Lambda_{\mu}(\xi)\widehat{V}(\xi,t)].
\end{equation}
\end{lemma}

\begin{proof} We calculate
\begin{equation*}
\begin{split}
&\widehat{\Gamma_lU}(\xi,t)=\mathcal{F}\{x_l\partial_tU+t\partial_lU\}(\xi,t)\\
&=i(\partial_{\xi_l}\widehat{\mathcal{N}})(\xi,t)+\partial_{\xi_l}[\Lambda_{\mu}(\xi)\widehat{U}(\xi,t)]+it\xi_l\widehat{U}(\xi,t)\\
&=i(\partial_{\xi_l}\widehat{\mathcal{N}})(\xi,t)+e^{-it\Lambda_{\mu}(\xi)}\partial_{\xi_l}[\Lambda_{\mu}(\xi)\widehat{V}(\xi,t)]-it(\partial_{\xi_l}\Lambda_{\mu})(\xi)e^{-it\Lambda_{\mu}(\xi)}\Lambda_{\mu}(\xi)\widehat{V}(\xi,t)+it\xi_l\widehat{U}(\xi,t).
\end{split}
\end{equation*}
This gives \eqref{wer1} since $(\partial_{\xi_l}\Lambda_{\mu})(\xi)\Lambda_\mu(\xi)=\xi_l$.
\end{proof}

Our main result in this section is the following proposition:

\begin{proposition}\label{DiEs1}
With the hypothesis in Proposition \ref{bootstrap}, for any $t\in[0,T]$, $\mathcal{L}\in\mathcal{V}_n^q$, $n\leq 2$, $l\in\{1,2,3\}$, and $k\in\mathbb{Z}$ we have
\begin{equation}\label{wer2}
2^{k/2}(2^{k^-}\langle t\rangle)^{\ga}\|P_k(x_lV^{\mathcal{L}h_{\al\be}})(t)\|_{L^2}+2^{k^+}\|P_k(x_lV^{\mathcal{L}\psi})(t)\|_{L^2}\lesssim\varep_0\langle t\rangle^{H(q+1,n+1)\delta}2^{-N(n+1)k^+}.
\end{equation}
\end{proposition}

\begin{proof} The identities \eqref{wer1} and \eqref{roco2} for $U^{\mathcal{L}h_{\al\be}}$ give
\begin{equation*}
\widehat{\Gamma_lU^{\mathcal{L}h_{\al\be}}}(\xi,t)=i(\partial_{\xi_l}\widehat{\LL\mathcal{N}^h_{\al\be}})(\xi,t)+e^{-it\Lambda_{wa}(\xi)}\partial_{\xi_l}[\Lambda_{wa}(\xi)\widehat{V^{\mathcal{L}h_{\al\be}}}(\xi,t)],
\end{equation*}
for $l\in\{1,2,3\}$. Therefore
\begin{equation*}
e^{-it\Lambda_{wa}(\xi)}\Lambda_{wa}(\xi)(\partial_{\xi_l}\widehat{V^{\mathcal{L}h_{\al\be}}})(\xi)=\widehat{\Gamma_lU^{\mathcal{L}h_{\al\be}}}(\xi)-i(\partial_{\xi_l}\widehat{\LL\mathcal{N}^h_{\al\be}})(\xi)-e^{-it\Lambda_{wa}(\xi)}(\xi_l/|\xi|)\widehat{V^{\mathcal{L}h_{\al\be}}}(\xi).
\end{equation*}
We multiply all the terms by $2^{-k/2}(2^{k^-}\langle t\rangle)^{\ga}\varphi_k(\xi)$ and take $L^2$ norms to show that
\begin{equation}\label{wer3}
\begin{split}
2^{k/2}&(2^{k^-}\langle t\rangle)^{\ga}\|\varphi_k(\xi)(\partial_{\xi_l}\widehat{V^{\mathcal{L}h_{\al\be}}})(\xi)\|_{L^2_\xi}\lesssim 2^{-k/2}(2^{k^-}\langle t\rangle)^{\ga}\|\varphi_k(\xi)\widehat{\Gamma_lU^{\mathcal{L}h_{\al\be}}}(\xi)\|_{L^2_\xi}\\
&+2^{-k/2}(2^{k^-}\langle t\rangle)^{\ga}\|\varphi_k(\xi)(\partial_{\xi_l}\widehat{\mathcal{N}^{\mathcal{L}h_{\al\be}}})(\xi)\|_{L^2_\xi}+2^{-k/2}(2^{k^-}\langle t\rangle)^{\ga}\|\varphi_k(\xi)\widehat{V^{\mathcal{L}h_{\al\be}}}(\xi)\|_{L^2_\xi}.
\end{split}
\end{equation}
To control the first term in the left-hand side of \eqref{wer2} it suffices to prove that
\begin{equation}\label{wer3.1}
2^{-k/2}(2^{k^-}\langle t\rangle)^{\ga}\|\varphi_k(\xi)(\partial_{\xi_l}\widehat{\LL\mathcal{N}^h_{\al\be}})(\xi)\|_{L^2_\xi}\lesssim\varep_1^2\langle t\rangle^{H(q+1,n+1)\delta}2^{-N(n+1)k^+}
\end{equation}
and 
\begin{equation}\label{wer3.2}
\begin{split}
2^{-k/2}(2^{k^-}\langle t\rangle)^{\ga}\|\varphi_k(\xi)\widehat{\Gamma_lU^{\mathcal{L}h_{\al\be}}}(\xi)\|_{L^2_\xi}&\lesssim \varep_0\langle t\rangle^{H(q+1,n+1)\delta}2^{-N(n+1)k^+},\\
2^{-k/2}(2^{k^-}\langle t\rangle)^{\ga}\|\varphi_k(\xi)\widehat{V^{\mathcal{L}h_{\al\be}}}(\xi)\|_{L^2_\xi}&\lesssim \varep_0\langle t\rangle^{H(q+1,n+1)\delta}2^{-N(n+1)k^+}.
\end{split}
\end{equation}
The bounds \eqref{wer3.2} follow from the main improved energy estimates in Proposition \ref{EnergEst2}, and the commutation identities
\begin{equation}\label{wer3.5}
\Gamma_lU^{\mathcal{L}h_{\al\be}}-U^{\Gamma_l\mathcal{L}h_{\al\be}}=-i\partial_l\Lambda_{wa}^{-1}U^{\mathcal{L}h_{\al\be}}.
\end{equation}
The nonlinear estimates \eqref{wer3.1} follow from  \eqref{plk2} and the observation $H(q,n)+\widetilde{\ell}(n)\leq H(q+1,n+1)-4$ (see \eqref{fvc1.0}).

The inequality for the Klein--Gordon component in \eqref{wer2} follows similarly, using the identity \eqref{wer1} for $\mu=kg$, the improved energy estimates in Proposition \ref{EnergEst2}, and the nonlinear bounds
\begin{equation}\label{wer3.3}
\|\varphi_k(\xi)(\partial_{\xi_l}\widehat{\LL\mathcal{N}^{\psi}})(\xi)\|_{L^2_\xi}\lesssim\varep_1^2\langle t\rangle^{H(q+1,n+1)\delta}2^{-N(n+1)k^+}.
\end{equation}
These nonlinear bounds follow from \eqref{plk3}
\end{proof}

\section{$Z$-norm control of the Klein-Gordon field}\label{kgz}

In this section we prove the bounds \eqref{bootstrap3.4} for the Klein-Gordon field. We notice that, unlike the energy norms, the $Z$ norm of the Klein--Gordon profile $V^\psi$ is not allowed to grow slowly in time. Because of this we need to renormalize the profile $V^\psi$.

\subsection{Renormalization}\label{Nor1}

We start from the equation $\partial_tV^{\psi}=e^{it\Lambda_{kg}}\mathcal{N}^{\psi}$ for the profile $V^{\psi}$, where $\mathcal{N}^{\psi}$ is the nonlinearity defined in \eqref{zaq11.3}. To extract the nonlinear phase correction we examine only the quadratic part of the nonlinearity, which is (see \eqref{npsi2})
\begin{equation}\label{qnpsi}
\mathcal{N}^{\psi,2}:=-h_{00}(\Delta\psi-\psi)+2h_{0j}\partial_0\partial_j\psi-h_{jk}\partial_j\partial_k\psi.
\end{equation}
The formulas in the second line of \eqref{on5} show that
\begin{equation*}
-\widehat{\partial_j\partial_k\psi}(\rho)=\rho_j\rho_k\frac{i\widehat{U^{\psi,+}}(\rho)-i\widehat{U^{\psi,-}}(\rho)}{2\Lambda_{kg}(\rho)},\qquad \widehat{\partial_0\partial_j\psi}(\rho)=i\rho_j\frac{\widehat{U^{\psi,+}}(\rho)+\widehat{U^{\psi,-}}(\rho)}{2}.
\end{equation*}
Therefore, using also the identitities $\widehat{U^{\psi,\pm}}(\rho,t)=e^{\mp it\Lambda_{kg}(\rho)}\widehat{V^{\psi,\pm}}(\rho,t)$,
\begin{equation}\label{Nor2}
e^{it\Lambda_{kg}(\xi)}\widehat{\mathcal{N}^{\psi,2}}(\xi,t)=\frac{1}{(2\pi)^3}\sum_{\pm}\int_{\mathbb{R}^3}ie^{it\Lambda_{kg}(\xi)}e^{\mp it\Lambda_{kg}(\xi-\eta)}\widehat{V^{\psi,\pm}}(\xi-\eta,t)\mathfrak{q}_{kg,\pm}(\xi-\eta,\eta,t)\,d\eta
\end{equation} 
where
\begin{equation}\label{Nor3}
\mathfrak{q}_{kg,\pm}(\rho,\eta,t):=\pm\widehat{h_{00}}(\eta,t)\frac{\Lambda_{kg}(\rho)}{2}+\widehat{h_{0j}}(\eta,t)\rho_j\pm\widehat{h_{jk}}(\eta,t)\frac{\rho_j\rho_k}{2\Lambda_{kg}(\rho)}.
\end{equation}

We would like to eliminate the bilinear interaction between $h_{\al\be}$ and $V^{\psi,+}$ in the first line of \eqref{Nor2} corresponding to $|\eta|\ll 1$. To identify the main term we approximate, heuristically, 
\begin{equation*}
\begin{split}
&\frac{1}{(2\pi)^3}\int_{|\eta|\ll \langle t\rangle^{-1/2}}ie^{it\Lambda_{kg}(\xi)}e^{-it\Lambda_{kg}(\xi-\eta)}\widehat{V^{\psi,+}}(\xi-\eta,t)\mathfrak{q}_{kg,+}(\xi-\eta,\eta,t)\,d\eta\\
&\approx i\frac{\widehat{V^{\psi,+}}(\xi,t)}{(2\pi)^3}\int_{|\eta|\ll \langle t\rangle^{-1/2}}e^{it\eta\cdot\nabla\Lambda_{kg}(\xi)}\Big\{\widehat{h_{00}}(\eta,t)\frac{\Lambda_{kg}(\xi)}{2}+\widehat{h_{0j}}(\eta,t)\xi_j+\widehat{h_{jk}}(\eta,t)\frac{\xi_j\xi_k}{2\Lambda_{kg}(\xi)}\Big\}\,d\eta\\
&\approx i\widehat{V^{\psi}}(\xi,t)\Big\{h_{00}^{low}(t\xi/\Lambda_{kg}(\xi),t)\frac{\Lambda_{kg}(\xi)}{2}+h_{0j}^{low}(t\xi/\Lambda_{kg}(\xi),t)\xi_j+h_{jk}^{low}(t\xi/\Lambda_{kg}(\xi),t)\frac{\xi_j\xi_k}{2\Lambda_{kg}(\xi)}\Big\}
\end{split}
\end{equation*}
where $h_{\al\be}^{low}$ are suitable low-frequency components of $h_{\al\be}$.

We can now define our nonlinear phase correction and nonlinear Klein-Gordon profile precisely. For any $h\in\{h_{\al\be}\}$ we define the low frequency component $h^{low}$ by
\begin{equation}\label{Nor4.1}
\widehat{h^{low}}(\rho,s):=\varphi_{\leq 0}(\langle s\rangle^{p_0}\rho)\widehat{h}(\rho,s),\qquad p_0:=0.68.
\end{equation}  
The reason for this choice of $p_0$, slightly bigger than $2/3$, will become clear later, in the proof of Lemma \ref{BootstrapZ4}. Then we define the Klein--Gordon phase correction
\begin{equation}\label{Nor4}
\begin{split}
\Theta_{kg}(\xi,t):=\int_0^t\Big\{&h_{00}^{low}(s\xi/\Lambda_{kg}(\xi),s)\frac{\Lambda_{kg}(\xi)}{2}\\
&+h_{0j}^{low}(s\xi/\Lambda_{kg}(\xi),s)\xi_j+h_{jk}^{low}(s\xi/\Lambda_{kg}(\xi),s)\frac{\xi_j\xi_k}{2\Lambda_{kg}(\xi)}\Big\}\,ds.
\end{split}
\end{equation}
Finally, we define the nonlinear (modified) Klein--Gordon profile $V_\ast^{\psi}$ by
\begin{equation}\label{Nor5}
\widehat{V^{\psi}_\ast}(\xi,t):=e^{-i\Theta_{kg}(\xi,t)}\widehat{V^{\psi}}(\xi,t).
\end{equation}

We notice that the functions $h^{low}$  are real-valued, thus $\Theta_{kg}$ is real-valued. Let $h^{high}_{\al\be}:=h_{\al\be}-h_{\al\be}^{low}$ and recall the definitions \eqref{Nor3}. For $X\in\{low,\,high\}$ let
\begin{equation}\label{Nor3.1}
\mathfrak{q}_{kg,\pm}^X(\rho,\eta,t):=\pm\widehat{h_{00}^X}(\eta,t)\frac{\Lambda_{kg}(\rho)}{2}+\widehat{h_{0j}^X}(\eta,t)\rho_j\pm\widehat{h_{jk}^X}(\eta,t)\frac{\rho_j\rho_k}{2\Lambda_{kg}(\rho)}.
\end{equation}

The formula \eqref{Nor2} and the equation $\partial_tV^{\psi}=e^{it\Lambda_{kg}}\mathcal{N}^{\psi}$ show that
\begin{equation}\label{Nor7}
\begin{split}
\partial_t\widehat{V^{\psi}_\ast}(\xi,t)&=e^{-i\Theta_{kg}(\xi,t)}\{\partial_t\widehat{V^{\psi}}(\xi,t)-i\widehat{V^{\psi}}(\xi,t)\dot{\Theta}_{kg}(\xi,t)\}=\sum_{a=1}^4\mathcal{R}_a^{\psi}(\xi,t),
\end{split}
\end{equation}
where
\begin{equation}\label{Nor9}
\mathcal{R}_1^{\psi}(\xi,t):=\frac{e^{-i\Theta_{kg}(\xi,t)}}{(2\pi)^3}\int_{\mathbb{R}^3}ie^{it\Lambda_{kg}(\xi)}e^{it\Lambda_{kg}(\xi-\eta)}\widehat{V^{\psi,-}}(\xi-\eta,t)\mathfrak{q}^{low}_{kg,-}(\xi-\eta,\eta,t)\,d\eta,
\end{equation}
\begin{equation}\label{Nor8}
\begin{split}
\mathcal{R}_2^{\psi}(\xi,t):=\frac{e^{-i\Theta_{kg}(\xi,t)}}{(2\pi)^3}\int_{\mathbb{R}^3}i\big\{&e^{it(\Lambda_{kg}(\xi)-\Lambda_{kg}(\xi-\eta))}\widehat{V^{\psi}}(\xi-\eta,t)\mathfrak{q}_{kg,+}^{low}(\xi-\eta,\eta,t)\\
&-e^{it(\xi\cdot\eta)/\Lambda_{kg}(\xi)}\widehat{V^{\psi}}(\xi,t)\mathfrak{q}^{low}_{kg,+}(\xi,\eta,t)\big\}\,d\eta,
\end{split}
\end{equation}
\begin{equation}\label{Nor10}
\mathcal{R}_3^{\psi}(\xi,t):=\frac{e^{-i\Theta_{kg}(\xi,t)}}{(2\pi)^3}\sum_{\iota\in\{+,-\}}\int_{\mathbb{R}^3}ie^{it\Lambda_{kg}(\xi)}e^{-it\Lambda_{kg,\iota}(\xi-\eta)}\widehat{V^{\psi,\iota}}(\xi-\eta,t)\mathfrak{q}^{high}_{kg,\iota}(\xi-\eta,\eta,t)\,d\eta,
\end{equation}
and
\begin{equation}\label{Nor10.1}
\mathcal{R}_4^{\psi}(\xi,t):=e^{-i\Theta_{kg}(\xi,t)}e^{it\Lambda_{kg}(\xi)}[\widehat{\mathcal{N}^\psi}(\xi,t)-\widehat{\mathcal{N}^{\psi,2}}(\xi,t)].
\end{equation}

\subsection{Improved control} In the rest of this section we prove our main $Z$-norm estimate for the profile $V^{\psi}$.

\begin{proposition}\label{ZKGim}
We have, for any $t\in[0,T]$,
\begin{equation}\label{dor1}
\|V^{\psi}(t)\|_{Z_{kg}}\lesssim\varep_0.
\end{equation}
\end{proposition}

Since $|\widehat{V^{\psi}}(\xi,t)|=|\widehat{V^{\psi}_\ast}(\xi,t)|$, in view of the definition \eqref{sec5.1} it suffices to prove that
\begin{equation}\label{Nor40}
\|\varphi_k(\xi)\{\widehat{V^{\psi}_\ast}(\xi,t_2)-\widehat{V^{\psi}_\ast}(\xi,t_1)\}\|_{L^\infty_\xi}\lesssim \varep_02^{-\delta m/2}2^{-k^-/2+\kappa k^-}2^{-N_0k^+}
\end{equation}
for any $k\in\mathbb{Z}$, $m\geq 1$, and $t_1,t_2\in[2^m-2,2^{m+1}]\cap [0,T]$. We prove these bounds in several steps. We start with the contribution of very low frequencies.

\begin{lemma}\label{transfer} 
The bounds \eqref{Nor40} hold if $k\leq -\kappa m$ or if $k\geq \delta'm-10$.
\end{lemma}

\begin{proof} {\bf{Step 1.}} We start with the case of large $k\geq \delta'm -10$. Notice that 
\begin{equation*}
\|Q_{j,k}V^\psi(t)\|_{L^2}\lesssim \varep_0\langle t\rangle^{\delta'}2^{-N(0)k^+},\qquad 2^{j+k}\|Q_{j,k}V^\psi(t)\|_{H^{0,1}_\Omega}\lesssim \varep_0\langle t\rangle^{\delta'}2^{-N(2)k^+},
\end{equation*}
for any $j\geq -k^-$, due to Propositions \ref{EnergEst2} and \ref{DiEs1}. Using now \eqref{consu2.1} we have
\begin{equation*}
\|\widehat{P_kV^\psi}(t)\|_{L^\infty}\lesssim 2^{-3k/2}\varep_0\langle t\rangle^{\delta'}\cdot 2^{-N(0)k^+(1-\delta)/4}2^{-N(2)k^+(3+\delta)/4}\lesssim \varep_0\langle t\rangle^{\delta'}2^{-3k/2}2^{-N_0k^+}2^{-dk^+},
\end{equation*}
for any $t\in[0,T]$ and $k\in\mathbb{Z}$ (recall \eqref{fvc1.1}). The bounds \eqref{Nor40} follow if $2^k\gtrsim 2^{\delta'm}$.

{\bf{Step 2.}} It remains to show that if  $k\leq -\kappa m$ and $t\in[2^m-2,2^{m+1}]\cap [0,T]$ then
\begin{equation}\label{nar1}
\|\varphi_k(\xi)\widehat{V^{\psi}}(\xi,t)\|_{L^\infty_\xi}\lesssim \varep_02^{-\delta m/2}2^{-k/2+\kappa k}.
\end{equation} It follows from Proposition \ref{DiEs1} that
\begin{equation*}
2^{k^+}\|\varphi_{k}(\xi)(\partial_{\xi_l}\widehat{V^{\mathcal{L}\psi}})(\xi,t)\|_{L^2_\xi}\lesssim\varep_0\langle t\rangle^{H(q+1,n+1)\delta}2^{-N_0k^++(n+1)dk^+},
\end{equation*}
for any $t\in[0,T]$, $k\in\mathbb{Z}$, $l\in\{1,2,3\}$, and $\mathcal{L}\in\mathcal{V}_n^q$, $n\in [0,2]$. Using Lemma \ref{hyt1} (i), we have
\begin{equation}\label{Nor41.2}
\sup_{j\geq -k^-}2^j\|Q_{jk}V^{\mathcal{L}\psi}(t)\|_{L^2}\lesssim\varep_0\langle t\rangle^{H(q+1,n+1)\delta}2^{-N_0k^+-k^++(n+1)dk^+}.
\end{equation}
We use now \eqref{Nor41.2} and \eqref{consu2} to estimate
\begin{equation*}
\|\widehat{P_kV^{\psi}}(t)\|_{L^\infty}\lesssim 2^{-3k/2}\big\{\sup_{j\geq -k^-}\|Q_{jk}V^{\psi}(t)\|_{H^{0,1}_\Omega}\big\}^{(1-\delta)/2}\big\{\varep_0\langle t\rangle^{H(1,2)\delta}2^k\big\}^{(1+\delta)/2}.
\end{equation*}
Therefore, recalling that $\kappa^2\geq 4\delta'$, for \eqref{nar1} it suffices to prove that
\begin{equation*}
\|P_kV^{\psi}(t)\|_{H^{0,1}_\Omega}\lesssim\varep_0\langle t\rangle^{2\delta'}2^{k+10\kappa k}
\end{equation*}
if $k\leq-\kappa m$ and  $t\in[2^m-2,2^{m+1}]\cap [0,T]$. In view of \eqref{hyt3.1}, for this it suffices to prove that
\begin{equation}\label{nar4}
\|P_lV^{\Omega\psi}(t)\|_{L^2}+\sum_{a=1}^3\|\varphi_l(\xi)(\partial_{\xi_a}\widehat{V^{\Omega\psi}})(\xi,t)\|_{L^2}\lesssim\varep_0\langle t\rangle^{2\delta'}2^{10\kappa l},
\end{equation}
for any $l\in\mathbb{Z}$ and $t\in[0,T]$, where $\Omega\in\{Id,\Omega_{23},\Omega_{31},\Omega_{12}\}$.

The bound on the first term in the left-hand side of \eqref{nar4} follows from \eqref{Nor41.2}. To bound the remaining terms we use the identities \eqref{wer1}. For \eqref{nar4} it suffices to prove that
\begin{equation}\label{nar5}
\|P_l\Gamma_aU^{\Omega\psi}(t)\|_{L^2}+\|\varphi_l(\xi)(\partial_{\xi_a}\widehat{\Omega\mathcal{N}^{\psi}})(\xi,t)\|_{L^2}\lesssim\varep_0\langle t\rangle^{2\delta'}2^{10\kappa l},
\end{equation}
for any $l\in\mathbb{Z}$, $t\in[0,T]$, and $a\in\{1,2,3\}$. The term $\|P_l\Gamma_aU^{\Omega\psi}(t)\|_{L^2}$ is bounded as claimed due to \eqref{Nor41.2} and the commutation identities \eqref{wer3.5} (with $\Lambda_{wa}$ replaced by $\Lambda_{kg}$ and $h_{\al\be}$ replaced by $\psi$). Therefore, it remains to prove that
\begin{equation}\label{nar6}
\begin{split}
\|P_k{\Omega\mathcal{N}^\psi}(\xi,t)\|_{L^2}&\lesssim\varep_0\langle t\rangle^{2\delta'}2^{10\kappa k}\min\{\langle t\rangle^{-1},2^{k^-}\}2^{-2k^+}\\
\|\varphi_k(\xi)(\partial_{\xi_a}\widehat{\Omega\mathcal{N}^{\psi}})(\xi,t)\|_{L^2}&\lesssim\varep_0\langle t\rangle^{2\delta'}2^{10\kappa k}
\end{split}
\end{equation}
for any $k\in\mathbb{Z}$, $t\in[0,T]$, $\Omega\in\{Id,\Omega_{23},\Omega_{31},\Omega_{12}\}$, and $a\in\{1,2,3\}$.

{\bf{Step 3.}} The bounds \eqref{nar6} are similar to the bounds in Lemmas \ref{dtv8.2} and \ref{dtv8.0}. We may assume $k\leq 0$ and the only issue is to gain the factors $2^{10\kappa k}$ and we are allowed to lose small powers $\langle t\rangle^{2\delta'}$. Notice that the cubic components $\Omega\,\mathcal{N}^{\psi,\geq 3}$ satisfy stronger bounds (this follows from \eqref{abc31.00}--\eqref{abc31.300} if $2^k\gtrsim\langle t\rangle^{-1}$ and the  $L^2$ estimates \eqref{gb5} and \eqref{cnb2} if $2^k\lesssim\langle t\rangle^{-1}$). After these reductions, with $I$ as in \eqref{mults}--\eqref{abc36.1}, for \eqref{nar6} it suffices to prove that
\begin{equation}\label{nar7}
\sum_{(k_1,k_2)\in\mathcal{X}_k}2^{k_2^+-k_1}\|P_kI[P_{k_1}U^{\mathcal{L}_1h,\iota_1},P_{k_2}U^{\mathcal{L}_2\psi,\iota_2}](t)\|_{L^2}\lesssim \varep_1^2\langle t\rangle^{2\delta'}2^{10\kappa k}\min(\langle t\rangle^{-1},2^{k})
\end{equation}
and
\begin{equation}\label{nar8}
\sum_{(k_1,k_2)\in\mathcal{X}_k}2^{k_2^+-k_1}\|\varphi_k(\xi)(\partial_{\xi_a}\mathcal{F}\{I[P_{k_1}U^{\mathcal{L}_1h,\iota_1},P_{k_2}U^{\mathcal{L}_2\psi,\iota_2}]\})(\xi,t)\|_{L^2_\xi}\lesssim\varep_1^2\langle t\rangle^{2\delta'}2^{10\kappa k},
\end{equation}
for any $k\leq 0$, $h\in\{h_{\al\be}\}$, $\iota_1,\iota_2\in\{+,-\}$, $a\in\{1,2,3\}$, $\mathcal{L}_1\in\mathcal{V}_{n_1}^0$, $\mathcal{L}_2\in\mathcal{V}_{n_2}^0$, $n_1+n_2\leq 1$.

{\bf{Substep 3.1.}} We prove first \eqref{nar7}. These bounds easily follow from \eqref{box82} when $2^{k}\lesssim\langle t\rangle^{-1}$. On the other hand, if $2^{k}\geq \langle t\rangle^{-1}$ then we estimate, using \eqref{cnb2} and \eqref{wws1}, 
\begin{equation*}
\begin{split}
2^{k_2^+-k_1}\|P_kI[P_{k_1}U^{\mathcal{L}_1h,\iota_1},P_{k_2}U^{\mathcal{L}_2\psi,\iota_2}](t)\|_{L^2}&\lesssim 2^{k_2^+-k_1}\|P_{k_1}U^{\mathcal{L}_1h}(t)\|_{L^\infty}\|P_{k_2}U^{\mathcal{L}_2\psi}(t)\|_{L^2}\\
&\lesssim\varep_1^22^{k_2^-}\langle t\rangle^{-1+\delta'}\min(1,2^{k_1^-}\langle t\rangle)^{1/2}2^{-4(k_1^++k_2^+)}.
\end{split}
\end{equation*}
This suffices to bound the contribution of the pairs $(k_1,k_2)\in\mathcal{X}_k$ in \eqref{nar7} for which $2^{k_2}\lesssim 2^{k/10}$. For the remaining pairs we have $\min(k_1,k_2)\geq k/10+10$, thus $|k_1-k_2|\leq 4$. Let $J_1$ denote the largest integer such that $2^{J_1}\leq\langle t\rangle(1+2^{k_1}\langle t\rangle)^{-\delta}$ and decompose $P_{k_1}U^{\mathcal{L}_1h,\iota_1}(t)=U_{\leq J_1,k_1}^{\mathcal{L}_1h,\iota_1}(t)+U_{>J_1,k_1}^{\mathcal{L}_1h,\iota_1}(t)$, as in \eqref{box15}. We use first Lemma \ref{Lembil2}, therefore
\begin{equation*}
\begin{split}
2^{k_2^+-k_1}&\big\|P_kI[U_{\leq J_1,k_1}^{\mathcal{L}_1h,\iota_1}(t),P_{k_2}U^{\mathcal{L}_2\psi,\iota_2}(t)]\big\|_{L^2}\\
&\lesssim 2^{k_2^+-k_1}2^{k/2}(1+2^{k_1}\langle t\rangle)^{\delta}\langle t\rangle^{-1}\|Q_{\leq J_1,k_1}V^{\mathcal{L}_1h}(t)\|_{H^{0,1}_{\Omega}}\|P_{k_2}U^{\mathcal{L}_2\psi}(t)\|_{L^2}\\
&\lesssim\varep_1^22^{k/2}\langle t\rangle^{-1+\delta'}2^{-k_1^-/2}2^{-4(k_1^++k_2^+)}.
\end{split}
\end{equation*}
Moreover, using \eqref{vcx1.1} and \eqref{wws2},
\begin{equation*}
\begin{split}
2^{k_2^+-k_1}&\big\|P_kI[U_{>J_1,k_1}^{\mathcal{L}_1h,\iota_1}(t),P_{k_2}U^{\mathcal{L}_2\psi,\iota_2}(t)]\big\|_{L^2}\lesssim 2^{k_2^+-k_1}\|U_{>J_1,k_1}^{\mathcal{L}_1h}(t)\|_{L^2}\|P_{k_2}U^{\mathcal{L}_2\psi}(t)\|_{L^\infty}\\
&\lesssim\varep_1^2\langle t\rangle^{-1+\delta'}2^{-3k_1^-/2}2^{-J_1}2^{-4(k_1^++k_2^+)}.
\end{split}
\end{equation*}
Therefore, recalling that $2^{k}\in[\langle t\rangle^{-1},1]$ and $\min(k_1,k_2)\geq k/10+10$, we have
\begin{equation*}
2^{k_2^+-k_1}\big\|P_kI[U_{\leq J_1,k_1}^{\mathcal{L}_1h,\iota_1}(t),P_{k_2}U^{\mathcal{L}_2\psi,\iota_2}(t)]\big\|_{L^2}\lesssim\varep_1^22^{k/2}\langle t\rangle^{-1+\delta'}2^{-k_1^-/2}2^{-4(k_1^++k_2^+)}.
\end{equation*}
This suffices to complete the proof of \eqref{nar7}.

{\bf{Substep 3.2.}} To prove \eqref{nar8} we write $U^{\mathcal{L}_2\psi,\iota_2}=e^{-it\Lambda_{kg,\iota_2}}V^{\mathcal{L}_2\psi,\iota_2}$ and examine the formula \eqref{abc36.1}. We make the change of variables $\eta\to\xi-\eta$ and notice that the $\partial_{\xi_a}$ derivative can hit either the multiplier $m(\eta,\xi-\eta)$, or the phase $e^{-it\Lambda_{kg,\iota_2}(\xi-\eta)}$, or the profile $\widehat{P_{k_2}V^{\mathcal{L}_2\psi,\iota_2}}(\xi-\eta)$. In the first two cases, the derivative effectively corresponds to multiplying by factors $\lesssim \langle t\rangle$ or $\lesssim 2^{-k_2^-}$. The desired bounds are then consequences of the bounds \eqref{nar7} (in the case $2^{k_2^-}\lesssim \langle t\rangle^{-1}$ we need to apply \eqref{box82} again to control the corresponding contributions).

Finally, assume that the $\partial_{\xi_a}$ derivative hits the profile $\widehat{P_{k_2}V^{\mathcal{L}_2\psi,\iota_2}}(\xi-\eta)$. Letting (as in Lemma \ref{dtv8.2}) $\widehat{U^{\mathcal{L}_2\psi,\iota_2}_{\ast a,k_2}}(\xi,t)=e^{-it\Lambda_{kg,\iota_2}(\xi)}\partial_{\xi_a}\{\varphi_{k_2}\cdot \widehat{V^{\mathcal{L}_2\psi,\iota_2}}\}(\xi,t)$ it suffices to prove that 
\begin{equation}\label{nar99.2}
\sum_{(k_1,k_2)\in\mathcal{X}_k}2^{k_2^+-k_1}\|P_kI[P_{k_1}U^{\mathcal{L}_1h,\iota_1},U^{\mathcal{L}_2\psi,\iota_2}_{\ast a,k_2}](t)\|_{L^2}\lesssim\varep_1^2\langle t\rangle^{2\delta'}2^{10\kappa k}.
\end{equation}
This follows easily using the $L^2$ bounds \eqref{abc99.3} and \eqref{vcx1}.
\end{proof}

We return now to the main estimates \eqref{Nor40}. In the remaining range $-\kappa m\leq k\leq\delta'm -10$, they follow from the identity \eqref{Nor7} and the bounds (which are proved in Lemmas \ref{BootstrapZ6}--\ref{BootstrapZ5} below) 
\begin{equation}\label{Nor50}
\Big\|\varphi_k(\xi)\int_{t_1}^{t_2}\mathcal{R}_a^{\psi}(\xi,s)\,ds\Big\|_{L^\infty_\xi}\lesssim \varep_1^22^{-100\delta' m}
\end{equation}
for $a\in\{1,2,3,4\}$, $m\geq 100$, and $t_1\leq t_2\in[2^{m-1},2^{m+1}]\cap [0,T]$.

In some estimates we need to use integration by parts in time (normal forms). For $\mu\in\{(kg,+),(kg,-)\}$, $\nu\in\{(wa,+),(wa,-)\}$, and $s\in[0,T]$ we define the operators $T^{kg}_{\mu\nu}$ by
\begin{equation}\label{Nor53.1}
T^{kg}_{\mu\nu}[f,g](\xi,s):=\int_{\mathbb{R}^3}\frac{e^{is\Phi_{(kg,+)\mu\nu}(\xi,\eta)}}{\Phi_{(kg,+)\mu\nu}(\xi,\eta)}m(\xi-\eta,\eta)\widehat{f}(\xi-\eta,s)\widehat{g}(\eta,s)\,d\eta,
\end{equation}
where $\Phi_{(kg,+)\mu\nu}(\xi,\eta)=\Lambda_{kg}(\xi)-\Lambda_{\mu}(\xi-\eta)-\Lambda_{\nu}(\eta)$ (see \eqref{on9.2}) and $m\in\mathcal{M}$ (see \eqref{mults}).

\begin{lemma}\label{BootstrapZ6} The bounds \eqref{Nor50} hold for $m\geq 100$, $k\in [-\kappa m,\delta'm]$, and $a=4$.
\end{lemma}

\begin{proof} We use the bounds \eqref{consu2}, combined with either \eqref{wer4.1}--\eqref{plk2} or \eqref{abc31.00}--\eqref{abc31.300} (in both cases $n=1$).  It follows that
\begin{equation}\label{Nor78.1}
\begin{split}
\|\widehat{P_k\mathcal{N}^\psi}(t)\|_{L^\infty}&\lesssim\varep_1^2\langle t\rangle^{-1/2+\delta'}2^{-k^-}2^{-N(2)k^+-2k^+},\\
\|\widehat{P_k\mathcal{N}^{\psi,\geq 3}}(t)\|_{L^\infty}&\lesssim\varep_1^2\langle t\rangle^{-1.1+\delta'}2^{-k^-}2^{-N(2)k^+-2k^+},
\end{split}
\end{equation}
for any $t\in[0,T]$ and $k\in\mathbb{Z}$. The estimates \eqref{Nor50} follow easily from definitions.
\end{proof}

\begin{lemma}\label{BootstrapZ3} The bounds \eqref{Nor50} hold for $m\geq 100$, $k\in [-\kappa m,\delta'm]$, and $a=1$.
\end{lemma}

\begin{proof} We examine the formula \eqref{Nor9}, substitute $h=i\Lambda_{wa}^{-1}(U^{h,+}-U^{h,-})/2$, $h\in \{h_{\al\be}\}$ (see \eqref{on5}) and decompose the input functions dyadically in frequency. For $h\in \{h_{\al\be}\}$ let 
\begin{equation}\label{ulow}
\widehat{U_{low}^{h,\iota_2}}(\xi,s):=\varphi_{\leq 0}(\langle s\rangle^{p_0}\xi)\widehat{U^{h,\iota_2}}(\xi,s),\qquad \widehat{V_{low}^{h,\iota_2}}(\xi,s):=\varphi_{\leq 0}(\langle s\rangle^{p_0}\xi)\widehat{V^{h,\iota_2}}(\xi,s).
\end{equation} 
With $I$ defined as in \eqref{mults}--\eqref{abc36.1}, it suffices to prove that, for $\iota_2\in\{+,-\}$,
\begin{equation}\label{Nor52}
\begin{split}
\sum_{(k_1,k_2)\in\mathcal{X}_k}2^{k_1^+-k_2}\Big\|\varphi_k(\xi)\int_{t_1}^{t_2}e^{is\Lambda_{kg}(\xi)-i\Theta_{kg}(\xi,s)}&\mathcal{F}\{I[P_{k_1}U^{\psi,-},P_{k_2}U^{h,\iota_2}_{low}]\}(\xi,s)\,ds\Big\|_{L^\infty_\xi}\lesssim \varep_1^22^{-\kappa m}.
\end{split}
\end{equation}

We estimate first, using just \eqref{vcx1.2} and \eqref{vcx1},
\begin{equation}\label{Nor53.5}
\begin{split}
2^{k_1^+-k_2}\big\|\mathcal{F}\{I[P_{k_1}U^{\psi,-},&P_{k_2}U^{h,\iota_2}_{low}]\}(\xi,s)\big\|_{L^\infty_\xi}\lesssim 2^{k_1^+-k_2}\|\widehat{P_{k_1}U^{\psi}}(s)\|_{L^\infty}\|\widehat{P_{k_2}U^{h}_{low}}(s)\|_{L^1}\\
&\lesssim \varep_1^22^{-k_1^-/2+\kappa k_1^-}2^{-N_0k_1^++k_1^+}2^{k_2-\delta k_2}2^{\delta m}
\end{split}
\end{equation} 
This suffices to bound the contribution of the pairs $(k_1,k_2)$ for which $k_2\leq -1.1m-10$. It remains to prove that
\begin{equation}\label{Nor54}
\begin{split}
\Big|\int_{t_1}^{t_2}&\int_{\mathbb{R}^3}e^{is\Lambda_{kg}(\xi)-i\Theta_{kg}(\xi,s)}m(\xi-\eta,\eta)e^{is\Lambda_{kg}(\xi-\eta)}\widehat{P_{k_1}V^{\psi,-}}(\xi-\eta,s)\\
&\times e^{-is\Lambda_{wa,\iota_2}(\eta)}\widehat{P_{k_2}V^{h,\iota_2}_{low}}(\eta,s)\, d\eta ds\Big|\lesssim \varep_1^22^{-\kappa m-\delta'm}2^{k_2},
\end{split}
\end{equation}
for any $\xi$ with $|\xi|\in[2^{k_1-4},2^{k_1+4}]$, provided that
\begin{equation}\label{Nor55}
k_2\in[-1.1m-10,-p_0m+10],\qquad k_1\in[-\kappa m-10,\delta' m+10]. 
\end{equation}

To prove \eqref{Nor54} we integrate by parts in time. Letting $\sigma=(kg,+)$, $\mu=(kg,-),\,\nu=(wa,\iota_2)$, we notice that $\Phi_{\sigma\mu\nu}(\xi,\eta)\gtrsim 1$ in the support of the integral.{\footnote{Here it is important that $\mu\neq (kg,+)$, so the phase is nonresonant. The nonlinear correction \eqref{Nor5} was done precisely to weaken the corresponding resonant contribution of the profile $V^{kg,+}$.}} The left-hand side of \eqref{Nor54} is dominated by $C(I_{kg}(\xi)+II_{kg}(\xi)+III_{kg}(\xi))$, where, with $T^{kg}_{\mu\nu}$ defined as in \eqref{Nor53.1},
\begin{equation}\label{Nor56}
\begin{split}
I_{kg}(\xi)&:=\Big(1+\int_{t_1}^{t_2}|\dot{\Theta}_{kg}(\xi,s)|\,ds\Big)\sup_{s\in [t_1,t_2]}|T_{\mu\nu}^{kg}[P_{k_1}V^{\psi,-},P_{k_2}V^{h,\iota_2}_{low}](\xi,s)|,\\
II_{kg}(\xi)&:=\int_{t_1}^{t_2}|T_{\mu\nu}^{kg}[\partial_s(P_{k_1}V^{\psi,-}),P_{k_2}V^{h,\iota_2}_{low}](\xi,s)|\,ds,\\
III_{kg}(\xi)&:=\int_{t_1}^{t_2}|T_{\mu\nu}^{kg}[P_{k_1}V^{\psi,-},\partial_s(P_{k_2}V^{h,\iota_2}_{low})](\xi,s)|\,ds.
\end{split}
\end{equation}

As in \eqref{Nor53.5}, assuming $k_1,k_2$ as in \eqref{Nor55}, we estimate for any $s\in[t_1,t_2]$ 
\begin{equation*}
|T_{\mu\nu}^{kg}[P_{k_1}V^{\psi,-},P_{k_2}V^{h,\iota_2}_{low}](\xi,s)|\lesssim\varep_1^22^{-k_1^-/2}2^{2k_2}2^{4\delta m}.
\end{equation*}
The definition \eqref{Nor4} shows that $|\dot{\Theta}_{kg}(\xi,s)|\lesssim 2^{k_1^+}\sup_{\al,\be}\|h_{\al\be}\|_{L^\infty}\lesssim 2^{k_1^+}2^{-m+\delta'm}$. Therefore 
\begin{equation}\label{Nor56.1}
I_{kg}(\xi)\lesssim \varep_1^22^{-k_1^-/2}2^{2k_2}2^{2\delta' m}.
\end{equation}
Similarly using \eqref{Nor78.1} and the bounds $\|\widehat{P_{k_2}V^{h,\iota_2}_{low}}(s)\|_{L^1}\lesssim\varep_12^{2k_2}2^{4\delta m}$ it follows that
\begin{equation*}
II_{kg}(\xi)\lesssim \varep_1^22^{m/2+2\delta' m}2^{-k_1^-}2^{2k_2}.
\end{equation*}
Finally, using also \eqref{wer4.0},
\begin{equation}\label{Nor56.2}
\|\mathcal{F}\{\partial_s(P_{k_2}V^{h,\iota_2}_{low})(s)\}\|_{L^1}\lesssim 2^{3k_2/2}\|\partial_s(P_{k_2}V^{h}_{low})(s)\|_{L^2}\lesssim \varep_12^{-m+\delta'm/2}2^{2k_2},
\end{equation}
and it follows that $III_{kg}(\xi)$ can be bounded as in \eqref{Nor56.1}. Therefore
\begin{equation*}
I_{kg}(\xi)+II_{kg}(\xi)+III_{kg}(\xi)\lesssim \varep_1^22^{m/2+2\delta' m}2^{-k_1^-}2^{2k_2}.
\end{equation*}
The desired bounds \eqref{Nor54} follow, recalling that $p_0=0.68$ and the restrictions \eqref{Nor55}.
\end{proof}

\begin{lemma}\label{BootstrapZ4} The bounds \eqref{Nor50} hold for $m\geq 100$, $k\in [-\kappa m,\delta'm]$, and $a=2$.
\end{lemma}

\begin{proof} We decompose $V^{\psi}=\sum_{(k_1,j_1)\in\mathcal{J}}V_{j_1,k_1}^{\psi,+}$ as in \eqref{on11.3} and examine the definition \eqref{Nor8}. For \eqref{Nor50} it suffices to prove that
\begin{equation*}
\begin{split}
\varphi_k(\xi)&\Big|\int_{\mathbb{R}^3}\big\{e^{it(\Lambda_{kg}(\xi)-\Lambda_{kg}(\xi-\eta))}\widehat{V_{j_1,k_1}^{\psi,+}}(\xi-\eta,t)\mathfrak{q}_{kg,+}^{low}(\xi-\eta,\eta,t)\\
&-e^{it(\xi\cdot\eta)/\Lambda_{kg}(\xi)}\widehat{V^{\psi,+}_{j_1,k_1}}(\xi,t)\mathfrak{q}^{low}_{kg,+}(\xi,\eta,t)\big\}\,d\eta\Big|\lesssim \varep_1^22^{-m-\kappa m}2^{-\delta j_1},
\end{split}
\end{equation*}
provided that $|k_1-k|\leq 10$ and $t\in[2^{m-1},2^{m+1}]$. Using also the definitions \eqref{Nor3.1} it suffices to prove that for any multiplier $m\in\mathcal{M}_0$ (see \eqref{mults0}) and $\al,\be\in\{0,1,2,3\}$ we have
\begin{equation}\label{Nor61}
\begin{split}
\varphi_k(\xi)&\Big|\int_{\mathbb{R}^3}\widehat{h_{\al\be}^{low}}(\eta,t)\big\{e^{it(\Lambda_{kg}(\xi)-\Lambda_{kg}(\xi-\eta))}\widehat{V_{j_1,k_1}^{\psi,+}}(\xi-\eta,t)m(\xi-\eta)\langle\xi-\eta\rangle\\
&-e^{it(\xi\cdot\eta)/\Lambda_{kg}(\xi)}\widehat{V^{\psi,+}_{j_1,k_1}}(\xi,t)m(\xi)\langle\xi\rangle\big\}\,d\eta\Big|\lesssim \varep_1^22^{-m-\kappa m}2^{-\delta j_1}.
\end{split}
\end{equation}

Recall that $\|\widehat{V_{j_1,k_1}^{\psi,+}}(t)\|_{L^\infty}\lesssim \varep_12^{-k_1}2^{-j_1/2+\delta j_1/2}2^{\delta' m}2^{-4k_1^+}$, see \eqref{vcx1.15}.
Thus, without using the cancellation of the two terms in the integral, the left-hand side of \eqref{Nor61} is bounded by
\begin{equation*}
C\varep_12^{-j_1/2+\delta j_1}2^{\delta' m}2^{-k_1^-}\|\widehat{h_{\al\be}^{low}}(t)\|_{L^1}\lesssim \varep_1^22^{-j_1/2+\delta j_1}2^{2\delta' m}2^{-k_1^-}2^{-p_0m}.
\end{equation*}
This suffices to prove \eqref{Nor61} when $j_1$ is large, i.e. $2^{j_1/2}\gtrsim 2^{(1.01-p_0)m}$. 

On the other hand, if $2^m\gg 1$ and $j_1/2\leq (1.01-p_0)m=0.33 m$ then we estimate
\begin{equation}\label{Nor65}
\begin{split}
\big|e^{it(\Lambda_{kg}(\xi)-\Lambda_{kg}(\xi-\eta))}-e^{it(\xi\cdot\eta)/\Lambda_{kg}(\xi)}\big|&\lesssim 2^{-2p_0m+m},\\
\big|\widehat{V_{j_1,k_1}^{\psi,+}}(\xi-\eta,t)m(\xi-\eta)\langle\xi-\eta\rangle-\widehat{V_{j_1,k_1}^{\psi,+}}(\xi,t)m(\xi)\langle\xi\rangle\big|&\lesssim \varep_12^{j_1/2}2^{2\delta' m}2^{-2k_1^-}2^{-p_0m},
\end{split}
\end{equation}
provided that $|\xi|\approx 2^k$ and $|\eta|\lesssim 2^{-p_0m}$. Indeed, the first bound follows from the observation that $\nabla\Lambda_{kg}(\xi)=\xi/\Lambda_{kg}(\xi)$. The second bound follows from \eqref{vcx1.15} once we notice that taking $\partial_\xi$ derivative of the localized profiles $\widehat{V_{j_1,k_1}^{\psi,+}}$ corresponds essentially to multiplication by a factor of $2^{j_1}$. If $j_1/2\leq 0.33 m$ it follows that the left-hand side of \eqref{Nor61} is bounded by
\begin{equation*}
C\varep_12^{-0.34 m}\|\widehat{h_{\al\be}^{low}}(t)\|_{L^1}\lesssim \varep_1^22^{-p_0m-0.34 m+\delta'm}.
\end{equation*}
This suffices to prove \eqref{Nor61} when $j_1/2\leq 0.33 m$, which completes the proof of the lemma.
\end{proof}

\begin{lemma}\label{BootstrapZ5}
The bounds \eqref{Nor50} hold for $m\geq 100$, $k\in[-\kappa m,\delta'm]$, and $a=3$.
\end{lemma}

\begin{proof} We examine the formula \eqref{Nor10}. Let $U_{high}^{h,\iota_2}:=U^{h,\iota_2}-U_{low}^{h,\iota_2}$, $V_{high}^{h,\iota_2}:=V^{h,\iota_2}-V_{low}^{h,\iota_2}$, see \eqref{ulow}. As in the proof of Lemma \ref{BootstrapZ3}, after simple reductions it suffices to prove that
\begin{equation}\label{Nor72}
\begin{split}
2^{k_1^+-k_2}\Big\|\varphi_k(\xi)\int_{t_1}^{t_2}e^{is\Lambda_{kg}(\xi)-i\Theta(\xi,s)}&\mathcal{F}\{I[P_{k_1}U^{\psi,\iota_1},P_{k_2}U^{h,\iota_2}_{high}]\}(\xi,s)\,ds\Big\|_{L^\infty_\xi}\lesssim \varep_1^22^{-101\delta' m}
\end{split}
\end{equation}
for $\iota_1,\iota_2\in\{+,-\}$, and $(k_1,k_2)\in\mathcal{X}_k$, $k_1,k_2\in [-p_0m-10,m/100]$.

As in the proof of Lemma \ref{BootstrapZ3} we integrate by parts in time to estimate
\begin{equation*}
\Big|\int_{t_1}^{t_2}e^{is\Lambda_{kg}(\xi)-i\Theta_{kg}(\xi,s)}\mathcal{F}\{I[P_{k_1}U^{\psi,\iota_1},P_{k_2}U^{h,\iota_2}_{high}]\}(\xi,s)\,ds\Big|\lesssim I'_{kg}(\xi)+II'_{kg}(\xi)+III'_{kg}(\xi),
\end{equation*}
where, with $\mu=(kg,\iota_1)$ and $\nu=(wa,\iota_2)$ and $T^{kg}_{\mu\nu}$ defined as in \eqref{Nor53.1},
\begin{equation*}
\begin{split}
I'_{kg}(\xi)&:=\Big(1+\int_{t_1}^{t_2}|\dot{\Theta}_{kg}(\xi,s)|\,ds\Big)\sup_{s\in[t_1,t_2]}|T_{\mu\nu}^{kg}[P_{k_1}V^{\psi,\iota_1},P_{k_2}V^{h,\iota_2}_{high}](\xi,s)|,\\
II'_{kg}(\xi)&:=\int_{t_1}^{t_2}|T_{\mu\nu}^{kg}[\partial_s(P_{k_1}V^{\psi,\iota_1}),P_{k_2}V^{h,\iota_2}_{high}](\xi,s)|\,ds,\\
III'_{kg}(\xi)&:=\int_{t_1}^{t_2}|T_{\mu\nu}^{kg}[P_{k_1}V^{\psi,\iota_1},\partial_s(P_{k_2}V^{h,\iota_2}_{high})](\xi,s)|\,ds.
\end{split}
\end{equation*}
Notice that $|\dot{\Theta}_{kg}(\xi,s)|\lesssim 2^{-m+4\delta' m}$, as a consequence of \eqref{wws1}. After possibly changing the multiplier $m$ in the definition \eqref{Nor53.1}, for \eqref{Nor72} it suffices to prove that
\begin{equation}\label{Nor75}
|\varphi_k(\xi)T_{\mu\nu}^{kg}[P_{k_1}V^{\psi,\iota_1},P_{k_2}V^{h,\iota_2}](\xi,s)|\lesssim \varep_1^22^{-110\delta' m}2^{k_2^-},
\end{equation}
\begin{equation}\label{Nor75.2}
2^m|\varphi_k(\xi)T_{\mu\nu}^{kg}[\partial_s(P_{k_1}V^{\psi,\iota_1}),P_{k_2}V^{h,\iota_2}](\xi,s)|\lesssim \varep_1^22^{-110\delta' m}2^{k_2^-},
\end{equation}
\begin{equation}\label{Nor75.3}
2^m|\varphi_k(\xi)T_{\mu\nu}^{kg}[P_{k_1}V^{\psi,\iota_1},\partial_s(P_{k_2}V^{h,\iota_2})](\xi,s)|\lesssim \varep_1^22^{-110\delta' m}2^{k_2^-},
\end{equation}
for any $s\in[2^{m-1},2^{m+1}]$, $k_1,k_2\in[-p_0m-10,m/100]$, $\mu=(kg,\iota_1)$, $\nu=(wa,\iota_2)$, $\iota_1,\iota_2\in\{+,-\}$.

{\bf{Step 1: proof of \eqref{Nor75}.}} If $k_1\leq -4\kappa m$ (so $k_2\geq -\kappa m-20$) then we can just use the $L^2$ bounds \eqref{vcx1} and \eqref{cnb2} on the two inputs, and Lemma \ref{pha2} (i). On the other hand, if $k_1\geq -4\kappa m$ then we decompose $P_{k_1}V^{\psi,\iota_1}=\sum_{j_1}V_{j_1,k_1}^{\psi,\iota_1}$ and $P_{k_2}V^{h,\iota_2}=\sum_{j_2}V_{j_2,k_2}^{h,\iota_2}$ as in \eqref{on11.3}. Let $\overline{k}:=\max(k,k_1,k_2)$ and recall that $|\Phi_{(kg,+)\mu\nu}(\xi,\eta)|\gtrsim 2^{k_2}2^{-2\overline{k}^+}$ in the support of the integrals defining the operators $T_{\mu\nu}^{kg}$ (see \eqref{pha3}). 

The contribution of the pairs $(V_{j_1,k_1}^{\psi,\iota_1},V_{j_2,k_2}^{h,\iota_2})$ for which $2^{\max(j_1,j_2)}\leq 2^{0.99m}2^{-6\overline{k}^+}$ is negligible,
\begin{equation}\label{Nor78}
|T^{kg}_{\mu\nu}[V_{j_1,k_1}^{\psi,\iota_1},V_{j_2,k_2}^{h,\iota_2}](\xi,s)|\lesssim \varep_1^22^{-2m}\,\,\,\,\text{ if }\,\,\,\,2^{\max(j_1,j_2)}\leq 2^{0.99m}2^{-6\overline{k}^+}.
\end{equation}
Indeed, this follows by integration by parts in $\eta$ (using Lemma \ref{tech5}), the bounds \eqref{Linfty3.4}, and the observation that the gradient of the phase admits a suitable lower bound $|\nabla_{\eta}\{s\Lambda_{kg,\iota_1}(\xi-\eta)+s\Lambda_{wa,\iota_2}(\eta)\}|\gtrsim \langle s\rangle 2^{-2k_1^+}$ in the support of the integral. On the other hand, we estimate
\begin{equation*}
\begin{split}
|T^{kg}_{\mu\nu}&[V_{j_1,k_1}^{\psi,\iota_1},V_{j_2,k_2}^{h,\iota_2}](\xi,s)|\lesssim 2^{-k_2}2^{2\overline{k}^+}\varep_1^22^{3k_2/2}\|\widehat{V_{j_1,k_1}^{\psi,\iota_1}}(s)\|_{L^\infty}\|\widehat{V_{j_2,k_2}^{h,\iota_2}}(s)\|_{L^2}\\
&\lesssim \varep_1^22^{\delta'm}2^{k_2/2}2^{2k_1^++2k_2^+}\cdot 2^{-k_1}2^{-j_1/2+\delta j_1}2^{-10k_1^+}2^{-j_2}2^{-k_2^-/2-4\delta k_2^-}2^{-10k_2^+}\\
&\lesssim \varep_1^22^{2\delta'm}2^{-k_1}2^{-j_1/2+\delta j_1}2^{-j_2}2^{-6\overline{k}^+},
\end{split}
\end{equation*}
using \eqref{vcx1.1}, \eqref{vcx1.15}, and Lemma \ref{pha2} (i). Recalling that $k_1\geq -4\kappa m$, this suffices to estimate the contribution of the pairs $(V_{j_1,k_1}^{\psi,\iota_1},V_{j_2,k_2}^{h,\iota_2})$ for which $2^{\max(j_1,j_2)}\geq 2^{0.99m}2^{-6\overline{k}^+}$, and the bound \eqref{Nor75} follows.

{\bf{Step 2: proof of \eqref{Nor75.2}.}} Notice that, for any $t\in[0,T]$,
\begin{equation}\label{loc1}
\partial_tV^{\psi,\iota_1}(t)=e^{it\Lambda_{kg,\iota_1}}\mathcal{N}^\psi(t)=e^{it\Lambda_{kg,\iota_1}}\mathcal{N}^{\psi,2}(t)+e^{it\Lambda_{kg,\iota_1}}\mathcal{N}^{\psi,\geq 3}(t).
\end{equation}
The contribution of the nonlinearity $\mathcal{N}^{\psi,\geq 3}$ can be bounded using \eqref{Nor78.1},
\begin{equation*}
\begin{split}
|\varphi_k(\xi)&T_{\mu\nu}^{kg}[e^{is\Lambda_{kg,\iota_1}}P_{k_1}\mathcal{N}^{\psi,\geq 3}(s),P_{k_2}V^{h,\iota_2}(s)](\xi)|\\
&\lesssim 2^{-k_2}2^{2\overline{k}^+}\|\widehat{P_{k_1}\mathcal{N}^{\psi,\geq 3}}(s)\|_{L^\infty}\|P_{k_2}V^{h}(s)\|_{L^2}2^{3\min(k_1,k_2)/2}\lesssim 2^{-1.09m}2^{k_2}2^{-\max(k_1,k_2)}2^{-6\overline{k}^+}.
\end{split}
\end{equation*}
This is better than the bounds \eqref{Nor75.2} since $\max(k_1,k_2)\geq k-10\geq -\kappa m-10$.

To bound the contribution of the quadratic components $\mathcal{N}^{\psi,2}$ we recall that $\mathcal{F}\{P_{k_1}\mathcal{N}^{\psi,2}\}(s)$ can be written as a sum of terms of the form
\begin{equation*}
\varphi_{k_1}(\gamma)\int_{\mathbb{R}^3}|\rho|^{-1}\langle\gamma-\rho\rangle m_3(\gamma-\rho)\widehat{U^{\psi,\iota_3}}(\gamma-\rho)\widehat{U^{h_4,\iota_4}}(\rho)\,d\rho,
\end{equation*}
where $\iota_3,\iota_4\in\{+,-\}$, $h_4\in\{h_{\al\be}\}$, and $m_3$ is a symbol as in \eqref{mults0} (see \eqref{npsi2}). We combine this with the formula \eqref{Nor53.1}. For \eqref{Nor75.2} it suffices to prove that, for any $\xi\in\mathbb{R}^3$,
\begin{equation}\label{Nor81}
\begin{split}
\Big|&\varphi_k(\xi)\int_{\mathbb{R}^3\times\mathbb{R}^3}\frac{\varphi_{k_1}(\xi-\eta)m(\xi-\eta,\eta)}{\Lambda_{kg}(\xi)-\Lambda_\mu(\xi-\eta)-\Lambda_\nu(\eta)}e^{-is\Lambda_\nu(\eta)}\widehat{P_{k_2}V^{h,\iota_2}}(\eta,s)\\
&\times m_3(\xi-\eta-\rho)\langle\xi-\eta-\rho\rangle|\rho|^{-1}\widehat{U^{\psi,\iota_3}}(\xi-\eta-\rho,s)\widehat{U^{h_4,\iota_4}}(\rho,s)\,d\eta d\rho\Big|\lesssim \varep_1^22^{k_2^-}2^{-1.005m},
\end{split}
\end{equation} 
provided that $\mu=(kg,\iota_1)$, $\nu=(wa,\iota_2)$, $\iota_1,\iota_2,\iota_3,\iota_4\in\{+,-\}$, and $k_1,k_2\in[-p_0m-10,m/10]$.

We decompose the solutions $U^{\psi,\iota_3}$, $U^{h_4,\iota_4}$, and $P_{k_2}V^{h,\iota_2}$ dyadically in frequency and space as in \eqref{on11.3}. Then we notice that the contribution when one of the parameters $j_3,k_3,j_4,k_4,j_2$ is large can be bounded using just $L^2$ estimates. For \eqref{Nor81} it suffices to prove that
\begin{equation}\label{Nor82}
2^{-k_2^-}2^{k_3^+-k_4}|\mathcal{C}_{kg}[e^{-is\Lambda_\theta}V_{j_3,k_3}^{\psi,\iota_3}(s),e^{-is\Lambda_\nu}V_{j_2,k_2}^{h,\iota_2}(s),e^{-is\Lambda_\vartheta}V^{h_4,\iota_4}_{j_4,k_4}(s)](\xi)|\lesssim \varep_1^32^{-1.01m}
\end{equation}
for any $k_2\in[-p_0m-10,m/10]$, $k_3,k_4\in[-2m,m/10]$, and $j_2,j_3,j_4\leq 2m$, where $\theta=(kg,\iota_3)$, $\vartheta=(wa,\iota_4)$, and, with $m\in\mathcal{M}$, $m_3,m_4\in\mathcal{M}_0$,
\begin{equation}\label{Nor83}
\begin{split}
\mathcal{C}_{kg}[f,g,h](\xi):=&\int_{\mathbb{R}^3\times\mathbb{R}^3}\frac{\varphi_k(\xi)\varphi_{k_1}(\xi-\eta)m(\xi-\eta,\eta)}{\Lambda_{kg}(\xi)-\Lambda_\mu(\xi-\eta)-\Lambda_\nu(\eta)}\\
&\times m_3(\xi-\eta-\rho)m_4(\rho)\cdot\widehat{f}(\xi-\eta-\rho)\widehat{g}(\eta)\widehat{h}(\rho)\,d\eta d\rho.
\end{split}
\end{equation}

{\bf{Substep 2.1.}} Assume first that
\begin{equation}\label{pol1}
j_3\geq 0.99 m-3k_3^+.
\end{equation}
Let $k^\ast:=\max\{k_2^+,k_3^+,k_4^+\}$. Let $Y$ denote the left-hand side of \eqref{Nor82}. Using Lemmas \ref{pha2} and \ref{L1easy} (i) we estimate
\begin{equation}\label{pol2}
\begin{split}
Y&\lesssim 2^{k^+_3-k_4}2^{-2k_2+6\max(k^+,k_2^+)}\|V_{j_3,k_3}^{\psi,\iota_3}(s)\|_{L^2}\|e^{-is\Lambda_\nu}V_{j_2,k_2}^{h,\iota_2}(s)\|_{L^\infty}\|V^{h_4,\iota_4}_{j_4,k_4}(s)\|_{L^2}\\
&\lesssim \varep_1^32^{4\delta'm}2^{-j_3}2^{-m}2^{-k_2^-}2^{-k_4^-/2}2^{-dk^\ast},
\end{split}
\end{equation}
where in the last line we used \eqref{wws1} and bounds from Lemma \ref{dtv2}. Since $2^{-k_2^-}\lesssim 2^{0.68m}$ and $j_3+3k_3^+\geq 0.99 m$, this suffices to prove \eqref{Nor82} when $k_4^-\geq -0.55m-10$. On the other hand, if $k_4\leq -0.55m-10$, then we estimate in the Fourier space. Using \eqref{pha3}, \eqref{vcx1.15}, and \eqref{vcx1}--\eqref{vcx1.1}
\begin{equation}\label{pol3}
\begin{split}
Y&\lesssim 2^{k_3^+-k_4}2^{-2k_2+6\max(k^+,k_2^+)}\|\widehat{V_{j_3,k_3}^{\psi,\iota_3}}(s)\|_{L^\infty}2^{3k_2/2}\|\widehat{V_{j_2,k_2}^{h,\iota_2}}(s)\|_{L^2}2^{3k_4/2}\|\widehat{V^{h_4,\iota_4}_{j_4,k_4}}(s)\|_{L^2}\\
&\lesssim \varep_1^32^{4\delta' m}2^{-j_3/2}2^{k_4}2^{-k_3}2^{-dk^\ast}.
\end{split}
\end{equation}
Since $k_4\leq -0.55m-10$, this suffices to prove \eqref{Nor82} when $k_3\geq -0.01m-10$. Finally, if $k_4\leq -0.55m-10$ and $k_3\leq -0.01m-10$ then $k_2\geq -\kappa m-10$ (due to the assumption $k\geq -\kappa m$) and a similar estimate gives
\begin{equation}\label{pol4}
\begin{split}
Y&\lesssim 2^{k_3^+-k_4}2^{-2k_2+6\max(k^+,k_2^+)}2^{3k_3/2}\|\widehat{V_{j_3,k_3}^{\psi,\iota_3}}(s)\|_{L^2}\|\widehat{V_{j_2,k_2}^{h,\iota_2}}(s)\|_{L^\infty}2^{3k_4/2}\|\widehat{V^{h_4,\iota_4}_{j_4,k_4}}(s)\|_{L^2}\\
&\lesssim \varep_1^32^{0.01 m}2^{-j_3}2^{k_4}2^{-dk^\ast}.
\end{split}
\end{equation} 
This completes the proof of \eqref{Nor82} when $j_3\geq 0.99 m-3k_3^+$.

{\bf{Substep 2.2.}} Assume now that
\begin{equation}\label{pol6}
j_3\leq 0.99 m-3k_3^+.
\end{equation}
We notice that the $\eta$ gradient of the phase $-s\Lambda_\theta(\xi-\eta-\rho)-s\Lambda_\nu(\eta)$ is $\gtrsim 2^m2^{-2k_3^+}$ in the support of the integral in \eqref{Nor83}. Similarly, the $\rho$ gradient of the phase $-s\Lambda_\theta(\xi-\eta-\rho)-s\Lambda_\vartheta(\rho)$ is $\gtrsim 2^m2^{-2k_3^+}$ in the support of the integral. Using Lemma \ref{tech5} (integration by parts in $\eta$ or $\rho$), the contribution is negligible unless
\begin{equation}\label{pol7}
j_2\geq 0.99 m-3k_3^+\quad\text{ and }\quad j_4\geq 0.99 m-3k_3^+.
\end{equation}

Given \eqref{pol7}, we estimate first, as in \eqref{pol3},
\begin{equation*}
\begin{split}
Y&\lesssim 2^{k_3^+-k_4}2^{-2k_2+6\max(k^+,k_2^+)}\|\widehat{V_{j_3,k_3}^{\psi,\iota_3}}(s)\|_{L^\infty}2^{3k_2/2}\|\widehat{V_{j_2,k_2}^{h,\iota_2}}(s)\|_{L^2}2^{3k_4/2}\|\widehat{V^{h_4,\iota_4}_{j_4,k_4}}(s)\|_{L^2}\\
&\lesssim \varep_1^32^{4\delta' m}2^{-j_2-j_4}2^{-k_2}2^{-k_3/2}2^{-dk^\ast}.
\end{split}
\end{equation*}
This suffices if $k_3\geq -0.2m-10$. On the other hand, if $k_3\leq -0.2m-10$ then we may assume that $\max\{k_2,k_4\}\geq -\kappa m-10$ (due to the assumption $k\geq -\kappa m$) and estimate as in \eqref{pol2},
\begin{equation*}
\begin{split}
Y&\lesssim 2^{k^+_3-k_4}2^{-2k_2+6\max(k^+,k_2^+)}\|e^{-is\Lambda_\theta}V_{j_3,k_3}^{\psi,\iota_3}(s)\|_{L^\infty}\|V_{j_2,k_2}^{h,\iota_2}(s)\|_{L^2}\|V^{h_4,\iota_4}_{j_4,k_4}(s)\|_{L^2}\\
&\lesssim \varep_1^32^{4\delta'm}2^{-m}2^{-j_2-j_4}2^{-5k_2/2}2^{-3k_4/2}2^{-dk^\ast}.
\end{split}
\end{equation*}
Since $2^{-k_2}\lesssim 2^{p_0m}$, this suffices to prove \eqref{Nor82} when $k_4\geq -0.1 m-10$. Finally, if $k_3,k_4\leq -0.1 m-10$ and $k_2\geq -\kappa m-10$ then we estimate as in \eqref{pol4}
\begin{equation*}
\begin{split}
Y&\lesssim 2^{k_3^+-k_4}2^{-2k_2+6\max(k^+,k_2^+)}2^{3k_3/2}\|\widehat{V_{j_3,k_3}^{\psi,\iota_3}}(s)\|_{L^2}\|\widehat{V_{j_2,k_2}^{h,\iota_2}}(s)\|_{L^\infty}2^{3k_4/2}\|\widehat{V^{h_4,\iota_4}_{j_4,k_4}}(s)\|_{L^2}\\
&\lesssim \varep_1^32^{0.01 m}2^{-j_4}2^{-j_2/2}2^{-dk^\ast},
\end{split}
\end{equation*} 
which suffices. This completes the proof of the the bounds \eqref{Nor75.2}.

{\bf{Step 3: proof of \eqref{Nor75.3}.}} Notice that, for any $t\in[0,T]$ and $h\in\{h_{\al\be}\}$,
\begin{equation}\label{loc2}
\partial_tV^{h,\iota_2}(t)=e^{it\Lambda_{wa,\iota_2}}\mathcal{N}^h(t)=e^{it\Lambda_{wa,\iota_2}}\mathcal{N}^{h,2}(t)+e^{it\Lambda_{wa,\iota_2}}\mathcal{N}^{h,\geq 3}(t).
\end{equation}
If $k_1\leq -0.01m-10$ then we may assume that $k_2\geq -\kappa m-10$ (due to the assumption $k\geq -\kappa m$) and estimate the left-hand side of \eqref{Nor75.3} using just \eqref{pha3}, \eqref{cnb2}, and \eqref{wer4.0},
\begin{equation*}
C2^m2^{-k_2}2^{6\overline{k}^+}\|\widehat{P_{k_1}V^\psi}\|_{L^2}\|\widehat{P_{k_2}\mathcal{N}^h}\|_{L^2}\lesssim \varep_1^22^{k_1^-}2^{2\kappa m},
\end{equation*}
which suffices. On the other hand, if $k_1\geq -0.01m-10$ then we decompose $P_{k_1}V^{\psi,\iota_1}=\sum_{j_1\geq -k_1^-}V^{\psi,\iota_1}_{j_1,k_1}$, and notice that the contribution of the localized profiles for which $j_1\geq 0.1m$ can also be bounded in a similar way, using \eqref{vcx1.15} and estimating $\widehat{P_{k_2}\mathcal{N}^h}$ in $L^1$ to gain a factor of $2^{3k_2/2}$. The contribution of the cubic and higher order nonlinearity can also be bounded in the same way, using the stronger estimates \eqref{abc3.00}. It remains to prove that
\begin{equation}\label{Nor91}
2^{-k_2^-}|\varphi_k(\xi)T_{\mu\nu}^{kg}[V^{\psi,\iota_1}_{j_1,k_1}(s),e^{is\Lambda_{wa,\iota_2}}P_{k_2}\mathcal{N}^{h,2}(s)](\xi)|\lesssim \varep_1^22^{-1.005m},
\end{equation}
for any $s\in[2^{m-1},2^{m+1}]$, $k_2\in[-p_0m-10,m/10]$, $k_1\geq -0.01m-10$, and $j_1\leq 0.1 m$.

We examine now the quadratic nonlinearities $\mathcal{N}^{h,2}_{\al\be}$ in \eqref{sac1.1}. They contain two main types of terms: bilinear interactions of the metric components and bilinear interactions of the Klein-Gordon field. We define the trilinear operators
\begin{equation}\label{Nor92}
\begin{split}
\mathcal{C}'_{kg}[f,g,h](\xi):=&\int_{\mathbb{R}^3\times\mathbb{R}^3}\frac{\varphi_k(\xi)\varphi_{k_2}(\eta)m(\xi-\eta,\eta)}{\Lambda_{kg}(\xi)-\Lambda_\mu(\xi-\eta)-\Lambda_\nu(\eta)}\\
&\times m_3(\eta-\rho)m_4(\rho)\cdot\widehat{f}(\xi-\eta)\widehat{g}(\eta-\rho)\widehat{h}(\rho)\,d\eta d\rho,
\end{split}
\end{equation} 
where $m\in\mathcal{M}$, $m_3,m_4\in\mathcal{M}_0$. For \eqref{Nor91} it suffices to prove that
\begin{equation}\label{Nor94}
\begin{split}
2^{-k_2^-}2^{|k_3-k_4|}|\mathcal{C}'_{kg}[e^{-is\Lambda_{\mu}}V^{\psi,\iota_1}_{j_1,k_1}(s),e^{-is\Lambda_{wa,\iota_3}}V^{h_3,\iota_3}_{j_3,k_3}(s),e^{-is\Lambda_{wa,\iota_4}}V^{h_4,\iota_4}_{j_4,k_4}(s)](\xi)|\\
\lesssim \varep_1^22^{-1.01m}2^{-\gamma(j_3+|k_3|+j_4+|k_4|)}
\end{split}
\end{equation}
and
\begin{equation}\label{Nor95}
\begin{split}
2^{-k_2^-}|\mathcal{C}'_{kg}[e^{-is\Lambda_{\mu}}V^{\psi,\iota_1}_{j_1,k_1}(s),e^{-is\Lambda_{kg,\iota_3}}V^{\psi,\iota_3}_{j_3,k_3}(s),e^{-is\Lambda_{kg,\iota_4}}V^{\psi,\iota_4}_{j_4,k_4}(s)](\xi)|\\
\lesssim \varep_1^22^{-1.01m}2^{-\gamma(j_3+|k_3|+j_4+|k_4|)},
\end{split}
\end{equation}
where $s,k_1,k_2,j_1$ are as in \eqref{Nor91}, $h_3,h_4\in\{h_{\al\be}\}$, and $(k_3,j_3),(k_4,j_4)\in\mathcal{J}$.

{\bf{Substep 3.1: proof of \eqref{Nor94}.}} Using just $L^2$ estimates, we may assume that $k_4\leq k_3\leq m/10+10$. Notice that the $\eta$ gradient of the phase $-s\Lambda_{kg,\iota_1}(\xi-\eta)-s\Lambda_{wa,\iota_3}(\eta-\rho)$ is $\gtrsim 2^m2^{-2k_1^+}$ in the support of the integral. Therefore, using integration by parts in $\eta$ (Lemma \ref{tech5}),  the integral is negligible if $j_3\leq 0.99 m-3k_1^+$. Similarly, by making the change of variables $\rho\to\eta-\rho$, the integral is negligible if $j_4\leq 0.99 m-3k_1^+$. Finally, if $j_3\geq 0.99 m-3k_1^+$ and $j_4\geq 0.99 m-3k_1^+$ then we use \eqref{pha3} to estimate the left-hand side of \eqref{Nor94} by
\begin{equation*}
\begin{split}
C2^{-k_2}2^{k_3-k_4}2^{-k_2+6\max(k^+,k_1^+)}&\|\widehat{V^{\psi,\iota_1}_{j_1,k_1}}(s)\|_{L^\infty}2^{3k_2/2}\|\widehat{V^{h_3,\iota_3}_{j_3,k_3}}(s)\|_{L^2}2^{3k_4/2}\|\widehat{V^{h_4,\iota_4}_{j_4,k_4}}(s)\|_{L^2}\\
&\lesssim \varep_1^32^{0.01m}2^{-k_2/2}2^{-k_1^-}2^{-j_3}2^{-j_4(1-\delta')}2^{-d\max\{k_1^+,k_3^+\}}.
\end{split}
\end{equation*}
This suffices to prove \eqref{Nor94}.

{\bf{Substep 3.2: proof of \eqref{Nor95}.}} Using Lemma \ref{pha2} (ii) and \eqref{wws2.5}, we estimate the left-hand side of \eqref{Nor95} by
\begin{equation*}
\begin{split}
C2^{-2k_2+6\max(k^+,k_1^+)}\|e^{-is\Lambda_{\mu}}V^{\psi,\iota_1}_{j_1,k_1}(s)\|_{L^\infty}\|V^{\psi,\iota_3}_{j_3,k_3}(s)\|_{L^2}\|V^{\psi,\iota_4}_{j_4,k_4}(s)\|_{L^2}\\
\lesssim \varep_1^32^{-1.49m}2^{-j_3}2^{-j_4}2^{-2k_2}2^{-k_3^+-k_4^+}.
\end{split}
\end{equation*} 
This suffices if $2^{-0.48m}2^{-j_3}2^{-j_4}2^{-2k_2}\lesssim 2^{-\gamma(j_3+j_4)}$. Otherwise, if $(1-\gamma)(j_3+j_4)+0.48 m\leq -2k_2-120$ then we may assume that $j_3\leq j_4$ (so $j_3\leq 0.45 m-40$ since $k_2\geq -0.68m-10$) and use \eqref{wws2.5} again to estimate the left-hand side of \eqref{Nor95} by
\begin{equation*}
\begin{split}
C2^{-2k_2+6\max(k^+,k_1^+)}&\|\widehat{V^{\psi,\iota_1}_{j_1,k_1}}(s)\|_{L^\infty}2^{3k_2/2}\|e^{-is\Lambda_{kg,\iota_3}}V^{\psi,\iota_3}_{j_3,k_3}(s)\|_{L^\infty}\|V^{\psi,\iota_4}_{j_4,k_4}(s)\|_{L^2}\\
&\lesssim \varep_1^32^{-1.49m}2^{-k_2^-/2}2^{-k_3^-/2}2^{-j_4}.
\end{split}
\end{equation*}
The bounds \eqref{Nor95} follow since $2^{-j_4}2^{-k_3^-/2}\lesssim 1$. This completes the proof of the lemma.
\end{proof}

\section{$Z$-norm control of the metric components}\label{waz}

In this section we prove the bounds \eqref{bootstrap3.4} for the metric components.

\begin{proposition}\label{ZWAim}
With the hypothesis of Proposition \ref{bootstrap}, for any $t\in[0,T]$, $a,b\in\{1,2,3\}$ and $\al,\be\in\{0,1,2,3\}$ we have
\begin{equation}\label{ben1}
\|V^{F}(t)\|_{Z_{wa}}+\|V^{\omega_a}(t)\|_{Z_{wa}}+\|V^{\vartheta_{ab}}(t)\|_{Z_{wa}}+\langle t\rangle^{-\delta}\|V^{h_{\al\be}}(t)\|_{Z_{wa}}\lesssim\varep_0.
\end{equation}
\end{proposition}

The rest of the section is concerned with the proof of this proposition. As in the previous section, we need to first renormalize the profiles $V^{h_{\al\be}}$. The nonlinear phase correction is determined only by the quadratic quasilinear components of the nonlinearity
\begin{equation}\label{bax1}
\mathcal{Q}_{\al\be}^2:=\big\{-h_{00}\Delta+2h_{0j}\partial_0\partial_j-h_{jk}\partial_j\partial_k\big\}h_{\al\be},
\end{equation}
see \eqref{sac1.2}. The point of the renormalization is to weaken some of the resonant bilinear interactions corresponding to very low frequencies of the metric components.

We define our nonlinear phase correction and the nonlinear profiles associated to the metric components as in subsection \ref{Nor1}. As in \eqref{Nor4.1}, for any $h\in\{h_{\mu\nu}\}$ we define the low frequency component $h^{low}$ by $\widehat{h^{low}}(\rho,s)=\varphi_{\leq 0}(\langle s\rangle^{p_0}\rho)\widehat{h}(\rho,s)$, $p_0=0.68$. Then we define the correction
\begin{equation}\label{bax3}
\begin{split}
\Theta_{wa}(\xi,t):=\int_0^t\Big\{&h_{00}^{low}(s\xi/\Lambda_{wa}(\xi),s)\frac{\Lambda_{wa}(\xi)}{2}\\
&+h_{0j}^{low}(s\xi/\Lambda_{wa}(\xi),s)\xi_j+h_{jk}^{low}(s\xi/\Lambda_{wa}(\xi),s)\frac{\xi_j\xi_k}{2\Lambda_{wa}(\xi)}\Big\}\,ds.
\end{split}
\end{equation}
Finally, we define the nonlinear (modified) profiles of the metric components by
\begin{equation}\label{bax4}
\widehat{V^{G}_\ast}(\xi,t):=e^{-i\Theta_{wa}(\xi,t)}\widehat{V^{G}}(\xi,t),\qquad G\in\{h_{\al\be}, F, \omega_a,\vartheta_{ab}\}.
\end{equation}

We notice that the functions $h^{low}$  are real-valued, thus $\Theta_{wa}$ is real-valued. Let $h^{high}_{\al\be}:=h_{\al\be}-h_{\al\be}^{low}$. For $X\in\{low,\,high\}$ let
\begin{equation}\label{bax5}
\mathfrak{q}_{wa,\pm}^X(\rho,\eta,t):=\pm\widehat{h_{00}^X}(\eta,t)\frac{\Lambda_{wa}(\rho)}{2}+\widehat{h_{0j}^X}(\eta,t)\rho_j\pm\widehat{h_{jk}^X}(\eta,t)\frac{\rho_j\rho_k}{2\Lambda_{wa}(\rho)}.
\end{equation}

To derive our main transport equations we start from the formulas 
\begin{equation}\label{bax4.1}
\partial_tV^{h_{\al\be}}=e^{it\Lambda_{wa}}\mathcal{N}^{h}_{\al\be}=e^{it\Lambda_{wa}}\mathcal{Q}^2_{\al\be}+e^{it\Lambda_{wa}}\mathcal{KG}^2_{\al\be}+e^{it\Lambda_{wa}}\mathcal{S}^2_{\al\be}+e^{it\Lambda_{wa}}\mathcal{N}_{\al\be}^{h,\geq 3},
\end{equation}
see \eqref{sac1.1}. The formulas in the first line of \eqref{on5} show that, for $h\in\{h_{\al\be}\}$,
\begin{equation*}
-\widehat{\partial_j\partial_kh}(\rho)=\rho_j\rho_k\frac{i\widehat{U^{h,+}}(\rho)-i\widehat{U^{h,-}}(\rho)}{2\Lambda_{wa}(\rho)},\qquad \widehat{\partial_0\partial_jh}(\rho)=i\rho_j\frac{\widehat{U^{h,+}}(\rho)+\widehat{U^{h,-}}(\rho)}{2}.
\end{equation*}
Therefore, using \eqref{bax1}, \eqref{bax5}, and the identitities $\widehat{U^{h,\pm}}(\rho,t)=e^{\mp it\Lambda_{wa}(\rho)}\widehat{V^{h,\pm}}(\rho,t)$, we have
\begin{equation}\label{bax4.2}
\widehat{\mathcal{Q}_{\al\be}^2}(\xi,t)=\frac{1}{(2\pi)^3}\sum_{\pm}\int_{\mathbb{R}^3}ie^{\mp it\Lambda_{wa}(\xi-\eta)}\widehat{V^{h_{\al\be},\pm}}(\xi-\eta,t)\mathfrak{q}_{wa,\pm}(\xi-\eta,\eta,t)\,d\eta,
\end{equation}
where $\mathfrak{q}_{wa,\pm}=\mathfrak{q}^{low}_{wa,\pm}+\mathfrak{q}^{high}_{wa,\pm}$. Finally, we notice that
\begin{equation}\label{bax4.3}
\dot{\Theta}_{wa}(\xi,t)=\frac{1}{(2\pi)^3}\int_{\mathbb{R}^3}\mathfrak{q}^{low}_{wa,+}(\xi,\eta,t)e^{i\eta\cdot t\xi/\Lambda_{wa}(\xi)}\,d\eta.
\end{equation} 

Combining \eqref{bax4.1}--\eqref{bax4.3} we derive our main equations for the modified profiles $V^{h_{\al\be}}_\ast$,
\begin{equation}\label{bax6}
\begin{split}
\partial_t\widehat{V^{h_{\al\be}}_\ast}(\xi,t)&=e^{-i\Theta_{wa}(\xi,t)}\{\partial_t\widehat{V^{h_{\al\be}}}(\xi,t)-i\widehat{V^{h_{\al\be}}}(\xi,t)\dot{\Theta}_{wa}(\xi,t)\}=\sum_{a=1}^6\mathcal{R}_a^{h_{\al\be}}(\xi,t),
\end{split}
\end{equation}
where
\begin{equation}\label{bax7}
\mathcal{R}_1^{h_{\al\be}}(\xi,t):=\frac{e^{-i\Theta_{wa}(\xi,t)}}{(2\pi)^3}\int_{\mathbb{R}^3}ie^{it\Lambda_{wa}(\xi)}e^{it\Lambda_{wa}(\xi-\eta)}\widehat{V^{h_{\al\be},-}}(\xi-\eta,t)\mathfrak{q}^{low}_{wa,-}(\xi-\eta,\eta,t)\,d\eta,
\end{equation}
\begin{equation}\label{bax8}
\begin{split}
\mathcal{R}_2^{h_{\al\be}}(\xi,t):=\frac{e^{-i\Theta_{wa}(\xi,t)}}{(2\pi)^3}\int_{\mathbb{R}^3}i\big\{&e^{it(\Lambda_{wa}(\xi)-\Lambda_{wa}(\xi-\eta))}\widehat{V^{h_{\al\be}}}(\xi-\eta,t)\mathfrak{q}_{wa,+}^{low}(\xi-\eta,\eta,t)\\
&-e^{it(\xi\cdot\eta)/\Lambda_{wa}(\xi)}\widehat{V^{h_{\al\be}}}(\xi,t)\mathfrak{q}^{low}_{wa,+}(\xi,\eta,t)\big\}\,d\eta,
\end{split}
\end{equation}
\begin{equation}\label{bax9}
\mathcal{R}_3^{h_{\al\be}}(\xi,t):=\frac{e^{-i\Theta_{wa}(\xi,t)}}{(2\pi)^3}\sum_{\iota\in\{+,-\}}\int_{\mathbb{R}^3}ie^{it\Lambda_{wa}(\xi)}e^{-it\Lambda_{wa,\iota}(\xi-\eta)}\widehat{V^{h_{\al\be},\iota}}(\xi-\eta,t)\mathfrak{q}^{high}_{wa,\iota}(\xi-\eta,\eta,t)\,d\eta,
\end{equation}
\begin{equation}\label{bax10.1}
\mathcal{R}_4^{h_{\al\be}}(\xi,t):=e^{-i\Theta_{wa}(\xi,t)}e^{it\Lambda_{wa}(\xi)}\widehat{\mathcal{KG}^2_{\al\be}}(\xi,t),
\end{equation}
\begin{equation}\label{bax10.2}
\mathcal{R}_5^{h_{\al\be}}(\xi,t):=e^{-i\Theta_{wa}(\xi,t)}e^{it\Lambda_{wa}(\xi)}\widehat{\mathcal{S}^2_{\al\be}}(\xi,t),
\end{equation}
and
\begin{equation}\label{bax10.3}
\mathcal{R}_6^{h_{\al\be}}(\xi,t):=e^{-i\Theta_{wa}(\xi,t)}e^{it\Lambda_{wa}(\xi)}\widehat{\mathcal{N}^{h,\geq 3}_{\al\be}}(\xi,t).
\end{equation}

\subsection{The first reduction}\label{WaRed} We return now to the proof of Proposition \ref{ZWAim}. Since $|\widehat{V^{G}}(\xi,t)|=|\widehat{V^{G}_\ast}(\xi,t)|$, in view of the definition \eqref{bax4} it suffices to prove that
\begin{equation}\label{bax12}
\begin{split}
\|\varphi_k(\xi)\{\widehat{V^{h_{\al\be}}_\ast}(\xi,t_2)-\widehat{V^{h_{\al\be}}_\ast}(\xi,t_1)\}\|_{L^\infty_\xi}&\lesssim \varep_02^{\delta m}2^{-k^--\kappa k^-}2^{-N_0k^+},\\
\|\varphi_k(\xi)\{\widehat{V^{H}_\ast}(\xi,t_2)-\widehat{V^{H}_\ast}(\xi,t_1)\}\|_{L^\infty_\xi}&\lesssim \varep_02^{-\delta m/2}2^{-k^--\kappa k^-}2^{-N_0k^+},
\end{split}
\end{equation}
for any $H\in\{F,\omega_a,\vartheta_{ab}\}$, $k\in\mathbb{Z}$, $m\geq 1$, and $t_1,t_2\in[2^m-2,2^{m+1}]\cap [0,T]$. As before, we show first that the bounds \eqref{bax12} hold if $k$ is too small or if $k$ is too large (relative to $m$). 

\begin{lemma}\label{transferwa} 
The bounds \eqref{bax12} hold if $k\leq -(\delta'/\kappa) m$ or if $k\geq \delta'm-10$.
\end{lemma}

\begin{proof} As in Lemma \ref{transfer}, we use Propositions \ref{EnergEst2} and \ref{DiEs1}, and the inequalities \eqref{consu2.1}. Thus
\begin{equation*}
\begin{split}
\|\widehat{P_kV^{h_{\al\be}}}(t)\|_{L^\infty}&\lesssim \varep_02^{-k-\delta k^-/2}\langle t\rangle^{\delta'/2}2^{-N(0)k^+(1-\delta)/4}2^{-N(2)k^+(3+\delta)/4}\\
&\lesssim \varep_02^{-k^--\delta k^-/2}\langle t\rangle^{\delta'/2}2^{-N_0k^+-dk^+}.
\end{split}
\end{equation*} 
This suffices to prove \eqref{bax12} if $k\leq -(\delta'/\kappa) m$ or if $k\geq \delta'm-10$, as claimed.
\end{proof}

In the remaining range $k\in[- (\delta'/\kappa)m,\delta'm-10]$, we use the identities \eqref{bax6}--\eqref{bax10.3}, so
\begin{equation*}
\widehat{V^{h_{\al\be}}_\ast}(\xi,t_2)-\widehat{V^{h_{\al\be}}_\ast}(\xi,t_1)=\sum_{a=1}^6\int_{t_1}^{t_2}\mathcal{R}_a^{h_{\al\be}}(\xi,s)\,ds.
\end{equation*}
We analyze the contributions of the nonlinear terms $\mathcal{R}_a^{h_{\al\be}}$ separately, and prove the bounds \eqref{bax12}. In fact, in all cases except for the semilinear wave interactions in the terms $\mathcal{R}_5^{h_{\al\be}}$ we can prove stronger bounds, namely
\begin{equation}\label{bax13}
\Big|\varphi_k(\xi)\int_{t_1}^{t_2}\mathcal{R}_a^{h_{\al\be}}(\xi,s)\,ds\Big|\lesssim\varep_1^22^{-\delta m/2}2^{-k^-}2^{-N_0k^+},
\end{equation}
with a gain of a factor of $2^{-\delta m/2}$ instead of a loss. See Lemmas \ref{BAX100}, \ref{BAX101}, \ref{BAX102}, \ref{BAX103}, and \ref{TAR100}. In the case of semilinear wave interactions ($a=5$) we prove slightly weaker bounds, which still suffice for \eqref{bax12}. See Lemma \ref{TAR101}.

\subsection{The nonlinear terms $\mathcal{R}_a^{h_{\al\be}}$, $a\in\{1,2,4,6\}$} In this subsection we consider some of the easier cases, when we can prove the stronger bounds \eqref{bax13}. As in \eqref{Nor53.1}, for $\mu,\nu\in\{(kg,+),(kg,-)\}$ or $\mu,\nu\in\{(wa,+),(wa,-)\}$, and $s\in[0,T]$ we define the operators $T^{wa}_{\mu\nu}$ by
\begin{equation}\label{bax13.1}
T^{wa}_{\mu\nu}[f,g](\xi,s):=\int_{\mathbb{R}^3}\frac{e^{is\Phi_{(wa,+)\mu\nu}(\xi,\eta)}}{\Phi_{(wa,+)\mu\nu}(\xi,\eta)}m(\xi-\eta,\eta)\widehat{f}(\xi-\eta,s)\widehat{g}(\eta,s)\,d\eta,
\end{equation}
where $\Phi_{(wa,+)\mu\nu}(\xi,\eta)=\Lambda_{wa}(\xi)-\Lambda_{\mu}(\xi-\eta)-\Lambda_{\nu}(\eta)$ (see \eqref{on9.2}) and $m\in\mathcal{M}$ (see \eqref{mults}).

\begin{lemma}\label{BAX100} 
The bounds \eqref{bax13} hold for $m\geq 100$, $k\in[-\kappa m,\delta'm]$, and $a=6$.
\end{lemma}

\begin{proof} We use \eqref{consu2}, combined with either \eqref{wer4.0}--\eqref{plk2} or \eqref{abc3.00}--\eqref{abc6.00}. Thus
\begin{equation}\label{bax20}
\begin{split}
\|\widehat{P_k\mathcal{N}^{h,2}_{\al\be}}\|_{L^\infty}&\lesssim \varep_1^2\langle t\rangle^{-1/2+\delta'}2^{-k^-/2}2^{-N(2)k^+},\\
\|\widehat{P_k\mathcal{N}^{h,\geq 3}_{\al\be}}\|_{L^\infty}&\lesssim \varep_1^2\langle t\rangle^{-5/4+6\delta'}2^{-3k^-/4}2^{-N(2)k^+}.
\end{split}
\end{equation}
The bounds in the second line suffice to prove \eqref{bax13} for $a=6$.
\end{proof}

\begin{lemma}\label{BAX101} 
The bounds \eqref{bax13} hold for $m\geq 100$, $k\in[-\kappa m,\delta'm]$, and $a=1$.
\end{lemma}

\begin{proof} This is similar to the proof of Lemma \ref{BootstrapZ3} and the main point is that the interaction is nonresonant so we can integrate by parts in time. We use the formula \eqref{bax5}, substitute $h=i\Lambda_{wa}^{-1}(U^{h,+}-U^{h,-})/2$, $h\in \{h_{\al\be}\}$ (see \eqref{on5}) and decompose all the input functions dyadically in frequency. With $U_{low}^{h,\pm}$ and $V_{low}^{h,\pm}$ defined as in \eqref {ulow}, it suffices to prove that
\begin{equation}\label{bax21}
\begin{split}
\sum_{(k_1,k_2)\in\mathcal{X}_k}2^{k_1-k_2}\Big\|\varphi_k(\xi)\int_{t_1}^{t_2}e^{is\Lambda_{wa}(\xi)-i\Theta_{wa}(\xi,s)}&\mathcal{F}\{I[P_{k_1}U^{h_1,-},P_{k_2}U^{h_2,\iota_2}_{low}]\}(\xi,s)\,ds\Big\|_{L^\infty_\xi}\lesssim \varep_1^22^{-\kappa m},
\end{split}
\end{equation}
for $h_1,h_2\in\{h_{\al\be}\}$ and $\iota_2\in\{+,-\}$.

We estimate first, using just \eqref{vcx1.2} and \eqref{vcx1},
\begin{equation}\label{bax21.1}
\begin{split}
\big\|\mathcal{F}\{I[P_{k_1}U^{h_1,-},&P_{k_2}U^{h_2,\iota_2}_{low}]\}(\xi,s)\big\|_{L^\infty_\xi}\lesssim \|\widehat{P_{k_1}U^{h_1,-}}\|_{L^\infty}\|\widehat{P_{k_2}U^{h_2,\iota_2}_{low}}\|_{L^1}\\
&\lesssim \varep_1^22^{-k_1^--\kappa k_1^-}2^{2k_2-\delta k_2}2^{2\delta m}.
\end{split}
\end{equation} 
This suffices to control the contribution of the pairs $(k_1,k_2)$ for which $k_2\leq -1.01m$. After this reduction it remains to prove that
\begin{equation}\label{bax22}
\begin{split}
\Big|\int_{t_1}^{t_2}&\int_{\mathbb{R}^3}e^{is\Lambda_{wa}(\xi)-i\Theta_{wa}(\xi,s)}m(\xi-\eta,\eta)e^{is\Lambda_{wa}(\xi-\eta)}\widehat{P_{k_1}V^{h_1,-}}(\xi-\eta,s)\\
&\times e^{-is\Lambda_{wa,\iota_2}(\eta)}\widehat{P_{k_2}V^{h_2,\iota_2}_{low}}(\eta,s)\, d\eta ds\Big|\lesssim \varep_1^22^{-0.01m}2^{k_2},
\end{split}
\end{equation}
for any $\xi$ with $|\xi|\in[2^{k_1-4},2^{k_1+4}]$, provided that
\begin{equation}\label{bax23}
k_2\in[-1.01m,-p_0m+10],\qquad k_1\in[-\kappa m-10,\delta'm+10]. 
\end{equation}

To prove \eqref{bax22} we integrate by parts in time. Notice that $\Lambda_{wa}(\xi)+\Lambda_{wa}(\xi-\eta)-\Lambda_{wa,\iota_2}(\eta)\gtrsim 2^{k_1^-}$ in the support of the integral. The left-hand side of \eqref{bax22} is dominated by $C(I_{wa}+II_{wa}+III_{wa})(\xi)$, where, with $\mu=(wa,-)$ and $\nu=(wa,\iota_2)$ and $T^{wa}_{\mu\nu}$ defined as in \eqref{bax13.1},
\begin{equation}\label{bax25}
\begin{split}
I_{wa}(\xi)&:=\Big(1+\int_{t_1}^{t_2}|\dot{\Theta}_{wa}(\xi,s)|\,ds\Big)\sup_{s\in[t_1,t_2]}|T_{\mu\nu}^{wa}[P_{k_1}V^{h_1,-},P_{k_2}V^{h_2,\iota_2}_{low}](\xi,s)|,\\
II_{wa}(\xi)&:=\int_{t_1}^{t_2}|T_{\mu\nu}^{wa}[\partial_s(P_{k_1}V^{h_1,-}),P_{k_2}V^{h_2,\iota_2}_{low}](\xi,s)|\,ds,\\
III_{wa}(\xi)&:=\int_{t_1}^{t_2}|T_{\mu\nu}^{wa}[P_{k_1}V^{h_1,-},\partial_s(P_{k_2}V^{h_2,\iota_2}_{low})](\xi,s)|\,ds.
\end{split}
\end{equation}

As in \eqref{bax21.1}, assuming $k_1,k_2$ as in \eqref{bax23}, we estimate for any $s\in[t_1,t_2]$ 
\begin{equation*}
|T_{\mu\nu}^{wa}[P_{k_1}V^{h_1,-},P_{k_2}V^{h_2,\iota_2}_{low}](\xi,s)|\lesssim\varep_1^22^{-2k_1^--\kappa k_1^-}2^{2k_2}2^{4\delta m}.
\end{equation*}
The definition \eqref{bax3} shows that $|\dot{\Theta}_{wa}(\xi,s)|\lesssim 2^{k_1^+}\sup_{\al,\be}\|h_{\al\be}(s)\|_{L^\infty}\lesssim 2^{-m+2\delta'm}$. Therefore 
\begin{equation}\label{bax26}
I_{wa}(\xi)\lesssim \varep_1^22^{2k_2}2^{0.01m}.
\end{equation}
Similarly using \eqref{bax20} and the bounds $\|\widehat{P_{k_2}V^{h_2,\iota_2}_{low}}(s)\|_{L^1}\lesssim\varep_12^{2k_2}2^{4\delta m}$ it follows that
\begin{equation*}
II_{wa}(\xi)\lesssim \varep_1^22^{m/2}2^{2k_2}2^{2\kappa m}.
\end{equation*}
Finally, using \eqref{Nor56.2}, we see that $III_{wa}(\xi)\lesssim \varep_1^22^{2k_2}2^{0.01m}$, as in \eqref{bax26}. The desired conclusion \eqref{bax22} follows.
\end{proof}

\begin{lemma}\label{BAX102} 
The bounds \eqref{bax13} hold for $m\geq 100$, $k\in[-\kappa m,\delta'm]$, and $a=2$.
\end{lemma} 

\begin{proof} This is similar to the proof of Lemma \ref{BootstrapZ4}. With $h_1=h_{\al\be}$ we decompose $V^{h_1}=\sum_{(k_1,j_1)\in\mathcal{J}}V_{j_1,k_1}^{h_1,+}$ as in \eqref{on11.3} and examine the definition \eqref{bax8}. It suffices to prove that
\begin{equation*}
\begin{split}
\varphi_k(\xi)&\Big|\int_{\mathbb{R}^3}\big\{e^{it(\Lambda_{wa}(\xi)-\Lambda_{wa}(\xi-\eta))}\widehat{V_{j_1,k_1}^{h_1,+}}(\xi-\eta,t)\mathfrak{q}_{wa,+}^{low}(\xi-\eta,\eta,t)\\
&-e^{it(\xi\cdot\eta)/\Lambda_{wa}(\xi)}\widehat{V^{h_1,+}_{j_1,k_1}}(\xi,t)\mathfrak{q}^{low}_{wa,+}(\xi,\eta,t)\big\}\,d\eta\Big|\lesssim \varep_1^22^{-m-\kappa m}2^{-\delta j_1},
\end{split}
\end{equation*}
provided that $|k_1-k|\leq 10$ and $t\in[2^{m-1},2^{m+1}]$. Using also the definitions \eqref{bax5} it suffices to prove that for any multiplier $m\in\mathcal{M}_0$ (see \eqref{mults0}) and $\al,\be\in\{0,1,2,3\}$ we have
\begin{equation}\label{bax28}
\begin{split}
\varphi_k(\xi)&\Big|\int_{\mathbb{R}^3}\widehat{h_{\al\be}^{low}}(\eta,t)\big\{e^{it(\Lambda_{wa}(\xi)-\Lambda_{wa}(\xi-\eta))}\widehat{V_{j_1,k_1}^{h_1,+}}(\xi-\eta,t)m(\xi-\eta)|\xi-\eta|\\
&-e^{it(\xi\cdot\eta)/\Lambda_{wa}(\xi)}\widehat{V^{h_1,+}_{j_1,k_1}}(\xi,t)m(\xi)|\xi|\big\}\,d\eta\Big|\lesssim \varep_1^22^{-m-\kappa m}2^{-\delta j_1}.
\end{split}
\end{equation}

Recall that $\|\widehat{V_{j_1,k_1}^{h_1,+}}(t)\|_{L^\infty}\lesssim \varep_12^{-3k_1/2}2^{-j_1/2+\delta j_1}2^{\delta'm}2^{-4k_1^+}$, see \eqref{vcx1.15}. Thus, without using the cancellation of the two terms in the integral, the left-hand side of \eqref{bax28} is bounded by
\begin{equation*}
C\varep_12^{-j_1/2+\delta j_1}2^{2\kappa m}\|\widehat{h_{\al\be}^{low}}(t)\|_{L^1}\lesssim \varep_1^22^{-j_1/2+\delta j_1}2^{3\kappa m}2^{-p_0m}.
\end{equation*}
This suffices to prove \eqref{Nor61} when $j_1$ is large, i.e. $2^{j_1/2}\gtrsim 2^{(1.01-p_0)m}$. 

On the other hand, if $2^m\gg 1$ and $j_1/2\leq (1.01-p_0)m=0.33 m$ then we estimate
\begin{equation}\label{bax30}
\begin{split}
\big|e^{it(\Lambda_{wa}(\xi)-\Lambda_{wa}(\xi-\eta))}-e^{it(\xi\cdot\eta)/\Lambda_{wa}(\xi)}\big|&\lesssim 2^{-2p_0m+m}2^{2\kappa m},\\
\big|\widehat{V_{j_1,k_1}^{h_1,+}}(\xi-\eta,t)m(\xi-\eta)|\xi-\eta|-\widehat{V_{j_1,k_1}^{h_1,+}}(\xi,t)m(\xi)|\xi|\big|&\lesssim \varep_12^{j_1/2}2^{3\kappa m}2^{-p_0m},
\end{split}
\end{equation}
provided that $|\xi|\approx 2^k$ and $|\eta|\lesssim 2^{-p_0m}$. Indeed, the first bound follows from the observation that $\nabla\Lambda_{wa}(\xi)=\xi/\Lambda_{wa}(\xi)$. The second bound follows from \eqref{vcx1.15} once we notice that taking $\partial_\xi$ derivative of the localized profiles $\widehat{V_{j_1,k_1}^{h_1,+}}$ corresponds essentially to multiplication by a factor of $2^{j_1}$. If $j_1/2\leq 0.33 m$ it follows that the left-hand side of \eqref{bax28} is bounded by
\begin{equation*}
C\varep_12^{-0.34 m}\|\widehat{h_{\al\be}^{low}}(t)\|_{L^1}\lesssim \varep_1^22^{-p_0m-0.34 m+\kappa m}.
\end{equation*}
This suffices to prove \eqref{bax28} when $j_1/2\leq 0.33 m$, which completes the proof of the lemma.
\end{proof}

\begin{lemma}\label{BAX103} 
The bounds \eqref{bax13} hold for $m\geq 100$, $k\in[-\kappa m,\delta'm]$, and $a=4$.
\end{lemma}

\begin{proof} Here we analyze bilinear interactions of the Klein-Gordon field with itself. The important observation is that these interactions are still nonresonant, due to Lemma \ref{pha2}. Recall that $\mathcal{KG}^2_{\al\be}=2\partial_\al\psi\partial_\be\psi+\psi^2m_{\al\be}$. We express $\psi,\partial_\al\psi$ in terms of the normalized profiles $U^{\psi,\pm}$ and decompose dyadically in frequency. It suffices to prove that
\begin{equation}\label{bax41}
\begin{split}
\Big\|\varphi_k(\xi)\int_{t_1}^{t_2}e^{is\Lambda_{wa}(\xi)-i\Theta_{wa}(\xi,s)}&\mathcal{F}\{I[P_{k_1}U^{\psi,\iota_1},P_{k_2}U^{\psi,\iota_2}]\}(\xi,s)\,ds\Big\|_{L^\infty_\xi}\lesssim \varep_1^22^{-50\delta' m},
\end{split}
\end{equation}
for any $(k_1,k_2)\in\mathcal{X}_k$ and $\iota_1,\iota_2\in\{+,-\}$.

The estimates \eqref{bax41} follow easily using just the $L^2$ bounds \eqref{cnb2} when $\min(k_1,k_2)\leq-3m/5$ or $\max(k_1,k_2)\geq m/20$. In the remaining range  we integrate by parts in time. As before, notice that $|\dot{\Theta}_{wa}(\xi,s)|\lesssim 2^{-m+2\delta'm}$. With $T_{\mu\nu}^{wa}$ as in \eqref{bax13.1}, it suffices to prove that
\begin{equation}\label{bax42}
\big|\varphi_k(\xi)T^{wa}_{\mu\nu}[P_{k_1}V^{\psi,\iota_1},P_{k_2}V^{\psi,\iota_2}](\xi,s)\big|\lesssim \varep_1^22^{-60 \delta' m},
\end{equation} 
and
\begin{equation}\label{bax42.5}
2^m\big|\varphi_k(\xi)T^{wa}_{\mu\nu}[\partial_s(P_{k_1}V^{\psi,\iota_1}),P_{k_2}V^{\psi,\iota_2}](\xi,s)\big|\lesssim \varep_1^22^{-60\delta' m},
\end{equation}
where $\mu=(kg,\iota_1)$, $\nu=(kg,\iota_2)$, $s\in[2^{m-1},2^{m+1}]$, and $k_1,k_2\in[-3m/5,m/20]$.

{\bf{Step 1: proof of \eqref{bax42}.}} We decompose $P_{k_1}V^{\psi,\iota_1}=\sum_{j_1\geq -k_1^-}V^{\psi,\iota_1}_{j_1,k_1}$ and $P_{k_2}V^{\psi,\iota_2}=\sum_{j_2\geq -k_2^-}V^{\psi,\iota_2}_{j_2,k_2}$ and in \eqref{on11.3}. The contribution of the pairs $(V^{\psi,\iota_1}_{j_1,k_1},V^{\psi,\iota_2}_{j_2,k_2})$ with $\max(j_1,j_2)\geq 0.01 m$ can be estimated easily, using the observation that $|\Phi_{(wa,+)\mu\nu}(\xi,\eta)|\gtrsim 2^{k^-}2^{-2\max(k_1^+,k_2^+)}$ in the support of the integral (see Lemma \ref{pha2} (i)) and the $L^2$ bounds \eqref{vcx1.1}.

On the other hand, if $\max(j_1,j_2)\leq 0.1 m$ then we have to show that
\begin{equation}\label{bax45}
\begin{split}
\Big|\varphi_k(\xi)\int_{\mathbb{R}^3}\frac{e^{-is[\Lambda_{\mu}(\xi-\eta)+\Lambda_{\nu}(\eta)]}}{\Phi_{(wa,+)\mu\nu}(\xi,\eta)}m(\xi-\eta,\eta)\widehat{V^{\psi,\iota_1}_{j_1,k_1}}(\xi-\eta,s)\widehat{V^{\psi,\iota_2}_{j_2,k_2}}(\eta,s)\,d\eta\Big|\lesssim \varep_1^22^{-\kappa m}.
\end{split}
\end{equation}
We observe that
\begin{equation}\label{bax44}
|\nabla\Lambda_{kg}(x)-\nabla\Lambda_{kg}(y)|\gtrsim |x-y|/(1+|x|^4+|y|^4)\qquad \text{ for any } x,y\in\mathbb{R}^3.
\end{equation}  
Therefore, using integration by parts in $\eta$ (Lemma \ref{tech5}), the left-hand side of \eqref{bax45} is $\lesssim C\varep_1^22^{-2m}$ if $\mu=-\nu$. On the other hand, if $\mu=\nu$ then the only space-resonant point (where the gradient of the phase vanishes)  is $\eta=\xi/2$ and we insert cutoff functions of the form $\varphi_{\leq 0}(2^{0.4m}(\eta-\xi/2))$ and $\varphi_{>1}(2^{0.4m}(\eta-\xi/2))$. Then we estimate the integral corresponding to $|\eta-\xi/2|\lesssim 2^{-0.4m}$ by $C\varep_1^22^{-1.1m}$, by placing the profiles $\mathcal{F}V^{\psi,\iota_1}_{j_1,k_1}$ and $\mathcal{F}V^{\psi,\iota_2}_{j_2,k_2}$ in $L^\infty$. Finally, we estimate the integral corresponding to $|\eta-\xi/2|\gtrsim 2^{-0.4m}$ by $C\varep_1^22^{-2m}$, using integration by parts in $\eta$ and \eqref{bax44}. This completes the proof of \eqref{bax45}.

{\bf{Step 2: proof of \eqref{bax42.5}.}} Recall that $\partial_tV^{\psi,\iota_1}(t)=e^{it\Lambda_{kg,\iota_1}}\mathcal{N}^{\psi,2}(t)+e^{it\Lambda_{kg,\iota_1}}\mathcal{N}^{\psi,\geq 3}(t)$, see \eqref{loc1}. The contribution of the cubic and higher order nonlinearity $\mathcal{N}^{\psi,\geq 3}$ is easy to estimate, using just the $L^2$ bounds \eqref{abc31.00} and the lower bounds $|\Phi_{(wa,+)\mu\nu}(\xi,\eta)|\gtrsim 2^{k^-}2^{-2\max(k_1^+,k_2^+)}$, which hold in the support of the operator.

To bound the contribution of $\mathcal{N}^{\psi,2}$ we define the trilinear operators
\begin{equation}\label{bax50}
\begin{split}
\mathcal{C}_{wa}[f,g,h](\xi):=&\int_{\mathbb{R}^3\times\mathbb{R}^3}\frac{\varphi_k(\xi)\varphi_{k_1}(\xi-\eta)m(\xi-\eta,\eta)}{\Lambda_{wa}(\xi)-\Lambda_\mu(\xi-\eta)-\Lambda_\nu(\eta)}\\
&\times m_3(\xi-\eta-\rho)m_4(\rho)\cdot\widehat{f}(\xi-\eta-\rho)\widehat{g}(\eta)\widehat{h}(\rho)\,d\eta d\rho,
\end{split}
\end{equation}
where $m\in\mathcal{M}$, $m_3,m_4\in\mathcal{M}_0$. For \eqref{bax42.5} it remains to prove that
\begin{equation}\label{bax51}
2^m\sum_{k_3,k_4\in\mathbb{Z}}2^{k_3^+-k_4}\big|\mathcal{C}_{wa}[P_{k_3}U^{\psi,\iota_3},P_{k_2}U^{\psi,\iota_2},P_{k_4}U^{h,\iota_4}](\xi,s)\big|\lesssim \varep_1^22^{-60\delta' m},
\end{equation}
where $\mu=(kg,\iota_1)$, $\nu=(kg,\iota_2)$, $s\in[2^{m-1},2^{m+1}]$, $h\in\{h_{\al\be}\}$, and $k_1,k_2\in[-3m/5,m/20]$.

Using Lemma \ref{pha2} (ii) and the bounds \eqref{cnb2} and \eqref{wws1} we estimate
\begin{equation}\label{bax55}
\begin{split}
\big|\mathcal{C}_{wa}[P_{k_3}U^{\psi,\iota_3},&P_{k_2}U^{\psi,\iota_2},P_{k_4}U^{h,\iota_4}](\xi,s)\big|\lesssim 2^{-k+6\max(k^+,k_2^+)}\|P_{k_3}U^{\psi}\|_{L^2}\|P_{k_2}U^{\psi}\|_{L^2}\|P_{k_4}U^{h}\|_{L^\infty}\\
&\lesssim \varep_1^32^{2\delta'm}2^{k_2^-}2^{k_3^-}2^{k_4^-(1-\delta)}\min(2^{-m},2^{k_4^-})2^{-20\max(k_2^+,k_3^+,k_4^+)}.
\end{split}
\end{equation}
This suffices to bound the contribution of the triplets $(k_2,k_3,k_4)$ for which $k_4\leq -1.01m$, or $k_3\leq -0.01 m$, or $\max(k_2,k_3,k_4)\geq 10\delta'm$. In the remaining range we further decompose $P_{k_3}U^{\psi,\iota_3}=\sum_{j_3}e^{-it\Lambda_{kg,\iota_3}}V^{\psi,\iota_3}_{j_3,k_3}$ and $P_{k_2}U^{\psi,\iota_2}=\sum_{j_2}e^{-it\Lambda_{kg,\iota_2}}V^{\psi,\iota_2}_{j_2,k_2}$. Notice that the contribution of the pairs $(j_2,j_3)$ for which $\max(j_2,j_3)\geq 0.1m$ can be suitably bounded, using an estimate similar to \eqref{bax55}. For \eqref{bax51} it remains to prove that
\begin{equation}\label{bax56}
\big|\mathcal{C}_{wa}[e^{-it\Lambda_{kg,\iota_3}}V^{\psi,\iota_3}_{j_3,k_3},e^{-it\Lambda_{kg,\iota_2}}V^{\psi,\iota_2}_{j_2,k_2},P_{k_4}U^{h,\iota_4}](\xi,s)\big|\lesssim \varep_1^32^{-1.01m}2^{k_4},
\end{equation} 
provided that $k_2,k_3\in[-0.1m,10\delta'm]$, $k_4\in[-1.01 m,10\delta'm]$, and $j_2,j_3\leq 0.1m$.

To prove this, we insert cutoff functions of the form $\varphi_{\leq 0}(2^{0.35m}(\rho-\xi))$ and $\varphi_{>1}(2^{0.35m}(\rho-\xi))$ in the integral in \eqref{bax50}. The contribution of the integral when $|\rho-\xi|\lesssim 2^{-0.35m}$ (which is nontrivial only if $2^{k_4}\gtrsim 2^{-\kappa m}$) is bounded as claimed by estimating in the Fourier space, with $\widehat{V^{\psi,\iota_3}_{j_3,k_3}}$ and $\widehat{V^{\psi,\iota_2}_{j_2,k_2}}$ estimated $L^2$ and $\widehat{P_{k_4}U^{h,\iota_4}}$ estimated in $L^\infty$. 

On the other hand, the integral when $|\rho-\xi|\gtrsim 2^{-0.35m}$ can be estimated as in the proof of \eqref{bax45}. Using \eqref{bax44}, the $\eta$ integral is bounded by $C\varep_1^22^{-1.1m}$, for any $\xi,\rho\in\mathbb{R}^3$. Then we notice that $\|\widehat{P_{k_4}U^{h,\iota_4}}(\rho,s)\|_{L^1_\rho}\lesssim\varep_12^{2k_4}2^{\delta'm}$, and the desired conclusion \eqref{bax56} follows. This completes the proof of the lemma.
\end{proof} 

\subsection{Localized bilinear wave interactions} In this subsection we start analyzing the remaining cases, where we have bilinear interactions of the metric components. These cases are more difficult because of the presence of time-resonant frequencies (parallel bilinear interactions), which prevent direct integration by parts in time. 

For $b\in\mathbb{Z}$, $\xi\in\mathbb{R}^3$, and multipliers $m\in\mathcal{M}$ we define the bilinear operators
\begin{equation}\label{baz1}
J_b[f,g](\xi)=J_{b;\iota_1\iota_2}[f,g](\xi):=\int_{\mathbb{R}^3}m(\xi-\eta,\eta)\widehat{f}(\xi-\eta)\widehat{g}(\eta)\varphi_b\big(\Xi_{\iota_1\iota_2}(\xi-\eta,\eta)\big)\,d\eta,
\end{equation}
for $\iota_1,\iota_2\in\{+,-\}$ (see the definition \eqref{par1}). As in the proof of Lemma \ref{PhaWave}, we remark that an expression of the form $\varphi_b\big(\Xi_{\iota_1\iota_2}(\xi-\eta,\eta)\big)\cdot\varphi_k(\xi)\varphi_{k_1}(\xi-\eta)\varphi_{k_2}(\eta)$ can be nontrivial only if either $b\geq -20$ or $2^b\lesssim 2^{k-\max(k_1,k_2)}$.   

We start with a lemma.

\begin{lemma}\label{gene1} Assume $m\geq 10$, $t\in[2^{m-1},2^{m+1}]\cap[0,T]$, $l,l_1,l_2\in\mathbb{Z}\cap [-m,m/5+10]$. 

(i) If $b\in[-m,2]$, $(l_1,j_1),(l_2,j_2)\in\mathcal{J}$, and $n\in L^\infty(\mathbb{R}^3\times\mathbb{R}^3)$ then
\begin{equation}\label{gene2}
\begin{split}
\Big|\varphi_l(\xi)&\int_{\mathbb{R}^3}n(\xi,\eta)\widehat{V^{h_1,\iota_1}_{j_1,l_1}}(\xi-\eta,t)\widehat{V^{h_2,\iota_2}_{j_2,l_2}}(\eta,t)\varphi_{\leq b}(\Xi_{\iota_1\iota_2}(\xi-\eta,\eta))\,d\eta\Big|\\
&\lesssim \varep_1^2\|n\|_{L^\infty}\cdot 2^{\delta'm}\min(2^{-2\max(l_1,l_2)},2^{2b-2l})2^{l_1+l_2}2^{-\max(j_1,j_2)}2^{-20(l_1^++l_2^+)},
\end{split}
\end{equation}
for any $\xi\in\mathbb{R}^3$, $h_1,h_2\in\{h_{\al\be}\}$, and $\iota_1,\iota_2\in\{+,-\}$. 

(ii) If $b\geq (-m+l-l_1-l_2)/2+\delta m/8$ then, for any $\xi\in\mathbb{R}^3$, $h_1,h_2\in\{h_{\al\be}\}$, and $\iota_1,\iota_2\in\{+,-\}$,
\begin{equation}\label{gene3}
\begin{split}
\sum_{j_1\geq-l_1^-,\,j_2\geq -l_2^-}&\big|\varphi_l(\xi)J_b[U^{h_1,\iota_1}_{j_1,l_1}(t),U^{h_2,\iota_2}_{j_2,l_2}(t)]\big|\\
&\lesssim \varep_1^22^{2\delta'm-m}\min(2^{-2\max(l_1,l_2)},2^{2b-2l})2^{-b+l_1+l_2}2^{-18(l_1^++l_2^+)}.
\end{split}
\end{equation}
As a consequence, for any $\xi\in\mathbb{R}^3$,
\begin{equation}\label{gene3x}
\sum_{b\leq 2}\sum_{j_1\geq-l_1^-,\,j_2\geq -l_2^-}\big|\varphi_l(\xi)J_b[U^{h_1,\iota_1}_{j_1,l_1}(t),U^{h_2,\iota_2}_{j_2,l_2}(t)]\big|\lesssim \varep_1^22^{2\delta'm-m}2^{\min(l_1,l_2)-l}2^{-18(l_1^++l_2^+)}.
\end{equation}
\end{lemma}

\begin{proof} (i) We may assume that $\|n\|_{L^\infty}=1$. Without loss of generality, by rotation, we may also assume that $j_2\geq j_1$ and $\xi=(\xi_1,0,0)$, $\xi_1\in[2^{l-1},2^{l+1}]$.  We notice that the $\eta$ integral in the left-hand side is supported in the set $\mathcal{R}_{\leq b;\xi}:=\big\{|\eta|\approx 2^{l_2},\,|\xi-\eta|\approx 2^{l_1},\,\sqrt{\eta_2^2+\eta_3^2}\lesssim X:=\min(2^{l_1},2^{l_2},2^{b+l_1+l_2-l})\big\}$ (the last bound on $\sqrt{\eta_2^2+\eta_3^2}$ holds when $b\leq -20$, see \eqref{par9.1}).

With $\mathbf{e}_1:=(1,0,0)$, we estimate the left-hand side of \eqref{gene2} by
\begin{equation*}
\begin{split}
C\|\widehat{V^{h_1,\iota_1}_{j_1,l_1}}&(t)\|_{L^\infty}\int_{\mathbb{R}^3}|\widehat{V^{h_2,\iota_2}_{j_2,l_2}}(\eta,t)|\mathbf{1}_{\mathcal{R}_{\leq b;\xi}}(\eta)\,d\eta\\
&\lesssim\|\widehat{V^{h_1,\iota_1}_{j_1,l_1}}(t)\|_{L^\infty}\int_{[0,\infty)\times{\mathbb{S}}^2,\,|\theta-\mathbf{e}_1|\lesssim 2^{-l_2}X}|\widehat{V^{h_2,\iota_2}_{j_2,l_2}}(r\theta,t)|\,r^2 dr d\theta.
\end{split}
\end{equation*}
Using \eqref{vcx1.3*}, \eqref{Linfty3.34} (with $p=1/\delta$ large) and \eqref{vcx1.1}, we can further estimate the right-hand side of the expression above by
\begin{equation*}
\begin{split}
C\varep_12^{\delta'm/2}2^{-22l_1^+}2^{-l_1^-}&\|\widehat{V^{h_2,\iota_2}_{j_2,l_2}}(r\theta,s)\|_{L^2(r^2dr)L^p_\theta}\cdot (2^{-l_2}X)^{2/p'}2^{3l_2/2}\\
&\lesssim\varep_1^22^{\delta'm}2^{-20(l_1^++l_2^+)}2^{-j_2}2^{-l_1-l_2}\min\{2^{\min(l_1,l_2)},2^{b+l_1+l_2-l}\}^2.
\end{split}
\end{equation*}
The bound \eqref{gene2} follows. 

(ii) Notice that \eqref{gene3x} follows from \eqref{gene2} and \eqref{gene3} by summation over $b$. To prove \eqref{gene3} we notice first that the contribution of the pairs $(j_1,j_2)$ for which $\max(j_1,j_2)\geq m+b-\delta m$ is bounded as claimed, as a consequence of \eqref{gene2}. 

We claim that the contribution of the remaining pairs $(j_1,j_2)$ is negligible, i.e.
\begin{equation}\label{baz3.1}
\big|\varphi_l(\xi)J_b[e^{-it\Lambda_{wa,\iota_1}}V^{h_1,\iota_1}_{j_1,l_1}(t),e^{-it\Lambda_{wa,\iota_2}}V^{h_2,\iota_2}_{j_2,l_2}(t)]\lesssim \varep_1^22^{-4m}2^{-20(l_1^++l_2^+)},
\end{equation}
if $\max(j_1,j_2)\leq m+b-\delta m$. For this we would like to use Lemma \ref{tech5}. Notice that, in the support of the integral, we always have the lower bounds
\begin{equation}\label{baz4}
\big|\nabla_\eta[\Lambda_{wa,\iota_1}(\xi-\eta)+\Lambda_{wa,\iota_2}(\eta)]\big|=\big|\Xi_{\iota_1\iota_2}(\xi-\eta,\eta)\big|\gtrsim 2^{b}.
\end{equation}

We would like to use Lemma \ref{tech5} with $K\approx 2^{m+b}$ and $\eps\approx 2^{\delta m/8}/K$. As in the proof of Lemma \ref{PhaWave} let $H_{b;\xi}(\eta)=2^{-2b}|\Xi_{\iota_1\iota_2}(\xi-\eta,\eta)|^2$ such that $\varphi_b(\Xi_{\iota_1\iota_2}(\xi-\eta,\eta))=\varphi''_0(H_{b;\xi}(\eta))$, $\varphi''_0(x):=\mathbf{1}_{[0,\infty)}(x)\varphi_0(\sqrt{x})$. If $2^b\gtrsim 1$ then the function $H_{b;\xi}$ satisfies differential bounds of the form $|D^\al_{\eta}H_{b;\xi}(\eta)|\lesssim_{|\alpha|}2^{|\alpha|\max(-l_1,-l_2)}$ in the support of the integral, for all multi-indices $\alpha\in\mathbb{Z}_+^3$, and the desired conclusion \eqref{baz3.1} follows from Lemma \ref{tech5}. 

Assume now that $b\leq -20$ (so $2^b\lesssim 2^{k-\max(k_1,k_2)}$) and, as before, $\xi=(\xi_1,0,0)$, $\xi_1>0$. The formula \eqref{par9.1} shows that
\begin{equation}\label{baz2.1}
|D_{\eta}^{\al}H_{b;\xi}(\eta)|\lesssim_{|\al|} 2^{-|\al|(b+l_1+l_2-l)},\qquad\al\in\mathbb{Z}_+^3
\end{equation}
for $\eta\in\mathcal{R}_{b;\xi}=\big\{|\eta|\approx 2^{l_2},\,|\xi-\eta|\approx 2^{l_1},\,\sqrt{\eta_2^2+\eta_3^2}\approx 2^{b+l_1+l_2-l}\big\}$. Notice that $\eps 2^{\max(j_1,j_2)}+\eps 2^{-(b+l_1+l_2-l)}\lesssim 2^{-\delta m/8}$, due to the assumptions $2b+m+l_1+l_2-l\geq\delta m/4$ and $2^{\max(j_1,j_2)}\lesssim K2^{-\delta m}$. The desired bounds \eqref{baz3.1} would follow from Lemma \ref{tech5} if we could verify the second bound in \eqref{ln2}. With $f:=K^{-1}s[\Lambda_{wa,\iota_1}(\xi-\eta)+\Lambda_{wa,\iota_2}(\eta)]$, we always have 
\begin{equation}\label{baz2.2}
|D^\al f(\eta)|\lesssim_{|\al|} K^{-1}2^m2^{-(|\al|-1)\min(l_1,l_2)}
\end{equation} 
in the support of the integral. Since $2^{-\min(l_1,l_2)}\lesssim K2^{-\delta m/8}$, it suffices to verify the bounds \eqref{ln2} when $|\al|=2$, i.e. $K^{-1}2^m2^{-\min(l_1,l_2)}\lesssim K2^{-\delta m/8}$. In view of the assumption $2b\geq -m-l_1-l_2+l+\delta m$, this holds when $2^{l-\max(l_1,l_2)}\gtrsim 1$, and the desired conclusion \eqref{baz3.1} follows in this case.

If $l\leq \max(l_1,l_2)-40$ then we need to be slightly more careful with the estimates \eqref{baz2.2}. Since $b\leq -20$, we may assume that $b\leq l-\max(l_1,l_2)+10$. We may also assume that $\iota_1=-\iota_2$, since otherwise $\varphi_b(\Xi_{\iota_1\iota_2}(\xi-\eta,\eta))\equiv 0$. We define $K,\eps,f$ as before and notice that the bounds \eqref{baz2.2} can be improved to
\begin{equation*}
|D^\al f(\eta)|\lesssim_{|\al|} K^{-1}2^m2^{l-\max(l_1,l_2)}2^{-(|\al|-1)\max(l_1,l_2)},
\end{equation*}
in the support of the integral. This suffices to verify the bounds \eqref{ln2} in Lemma \ref{tech5} in the remaining case, and completes the proof of \eqref{baz3.1}. 
\end{proof}

Our main result in this subsection is the following lemma, in which we show that the contribution of non-parallel wave interactions is suitably small.

\begin{lemma}\label{LemBaz1}
Assume that $m\geq 100$, $t_1,t_2\in[2^{m-1},2^{m+1}]\cap [0,T]$, $k,k_1,k_2\in\mathbb{Z}$, $k\in[-\kappa m/4,\delta'm]$, $-p_0m-10\leq k_2\leq k_1\leq m/10$, and $q\geq (-m+k-k_1-k_2)/2+\delta m/8$. Then
\begin{equation}\label{baz2}
\begin{split}
\Big\|\varphi_k(\xi)\int_{t_1}^{t_2}e^{is\Lambda_{wa}(\xi)-i\Theta_{wa}(\xi,s)}J_q[P_{k_1}U^{h_1,\iota_1},P_{k_2}U^{h_2,\iota_2}](\xi,s)\,ds\Big\|_{L^\infty_\xi}\lesssim \varep_1^22^{-0.001m}2^{k_2^-},
\end{split}
\end{equation}
for any $h_1,h_2\in\{h_{\al\be}\}$ and $\iota_1,\iota_2\in\{+,-\}$.
\end{lemma}

\begin{proof} We notice that the desired bounds follow directly from \eqref{gene3} if $q\leq -0.002m+10$ (recall that $k,k_1\geq -\kappa m/4-10$). On the other hand, if
\begin{equation}\label{baz8}
q\in [-0.002m+10,2]
\end{equation}
then we integrate by parts in time. Notice that
\begin{equation}\label{baz9}
\big|\Lambda_{wa}(\xi)-\Lambda_{wa,\iota_1}(\xi-\eta)-\Lambda_{wa,\iota_2}(\eta)\big|^{-1}\lesssim 2^{-2q-k_2}
\end{equation}
in the support of the integral, as a consequence of \eqref{par5}. We define the operators $T^{wa}_{\mu\nu;q}$ by
\begin{equation}\label{baz10}
T^{wa}_{\mu\nu;q}[f,g](\xi,s):=\int_{\mathbb{R}^3}\frac{e^{is\Phi_{(wa,+)\mu\nu}(\xi,\eta)}}{\Phi_{(wa,+)\mu\nu}(\xi,\eta)}m(\xi-\eta,\eta)\widehat{f}(\xi-\eta,s)\widehat{g}(\eta,s)\varphi_q\big(\Xi_{\iota_1\iota_2}(\xi-\eta,\eta)\big)\,d\eta,
\end{equation}
where $\mu=(wa,\iota_1)$, $\nu=(wa,\iota_2)$, and $\Phi_{(wa,+)\mu\nu}(\xi,\eta)=\Lambda_{wa}(\xi)-\Lambda_{\mu}(\xi-\eta)-\Lambda_{\nu}(\eta)$. As in Lemma \ref{BAX101}, we integrate by parts in time and recall that $|\dot{\Theta}_{wa}(\xi,s)|\lesssim 2^{-m+2\delta'm}$. For \eqref{baz2} it suffices to prove that, for any $s\in[2^{m-1},2^{m+1}]$,
\begin{equation}\label{baz20}
|\varphi_k(\xi)T_{\mu\nu;q}^{wa}[P_{k_1}V^{h_1,\iota_1},P_{k_2}V^{h_2,\iota_2}](\xi,s)|\lesssim \varep_1^22^{-0.002m}2^{k_2^-},
\end{equation}
\begin{equation}\label{baz22}
2^m|\varphi_k(\xi)T_{\mu\nu;q}^{wa}[P_{k_1}V^{h_1,\iota_1},\partial_s(P_{k_2}V^{h_2,\iota_2})](\xi,s)|\lesssim \varep_1^22^{-0.002m}2^{k_2^-},
\end{equation}
\begin{equation}\label{baz21}
2^m|\varphi_k(\xi)T_{\mu\nu;q}^{wa}[\partial_s(P_{k_1}V^{h_1,\iota_1}),P_{k_2}V^{h_2,\iota_2}](\xi,s)|\lesssim \varep_1^22^{-0.002m}2^{k_2^-}.
\end{equation}

{\bf{Step 1: proof of \eqref{baz20}.}} We decompose $P_{k_1}V^{h_1,\iota_1}=\sum_{j_1}V^{h_1,\iota_1}_{j_1,k_1}$ and $P_{k_2}V^{h_2,\iota_2}=\sum_{j_2}V^{h_2,\iota_2}_{j_2,k_2}$. The contribution of the pairs $(j_1,j_2)$ for which $\max(j_1,j_2)\leq 0.99m$ is negligible, due to Lemma \ref{tech5}, the assumptions $q\geq -0.002m-10$, $k_2\geq -p_0m-10$, and the bounds \eqref{baz4}, \eqref{baz9}. We estimate also
\begin{equation*}
\begin{split}
|\varphi_k(\xi)T_{\mu\nu;q}^{wa}[V^{h_1,\iota_1}_{j_1,k_1},&V^{h_2,\iota_2}_{j_2,k_2}](\xi,s)|\lesssim 2^{-2q-k_2}\|\widehat{V^{h_1,\iota_1}_{j_1,k_1}}(s)\|_{L^\infty}2^{3k_2}\|\widehat{V^{h_2,\iota_2}_{j_2,k_2}}(s)\|_{L^\infty}\\
&\lesssim\varep_1^22^{5\delta'm}2^{-4k_1^+}2^{2k_2-2q}(2^{-3k_2/2}2^{-j_2/2+\delta j_2})(2^{-3k_1/2}2^{-j_1/2+\delta j_1})\\
&\lesssim \varep_1^22^{0.01m}2^{-2q}2^{k_2^-/2}2^{-(j_1+j_2)/2+\delta(j_1+j_2)}
\end{split}
\end{equation*}
using \eqref{vcx1.15}. Since $2^{-2q}\lesssim 2^{0.01m}$ (due to \eqref{baz8}) and $2^{-k_2^-}\lesssim 2^{p_0m}$, this suffices to control the contribution of the pairs $(j_1,j_2)$ with $j_1+j_2\geq 0.99m$. The bounds \eqref{baz20} follow.

{\bf{Step 2: proof of \eqref{baz22}.}} Recall \eqref{loc2}. The contribution of the cubic terms $\mathcal{N}^{h,\geq 3}(s)$ can be estimated easily, using the bounds in the second line of \eqref{bax20}. The quadratic nonlinearities $\mathcal{N}^{h,2}$ contain two main types of terms: bilinear interactions of the metric components and bilinear interactions of the Klein-Gordon field (see \eqref{sac1.1}). The desired bounds follow from \eqref{baz31}--\eqref{baz32} in Lemma \ref{LemBaz2} below. 

{\bf{Step 3: proof of \eqref{baz21}.}} As before, the contribution of the cubic and higher order nonlinearities $\mathcal{N}^{h,\geq 3}(s)$ can be estimated using the bounds in the second line of \eqref{bax20}. The quadratic nonlinearities $\mathcal{N}^{h,2}$ can be estimated using the change of variables $\eta\to\xi-\eta$ and the bounds \eqref{baz33}--\eqref{baz34} in Lemma \ref{LemBaz2} below.
\end{proof} 

We estimate now the trilinear operators arising in the proof of the previous lemma:   

\begin{lemma}\label{LemBaz2} For $m_3,m_4\in\mathcal{M}_0$, $m\in\mathcal{M}$, (see \eqref{mults})  we define the trilinear operators
\begin{equation}\label{baz30}
\begin{split}
\mathcal{C}^{q,l}_{wa}[f,g,h](\xi):=&\int_{\mathbb{R}^3\times\mathbb{R}^3}\frac{\varphi_k(\xi)m(\xi-\eta,\eta)}{\Lambda_{wa}(\xi)-\Lambda_\mu(\xi-\eta)-\Lambda_\nu(\eta)}\varphi_q\big(\Xi_{\iota_1\iota_2}(\xi-\eta,\eta)\big)\varphi_{l}(\eta)\\
&\times m_3(\eta-\rho)m_4(\rho)\cdot\widehat{f}(\xi-\eta)\widehat{g}(\eta-\rho)\widehat{h}(\rho)\,d\eta d\rho,
\end{split}
\end{equation} 
where $q,l\in\mathbb{Z}$, and $\mu=(wa,\iota_1),\nu=(wa,\iota_2)$. Assume that $m\geq 100$, $s\in[2^{m-1},2^{m+1}]\cap [0,T]$, $k,k_1,k_2\in\mathbb{Z}$, $k\in[-\kappa m/4,\delta'm]$, $-p_0m-10\leq k_2\leq k_1\leq m/10$, and $q\geq -0.002 m+10$. Then 
\begin{equation}\label{baz31}
2^{-k_2^-}2^{|k_3-k_4|}\big|\mathcal{C}^{q,k_2}_{wa}[U^{h_1,\iota_1}_{j_1,k_1}(s),U^{h_3,\iota_3}_{j_3,k_3}(s),U^{h_4,\iota_4}_{j_4,k_4}(s)](\xi)\big|\lesssim \varep_1^32^{-1.002m},
\end{equation}
\begin{equation}\label{baz32}
2^{-k_2^-}\big|\mathcal{C}^{q,k_2}_{wa}[U^{h_1,\iota_1}_{j_1,k_1}(s),U^{\psi,\iota_3}_{j_3,k_3}(s),U^{\psi,\iota_4}_{j_4,k_4}(s)](\xi)\big|\lesssim \varep_1^32^{-1.002m},
\end{equation}
\begin{equation}\label{baz33}
2^{-k_2^-}2^{|k_3-k_4|}\big|\mathcal{C}^{q,k_1}_{wa}[U^{h_2,\iota_2}_{j_2,k_2}(s),U^{h_3,\iota_3}_{j_3,k_3}(s),U^{h_4,\iota_4}_{j_4,k_4}(s)](\xi)\big|\lesssim \varep_1^32^{-1.002m},
\end{equation}
and
\begin{equation}\label{baz34}
2^{-k_2^-}\big|\mathcal{C}^{q,k_1}_{wa}[U^{h_2,\iota_2}_{j_2,k_2}(s),U^{\psi,\iota_3}_{j_3,k_3}(s),U^{\psi,\iota_4}_{j_4,k_4}(s)](\xi)\big|\lesssim \varep_1^32^{-1.002m},
\end{equation}
for any $(k_1,j_1),(k_2,j_2),(k_3,j_3),(k_4,j_4)\in\mathcal{J}$ and $h_1,h_2,h_3,h_4\in\{h_{\al\be}\}$.
\end{lemma}

\begin{proof} Let $Y_1,Y_2,Y_3,Y_4$ denote the expressions in the left-hand sides of \eqref{baz31}, \eqref{baz32}, \eqref{baz33}, and \eqref{baz34} respectively. We may assume that $\xi=(\xi_1,0,0)$, $\xi_1>0$. We remark that the bounds \eqref{baz31} and \eqref{baz33} are different because $k_2$ can be very small, $k_2\geq -p_0m-10$, but $k_1$ cannot be so small, $k_1\geq -\kappa m/4-4$. The same remark applies to the bounds \eqref{baz32} and \eqref{baz34}.

{\bf{Step 1: proof of \eqref{baz31}.}} We may assume that $k,k_1\geq-\kappa m/4-10$. We estimate first, using \eqref{baz9}, \eqref{gene3x} (or \eqref{vcx1} if $\min(k_3,k_4)\leq -m+\delta m$), and \eqref{vcx1.15}
\begin{equation}\label{baz35}
\begin{split}
Y_1&\lesssim 2^{-k_2^-}2^{|k_3-k_4|}2^{-2q-k_2}\|\widehat{U^{h_1,\iota_1}_{j_1,k_1}}(s)\|_{L^\infty}\cdot \varep_1^22^{3k_2}2^{2\delta'm-m}2^{-k_2+\min\{k_3,k_4\}}2^{-10\max\{k_3^+,k_4^+\}}\\
&\lesssim \varep_1^32^{-0.995 m}2^{-j_1/2+\delta j_1}2^{\max\{k_3^-,k_4^-\}}2^{-8\max\{k_3^+,k_4^+\}}2^{-8k_1^+}.
\end{split}
\end{equation} 
This suffices to prove the desired bounds unless
\begin{equation}\label{baz35.3}
j_1\leq 0.02 m,\qquad k_1\in[-\kappa m/4-10,0.01m],\qquad\max\{k_3,k_4\}\in[-0.01m,0.01m].
\end{equation}
On the other hand, if the inequalities in \eqref{baz35.3} hold, then we analyze several subcases.

{\bf{Substep 1.1.}} Assume first that the inequalities in \eqref{baz35.3} hold and, in addition, 
\begin{equation}\label{baz36}
\min\{k_3,k_4\}\leq -0.03 m-30.
\end{equation}
By symmetry we may assume that $k_4\leq k_3$, therefore $|k_2-k_3|\leq 4$. We fix $\rho$ with $|\rho|\leq 2^{k_4+2}$ and estimate the $\eta$ integral by $C2^{-2q-k_2}2^{-0.99m}2^{-4k_1^+}2^{-4k_3^+}$ using Lemma \ref{gene1} with $2^l\approx 2^k$, $2^{l_1}\approx 2^{k_1}$, $2^{l_2}\approx 2^{k_3}$. Thus
\begin{equation*}
Y_1\lesssim 2^{-k_2^-}2^{k_3-k_4}2^{-2q-k_2}2^{-0.98m}2^{-4k_1^+}2^{-4k_3^+}2^{3k_4}
\end{equation*}
The conclusion follows in this case.

{\bf{Substep 1.2.}} Assume now that the inequalities in \eqref{baz35.3} hold and, in addition, 
\begin{equation}\label{baz40.2}
\min\{k_3,k_4\}\geq -0.03 m-30\qquad\text{ and }\qquad \max\{j_3,j_4\}\geq 0.97 m+k_2-100.
\end{equation}
We estimate first in the physical space, using \eqref{ener62}, \eqref{par6}, and \eqref{wws1},
\begin{equation}\label{baz40.11}
\begin{split}
Y_1&\lesssim 2^{-k_2^-}2^{|k_3-k_4|}2^{-2q-k_2}\|U^{h_1,\iota_1}_{j_1,k_1}(s)\|_{L^\infty}\|U^{h_3,\iota_3}_{j_3,k_3}(s)\|_{L^2}\|U^{h_4,\iota_4}_{j_4,k_4}(s)\|_{L^2}\\
&\lesssim \varep_1^32^{-0.95m}2^{-2k_2^-}2^{-(1-\delta)(j_3+j_4)}2^{-10k_1^+},
\end{split}
\end{equation}
which suffices if $k_2\geq -0.3m-30$.

On the other hand, if $k_2\leq -0.3m-30$ then we may assume that $|k_3-k_4|\leq 4$ and, by symmetry, $j_3\leq j_4$. Notice that the $\rho$ integral is bounded by $C2^{2\delta'm}2^{-j_3}2^{-k_3/2}2^{-j_4}2^{-k_4/2}2^{-4k_3^+}$, using just $L^2$ estimates on the two components. This would suffice if we had slightly stronger bounds on $j_3,j_4$, such as  $j_4\geq 1.03m+k_2-100$ or $j_3\geq 0.1m$. In the remaining case, when $\big|j_4-m-k_2+100\big|\leq 0.03m$ and $j_3\leq 0.1 m$, we need to gain by integration by parts in $\eta$. 

Recall that $\xi=(\xi_1,0,0)$ and insert cutoff functions of the form $\varphi_{\leq n}(\rho_2,\rho_3)$ and $\varphi_{>n}(\rho_2,\rho_3)$, where $n:=-0.15m+100$. More precisely, for $\ast\in\{\leq n,>n\}$ we define
\begin{equation*}
\begin{split}
&G_\ast(\xi):=2^{-k_2^-}\int_{\mathbb{R}^3\times\mathbb{R}^3}\frac{\varphi_k(\xi)m(\xi-\eta,\eta)\varphi_{k_2}(\eta)}{\Lambda_{wa}(\xi)-\Lambda_\mu(\xi-\eta)-\Lambda_\nu(\eta)}\varphi_q(\Xi_{\iota_1\iota_2}(\xi-\eta,\eta))\varphi_\ast(\rho_2,\rho_3)m_3(\eta-\rho)m_4(\rho)\\
&\times e^{-is[\Lambda_{wa,\iota_1}(\xi-\eta)+\Lambda_{wa,\iota_3}(\eta-\rho)+\Lambda_{wa,\iota_4}(\rho)]}\widehat{V^{h_1,\iota_1}_{j_1,k_1}}(\xi-\eta,s)\widehat{V^{h_3,\iota_3}_{j_3,k_3}}(\eta-\rho,s)\widehat{V^{h_4,\iota_4}_{j_4,k_4}}(\rho,s)\,d\eta d\rho.
\end{split}
\end{equation*} 
Notice that the $\eta$ derivative of the phase is bounded from below
\begin{equation*}
\begin{split}
\big|\nabla_{\eta}[&\Lambda_{wa,\iota_1}(\xi-\eta)+\Lambda_{wa,\iota_3}(\eta-\rho)]\big|\geq \frac{|(\rho_2,\rho_3)|}{|\eta-\rho|}-|\eta|\Big(\frac{1}{|\xi-\eta|}+\frac{1}{|\eta-\rho|}\Big)\gtrsim 2^{n-k_3}
\end{split}
\end{equation*} 
in the support of the integral defining $G_{> n}$. Thus $|G_{>n}(\xi)|\lesssim 2^{-2m}$, using integration by parts in $\eta$ with Lemma \ref{tech5} (recall that $j_1,j_3\leq 0.1m$). 

Finally we estimate $|G_{\leq n}(\xi)|$ as in the proof of Lemma \ref{gene1} (i), with $p=1/\delta$,
\begin{equation*}
\begin{split}
|G_{\leq n}(\xi)|&\lesssim 2^{-k_2}2^{-2q-k_2}2^{3k_2}\|\widehat{V^{h_1,\iota_1}_{j_1,k_1}}(s)\|_{L^\infty}\|\widehat{V^{h_3,\iota_3}_{j_3,k_3}}(s)\|_{L^\infty}\int_{\mathbb{R}^3}\varphi_{\leq n}(\rho_2,\rho_3)|\widehat{V^{h_4,\iota_4}_{j_4,k_4}}(\rho,s)|\,d\rho\\
&\lesssim \varep_1^22^{k_2}2^{0.04m}\|\widehat{V^{h_4,\iota_4}_{j_4,k_4}}(r\theta,s)\|_{L^2(r^2dr)L^p_\theta}\cdot 2^{(2n-2k_4)/p'}2^{3k_4/2}\\
&\lesssim\varep_1^32^{0.1m}2^{2n}2^{k_2}2^{-j_4}.
\end{split}
\end{equation*}
The desired conclusion follows since $2^{-j_4+k_2}\lesssim 2^{-0.97m}$, see \eqref{baz40.2}, and $2^{2n}\lesssim 2^{-0.3m}$.

{\bf{Substep 1.3.}} Assume now that the inequalities in \eqref{baz35.3} hold and, in addition, 
\begin{equation}\label{baz40}
\min\{k_3,k_4\}\geq -0.03m-30\qquad\text{ and }\qquad \max\{j_3,j_4\}\leq 0.97m+k_2-100.
\end{equation}
By symmetry, we may assume that $k_4\leq k_3$. We insert first cutoff functions of the form $\varphi_{q'}\big(\Xi_{\iota_3\iota_4}(\eta-\rho,\rho)\big)$ in the integral \eqref{baz30}. Using \eqref{baz3.1} (with $2^{l_1}\approx 2^{k_3}$, $2^{l_2}\approx 2^{k_4}$, $2^l\approx 2^{k_2}$, $b=q'$) the contribution is negligible if $q'\geq k_2-0.025m$. On the other hand, if 
\begin{equation}\label{baz40.12}
q'_0\leq q'\leq k_2-0.025m,\quad\text{ where }\quad q'_0:=-m/2+k_2/2+0.03m,
\end{equation}
then the contribution of the $\rho$ integral is bounded by $C\varep_1^22^{2\delta'm-1.02m}2^{k_3+k_4-k_2}2^{-4k_3^+-4k_4^+}$ (due to \eqref{gene3}) for any $\eta$. The desired bounds then follow, once we notice that the $\eta$ integral gains a factor of $2^{3k_2}$. 

To bound the contribution of $\varphi_{\leq q'_0}\big(\Xi_{\iota_3\iota_4}(\eta-\rho,\rho)\big)$, we further insert cutoff functions of the form $\varphi_{\leq n'}(\rho')$ and $\varphi_{>n'}(\rho')$, where $n':=-0.05m+k_4+100$ and $\rho'=(\rho_2,\rho_3)$. More precisely, as before, for $\ast\in\{\leq n',>n'\}$ we define
\begin{equation*}
\begin{split}
G'_\ast(\xi)&:=2^{-k_2^-}\int_{\mathbb{R}^3\times\mathbb{R}^3}\frac{\varphi_k(\xi)m(\xi-\eta,\eta)\varphi_{k_2}(\eta)\varphi_q(\Xi_{\iota_1\iota_2}(\xi-\eta,\eta))}{\Lambda_{wa}(\xi)-\Lambda_\mu(\xi-\eta)-\Lambda_\nu(\eta)}\\
&\times\varphi_{\leq q'_0}\big(\Xi_{\iota_3\iota_4}(\eta-\rho,\rho)\big)\varphi_\ast(\rho')m_3(\eta-\rho)m_4(\rho)\\
&\times e^{-is[\Lambda_{wa,\iota_1}(\xi-\eta)+\Lambda_{wa,\iota_3}(\eta-\rho)+\Lambda_{wa,\iota_4}(\rho)]}\widehat{V^{h_1,\iota_1}_{j_1,k_1}}(\xi-\eta,s)\widehat{V^{h_3,\iota_3}_{j_3,k_3}}(\eta-\rho,s)\widehat{V^{h_4,\iota_4}_{j_4,k_4}}(\rho,s)\,d\eta d\rho.
\end{split}
\end{equation*}
For \eqref{baz31} it remains to prove that
\begin{equation}\label{baz40.8}
2^{k_3-k_4}|G'_{\leq n'}(\xi)|+2^{k_3-k_4}|G'_{>n'}(\xi)|\lesssim \varep_1^32^{-1.01m}.
\end{equation}

Notice that the integral in the definition of $G'_{\leq n'}$ is supported in the set 
\begin{equation*}
\{(\eta,\rho):|\rho'|\lesssim 2^{n'},\,\widetilde{\Xi}(\eta,\rho)\lesssim 2^{k_3-k_2}2^{q'_0}\},
\end{equation*}
due to \eqref{par73} and the assumption $\Xi_{\iota_3\iota_4}(\eta-\rho,\rho)\lesssim 2^{q'_0}$. Therefore, using also Lemma \ref{box100.5},
\begin{equation}\label{baz40.3}
\begin{split}
|G'_{\leq n'}(\xi)|&\lesssim 2^{-k_2}2^{2q-k_2}\cdot 2^{2n'}2^{k_4}2^{k_2}(2^{k_3}2^{q'_0})^2\cdot \varep_1^32^{-k_1-k_3-k_4}2^{-8k_1^+-8k_3^+}2^{2\delta'm}\\
&\lesssim \varep_1^32^{-1.02m} 2^{2k_4}2^{-4k_3^+}.
\end{split}
\end{equation}

Finally we have to bound the functions $G'_{>n'}$. This can be done as in \eqref{baz40.3} if $j_3+j_4\geq 0.2 m$, since the gain of $2^{2n'}$ can be replaced by a gain of $2^{-j_3/2}2^{-j_4/2}$ coming from Lemma \ref{box100.5}. On the other hand, if $j_3+j_4\leq 0.2 m$ then we claim that 
\begin{equation}\label{baz40.4}
|G'_{>n'}(\xi)|\lesssim\varep_1^32^{-2m}.
\end{equation}
To see this we use integration by parts in $\eta$. We show that
\begin{equation}\label{baz40.5}
|\nabla_\eta\{\Lambda_{wa,\iota_1}(\xi-\eta)+\Lambda_{wa,\iota_3}(\eta-\rho)\}|\gtrsim 2^{-0.46m}2^{-k_2/2}
\end{equation}
in the support of the integral defining $G'_{>n'}(\xi)$. In view of Lemma \ref{tech5} (with $K\approx 2^{0.54m}2^{-k_2/2}$, $\epsilon=K^{-1}2^{\delta m}$), this would clearly suffice to prove \eqref{baz40.4}. 

To prove \eqref{baz40.5}, assume for contradiction that it fails, so $\Xi_{\iota_1\iota_3}(\xi-\eta,\eta-\rho)\leq 2^{-0.46m-k_2/2}$ for some $\eta,\rho$ in the support of the integral defining $G'_{>n'}(\xi)$. Since $\Xi_{\iota_3\iota_4}(\eta-\rho,\rho)\leq 2^{-0.47m+k_2/2}$, it follows from \eqref{par73} that $\widetilde{\Xi}(\eta,\rho)\leq 2^{-0.47m-k_2/2+k_3+4}$ and it follows from \eqref{par72}--\eqref{par72.4} that $\widetilde{\Xi}(\xi-\eta,\rho)\leq 2^{-0.46m-k_2/2+4}$. Therefore, using again \eqref{par72}--\eqref{par72.4} we have $\widetilde{\Xi}(\xi,\rho)\leq 2^{-0.46m-k_2/2+8}$, in contradiction with the assumption $|\rho'|\geq 2^{n'-4}\geq 2^{-0.05m+k_4+90}$ (recall that $-k_2/2\leq p_0m/2+10\leq 0.34m+10)$.  This completes the proof of \eqref{baz40.4}.

{\bf{Step 2: proof of \eqref{baz32}.}} An estimate similar to \eqref{baz35} still holds, using Lemma \ref{abc2.2} instead of \eqref{gene3x}. This proves the desired conclusion when $j_1\geq 0.1m$. On the other hand, if $j_1\leq 0.1 m$ then we notice that $|\nabla_\eta[s\Lambda_{wa,\iota_1}(\xi-\eta)+s\Lambda_{kg,\iota_3}(\eta-\rho)]|\gtrsim 2^m2^{-2k_3^+}$ in the support of the integral. Therefore, using integration by parts in $\eta$ with Lemma \ref{tech5}, $Y_2$ is negligible if $j_3\leq 0.95 m-3k_3^+$. Similarly, after making the change of variables $\rho\to\eta-\rho$, $Y_2$ is negligible if $j_4\leq 0.95 m-3k_4^+$. Finally, if
\begin{equation*}
j_1\leq 0.1 m,\qquad j_3\geq 0.95 m-3k_3^+,\qquad j_4\geq 0.95 m-3k_4^+
\end{equation*} 
then we estimate
\begin{equation*}
Y_2\lesssim 2^{-k_2^-}2^{-2q-k_2}\|\widehat{V^{h_1,\iota_1}_{j_1,k_1}}(s)\|_{L^\infty}2^{3k_2}\|V^{\psi,\iota_3}_{j_3,k_3}(s)\|_{L^2}\|V^{\psi,\iota_4}_{j_4,k_4}(s)\|_{L^2}\lesssim \varep_1^32^{0.1m}2^{-j_3}2^{-j_4}2^{-6k_3^+-6k_4^+}.
\end{equation*}
The bounds \eqref{baz32} follow.

{\bf{Step 3: proof of \eqref{baz33}.}} We may assume that $k_4\leq k_3$, thus $k,k_1,k_3\geq-\kappa m/4-10$ and $k\leq k_1+6$ and $k_1\leq k_3+6$. The main frequency parameters in the proof are $k_2$ and $k_4$. We estimate first, using \eqref{baz9} and \eqref{gene3x} (or \eqref{vcx1} if $k_4\leq -m+\delta m$)
\begin{equation}\label{baz51.6}
\begin{split}
Y_3&\lesssim 2^{-k_2^-}2^{k_3-k_4}2^{-2q-k_2}2^{3k_2/2}\|\widehat{U^{h_2,\iota_2}_{j_2,k_2}}(s)\|_{2}\cdot \varep_1^22^{2\delta'm-m}2^{k_4-k_1}2^{-10k_3^+}\\
&\lesssim \varep_1^32^{-0.995 m}2^{-j_2-k_2^-}2^{-8k_3^+}.
\end{split}
\end{equation} 
This suffices unless $j_2+k_2^-\leq 0.01m$. To deal with this case we analyze several subcases.

{\bf{Substep 3.1.}} Assume first that 
\begin{equation}\label{baz51}
j_2+k_2^-\leq 0.01m\qquad\text{ and }\qquad k_4\leq k_2-0.03 m-30.
\end{equation}
This is similar to the case analyzed in \eqref{baz36}. We fix $\rho$ with $|\rho|\leq 2^{k_4+2}$ and estimate the $\eta$ integral by $C2^{-2q-k_2}2^{-0.99m}2^{k_2-k}2^{-4k_3^+}$ using Lemma \ref{gene1} with $2^l\approx 2^k$, $2^{l_1}\approx 2^{k_2}$, $2^{l_2}\approx 2^{k_3}$. Thus
\begin{equation*}
Y_2\lesssim 2^{-k_2^-}2^{k_3-k_4}2^{-2q-k_2}2^{-0.98m}2^{k_2-k}2^{-4k_3^+}2^{2k_4}
\end{equation*}
The desired conclusion follows since $2^{k_4-k_2}\lesssim 2^{-0.03m}$, see \eqref{baz51}.

{\bf{Substep 3.2.}} Assume now that 
\begin{equation}\label{baz55}
j_2+k_2^-\leq 0.01m,\qquad k_4\geq k_2-0.03m-30,\qquad j_3\geq 0.97m-100.
\end{equation}
As in \eqref{baz40.11}, we estimate in the physical space, using \eqref{ener62}, \eqref{par6}, and \eqref{wws1},
\begin{equation}\label{baz55.1}
\begin{split}
Y_2&\lesssim 2^{-k_2^-}2^{k_3-k_4}2^{-2q-k_2}\|U^{h_2,\iota_2}_{j_2,k_2}(s)\|_{L^\infty}\|U^{h_3,\iota_3}_{j_3,k_3}(s)\|_{L^2}\|U^{h_4,\iota_4}_{j_4,k_4}(s)\|_{L^2}\\
&\lesssim \varep_1^32^{-0.99m}2^{-k_2}2^{-3k_4/2}2^{-j_3}2^{-j_4}2^{-4k_3^+}.
\end{split}
\end{equation}
Since $2^{-k_2}\lesssim 2^{p_0m}=2^{0.68m}$ and $2^{-k_4-j_4}\lesssim 1$, this suffices to prove \eqref{baz33} when $-k_4\leq 0.55 m$. On the other hand, if $k_4\leq-0.55 m$ then we can bound simply, using \eqref{vcx1.15} and \eqref{vcx1},
\begin{equation*}
\begin{split}
Y_2&\lesssim 2^{-k_2^-}2^{k_3-k_4}2^{-2q-k_2}2^{3k_2/2}\|\widehat{U^{h_2,\iota_2}_{j_2,k_2}}(s)\|_{L^2}\|\widehat{U^{h_3,\iota_3}_{j_3,k_3}}(s)\|_{L^\infty}2^{3k_4/2}\|\widehat{U^{h_4,\iota_4}_{j_4,k_4}}(s)\|_{L^2}\\
&\lesssim \varep_1^32^{0.01m}2^{k_4}2^{-j_3/2+\delta j_3}2^{-4k_3^+}.
\end{split}
\end{equation*}
Given \eqref{baz55}, this suffices to prove \eqref{baz33} if $k_4\leq-0.55 m$. 

{\bf{Substep 3.3.}} Assume now that
\begin{equation}\label{baz58}
j_2+k_2^-\leq 0.01m,\qquad k_4\geq k_2-0.03m-30,\qquad j_4\geq 0.97m-100.
\end{equation}
The bounds \eqref{baz55.1} still hold, but they only suffice to prove \eqref{baz33} when $2^{j_3+k_2+3k_4/2}\gtrsim 2^{-0.95m}$. This holds if $2^{j_3+k_2}\gtrsim 2^{0.12m}$ (because $2^{k_4}\gtrsim 2^{-0.71m}$) or when $2^{k_4}\gtrsim 2^{-0.18m}$ (because $2^{j_3+k_2}\gtrsim 2^{-0.68m}$. It remains to prove \eqref{baz33} in the case when
\begin{equation}\label{baz58.1}
j_3+k_2\leq 0.12m\qquad\text{ and }\qquad k_4\leq -0.18m-50.
\end{equation}

Assuming that both \eqref{baz58} and \eqref{baz58.1} hold, we consider the operator defined by $\eta$ integration first. More precisely, with $n:=-0.05m+50$ and $\ast\in\{\leq n,>n\}$ we define
\begin{equation}\label{baz59}
\begin{split}
&H_{\ast}(\xi,\rho):=2^{-k_2^-}\int_{\mathbb{R}^3}\frac{\varphi_k(\xi)m(\theta,\xi-\theta)\varphi_{k_1}(\xi-\theta)}{\Lambda_{wa}(\xi)-\Lambda_\mu(\theta)-\Lambda_\nu(\xi-\theta)}\varphi_q(\Xi_{\iota_1\iota_2}(\theta,\xi-\theta))\varphi_{\ast}(\Xi_{\iota_2\iota_3}(\theta,\xi-\rho-\theta))\\
&\times m_3(\xi-\rho-\theta)e^{-is[\Lambda_{wa,\iota_2}(\theta)+\Lambda_{wa,\iota_3}(\xi-\rho-\theta)]}\widehat{V^{h_2,\iota_2}_{j_2,k_2}}(\theta,s)\widehat{V^{h_3,\iota_3}_{j_3,k_3}}(\xi-\rho-\theta,s)\,d\theta.
\end{split}
\end{equation} 
This corresponds to the $\eta$ integral in the operator in \eqref{baz30}, after making the change of variables $\eta\to\xi-\theta$ and inserting angular cutoff functions of the form $\varphi_{\ast}(\Xi_{\iota_2\iota_3}(\theta,\xi-\rho-\theta))$. For \eqref{baz33} it suffices to prove that, for $\ast\in\{\leq n,>n\}$ and $\xi\in\mathbb{R}^3$,
\begin{equation}\label{baz59.1}
2^{k_3-k_4}\int_{\mathbb{R}^3}|H_{\ast}(\xi,\rho)|\big|\widehat{U^{h_4,\iota_4}_{j_4,k_4}}(\rho,s)\big|\,d\rho\lesssim \varep_1^32^{-1.01m}.
\end{equation} 

Using just \eqref{gene2} (with $\xi\to\xi-\rho$, $2^l\approx 2^{k}$, $2^{l_1}\approx 2^{k_2}$, $2^{l_2}\approx 2^{k_3}$), we bound
\begin{equation*}
|H_{\leq n}(\xi,\rho)|\lesssim 2^{-k_2^-}\cdot\varep_1^2 2^{-2q-k_2}2^{2n-2k}2^{2k_2^-}2^{-8k_3^+}\lesssim \varep_1^22^{-0.09m}2^{-4k_3^+}.
\end{equation*}
The bound \eqref{baz59.1} for $H_{\leq n}$ follows since $\big\|\widehat{U^{h_4,\iota_4}_{j_4,k_4}}(\rho,s)\big\|_{L^1_\rho}\lesssim \varep_12^{-j_4}2^{k_4}2^{\delta'm}$ and $2^{-j_4}\lesssim 2^{-0.97m}$.

On the other hand, we claim that $H_{>n}$ is negligible, i. e. $|H_{>n}(\xi,\rho)|\lesssim \varep_1^22^{-4m}2^{-4k_3^+}$. This follows by integration by parts in $\theta$ using Lemma \ref{tech5} (as in the proof of \eqref{baz3.1}), once we notice that the $\theta$ gradient of the phase is bounded from below by $c2^{n}2^m\gtrsim 2^{0.9m}$ in the support of the integral, and recall that $2^{\max(j_2,j_3)}\lesssim 2^{0.8m}$ (due to \eqref{baz58}--\eqref{baz58.1}). This completes the proof of \eqref{baz33} when \eqref{baz58} holds.

{\bf{Substep 3.4.}} Finally, assume that 
\begin{equation}\label{baz52}
j_2+k_2^-\leq 0.01m,\qquad k_4\geq k_2-0.03m-30,\qquad \max(j_3,j_4)\leq 0.97m-100.
\end{equation}
This is similar to the case analyzed in \eqref{baz40}. We insert the cutoff functions $\varphi_{q'}\big(\Xi_{\iota_3\iota_4}(\eta-\rho,\rho)\big)$ in the integral \eqref{baz30}. Using \eqref{baz3.1} (with $l_1=k_3$, $l_2=k_4$, $l=k_1$, $b=q'$) the contribution is negligible if $q'\geq -0.025m$. On the other hand, if
\begin{equation}\label{baz52.2}
q'_1\leq q'\leq -0.025m,\quad\text{ where }\quad q'_1:=-m/2-k_4/2+0.01 m
\end{equation}
then the contribution of the $\rho$ integral is bounded by $C\varep_1^22^{-0.02m-m}2^{-2k_1}2^{k_3+k_4}2^{-4k_3^+}$ (due to \eqref{gene3}) for any $\eta$. The desired conclusion then follows, once we notice that the $\eta$ integral gains a factor of $2^{2k_2}2^{\delta'm}$.

To bound the contribution of $\varphi_{\leq q'_1}\big(\Xi_{\iota_3\iota_4}(\eta-\rho,\rho)\big)$, we make the change of variables $\eta\to\xi-\theta$ and insert angular cutoff functions of the form $\varphi_{\ast}(\Xi_{\iota_2\iota_3}(\theta,\xi-\rho-\theta))$. More precisely we define
\begin{equation*}
\begin{split}
H'_{\ast}(\xi,\rho):=&2^{-k_2^-}\int_{\mathbb{R}^3}\frac{\varphi_k(\xi)m(\theta,\xi-\theta)\varphi_{k_1}(\xi-\theta)}{\Lambda_{wa}(\xi)-\Lambda_\mu(\theta)-\Lambda_\nu(\xi-\theta)}\varphi_q(\Xi_{\iota_1\iota_2}(\theta,\xi-\theta))\\
&\times \varphi_{\leq q'_1}\big(\Xi_{\iota_3\iota_4}(\xi-\theta-\rho,\rho)\big)\varphi_{\ast}(\Xi_{\iota_2\iota_3}(\theta,\xi-\rho-\theta))\\
&\times m_3(\xi-\rho-\theta)e^{-is[\Lambda_{wa,\iota_2}(\theta)+\Lambda_{wa,\iota_3}(\xi-\rho-\theta)]}\widehat{V^{h_2,\iota_2}_{j_2,k_2}}(\theta,s)\widehat{V^{h_3,\iota_3}_{j_3,k_3}}(\xi-\rho-\theta,s)\,d\theta.
\end{split}
\end{equation*}
where $n:=-0.05m+50$ as in \eqref{baz59} and $\ast\in\{\leq n,>n\}$.
For \eqref{baz31} it remains to prove that
\begin{equation}\label{baz52.8}
2^{k_3-k_4}\int_{\mathbb{R}^3}|H'_{\ast}(\xi,\rho)|\big|\widehat{U^{h_4,\iota_4}_{j_4,k_4}}(\rho,s)\big|\,d\rho\lesssim \varep_1^32^{-1.01m}.
\end{equation}

We consider first the contribution of $|H'_{\leq n}(\xi,\rho)|$. We use the restrictions 
\begin{equation}\label{baz52.9}
\widetilde{\Xi}(\xi-\theta-\rho,\rho)\lesssim 2^{q'_1}\qquad\text{ and }\qquad\widetilde{\Xi}(\xi-\theta-\rho,\theta)\lesssim 2^{n},
\end{equation}
which hold in the support of the defining integrals. Therefore $\widetilde{\Xi}(\xi-\theta,\rho)\lesssim 2^{q'_1}2^{k_3-k_1}$, using \eqref{par73}. Moreover, since $2^{q'_1}\lesssim 2^{n}$, we have $\widetilde{\Xi}(\rho,\theta)\lesssim 2^n$ (using \eqref{par72.4}). In addition $\widetilde{\Xi}(\xi-\theta,\theta)\lesssim 2^n2^{k_3-k_1}$ and then $\widetilde{\Xi}(\xi,\theta)\lesssim 2^n2^{k_3-k}$ (using \eqref{par73}). Therefore the support of the $(\theta,\rho)$ integral is included in the set $\{(\theta,\rho):\widetilde{\Xi}(\xi,\theta)\lesssim 2^n2^{k_3-k},\,\widetilde{\Xi}(\xi-\theta,\rho)\lesssim 2^{q'_1}2^{k_3-k_1}\}$. Thus
\begin{equation*}
\begin{split}
2^{k_3-k_4}\int_{\mathbb{R}^3}|H'_{\leq n}(\xi,\rho)|\big|\widehat{U^{h_4,\iota_4}_{j_4,k_4}}(\rho,s)\big|\,d\rho&\lesssim 2^{k_3-k_4}2^{-2q-2k_2}2^{2k_3^+}\|\widehat{V^{h_2,\iota_2}_{j_2,k_2}}(s)\|_{L^\infty}\|\widehat{V^{h_3,\iota_3}_{j_3,k_3}}(s)\|_{L^\infty}\\
&\times\|\widehat{U^{h_4,\iota_4}_{j_4,k_4}}(s)\|_{L^\infty}\cdot 2^{3k_4}(2^{q'_1}2^{k_3-k_1})^22^{3k_2}(2^{n}2^{k_3-k})^2\\
&\lesssim \varep_1^32^{2n}2^{-0.95m}2^{-8k_3^+},
\end{split}
\end{equation*}
using Lemma \ref{box100.5} and \eqref{baz52.2} in the last inequality. This gives \eqref{baz52.8} when $\ast=\leq n$.

The same argument also gives the desired bounds \eqref{baz52.8} when $\ast=>n$ if $j_3\geq 0.2 n$ or if $k_3\geq 0.01m$. In the remaining case, when $j_3\leq 0.2m$ and $k_3\in[-\kappa m/4-10,0.01m]$ we can integrate by parts in $\theta$, using Lemma \ref{tech5}, to see that the contribution is negligible, $|H'_{\ast}(\xi,\rho)|\lesssim \varep_1^22^{-4m}$. This completes the proof of \eqref{baz52.8}.

{\bf{Step 4: proof of \eqref{baz34}.}} This is similar to the proof of \eqref{baz32}. Using Lemma \ref{abc2.2} instead of \eqref{gene3x} and estimating as in \eqref{baz51.6}, the bounds \eqref{baz34} follow when $j_2+k_2^-\geq 0.01m$ or when $\max\{k_3,k_4\}\geq 0.01m$. On the other hand, if $k_3\leq 0.01m$ then we notice that $|\nabla_\eta[s\Lambda_{wa,\iota_2}(\xi-\eta)+s\Lambda_{kg,\iota_3}(\eta-\rho)]|\gtrsim 2^{0.95m}$ in the support of the integral. Therefore, using Lemma \ref{tech5}, $Y_4$ is negligible if $j_3\leq 0.9m$. Similarly, $Y_4$ is negligible if $j_4\leq 0.9m$. Finally, if $j_2+k_2^-\leq 0.01 m$ and $\min(j_3,j_4)\geq 0.9m$ then we estimate
\begin{equation*}
Y_4\lesssim 2^{-k_2^-}2^{-2q-k_2}\|\widehat{U^{h_2,\iota_2}_{j_2,k_2}}(s)\|_{L^\infty}2^{3k_2}\|U^{\psi,\iota_3}_{j_3,k_3}(s)\|_{L^2}\|U^{\psi,\iota_4}_{j_4,k_4}(s)\|_{L^2}\lesssim \varep_1^32^{0.1m}2^{-j_3}2^{-j_4}.
\end{equation*}
The bounds \eqref{baz34} follow in this last case as well.
\end{proof}

We conclude this subsection with an estimate on certain cubic expressions.

\begin{lemma}\label{SAR100}
Assume that $m\geq 100$, $t\in[2^{m-1},2^{m+1}]\cap [0,T]$, $k,k_1,k_2\in\mathbb{Z}$, $k\in[-\kappa m,\delta'm]$, $k_1,k_2\in[-p_0m-10,m/10+10]$, $h_1\in\{h_{\al\be}\}$, $\iota_1,\iota_2\in\{+,-\}$, and $q\in[-m,2]$. Then
\begin{equation}\label{sar1}
\begin{split}
\varphi_k(\xi)\int_{\mathbb{R}^3}|\widehat{P_{k_1}U^{h_1,\iota_1}}(\xi-\eta,t)|&|\widehat{P_{k_2}H}(\eta,t)|\varphi_{\leq q}\big(\Xi_{\iota_1\iota_2}(\xi-\eta,\eta)\big)\,d\eta\\
&\lesssim \varep_1^32^{2q}2^{-m+2\kappa m}2^{\min(k_1^-,k_2^-)}2^{-10(k_1^++k_2^+)}
\end{split}
\end{equation}
for any $\xi\in\mathbb{R}^3$, and $H=|\nabla|^{-1}\mathcal{N}$, $\mathcal{N}\in\mathcal{QU}$ (see \eqref{sho3.9}).
\end{lemma}

\begin{proof} It follows from \eqref{sho4} and \eqref{Linfty3.34} that 
\begin{equation}\label{sar2}
\|\widehat{P_{k_2}H}(r\theta,t)\|_{L^2(r^2dr)L^p_\theta}\lesssim\varep_1^22^{2\delta'm-m}2^{-k_2/2}2^{-15k_2^+},
\end{equation}
for $p=1/\delta$. The proof of \eqref{sar1} is similar to the proof of Lemma \ref{gene1} (i). We may assume that $\xi=(\xi_1,0,0)$, $\xi_1\in[2^{k-1},2^{k+1}]$.  The $\eta$ integral in the left-hand side is supported in the set $\big\{|\eta|\approx 2^{k_2},\,|\xi-\eta|\approx 2^{k_1},\,\sqrt{\eta_2^2+\eta_3^2}\lesssim 2^{q+k_1+k_2-k}\big\}$. With $\mathbf{e}_1:=(1,0,0)$, we estimate the left-hand side of \eqref{sar1} by
\begin{equation*}
\begin{split}
C\|\widehat{P_{k_1}U^{h_1,\iota_1}}(t)\|_{L^\infty}&\int_{[0,\infty)\times{\mathbb{S}}^2,\,|\theta-\mathbf{e}_1|\lesssim 2^{q+k_1-k}}|\widehat{P_{k_2}H}(r\theta,t)|\,r^2 dr d\theta\\
&\lesssim \|\widehat{P_{k_1}U^{h_1,\iota_1}}(t)\|_{L^\infty}\|\widehat{P_{k_2}H}(r\theta,t)\|_{L^2(r^2dr)L^p_\theta}\cdot 2^{2(q+k_1-k)/p'}2^{3k_2/2}\\
&\lesssim \varep_12^{-k_1}2^{\delta'm}2^{-20k_1^+}\cdot \varep_1^22^{3\delta'm-m}2^{-k_2/2}2^{-15k_2^+}2^{2q+2k_1-2k}2^{3k_2/2}.
\end{split}
\end{equation*}
The desired conclusion follows once we recall that $k\geq-\kappa m$. 
\end{proof}

\subsection{The nonlinear terms $\mathcal{R}_3^{h_{\al\be}}$} We are now ready to bound the quasilinear contributions defined in \eqref{bax9}, for which we can still prove strong bounds.

\begin{lemma}\label{TAR100} 
The bounds \eqref{bax13} hold for $m\geq 100$, $k\in[-\kappa m,\delta'm]$, and $a=3$.
\end{lemma}

\begin{proof} The formulas \eqref{bax5} and \eqref{zaq5} show that
\begin{equation*}
\begin{split}
&\mathfrak{q}_{wa,\iota_1}(\theta,\eta)=\iota_1\widehat{h_{00}}(\eta)\frac{\Lambda_{wa}(\theta)}{2}+\widehat{h_{0j}}(\eta)\theta_j+\iota_1\widehat{h_{jk}}(\eta)\frac{\theta_j\theta_k}{2\Lambda_{wa}(\theta)}\\
&=\widehat{F}(\eta)m_{F,\iota_1}(\theta,\eta)+\widehat{\underline{F}}(\eta)m_{\underline{F},\iota_1}(\theta,\eta)+\widehat{\rho}(\eta)m_{\rho,\iota_1}(\theta,\eta)+\widehat{\omega_l}(\eta)m_{\omega_l,\iota_1}(\theta,\eta)\\
&+\widehat{\Omega_a}(\eta)m_{\Omega_a,\iota_1}(\theta,\eta)+\widehat{\vartheta_{ab}}(\eta)m_{\vartheta_{ab},\iota_1}(\theta,\eta),
\end{split}
\end{equation*}
where
\begin{equation}\label{tar1}
\begin{split}
&m_{F,\iota_1}(\theta,\eta):=\iota_1\frac{|\eta|^2|\theta|^2-\eta_j\eta_k\theta_j\theta_k}{2|\eta|^2|\theta|},\qquad m_{\underline{F},\iota_1}(\theta,\eta):=\iota_1\frac{|\eta|^2|\theta|^2+\eta_j\eta_k\theta_j\theta_k}{2|\eta|^2|\theta|},\\
&m_{\rho,\iota_1}(\theta,\eta):=-i\frac{\eta_j\theta_j}{|\eta|},\qquad m_{\omega_l,\iota_1}(\theta,\eta):=i\frac{\in_{jkl}\eta_k\theta_j}{|\eta|},\\
&m_{\Omega_a,\iota_1}(\theta,\eta):=\iota_1\frac{(\in_{kla}\eta_j\eta_l+\in_{jla}\eta_k\eta_l)\theta_j\theta_k}{2|\eta|^2|\theta|},\qquad m_{\vartheta_{ab},\iota_1}(\theta,\eta):=-\iota_1\frac{\in_{jpa}\in_{kqb}\eta_p\eta_q\theta_j\theta_k}{2|\eta|^2|\theta|}.
\end{split}
\end{equation}
We rewrite $\mathfrak{q}_{wa,\iota_1}$ in terms of the normalized solutions $U^{G,\iota_2}$, using \eqref{on5},
\begin{equation}\label{tar2}
\mathfrak{q}_{wa,\iota_1}(\theta,\eta)=\sum_{\iota_2\in\{+,-\}}\sum_{G\in\{F,\underline{F},\rho,\omega_a,\Omega_a,\vartheta_{ab}\}}\iota_2\frac{i\widehat{U^{G,\iota_2}}(\eta)}{2|\eta|}m_{G,\iota_1}(\theta,\eta).
\end{equation}

We substitute this formula into the definition \eqref{bax9} and decompose dyadically. Let
\begin{equation}\label{tar3}
\begin{split}
\mathcal{R}^{h_{\al\be};\iota_1,\iota_2}_{3;k_1,k_2,q}(\xi,t)&:=\frac{-e^{-i\Theta_{wa}(\xi,t)}}{(2\pi)^3}\int_{\mathbb{R}^3}e^{it\Lambda_{wa}(\xi)}\varphi_q(\Xi_{\iota_1\iota_2}(\xi-\eta,\eta))\widehat{P_{k_1}U^{h_{\al\be},\iota_1}}(\xi-\eta,t)\\
&\times\sum_{G\in\{F,\underline{F},\rho,\omega_a,\Omega_a,\vartheta_{ab}\}}\iota_2\widehat{P_{k_2}U^{G,\iota_2}}(\eta,t)\frac{m_{G,\iota_1}(\xi-\eta,\eta)}{2|\eta|}\varphi_{\geq 1}(\langle t\rangle^{p_0}\eta)\,d\eta,
\end{split}
\end{equation}
and notice that
\begin{equation*}
\varphi_k(\xi)\mathcal{R}^{h_{\al\be}}_{3}(\xi,t)=\sum_{\iota_1,\iota_2\in\{+,-\}}\sum_{(k_1,k_2)\in\mathcal{X}_k}\sum_{q\leq 4}\varphi_k(\xi)\mathcal{R}^{h_{\al\be};\iota_1,\iota_2}_{3;k_1,k_2,q}(\xi,t).
\end{equation*}
For \eqref{bax13} it suffices to prove that if $m\geq 100$, $k\in[-\kappa m,\delta'm]$, and $\xi\in\mathbb{R}^3$, then
\begin{equation}\label{tar4}
\Big|\sum_{\iota_1,\iota_2\in\{+,-\}}\sum_{(k_1,k_2)\in\mathcal{X}_k}\sum_{q\leq 4}\varphi_k(\xi)\int_{t_1}^{t_2}\mathcal{R}^{h_{\al\be};\iota_1,\iota_2}_{3;k_1,k_2,q}(\xi,s)\,ds\Big|\lesssim\varep_1^22^{-\delta m/2}2^{-N_0k^+}.
\end{equation}

Notice that all the multipliers $m_{G,\iota_1}$ are sums of symbols of the form $|\theta|m_1(\theta)m_2(\eta)$, with $m_1,m_2\in\mathcal{M}_0$ (see definition \eqref{mults0}). We may assume that $k_2\geq -p_0m-10$ in the sum in \eqref{tar4}. Using the bounds $\|\widehat{P_{k_1}U^{h_{\al\be},\iota_1}}\|_{L^1}\lesssim \varep_12^{2k_1-\delta k_1}2^{\delta'm}$, it is easy to see that the contribution in the sum corresponding to the pairs $(k_1,k_2)$ for which $k_1\leq -p_0m-10$ (thus $|k_2-k|\leq 4$) is bounded by $C2^{-0.1m}$. Also, the contribution of the pairs $(k_1,k_2)$ for which $\max(k_1,k_2)\geq m/10$ is bounded by $C2^{-0.1m}$, using just $L^2$ estimates on the profiles. Finally, the contribution of the indices $q$ with $q\geq q_0:=(-m+k-k_1-k_2)/2+\delta m/8$ is bounded as claimed, due to Lemma \ref{LemBaz1}. After these reductions, it suffices to prove that 
\begin{equation}\label{tar8}
\begin{split}
\Big|\varphi_k(\xi)\sum_{G\in\{F,\underline{F},\rho,\omega_a,\Omega_a,\vartheta_{ab}\}}&\int_{\mathbb{R}^3}\varphi_{\leq q_0}(\Xi_{\iota_1\iota_2}(\xi-\eta,\eta))\widehat{P_{k_1}U^{h_{\al\be},\iota_1}}(\xi-\eta,t)\widehat{P_{k_2}U^{G,\iota_2}}(\eta,t)\\
&\times\frac{m_{G,\iota_1}(\xi-\eta,\eta)}{2|\eta|}\varphi_{\geq 1}(\langle t\rangle^{p_0}\eta)\,d\eta\Big|\lesssim\varep_1^22^{-\delta m-m}2^{-N_0k^+},
\end{split}
\end{equation}
for any $t\in [2^{m-1},2^{m+1}]$, $\iota_1,\iota_2\in\{+,-\}$, and $k_1,k_2\in[-p_0m-10,m/10]$.

To prove \eqref{tar8} we use the fact that the multipliers $m_{G,\iota_1}$ satisfy certain null structure conditions. Indeed, notice that
\begin{equation}\label{tar9}
|v_aw_b-v_bw_a|^2+\big(|v|^2|w|^2-|(v\cdot w)|^2\big)\lesssim |v||w|\big(|v||w|-|(v\cdot w)|\big)\lesssim |v|^2|w|^2\Xi_{\iota}(v,w)^2
\end{equation} 
for any $\iota\in\{+,-\}$, $a,b\in\{1,2,3\}$, and $v,w\in\mathbb{R}^3$. Therefore, using the definitions \eqref{tar1},
\begin{equation}\label{tar10}
|m_{G,\iota_1}(\xi-\eta,\eta)|\lesssim 2^{k_1}2^{q_0}
\end{equation}
in the support of the integral in \eqref{tar8}, for $G\in\{F,\omega_a,\Omega_a,\vartheta_{ab}\}$. Using \eqref{gene2}, it follows that the integrals in \eqref{tar8} corresponding to the functions $F,\omega_a,\Omega_a,\vartheta_{ab}$ are bounded by
\begin{equation*}
C\varep_1^22^{\delta'm}2^{k_1+k_2}\cdot 2^{\min(k_1,k_2)}2^{2q_0-2k}\cdot 2^{k_1-k_2}2^{q_0}2^{-4(k_1^++k_2^+)}\lesssim \varep_1^22^{2\kappa m}2^{-3m/2}2^{-\min(k_1,k_2)/2},
\end{equation*}
which suffices.

Clearly, the symbols $m_{\rho,\pm}$ and $m_{\underline{F},\pm}$ in \eqref{tar1} do not satisfy favorable null structure bounds like \eqref{tar10}. However, we can use the identities \eqref{zaq11.2}--\eqref{zaq11} to combine the $\rho$ and $\underline{F}$ terms and extract a cancellation. Indeed, notice that
\begin{equation}\label{tar12}
\begin{split}
U^{\rho,+}&=\partial_t\rho-i|\nabla|\rho\\
&=\partial_t(R_0\underline{F}+R_0\tau+|\nabla|^{-1}E_0^{\geq 2})-i|\nabla|(R_0\underline{F}+R_0\tau+|\nabla|^{-1}E_0^{\geq 2})\\
&=-|\nabla|\underline{F}+|\nabla|^{-1}\mathcal{N}^{\underline{F}}-|\nabla|\tau+|\nabla|^{-1}\mathcal{N}^{\tau}+|\nabla|^{-1}\partial_0E_0^{\geq 2}-i\partial_t\underline{F}-i\partial_t\tau-iE_0^{\geq 2}\\
&=-iU^{\underline{F},+}+\big\{-iU^{\tau,+}+|\nabla|^{-1}\mathcal{N}^{\underline{F}}+|\nabla|^{-1}\mathcal{N}^{\tau}+|\nabla|^{-1}\partial_0E_0^{\geq 2}-iE_0^{\geq 2}\big\}.
\end{split}
\end{equation}
Therefore $U^{\rho,\pm}=\mp iU^{\underline{F},\pm}+H^{\pm}$, where the functions $H^+$ and $H^-=\overline{H^+}$ satisfy the bounds \eqref{sar1} of Lemma \ref{SAR100}. 

We combine the contributions of $U^{\rho,\pm}$ and $U^{\uF,\pm}$ in the left-hand side of \eqref{tar8}. Notice that
\begin{equation*}
\begin{split}
\widehat{U^{\rho,\iota_2}}(\eta)&m_{\rho,\iota_1}(\xi-\eta,\eta)+\widehat{U^{\uF,\iota_2}}(\eta)m_{\uF,\iota_1}(\xi-\eta,\eta)\\
&=\widehat{U^{\uF,\iota_2}}(\eta)\widetilde{m}(\xi-\eta,\eta)+\widehat{H^{\iota_2}}(\eta)m_{\rho,\iota_1}(\xi-\eta,\eta),
\end{split}
\end{equation*}
as a consequence of \eqref{tar1} and \eqref{tar12}, where 
\begin{equation*}
\widetilde{m}(\theta,\eta):=-\iota_2\frac{\eta_j\theta_j}{|\eta|}+\iota_1\frac{|\eta|^2|\theta|^2+(\eta_j\theta_j)^2}{2|\eta|^2|\theta|}=\frac{\iota_1|\theta|}{8}|\Xi_{\iota_1\iota_2}(\theta,\eta)|^4.
\end{equation*}
The main point is that the combined symbols $\widetilde{m}$ satisfy favorable null-structure bounds of the form $|\widetilde{m}(\xi-\eta,\eta)|\lesssim 2^{k_1}2^{q_0}$, similar to \eqref{tar10} (see also \eqref{par2.1}), in the support of the integral in \eqref{tar8}. As before, this suffices to bound the corresponding contributions in \eqref{tar8}. Finally, the contributions of the functions $H^{\pm}$ to the left-hand side of \eqref{tar8} are bounded as claimed, as a consequence of \eqref{sar1}. This completes the proof of the lemma.
\end{proof}

\subsection{The nonlinear terms $\mathcal{R}_5^{h_{\al\be}}$} Finally, for the contribution of the semilinear quadratic terms we prove weaker bounds, but still sufficient to conclude the proof of \eqref{bax12}. 

\begin{lemma}\label{TAR101} 
Assume that $m\geq 100$, $t_1,t_2\in[2^{m-1},2^{m+1}]\cap [0,T]$, and $k\in[-\kappa m,\delta'm]$. Then, for any $\xi\in\mathbb{R}^3$,
\begin{equation}\label{tar14}
\Big|\varphi_k(\xi)\int_{t_1}^{t_2}\mathcal{R}_5^{h_{\al\be}}(\xi,s)\,ds\Big|\lesssim\varep_1^22^{\delta m/2}2^{-k^-}2^{-N_0k^+}.
\end{equation}
Moreover, for $G\in\{F,\Omega_a,\vartheta_{ab}\}$ we have
\begin{equation}\label{tar15}
\Big|\varphi_k(\xi)\int_{t_1}^{t_2}\mathcal{R}_5^{G}(\xi,s)\,ds\Big|\lesssim\varep_1^22^{-\delta m/2}2^{-k^-}2^{-N_0k^+},
\end{equation}
where (compare with the definitions \eqref{zaq2})
\begin{equation}\label{tar16}
\mathcal{R}^F_5:=\frac{1}{2}[\mathcal{R}_5^{h_{00}}+R_jR_k\mathcal{R}_5^{h_{jk}}],
\quad\mathcal{R}_5^{\Omega_a}:=\in_{akl}R_kR_m\mathcal{R}_5^{h_{lm}},
\quad\mathcal{R}_5^{\vartheta_{ab}}:=\in_{amp}\in_{bnq}R_mR_n\mathcal{R}_5^{h_{pq}}.
\end{equation}
\end{lemma}

\begin{proof} Recall that $\mathcal{R}_5^{h_{\al\be}}(\xi,t):=e^{-i\Theta_{wa}(\xi,t)}e^{it\Lambda_{wa}(\xi)}\widehat{\mathcal{S}^2_{\al\be}}(\xi,t)$ and $\mathcal{S}^2_{\al\be}=-(Q^2_{\al\be}+P^2_{\al\be})$, where the semilinear nonlinearities $Q^2_{\al\be}$ and $P^2_{\al\be}$ are explicitly given in \eqref{zaq20}--\eqref{zaq21}. We decompose $Q_{\al\be}^2+P_{\al\be}^2=\Pi^I_{\al\be}+\mathcal{C}^{I,\,\geq 3}_{\al\be}-\mathcal{S}^{I,2}_{\al\be}$ according to Lemma \ref{Pstructure} (there are no commutator remainders when $\LL=I$). The desired conclusions follow from Lemmas \ref{TAR102}--\ref{TAR104} below.
\end{proof}

We prove first strong bounds for the cubic terms and the semilinear null forms.

\begin{lemma}\label{TAR102} 
Assume that $m,k,t_1,t_2,\xi$ are as in Lemma \ref{TAR101} and $\mathcal{C}^{I,\,\geq 3}$ is a semilinear cubic remainder (see Definition \ref{nake4}). Then
\begin{equation}\label{tar21}
\Big|\varphi_k(\xi)\int_{t_1}^{t_2}e^{-i\Theta_{wa}(\xi,s)}e^{is\Lambda_{wa}(\xi)}\widehat{\mathcal{C}^{I,\,\geq 3}}(\xi,s)\,ds\Big|\lesssim\varep_1^22^{-\delta m/2}2^{-N_0k^+}.
\end{equation}
\end{lemma}

\begin{proof} We decompose dyadically in frequency and use first Lemma \ref{SAR100}. It remains to bound the contribution of the middle frequencies. Bilinear interactions of the quadratic and higher order expression are easy to bound, using just \eqref{sho4}. It remains to prove that
\begin{equation}\label{tar22}
\begin{split}
&\big|\varphi_k(\xi)\widehat{\mathcal{C}_{k_1,k_2}}(\xi,s)\big|\lesssim\varep_1^22^{-1.01m}\\
\text{ where }\quad&\widehat{\mathcal{C}_{k_1,k_2}}(\xi,s):=\int_{\mathbb{R}^3}m_1(\xi-\eta)\widehat{P_{k_1}U^{h_1,\iota_1}}(\xi-\eta,s)m_2(\eta)\widehat{P_{k_2}\mathcal{N}}(\eta,s)\,d\eta,
\end{split}
\end{equation}
for any $s\in[2^{m-1},2^{m+1}]$, $k_1,k_2\in[-4\kappa m,10\kappa m]$, $m_1,m_2\in\mathcal{M}_0$, $h_1\in\{h_{\al\be}\}$, $\iota_1\in\{+,-\}$, and $\mathcal{N}\in\mathcal{QU}$. Using \eqref{wws1} and Corollary \ref{SHO4} we estimate
\begin{equation*}
\begin{split}
&\|P_k\mathcal{L}\mathcal{C}_{k_1,k_2}(t)\|_{L^2}\lesssim 2^{-2m+0.01m},\qquad \sum_{l\in\{1,2,3\}}\|P_k(x_l\mathcal{L}\mathcal{C}_{k_1,k_2})(t)\|_{L^2}\lesssim 2^{-m+0.01m},
\end{split}
\end{equation*}
for any $\mathcal{L}=\Omega^\gamma$, $|\gamma|\leq 1$. The desired bounds \eqref{tar22} follow by interpolation, see \eqref{consu2}.
\end{proof}

\begin{lemma}\label{TAR103} 
Assume that $m,k,t_1,t_2,\xi$ are as in Lemma \ref{TAR101} and $\Pi$ is a semilinear null form (see Definition \ref{nullStr1}). Then
\begin{equation}\label{tar25}
\Big|\varphi_k(\xi)\int_{t_1}^{t_2}e^{-i\Theta_{wa}(\xi,s)}e^{is\Lambda_{wa}(\xi)}\widehat{\Pi}(\xi,s)\,ds\Big|\lesssim\varep_1^22^{-\delta m/2}2^{-N_0k^+}.
\end{equation}
\end{lemma}

\begin{proof} After decomposing dyadically in frequency, it remains to prove that
\begin{equation}\label{tar27}
\sum_{(k_1,k_2)\in\mathcal{X}_k}Y_{k_1,k_2}\lesssim \varep_1^22^{-\delta m/2}2^{-N_0k^+}
\end{equation}  
where, with $n\in\mathcal{M}^0_{\iota_1\iota_2}$ being a null multiplier as in Definition \ref{nullStr1}.
\begin{equation}\label{tar28}
Y_{k_1,k_2}:=\Big|\varphi_k(\xi)\int_{t_1}^{t_2}e^{is\Lambda_{wa}(\xi)-i\Theta_{wa}(\xi,s)}\int_{\mathbb{R}^3}\widehat{P_{k_1}U^{h_1,\iota_1}}(\xi-\eta,s)\widehat{P_{k_2}U^{h_2,\iota_2}}(\eta,s)n(\xi-\eta,\eta)\,d\eta ds\Big|,
\end{equation}

The sum over the pairs $(k_1,k_2)$ in \eqref{tar27} for which $\min(k_1,k_2)\leq -0.01m$ or $\max(k_1,k_2)\geq 0.01m$ can be bounded as claimed using \eqref{gene3x} and disregarding the cancellation in the multiplier $n$. On the other hand, if $|k_1|,|k_2|\leq 0.01m$ then we insert cutoff functions of the form $\varphi_{\leq q_0}(\Xi_{\iota_1\iota_2}(\xi-\eta,\eta))$ and $\varphi_{> q_0}(\Xi_{\iota_1\iota_2}(\xi-\eta,\eta))$, $q_0:=(-m+k-k_1-k_2)/2+\delta m/8$ in the integral in \eqref{tar28}. The contribution corresponding to $\varphi_{> q_0}(\Xi_{\iota_1\iota_2}(\xi-\eta,\eta))$ can be bounded using Lemma \ref{LemBaz1}. The contribution corresponding to $\varphi_{\leq q_0}(\Xi_{\iota_1\iota_2}(\xi-\eta,\eta))$ can be bounded using \eqref{gene2} and the null structure bound $|n_2(\xi-\eta,\eta)|\lesssim 2^{q_0}$ if $\Xi_{\iota_1\iota_2}(\xi-\eta,\eta)\lesssim 2^{q_0}$. 
\end{proof}

Finally, we bound the contributions of the terms $\HH_{\al\be}:=\mathcal{S}^{I,2}_{\al\be}$.

\begin{lemma}\label{TAR104} Assume that $m,k,t_1,t_2,\xi$ are as in Lemma \ref{TAR101}. Then 
\begin{equation}\label{tar31}
\Big|\varphi_k(\xi)\int_{t_1}^{t_2}e^{-i\Theta_{wa}(\xi,s)}e^{is\Lambda_{wa}(\xi)}\widehat{\mathcal{H}_{\al\be}}(\xi,s)\,ds\Big|\lesssim\varep_1^22^{\delta m/2}2^{-k^-}2^{-N_0k^+}.
\end{equation}
Moreover, for $G\in\{F,\Omega_a,\vartheta_{ab}\}$ we have
\begin{equation}\label{tar32}
\Big|\varphi_k(\xi)\int_{t_1}^{t_2}e^{-i\Theta_{wa}(\xi,s)}e^{is\Lambda_{wa}(\xi)}\widehat{\mathcal{H}^G}(\xi,s)\,ds\Big|\lesssim\varep_1^22^{-\delta m/2}2^{-N_0k^+},
\end{equation}
where
\begin{equation}\label{tar33}
\mathcal{H}^F:=\frac{1}{2}[\mathcal{H}_{00}+R_jR_k\mathcal{H}_{jk}],
\quad\mathcal{H}^{\Omega_a}:=\in_{akl}R_kR_m\mathcal{H}_{lm},
\quad\mathcal{H}^{\vartheta_{ab}}:=\in_{amp}\in_{bnq}R_mR_n\mathcal{H}_{pq}.
\end{equation}
\end{lemma}

\begin{proof} The bounds \eqref{tar32} follow as in Lemma \ref{TAR103} once we notice that $\mathcal{H}^F, \mathcal{H}^{\Omega_a}, \mathcal{H}^{\vartheta_{ab}}$ are defined by multipliers that suitable null structure bounds. Indeed, the identities \eqref{on5} show that the functions $\mathcal{H}^F$, $\mathcal{H}^{\Omega_a}$, $\mathcal{H}^{\vartheta_{ab}}$ are defined by the multipliers
\begin{equation*}
n^F_{\iota_1,\iota_2}(\theta,\eta)=C_{\iota_1\iota_2}\frac{(\theta\cdot\eta)^2}{|\theta|^2|\eta|^2}\Big(1-\iota_1\iota_2\frac{(\theta+\eta)\cdot\theta}{|\theta+\eta||\theta|}\frac{(\theta+\eta)\cdot\eta}{|\theta+\eta||\eta|}\Big),
\end{equation*}
\begin{equation*}
n^{\Omega_a}_{\iota_1,\iota_2}(\theta,\eta)=C_{\iota_1\iota_2}\frac{(\theta\cdot\eta)^2}{|\theta|^2|\eta|^2}\frac{(\theta+\eta)\cdot\eta}{|\theta+\eta||\eta|}\in_{akl}\frac{(\theta+\eta)_k\theta_l}{|\theta+\eta||\theta|},
\end{equation*}
\begin{equation*}
n^{\vartheta_{ab}}_{\iota_1,\iota_2}(\theta,\eta)=C_{\iota_1\iota_2}\frac{(\theta\cdot\eta)^2}{|\theta|^2|\eta|^2}\in_{amp}\in_{bnq}\frac{(\theta+\eta)_m\theta_p}{|\theta+\eta||\theta|}\frac{(\theta+\eta)_n\eta_q}{|\theta+\eta||\eta|},
\end{equation*}
where $\iota_1,\iota_2\in\{+,-\}$. We observe that these multipliers satisfy null structure bounds as before,  $|n^\ast_{\iota_1\iota_2}(\xi-\eta,\eta)|\lesssim 2^{q_0}$ if $\Xi_{\iota_1\iota_2}(\xi-\eta,\eta)\lesssim 2^{q_0}$ (see also \eqref{tar9}). 

We prove now \eqref{tar31}. Using the definitions and decomposing dyadically in frequency, it suffices to show that
\begin{equation}\label{tar37}
\sum_{(k_1,k_2)\in\mathcal{X}_k}Y'_{k_1,k_2}\lesssim \varep_1^22^{\delta m/2}2^{-k^-}2^{-N_0k^+},
\end{equation}  
where, with $\vartheta_1,\vartheta_2\in\{\vartheta_{ab}\}$, $\iota_1,\iota_2\in\{+,-\}$, and $m_1,m_2\in\mathcal{M}$,
\begin{equation}\label{tar38}
\begin{split}
Y'_{k_1,k_2}:=\Big|\varphi_k(\xi)\int_{t_1}^{t_2}&e^{is\Lambda_{wa}(\xi)-i\Theta_{wa}(\xi,s)}\int_{\mathbb{R}^3}m_1(\xi-\eta)m_2(\eta)\\
&\times\widehat{P_{k_1}U^{\vartheta_1,\iota_1}}(\xi-\eta,s)\widehat{P_{k_2}U^{\vartheta_2,\iota_2}}(\eta,s)\,d\eta ds\Big|.
\end{split}
\end{equation}

As in Lemma \ref{TAR103}, the sum over the pairs $(k_1,k_2)$ in \eqref{tar37} for which $\min(k_1,k_2)\leq -0.01m$ or $\max(k_1,k_2)\geq 0.01m$ can be bounded as claimed using \eqref{gene3x}. On the other hand, if $|k_1|,|k_2|\leq 0.01m$ then we insert cutoff functions of the form $\varphi_{\leq q_0}(\Xi_{\iota_1\iota_2}(\xi-\eta,\eta))$ and $\varphi_{> q_0}(\Xi_{\iota_1\iota_2}(\xi-\eta,\eta))$, $q_0:=(-m+k-k_1-k_2)/2+\delta m/8$ in the integral in \eqref{tar38}. The contribution corresponding to $\varphi_{> q_0}(\Xi_{\iota_1\iota_2}(\xi-\eta,\eta))$ can be bounded using Lemma \ref{LemBaz1}. Finally, the contribution corresponding to $\varphi_{\leq q_0}(\Xi_{\iota_1\iota_2}(\xi-\eta,\eta))$ can be bounded, using the estimates in the first line of \eqref{vcx1.2} and the fact that the support of the $\eta$-integral has volume $\lesssim 2^{\min(k_1,k_2)}2^{2q_0+2k_1+2k_2-2k}$ (see the sets $\mathcal{R}_{\leq b;\xi}$ defined in Lemma \ref{gene1}), by
\begin{equation*}
\begin{split}
C2^m2^{\min(k_1,k_2)}&2^{2q_0+2k_1+2k_2-2k}\cdot\varep_1^22^{-k_1^--\kappa k_1^-}2^{-N_0k_1^+}2^{-k_2^--\kappa k_2^-}2^{-N_0k_2^+}\\
&\lesssim \varep_1^22^{\delta m/4}2^{-k^-}(2^{\min(k_1,k_2)}2^{-\kappa k_1^-}2^{-\kappa k_2^-})2^{-k^++k_1^++k_2^+}2^{-N_0k_1^+}2^{-N_0k_2^+}.
\end{split}
\end{equation*}
This suffices to bound the remaining contributions, as claimed in \eqref{tar37}. Notice that this last estimate is tighter than before, and it relies on the strong bounds in \eqref{vcx1.2}, without $2^{C\delta m}$ losses, and on the weak null structure of the nonlinearity.
\end{proof}

\chapter{The main theorems}\label{Chapter6}

\section{Global regularity and asymptotics}\label{Proof1}

In this section we prove our first set of main theorems. All of the theorems below rely on Proposition \ref{bootstrap}, and some of the ingredients in its proof.

\subsection{Global regularity} We start with a quantitative global regularity result:

\begin{theorem}\label{winr} Assume that $(\overline{g}_{ij},k_{ij},\psi_0,\psi_1)$ is an initial data set on $\Sigma_0=\{(x,t)\in\mathbb{R}^4:t=0\}$ that satisfies the smallness conditions \eqref{asum1} and the constraint equations \eqref{asum2}. 

(i) Then the reduced Einstein-Klein-Gordon system
\begin{equation}\label{winr0.1}
\begin{split}
\widetilde{\square}_\g \g_{\al\be}+2\partial_\al\psi\partial_\be\psi+\psi^2\g_{\al\be}-F^{\geq 2}_{\al\be}(\g,\partial \g)&=0,\\
\widetilde{\square}_\g\psi-\psi&=0,
\end{split}
\end{equation}
admits a unique global solution $(\g,\psi)$ in $M=\{(x,t)\in\mathbb{R}^4:t\geq 0\}$, with initial data $(\overline{g}_{ij},k_{ij},\psi_0,\psi_1)$ on $\Sigma_0$ (as described in \eqref{data}) and satisfying $\|\g-m\|_{C^4(M)}+\|\psi\|_{C^4(M)}\lesssim\varep_0$. The solution also satisfies the harmonic gauge conditions in $M$
\begin{equation}\label{winr0.2}
0=\mathbf{\Gamma}_\mu=\g^{\al\be}\partial_\al \g_{\be\mu}-(1/2)\g^{\al\be}\partial_\mu \g_{\al\be},\qquad\mu\in\{0,1,2,3\}.
\end{equation}

(ii) Let $h_{\al\be}=\g_{\al\be}-m_{\al\be}$ as in \eqref{zaq1} and define the functions $U^\ast,V^\ast$ as in \eqref{variables4}--\eqref{variables4L}. For any $t\in[0,\infty)$, $\al,\be\in\{0,1,2,3\}$, and $a,b\in\{1,2,3\}$ we have
\begin{equation}\label{winr1}
\sup_{n\leq 3,\,\mathcal{L}\in\mathcal{V}_{n}^q}\langle t\rangle^{-H(q,n)\delta}\big\{\|(\langle t\rangle|\nabla|_{\leq 1})^{\gamma}|\nabla|^{-1/2}U^{\mathcal{L}h_{\al\be}}(t)\|_{H^{N(n)}}+\|U^{\mathcal{L}\psi}(t)\|_{H^{N(n)}}\big\}\lesssim\varep_0,
\end{equation}
\begin{equation}\label{winr2}
\begin{split}
\sup_{n\leq 2,\,\mathcal{L}\in\mathcal{V}_{n}^q}&\sup_{k\in\mathbb{Z},\,l\in\{1,2,3\}}\,2^{N(n+1)k^+}\langle t\rangle^{-H(q+1,n+1)\delta}\\
&\big\{2^{k/2}(2^{k^-}\langle t\rangle)^{\gamma}\|P_k(x_lV^{\mathcal{L}h_{\al\be}})(t)\|_{L^2}+2^{k^+}\|P_k(x_lV^{\mathcal{L}\psi})(t)\|_{L^2}\big\}\lesssim\varep_0,
\end{split}
\end{equation}
and
\begin{equation}\label{winr3}
\|V^{F}(t)\|_{Z_{wa}}+\|V^{\omega_a}(t)\|_{Z_{wa}}+\|V^{\vartheta_{ab}}(t)\|_{Z_{wa}}+\langle t\rangle^{-\delta}\|V^{h_{\al\be}}(t)\|_{Z_{wa}}+\|V^{\psi}(t)\|_{Z_{kg}}\lesssim\varep_0.
\end{equation}
where $H(q,n)$ is defined as in \eqref{fvc1.0}.

(iii) For any $k\in\mathbb{Z}$, $t\in[0,\infty)$, and $\mathcal{L}\in\mathcal{V}_n^q$, $n\leq 2$, the functions $\LL h_{\al\be}$ and $\LL\psi$ satisfy the pointwise decay bounds
\begin{equation}\label{winr4}
\|P_kU^{\mathcal{L}h_{\al\be}}(t)\|_{L^\infty}\lesssim \varep_0\langle t\rangle^{-1+\delta'/2}2^{k^-}2^{-N(n+1)k^++2k^+}\min\{1,\langle t\rangle 2^{k^-}\}^{1/2}
\end{equation}
and
\begin{equation}\label{winr5}
\|P_kU^{\mathcal{L}\psi}(t)\|_{L^\infty}\lesssim \varep_0\langle t\rangle^{-1+\delta'/2}2^{k^-/2}2^{-N(n+1)k^++2k^+}\min\{1,\langle t\rangle 2^{2k^-}\}^{1/2}.
\end{equation}
\end{theorem}

\begin{proof} {\bf{Step 1.}} We prove first suitable bounds on the functions $h_{\al\be}$ and $\psi$ for $t\in[0,2]$.  Indeed, notice first that, at time $t=0$,
\begin{equation}\label{IniDat6}
\begin{split}
&\sum_{|\beta'|\leq |\beta|+\beta_0-1\leq n}\||\nabla|^{-1/2}|\nabla|_{\leq 1}^\gamma(x^{\beta'}\partial^{\beta}\partial_0^{\beta_0}h_{\al\be})(0)\|_{H^{N(n)}}\lesssim\varep_0,\\
&\sum_{|\beta'|,|\beta|+\beta_0-1\leq n}\|(x^{\beta'}\partial^{\beta}\partial_0^{\beta_0}\psi)(0)\|_{H^{N(n)}}\lesssim\varep_0,
\end{split}
\end{equation}
for any $n\in[0,3]$ and $\beta_0\in\{0,1\}$, where $x^{\beta'}=x_1^{\beta'_1}x_2^{\beta'_2}x_3^{\beta'_3}$ and $\partial^{\beta}=\partial_1^{\beta_1}\partial_2^{\beta_2}\partial_3^{\beta_3}$. These bounds follow directly from \eqref{data} and \eqref{asum1}, by passing to the Fourier space.

We can now construct the functions $h_{\al\be}$ and $\psi$ by solving that coupled equations
\begin{equation}\label{winr9}
(\partial_0^2-\Delta)h_{\al\be}=\mathcal{N}^h_{\al\be},\qquad (\partial_0^2-\Delta+1)\psi=\mathcal{N}^\psi,
\end{equation}
see Proposition \ref{LemEKG}. Using standard energy estimates, similar to \eqref{EnergyIncrement1}-\eqref{EnergyIncrement1.1}, the solutions $h_{\al\be},\psi$ are well defined $C^2$ functions on $\mathbb{R}^3\times [0,2]$ and satisfy bounds similar to \eqref{IniDat6},
\begin{equation}\label{winr10}
\begin{split}
&\sum_{|\beta'|\leq |\beta|+\beta_0-1\leq n}\||\nabla|^{-1/2}|\nabla|_{\leq 1}^\gamma(x^{\beta'}\partial^{\beta}\partial_0^{\beta_0}h_{\al\be})(t)\|_{H^{N(n)}}\lesssim\varep_0,\\
&\sum_{|\beta'|,|\beta|+\beta_0-1\leq n}\|(x^{\beta'}\partial^{\beta}\partial_0^{\beta_0}\psi)(t)\|_{H^{N(n)}}\lesssim\varep_0,
\end{split}
\end{equation}
for any $t\in[0,2]$, $n\in[0,3]$, and $\beta_0\in\{0,1\}$. 

We would like to show now that the estimates \eqref{winr10} hold for all $\beta_0\in[0,3]$. Indeed, using the equations \eqref{winr9}, we can replace $\partial_0^2 h_{\al\be}$ with $\Delta h_{\al\be}+\mathcal{N}^h_{\al\be}$ and $\partial_0^2 \psi$ with $\Delta\psi-\psi+\mathcal{N}^\psi$. Simple product estimates using just \eqref{winr10} and the formulas in Proposition \ref{LemEKG} show that
\begin{equation*}
\begin{split}
&\sum_{|\beta'|\leq |\beta|+\beta_0-1\leq n,\,\beta_0\in\{2,3\}}\||\nabla|^{-1/2}|\nabla|_{\leq 1}^\gamma(x^{\beta'}\partial^{\beta}\partial_0^{\beta_0-2}\mathcal{N}^h_{\al\be})(t)\|_{H^{N(n)}}\lesssim \varep_0,\\
&\sum_{|\beta'|,|\beta|+\beta_0-1\leq n,\,\beta_0\in\{2,3\}}\|(x^{\beta'}\partial^{\beta}\partial_0^{\beta_0-2}\mathcal{N}^{\psi})(t)\|_{H^{N(n)}}\lesssim\varep_0.
\end{split}
\end{equation*} 
Therefore the bounds \eqref{winr10} hold for all $t\in[0,2]$ and any $n,\beta_0\in[0,3]$.

{\bf{Step 2.}} We prove now the bounds \eqref{winr1}--\eqref{winr3} for $t\in[0,2]$. Using \eqref{winr10} we have
\begin{equation}\label{winr11}
\begin{split}
&\||\nabla|^{-1/2}|\nabla|_{\leq 1}^\gamma(\nabla_{x,t}\LL h_{\al\be})(t)\|_{H^{N(n)}}\lesssim\varep_0,\\
&\|\LL \psi(t)\|_{H^{N(n)}}+\|\nabla_{x,t}\LL \psi(t)\|_{H^{N(n)}}\lesssim\varep_0,
\end{split}
\end{equation}
for any $t\in[0,2]$, $n\in[0,3]$, and $\LL\in\mathcal{V}^q_n$. The energy bounds \eqref{winr1} follow. The weighted bounds \eqref{winr2} then follow using Lemma \ref{ident} and \eqref{winr10} (for the nonlinear estimates), as in the proof on Proposition \ref{DiEs1}.

We prove now the $Z$ norm bounds \eqref{winr3}. It follows from \eqref{winr2} and \eqref{hyt3.1} that
\begin{equation}\label{IniDat2}
2^{N(2)k^+}2^{k/2}2^{\gamma k^-}2^{j}\|Q_{j,k}V^{\Omega h_{\al\be}}(t)\|_{L^2}\lesssim\varep_0,
\end{equation}
for any $t,\in[0,2]$, $k\in\mathbb{Z}$, and $\Omega\in\{\Omega_{23},\Omega_{31},\Omega_{12}\}$. Moreover
\begin{equation*}
2^{N(0)k^+}2^{-k/2}2^{\gamma k^-}\|Q_{j,k}V^{h_{\al\be}}(t)\|_{L^2}\lesssim\varep_0,
\end{equation*}
due to \eqref{winr1}. Using \eqref{consu2.1} it follows that $\|\widehat{P_k V^{h_{\al\be}}}\|_{L^\infty}\lesssim\varep_02^{-k-\gamma k^-}2^{-N_0k^+}$ for any $k\in\mathbb{Z}$, thus $\|V^{h_{\al\be}}(t)\|_{Z_{wa}}\lesssim\varep_0$ as desired.

Similarly, using \eqref{winr10} for any $t\in[0,2]$ we have
\begin{equation}\label{ku50}
\|U^{\psi}(t)\|_{H^{N(0)}}+\|\langle x\rangle^{2}U^{\psi}(t)\|_{H^{N(2)}}\lesssim\varep_0.
\end{equation}
In particular $\|P_kU^{\psi}(t)\|_{L^1}\lesssim\varep_0$, which gives $\|\widehat{P_kU^{\psi}}(t)\|_{L^\infty}\lesssim\varep_0$ for any $k\in\mathbb{Z}$. This suffices for $k\leq 0$. On the other hand, if $k\geq 0$ then it follows from \eqref{ku50} that
\begin{equation*}
\|P_kU^{\psi}(t)\|_{L^2}\lesssim\varep_02^{-N(0) k^+},\qquad \|\,|x|^2P_kU^{\psi}(t)\|_{L^2}\lesssim \varep_02^{-N(2) k^+}.
\end{equation*}
Thus $\|P_kU^{\psi}(t)\|_{L^1}\lesssim\varep_02^{-(N(0)+3N(2))/4 k^+}$, which gives the desired control $\|V^{\psi}(t)\|_{Z_{kg}}\lesssim\varep_0$. 

{\bf{Step 3.}} To summarize, given suitable initial data we constructed the solution $(h_{\al\be},\psi)$ of the system \eqref{winr9} on the time interval $[0,2]$ satisfying the bounds \eqref{winr10}. Letting $\g_{\al\be}=m_{\al\be}+h_{\al\be}$, the metric $\g$ (which is close to the Minkowski metric $m$) and the field $\psi$ satisfy the reduced Einstein-Klein-Gordon system \eqref{winr0.1} in $\mathbb{R}^3\times[0,2]$. The harmonic gauge condition \eqref{winr0.2} holds at time $t=0$ (due to the constraint equations \eqref{asum2}), therefore it holds in $\mathbb{R}^3\times[0,2]$ due to the reduced wave equations \eqref{tr6}. 

We apply now the main Proposition \ref{bootstrap}. A standard continuity argument shows that the solution $(h_{\al\be},\psi)$ can be extended globally in time, and satisfies the bootstrap bounds \eqref{winr1}--\eqref{winr3} and the pointwise smallness bounds $\|h_{\al\be}(t)\|_{L^\infty}\lesssim\varep_0$ for all $t\geq 0$. These pointwise bounds follow from \eqref{Linfty1} and \eqref{winr2} (compare with the proof of \eqref{wws1}) and are needed to justify that the metric $\g_{\al\be}$ is Lorentzian, so $\g^{\al\be}$ are well defined.

Finally, the bounds \eqref{winr4}--\eqref{winr5} are similar to the bounds \eqref{wws1} and \eqref{wws2}.  The weighted profile bounds \eqref{winr2} and the estimates \eqref{hyt3.1} show that
\begin{equation}\label{winr23}
2^{k/2}(2^{k^-}\langle t\rangle)^{\gamma}2^j\|Q_{j,k}V^{\mathcal{L}h}(t)\|_{L^2}+2^{k^+}2^j\|Q_{j,k}V^{\mathcal{L}\psi}(t)\|_{L^2}\lesssim \varep_02^{-N(n+1)k^+}\langle t\rangle^{H(q+1,n+1)\delta},
\end{equation}
for any $t\in[0,\infty)$, $(k,j)\in\mathcal{J}$, and  $\mathcal{L}\in\mathcal{V}_n^q$, $n\leq 2$. The desired bounds \eqref{winr4}--\eqref{winr5} follow from the linear estimates in Lemma \ref{LinEstLem}.
\end{proof}

\subsection{Decay of the metric components and the Klein-Gordon field} We prove now several estimates in the physical space. We introduce the tensor-fields 
\begin{equation}\label{winr20}
L:=\partial_t+\partial_r,\qquad\underline{L}:=\partial_t-\partial_r,\qquad \Pi^{\al\be}:=r^{-2}\big[\Omega_{12}^\al\Omega_{12}^\beta+\Omega_{23}^\al\Omega_{23}^\beta+\Omega_{31}^\al\Omega_{31}^\beta\big],
\end{equation}
where $r:=|x|$ and $\partial_r:=|x|^{-1}x^j\partial_j$. Notice that
\begin{equation}\label{Complete1}
m^{\alpha\beta}=-\frac{1}{2}\big\{L^\alpha\underline{L}^\beta+\underline{L}^\alpha L^\beta\big\}+\Pi^{\al\be}.
\end{equation}

Given a vector-field $V$ we define the (Minkowski) derivative operator $\partial_V:=V^\alpha\partial_\alpha$. Let
\begin{equation}\label{winr20.4}
\mathcal{T}:=\{L,r^{-1}\Omega_{12},r^{-1}\Omega_{23},r^{-1}\Omega_{31}\}
\end{equation}
denote the set of ``good" vector-fields, tangential to the (Minkowski) light cones. For $n\in\{0,1,2\}$ and $p\leq 6$ we define also the sets of differentiated metric components
\begin{equation}\label{winr25.6}
\mathcal{H}_{n,p}:=\big\{\partial_1^{a_1}\partial_2^{a_2}\partial_3^{a_3}\LL h_{\alpha\beta}:a_1+a_2+a_3\leq p,\,\LL\in\mathcal{V}^n_n,\,\al,\be\in\{0,1,2,3\}\big\}.
\end{equation}
We show first that the metric components and their derivatives have suitable decay,
\begin{equation}\label{winr25.7}
\vert h(x,t)\vert+\langle t+r\rangle\vert \partial_Vh(x,t)\vert+\langle t-r\rangle\vert\nabla_{x,t}h(x,t)\vert\lesssim \varepsilon_0\langle t+r\rangle^{2\delta'-1},
\end{equation}
where $V\in\mathcal{T}$ is a good vector-field and $h\in\{h_{\al\be}\}$. More precisely:

\begin{theorem}\label{Precisedmetric} Assume that $(\g,\psi)$ is a global solution of the Einstein-Klein-Gordon system as given by Theorem \ref{winr}.

(i) For any $H\in\mathcal{H}_{2,6}$ we have
\begin{equation}\label{winr25}
\vert H(x,t)\vert+\vert\partial_0H(x,t)\vert\lesssim \varepsilon_0\langle t+r\rangle^{\delta'/2-1}.
\end{equation}

(ii) Moreover, if $H'\in\mathcal{H}_{1,4}$, $H''\in\mathcal{H}_{0,3}$, $a\in\{1,2,3\}$, and $\Omega\in\{\Omega_{12},\,\Omega_{23},\,\Omega_{31}\}$ then
\begin{equation}\label{ImprovedBoundsInside}
\begin{split}
\langle t+r\rangle\vert \partial_{r^{-1}\Omega} H'(x,t)\vert+\langle r\rangle\vert \partial_{L} H'(x,t)\vert\vert&\lesssim \varepsilon_0\langle t+r\rangle^{2\delta'-1},\\
\langle t- r\rangle\vert\partial_aH'(x,t)\vert+\min(\langle r\rangle,\langle t- r\rangle)\vert\partial_0H'(x,t)\vert&\lesssim \varepsilon_0\langle t+r\rangle^{2\delta'-1},\\
\langle t+r\rangle\vert\partial_L H''(x,t)\vert+\langle t- r\rangle\vert\partial_0H''(x,t)\vert&\lesssim \varepsilon_0\langle t+r\rangle^{2\delta'-1}.
\end{split}
\end{equation}

(iii) The scalar field decays faster but with limited improvement: for $\Psi=\partial_1^{a_1}\partial_2^{a_2}\partial_3^{a_3}\LL_1\psi$ for some $\LL_1\in\mathcal{V}^1_1$ and $a_1+a_2+a_3\leq 4$ we have
\begin{equation}\label{winr25.4}
\begin{split}
\vert \Psi(x,t)\vert+\vert\partial_0\Psi(x,t)+\vert\langle\nabla\rangle\Psi(x,t)\vert&\lesssim\varepsilon_0\langle t+r\rangle^{\delta'/2-1}\langle r\rangle^{-1/2},\\
\vert \partial_b\Psi(x,t)\vert&\lesssim\varepsilon_0\langle t+r\rangle^{\delta'/2-3/2},\qquad b\in\{1,2,3\}.
\end{split}
\end{equation}
\end{theorem}

\begin{proof} {\bf{Step 1.}} We prove first the bounds \eqref{winr25}, using the profile bounds \eqref{winr23} and linear estimates. With $\LL_2\in\mathcal{V}_2^2$ and $H=\partial^a\LL_2 h_{\al\be}\in\mathcal{H}_{2,6}$, $\partial^a=\partial_1^{a_1}\partial_2^{a_2}\partial_3^{a_3}$, we have
\begin{equation}\label{winr28}
|P_kH(x,t)|+|P_k\partial_0H(x,t)|\lesssim\sum_{R\in\{|\nabla|^{-1},\mathrm{Id}\}}\sum_{j\geq -k^-}\big|\big(\partial^a RP'_k(e^{-it\Lambda_{wa}}Q_{j,k}V^{\LL_2 h_{\al\be}})\big)(x,t)\big|
\end{equation}
for any $k\in\mathbb{Z}$ and $(x,t)\in M$. Using \eqref{Linfty1} and \eqref{winr23} we thus estimate
\begin{equation}\label{winr29}
\begin{split}
\|P_kH(t)\|_{L^\infty}+\|P_k\partial_0H(t)\|_{L^\infty}&\lesssim 2^{7k^+}2^{k/2}\sum_{j\geq -k^{-}}\min(1,2^{j}\langle t\rangle^{-1})\|Q_{j,k}V^{\LL_2 h_{\al\be}}(t)\|_{L^2}\\
&\lesssim \varep_0\langle t\rangle^{-1+H(3,3)\delta+\delta}2^{-2k^+}\min(1,2^{k^-}\langle t\rangle)^{1-\delta}.
\end{split}
\end{equation}

The bounds \eqref{winr25} follow from \eqref{winr29} if $r\lesssim \langle t\rangle$. On the other hand, if $r=|x|\geq 4\langle t\rangle$ then we still use \eqref{winr28} first and notice that the contribution of the pairs $(k,j)$ with $2^k\geq 2^{-10}r^{1/2}$ or $2^j\geq  2^{-10}r^{1-\delta}$ can still be bounded as in \eqref{winr29}. On the other hand, if $2^k,2^j\leq 2^{-10}|x|^{1-\delta}$ and $m\in\mathcal{M}_0$ then we have rapid decay, 
\begin{equation}\label{winr30}
\Big|\int_{\mathbb{R}^3}e^{-it|\xi|}e^{ix\cdot\xi}m(\xi)\widehat{V^{\LL_2 h_{\al\be}}_{j,k}}(\xi,t)\,d\xi\Big|\lesssim \varep_0|x|^{-10}2^{-10k^+},
\end{equation}
using integration by parts in $\xi$ (Lemma \ref{tech5}). The desired bounds \eqref{winr25} follow. 

{\bf{Step 2.}} We consider now the bounds \eqref{ImprovedBoundsInside}, which we prove in several stages (more precisely, the bounds \eqref{ImprovedBoundsInside} follow from \eqref{winr35}, \eqref{winr32}, \eqref{winr37}, and \eqref{winr39}). 

We prove first that if $H'\in\mathcal{H}_{1,4}$, $r=|x|\geq 2\langle t\rangle$, and $\mu\in\{0,1,2,3\}$ then
\begin{equation}\label{winr35}
\vert\partial_\mu H'(x,t)\vert\lesssim \varepsilon_0 r^{\delta'-2}.
\end{equation}
Indeed, as in \eqref{winr28} we estimate
\begin{equation}\label{winr36}
|\partial_\mu H'(x,t)|\leq\sum_{(k,j)\in\mathcal{J}}\sum_{R\in\{R_a,\mathrm{Id}\}}\big|\big(R\partial_1^{a_1}\partial_2^{a_2}\partial_3^{a_3}P'_k(e^{-it\Lambda_{wa}}Q_{j,k}V^{\LL_1 h_{\al\be}})\big)(x,t)\big|.
\end{equation}
Using now \eqref{winr40}, for any $k\in\mathbb{Z}$ we have
\begin{equation*}
|x|\big|P_k(\partial_\mu H')(x,t)\big|\lesssim (1+2^k|x|)^\delta 2^{k/2}2^{4k^+}\sum_{j\geq -k^-}\|Q_{j,k}V^{\LL_1 h_{\al\be}}(t)\|_{H^{0,1}_\Omega}.
\end{equation*}
The sum over the pairs $(k,j)\in\mathcal{J}$ with $2^j\gtrsim r^{1-\delta}$ is controlled as claimed using \eqref{winr23}. On the other hand, since $r\geq 2\langle t\rangle$, the sum over the pairs  $(k,j)\in\mathcal{J}$ with $2^j\leq 2^{-10} r^{1-\delta}$ is negligible, as in \eqref{winr30}. The bounds \eqref{winr35} follow.

We show now that if $H'\in\mathcal{H}_{1,4}$ and $r\leq 4\langle t\rangle$ then
\begin{equation}\label{winr32}
\vert r^{-1}\Omega_{ab} H'(x,t)\vert\lesssim \varepsilon_0\langle t\rangle^{\delta'-2}.
\end{equation}
Indeed, this follows using the identities
\begin{equation}\label{winr31}
t\Omega_{ab}=x_a\Gamma_b-x_b\Gamma_a,
\end{equation}
and the bounds \eqref{winr25} applied to the functions $\Gamma H'$.

We prove now that if $H'\in\mathcal{H}_{1,4}$ and $r\leq 4\langle t\rangle$ then
\begin{equation}\label{winr36.5}
|\partial_LH'(x,t)|\lesssim\varep_0\langle r\rangle^{-1}\langle t\rangle^{-1+\delta'}.
\end{equation}
This follows from \eqref{winr25} if $r\lesssim 1$. On the other hand, if $2^{4}\leq r\leq 4\langle t\rangle$ then we write
\begin{equation}\label{winr36.7}
|x|\partial_LH'(x,t)=\frac{1}{2}\big\{|x|\partial^a(U^{\LL_1h_{\al\be}}+\overline{U^{\LL_1h_{\al\be}}})(x,t)+ix^bR_b\partial^a(U^{\LL_1h_{\al\be}}-\overline{U^{\LL_1h_{\al\be}}})(x,t)\big\},
\end{equation}
where $\partial^a=\partial_1^{a_1}\partial_2^{a_2}\partial_3^{a_3}$, $\mathcal{L}_1\in\mathcal{V}_1^1$, $H'=\partial^a\mathcal{L}_1h_{\al\be}$. Therefore, using the bounds \eqref{winr41} below,
\begin{equation}\label{winr36.8}
\begin{split}
\big||x|\partial_LH'(x,t)\big|&\leq \big|(|x|+ix^bR_b)\partial^a U^{\LL_1h_{\al\be}}(x,t)\big|\\
&\lesssim\sum_{2^k\leq t^{-1}}t2^{3k/2}\|P_kU^{\LL_1h_{\al\be}}(t)\|_{L^2}\\
&+\sum_{2^k\geq t^{-1},\,(k,j)\in\mathcal{J}}2^{4k^+}(2^{k}t)^{2\delta}2^{k/2}\min(1,2^jt^{-1})\|Q_{j,k}V^{\LL_1h_{\al\be}}(t)\|_{H^{0,1}_\Omega}.
\end{split}
\end{equation}
The desired bounds \eqref{winr36.5} now follow from \eqref{winr1} and \eqref{winr23}.

We prove now that if $H'\in\mathcal{H}_{1,4}$, $r\leq 4\langle t\rangle$, and $a\in\{1,2,3\}$ then
\begin{equation}\label{winr37}
\vert\partial_a H'(x,t)\vert\lesssim \varepsilon_0 \langle t\rangle^{\delta'-1}\langle t-r\rangle^{-1}.
\end{equation}
Indeed, for this we use the identity
\begin{equation}\label{winr37.2}
(r-t)\partial_r=rL-r^{-1}x^b\Gamma_b.
\end{equation}
Using now \eqref{winr36.5} and \eqref{winr25} it follows that $|t-r||\partial_rH'(x,t)|\lesssim  \varepsilon_0 \langle t\rangle^{\delta'-1}$. Using also \eqref{winr32} it follows that $|t-r||\partial_aH'(x,t)|\lesssim  \varepsilon_0 \langle t\rangle^{\delta'-1}$, $a\in\{1,2,3\}$. The bounds \eqref{winr37} follow using also \eqref{winr25} in the case $|t-r|\lesssim 1$. 

{\bf{Step 3.}} Finally, we prove that if $H''\in\mathcal{H}_{0,3}$ and $r\leq 4\langle t\rangle$ then 
\begin{equation}\label{winr39}
|\partial_LH''(x,t)|\lesssim\varep_0\langle t\rangle^{-2+2\delta'}.
\end{equation}
This follows from \eqref{winr36.5} $t\lesssim 1$ or if $r\gtrsim t$. 

On the other hand, to prove the bounds when $t\geq 2^{8}$ and $r\leq t/8$ we notice that
\begin{equation}\label{winr39.1}
\partial_{\underline{L}}\partial_L=\partial_0^2-\partial_r^2=\partial_0^2-\Delta+(2/r)\partial_r+r^{-2}(\Omega_{12}^2+\Omega_{23}^2+\Omega_{31}^2),
\end{equation}
where $\underline{L}=\partial_t-\partial_r$. We apply $\partial_{\underline{L}}$ to $\partial_LH''$ and use the wave equations \eqref{zaq11.1}. Thus
\begin{equation}\label{winr39.2}
\partial_{\underline{L}}\partial_LH''=W'':=\partial^a\mathcal{N}^h_{\al\be}+(2/r)\partial_r(\partial^ah_{\al\be})+r^{-2}(\Omega_{12}^2+\Omega_{23}^2+\Omega_{31}^2)(\partial^ah_{\al\be}).
\end{equation}

We prove now that if $s\geq 2^6$ and $|y|\leq s/2$ then
\begin{equation}\label{winr39.4}
|W''(y,s)|\lesssim \varep_0|y|^{-1}\langle s\rangle^{-2+3\delta'/2}.
\end{equation}
Indeed, for the nonlinear terms in $\partial^a\mathcal{N}^h_{\al\be}$ this follows from the formula \eqref{sac1} and the bounds \eqref{winr25}, \eqref{winr25.4} (which is proved below), \eqref{winr36.5}, and \eqref{winr37}. For the term $(2/r)\partial_r(\partial^ah_{\al\be})$ these bounds follow directly from \eqref{winr37}, while for the term $r^{-2}(\Omega_{12}^2+\Omega_{23}^2+\Omega_{31}^2)(\partial^ah_{\al\be})$ the bounds \eqref{winr39.4} follow using \eqref{winr31} and \eqref{winr25}.

We can now use the identity $\partial_{\underline{L}}\partial_LH''=W''$ and integrate along the vector-field $\underline{L}$ to complete the proof of \eqref{winr39}. Indeed, for any function $F$ and $\lambda\in [0,t/2]$ we have 
\begin{equation*}
F(x,t)-F\big(x+\lambda \frac{x}{|x|},t-\lambda\big)=\int_0^\lambda-\frac{d}{ds}F\big(x+s \frac{x}{|x|},t-s\big)\,ds=\int_0^\lambda (\partial_{\underline{L}}F)\big(x+s \frac{x}{|x|},t-s\big)\,ds.
\end{equation*}
We apply this with $F=\partial_LH''$ and $\lambda=t/8$. Since $|F(x+\lambda x/|x|,t-\lambda)|\lesssim \varep_0t^{-2+\delta'}$ (due to \eqref{winr36.5}) and $|(\partial_{\underline{L}}F)(x+sx/|x|,t-s)|\lesssim \varep_0(s+|x|)^{-1}t^{-2+3\delta'/2}$ for any $s\in[0,\lambda]$ (due to \eqref{winr39.4}), it follows that
\begin{equation}\label{winr39.6}
|\partial_LH''(x,t)|\lesssim \varep_0t^{-2+3\delta'/2}\ln(t/|x|)\qquad\text{ for }t\geq 2^{8}\text{ and }|x|\leq t/8.
\end{equation}

The bounds \eqref{winr39} follow if $|x|\geq t^{-4}$. Moreover, if $|x'|\leq 2t^{-4}$ then $|\partial_rH''(x',t)|\lesssim \varep_0t^{-2+\delta'}$ (due to \eqref{winr37}) and $|\partial_a\partial_0H''(x',t)|\lesssim \varep_0t^{-1+\delta'}$, $a\in\{1,2,3\}$ (due to \eqref{winr25}). In view of \eqref{winr39.6} it follows that $|\partial_0H''(x,t)|\lesssim \varep_0t^{-2+2\delta'}$ if $|x|\leq 2t^{-4}$, and the desired bounds \eqref{winr39} follow for all $x$ with $|x|\leq t/8$.

{\bf{Step 4.}} We prove now the bounds \eqref{winr25.4} on the scalar field. These bounds follow in the same way as \eqref{winr35} if $r\geq 2\langle t\rangle$ or if $r+t\lesssim 1$.

It remains to consider the case $t\geq 8$ and $r\leq 4t$. Notice that $t\partial_b\Psi=\Gamma_b\Psi-x_b\partial_0\Psi$, so the bounds $|\partial_b\Psi(x,t)|\lesssim \varep_0t^{-3/2+\delta'/2}$ in \eqref{winr25.4} follow from the bounds in the first line (for $\partial_0\psi$), and the $L^\infty$ estimates \eqref{wws2} (for $\Gamma_b\psi$). To summarize, it remains to prove that
\begin{equation}\label{winv1}
\vert \Psi(x,t)\vert+\vert\partial_0\Psi(x,t)\vert+\vert\langle\nabla\rangle\Psi(x,t)\vert\lesssim\varepsilon_0t^{\delta'/2-1}\langle r\rangle^{-1/2}\qquad\text{ if }t\geq 8\text{ and }r\leq 4t.
\end{equation}

These bounds follow from \eqref{wws2} if $|r|\lesssim 1$. On the other hand, if $r\geq 1$ then we estimate
\begin{equation}\label{winv2}
\begin{split}
|P_k\Psi(x,t)|+&|P_k\partial_0\Psi(x,t)|+\vert P_k\langle\nabla\rangle\Psi(x,t)\vert\\
&\lesssim\sum_{R\in\{\Lambda_{kg}^{-1},\mathrm{Id}\}}\sum_{j\geq -k^-}\big|\big(\partial^a RP'_k(e^{-it\Lambda_{kg}}Q_{j,k}V^{\LL_1\psi})\big)(x,t)\big|,
\end{split}
\end{equation}
where $\Psi=\partial^a\LL_1 \psi$, $\LL_1\in\mathcal{V}_1^1$. Using \eqref{winr23} we estimate
\begin{equation}\label{winv3}
\begin{split}
\big|\big(\partial^a RP'_k(e^{-it\Lambda_{kg}}Q_{j,k}V^{\LL_1\psi})\big)(x,t)\big|&\lesssim 2^{4k^+}2^{3k/2}\|Q_{j,k}V^{\LL_1\psi}(t)\|_{L^2}\\
&\lesssim \varep_02^{-5k^+}2^{3k/2}2^{-j}t^{\delta'/4}.
\end{split}
\end{equation}

If $2^k\lesssim t^{-1/2}$ then we use \eqref{winv3} to estimate the contribution of the pairs $(k,j)$ with $2^j\geq r^{1-\delta}2^{-10}$. On the other hand, if $2^j\leq r^{1-\delta}2^{-10}$ then we have rapid decay,
\begin{equation}\label{winv4}
\big|\big(\partial^a RP'_k(e^{-it\Lambda_{kg}}Q_{j,k}V^{\LL_1\psi})\big)(x,t)\big|\lesssim \varep_02^{3k/2}r^{-4},
\end{equation}
using integration by parts in the Fourier space (Lemma \ref{tech5}). Therefore, using \eqref{winv2},
\begin{equation}\label{winv5}
|P_k\Psi(x,t)|+|P_k\partial_0\Psi(x,t)|\lesssim \varep_0(2^{2k}t)r^{-1/2}t^{\delta'/3-1}.
\end{equation}

Assume now that $2^{k-40}\geq t^{-1/2}$. Using \eqref{Linfty3.6} and \eqref{winr23} we have
\begin{equation}\label{winv9}
\begin{split}
\big|\big(\partial^a RP'_k(e^{-it\Lambda_{kg}}Q_{j,k}V^{\LL_1\psi})\big)(x,t)\big|&\lesssim 2^{5k^+}t^{-3/2}2^{j/2-k^-}t^\delta\|Q_{j,k}V^{\LL_1\psi}(t)\|_{H^{0,1}_\Omega}\\
&\lesssim\varep_02^{-2k^+}t^{-3/2+\delta'/3}2^{-k/2}
\end{split}
\end{equation}
provided that $2^j\leq 2^{k^--20}t$. Moreover, if $r^{1-\delta}\geq t2^{k+20}$ and $2^j\leq 2^{k^--20}t$ then 
\begin{equation*}
\big|\big(\partial^a RP'_k(e^{-it\Lambda_{kg}}Q_{j,k}V^{\LL_1\psi})\big)(x,t)\big|\lesssim \varep_02^{-2k^+}r^{-6}
\end{equation*}
using integration by parts in the Fourier space (Lemma \ref{tech5}). Therefore, using \eqref{winv9},
\begin{equation}\label{winv7}
\sum_{2^j\in[2^{-k^-},2^{k^--20}t]}\big|\big(\partial^a RP'_k(e^{-it\Lambda_{kg}}Q_{j,k}V^{\LL_1\psi})\big)(x,t)\big|\lesssim \varep_02^{-2k^+}t^{-1+2\delta'/5}r^{-1/2}.
\end{equation}

To bound the contribution of the sum over $j$ large with $2^j\geq 2^{k^--20}t$ we notice that
\begin{equation}\label{winv8}
\begin{split}
\big|\big(\partial^a RP'_k(e^{-it\Lambda_{kg}}Q_{j,k}V^{\LL_1\psi})\big)(x,t)\big|&\lesssim r^{-1}(1+2^kr)^\delta 2^{k/2}2^{4k^+}\|Q_{j,k}V^{\LL_1\psi}\|_{H^{0,1}_\Omega}\\
&\lesssim \varep_0r^{-1}2^{-4k^+}2^{k/2}2^{-j}t^{\delta'/2-4\delta}.
\end{split}
\end{equation}
as a consequence of \eqref{winr40} and \eqref{winr23}. We use now \eqref{winv3} if $r\leq 2^{-k}$ and \eqref{winv8} if $r\geq 2^{-k}$ to conclude that
\begin{equation*}
\sum_{2^j\geq 2^{k^--20}t}\big|\big(\partial^a RP'_k(e^{-it\Lambda_{kg}}Q_{j,k}V^{\LL_1\psi})\big)(x,t)\big|\lesssim \varep_02^{-2k^+}t^{-1+\delta'/2-4\delta}r^{-1/2}.
\end{equation*}

Using also \eqref{winv7} and \eqref{winv2}, if $2^{k-40}\geq t^{-1/2}$ we have
\begin{equation*}
|P_k\Psi(x,t)|+|P_k\partial_0\Psi(x,t)|+\vert P_k\langle\nabla\rangle\Psi(x,t)\vert\lesssim \varep_02^{-2k^+}t^{-1+\delta'/2-4\delta}r^{-1/2}.
\end{equation*}
The bounds \eqref{winv1} follow if $r\geq 1$ by summation over $k$, using also \eqref{winv5}.
\end{proof}

Using also the harmonic gauge condition \eqref{winr0.2} we can prove some additional bounds on the derivatives of the metric $h_{\al\be}$. More precisely:

\begin{lemma}\label{extrader} With $(\g,\psi)$ as in Theorem \ref{winr} and $\Pi^{\al\be}$ defined as in \eqref{winr20}, we have the additional bounds
\begin{equation}\label{extrader1}
\big|\underline{L}^\mu V^\al L^\be\partial_\mu (\partial^ah_{\al\be})(x,t)\big|+\big|\underline{L}^\mu\Pi^{\al\be}\partial_\mu (\partial^ah_{\al\be})(x,t)\big|\lesssim \varep_0\langle t+r\rangle^{-2+3\delta'}
\end{equation}
for any $(x,t)\in M$, $V\in\mathcal{T}$, and $\partial^a=\partial_1^{a_1}\partial_2^{a_2}\partial_3^{a_3}$, $a_1+a_2+a_3\leq 3$.
\end{lemma}

\begin{proof} We write the harmonic gauge condition in the form
\begin{equation*}
m^{\al\be}\partial_\al h_{\be\mu}-(1/2)m^{\al\be}\partial_\mu h_{\al\be}=-g_{\geq 1}^{\al\be}\partial_\al h_{\be\mu}+(1/2)g_{\geq 1}^{\al\be}\partial_\mu h_{\al\be}.
\end{equation*}
Using \eqref{winr25}--\eqref{ImprovedBoundsInside}, it follows that, for $\mu\in\{0,1,2,3\}$,
\begin{equation}\label{extrader5}
m^{\al\be}\big\{\partial_\al (\partial^ah_{\be\mu})-(1/2)\partial_\mu (\partial^ah_{\al\be})\big\}=O(\varep_0\langle t+r\rangle^{-2+3\delta'}),
\end{equation}
where $f=O(g)$ means $|f(x,t)|\lesssim g(x,t)$ for all $(x,t)\in M$. We use now the formula \eqref{Complete1}, and eliminate some of the terms using \eqref{ImprovedBoundsInside}, to conclude that
\begin{equation}\label{extrader6}
L^\be\underline{L}^\al\partial_\al (\partial^ah_{\be\mu})-\frac{1}{2}\big\{L^\alpha\underline{L}^\beta+\underline{L}^\alpha L^\beta\big\}\partial_\mu (\partial^ah_{\al\be})+\Pi^{\al\be}\partial_\mu (\partial^ah_{\al\be})=O(\varep_0\langle t+r\rangle^{-2+3\delta'}),
\end{equation}
for $\mu\in\{0,1,2,3\}$. The desired conclusions in \eqref{extrader1} follow by multiplying with either $V^\mu$, $V\in\mathcal{T}$, or with $\underline{L}^\mu$. 
\end{proof}

We prove now almost sharp bounds on two derivatives of the metric tensor, in the region $\{|x|\gtrsim t\gtrsim 1\}$. These bounds are used later to prove weak peeling estimates.

\begin{lemma}\label{twoder0}
Assume that $(\g,\psi)$ is a global solution of the Einstein-Klein-Gordon system as given by Theorem \ref{winr}. If $V_1,V_2\in\mathcal{T}$ and $H\in\mathcal{H}_{0,3}$ (see \eqref{winr20.4}--\eqref{winr25.6}) then
\begin{equation}\label{twoder1}
\langle t-r\rangle^2|\partial_{\underline{L}}^2H(x,t)|+\langle t-r\rangle\langle r\rangle|\partial_{\underline{L}}\partial_{V_1}H(x,t)|+\langle r\rangle^2|\partial_{V_2}\partial_{V_1}H(x,t)|\lesssim \varep_0\langle r\rangle^{-1+3\delta'},
\end{equation}
for any $(x,t)\in M':=\{(x,t)\in M:\,t\geq 1\text{ and }|x|\geq 2^{-8}t\}$. 
\end{lemma}

\begin{proof} {\bf{Step 1.}} The bounds are easy when either $V_1$ or $V_2$ is a rotation vector-field. Indeed, if $V_1=r^{-1}\Omega_{ab}$ then $\partial_{V_1}H$ is a sum of functions of the form $r^{-1}H'$, $H'\in\mathcal{H}_{1,3}$, so 
\begin{equation}\label{twoder2}
\langle t-r\rangle\langle r\rangle|\partial_{\underline{L}}\partial_{V_1}H(x,t)|+\langle r\rangle^2|\partial_{V_2}\partial_{V_1}H(x,t)|\lesssim \varep_0\langle r\rangle^{-1+2\delta'},
\end{equation}
as a consequence of \eqref{ImprovedBoundsInside}. Moreover, if $V_2=r^{-1}\Omega_{ab}$ then the commutators $[V_2,\underline{L}]$ and $[V_2,V_1]$ are sums of vector-fields of the form $r^{-1}W$, $W\in\mathcal{T}$. Therefore, using also \eqref{twoder2},
\begin{equation}\label{twoder3}
\langle t-r\rangle\langle r\rangle|\partial_{V_2}\partial_{\underline{L}}H(x,t)|+\langle r\rangle^2|\partial_{V_2}\partial_{V_1}H(x,t)|\lesssim \varep_0\langle r\rangle^{-1+2\delta'},
\end{equation}

We show now that 
\begin{equation}\label{twoder4}
\langle t-r\rangle\langle r\rangle|\partial_{\underline{L}}\partial_{L}H(x,t)|\lesssim \varep_0\langle r\rangle^{-1+2\delta'}.
\end{equation}
Recalling that $H=\partial^ah_{\al\be}$ and using the formula \eqref{winr39.2}, we have
\begin{equation}\label{twoder5}
\partial_{\underline{L}}\partial_LH=\partial^a\mathcal{N}^h_{\al\be}+(2/r)\partial_r(\partial^ah_{\al\be})+r^{-2}(\Omega_{12}^2+\Omega_{23}^2+\Omega_{31}^2)(\partial^ah_{\al\be}).
\end{equation}
As in the proof of \eqref{winr39.4}, it is easy to see that
\begin{equation*}
\begin{split}
(2/r)|\partial_r(\partial^ah_{\al\be})(x,t)|&\lesssim \langle t-r\rangle^{-1}\langle r\rangle^{-2+2\delta'},\\
r^{-2}|(\Omega_{12}^2+\Omega_{23}^2+\Omega_{31}^2)(\partial^ah_{\al\be})(x,t)|&\lesssim\langle r\rangle^{-3+2\delta'},
\end{split}
\end{equation*}
using \eqref{ImprovedBoundsInside} and \eqref{winr25}. The nonlinearity $\partial^a\mathcal{N}^h_{\al\be}(x,t)$ can be estimated using the formula \eqref{sac1} and the bounds \eqref{winr25}--\eqref{winr25.4}. The desired bounds \eqref{twoder4} follow.

{\bf{Step 2.}} For \eqref{twoder1} it remains to prove that
\begin{equation}\label{twoder7}
\langle t-r\rangle^2|\partial_{\underline{L}}^2H(x,t)|+\langle r\rangle^2|\partial_{L}^2H(x,t)|\lesssim \varep_0\langle r\rangle^{-1+3\delta'},
\end{equation}
for any $(x,t)\in M'$. We define the vector-field
\begin{equation}\label{twoder8}
\widetilde{\Gamma}:=r^{-1}x^b\Gamma_b=r\partial_t+t\partial_r
\end{equation}
in $M'$, and notice that
\begin{equation}\label{twoder9}
\widetilde{\Gamma}=(1/2)[(r+t)L+(r-t)\underline{L}].
\end{equation}
Moreover, using \eqref{ImprovedBoundsInside},
\begin{equation}\label{twoder10}
|\partial_L(\widetilde{\Gamma}H)(x,t)|\lesssim \varep_0\langle r\rangle^{-2+2\delta'},\qquad |\partial_{\underline{L}}(\widetilde{\Gamma}H)(x,t)|\lesssim \varep_0\langle r\rangle^{-1+2\delta'}\langle t-r\rangle^{-1}.
\end{equation}
Using \eqref{twoder9} it follows that $|\partial_L[(r+t)\partial_L+(r-t)\partial_{\underline{L}}] H(x,t)|\lesssim \varep_0\langle r\rangle^{-2+2\delta'}$, thus
\begin{equation*}
|(r+t)\partial_L^2H(x,t)|\lesssim \varep_0\langle r\rangle^{-2+2\delta'}+|\partial_LH(x,t)|+\langle t-r\rangle |\partial_L\partial_{\underline{L}}H(x,t)|.
\end{equation*}
The bound on $|\partial_{L}^2H(x,t)|$ in \eqref{twoder7} follows using also \eqref{ImprovedBoundsInside} and \eqref{twoder4}. 

Similarly, using \eqref{twoder9}--\eqref{twoder10} we have 
\begin{equation*}
|\partial_{\underline{L}}[(r-t)\partial_{\underline{L}}+(r+t)\partial_L] H(x,t)|\lesssim \varep_0\langle r\rangle^{-1+2\delta'}\langle t-r\rangle^{-1},
\end{equation*}
therefore
\begin{equation*}
|(r-t)\partial_{\underline{L}}^2H(x,t)|\lesssim \varep_0\langle r\rangle^{-1+2\delta'}\langle t-r\rangle^{-1}+|\partial_{\underline{L}}H(x,t)|+\langle r\rangle |\partial_{\underline{L}}\partial_LH(x,t)|.
\end{equation*}
The bound on $|\partial_{\underline{L}}^2H(x,t)|$ in \eqref{twoder7} follows using also \eqref{ImprovedBoundsInside} and \eqref{twoder4} when $|t-r|\geq 1$, or using \eqref{winr25} when $|t-r|\leq 1$. This completes the proof of the lemma.
\end{proof}

We prove now the additional linear estimates we used in Theorem \ref{Precisedmetric}.

\begin{lemma}\label{addiline}
(i) For any $f$ in $L^2(\mathbb{R}^3)$, $x\in\mathbb{R}^3$, and $k\in\mathbb{Z}$ we have
\begin{equation}\label{winr40}
\big||x|P_kf(x)\big|\lesssim (1+2^{k}|x|)^{\delta}2^{k/2}\|P_kf\|_{H^{0,1}_\Omega}.
\end{equation}

(ii) In addition, if $t\geq 1$, $|x|\leq 8t$, $(k,j)\in\mathcal{J}$ then
\begin{equation}\label{winr41}
\big|(|x|+ix^aR_a)T_m(e^{-it\Lambda_{wa}}f_{j,k})(x)\big|\lesssim (1+2^{k}t)^{2\delta}2^{k/2}\min(1,2^jt^{-1})\|f^\ast_{j,k}\|_{H^{0,1}_\Omega}.
\end{equation}
Here $f_{j,k}=P'_kQ_{j,k}f$, $f^\ast_{j,k}=Q_{j,k}f$ are as in \eqref{defin}--\eqref{defin2}, $R_a=|\nabla|^{-1}\partial_a$ are the Riesz transforms, and the linear operators $T_m$ are defined by $T_mg=\mathcal{F}^{-1}(m\cdot \widehat{g})$, where $m\in\mathcal{M}_0$. 

Finally, if $|x|\in[2^{-20}t,8t]$ and $2^j\leq t(1+2^kt)^{-\delta}$ then
\begin{equation}\label{winr41.2}
\big|(|x|+ix^aR_a)T_m(e^{-it\Lambda_{wa}}f_{j,k})(x)\big|\lesssim (1+2^{k}t)^{2\delta}2^{-k/2}t^{-1}\|f_{j,k}\|_{H^{0,2}_\Omega}.
\end{equation}
\end{lemma}

\begin{proof} (i) Clearly $\|P_kf\|_{L^\infty}\lesssim 2^{3k/2}\|P_kf\|_{L^2}$, so the bounds \eqref{winr40} follow if $|x|\lesssim 2^{-k}$. On the other hand, if $|x|\geq 2^{-k+20}$ then we estimate
\begin{equation}\label{winr41.5}
|P_kf(x)|\lesssim \int_{\mathbb{R}^3}|P_kf(y)|\cdot 2^{3k}(1+2^{k}|x-y|)^{-8}\,dy.
\end{equation}
Using the Sobolev embedding along the spheres $\mathbb{S}^2$, for any $g\in H^{0,1}_\Omega$ and $p\in[2,\infty)$ we have
\begin{equation}\label{winr43}
\big\Vert g(r\theta)\,\big\Vert_{L^2(r^2dr)L^p_\theta}
\lesssim_p \sum_{m_1+m_2+m_3\leq 1}\Vert \Omega_{23}^{m_1}\Omega_{31}^{m_2}\Omega_{12}^{m_3}g\Vert_{L^2}\lesssim_p \Vert g\Vert_{H^{0,1}_\Omega}.
\end{equation}
Therefore, for $x\in\mathbb{R}^3$ with $|x|\geq 2^{-k+20}$ we can estimate the right-hand side of \eqref{winr41.5} by
\begin{equation*}
\begin{split}
C\|P_kf(r\theta)\|_{L^2(r^2dr)L^p_\theta}\|2^{3k}(1+2^k|x-r\theta|)^{-8}\|_{L^2(r^2dr)L^{p'}_\theta}\lesssim_p \|P_kf\|_{H^{0,1}_\Omega}\cdot 2^{3k}2^{-k/2}|x|(2^k|x|)^{-2/p'}.
\end{split}
\end{equation*}
The desired bounds \eqref{winr40} follow in this case as well, by taking $p$ sufficiently large such that $1-1/p'\leq \delta/4$.

(ii) We prove now the bounds \eqref{winr41}. Notice that we may assume $2^j\leq t(1+2^k t)^{-\delta}$, since otherwise the bounds follow from \eqref{winr40}. We write
\begin{equation}\label{winr44}
(|x|+ix^aR_a)T_m(e^{-it\Lambda_{wa}}f_{j,k})(x)=C\int_{\mathbb{R}^3}(|x|-x\cdot\xi/|\xi|)m(\xi)e^{ix\cdot\xi}e^{-it|\xi|}\widehat{f_{j,k}}(\xi)\,d\xi.
\end{equation}
We may also assume $|x|\in[t/2,2t]$ (otherwise we have rapid decay using integration by parts). The conclusion follows from \eqref{winr41.2} and the observation that $\|f_{j,k}\|_{H^{0,2}_\Omega}\lesssim 2^{j+k}\|f^\ast_{j,k}\|_{H^{0,1}_\Omega}$.

It remains to prove the bounds \eqref{winr41.2}. We still use the formula \eqref{winr44}, and notice that the desired bounds follow easily if $2^k t\lesssim 1$, so we may assume that $2^kt\geq 2^{50}$. By rotation invariance,  we may assume $x=(x_1,0,0)$, $x_1\in[2^{-20}t,8t]$. Then we estimate the right-hand side of \eqref{winr44} by $Ct\sum_{b,c\in[0,2^{-20}(2^kt)^{1/2}]}|J_{b,c}|+R$, where
\begin{equation}\label{tric6}
\begin{split}
&J_{b,c}:=\int_{\mathbb{R}^3}\widehat{f_{j,k}}(\xi)\varphi_{[k-4,k+4]}(\xi)m(\xi)(1-\xi_1/|\xi|)\mathbf{1}_+(\xi_1)e^{ix_1\xi_1-it|\xi|}\psi_{b,c}(\xi)\,d\xi,\\
&\psi_{b,c}(\xi):=\varphi_b^{[0,\infty)}(\xi_2/2^{\lambda})\varphi_c^{[0,\infty)}(\xi_3/2^{\lambda}),\qquad 2^{\lambda}:=t^{-1/2}2^{k/2},
\end{split}
\end{equation}
and $R$ is an acceptable remainder that can be estimated using Lemma \ref{tech5}.

To prove \eqref{winr41.2} it suffices to show that for any $b,c\in[0,2^{-20}(2^kt)^{1/2}]$ we have
\begin{equation}\label{tric9}
|J_{b,c}|\lesssim t^{-2}2^{-k/2}(t 2^{k})^{\delta}\|f_{j,k}\|_{H^{0,2}_\Omega}.
\end{equation}
Notice that
\begin{equation}\label{tric9.02}
m(\xi)(1-\xi_1/|\xi|)\mathbf{1}_+(\xi_1)=\frac{m(\xi)}{1+\xi_1/|\xi|}\frac{\xi_2^2+\xi_3^2}{|\xi|^2}\mathbf{1}_+(\xi_1)=m'(\xi)\frac{\xi_2^2+\xi_3^2}{|\xi|^2}\mathbf{1}_+(\xi_1),
\end{equation}
in the support of the integrals defining $J_{b,c}$, for some suitable symbol $m'\in\mathcal{M}_0$.

We estimate first $|J_{0,0}|$. For any $p\in[2,\infty)$, using also \eqref{winr43} we have
\begin{equation}\label{tric9.2}
|J_{0,0}|\lesssim \|\widehat{f_{j,k}}(r\theta)\|_{L^2(r^2dr)L^p_\theta}(2^{\lambda-k})^{2/p'} 2^{3k/2}2^{2\lambda-2k} \lesssim_p \|f_{j,k}\|_{H^{0,1}_\Omega}\cdot t^{-2}2^{k/2}(t2^{k})^{1/p},
\end{equation} 
where the factor $(2^{\lambda-k})^{2/p'} 2^{3k/2}$ is due to the $L^2(r^2dr)L^{p'}_\theta$ norm of the support of the integral, and the factor $2^{2\lambda-2k}$ is due to the null factor $(\xi_2^2+\xi_3^2)/|\xi|^2$ in \eqref{tric9.02}. This is consistent with the desired bound \eqref{tric9}, by taking $p$ large enough.

To prove \eqref{tric9} when $(b,c)\neq (0,0)$ we may assume without loss of generality that $b\geq c$, so $b\geq\max(c,1)$. We integrate by parts in the integral in \eqref{tric6}, up to eight times, using the rotation vector-field $\Omega_{12}=\xi_1\partial_{\xi_2}-\xi_2\partial_{\xi_1}$. Since $\Omega_{12}\{x_1\xi_1-t|\xi|\}=-\xi_2x_1$, every integration by parts gains a factor of $t2^{\lambda+b}\approx t^{1/2}2^{k/2+b}$ and loses a factor $t^{1/2}2^{k/2}$. If $\Omega_{12}$ hits the function $\widehat{f_{j,k}}$ twice then we stop integrating by parts and bound the integral by estimating $\Omega_{12}^2\widehat{g_{j,k}}$ in $L^2$. As in \eqref{tric9.2} it follows that
\begin{equation*}
\begin{split}
|J_{b,c}|&\lesssim \big\{\|\widehat{f_{j,k}}(r\theta)\|_{L^2(r^2dr)L^p_\theta}+\|\widehat{\Omega_{12}f_{j,k}}(r\theta)\|_{L^2(r^2dr)L^p_\theta}\big\}(2^{\lambda-k})^{2/p'}2^{3k/2}2^{2\lambda-2k}2^{-b}\\
&+\|\Omega_{12}^2\widehat{f_{j,k}}\|_{L^2}(2^{\lambda+b}2^{k/2})(t2^{\lambda+b})^{-2}2^{2\lambda+2b-2k},
\end{split}
\end{equation*}
which gives the desired bound \eqref{tric9}. This completes the proof of the lemma.
\end{proof}

\subsection{Null and time-like geodesics} We consider now the future-directed causal geo\-de\-sics in our space-time $M$, and prove that they extend forever (in the affine parametrization) and become asymptotically parallel to the geodesics of the Minkowski space. More precisely:

\begin{theorem}\label{winr100} With $(\g,\psi)$ as in Theorem \ref{winr}, assume $p=(p^0,p^1,p^2,p^3)$ is a point in $M$ and $v=v^\al\partial_\al$ is a null or time-like vector at $p$, normalized with $v^0=1$. Then there is a unique affinely parametrized global geodesic curve $\gamma:[0,\infty)\to M$ with
\begin{equation}\label{winr101}
\gamma(0)=p=(p^0,p^1,p^2,p^3),\qquad\dot{\gamma}(0)=v=(v^0,v^1,v^2,v^3).
\end{equation}
Moreover, there is a vector $v_\infty=(v^0_\infty,v^1_\infty,v^2_\infty,v^3_\infty)$ such that, for any $s\in[0,\infty)$,
\begin{equation}\label{winr101.4}
|\dot{\gamma}(s)-v_\infty|\lesssim\varep_0(1+s)^{-1+6\delta'}\quad\text{ and }\quad m_{\al\be}v^\al_{\infty}v^\be_\infty=\g_{\al\be}(p)v^\al v^\be.
\end{equation}
The implicit constant in \eqref{winr101.4} is independent of $p$. As a consequence
\begin{equation}\label{winr101.6}
|\gamma(s)-v_\infty s-p|\lesssim\varep_0(1+s)^{6\delta'}\qquad\text{ for any }s\in[0,\infty).
\end{equation}

\end{theorem}

\begin{proof}  The proof uses only Theorem \ref{Precisedmetric}, Lemma \ref{extrader}, and the definitions. 

{\bf{Step 1.}} Assume that $T>0$ and $\gamma=(\gamma^0,\gamma^1,\gamma^2,\gamma^3):[0,T)\to M$ is a $C^4$ curve satisfying the geodesic equation
\begin{equation}\label{winr102}
\ddot{\gamma}^\mu+\dot{\gamma}^\al\dot{\gamma}^\be{\mathbf{\Gamma}^\mu}_{\al\be}=0,\qquad\mu\in\{0,1,2,3\},
\end{equation}
with initial data \eqref{winr101}. Let $V(s)=(V^0,V^1,V^2,V^3)(s):=\dot{\gamma}(s)$. The identity \eqref{winr102} implies the norm conservation identity
\begin{equation}\label{winr103}
{V}^\al(s){V}^\be(s) \mathbf{g}_{\al\be}(\gamma(s))=\mathrm{constant}\qquad\text{ for }s\in[0,T),
\end{equation}
as well as the general identity
\begin{equation}\label{winr104}
\frac{d}{ds}\big\{V^\be \g_{\al\be}\big\}=(1/2)V^\mu V^\nu\partial_\al h_{\mu\nu},\qquad\alpha\in\{0,1,2,3\}.
\end{equation}

Using first \eqref{winr103} it follows that $V^\al(s)V^\be(s) \g_{\al\be}\leq 0$ for all $s\in[0,T)$. In particular
\begin{equation}\label{winr106}
|V'(s)|\leq 1.1 |V^0(s)|\qquad\text{ for any }s\in[0,T),
\end{equation}
where $V':=(V^1,V^2,V^3)$. Since $V^0(0)=1$, we have $V^0(s)>0$ for all $s\in[0,T)$. 

{\bf{Step 2.}} The main idea is to prove that, for any $s\in[0,T)$ and $\al\in\{0,1,2,3\}$,
\begin{equation}\label{winr110}
V^0(s)\in[3/4,4/3]\qquad\text{ and }\qquad \int_0^s|\partial_\al h_{\mu\nu}(\gamma(u))V^\mu(u)V^\nu(u)|\,du\lesssim\varep_0.
\end{equation}
We use again a bootstrap argument. Assume that, for some $T'<T$, the weaker inequalities 
\begin{equation}\label{winr110.5}
V^0(s)\in[2/3,3/2]\qquad\text{ and }\qquad \int_0^s|\partial_\al h_{\mu\nu}(\gamma(u))V^\mu(u)V^\nu(u)|\,du\leq\varep_1,
\end{equation}
hold for any $s\in[0,T')$, where $\varep_1=\varep_0^{2/3}$ as before. It suffices to show that the stronger inequalities \eqref{winr110} hold for any $s\in[0,T')$, under the bootstrap assumption \eqref{winr110.5}.

It follows from \eqref{winr110.5} and \eqref{winr106} that, for any $s\in[0,T')$ and $\mu\in\{0,1,2,3\}$,
\begin{equation}\label{winr110.6}
|V^\mu(s)|\leq 2,\qquad \gamma^0(s)-p^0\in[2s/3,3s/2].
\end{equation}
We apply now \eqref{winr104} with $\alpha=0$. Let
\begin{equation}\label{winr108}
A(s):=-V^\beta(s) \g_{0\beta}(\gamma(s))=V^0(s)-V^\beta(s) h_{0\beta}(\gamma(s)).
\end{equation}
Using \eqref{winr104} we calculate
\begin{equation*}
\partial_sA=-(1/2)V^\mu V^\nu\partial_0 h_{\mu\nu}.
\end{equation*}
Using the bootstrap assumption \eqref{winr110.5} it follows that $|A(s)-A(0)|\leq \varep_1$ for any $s\in [0,T')$. Thus $|V^0(s)-1|\lesssim\varep_1$ (see \eqref{winr108}), and the bounds $V^0(s)\in[3/4,4/3]$ in \eqref{winr110} follow.

To prove the second bounds in \eqref{winr110} we would like to use \eqref{extrader1}, but for this we need to link the vectors $V(s)$ and $L(\gamma(s))$. We define the function $B:[0,T')\to\mathbb{R}$ by
\begin{equation}\label{winr120}
B(s):=1+\gamma^0(s)+(1+\gamma^0(s))^{2\delta'}-\big\{1+(\gamma^1(s))^2+(\gamma^2(s))^2+(\gamma^3(s))^2\big\}^{1/2}.
\end{equation}
One should think of $B$ as a slight modification (for the purpose of making it increasing along the geodesic curve $\gamma$) of the function $\gamma^0(s)-|\gamma'(s)|$, where $\gamma'(s):=(\gamma^1(s),\gamma^2(s),\gamma^3(s))$ and $|\gamma'(s)|:=[(\gamma^1(s))^2+(\gamma^2(s))^2+(\gamma^3(s))^2]^{1/2}$. We calculate
\begin{equation}\label{winr121}
(\partial_sB)(s)=V^0(s)+2\delta'V^0(s)(1+\gamma^0(s))^{2\delta'-1}-\frac{(\gamma^1V^1)(s)+(\gamma^2V^2)(s)+(\gamma^3V^3)(s)}{\big\{1+(\gamma^1(s))^2+(\gamma^2(s))^2+(\gamma^3(s))^2\big\}^{1/2}}.
\end{equation}
Notice that
\begin{equation}\label{winr125}
\sqrt{1+(B(s))^2}\lesssim \langle \gamma^0(s)-|\gamma'(s)|\rangle (1+\gamma^0(s))^{2\delta'}.
\end{equation}
Moreover, since the vector $V(s)$ is time-like or null and using \eqref{winr25}, we have 
\begin{equation*}
0\leq V^0(s)-|V'(s)|+C\varep_0V^0(s)(1+\gamma^0(s)+|\gamma'(s)|)^{-1+\delta'}
\end{equation*} 
for some constant $C\geq 1$, where $V'(s)=(V^1(s),V^2(s),V^3(s))$. Therefore
\begin{equation}\label{winr126}
(\partial_sB)(s)\geq \big|V^0(s)-|V'(s)|\big|+\Big\{|V'(s)|-\frac{(\gamma'\cdot V')(s)}{[1+|\gamma'(s)|^2]^{1/2}}\Big\}.
\end{equation}

We would like now to express the vector $V(s)$ in terms of the good vectors at the point $\gamma(s)$. More precisely, for any $s\in [0,T')$ we would like to decompose
\begin{equation}\label{winr127}
V(s)=|V'(s)|L(\gamma(s))+H'(\gamma(s))+E(\gamma(s)),
\end{equation}
where $H'=H^1\partial_1+H^2\partial_2+H^3\partial_3$ is a horizontal vector tangential to the sphere,
\begin{equation}\label{winr128}
|H^a(\gamma(s))|\lesssim \sqrt{(\partial_sB)(s)},\qquad H'(\gamma(s))\cdot \gamma'(s)=0,
\end{equation}
and $E=E^0\partial_0+E^1\partial_1+E^2\partial_2+E^3\partial_3$ is an error term,
\begin{equation}\label{winr129}
|E^\mu(\gamma(s))|\lesssim (\partial_sB)(s),\qquad \mu\in\{0,1,2,3\}.
\end{equation}
Indeed, one can simply take $H'=0$ if $|\gamma'(s)|\leq 1$, since in this case $|(\partial_sB)(s)|\gtrsim 1$ (see \eqref{winr126}). On the other hand, if $|\gamma'(s)|\geq 1$ then we use the basic decomposition
\begin{equation*}
V=V^0\partial_0+(V'\cdot\partial_r)\partial_r+H'=|V'|L+H'+\{(V^0-|V'|)\partial_0+[(V'\cdot\partial_r)-|V'|]\partial_r\},
\end{equation*}
where $H'=H^a\partial_a$, $H'\cdot\partial_r=0$. Since $\partial_r=|\gamma'|^{-1}\gamma'$, $|H'|^2=|V'|^2-(V'\cdot\partial_r)^2$ and $\partial_sB\geq |V^0-|V'||+(|V'|-V'\cdot\partial_r)$ (see \eqref{winr126}), the desired conclusions \eqref{winr127}--\eqref{winr129} follow.

We show now that, for any $s\in[0,T')$ and $\al\in\{0,1,2,3\}$,
\begin{equation}\label{winr130}
|\partial_\al h_{\mu\nu}(\gamma(s))V^\mu(s)V^\nu(s)|\lesssim \varep_0(1+|\gamma(s)|)^{-2+5\delta'}+\varep_0(1+|\gamma(s)|)^{-1+5\delta'}\frac{(\partial_sB)(s)}{(1+B(s)^2)^{1/2}}.
\end{equation}
Indeed, in view of \eqref{winr127} and \eqref{winr110} the left-hand side of \eqref{winr130} is bounded by
\begin{equation*}
\begin{split}
C\{&|[\partial_\al h_{\mu\nu}L^\mu L^\nu](\gamma(s))|+|[\partial_\al h_{\mu\nu}L^\mu (H')^\nu](\gamma(s))|\\
&+|[\partial_\al h_{\mu\nu}(H')^\mu (H')^\nu](\gamma(s))|+|[\partial_\al h_{\mu\nu}V^\mu E^\nu](\gamma(s))|+|[\partial_\al h_{\mu\nu}E^\mu E^\nu](\gamma(s))|.
\end{split}
\end{equation*}
The terms in the first line are bounded by $C\varep_0(1+|\gamma(s)|)^{-2+4\delta'}$, due to \eqref{ImprovedBoundsInside} and \eqref{extrader1}. The terms in the second line are bounded by $C\varep_0(1+|\gamma(s)|)^{-1+2\delta'}\langle \gamma^0(s)-|\gamma'(s)|\rangle^{-1}(\partial_sB)(s)$, due to \eqref{ImprovedBoundsInside} and \eqref{winr128}--\eqref{winr129}. The desired bounds \eqref{winr130} follow using also \eqref{winr125}.

Finally, we can complete the proof of the second estimate in \eqref{winr110}. Since $1+|\gamma(s)|\approx 1+s+p^0+|p'|\approx 1+s+|p|$ (due to \eqref{winr110.6}), it follows from \eqref{winr130} that
\begin{equation}\label{winr131}
\begin{split}
\int_{t_1}^{t_2}|\partial_\al &h_{\mu\nu}(\gamma(s))V^\mu(s)V^\nu(s)|\,ds\lesssim \int_{t_1}^{t_2}\varep_0(1+s+|p|)^{-2+5\delta'}\,ds\\
&+\int_{t_1}^{t_2}\varep_0(1+s+|p|)^{-1+5\delta'}\frac{d}{ds}\log[\sqrt{1+B(s)^2}+B(s)]\,ds\\
&\lesssim\varep_0(1+t_1+|p|)^{-1+6\delta'}
\end{split}
\end{equation}
for any $t_1\leq t_2\in[0,T')$, where we used integration by parts and the observation that $|B(s)|\lesssim 1+|\gamma(s)|\lesssim  1+s+|p|$. The conclusion \eqref{winr110} follows by setting $t_1=0$.

{\bf{Step 3.}} We can now prove the conclusions of the theorem. In view of \eqref{winr110}, \eqref{winr110.6} and the standard existence theory of solutions of ODEs, the geodesic curve $\gamma$ extends for all values of $s\in[0,\infty)$ as a smooth solution of the equation \eqref{winr102}. Thus one can take $T=\infty$, and the inequalities \eqref{winr110} are satisfied for all $s\in[0,\infty)$. 

We apply  \eqref{winr104} for all $\al\in\{0,1,2,3\}$ and integrate between times $s_1<s_2\in[0,\infty)$. Using \eqref{winr131} we have
\begin{equation*}
\big|V^\be(s_2) \g_{\al\be}(\gamma(s_2))-V^\be(s_1) \g_{\al\be}(\gamma(s_1))\big|\lesssim \varep_0(1+s_1+|p|)^{-1+6\delta'}.
\end{equation*}
Since $|h_{\al\be}(\gamma(s_2))|+|h_{\al\be}(\gamma(s_1))|\lesssim \varep_0(1+s_1+|p|)^{-1+\delta'}$ it follows that
\begin{equation}\label{winr135}
\big|V^\al(s_2) -V^\al(s_1) \big|\lesssim \varep_0(1+s_1+|p|)^{-1+6\delta'}
\end{equation}
for any $s_1<s_2\in[0,\infty)$\ and $\al\in\{0,1,2,3\}$. In particular, $v_\infty:=\lim_{s\to\infty}V(s)$ exists. Notice that \eqref{winr101.4} follows from \eqref{winr135} and \eqref{winr103}. The bounds in \eqref{winr101.6} then follow by integrating the bounds $|\dot{\gamma}(u)-v_\infty|\lesssim\varep_0(1+u)^{-1+6\delta'}$ from $0$ to $s$.
\end{proof}

\section{Weak peeling estimates and the ADM energy}\label{SecPeeling}

\subsection{Peeling estimates} In this section we prove weak peeling estimates for the Riemann tensor of our space-time. The Riemann tensor $\mathbf{R}$ satisfies the symmetry properties
\begin{equation}\label{Peel1}
\begin{split}
&\mathbf{R}_{\al\be\mu\nu}=-\mathbf{R}_{\be\al\mu\nu}=-\mathbf{R}_{\al\be\nu\mu}=\mathbf{R}_{\mu\nu\al\be},\\
&\mathbf{R}_{\al\be\mu\nu}+\mathbf{R}_{\be\mu\al\nu}+\mathbf{R}_{\mu\al\be\nu}=0.
\end{split}
\end{equation}
In our case, we also have the Einstein-field equations
\begin{equation}\label{Peel2}
\g^{\mu\nu}\mathbf{R}_{\al\mu\be\nu}=\mathbf{R}_{\al\be}=\D_\al\psi\D_\be\psi+(\psi^2/2)\g_{\al\be}.
\end{equation}

The rates decay of the components of the Riemann tensor are mainly determined by their {\it{signatures}}. To define this, we use (Minkowski) frames $(L,\underline{L},e_a)$, where $L,\underline{L}$ are as in \eqref{winr20} and $e_a\in\mathcal{T}_h:=\{r^{-1}\Omega_{12},r^{-1}\Omega_{23},r^{-1}\Omega_{31}\}$. We assign signature $+1$ to the vector-field $L$, $-1$ to the vector-field $\underline{L}$, and $0$ to the horizontal vector-fields in $\mathcal{T}_h$. With $e_1,e_2,e_3,e_4\in\mathcal{T}_h$ we define $\mathrm{Sig}(a)$ as the set of components of the Riemann tensor of total signature $a$, so
\begin{equation}\label{Peel3}
\begin{split}
\mathrm{Sig}(-2)&:=\{\mathbf{R}(\underline{L},e_1,\underline{L},e_2)\},\\
\mathrm{Sig}(2)&:=\{\mathbf{R}(L,e_1,L,e_2)\},\\
\mathrm{Sig}(-1)&:=\{\mathbf{R}(\underline{L},e_1,e_2,e_3),\mathbf{R}(\underline{L},L,\underline{L},e_1)\},\\
\mathrm{Sig}(1)&:=\{\mathbf{R}(L,e_1,e_2,e_3),\mathbf{R}(L,\underline{L},L,e_1)\},\\
\mathrm{Sig}(0)&:=\{\mathbf{R}(e_1,e_2,e_3,e_4),\mathbf{R}(L,\underline{L},e_1,e_2),\mathbf{R}(L,e_1,\underline{L},e_2),\mathbf{R}(L,\underline{L},L,\underline{L})\}.
\end{split}
\end{equation}
These components capture the entire curvature tensor, due to the symmetries \eqref{Peel1}. 

Notice that we define our decomposition in terms of the Minkowski null pair $(L,\underline{L})$ instead of more canonical null frames (or tetrads) adapted to the metric $\g$ (see, for example, \cite{ChKl}, \cite{KlNi}, \cite{KlNi2}). However, as we explain below, the weak peeling estimates are invariant under natural changes of the frame, and the rate of decay depends only on the signature of the component, except for the components $\mathbf{R}(L,e_1,\underline{L},e_2)$. More precisely:

\begin{theorem}\label{PeelingR}
Assume that  $(\g,\psi)$ as in Theorem \ref{winr}, $(x,t)\in M'=\{(x,t)\in M:\,t\geq 1\text{ and }|x|\geq 2^{-8}t\}$, and $\Psi_{(a)}\in\mathrm{Sig}(a')$ for $a\leq a'\in\{-2,-1,0,1,2\}$. Then
\begin{equation}\label{Peel4}
\begin{split}
|\Psi_{(-2)}(x,t)|&\lesssim \varep_0\langle r\rangle^{3\delta'-1}\langle t-r\rangle^{-2},\\
|\Psi_{(-1)}(x,t)|&\lesssim \varep_0\langle r\rangle^{5\delta'-2}\langle t-r\rangle^{-1},\\
|\Psi_{(2)}(x,t)|+|\Psi_{(1)}(x,t)|&\lesssim \varep_0\langle r\rangle^{7\delta'-3}.
\end{split}
\end{equation}
Moreover, if $\Psi_{(0)}^1\in\mathrm{Sig}_1(0):=\{\mathbf{R}(e_1,e_2,e_3,e_4),\mathbf{R}(L,\underline{L},e_1,e_2),\mathbf{R}(L,\underline{L},L,\underline{L})\}$ then
\begin{equation}\label{Peel5}
|\Psi_{(0)}^1(x,t)|\lesssim \varep_0\langle r\rangle^{7\delta'-3}.
\end{equation}
\end{theorem}

\begin{remark}
(i) In view of \eqref{Peel2}, we notice that the Ricci components decay at most cubically, $\mathbf{R}_{\alpha\beta}=O(\langle r\rangle^{-3+})$ in $M'$. As a result we do not expect uniform estimates of order better than cubic for any components of the Riemann curvature tensor, so the weak peeling estimates in Theorem \ref{PeelingR} are optimal in this sense, at least up to $r^{C\delta'}$ losses. In fact, the almost cubic decay is also formally consistent with the weak peeling estimates of Klainerman--Nicolo \cite[Theorem 1.2 (b)]{KlNi2} in the setting of our more general metrics.

(ii) We show in Proposition \ref{PeelingR2} below that our weak peeling estimates are invariant under suitable changes of the frame. Moreover, we also show that the natural $\langle r\rangle^{-3+}$ decay of the signature 0 components $\mathbf{R}(L,e_1,\underline{L},e_2)$, which is missing in the estimates \eqref{Peel4}--\eqref{Peel5}, can be restored if one works with a null vector-field $L$.
\end{remark}

\begin{proof} The estimates follow from the formulas \eqref{Peel1}--\eqref{Peel3}, the bounds on first and second order derivatives of $h_{\al\be}$ in Theorem \ref{Precisedmetric} and Lemmas \ref{extrader}--\ref{twoder0}, and the general identity
\begin{equation}\label{Peel6}
\begin{split}
&{\bf R}_{\alpha\beta\mu\nu}=-\partial_\alpha{\bf\Gamma}_{\nu\beta\mu}+\partial_\beta{\bf\Gamma}_{\nu\alpha\mu}+\g^{\rho\lambda}{\bf\Gamma}_{\rho\beta\mu}{\bf\Gamma}_{\lambda\alpha\nu}-\g^{\rho\lambda}{\bf\Gamma}_{\rho\alpha\mu}{\bf\Gamma}_{\lambda\beta\nu}\\
&=\frac{1}{2}\big[\partial_\al\partial_\nu h_{\be\mu}+\partial_\be\partial_\mu h_{\al\nu}-\partial_\be\partial_\nu h_{\al\mu}-\partial_\al\partial_\mu h_{\be\nu}\big]+\g^{\rho\lambda}\big[{\bf\Gamma}_{\rho\beta\mu}{\bf\Gamma}_{\lambda\alpha\nu}-{\bf\Gamma}_{\rho\alpha\mu}{\bf\Gamma}_{\lambda\beta\nu}\big].
\end{split}
\end{equation}

{\bf{Step 1.}} In view of \eqref{ImprovedBoundsInside}--\eqref{winr25.4} and \eqref{twoder1} we have the general bounds
\begin{equation}\label{Peel7}
\begin{split}
|\mathbf{\Gamma}_{\al\be\mu}(x,t)|&\lesssim\varep_0\langle t-r\rangle^{-1}\langle r\rangle^{-1+2\delta'},\\
|\partial_\al\partial_\be h_{\mu\nu}(x,t)|&\lesssim \varep_0\langle t-r\rangle^{-2}\langle r\rangle^{-1+3\delta'},\\
|\mathbf{R}_{\al\be}(x,t)|&\lesssim \varep_0\langle r\rangle^{-3+2\delta'},
\end{split}
\end{equation}
for any $(x,t)\in M'$ and $\al,\be,\mu,\nu\in\{0,1,2,3\}$. Also, if $V_1,V_2,V_3\in \mathcal{T}$ then, using also \eqref{extrader1}
\begin{equation}\label{Peel8}
\begin{split}
|V_1^\al V_2^\be &V_3^\mu\mathbf{\Gamma}_{\al\be\mu}(x,t)|+|\underline{L}^\al L^\be V_1^\mu\mathbf{\Gamma}_{\al\be\mu}(x,t)|\\
&+| L^\al V_1^\be\underline{L}^\mu\mathbf{\Gamma}_{\al\be\mu}(x,t)|+|V_1^\al\underline{L}^\be L^\mu\mathbf{\Gamma}_{\al\be\mu}(x,t)|\lesssim\varep_0\langle r\rangle^{-2+3\delta'}
\end{split}
\end{equation}
and
\begin{equation}\label{Peel8.5}
\begin{split}
|V_1^\al V_2^\be\partial_\al\partial_\be h_{\mu\nu}(x,t)|&\lesssim \varep_0\langle r\rangle^{-3+3\delta'},\\
|V_1^\al \underline{L}^\be\partial_\al\partial_\be h_{\mu\nu}(x,t)||&\lesssim \varep_0\langle t-r\rangle^{-1}\langle r\rangle^{-2+3\delta'}.
\end{split}
\end{equation}

We can now prove the bounds \eqref{Peel4}. The bounds on $\Psi_{(-2)}$ and $\Psi_{(2)}$ follow, since all the terms in the right-hand side of \eqref{Peel6} have suitable decay (using just \eqref{Peel7} for $\Psi_{(-2)})$ and \eqref{Peel8}--\eqref{Peel8.5} for $\Psi_{(2)}$). 

If $\Psi_{(-1)}$ is one of the curvature components in $\mathrm{Sig}(a')$, $a'\geq -1$, containing at most one vector-field $\underline{L}$ then the desired bounds follow using \eqref{Peel7} and \eqref{Peel8.5}. On the other hand, if $\Psi_{(-1)}=\mathbf{R}(\underline{L},L,\underline{L},V)$, $V\in\mathcal{T}$, then we use \eqref{Complete1} and the Einstein equations \eqref{Peel2}, thus
\begin{equation}\label{Peel10}
\begin{split}
-\mathbf{R}(\underline{L},V)&=m^{\al\be}\mathbf{R}(\underline{L},\partial_\al,\partial_\be,V)+g_{\geq 1}^{\al\be}\mathbf{R}(\underline{L},\partial_\al,\partial_\be,V)\\
&=-(1/2)\mathbf{R}(\underline{L},L,\underline{L},e_1)+\Pi^{\al\be}\mathbf{R}(\underline{L},\partial_\al,\partial_\be,e_1)+g_{\geq 1}^{\al\be}\mathbf{R}(\underline{L},\partial_\al,\partial_\be,e_1).
\end{split}
\end{equation}
Using \eqref{winr25} and \eqref{Peel7} it follows that $\mathbf{R}(\underline{L},L,\underline{L},V)=O(\varep_0\langle t-r\rangle^{-1}\langle r\rangle^{-2+5\delta'})$, as desired.

If $\Psi_{(1)}=\mathbf{R}(L,e_1,e_2,e_3)$ then the desired bounds in \eqref{Peel4} follow from \eqref{Peel8}--\eqref{Peel8.5}. On the other hand, if $\Psi_{1}=\mathbf{R}(L,\underline{L},L,e_1)$ then we use again the Einstein equations as in \eqref{Peel10}, and estimate
\begin{equation*}
|\mathbf{R}(L,\underline{L},L,e_1)|\lesssim |\mathbf{R}(L,e_1)|+|g_{\geq 1}^{\al\be}\mathbf{R}(L,\partial_\al,\partial_\be,e_1)|+\big|\Pi^{\al\be}\mathbf{R}(L,\partial_\al,\partial_\be,e_1)\big|.
\end{equation*}
The desired bounds follow using \eqref{Peel7}, \eqref{Peel6}, \eqref{winr25}, and the definition of $\Pi$ in \eqref{winr20}.

{\bf{Step 2.}} We bound now the components of signature 0.  Clearly, using just \eqref{Peel7}--\eqref{Peel8.5},
\begin{equation}\label{Peel12}
|\mathbf{R}(e_1,e_2,e_3,e_4)(x,t)|\lesssim\varep_0\langle r\rangle^{-3+5\delta'}.
\end{equation}

We prove now that for any $(x,t)\in M'$
\begin{equation}\label{Peel13}
|\mathbf{R}(L,\underline{L},e_1,e_2)(x,t)|\lesssim\varep_0\langle r\rangle^{-3+5\delta'}.
\end{equation}
The quadratic terms involving the connection coefficients can be bounded easily,
\begin{equation}\label{Peel14}
\begin{split}
|L^\al&\underline{L}^\be e_1^\mu e_2^\nu\g^{\rho\lambda}{\bf\Gamma}_{\rho\beta\mu}{\bf\Gamma}_{\lambda\alpha\nu}|\lesssim |L^\al\underline{L}^\be e_1^\mu e_2^\nu g_{\geq 1}^{\rho\lambda}{\bf\Gamma}_{\rho\beta\mu}{\bf\Gamma}_{\lambda\alpha\nu}|+|L^\al\underline{L}^\be e_1^\mu e_2^\nu m^{\rho\lambda}{\bf\Gamma}_{\rho\beta\mu}{\bf\Gamma}_{\lambda\alpha\nu}|\\
&\lesssim \varep_0\langle r \rangle^{-3+5\delta'}+|L^\al\underline{L}^\be e_1^\mu e_2^\nu \Pi^{\rho\lambda}{\bf\Gamma}_{\rho\beta\mu}{\bf\Gamma}_{\lambda\alpha\nu}|+|L^\al\underline{L}^\be e_1^\mu e_2^\nu (L^\rho\underline{L}^\lambda+\underline{L}^\rho L^\lambda){\bf\Gamma}_{\rho\beta\mu}{\bf\Gamma}_{\lambda\alpha\nu}|\\
&\lesssim \varep_0\langle r \rangle^{-3+5\delta'},
\end{split}
\end{equation}
for any $e_1,e_2\in\mathcal{T}_h$, using \eqref{winr25}, \eqref{Peel7}--\eqref{Peel8}, and \eqref{Complete1}. Moreover, we also have
\begin{equation}\label{Peel15}
|\underline{L}^\al e_1^\be L^\nu  e_2^\mu\partial_\al\partial_\be h_{\mu\nu}(x,t)|\lesssim \varep_0\langle r \rangle^{-3+5\delta'}.
\end{equation}
Indeed, to see this we start with the harmonic gauge condition \eqref{winr0.2} and write it in the form
\begin{equation*}
m^{\al\be}\partial_\nu\partial_\al h_{\be\mu}-(1/2)m^{\al\be}\partial_\nu\partial_\mu h_{\al\be}=-\partial_\nu(g_{\geq 1}^{\al\be}\partial_\al h_{\be\mu})+(1/2)\partial_\nu(g_{\geq 1}^{\al\be}\partial_\mu h_{\al\be}).
\end{equation*}
Using now \eqref{Complete1} and multiplying by $2e_1^\nu e_2^\mu$ we have
\begin{equation*}
\begin{split}
-(L^\al\underline{L}^\be+\underline{L}^\al L^\be)&e_1^\nu e_2^\mu\partial_\nu\partial_\al h_{\be\mu}+2\Pi^{\al\be}e_1^\nu e_2^\mu\partial_\nu\partial_\al h_{\be\mu}\\
&=m^{\al\be}e_1^\nu e_2^\mu\partial_\nu\partial_\mu h_{\al\be}-2e_1^\nu e_2^\mu\partial_\nu(g_{\geq 1}^{\al\be}\partial_\al h_{\be\mu})+e_1^\nu e_2^\mu\partial_\nu(g_{\geq 1}^{\al\be}\partial_\mu h_{\al\be}).
\end{split}
\end{equation*}
Thus
\begin{equation*}
\begin{split}
|\underline{L}^\al L^\be &e_1^\nu e_2^\mu\partial_\nu\partial_\al h_{\be\mu}|\lesssim |L^\al\underline{L}^\be e_1^\nu e_2^\mu\partial_\nu\partial_\al h_{\be\mu}|+|\Pi^{\al\be}e_1^\nu e_2^\mu\partial_\nu\partial_\al h_{\be\mu}|+|m^{\al\be}e_1^\nu e_2^\mu\partial_\nu\partial_\mu h_{\al\be}|\\
&+|e_1^\nu e_2^\mu\partial_\nu g_{\geq 1}^{\al\be}\partial_\al h_{\be\mu}|+|e_1^\nu e_2^\mu g_{\geq 1}^{\al\be}\partial_\nu\partial_\al h_{\be\mu}|+|e_1^\nu e_2^\mu\partial_\nu g_{\geq 1}^{\al\be}\partial_\mu h_{\al\be}|+|e_1^\nu e_2^\mu g_{\geq 1}^{\al\be}\partial_\nu\partial_\mu h_{\al\be}|,
\end{split}
\end{equation*}
and the bounds \eqref{Peel15} follow using \eqref{twoder1}, \eqref{ImprovedBoundsInside}, and \eqref{winr25}. The desired estimates \eqref{Peel13} now follow from \eqref{Peel14}--\eqref{Peel15}, the formula \eqref{Peel6}, and the bounds \eqref{twoder1}.

Finally, we prove now that for any $(x,t)\in M'$
\begin{equation}\label{Peel20}
|\mathbf{R}(L,\underline{L},L,\underline{L})(x,t)|\lesssim\varep_0\langle r\rangle^{-3+7\delta'}.
\end{equation}
Indeed, we start with the Einstein-field equations \eqref{Peel2}, thus
\begin{equation*}
\begin{split}
-\mathbf{R}(L,\underline{L})&=m^{\al\be}\mathbf{R}(L,\partial_\al,\partial_\be,\underline{L})+g_{\geq 1}^{\al\be}\mathbf{R}(L,\partial_\al,\partial_\be,\underline{L})\\
&=-(1/2)\mathbf{R}(L,\underline{L},L,\underline{L})+\Pi^{\al\be}\mathbf{R}(L,\partial_\al,\partial_\be,\underline{L})+g_{\geq 1}^{\al\be}\mathbf{R}(L,\partial_\al,\partial_\be,\underline{L}).
\end{split}
\end{equation*}
We use \eqref{Peel7} and the curvature bounds \eqref{Peel4} proved earlier to estimate the error terms in the identity above. For \eqref{Peel20} it suffices to prove that
\begin{equation}\label{Peel21}
\big|\Pi^{\mu\nu}\mathbf{R}(\partial_\mu,L,\partial_\nu,\underline{L})(x,t)\big|\lesssim\varep_0\langle r\rangle^{-3+7\delta'}.
\end{equation}

To prove \eqref{Peel21} we apply again the Einstein equations and \eqref{Complete1}, so
\begin{equation*}
\begin{split}
\Pi^{\mu\nu}\mathbf{R}_{\mu\nu}&=\Pi^{\mu\nu}m^{\al\be}\mathbf{R}(\partial_\mu,\partial_\al,\partial_\nu,\partial_\be)+\Pi^{\mu\nu}g_{\geq 1}^{\al\be}\mathbf{R}(\partial_\mu,\partial_\al,\partial_\nu,\partial_\be)\\
&=-\Pi^{\mu\nu}\mathbf{R}(\partial_\mu,L,\partial_\nu,\underline{L})+\Pi^{\mu\nu}\Pi^{\al\be}\mathbf{R}(\partial_\mu,\partial_\al,\partial_\nu,\partial_\be)+\Pi^{\mu\nu}g_{\geq 1}^{\al\be}\mathbf{R}(\partial_\mu,\partial_\al,\partial_\nu,\partial_\be).
\end{split}
\end{equation*}
In view of \eqref{Peel7}, \eqref{Peel12}, and \eqref{winr25}, for \eqref{Peel21} it suffices to prove that
\begin{equation}\label{Peel22}
\big|\Pi^{\mu\nu}\mathbf{R}(\partial_\mu,W_1,\partial_\nu,W_2)(x,t)\big|\lesssim\varep_0\langle r\rangle^{-2+6\delta'},
\end{equation}
for any $W_1,W_2\in\{L,\underline{L},r^{-1}\Omega_{12},r^{-1}\Omega_{23},r^{-1}\Omega_{31}\}$. 

The bounds \eqref{Peel22} follow from \eqref{Peel4} and \eqref{Peel12}, unless $W_1=W_2=\underline{L}$. In this case we use again the Einstein-field equations and \eqref{Complete1}, so
\begin{equation*}
\begin{split}
\mathbf{R}(\underline{L},\underline{L})&=m^{\al\be}\mathbf{R}(\underline{L},\partial_\al,\underline{L},\partial_\be)+g_{\geq 1}^{\al\be}\mathbf{R}(\underline{L},\partial_\al,\underline{L},\partial_\be)\\
&=\Pi^{\al\be}\mathbf{R}(\underline{L},\partial_\al,\underline{L},\partial_\be)+g_{\geq 1}^{\al\be}\mathbf{R}(\underline{L},\partial_\al,\underline{L},\partial_\be).
\end{split}
\end{equation*}
The desired bounds \eqref{Peel22} follow, which completes the proof of the theorem.
\end{proof}

We would like to show that our weak peeling estimates are invariant under natural changes of frames. Assume $L',\underline{L}',r^{-1}\Omega'_{12},r^{-1}\Omega'_{23},r^{-1}\Omega'_{31}$ are vector-fields in $M'$, which are small perturbations of the corresponding Minkowski vector-fields, i.e.
\begin{equation}\label{Peel29}
\begin{split}
|(L-L')(x,t)|&+|(\underline{L}-\underline{L}')(x,t)|+|(\Omega_{12}-\Omega'_{12})(x,t)|\\
&+|(\Omega_{23}-\Omega'_{23})(x,t)|+|(\Omega_{31}-\Omega'_{31})(x,t)|\leq \langle r\rangle^{-1+c_0},
\end{split}
\end{equation} 
for any $(x,t)\in M'$, where $c_0\in[2\delta',1/10]$. For $e'_1,e'_2,e'_3,e'_4\in\mathcal{T}'_h:=\{r^{-1}\Omega'_{12},r^{-1}\Omega'_{23},r^{-1}\Omega'_{31}\}$ we define the associated sets of curvature components as in \eqref{Peel3},
\begin{equation}\label{Peel30}
\begin{split}
\mathrm{Sig}'(-2)&:=\{\mathbf{R}(\underline{L}',e'_1,\underline{L}',e'_2)\},\\
\mathrm{Sig}'(2)&:=\{\mathbf{R}(L',e'_1,L',e'_2)\},\\
\mathrm{Sig}'(-1)&:=\{\mathbf{R}(\underline{L}',e'_1,e'_2,e'_3),\mathbf{R}(\underline{L}',L',\underline{L}',e'_1)\},\\
\mathrm{Sig}'(1)&:=\{\mathbf{R}(L',e'_1,e'_2,e'_3),\mathbf{R}(L',\underline{L}',L',e'_1)\},\\
\mathrm{Sig}'(0)&:=\mathrm{Sig}'_1(0)\cup\mathrm{Sig}'_2(0),\\
\mathrm{Sig}'_1(0)&:=\{\mathbf{R}(e'_1,e'_2,e'_3,e'_4),\mathbf{R}(L',\underline{L}',e'_1,e'_2),\mathbf{R}(L',\underline{L}',L',\underline{L}')\},\\
\mathrm{Sig}'_2(0)&:=\{\mathbf{R}(L',e'_1,\underline{L}',e'_2)\}.
\end{split}
\end{equation}

\begin{proposition}\label{PeelingR2}
Assume that $\Psi'_{(a)}\in\mathrm{Sig}'(a')$ for $a\leq a'\in\{-2,-1,0,1,2\}$, and $\Psi'^1_{(0)}\in\mathrm{Sig}'_1(0)$. Then $\Psi'_{(a)}, \Psi'^1_{(0)}$ satisfy similar bounds as in \eqref{Peel4}--\eqref{Peel5},
\begin{equation}\label{Peel31}
\begin{split}
|\Psi'_{(-2)}(x,t)|&\lesssim \varep_0\langle r\rangle^{-1+3\delta'}\langle t-r\rangle^{-2},\\
|\Psi'_{(-1)}(x,t)|&\lesssim \varep_0\langle r\rangle^{-2+3\delta'+c_0}\langle t-r\rangle^{-1},\\
|\Psi'_{(2)}(x,t)|+|\Psi'_{(1)}(x,t)|+|\Psi'^1_{(0)}(x,t)|&\lesssim \varep_0\langle r\rangle^{-3+3\delta'+2c_0}.
\end{split}
\end{equation}
Moreover, if the vector-field $L'$ satisfies the almost null bounds 
\begin{equation}\label{Peel31.3}
|\g(L',L')(x,t)|\lesssim \langle r\rangle^{-2+2c_0},\qquad \text{ for any }(x,t)\in M',
\end{equation}
(see Lemma \ref{AlmostOpticalLem} below for the construction of such vector-fields associated to almost optical functions) and $\Psi'^2_{(0)}\in\mathrm{Sig}'_2(0)$, then we have the additional bounds
\begin{equation}\label{Peel31.5}
|\Psi'^2_{(0)}(x,t)|\lesssim \varep_0\langle r\rangle^{-3+3\delta'+2c_0}.
\end{equation}

\end{proposition}

\begin{proof} {\bf{Step 1.}} We prove first the bounds \eqref{Peel31}. Assume that $W_1,\ldots,W_4$ are vector-fields in $\{L,\underline{L},r^{-1}\Omega_{ab}\}$ and let $W'_1,\ldots,W'_4$ denote their corresponding perturbations. In view of Theorem \ref{PeelingR} all the components of the curvature tensor in the Minkowski null frame are bounded by $C\varep_0\langle r\rangle^{-1+3\delta'}\langle t-r\rangle^{-2}$. Therefore we can estimate
\begin{equation}\label{Peel33}
\begin{split}
\big|&\mathbf{R}(W'_1,W'_2,W'_3,W'_4)-\mathbf{R}(W_1,W_2,W_3,W_4)\big|\lesssim \big|\mathbf{R}(W'_1-W_1,W_2,W_3,W_4)\big|\\
&+\big|\mathbf{R}(W_1,W'_2-W_2,W_3,W_4)\big|+\big|\mathbf{R}(W_1,W_2,W'_3-W_3,W_4)\big|\\
&+\big|\mathbf{R}(W_1,W_2,W_3,W'_4-W_4)\big|+\varep_0\langle r\rangle^{-2+2c_0}\langle r\rangle^{-1+3\delta'}\langle t-r\rangle^{-2}.
\end{split}
\end{equation}
The last remainder term in the right-hand side of \eqref{Peel33} is compatible with all the desired estimates \eqref{Peel31}. The other four curvature terms in the right-hand side are all bounded by $C\varep_0\langle r\rangle^{-1+c_0}\langle r\rangle^{-1+3\delta'}\langle t-r\rangle^{-2}$, which still suffices to prove the estimates on $|\Psi'_{(-2)}|$ and $|\Psi'_{(-1)}|$ in the first two lines of \eqref{Peel31}.

On the other hand, we can also use \eqref{Peel33} to estimate the terms $|\Psi'_{(2)}|$ and $|\Psi'_{(1)}|$ in the last line of \eqref{Peel31} by
\begin{equation}\label{Peel34}
C\varep_0\langle r\rangle^{-1+c_0}\langle r\rangle^{-2+5\delta'}\langle t-r\rangle^{-1}.
\end{equation}
This is because the four curvature terms in the right-hand of \eqref{Peel33} contain only components of signature $\geq -1$, since the change of one vector-field can reduce the signature by at most $2$. The remaining term $|\Psi'^1_{(0)}|$ is also bounded by the expression \eqref{Peel34}, because in the case of the components in $\mathrm{Sig}_1(0)$ the change of one vector-field can only reduce the total signature by $1$.{\footnote{We remark that this general argument fails for the components $\mathbf{R}(L,e_1,\underline{L},e_2)$ in $\mathrm{Sig_2(0)}$: even if such a component is $O(\langle r\rangle^{-3+})$ for a choice of frame, this bound is not invariant under a small change of the frame (satisfying \eqref{Peel29}) because the replacement of the vector $L$ in the first position would bring in errors of the form $\langle r\rangle^{-1+}\mathbf{R}(\underline{L},e_1,\underline{L},e_2)$, which are too big in the wave region.}} The desired bounds \eqref{Peel31} follow.

{\bf{Step 2.}} We prove now the bounds \eqref{Peel31.5}. Notice first that we may replace $\underline{L}'$, $e'_1$, and $e'_2$ with $\underline{L}$, $e_1$, and $e_2$ respectively, at the expense of components of signature $\geq -1$, thus bounded by the expression in \eqref{Peel34}. We may also assume that $L'=L+\rho\underline{L}$, where $|\rho(x,t)|\lesssim \langle r\rangle^{-1+c_0}$. Using \eqref{Peel1} we have
\begin{equation*}
\mathbf{R}(L',e_1,\underline{L},e_2)-\mathbf{R}(\underline{L},e_1,L',e_2)=\mathbf{R}(L',\underline{L},e_1,e_2).
\end{equation*}
Using the Einstein-field equations \eqref{Peel2} and \eqref{Complete1} we have
\begin{equation*}
\begin{split}
&\mathbf{R}(e_1,e_2)=m^{\al\be}\mathbf{R}(\partial_\al,e_1,\partial_\be,e_2)+g_{\geq 1}^{\al\be}\mathbf{R}(\partial_\al,e_1,\partial_\be,e_2)\\
&=-(1/2)[\mathbf{R}(L,e_1,\underline{L},e_2)+\mathbf{R}(\underline{L},e_1,L,e_2)]+\Pi^{\al\be}\mathbf{R}(\partial_\al,e_1,\partial_\be,e_2)+g_{\geq 1}^{\al\be}\mathbf{R}(\partial_\al,e_1,\partial_\be,e_2).
\end{split}
\end{equation*}
In view of \eqref{Peel31} and \eqref{Peel7}, we have
\begin{equation*}
|\mathbf{R}(L',\underline{L},e_1,e_2)(x,t)|+|\mathbf{R}(e_1,e_2)(x,t)|+|\Pi^{\al\be}\mathbf{R}(\partial_\al,e_1,\partial_\be,e_2)(x,t)|\lesssim \varep_0\langle r\rangle^{7\delta'-3}.
\end{equation*}
Combining the last three equations and recalling that $L=L'-\rho\underline{L}$ we have
\begin{equation*}
\begin{split}
\mathbf{R}(L',e_1,\underline{L},e_2)-\mathbf{R}(\underline{L},e_1,L',e_2)&=O(\varep_0\langle r\rangle^{7\delta'-3}),\\
[\mathbf{R}(L',e_1,\underline{L},e_2)+\mathbf{R}(\underline{L},e_1,L',e_2)]-2\rho\mathbf{R}(\underline{L},e_1,\underline{L},e_2)-2g_{\geq 1}^{\al\be}\mathbf{R}(\partial_\al,e_1,\partial_\be,e_2)&=O(\varep_0\langle r\rangle^{7\delta'-3}).
\end{split}
\end{equation*}
Therefore, to prove \eqref{Peel31.5} it suffices to show that
\begin{equation}\label{Peel40}
\rho\mathbf{R}(\underline{L},e_1,\underline{L},e_2)+g_{\geq 1}^{\al\be}\mathbf{R}(\partial_\al,e_1,\partial_\be,e_2)=O(\varep_0\langle r\rangle^{-3+3\delta'+2c_0}).
\end{equation}

The contributions of the curvature components in the term $\mathbf{R}(\partial_\al,e_1,\partial_\be,e_2)$ in \eqref{Peel40} are all bounded by the expression in \eqref{Peel34}, with the exception of the component $\mathbf{R}(\underline{L},e_1,\underline{L},e_2)$ which has signature $-2$.  Therefore, if we decompose $\partial_\al=A_\al\underline{L}+B_\al L+W$, where $W\cdot\partial_0=W\cdot\partial_r=0$, then the left-hand side \eqref{Peel40} is
\begin{equation}\label{Peel41}
\big[\rho+g_{\geq 1}^{\al\be}A_\al A_\be\big]\mathbf{R}(\underline{L},e_1,\underline{L},e_2)+O(\varep_0\langle r\rangle^{-3+3\delta'+2c_0}).
\end{equation}
Moreover, $A_\al=-(1/2)\partial_\al u^0$, where $u^0(x,t):=|x|-t$. Therefore
\begin{equation*}
\rho+g_{\geq 1}^{\al\be}A_\al A_\be=\rho+(1/4)g_{\geq 1}^{\al\be}\partial_\al u^0\partial_\be u^0=\rho-(1/4)L^\al L^\be h_{\al\be}+O(\langle r\rangle^{-2+2\delta'}),
\end{equation*}
where the last identity follow from the explicit formulas \eqref{zaq22}. Since $L=L'-\rho\underline{L}$, $\rho=O(\langle r\rangle^{-1+c_0})$, and $\g(L',\underline{L})=-2+O(\langle r\rangle^{-1+c_0})$ we have
\begin{equation*}
4\rho=\g(L,L)-\g(L',L')+O(\langle r\rangle^{-2+2c_0}).
\end{equation*}
The last two identities and the assumption \eqref{Peel31.3} show that $\rho+g_{\geq 1}^{\al\be}A_\al A_\be=O(\langle r\rangle^{-2+2c_0})$. The desired conclusion \eqref{Peel40} follows using also \eqref{Peel41}.
\end{proof}

\subsection{The ADM energy} The ADM energy measures the total deviation of our space-time from the Minkowski solution. In our asymptotically flat case it is calculated by integrating on large spheres inside the surfaces $\Sigma_t=\{(x,t)\in M:\,x\in\mathbb{R}^3\}$, according to the formula
\begin{equation}\label{ADM1}
E_{ADM}(t):=\frac{1}{16\pi}\lim_{R\to\infty}\int_{S_{R,t}}(\partial_j\g_{nj}-\partial_n\g_{jj})\frac{x^n}{|x|}\,dx,
\end{equation}
where the integration is over large (Euclidean) spheres $S_{R,t}\subseteq\Sigma_t$ of radius $R$. Using Stokes theorem and the definitions \eqref{zaq21.1}, we can rewrite 
\begin{equation}\label{ADM2}
E_{ADM}(t)=\frac{1}{16\pi}\lim_{R\to\infty}\int_{|x|\leq R}-2\Delta\tau(x,t)\,dx.
\end{equation}
We analyze first the density function $-2\Delta\tau$.

\begin{lemma}\label{ADM20}
We can decompose
\begin{equation}\label{ADM21}
\begin{split}
-2\Delta\tau&=-\delta_{jk}P^2_{jk}+\{(\partial_0\psi)^2+\psi^2+\partial_j\psi\partial_j\psi\}+O^1+\partial_jO^2_j,
\end{split}
\end{equation}
where $P^2_{jk}$ are defined as in \eqref{zaq21}, 
\begin{equation}\label{ADM22}
\|O^1(t)\|_{L^1(\mathbb{R}^3)}\lesssim\varep_0^2\langle t\rangle^{-\kappa}\qquad\text{ for any }t\geq 0,
\end{equation}
and
\begin{equation}\label{ADM21.5}
\begin{split}
O_j^{2}:&=-h_{0j}\partial_0h_{00}+h_{00}\partial_0 h_{0j}-h_{0n}\partial_0 h_{nj}+h_{nj}\partial_0h_{0n}-h_{kn}\partial_jh_{kn}+h_{0n}\partial_j h_{0n}\\
&-h_{0n}\partial_n h_{0j}+h_{kn}\partial_n h_{kj}-h_{00}\partial_j(\tau+\underline{F})+2h_{0j}\partial_0(\tau+\underline{F})-h_{nj}\partial_n(\tau+\underline{F}).
\end{split}
\end{equation}
\end{lemma}

\begin{proof} Recall the identity \eqref{zaq11}
\begin{equation}\label{ADM5}
-2\Delta\tau=\partial_\al E_\al^{\geq 2}+\underline{F}(\mathcal{N}^h)+\tau(\mathcal{N}^h)=\partial_\al E_\al^{\geq 2}+\frac{1}{2}\mathcal{N}^h_{00}+\frac{1}{2}\delta_{jk}\mathcal{N}^h_{jk}.
\end{equation}

{\bf{Step 1.}} We use the formulas \eqref{sac1}--\eqref{zaq21} and \eqref{zaq5.1} and identify first the $L^1$ errors. Indeed, all the cubic and higher order terms are bounded by $C\varep_0^2\langle t-|x|\rangle^{-1}\langle t+|x|\rangle^{-3+\kappa}$ or by $C\varep_0^2\langle x\rangle^{-1}\langle t+|x|\rangle^{-3+\kappa}$, due to \eqref{winr25}--\eqref{winr25.4}, thus acceptable $L^1$ errors. The semilinear quadratic null forms in $Q^2_{\al\be}$ (see \eqref{zaq20}) are also acceptable $L^1$ errors because
\begin{equation}\label{ADM23}
\big|(m^{\al\be}\partial_\al h_1\partial_\be h_2)(x,t)\big|+\big|(\partial_\mu h_1\partial_\nu h_2-\partial_\nu h_1\partial_\mu h_2)(x,t)\big|\lesssim\varep_0^2\langle t-|x|\rangle^{-1}\langle t+|x|\rangle^{-3+\kappa},
\end{equation}
for any $\mu,\nu\in\{0,1,2,3\}$ and $h_1,h_2\in\{h_{\al\be}\}$. This follows from \eqref{ImprovedBoundsInside} and the observation
\begin{equation}\label{ADM23.5}
\partial_0=-\partial_r+L\qquad\text{ and }\qquad\partial_j=(x_j/r)\partial_r+\text{ sum of good vector-fields in }\mathcal{T}.
\end{equation}

The Klein-Gordon contributions coming from \eqref{sac1.5} are the second term in the right-hand side of \eqref{ADM21}. It remains to analyze the quadratic semilinear and quasilinear terms involving the metric components, which are
\begin{equation}\label{ADM24}
\frac{1}{2}(\mathcal{Q}^2_{00}+\mathcal{Q}^2_{jj})-\frac{1}{2}(P^2_{00}+P^2_{jj})+\partial_0E_0^2+\partial_jE_j^2,
\end{equation}
where, using \eqref{zaq5.1} and \eqref{zaq22},
\begin{equation}\label{ADM25}
\begin{split}
E^{2}_\mu&:=-g_1^{\al\be}\partial_\al h_{\be\mu}+\frac{1}{2}g_1^{\al\be}\partial_\mu h_{\al\be}=h_{00}\partial_0 h_{0\mu}\\
&-h_{0k}(\partial_0 h_{k\mu}+\partial_k h_{0\mu})+h_{kn}\partial_k h_{n\mu}-\frac{1}{2}h_{00}\partial_\mu h_{00}+h_{0k}\partial_\mu h_{0k}-\frac{1}{2}h_{kn}\partial_\mu h_{kn}.
\end{split}
\end{equation}
Notice that all the terms in \eqref{ADM24} are either of the form $\partial h\cdot \partial h$ or $h\cdot\partial^2h$. To extract acceptable $L^1$ errors, we notice that, for any $\mu,\nu\in\{0,1,2,3\}$, $h\in\{h_{\al\be}\}$, and $V\in\mathcal{T}$,
\begin{equation}\label{ADM25.5}
\|h(t)\|_{L^p}\lesssim\varep_0\langle t\rangle^{-\kappa}\quad\text{ if }p\geq 3+4\kappa,
\end{equation}
\begin{equation}\label{ADM25.6}
\|\partial_\mu h(t)\|_{L^p}+\|\partial_\mu\partial_\nu h(t)\|_{L^p}\lesssim\varep_0\langle t\rangle^{-\kappa}\quad\text{ if }p\geq 2+4\kappa,
\end{equation}
\begin{equation}\label{ADM25.8}
\|(1-\chi_t)\cdot\partial_\mu h(t)\|_{L^p}+\|(1-\chi_t)\cdot\partial_\mu\partial_\nu h(t)\|_{L^p}+\|V^\mu\partial_\mu h(t)\|_{L^p}\lesssim\varep_0\langle t\rangle^{-\kappa}\quad\text{ if }p \geq 3/2+4\kappa,
\end{equation}
where $\chi_t$ is a smooth characteristic functions of the wave region, for example $\chi_t\equiv 0$ if $t\leq 8$ and $\chi_t(x)=\varphi_{\leq -4}((t-|x|)/t)$ if $t\geq 8$. These bounds follow from Theorem \ref{Precisedmetric}. Moreover
\begin{equation}\label{ADM25.9}
\begin{split}
\|(\partial_{\al}&h_1\cdot \partial_\be h_2)(t)\|_{L^p}+\|(h_1\cdot \partial_{\al}\partial_\be h_2)(t)\|_{L^p}\\
&+\|\mathcal{N}^h_{\al\be}(t)\|_{L^p}+\|E^{\geq 2}_{\al}(t)\|_{L^p}+\|\partial_\mu E^{\geq 2}_{\al}(t)\|_{L^p}\lesssim\varep_0\langle t\rangle^{-\kappa},
\end{split}
\end{equation}
for any $p\geq 1+2\kappa$. In particular, all the terms in \eqref{ADM24} are barely missing to being acceptable $L^1$ errors. All the semilinear terms that contain a good derivative are acceptable $L^1$ errors.

{\bf{Step 2.}} We analyze now the quadratic expressions in \eqref{ADM24}. Since $(\partial_0^2-\Delta)h_{\al\be}=\mathcal{N}^h_{\al\be}$, we can use \eqref{ADM25.9} to write
\begin{equation*}
\begin{split}
\partial_0E_0^2&\sim
\frac{1}{2}\partial_0h_{00}\partial_0h_{00}+\frac{1}{2}h_{00}\Delta h_{00}-\partial_0h_{0k}\partial_kh_{00}-h_{0k}\partial_0\partial_kh_{00}\\
&+\partial_0h_{kn}\partial_kh_{n0}+h_{kn}\partial_0\partial_kh_{n0}-\frac{1}{2}\partial_0h_{kn}\partial_0h_{kn}-\frac{1}{2}h_{kn}\Delta h_{kn}.
\end{split}
\end{equation*}
where in this proof $F\sim G$ means $\|F-G\|_{L^1}\lesssim\varep_0^2\langle t\rangle^{-\kappa}$. Using also \eqref{ADM23} we get
\begin{equation*}
\partial_0E_0^2\sim\frac{1}{2}\partial_j(h_{00}\partial_jh_{00})-\partial_k(h_{0k}\partial_0h_{00})+\partial_k(h_{nk}\partial_0h_{n0})-\frac{1}{2}\partial_j(h_{kn}\partial_jh_{kn})
\end{equation*}
Therefore, using again \eqref{ADM25}, we have $\partial_0E_0^2+\partial_jE_j^2\sim\partial_jO_j^{2,1}$, where
\begin{equation}\label{ADM28}
\begin{split}
O_j^{2,1}&:=-h_{0j}\partial_0h_{00}+h_{nj}\partial_0h_{n0}-h_{kn}\partial_jh_{kn}+h_{00}\partial_0 h_{0j}\\
&-h_{0n}\partial_0 h_{nj}-h_{0n}\partial_n h_{0j}+h_{kn}\partial_k h_{nj}+h_{0n}\partial_j h_{0n}.
\end{split}
\end{equation}

We examine now the terms $\mathcal{Q}^2_{00}$ and $\mathcal{Q}^2_{jj}$. Using the identities \eqref{sac1.2} we write
\begin{equation}\label{ADM28.5}
\frac{1}{2}\big(\mathcal{Q}^2_{00}+\mathcal{Q}^2_{nn}\big)=\partial_jO_j^{2,2}+\big\{\partial_j h_{00}\partial_j(\tau+\underline{F})-2\partial_jh_{0j}\partial_0(\tau+\underline{F})+\partial_n h_{nj}\partial_j(\tau+\underline{F})\big\},
\end{equation}
where $(1/2)(h_{00}+\delta_{jk}h_{jk})=\tau+\underline{F}$ (see \eqref{zaq21} and \eqref{zaq21.1}) and
\begin{equation}\label{ADM29}
O_j^{2,2}:=-h_{00}\partial_j(\tau+\underline{F})+2h_{0j}\partial_0(\tau+\underline{F})-h_{nj}\partial_n(\tau+\underline{F}).
\end{equation}
The desired formula \eqref{ADM21.5} follows from \eqref{ADM28} and \eqref{ADM29}.

{\bf{Step 3.}} We identify now the contribution of the semilinear terms. We show first that the semilinear terms in the bracket in the right-hand side of \eqref{ADM28.5} are acceptable $L^1$ errors. Indeed, using \eqref{ADM23.5} and \eqref{ADM25.6}--\eqref{ADM25.8}, we may replace $\partial_j$ with $-(x_j/r)\partial_0$, at the expense of acceptable errors. The semilinear expression in \eqref{ADM28.5} is 
\begin{equation*}
\begin{split}
&\sim\big\{\partial_0 h_{00}\partial_0(\tau+\underline{F})+2(x_j/r)\partial_0h_{0j}\partial_0(\tau+\underline{F})+(x_nx_j/r^2)\partial_0 h_{nj}\partial_0(\tau+\underline{F})\big\}\\
&\sim\partial_0(\tau+\underline{F})\cdot (L^\al L^\be\partial_0h_{\al\be}).
\end{split}
\end{equation*}
This is an acceptable $L^1$ error due to Lemma \ref{extrader} and \eqref{ImprovedBoundsInside}.

Finally, we examine \eqref{zaq21} and notice that $P_{00}\sim P_{jj}$, due to \eqref{ADM23.5} and \eqref{ADM25.6}--\eqref{ADM25.8}. The contribution of $-(1/2)(P^2_{00}+P^2_{jj})$ is the first term in the right-hand side of \eqref{ADM21}.
\end{proof}

We prove now that the ADM energy is well defined, conserved in time, non-negative, and can be linked to the scattering data of our space-time. More precisely:

\begin{proposition}\label{ADM4}
We have $\Delta\tau(t)\in L^1(\mathbb{R}^3)$ and
\begin{equation}\label{ADM4.5}
E_{ADM}(t)=\frac{1}{16\pi}\int_{\mathbb{R}^3}-2\Delta\tau(x,t)\,dx=E_{ADM}
\end{equation}
does not depend on $t\in[0,\infty)$. Moreover, for any $t\geq 0$,
\begin{equation}\label{ADM4.6}
E_{ADM}=\frac{1}{16\pi}\int_{\mathbb{R}^3}\big\{|U^\psi(t)|^2+(1/4)\sum_{m,n\in\{1,2,3\}}|U^{\vartheta_{mn}}(t)|^2\big\}\,dx+O(\varep_0^2\langle t\rangle^{-\kappa}).
\end{equation}
In particular, recalling the scattering profiles $V^\psi_\infty$ and $V^G_\infty$ from \eqref{Alop10.1}, we have
\begin{equation}\label{ADM4.7}
E_{ADM}=\frac{1}{16\pi}\|V^\psi_\infty\|^2+\frac{1}{64\pi}\sum_{m,n\in\{1,2,3\}}\|V^{\vartheta_{mn}}_\infty\|_{L^2}^2.
\end{equation}
\end{proposition}

\begin{proof} {\bf{Step 1.}} It follows from \eqref{ADM21}, \eqref{winr25}, \eqref{ImprovedBoundsInside}, and \eqref{twoder1} that
\begin{equation}\label{ADM8}
\|\partial_\mu E_\nu^{\geq 2}(t)\|_{(L^1\cap L^2)(\mathbb{R}^3)}+\|\mathcal{N}^h_{\mu\nu}\|_{(L^1\cap L^2)(\mathbb{R}^3)}\lesssim\varep_0\langle t\rangle^4,\qquad t\in[0,\infty),
\end{equation}
for any $\mu,\nu\in\{0,1,2,3\}$. In particular $\Delta\tau\in L^1(\mathbb{R}^3)$ (due to \eqref{zaq11}),  and the identity \eqref{ADM4.5} holds. To prove that the energy is constant we estimate, for any $t_1\leq t_2\in[0,\infty)$ and $R$ large,
\begin{equation*}
\begin{split}
\Big|&\int_{\mathbb{R}^3}\Delta\tau(x,t_2)\varphi_{\leq 0}(|x|/R)\,dx-\int_{\mathbb{R}^3}\Delta\tau(x,t_1)\varphi_{\leq 0}(|x|/R)\,dx\Big|\\
&=\Big|\int_{t_1}^{t_2}\int_{\mathbb{R}^3}\partial_0\partial_j\tau(x,s)\cdot\partial_j\varphi_{\leq 0}(|x|/R)\,dx\Big|\lesssim |t_2-t_1|R^{-1/4}\sup_{s\in[t_1,t_2]}\sum_{j\in\{1,2,3\}}\|\partial_0\partial_j\tau(s)\|_{L^{4/3}}.
\end{split}
\end{equation*}
Using now \eqref{zaq11} and \eqref{ADM8} it follows that
\begin{equation*}
\lim_{R\to\infty}\Big|\int_{\mathbb{R}^3}\Delta\tau(x,t_2)\varphi_{\leq 0}(|x|/R)\,dx-\int_{\mathbb{R}^3}\Delta\tau(x,t_1)\varphi_{\leq 0}(|x|/R)\,dx\Big|=0,
\end{equation*}
for any $t_1\leq t_2\in[0,\infty)$, so the function $M_{ADM}$ is constant in time.

{\bf{Step 2.}} We prove now the identity \eqref{ADM4.6}. We start from \eqref{ADM21} and notice that the contribution of the Klein-Gordon field $\psi$ is given by the integral of $|U^\psi(t)|^2$, as claimed. The divergence term $\partial_jO^2_j$ does not contribute, because $\|O_j^2(t)\|_{L^{4/3}}\lesssim \varep_0^2$. 

Finally, to calculate the contribution of $\delta_{jk}P_{jk}$ we would like to use Lemma \ref{Pstructure}. To apply it, we use the bounds \eqref{ADM40} proved below. In particular, using also \eqref{ADM23.5} we have
\begin{equation}\label{ADM43}
\big\|(m^{\al\be}\partial_\al Rh\cdot\partial_\be R'h')(t)\big\|_{L^1}+\big\|(\partial_\mu Rh\cdot\partial_\nu R'h'-\partial_\nu Rh\cdot\partial_\mu R'h')(t)\big\|_{L^1}\lesssim\varep_0^2\langle t\rangle^{-\kappa},
\end{equation}
for any $h,h'\in\{h_{\al\be}\}$, $\mu,\nu\in\{0,1,2,3\}$, and any compounded Riesz transforms $R,R'$.

We can now use the calculations in Lemma \ref{Pstructure} with $\LL=\LL_1=\LL_2=Id$. The cubic and higher order terms are all acceptable $L^1$ errors, due to \eqref{ADM25.9}. Also, most of the quadratic terms in \eqref{Pstruc3.5} can be estimated using \eqref{ADM23}, for example
\begin{equation*}
R_jR_k\partial_nG_1\cdot R_jR_k\partial_nG_2=\partial_jR_kR_nG_1\cdot \partial_nR_jR_kG_2\sim \partial_nR_kR_nG_1\cdot \partial_jR_jR_kG_2=\partial_kG_1\cdot \partial_kG_2,
\end{equation*}
for any $G_1,G_2\in\{F,\underline{F},\omega_a,\vartheta_{ab}\}$, where $f\sim g$ means $\|f-g\|_{L^1}\lesssim\varep_0^2\langle t\rangle^{-\kappa}$ as before. The terms containing $\tau$ are all acceptable $L^1$ errors, due to \eqref{zaq11} and \eqref{ADM40} below. The only remaining terms are $(1/2)\partial_j\va_{mn}\cdot\partial_j\va_{mn}\sim(1/2)\partial_0\va_{mn}\cdot\partial_0\va_{mn}\sim(1/4)(\partial_j\va_{mn}\cdot\partial_j\va_{mn}+\partial_0\va_{mn}\cdot\partial_0\va_{mn})$, coming from $-\delta_{jk}A^{\va\va}_{jk}$, which lead to the identity \eqref{ADM4.6}.
\end{proof}

We collect now some bounds on the Riesz transforms of the metric components, which are used in Proposition \ref{ADM4} and Theorem \ref{BondThm} below.

\begin{lemma}\label{DecayX}
Assume that $R$ is a compounded Riesz transform $R=R_1^{a_1}R_2^{a_2}R_3^{a_3}$, $a_1+a_2+a_3\leq 6$, $h\in\{h_{\al\be}\}$, $\mathcal{N}^h\in\{\mathcal{N}^h_{\al\be}\}$, and $V\in\mathcal{T}$ is a good vector-field. Then
\begin{equation}\label{ADM40}
\begin{split}
\|\partial_\mu Rh(t)\|_{L^p}&\lesssim\varep_0\langle t\rangle^{-\kappa}\quad\text{ if }p\geq 2+4\kappa,\\
\|V^\mu\partial_\mu Rh(t)\|_{L^q}&\lesssim\varep_0\langle t\rangle^{-\kappa}\quad\text{ if }q \geq 3/2+4\kappa,\\
\big\||\nabla|^{-1}R\mathcal{N}^h(t)\big\|_{L^q}+\big\||\nabla|^{-1}R\partial_\mu E^{\geq 2}_{\al}(t)\big\|_{L^q}&\lesssim\varep_0\langle t\rangle^{-\kappa}\quad\text{ if }q \geq 3/2+4\kappa,
\end{split}
\end{equation}
for any $\mu,\alpha\in\{0,1,2,3\}$. In addition, if $t\leq 1$ or if $t\geq 1$ and $|x|\in[2^{-10}t,2^{10}t]$ then
\begin{equation}\label{ADM40.5}
\begin{split}
|Rh(x,t)|+|\partial_\mu Rh(x,t)|&\lesssim\varep_0\langle t\rangle^{-1+\kappa},\\
|V^\mu\partial_\mu Rh(x,t)|&\lesssim\varep_0\langle t\rangle^{-4/3+\kappa},\\
\big||\nabla|^{-1}R\mathcal{N}^h(x,t)\big|+\big||\nabla|^{-1}RE^{\geq 2}_{\al}(x,t)\big|+\big||\nabla|^{-1}R\partial_\mu E^{\geq 2}_{\al}(x,t)\big|&\lesssim\varep_0\langle t\rangle^{-4/3+\kappa}.
\end{split}
\end{equation}

\end{lemma}

\begin{proof} The bounds in the first line of \eqref{ADM40} follow directly from \eqref{ADM25.6}. The bounds in the third line follow from \eqref{ADM25.9} and the Hardy--Littlewood--Sobolev inequality. To estimate $\|V^\mu\partial_\mu Rh(t)\|_{L^q}$ we notice that the contribution of $(1-\chi_t)\partial_\mu h(t)$ is bounded easily, due to the estimate on the first term in the left-hand side of \eqref{ADM25.8}. Moreover
\begin{equation}\label{ADM41}
\|V^\mu R(\chi_t\partial_\mu h(t))\|_{L^q}\lesssim\|V^\mu \chi_t\partial_\mu h(t)\|_{L^q}+\big\|[V^\mu,R](\chi_t\partial_\mu h(t))\big\|_{L^q}.
\end{equation}
The first term is bounded by $C\varep_0\langle t\rangle^{-\kappa}$, due to \eqref{ADM25.8}. The second term is a Calder\'{o}n commutator that can be estimated easily, with $g(x):=\chi_t(x)\partial_\mu h(x,t)$, 
\begin{equation}\label{ADM41.5}
\big |[V^\mu,R](g)(x)\big |\lesssim\int_{\mathbb{R}^3}|g(y)||x-y|^{-3}|V^\mu(x)-V^\mu(y)|\,dy\lesssim\int_{\mathbb{R}^3}|g(y)|t^{-1}|x-y|^{-2}\,dy.
\end{equation}
Therefore, recalling that $|g(x)|\lesssim \varep_0\chi_t(x)\langle t-|x|\rangle^{-1}t^{2\delta'-1}$ (see \eqref{ImprovedBoundsInside}), we have
\begin{equation*}
\big\|[V^\mu,R](\chi_t\partial_\mu h(t))\big\|_{L^q}\lesssim t^{-1}\|g\|_{L^{p_1}}\lesssim\varep_0t^{-1}t^{2\delta'-1}t^{2/p_1},
\end{equation*}
where $1/p_1=1/q+1/3$ (by fractional integration). The desired bounds \eqref{ADM40} follow using also \eqref{ADM41}. 

We notice now that the pointwise bounds \eqref{ADM40.5} follow easily from \eqref{winr1} if $t \lesssim 1$. On the other hand, if $t\geq 8$ then we fix $K_0$ the largest integer such that $2^{K_0}\leq \langle t\rangle^{8}$ and notice that the contribution of low or high frequencies $P_{\leq -K_0-1}h+P_{\geq K_0+1}h$ is suitably bounded due to \eqref{winr1}. For the medium frequencies we use the general estimate
\begin{equation}\label{ADM50}
\|RP_{[-K_0,K_0]}f\|_{L^\infty}\lesssim \log (2+t)\|P_{[-K_0,K_0]}f\|_{L^\infty}.
\end{equation}
The bounds in the first line of \eqref{ADM40.5} follow from \eqref{winr25} and \eqref{ImprovedBoundsInside}. 

To prove the bounds in the second line we estimate
\begin{equation*}
|V^\mu R(P_{[-K_0,K_0]}\partial_\mu h)(x,t)|\lesssim\|RP_{[-K_0,K_0]}(V^\mu \partial_\mu h)(t)\|_{L^\infty}+\big|[V^\mu,RP_{[-K_0,K_0]}](\partial_\mu h)(x,t)\big|.
\end{equation*}
The first term in the right-hand side is bounded by $\varep_0 t^{-2+4\delta'}$, due to \eqref{ImprovedBoundsInside} and \eqref{ADM50}. Recalling that $|x|\approx t$, the second term can be bounded as in \eqref{ADM41.5}
\begin{equation*}
\big|[V^\mu,RP_{[-K_0,K_0]}](\partial_\mu h)(x,t)\big|\lesssim\int_{\mathbb{R}^3}|\partial_\mu h(y)|t^{-1}|x-y|^{-2}\,dy\lesssim\varep_0\int_{\mathbb{R}^3}\frac{\langle t+|y|\rangle^{4\delta'-1}}{\langle t-|y|\rangle t|x-y|^{2}}\,dy,
\end{equation*}
where we used \eqref{ImprovedBoundsInside} for the second estimate. The last integral in the inequalities above is bounded by $Ct^{-4/3+5\delta'}$, which suffices to prove the estimates in the second line of \eqref{ADM40.5}.

The estimates in the third line follow from H\"{o}lder's inequality once we notice that the functions $\mathcal{N}^h$, $E^{\geq 2}_{\al}$, and $\partial_\mu E^{\geq 2}_{\al}$ are all bounded by $C\varep_0^2t^{-4/3+5\delta'}$ in $L^p$ for all $p\geq 3-\delta'$.
\end{proof}

\subsection{The linear momentum} With $\Sigma_t$ as before, let $N$ denote its associated future unit normal vector-field. We define the second fundamental form 
\begin{equation}\label{linmo1}
k_{ab}:=-\g(\D_{\partial_a}N,\partial_b)=\g(N,\D_{\partial_a}\partial_b)=N^\al\mathbf{\Gamma}_{\al a b},\qquad a,b\in\{1,2,3\}.
\end{equation} 
Let $\overline{g}_{jk}=\mathbf{g}_{jk}$ denote the induced (Riemannian) metric on $\Sigma_t$. Our main result in this section is the following:

\begin{proposition}\label{linmo2+}
We define the linear momentum
\begin{equation}\label{linmo2}
\mathbf{p}_a(t):=\frac{1}{8\pi}\lim_{R\to\infty}\int_{S_{R,t}}\pi_{ab}\frac{x^b}{|x|}\,dx,\qquad \pi_{ab}:=k_{ab}-\mathrm{tr}k\overg_{ab},
\end{equation}
where $S_{R,t}\subseteq \Sigma_t$ denotes the sphere of radius $R$ as before. Then the functions $\mathbf{p}_a$, $a\in\{1,2,3\}$ are well defined and independent of $t\in[0,\infty)$. Moreover, for any $t\geq 0$,
\begin{equation}\label{linmo4}
\mathbf{p}_a=-\frac{1}{16\pi}\int_{\mathbb{R}^3}\big\{2\partial_0\psi\partial_a\psi+(1/2)\sum_{m,n\in\{1,2,3\}}\partial_0\va_{mn}\partial_a\va_{mn}\big\}(t)\,dx+O(\varep_0^2\langle t\rangle^{-\kappa}).
\end{equation}
In particular 
\begin{equation}\label{linmo4.5}
\sum_{a\in\{1,2,3\}}\mathbf{p}_a^2\leq E_{ADM}^2,
\end{equation}
so the ADM mass $M_{ADM}:=\big(E_{ADM}^2-\sum_{a\in\{1,2,3\}}\mathbf{p}_a^2\big)^{1/2}\geq 0$ is well defined.
\end{proposition}

\begin{proof} This is similar to the proof of Lemma \ref{ADM20} and Proposition \ref{ADM4}. 

{\bf{Step 1.}}  Since $t$ is fixed, the quadratic and higher order terms do not contribute to the integral as $R\to\infty$, so we may redefine
\begin{equation}\label{linmo9}
\mathbf{p}_a(t)=\frac{1}{8\pi}\lim_{R\to\infty}\int_{S_{R,t}}\pi^1_{ab}\frac{x^b}{|x|}\,dx,\qquad \pi^1_{ab}:=k^1_{ab}-\delta_{ab}\delta_{jk}k^1_{jk}, 
\end{equation}
where $k^1_{ab}:=\mathbf{\Gamma}_{0ab}$ is the linear part of the second fundamental form $k$. We calculate
\begin{equation*}
\begin{split}
2k^1_{ab}&=\partial_ah_{0b}+\partial_bh_{0a}-\partial_0h_{ab}=R_aR_b(-2|\nabla|\rho-\partial_0F+\partial_0\uF)\\
&+(\in_{akl}R_bR_k+\in_{bkl}R_aR_k)(|\nabla|\omega_l+\partial_0\Omega_l)-\in_{apm}\in_{bqn}\partial_0R_pR_q\va_{mn},
\end{split}
\end{equation*}
using \eqref{zaq5l}. Recall that $R_mR_n\va_{mn}=0$ and $\delta_{mn}\va_{mn}=-2\tau$ (see \eqref{zaq2l} and \eqref{zaq21.1}). Thus
\begin{equation}\label{linmo11}
2\partial_b\pi^1_{ab}=-\in_{akl}|\nabla|R_k(|\nabla|\omega_l+\partial_0\Omega_l)+2\partial_a\partial_0\tau.
\end{equation}
Using \eqref{zaq2l} and then  \eqref{zaq5.1} we calculate
\begin{equation*}
\begin{split}
-\in_{akl}|\nabla|R_k(|\nabla|\omega_l+\partial_0\Omega_l)&=(\delta_{km}\delta_{an}-\delta_{kn}\delta_{am})R_kR_m\big (|\nabla|^2h_{0n}+\partial_0\partial_bh_{nb}\big)\\
&=-(\delta_{an}+R_aR_n)\big (|\nabla|^2h_{0n}+\partial_0\partial_bh_{nb}\big)\\
&=-(\delta_{an}+R_aR_n)\big (-\Delta h_{0n}+\partial_0^2h_{0n}+\partial_0E_n^{\geq 2}\big)\\
&=-(\delta_{an}+R_aR_n)\big (\mathcal{N}^h_{0n}+\partial_0E_n^{\geq 2}\big).
\end{split}
\end{equation*}
Using \eqref{zaq11} and the last two identities we calculate (after cancelling four terms)
\begin{equation}\label{linmo12}
2\partial_b\pi^1_{ab}=-\mathcal{N}^h_{0a}-\partial_aE_0^{\geq 2}-\partial_0E_a^{\geq 2}.
\end{equation}

{\bf{Step 2.}} Using Stokes theorem, \eqref{ADM8}, and \eqref{linmo12}, the limit in \eqref{linmo9} exists and
\begin{equation}\label{linmo13}
\mathbf{p}_a(t)=\frac{1}{16\pi}\int_{\mathbb{R}^3}2\partial_b\pi^1_{ab}\,dx=-\frac{1}{16\pi}\int_{\mathbb{R}^3}\{\mathcal{N}^h_{0a}+\partial_aE_0^{\geq 2}+\partial_0E_a^{\geq 2}\}\,dx.
\end{equation}
Moreover, using \eqref{linmo11}, 
\begin{equation*}
\partial_0(2\partial_b\pi^1_{ab})=-\in_{akl}\partial_k(|\nabla|\partial_0\omega_l+\partial_0^2\Omega_l)+2\partial_a\partial_0^2\tau.
\end{equation*}
In view of \eqref{zaq11.2}--\eqref{zaq11} and \eqref{ADM8}, we have 
\begin{equation*}
\big\|(|\nabla|\partial_0\omega_l+\partial_0^2\Omega_l)(t)\big\|_{L^{4/3}}+\big\|\partial_0^2\tau(t)\|_{L^{4/3}}\lesssim\varep_0\langle t\rangle^4,
\end{equation*}
and the same argument as in the proof of Proposition \ref{ADM4} shows that the components $\mathbf{p}_a(t)$ are constant in time.

To prove the identity \eqref{linmo4} we we use the formula \eqref{linmo13}, extract the time decaying components, and then let $t\to\infty$, just like in the proofs of Lemma \ref{ADM20} and Proposition \ref{ADM4}. Indeed, all the cubic and higher order terms and all the quadratic null terms lead to time decaying contributions, thus using \eqref{sac1}--\eqref{zaq21} and \eqref{zaq5.1},
\begin{equation*}
\mathcal{N}^h_{0a}+\partial_aE_0^{\geq 2}+\partial_0E_a^{\geq 2}\sim 2\partial_0\psi\partial_a\psi+\mathcal{Q}^2_{0a}-P_{0a}^2+\partial_0E_a^{2}+\partial_aE_0^{2},
\end{equation*}
where $F\sim G$ means $\|F-G\|_{L^1}\lesssim \varep_0^2\langle t\rangle^{-\kappa}$ as before. Also, as we know from the proof of Proposition \ref{ADM4}, derivative terms of the form $\partial(h\cdot\partial h)$ do not contribute to the integral in the limit $t\to\infty$. As in the proof of Lemma \ref{ADM20}, the terms $\mathcal{Q}^2_{0a}$, $\partial_0E_a^{2}$, and $\partial_aE_0^{2}$ are sums of derivatives and $L^1$ acceptable errors. The only terms that contribute in the limit are the terms $2\partial_0\psi\partial_a\psi$ and 
$(1/2)R_pR_q\partial_0\va_{mn}R_pR_q\partial_a\va_{mn}$ coming from $-P^2_{0a}$ after removing the $L^1$ acceptable errors (as in Proposition \ref{ADM4}). In view of \eqref{linmo13}, this leads to the desired formula \eqref{linmo4}. Finally, the inequality \eqref{linmo4.5} follows using also \eqref{ADM4.6} and letting $t\to\infty$.
\end{proof}

\section{Asymptotically optical functions and the Bondi energy}\label{OptBo}

Our final application concerns the construction of Bondi energy functions, with good monotonicity properties along null infinity.  We would like to thank Yakov Shlapentokh-Rothman for useful discussions on this topic.

\subsection{Almost optical functions and the Friedlander fields} In order to get precise information on the asymptotic behavior of the metric in the physical space we need to understand the bending of the light cones caused by the long-range effect of the nonlinearity (i.e. the modified scattering). 

In the Minkowski space, the outgoing light cones correspond to the level sets of the optical function $u^0=|x|-t$. In our case, the analogous objects we use are what we call {\it{almost (or asymptotically) optical functions}} $u$, which are close to $u^0$ but better adapted to the null geometry of our problem. Recall that the metrics we consider here have slow $O(\langle r\rangle^{-1+})$ decay at infinity, and we expect a nontrivial deviation that is not radially isotropic. 

We first define and construct a suitable class of asymptotically optical functions.

\begin{lemma}\label{AlmostOpticalLem}
There exists a $C^4$ almost optical function $u:M'\to\mathbb{R}$ satisfying the properties
\begin{equation}\label{AlmostOptical}
\begin{split}
u(x,t)&=\vert x\vert-t+u^{cor}(x,t),\qquad \g^{\alpha\beta}\partial_\alpha u\partial_\beta u=O(\varepsilon_0\langle r\rangle^{-2+6\delta'})
\end{split}
\end{equation}
and, for any $\mu\in\{0,1,2,3\}$,
\begin{equation}\label{EstimuGradu}
\begin{split}
&u^{cor}=O(\varepsilon_0\langle r\rangle^{3\delta'}),\qquad \partial_\mu u^{cor}=O(\varepsilon_0\langle r\rangle^{3\delta'-1}),\qquad\partial_\mu(L^\alpha\partial_\alpha u^{cor})=O(\varep_0\langle r\rangle^{3\delta'-2}).
\end{split}
\end{equation}
In addition, $u^{cor}$ is close to $\Theta_{wa}/|x|$ (see \eqref{bax3}) in the vicinity of the light cone,
\begin{equation}\label{uCorThetaClose}
\begin{split}
\Big\vert u^{cor}(x,t)-\frac{\Theta_{wa}(x,t)}{|x|}\Big\vert\lesssim \varep_0\langle r\rangle^{-1+3\delta'}(\langle r\rangle^{p_0}+\langle t-\vert x\vert \rangle),\quad \text{ if }(x,t)\in M',\,\big|t-|x|\big|\leq t/10.
\end{split}
\end{equation}
\end{lemma}

\begin{remark} The classical approach, see for example \cite{ChKl}, is to construct exact optical functions, satisfying the stronger identity $\g^{\alpha\beta}\partial_\alpha u\partial_\beta u=0$ instead of the approximate identity in \eqref{AlmostOptical}. We could do this too, but we prefer to work here with almost optical functions instead of exact optical functions because they are much easier to construct and their properties still suffice for our two main applications (the improved peeling estimates in Proposition \ref{PeelingR2} and the construction of the Bondi energy in Theorem \ref{BondThm}).
\end{remark}

\begin{proof} We define the function $H_L:M\to\mathbb{R}$, 
\begin{equation}\label{Alop1}
H_L:=\frac{1}{2}L^\al L^\be h_{\al\be}=-\frac{1}{2}g_{\geq 1}^{\al\be}\partial_\al u^0\partial_\be u^0+\frac{1}{2}g^{\al\be}_{\geq 2}\partial_\al u^0\partial_\be u^0,
\end{equation}
where the identity holds due to \eqref{zaq22}. Notice that
\begin{equation}\label{Alop1.2}
H_L=O(\varep_0\langle t+r\rangle^{-1+\delta'}),\qquad \partial_\mu H_L=O(\varep_0r^{-1}\langle t+r\rangle^{-1+3\delta'}),
\end{equation}
in $M$, as a consequence of \eqref{winr25}, \eqref{ImprovedBoundsInside}, and \eqref{extrader1}. We will define $u^{cor}$ such that
\begin{equation}\label{Alop1.4}
L^\al\partial_\al u^{cor}=H_L,
\end{equation}
in addition to the bounds in \eqref{EstimuGradu}--\eqref{uCorThetaClose}.

{\bf{Step 1.}} For $s,b\in[0,\infty)$ we define the projections $\Pi^-_b$ and $\Pi^+_b$
\begin{equation}\label{Alop2}
\begin{split}
&\Pi^-_bH_L(s):=L^\al L^\be\mathcal{F}^{-1}\{\varphi_{\leq 0}(\langle b\rangle^{p_0}\xi)\widehat{h_{\al\be}}(\xi,s)\},\\
&\Pi^+_bH_L(s):=L^\al L^\be\mathcal{F}^{-1}\{\varphi_{\geq 1}(\langle b\rangle^{p_0}\xi)\widehat{h_{\al\be}}(\xi,s)\},
\end{split}
\end{equation}
as in section \ref{waz}, where $p_0=0.68$. Then we define the correction $u^{cor}$ by integrating the low frequencies of $H_L$ from $0$ to $t$ and the high frequencies from $t$ to $\infty$. More precisely, let
\begin{equation}\label{Alop5}
\begin{split}
u^{cor}_1(x,t):=&-\int_{t}^\infty(\Pi^+_{s}H_{L})(x+(s-t)x/|x|,s)\,ds+\int_{0}^{t}(\Pi_s^-H_L)(sx/|x|,s)\,ds\\
&-\int_{t}^\infty\left\{(\Pi^-_{s}H_L)(x+(s-t)x/|x|,s)-(\Pi^-_sH_L)(sx/|x|,s)\right\}\,ds,
\end{split}
\end{equation}
\begin{equation}\label{Alop5.2}
\begin{split}
u^{cor}_2(x,t):=&-\int_{|x|}^\infty (\Pi^+_{s}H_{L})(sx/|x|,s+t-|x|)\,ds+\int_{0}^{|x|}(\Pi_s^-H_L)(sx/|x|,s)\,ds\\
&-\int_{|x|}^\infty\left\{(\Pi^-_{s}H_L)(sx/|x|,s+t-|x|)-(\Pi^-_sH_L)(sx/|x|,s)\right\}\,ds.
\end{split}
\end{equation}
We fix a smooth function $\chi_1:\mathbb{R}\to[0,1]$ supported in $(-\infty,2]$ and equal to $1$ in $(-\infty,1]$, let $\chi_2:=1-\chi_1$, and define
\begin{equation}\label{Alop5.4}
u^{cor}(x,t):=u^{cor}_1(x,t)\chi_1(|x|-t)+u^{cor}_2(x,t)\chi_2(|x|-t).
\end{equation}

Notice that, formally, one can rewrite the formula \eqref{Alop5} as
\begin{equation}\label{Alop5.5}
u^{cor}_1(x,t)\approx-\int_{t}^\infty H_{L}(x+(s-t)x/|x|,s)\,ds+\int_{0}^{\infty}(\Pi_s^-H_L)(sx/|x|,s)\,ds,
\end{equation}
which is consistent with the desired transport identity \eqref{Alop1.4}. However, the two infinite integrals in \eqref{Alop5.5} do not converge, and we need to reorganize them as in \eqref{Alop5} to achieve convergence.  A similar remark applies to the definition of $u^{cor}_2$ in \eqref{Alop5.2}.

We prove now the bounds \eqref{EstimuGradu}. Notice that, for any $b\geq 0$ and $(y,s)\in M$,
\begin{equation}\label{Alop6}
\begin{split}
\vert \Pi^-_bH_L(y,s)\vert+|y|\vert\nabla_{y,s}(\Pi^-_bH_L)(y,s)\vert&\lesssim\varepsilon_0 \langle |y|+s\rangle^{3\delta^\prime-1},\\
\langle y\rangle\vert (\Pi^+_{b}H_L)(y,s)\vert&\lesssim\varep_0\langle |y|+s\rangle^{3\delta^\prime-1}\langle b\rangle^{p_0},
\end{split}
\end{equation}
due to \eqref{Alop1.2}. Finally, for $(x,t)\in M'$ and $s\geq\max(\rho,0)$, $\rho\in[\max(t,|x|)-3,\max(t,|x|)+3]$, using the bounds in the first line of \eqref{Alop6} we have
\begin{equation}\label{Alop7}
\vert (\Pi^-_{b}H_L)(x+(s-\rho)x/|x|,s+t-\rho)-(\Pi^-_bH_L)(sx/|x|,s)\vert\lesssim \varep_0\langle t-|x|\rangle\langle s\rangle^{3\delta^\prime-2}.
\end{equation}
Using the definitions \eqref{Alop5}--\eqref{Alop5.4}, it follows that the functions $u^{cor,1},u^{cor,2}$ are well defined for any $(x,t)\in M'$, and moreover satisfy the estimates
\begin{equation}\label{MainBounds}
\chi_a(|x|-t)\vert u^{cor}_a(x,t)\vert\lesssim\varepsilon_0\langle x\rangle^{3\delta'},\qquad \chi_a(|x|-t)\vert \partial_\mu u^{cor}_a(x,t)\vert\lesssim\varepsilon_0\langle x\rangle^{3\delta^\prime-1},
\end{equation}
for $a\in\{1,2\}$ and $\mu\in\{0,1,2,3\}$. Using again the gradient bounds in \eqref{Alop6} and \eqref{Alop1.2}, we also have $\vert u^{cor}_1(x,t)-u^{cor}_2(x,t)\vert\lesssim\varepsilon_0\langle x\rangle^{3\delta^\prime-1}$ if $|x|-t\in[1,2]$. Finally,
\begin{equation*}
L^\al\partial_\al u^{cor}_1=H_L \text{ if }|x|-t\leq 2\qquad\text{ and }\qquad L^\al\partial_\al u^{cor}_2=H_L \text{ if }|x|-t\geq 1,
\end{equation*}
and the desired bounds in \eqref{EstimuGradu} follow.

{\bf{Step 2.}} We calculate now in $M'$
\begin{equation*}
\begin{split}
\g^{\alpha\beta}\partial_\alpha u\partial_\beta u&=(m^{\al\be}+g_{\geq 1}^{\al\be})\partial_\alpha (u^0+u^{cor})\partial_\beta (u^0+u^{cor})\\
&=2m^{\al\be}\partial_\alpha u^0\partial_\beta u^{cor}+g_{\geq 1}^{\al\be}\partial_\alpha u^0\partial_\beta u^0+O(\varep_0\langle r\rangle^{-2+6\delta'})\\
&=2L^\be\partial_\be u^{cor}-2H_L+O(\varep_0\langle r\rangle^{-2+6\delta'}),
\end{split}
\end{equation*}
using \eqref{Alop1} and \eqref{EstimuGradu}. The bounds in \eqref{AlmostOptical} follow using also \eqref{Alop1.4}.

Finally, to prove \eqref{uCorThetaClose} we examine the definition \eqref{bax3} and notice that
\begin{equation}\label{Alop9}
\frac{\Theta_{wa}(x,t)}{|x|}=\int_0^t(\Pi_s^-H_L)(sx/|x|,s)\,ds.
\end{equation}
Therefore, using \eqref{Alop6}, if $(x,t)\in M'$ and $t=|x|$
\begin{equation*}
\Big\vert u^{cor}(x,t)-\frac{\Theta_{wa}(x,t)}{|x|}\Big\vert\lesssim \varep_0\langle r\rangle^{-1+3\delta'+p_0}.
\end{equation*}
The desired estimates \eqref{uCorThetaClose} follow using also the bounds $\partial_r u^{cor}=O(\varepsilon_0\langle r\rangle^{3\delta'-1})$.
\end{proof}

We will prove now asymptotic formulas in the physical space for some of the metric components and for the Klein-Gordon field. These formulas will be used in the Bondi energy analysis in subsection \ref{BondiMass} below.

Recall the definitions
\begin{equation}\label{Alop10}
\begin{split}
\widehat{V_\ast^G}(\xi,t)&=\widehat{V^G}(\xi,t)e^{-i\Theta_{wa}(\xi,t)},\qquad G\in\{F,\omega_a,\vartheta_{ab}\},\\
\widehat{V_\ast^\psi}(\xi,t)&=\widehat{V^\psi}(\xi,t)e^{-i\Theta_{kg}(\xi,t)},
\end{split}
\end{equation}
see \eqref{bax4} and \eqref{Nor5}. It follows from \eqref{bax12} and \eqref{Nor40} that there are functions $V_\infty^G\in Z_{wa}$, $G\in\{F,\omega_a,\vartheta_{ab}\}$, and $V_\infty^\psi\in Z_{kg}$ such that, for any $t\geq 0$,
\begin{equation}\label{Alop10.1}
\sum_{G\in\{F,\omega_a,\vartheta_{ab}\}}\|V^G_\ast(t)-V^G_\infty\|_{Z_{wa}}+\|V^\psi_\ast(t)-V^\psi_\infty\|_{Z_{kg}}\lesssim\varep_0\langle t\rangle^{-\delta/2}.
\end{equation}

We define the smooth characteristic function $\chi_{wa}$ of the wave region by
\begin{equation}\label{Alop10.3}
\chi_{wa}:\{(x,t)\in M:t\geq 8\}\to [0,1],\qquad \chi_{wa}(x,t):=\varphi_{\leq 0}((|x|-t)/t^{0.4}).
\end{equation}
We define also the function
\begin{equation}\label{Alop10.4}
\nu_{kg}:\{(x,t)\in M: t\geq 8,\,|x|<t\}\to\mathbb{R}^3,\qquad\nu_{kg}(x,t):=\frac{x}{\sqrt{t^2-|x|^2}},
\end{equation}
such that $\nu_{kg}(x,t)$ is the critical point of the function $\xi\to x\cdot\xi-t\sqrt{1+|\xi|^2}$. We are now ready to state our main proposition describing the solutions in the physical space.

\begin{proposition}\label{Alop20}
For any $G\in \{F,\omega_a,\vartheta_{ab}\}$ and $t\in[8,\infty)$ we have
\begin{equation}\label{Alop21}
\begin{split}
U^G(x,t)&=\frac{-i\chi_{wa}(x,t)}{4\pi^2|x|}\int_0^\infty e^{i\rho u(x,t)}\varphi_{[-k_0,k_0]}(\rho)\widehat{V^G_\ast}(\rho x/|x|,t)\rho\,d\rho+U^G_{rem}(x,t),\\
U^{\psi}(x,t)&=\frac{\mathbf{1}_+(t-|x|)}{\sqrt{8\pi^3}e^{i\pi/4}}\frac{e^{-i\sqrt{t^2-|x|^2}}t}{(t^2-|x|^2)^{5/4}}e^{i\Theta_{kg}(\nu_{kg}(x,t),t)}\widehat{P_{[-k_0,k_0]}V^\psi_\ast}(\nu_{kg}(x,t),t)+U^\psi_{rem}(x,t),
\end{split}
\end{equation}
where $u$ is an almost optical function satisfying \eqref{AlmostOptical}--\eqref{uCorThetaClose} and $k_0$ denotes the smallest integer for which $2^{k_0}\geq t^{\delta'}$. The remainders $U^G_{rem}$ and $U^\psi_{rem}$ satisfy the $L^2$ bounds
\begin{equation}\label{Alop21.2}
\sum_{G\in\{F,\omega_a,\vartheta_{ab}\}}\|U^G_{rem}(t)\|_{L^2}+\|U^\psi_{rem}(t)\|_{L^2}\lesssim\varep_0t^{-\delta},\qquad\text{ for any }t\geq 8.
\end{equation}
\end{proposition}

\begin{proof} {\bf{Step 1.}} We prove first the conclusions concerning the variables $U^G$. We start from the formula
\begin{equation}\label{Alop22}
\widehat{U^G}(\xi,t)=e^{-it|\xi|}e^{i\Theta_{wa}(\xi,t)}\widehat{V^G_\ast}(\xi,t),
\end{equation}
and extract acceptable $L^2$ remainders until we reach the desired formula in \eqref{Alop21}.

We may assume that $t\gg 1$ and let $J_0$ denote the smallest integers for which $2^{J_0}\geq t^{1/3}$. We define
\begin{equation}\label{Alop23}
\begin{split}
V^G_{\ast,1}&:=(I-P_{[-k_0,k_0]})V^G_\ast+P_{[-k_0-2,k_0+2]}[\varphi_{\geq J_0+1}\cdot P_{[-k_0,k_0]}V^G_\ast],\\
V^G_{\ast,2}&:=P_{[-k_0-2,k_0+2]}[\varphi_{\leq J_0}\cdot P_{[-k_0,k_0]}V^G_\ast],
\end{split}
\end{equation}
and notice that $V^G_{\ast}:=V^G_{\ast,1}+V^G_{\ast,2}$. We show first that
\begin{equation}\label{Alop24}
\|V^G_{\ast,1}(t)\|_{L^2}\lesssim \varep_0t^{-\delta'/4}.
\end{equation}
Indeed, notice that $\|(I-P_{[-k_0,k_0]})V^G_\ast(t)\|_{L^2}\lesssim\varep_0t^{-\delta'/4}$, due to \eqref{winr1} and \eqref{Alop10}. To bound the remaining term we examine the definition \eqref{bax3} and notice that
\begin{equation}\label{Alop24.5}
\big|D^\al_\xi [e^{\pm i\Theta_{wa}(\xi,t)}]\big|\lesssim_\al t^{|\alpha|(1-p_0+2\delta')}t^{2\delta'}\qquad\text{ if }t^{-\delta'}\lesssim |\xi|\lesssim t^{\delta'}.
\end{equation}
Let $A_t$ denote the operator defined by the Fourier multiplier $\xi\to e^{- i\Theta_{wa}(\xi,t)}\varphi_{[-k_0-2,k_0+2]}(\xi)$. Notice that
\begin{equation*}
P_{[-k_0,k_0]}V^G_\ast=A_t[P_{[-k_0,k_0]}V^G]=A_t[\varphi_{\leq J_0-4}\cdot P_{[-k_0,k_0]}V^G]+A_t[\varphi_{\geq J_0-3}\cdot P_{[-k_0,k_0]}V^G].
\end{equation*}
If view of \eqref{Alop24.5}, the kernel of the operator $A_t$ decays rapidly if $|x|\gtrsim 2^{J_0}$, thus
\begin{equation}\label{Alop25.7}
\big\|\varphi_{\geq J_0+1}\cdot A_t[\varphi_{\leq J_0-4}\cdot P_{[-k_0,k_0]}V^G]\big\|_{L^2}\lesssim\varep_0 t^{-1}.
\end{equation}
Moreover, using \eqref{winr23},
\begin{equation}\label{Alop25.8}
\big\|A_t[\varphi_{\geq J_0-3}\cdot P_{[-k_0,k_0]}V^G]\big\|_{L^2}\lesssim\varep_0 2^{-J_0}t^{2\delta'}\lesssim\varep_0 t^{-1/3+2\delta'}.
\end{equation}
The bounds \eqref{Alop24} follow, using the definition. Therefore
\begin{equation}\label{Alop26}
\|U^G_{rem,1}(t)\|_{L^2}\lesssim\varep_0t^{-\delta'/4}\qquad\text{ where }\qquad U^G_{rem,1}(t):=\mathcal{F}^{-1}\{e^{-it|\xi|}e^{i\Theta_{wa}(\xi,t)}\widehat{V^G_{\ast,1}}(\xi,t)\}.
\end{equation}

We define now 
\begin{equation}\label{Alop28}
U^G_{rem,2}:=(1-\chi_{wa})\cdot \mathcal{F}^{-1}\{e^{-it|\xi|}e^{i\Theta_{wa}(\xi,t)}\widehat{V^G_{\ast,2}}(\xi,t)\}.
\end{equation}
Using integration by parts in $\xi$ (Lemma \ref{tech5}) and the formula \eqref{Alop23} we have rapid decay,
\begin{equation}\label{Alop29}
\|U^G_{rem,2}(t)\|_{L^2}\lesssim\varep_0t^{-1}.
\end{equation}

To estimate the main term $\chi_{wa}\cdot \mathcal{F}^{-1}\{e^{-it|\xi|}e^{i\Theta_{wa}(\xi,t)}\widehat{V^G_{\ast,2}}(\xi,t)\}$ we write it in the form
\begin{equation}\label{Alop30}
U_3^G(x,t):=\frac{\chi_{wa}(x,t)}{8\pi^3}\int_{\mathbb{R}^3}e^{ix\cdot\xi}e^{-it|\xi|}e^{i\Theta_{wa}(\xi,t)}\widehat{V^G_{\ast,2}}(\xi,t)\,d\xi.
\end{equation}
We can extract more remainders by inserting angular cutoffs. Notice that if we insert the factor $\varphi_{\geq 1}(t^{0.49}(x/|x|-\xi/|\xi|))$ in the integral above then the corresponding contribution is a rapidly decreasing $L^2$ remainder. Passing to polar coordinates $x=r\omega$, $\xi=\rho\theta$, $\omega,\theta\in\mathbb{S}^2$, it remains to estimate the integral
\begin{equation}\label{Alop31}
U_4^G(x,t):=\frac{\chi_{wa}(x,t)}{8\pi^3}\int_0^\infty\int_{\mathbb{S}^2}e^{ir\rho\omega\cdot\theta}e^{-it\rho}e^{i\rho\Theta_{wa}(\theta,t)}\widehat{V^G_{\ast,2}}(\rho\theta,t)\varphi_{\leq 0}(t^{0.49}(\omega-\theta))\,\rho^2 d\rho d\theta.
\end{equation}

For any $\omega,\theta\in\mathbb{S}^2$ and $\rho\in[2^{-10}t^{-\delta'},2^{10}t^{\delta'}]$ we have
\begin{equation}\label{Alop33}
|e^{i\rho\Theta_{wa}(\theta,t)}-e^{i\rho\Theta_{wa}(\omega,t)}|+|\widehat{V^G_{\ast,2}}(\rho\theta,t)-\widehat{V^G_{\ast,2}}(\rho\omega,t)|\lesssim \varep_0t^{4\delta'}|\theta-\omega|.
\end{equation}
Indeed, using the definitions we have $|\Omega_\theta[e^{\pm i\rho\Theta_{wa}(\theta,t)}]|\lesssim t^{2\delta'}$, where $\Omega_\theta$ is any of the rotation vector-fields in the variable $\theta\in\mathbb{S}^2$. The bounds \eqref{Alop33} follow using also \eqref{vcx1.3*}. Therefore we can further replace the angular variable $\theta$ with $\omega$ in two places in the integral in \eqref{Alop31}, at the expenses of acceptable errors. It remains to estimate the integral
\begin{equation}\label{Alop34}
U_5^G(x,t):=\frac{\chi_{wa}(x,t)}{8\pi^3}\int_0^\infty e^{-it\rho}e^{i\rho\Theta_{wa}(\omega,t)}\widehat{V^G_{\ast,2}}(\rho\omega,t)\Big\{\int_{\mathbb{S}^2}e^{ir\rho\omega\cdot\theta}\varphi_{\leq 0}(t^{0.49}(\omega-\theta))\,d\theta\Big\}\rho^2 d\rho.
\end{equation}

The integral over $\mathbb{S}^2$ in \eqref{Alop34} does not depend on $\omega$ and can be calculated explicitly. We may assume $\omega=(0,0,1)$ and use spherical coordinates to write this integral in the form 
\begin{equation*}
\begin{split}
2\pi\int_{0}^{\pi}e^{ir\rho\cos y}&\varphi_{\leq 0}(t^{0.49}2\sin(y/2))\,\sin y \,dy=8\pi\int_{0}^{1}e^{ir\rho(1-2z^2)}\varphi_{\leq 0}(t^{0.49}2z)z\,dz\\
&=\frac{2\pi}{r\rho}e^{ir\rho}\int_0^\infty e^{-i\al}\varphi_{\leq 0}\Big(\frac{2t^{0.49}}{\sqrt{2r\rho}}\sqrt\al\Big)\,d\al=\frac{-2\pi i}{r\rho}e^{ir\rho}+O(t^{-2}).
\end{split}
\end{equation*}
We substitute this into \eqref{Alop34}, and it remains to estimate the integral
\begin{equation}\label{Alop35}
U_6^G(x,t):=\frac{-i\chi_{wa}(x,t)}{4\pi^2r}\int_0^\infty e^{-it\rho}e^{i\rho \Theta_{wa}(\omega,t)}\widehat{V^G_{\ast,2}}(\rho\omega,t)e^{ir\rho}\rho\,d\rho.
\end{equation}

In view of \eqref{Alop25.7}--\eqref{Alop25.8} we can replace now the factor $\widehat{V^G_{\ast,2}}(\rho\omega,t)$ in \eqref{Alop35} with $\mathcal{F}\{P_{[-k_0,k_0]}V^G_\ast\}(\rho\omega,t)$. Then we replace $\Theta_{wa}(\omega,t)=\Theta_{wa}(x/|x|,t)$ with $u^{cor}(x,t)$, up to acceptable errors (due to \eqref{uCorThetaClose}). The desired formula in \eqref{Alop21} follows, since $u(x,t)=|x|-t+u^{cor}(x,t)$.

{\bf{Step 2.}} The proof for the Klein-Gordon variable $U^\psi$ is similar. We start from the formula
\begin{equation}\label{Alop40}
\widehat{U^\psi}(\xi,t)=e^{-it\sqrt{1+|\xi|^2}}e^{i\Theta_{kg}(\xi,t)}\widehat{V^\psi_\ast}(\xi,t),
\end{equation}
and extract acceptable $L^2$ remainders until we reach the desired formula in \eqref{Alop21}. We may assume $t\gg 1$, set $J_0$ such that $2^{J_0}\approx t^{1/3}$, and define, as in \eqref{Alop23},
\begin{equation}\label{Alop41}
\begin{split}
V^\psi_{\ast,1}&:=(I-P_{[-k_0,k_0]})V^\psi_\ast+P_{[-k_0-2,k_0+2]}[\varphi_{\geq J_0+1}\cdot P_{[-k_0,k_0]}V^\psi_\ast],\\
V^\psi_{\ast,2}&:=P_{[-k_0-2,k_0+2]}[\varphi_{\leq J_0}\cdot P_{[-k_0,k_0]}V^\psi_\ast].
\end{split}
\end{equation}
The definition \eqref{Nor4} shows that
\begin{equation}\label{Alop42}
\big|D^\al_\xi [e^{\pm i\Theta_{kg}(\xi,t)}]\big|\lesssim_\al t^{|\alpha|(1-p_0+2\delta')}t^{2\delta'}\qquad\text{ if }t^{-\delta'}\lesssim |\xi|\lesssim t^{\delta'}.
\end{equation}
As in \eqref{Alop24} it follows that $\|V^\psi_{\ast,1}(t)\|_{L^2}\lesssim \varep_0t^{-\delta'/4}$, so the contribution of $V^\psi_{\ast,1}$ is an acceptable remainder. The contribution of $V^\psi_{\ast,2}$ is also a remainder in the region $\{|x|\geq t\}$, due to Lemma \ref{tech5}. It remains to estimate the main contribution,
\begin{equation}\label{Alop30.3}
U_2^\psi(x,t):=\frac{\mathbf{1}_+(t-|x|)}{8\pi^3}\int_{\mathbb{R}^3}e^{ix\cdot\xi}e^{-it\sqrt{1+|\xi|^2}}e^{i\Theta_{kg}(\xi,t)}\widehat{V^\psi_{\ast,2}}(\xi,t)\,d\xi.
\end{equation}

We can insert a cutoff function of the form $\varphi_{\leq 0}(t^{0.49}(\xi-\nu_{kg}(x,t)))$ in the integral above to localize near the critical point, and then replace $e^{i\Theta_{kg}(\xi,t)}$ and $\widehat{V^\psi_{\ast,2}}(\xi,t)$ with $e^{i\Theta_{kg}(\nu_{kg}(x,t),t)}$ and $\widehat{V^\psi_{\ast,2}}(\nu_{kg}(x,t),t)$ respectively, at the expense of acceptable remainders (using \eqref{Alop42} and \eqref{vcx1.15} respectively). 

Therefore, the desired identity in \eqref{Alop21} holds if $\nu_{kg}(x,t)\geq 2^{k_0+4}$ (with the main term vanishing). On the other hand, if $\nu_{kg}(x,t)\leq 2^{k_0+4}$ (thus $t-|x|\gtrsim t^{1-4\delta'}$) then the remaining $\xi$ integral in \eqref{Alop30.3} can be estimated explicitly,
\begin{equation*}
\int_{\mathbb{R}^3}e^{ix\cdot\xi}e^{-it\sqrt{1+|\xi|^2}}\varphi_{\leq 0}(t^{0.49}(\xi-\nu_{kg}(x,t)))\,d\xi=e^{-i\sqrt{t^2-|x|^2}}\frac{e^{-i\pi/4}(2\pi)^{3/2}t}{(t^2-|x|^2)^{5/4}}+O(t^{-7/4}),
\end{equation*}
where the approximate identity follows from the standard stationary phase formula. After these reductions, the remaining main term is
\begin{equation*}
U_3^\psi(x,t):=\frac{\mathbf{1}_+(t-|x|)}{\sqrt{8\pi^3}}e^{i\Theta_{kg}(\nu_{kg}(x,t),t)}\widehat{V^\psi_{\ast,2}}(\nu_{kg}(x,t),t)\frac{e^{-i\sqrt{t^2-|x|^2}}e^{-i\pi/4}t}{(t^2-|x|^2)^{5/4}}.
\end{equation*}
The desired conclusion in \eqref{Alop21} follows using \eqref{vcx1.15}.
\end{proof}

\subsection{The Bondi energy}\label{BondiMass} We can define now a more refined concept of energy function. For this we fix $t\geq 1$, define the hypersurface $\Sigma_t:=\{(x,t)\in M:\,x\in\mathbb{R}^3\}$, and let $\overline{g}_{jk}=\mathbf{g}_{jk}$ denote the induced (Riemannian) metric on $\Sigma_t$. Let $\overline{g}^{jk}$ the inverse of the matrix $\overline{g}_{jk}$, $\overline{g}^{jk}\overline{g}_{jn}=\delta^k_n$, let $\overline{D}$ denote the covariant derivative on $\Sigma_t$ induced by the metric $\overline{g}$. Notice that
\begin{equation}\label{Bond1}
|\overline{g}^{jk}-\mathbf{g}^{jk}|\lesssim \varep_0^2\langle t+r\rangle^{-2+2\delta'},\qquad\overline{\Gamma}_{njk}=\mathbf{\Gamma}_{njk},
\end{equation}
for any $n,j,k\in\{1,2,3\}$. With $u$ an almost optical function as in \eqref{AlmostOptical} let
\begin{equation}\label{Bond2}
\mathbf{n}_j:=\partial_ju(\overline{g}^{ab}\partial_au\partial_bu)^{-1/2},
\end{equation}
denote the unit vector-field in $\Sigma'_t:=\{x\in\Sigma_t:\,|x|\geq 2^{-8}t\}$, normal to the level sets of the function $u$. In this section we use the metric $\overline{g}$ to raise and lower indices.

We fix a function $u$ as in Lemma \ref{AlmostOpticalLem}. For $R\in\mathbb{R}$ and $t$ large (say $t\geq 2|R|+10$) we define the modified spheres $S^u_{R,t}:=\{x\in\Sigma_t:\,u(x,t)=R\}$. We would like to define
\begin{equation}\label{Bond5}
E_{Bondi}(R):=\frac{1}{16\pi}\lim_{t\to\infty}\int_{S^u_{R,t}} W_j\mathbf{n}^j\,d\sigma,
\end{equation}
for a suitable vector-field $W$, where $d\sigma= d\sigma(\overline{g})$ is the surface measure induced by the metric $\overline{g}$. The main issue is the convergence as $t\to\infty$ of the integrals in the formula above, for $R$ fixed. For this we need to be careful both with the choice of surfaces of integration $S^u_{R,t}$ and the choice of the vector-field $W$ we integrate. 

Our main theorem in this section is the following:

\begin{theorem}\label{BondThm}
Let 
\begin{equation}\label{Bond7}
W_j:=\overg^{ab}(\partial_a h_{jb}-\partial_j h_{ab}).
\end{equation}
Then the limit in \eqref{Bond5} exists, and $E_{Bondi}:\mathbb{R}\to \mathbb{R}$ is a well-defined continuous and increasing function. Moreover
\begin{equation}\label{Bond7.1}
\begin{split}
\lim_{R\to-\infty}E_{Bondi}(R)&=E_{KG}:=\frac{1}{16\pi}\|V^\psi_\infty\|_{L^2}^2,\\
\lim_{R\to\infty}E_{Bondi}(R)&=E_{ADM}.
\end{split}
\end{equation}
\end{theorem}

\begin{proof} {\bf{Step 1.}} We decompose $W=W^1+W^{\geq 2}$,
\begin{equation}\label{Bond7.5}
W^1_j:=\delta^{ab}(\partial_a h_{jb}-\partial_j h_{ab}),\qquad W^{\geq 2}_j:=\overg^{ab}_{\geq 1}(\partial_a h_{jb}-\partial_j h_{ab}).
\end{equation}
To calculate the linear contribution of $W^1_j$ we use the divergence theorem
\begin{equation}\label{Bond8}
\int_{S^u_{R,t}} W^1_j\mathbf{n}^j\,d\sigma=\int_{B^u_{R,t}} \overline{D}^jW^1_j\,d\mu,
\end{equation}
where $B^u_{R,t}:=\{x\in\Sigma_t:\,u(x,t)\leq R\}$ is the ball of radius $R$. Then we calculate
\begin{equation}\label{Bond9}
\overline{D}^jW^1_j=\overg^{jk}\big\{\partial_kW^1_j-{\overline{\Gamma}^m}_{jk}W^1_m\big\}=\delta^{jk}\partial_kW^1_j+\overg^{jk}_1\partial_kW^1_j-\delta^{jk}{\overline{\Gamma}^m}_{jk}W^1_m+E_1,
\end{equation} 
where $E_1$ is a cubic and higher order term and $\overg^{jk}_1=g^{jk}_1=-h_{jk}$ is the linear part of the metric $\overg^{jk}$. As in the proof of Lemma \ref{ADM20}, the cubic term satisfies the bounds $\|E_1(t)\|_{L^1}\lesssim\varep_0^2t^{-\kappa}$. The point of the identity \eqref{Bond9} is that the linear part $\delta^{jk}\partial_kW^1_j=-2\Delta\tau$ satisfies the equations \eqref{ADM5}. We can therefore apply the results of Lemma \ref{ADM20}, and write
\begin{equation}\label{Bond10}
\begin{split}
&\overline{D}^jW^1_j\sim-\delta_{jk}P^2_{jk}+\{(\partial_0\psi)^2+\psi^2+\partial_j\psi\partial_j\psi\}+\partial_jO^2_j-h_{jk}\partial_jW^1_k-\delta^{jk}{\overline{\Gamma}^m}_{jk}W^1_m\\
&\sim-\delta_{jk}P^2_{jk}+\{(\partial_0\psi)^2+\psi^2+\partial_j\psi\partial_j\psi\}+\partial_j\{O^2_j-h_{jk}W^1_k\}+\{\partial_jh_{jk}W^1_k-{\overline{\Gamma}}_{mjj}W^1_m\}.
\end{split}
\end{equation} 
As in the previous section, $F\sim G$ means $\|F-G\|_{L^1}\lesssim\varep_0^2t^{-\kappa}$, and $O^2_j$ is defined in \eqref{ADM21.5}.

We show now that the last semilinear term in \eqref{Bond10} is an acceptable $L^1$ error. Indeed, using \eqref{zaq5.1} and the definitions,
\begin{equation*}
\partial_jh_{jk}W^1_k-{\overline{\Gamma}}_{mjj}W^1_m=\frac{1}{2}\partial_kh_{jj}(\partial_a h_{ka}-\partial_k h_{aa})=\frac{1}{2}\partial_k(2\tau-F+\underline{F})\cdot(\in_{klm}\partial_l\Omega_m-2\partial_k\tau).
\end{equation*}
Using now \eqref{ADM43} and \eqref{ADM40} we have
\begin{equation}\label{Bond11}
\partial_jh_{jk}W^1_k-{\overline{\Gamma}}_{mjj}W^1_m\sim 0.
\end{equation}
Moreover, since $\widetilde{O}^2_j:=O^2_j-h_{jk}W^1_k$ is quadratic, 
\begin{equation*}
\partial_j\{O^2_j-h_{jk}W^1_k\}\sim\overg^{kj}\overline{D}_k\widetilde{O}^2_j.
\end{equation*}
Using the divergence theorem, the formulas \eqref{Bond8}, \eqref{Bond10}, and the proof of Proposition \ref{ADM4}, we have
\begin{equation}\label{Bond15}
\begin{split}
\int_{S^u_{R,t}} W^1_j\mathbf{n}^j\,d\sigma&=\int_{B^u_{R,t}} \big\{[(\partial_0\psi)^2+\psi^2+\partial_j\psi\partial_j\psi](t)+(1/2)\sum_{m,n\in\{1,2,3\}}|\partial_0\vartheta_{mn}(t)|^2\big\}\,d\mu\\
&+\int_{S^u_{R,t}} \widetilde{O}^2_j\mathbf{n}^j\,d\sigma+O(\varep_0^2\langle t\rangle^{-\kappa}).
\end{split}
\end{equation}

{\bf{Step 2.}} We examine now the term $\widetilde{O}^2_j\mathbf{n}^j$ in the surface integral above. This is a quadratic term, thus generically bounded by $C\varep_0\langle t\rangle^{-2+}$, so its integral does not vanish as $t\to\infty$. However, many of its pieces have additional structure (such as good vector-fields), therefore satisfy slightly better estimates and do not contribute in the limit as $t\to\infty$. 

More precisely, in view of \eqref{ADM40.5}, \eqref{zaq11.2}, and \eqref{zaq11}, if $|x|\in[t/8,8t]$ then
\begin{equation}\label{Bond16}
\begin{split}
|R(|\nabla|\rho-\partial_0\underline{F})(x,t)|&+|R(\partial_0\rho+|\nabla|\underline{F})(x,t)|+|R(|\nabla|\Omega_j-\partial_0\omega_j)(x,t)|\\
&+|R(\partial_0\Omega_j+|\nabla|\omega_j)(x,t)|+|R\partial_\mu\tau(x,t)|\lesssim \varep_0t^{-5/4}
\end{split}
\end{equation}
for any $i\in\{1,2,3\}$, $\mu\in\{0,1,2,3\}$, and any compounded Riesz transform $R=R_1^{a_1}R_2^{a_2}R_3^{a_3}$, $a_1+a_2+a_3\leq 6$. Moreover, with $\widehat{x}_j:=x^j/|x|$, $j\in\{1,2,3\}$, we have
\begin{equation}\label{Bond17}
\begin{split}
t^{-1}|(\mathbf{n}^j-\widehat{x}_j)(x,t)|+|(\partial_0+\partial_r)RG(x,t)|+|(\partial_j+\widehat{x}_j\partial_0)RG(x,t)|\lesssim\varep_0t^{-5/4},
\end{split}
\end{equation}
for any $j\in\{1,2,3\}$ and $G\in\{F,\underline{F},\rho,\omega_a,\Omega_a,\vartheta_{ab}\}$, as a consequence of \eqref{ADM23.5}, \eqref{ADM40.5}, and Lemma \ref{AlmostOpticalLem}. 

We regroup now the terms in the formula \eqref{ADM21.5}, with respect to the nondifferentiated metric component. Using \eqref{zaq5l} and \eqref{Bond16}-\eqref{Bond17} we write
\begin{equation*}
\begin{split}
\mathbf{n}^j\{h_{00}\partial_0h_{0j}-h_{00}\partial_j(\tau+\underline{F})\}&\simeq h_{00}\widehat{x}_j(-R_j\partial_0\rho-\partial_j\underline{F}+\in_{jkl}\partial_0R_k\omega_l)\\
&\simeq h_{00}\widehat{x}_j\in_{jkl}|\nabla|R_k\Omega_l\simeq 0,
\end{split}
\end{equation*}
where in this section $f\simeq g$ means $|(f-g)(x,t)|\lesssim\varep_0^2 t^{-2-\kappa}$ for all $(x,t)$ with $t\geq 1$ and $|x|\in[t/4,4t]$. The other terms in $\mathbf{n}^j(O_j^2-h_{jk}W^1_k)$ can be simplified in a similar way, 
\begin{equation*}
\begin{split}
\mathbf{n}^j\{-h_{0j}\partial_0h_{00}+2h_{0j}\partial_0(\tau+\underline{F})\}&\simeq h_{0j}\widehat{x}_j\partial_0(-F+\underline{F})\simeq h_{0j}\partial_j(F-\underline{F}),
\end{split}
\end{equation*}
\begin{equation*}
\begin{split}
\mathbf{n}^j\{-h_{0n}&\partial_0h_{nj}+h_{0n}\partial_jh_{0n}-h_{0n}\partial_nh_{0j}\}\simeq h_{0n}\big\{\partial_jh_{nj}+\widehat{x}_j\in_{nab}\partial_jR_a\omega_b-\widehat{x}_j\in_{jab}\partial_nR_a\omega_b\big\}\\
&\simeq h_{0n}\big\{\partial_n(\uF-F)+\in_{nab}\partial_a\Omega_b-\in_{nab}\partial_0R_a\omega_b\big\}\simeq h_{0n}\partial_n(\uF-F),
\end{split}
\end{equation*}
\begin{equation*}
\begin{split}
\mathbf{n}^j\{h_{nj}&\partial_0h_{0n}-h_{nj}\partial_n(\tau+\uF)-h_{nj}W^1_n\}\simeq h_{nj}\big\{-\partial_jh_{0n}-\widehat{x}_j\partial_n\uF+\widehat{x}_j(\partial_nh_{aa}-\partial_ah_{an})\big\}\\
&\simeq h_{nj}\big\{\partial_jR_n\rho-\widehat{x}_j\partial_n\uF-\in_{nab}\partial_jR_a\omega_b-\widehat{x}_j\in_{nab}\partial_a\Omega_b\big\}\\
&\simeq h_{nj}\big\{-\widehat{x}_jR_n(\partial_0\rho+|\nabla|\uF)+\widehat{x}_j\in_{nab}R_a(\partial_0\omega_b-|\nabla|\Omega_b)\big\}\simeq 0.
\end{split}
\end{equation*}
Summing up these identities and using \eqref{zaq22} we get
\begin{equation}\label{Bond18}
\mathbf{n}^j\widetilde{O}^2_j\simeq \mathbf{n}^jh_{kn}(\partial_nh_{kj}-\partial_jh_{kn})\simeq -\mathbf{n}^j\overg_{\geq 1}^{kn}(\partial_nh_{kj}-\partial_jh_{kn}).
\end{equation}
Therefore, $\mathbf{n}^j\widetilde{O}^2_j+\mathbf{n}^jW^{\geq 2}_j\simeq 0$, and \eqref{Bond15} gives
\begin{equation}\label{Bond20}
\begin{split}
\int_{S^u_{R,t}} W_j\mathbf{n}^j\,d\sigma=\int_{B^u_{R,t}} \big\{[&(\partial_0\psi)^2+\psi^2+\partial_j\psi\partial_j\psi](t)\\
&+(1/2)\sum_{m,n\in\{1,2,3\}}|\partial_0\vartheta_{mn}(t)|^2\big\}\,d\mu+O(\varep_0^2\langle t\rangle^{-\kappa}).
\end{split}
\end{equation}

{\bf{Step 3.}} We fix now $R\in\mathbb{R}$ and let $t\to\infty$. At this stage, for the limit to exist it is important that the almost optical function $u$ has the properties stated in Lemma \ref{AlmostOpticalLem}. 

Recall the scattering profiles $V^\psi_\infty$ and $V^{\va_{ab}}_\infty$ defined in \eqref{Alop10}--\eqref{Alop10.1}. We show first that, for any $R\in\mathbb{R}$,
\begin{equation}\label{Bond21}
\lim_{t\to\infty}\int_{B^u_{R,t}} [(\partial_0\psi)^2+\psi^2+\partial_j\psi\partial_j\psi](x,t)\,d\mu=\|V^\psi_\infty\|_{L^2}^2.
\end{equation}

Indeed, we notice that
\begin{equation}\label{Bond22}
\int_{\mathbb{R}^3} [(\partial_0\psi)^2+\psi^2+\partial_j\psi\partial_j\psi](x,t)\,d x=\big\|U^\psi(t)\|_{L^2}^2=\|V_\ast^\psi(t)\|_{L^2}^2,
\end{equation}
and $d\mu=dx(1+O(\varep_0\langle t\rangle ^{-1+2\delta'}))$. Therefore, for \eqref{Bond21} it suffices to prove that
\begin{equation}\label{Bond23}
\lim_{t\to\infty}\int_{\mathbb{R}^3\setminus B^u_{R,t}} [(\partial_0\psi)^2+\psi^2+\partial_j\psi\partial_j\psi](x,t)\,dx=0.
\end{equation}
Recalling that $U^\psi(t)=\partial_0\psi(t)-i\Lambda_{kg}\psi(t)$ and $|u(x,t)-|x|+t|\lesssim\varep_0\langle x\rangle^{3\delta'}$ (see Lemma \ref{AlmostOpticalLem}), for \eqref{Bond23} it suffices to prove that, for operators $A\in\{I,\langle\nabla\rangle^{-1},\partial_j\langle\nabla\rangle^{-1}\}$ we have
\begin{equation}\label{Bond24}
\lim_{t\to\infty}\int_{|x|\geq t-t^{1/2}} |AU^\psi(x,t)|^2\,dx=0.
\end{equation}

This follows like in the proof of Proposition \ref{Alop20}. Indeed, with $k_0,J_0$ defined as the smallest integers for which $2^{k_0}\geq t^{\delta'}$ and $2^{J_0}\geq t^{1/3}$,  we have
\begin{equation*}
\big\|A(I-P_{[-k_0,k_0]})U^\psi(t)\big\|_{L^2}+\big\|AP_{[-k_0-2,k_0+2]}(\varphi_{\geq J_0+1}\cdot P_{[-k_0,k_0]}U^\psi)(t)\big\|_{L^2}\lesssim\varep_0t^{-\delta},
\end{equation*}
due to \eqref{winr23}. Moreover, the remaining component $AP_{[-k_0-2,k_0+2]}(\varphi_{\leq J_0}\cdot P_{[-k_0,k_0]}U^\psi)(t)$ is rapidly decreasing in the region $\{|x|\geq t-t^{1/2}\}$, so the desired limit \eqref{Bond23} follows.

{\bf{Step 4.}} To calculate the contribution of the metric components $|\partial_0\vartheta_{ab}|^2$ we use Proposition \ref{Alop20}. For $v\in\mathbb{R}$, $\theta\in\mathbb{S}^2$, and $t\geq 10$ we define
\begin{equation}\label{Bond31}
L_{ab}(v,\theta,t):=\Re\Big\{\frac{-i}{4\pi^2}\int_0^\infty e^{i\rho v}\varphi_{[-k_0,k_0]}(\rho)\widehat{V^{\vartheta_{ab}}_\ast}(\rho\theta,t)\rho\,d\rho\Big\}.
\end{equation}

We show first that, for any $R\in\mathbb{R}$, $a,b\in\{1,2,3\}$, and $t\geq (2+|R|)^{10}$  we have
\begin{equation}\label{Bond32}
\int_{B^u_{R,t}} |\partial_0\vartheta_{ab}(x,t)|^2\,d\mu=\int_{[-t^{0.4}/8,R]\times \mathbb{S}^2}|L_{ab}(v,\theta,t)|^2\,dv d\theta+O(\varep_0^2t^{-\delta}).
\end{equation}
Indeed, we can first replace the measure $d\mu$ by $dx$, at the expense of an acceptable error. Then, using \eqref{Alop24} and \eqref{Alop28}, we may assume that the integration is over the domain $D^u_{[-t^{0.4}/8,R],t}$ where $D^u_{[A,B],t}:=\{x\in\Sigma_t:\,u(x,t)\in[A,B]\}$, since the integration in the interior region $\{|x|\leq t-t^{0.4}/20\}$ of $|\partial_0\vartheta_{ab}(t)|^2$ produces an acceptable error. In addition, we may replace $\partial_0\vartheta_{ab}(x,t)$ with $|x|^{-1}L_{ab}(u(x,t),x/|x|,t)$ at the expense of acceptable errors, due to the first identity in \eqref{Alop21}. To summarize,
\begin{equation}\label{Bond33}
\int_{B^u_{R,t}} |\partial_0\vartheta_{ab}(x,t)|^2\,d\mu=\int_{D^u_{[-t^{0.4}/8,R],t}}|x|^{-2}|L_{ab}(u(x,t),x/|x|,t)|^2\,dx+O(\varep_0^2t^{-\delta}).
\end{equation}
We pass now to polar coordinates $x=r\theta$, and then make the change of variables $r\to v:=u(r\theta,t)$. The desired approximate identity \eqref{Bond32} follows using also \eqref{EstimuGradu}.

We apply now \eqref{Alop10.1} to conclude that
\begin{equation}\label{Bond34}
\lim_{t\to\infty}\int_{B^u_{R,t}} |\partial_0\vartheta_{ab}(x,t)|^2\,d\mu=\int_{(-\infty,R]\times \mathbb{S}^2}|L^\infty_{ab}(v,\theta)|^2\,dv d\theta
\end{equation}
where $L^\infty_{ab}(v,\theta)=\lim_{t\to\infty}L_{ab}(v,\theta,t)$ in $L^2(\mathbb{R}\times\mathbb{S}^2)$ is given by
\begin{equation}\label{Bond31.6}
L^\infty_{ab}(v,\theta):=\Re\Big\{\frac{-i}{4\pi^2}\int_0^\infty e^{i\rho v}\widehat{V^{\vartheta_{ab}}_\infty}(\rho\theta)\rho\,d\rho\Big\}.
\end{equation}
Combining this with \eqref{Bond20} and \eqref{Bond21}, we have
\begin{equation*}
\begin{split}
\lim_{t\to\infty}\int_{S^u_{R,t}} W_j\mathbf{n}^j\,d\sigma=\|V^\psi_\infty\|_{L^2}^2+\frac{1}{2}\sum_{m,n\in\{1,2,3\}}\int_{(-\infty,R]\times \mathbb{S}^2}|L^\infty_{mn}(v,\theta)|^2\,dv d\theta.
\end{split}
\end{equation*}
Recalling the definitions \eqref{Bond5} we have
\begin{equation}\label{Bond36}
E_{Bondi}(R)=\frac{1}{16\pi}\|V^\psi_\infty\|_{L^2}^2+\frac{1}{32\pi}\sum_{m,n\in\{1,2,3\}}\int_{(-\infty,R]\times \mathbb{S}^2}|L^\infty_{mn}(v,\theta)|^2\,dv d\theta,
\end{equation}
which is clearly well-defined, continuous, and increasing on $\mathbb{R}$. The limit as $R\to -\infty$ in \eqref{Bond7.1} follows since $L^\infty_{mn}\in L^2(\mathbb{R}\times\mathbb{S}^2)$. To prove the limit as $R\to \infty$ we use \eqref{ADM4.6} and let $t\to\infty$. Clearly $\lim_{t\to\infty}\|U^\psi(t)\|_{L^2}^2=\|V^\psi_\infty\|_{L^2}^2$ and, using \eqref{ADM43}, $\|U^{\va_{mn}}(t)\|_{L^2}^2=2\|\partial_0\va_{mn}(t)\|_{L^2}^2$. The desired limit as $R\to\infty$ in \eqref{Bond7.1} follows using also \eqref{Bond32}. 
\end{proof}

\subsection{The interior energy} We see that the total Klein-Gordon energy $E_{KG}$ defined in \eqref{Bond7.1} is part of the null Bondi energy $E_{Bondi}(R)$, for all $R\in\mathbb{R}$. This is consistent with the geometric intuition, because the matter travels at speeds lower than the speed of light, and accumulates at the future timelike infinity, not at null infinity. We show now that this Klein-Gordon energy can be further radiated by taking limits along suitable time-like cones.

\begin{proposition}\label{Bond40}
For $\alpha\in (0,1)$ let 
\begin{equation}\label{Bond41}
E_{i^+}(\alpha):=\frac{1}{16\pi}\lim_{t\to\infty}\int_{S_{\alpha t,t}}(\partial_j h_{nj}-\partial_n h_{jj})\frac{x^n}{|x|}\,dx,
\end{equation}
where the integration is over the Euclidean spheres $S_{\alpha t,t}\subseteq\Sigma_t$ of radius $\alpha t$. Then the limit in \eqref{Bond41} exists, and $E_{i^+}:(0,1)\to\mathbb{R}$ is a well-defined continuous and increasing function. Moreover, we have
\begin{equation}\label{Bond42}
\lim_{\alpha\to 0}E_{i^+}(\alpha)=0,\qquad \lim_{\alpha\to 1}E_{i^+}(\alpha)=E_{KG}.
\end{equation}
\end{proposition}

\begin{proof} We notice that the definition \eqref{Bond41} is similar to the definition of the ADM energy in \eqref{ADM1}. We could also use a more ``geometric" definition involving the vector-field $W$ (see \eqref{Bond7}), but this would make no difference here as $t\to\infty$, because the integrand $|\partial_j h_{nj}-\partial_n h_{jj}|$ is already bounded by $C_\alpha \langle t\rangle^{-2+4\delta'}$. 

{\bf{Step 1.}} We may assume that $t$ is large, say $t\geq \alpha^{-10}+(1-\alpha)^{-10}$. Using Stokes theorem and the definitions \eqref{zaq21.1}, we can rewrite 
\begin{equation}\label{Bond43}
\int_{S_{\alpha t,t}}(\partial_j h_{nj}-\partial_n h_{jj})\frac{x^n}{|x|}\,dx=\int_{|x|\leq \alpha t}-2\Delta\tau(x,t)\,dx.
\end{equation}
The density function $-2\Delta\tau$ was analyzed in Lemma \ref{ADM20}. The contributions the error terms $O^1$ and $\partial_jO^2_j$ decay as $t\to\infty$, and the contribution of the metric components $-\delta_{jk}P_{jk}^2$ also decays because $|\partial_\be h_{\mu\nu}(x,t)|\lesssim_\alpha\varep_0 \langle t\rangle^{-2+2\delta'}$ in the ball $\{|x|\leq \alpha t\}$ (due to \eqref{ImprovedBoundsInside}). Thus
\begin{equation}\label{Bond44}
\int_{|x|\leq \alpha t}-2\Delta\tau(x,t)\,dx=\int_{|x|\leq \alpha t}\{(\partial_0\psi)^2+\psi^2+\partial_j\psi\partial_j\psi\}(x,t)\,dx+O_\alpha(\varep_0^2t^{-\kappa}).
\end{equation}

{\bf{Step 2.}} To apply Proposition \ref{Alop20} we would like to show now that
\begin{equation}\label{Bond45}
\int_{|x|\leq \alpha t}\{(\partial_0\psi)^2+\psi^2+\partial_j\psi\partial_j\psi\}(x,t)\,dx=\int_{|x|\leq \alpha t}|U^\psi(x,t)|^2\,dx+O_\alpha(\varep_0^2t^{-\kappa}).
\end{equation}
Indeed, the real part of $U^\psi$ is $\partial_0\psi$, as needed. The imaginary part of $U^\psi$ is $-\langle \nabla\rangle\psi$, and some care is needed because the functions $\psi$ and $\partial_j\psi$ are connected to $\langle \nabla\rangle\psi$ by nonlocal operators. 

In view of \eqref{winr25.4}, we may replace the integral over the ball $\{|x|\leq \alpha t\}$ with a suitably smooth version, using the function $\chi_1\big(t^{-0.9}(|x|-\alpha t)\big)$, where $\chi_1:\mathbb{R}\to[0,1]$ is a smooth function supported in $(-\infty,2]$ and equal to $1$ in $(-\infty,1]$. For \eqref{Bond45} it suffices to prove that
\begin{equation}\label{Bond46}
\int_{\mathbb{R}^3}\chi_1\big(t^{-0.9}(|x|-\alpha t)\big)\{\psi^2+\partial_j\psi\partial_j\psi-(\langle \nabla\rangle\psi)^2\}(x,t)\,dx=O_\alpha(\varep_0^2t^{-\kappa}).
\end{equation}
To prove this we notice that
\begin{equation}\label{Bond46.1}
\int_{\mathbb{R}^3}\chi_1\big(t^{-0.9}(|x|-\alpha t)\big)G(x)\overline{G(x)}\,dx=C\int_{\mathbb{R}^3\times\mathbb{R}^3}\widehat{G}(\xi)\overline{\widehat{G}(\eta)} K_{1,t}(\xi-\eta)\,d\xi d\eta,
\end{equation}
where
\begin{equation}\label{Bond46.2}
K_{1,t}(\rho):=\int_{\mathbb{R}^3}\chi_1\big(t^{-0.9}(|x|-\alpha t)\big)e^{ix\cdot\rho}\,dx.
\end{equation}
We apply the identity \eqref{Bond46.1} for $G\in\{\langle \nabla\rangle\psi,\psi,\partial_j\psi\}$. For \eqref{Bond46} it suffices to prove that 
\begin{equation}\label{Bond46.3}
\Big|\int_{\mathbb{R}^3\times\mathbb{R}^3}\widehat{\psi}(\xi,t)\overline{\widehat{\psi}(\eta,t)} \big(\langle\xi\rangle\langle\eta\rangle-1-\xi_j\eta_j\big)K_{1,t}(\xi-\eta)\,d\xi d\eta\Big|\lesssim_{\alpha}\varep_0^2t^{-\kappa}.
\end{equation}
Using integration by parts it is easy to see that $|K_{1,t}(\rho)|\lesssim_\al t^3(1+t^{0.9}|\rho|)^{-10}$ (the rapid decay here is the main reason for replacing the characteristic function of the ball $\{|x|\leq \alpha t\}$ by the smooth approximation $\chi_1(t^{-0.9}(|x|-\alpha t))$). Moreover $\big|\langle\xi\rangle\langle\eta\rangle-1-\xi_j\eta_j\big|\lesssim\langle\xi\rangle\langle\eta\rangle|\xi-\eta|$, thus
\begin{equation*}
\big|\langle\xi\rangle\langle\eta\rangle-1-\xi_j\eta_j\big|\,|K_{1,t}(\xi-\eta)|\lesssim_\al\langle\xi\rangle\langle\eta\rangle t^{2.1}(1+t^{0.9}|\xi-\eta|)^{-9}.
\end{equation*}
Therefore, the left-hand side of \eqref{Bond46.3} is bounded by
\begin{equation*}
C_\alpha\|\widehat\psi(\xi,t)\langle\xi\rangle\|_{L^2}^2t^{-0.5}\lesssim_\al\varep_0^2t^{-0.5},
\end{equation*}
as desired. This completes the proof of \eqref{Bond45}.

{\bf{Step 3.}} We use now \eqref{Alop21}, pass to polar coordinates, and change variables to calculate
\begin{equation*}
\begin{split}
\int_{|x|\leq \alpha t}|U^\psi(x,t)|^2\,dx&=\frac{1}{8\pi^3}\int_{[0,\alpha t]\times\mathbb{S}^2}\frac{t^2r^2}{(t^2-r^2)^{5/2}}\Big|\widehat{P_{[-k_0,k_0]}V^\psi_\ast}\Big(\frac{r\theta}{\sqrt{t^2-r^2}},t\Big)\Big|^2\,dr d\theta+O_\alpha(\varep_0^2t^{-\delta})\\
&=\frac{1}{8\pi^3}\int_{[0,\alpha/\sqrt{1-\alpha^2}]\times\mathbb{S}^2}\rho^2\big|\widehat{P_{[-k_0,k_0]}V^\psi_\ast}\big(\rho\theta,t\big)\big|^2\,d\rho d\theta+O_\alpha(\varep_0^2t^{-\delta}).
\end{split}
\end{equation*}
Therefore, using also \eqref{Bond43}--\eqref{Bond45} and letting $t\to\infty$, we have
\begin{equation}\label{Bond47}
E_{i^+}(\alpha)=\frac{1}{16\pi}\frac{1}{8\pi^3}\int_{[0,\alpha/\sqrt{1-\alpha^2}]\times\mathbb{S}^2}\rho^2\big|\widehat{V^\psi_\infty}\big(\rho\theta\big)\big|^2\,d\rho d\theta,
\end{equation}
which is clearly a well-defined continuous and increasing function of $\alpha$ satisfying \eqref{Bond42} (recall that $E_{KG}=1/(16\pi)\|V^\psi_\infty\|_{L^2}^2=1/(16\pi)(8\pi^3)^{-1}\|\widehat{V^\psi_\infty}\|_{L^2}^2$).
\end{proof}

\end{document}